\documentclass[11pt,a4paper]{article}
\usepackage{authblk}
\usepackage[T1]{fontenc}
\usepackage[utf8]{inputenc}
\usepackage[british]{babel}
\usepackage{lmodern}
\usepackage[a4paper,margin=29mm]{geometry}
\usepackage{amsmath,amssymb,amsthm,mathtools}
\usepackage{booktabs,tabularx,array}
\usepackage{enumitem}
\usepackage{graphicx,float,caption}
\usepackage{microtype}
\usepackage{xcolor}
\usepackage[colorlinks=true,linkcolor=blue,citecolor=blue,urlcolor=blue]{hyperref}
\hypersetup{pdftitle={The Gleason--Kahane--Zelazko theorem: a gentle introduction},pdfauthor={Catalin Badea}}
\numberwithin{equation}{section}
\newtheorem{theorem}{Theorem}[section]
\newtheorem{proposition}[theorem]{Proposition}
\newtheorem{lemma}[theorem]{Lemma}
\newtheorem{corollary}[theorem]{Corollary}
\theoremstyle{definition}

\newtheorem{example}[theorem]{Example}
\theoremstyle{remark}
\newtheorem{remark}[theorem]{Remark}
\newcommand{\C}{\mathbb C}
\newcommand{\R}{\mathbb R}

\newcommand{\one}{\mathbf 1}
\newcommand{\inv}{A^{-1}}
\newcommand{\spA}{\sigma_A}
\newcommand{\GKZ}{Gleason--Kahane--\.{Z}elazko}
\newcommand{\norm}[1]{\left\lVert#1\right\rVert}

\DeclareMathOperator{\Rea}{Re}

\DeclareMathOperator{\image}{image}
\newcolumntype{Y}{>{\raggedright\arraybackslash}X}
\setlist[enumerate]{itemsep=3pt,topsep=5pt}
\title{The Gleason--Kahane--\.{Z}elazko theorem: a gentle introduction}
\author{Catalin Badea}
\affil{%
Department of Mathematics and Statistics, University of Reading, United Kingdom,\\
and UMR 8524 -- Laboratoire Paul Painlev\'e, Universit\'e de Lille, France\\
\texttt{c.badea@reading.ac.uk}, \texttt{cbadea@univ-lille.fr}\\
\url{https://www.reading.ac.uk/maths-and-stats/staff/catalin-badea}%
}

\date{}

\ifdefined\pdfobj\else
  \PackageError{gkz-pedago-15}{Please compile this self-contained file with pdfLaTeX}
    {The embedded dedication uses native pdfTeX image objects.}
\fi
\begingroup
\immediate\pdfobj stream attr{
  /Type /XObject /Subtype /Image
  /Width 1586 /Height 2264
  /ColorSpace /DeviceRGB /BitsPerComponent 8
  /Filter [/ASCIIHexDecode /DCTDecode]
}{
FFD8FFE000104A46494600010100000100010000FFDB0043000D090A0B0A080D0B0A0B0E0E0D0F13201513121213271C1E17202E2931302E292D2C33
3A4A3E333646372C2D405741464C4E525352323E5A615A50604A51524FFFDB0043010E0E0E131113261515264F352D354F4F4F4F4F4F4F4F4F4F4F4F
4F4F4F4F4F4F4F4F4F4F4F4F4F4F4F4F4F4F4F4F4F4F4F4F4F4F4F4F4F4F4F4F4F4F4F4F4F4FFFC000110808D8063203012200021101031101FFC400
1F0000010501010101010100000000000000000102030405060708090A0BFFC400B5100002010303020403050504040000017D010203000411051221
31410613516107227114328191A1082342B1C11552D1F02433627282090A161718191A25262728292A3435363738393A434445464748494A53545556
5758595A636465666768696A737475767778797A838485868788898A92939495969798999AA2A3A4A5A6A7A8A9AAB2B3B4B5B6B7B8B9BAC2C3C4C5C6
C7C8C9CAD2D3D4D5D6D7D8D9DAE1E2E3E4E5E6E7E8E9EAF1F2F3F4F5F6F7F8F9FAFFC4001F0100030101010101010101010000000000000102030405
060708090A0BFFC400B51100020102040403040705040400010277000102031104052131061241510761711322328108144291A1B1C109233352F015
6272D10A162434E125F11718191A262728292A35363738393A434445464748494A535455565758595A636465666768696A737475767778797A828384
85868788898A92939495969798999AA2A3A4A5A6A7A8A9AAB2B3B4B5B6B7B8B9BAC2C3C4C5C6C7C8C9CAD2D3D4D5D6D7D8D9DAE2E3E4E5E6E7E8E9EA
F2F3F4F5F6F7F8F9FAFFDA000C03010002110311003F00F4EA28A2800A28A2800A28A2800A28A2800A28A2800A28A2800A28A2800A28A2800A28A280
0A28A2800A28A2800A28A2800A28A2800A28A2800A28A2800A28A2800A28A2800A28A2800A28A2800A28A2800A28A2800A28A2800A28A2800A400280
00000E0014B4500148402304034B450014514500145145001451450014514500145145002606738191C66968A2800A400018000A5A2800A28A2800A4
C0CE7033EB4B45001451450014840241201C723DA968A0028A28A0028A28A0028A28A0028A28A0028A28A0028A28A0028A28A0028A28A0028A28A002
8A28A0028A28A0028A28A0028A28A0028A28A0028A28A0028A28A0028A28A0028A28A0028A28A0028A28A0028A28A0028A28A0028A28A0028A28A002
8A28A00AD2595B4A4B3C2A49EA4719FCAA748D6340A83007414EA295905D8514514C028A28A0028A28A0028A28A0028A28A0028A28A0028A28A0028A
28A0028A28A0028A28A0028A28A0028A28A0028A28A0028A28A0028A28A0028A28A0028A28A004000248039EBEF41008C100FD6968A0028A28A0028A
28A0028A28A0028A28A0028A28A0028A28A0028A28A0028A28A0028A28A0028A28A0028A28A0028A28A0028A28A0028A28A0028A28A0028A28A0028A
28A0028A28A0028A28A006ED5FEE8FCA8A75140051451400514514005145140051451400514514005145140051451400514514005145140051451400
514514005145140051451400514514005145140051451400514514005145140051451400514514005145140051451400514514005145140051451400
514514005145140051451400514514005145140051451400514514005145140051451400514514005145140051451400514514005145140051451400
514514005145140051451400514514005145140051451400514514005145140051451400514514005145140051451400514514005145140051451400
514514005145140051451400514514005145140051451400514514005145140051451400514514005145140051451400514514005145140051451400
514514005145140051451400514514005145140051451400514514005145140051451400514514005145140051451400514514005145140051451400
514514005145140051451400514514005145140051451400514514005145140051451400514514005145140051451400514514005145140051451400
514514005145140051451400514514005145140051451400514514005145140051451400514514005145140051451400514514005145140051451400
514514005145140051451400514514005145140051451400514514005145140051451400514514005145140051451400514514005145140051451400
514514005145140051451400514514005145140051451400514514005145140051451400514514005145140051451400514514005145140051451400
514514005145140051451400514514005145140051451400514514005145140051451400514514005145140051451400514514005145140051451400
514514005145729E32D7E5B08D2C6C08373371C7256802EEA9E2382CAE1618B6CAC0E1FDAB6A29165896443956008AE0AE7479ACBC37F6BBB65F3D8E
589EBC8E2BB3D1893A3D996183E52E7F2A486D17A8A28A620A28A2800A28AC6D435C8ADDDA1B54F3A65FFBE45006CD158DA1EA73DF34A970A372FCC0
8E31ED5B34005145140051484800924003A934C8A549A31244C191BA11DE8024A28A2800A28AA5A8EA56DA6DB996E640303217B9A00BB45721A6EA7A
86AFAEC60298ADE1F9CE3B8F4AEBE9260145145300A28A2800A28A2800A28A2800A28A648DB2367C6768271400ACCAAA5988007726ABCB7F6B16774A
A7033F2F35C4EA3A9EA1A8EB70D8A975477030BC015D545A65AD95A3B4C77E172598E314AE558B96B7F6F74EC90B12CA33C8ED566B9CF0C079EE2E6E
C28580FC91FBF3CD7494D09E8145145020A28A2800A28AA11EA90CBA91B288166504B37614017E8A28A0028A28A0028A28A0028A28A0028A28A0028A
28A0028A28A0028A28A0028A28A0028A28A0028A28A0028A28A0028A28A0028A28A0028A28A0028A28A0028A28A0028A28A0028A28A0028A28A0028A
28A0028A28A0028A28A0028A28A0028A28A0028A28A0028A28A0028A28A0028A28A0028A28A0028A28A0028A28A0028A28A0028A28A0028A28A0028A
28A0028A28A0028A28A0028A28A0028A28A0028A28A0028A28A0028A28A0028A28A0028A28A0028A28A0028A28A0028A28A0028A28A0028A28A0028A
28A0028A28A0028A28A0028A28A0028A28A0028A28A0028A28A0028A28A0028A28A0028A28A0028A28A0028A28A0028A28A0028A28A0028A28A0028A
28A0028A28A0028A28A0028A28A0028A28A0028A28A0028A28A0028A28A0028A28A0028A28A0028A28A0028A28A0028A28A0028A28A00CED6B548749
D3A4B9948DC0108BDD8D725E13D267D4F557D735142431CAABF383506BB72DE21F13A69F6E77C103856F407EB5DF5A5B47696C9044A02A0C7D7DE96E
3D9181E31C4D1DA5A06F9E5900C7E23FC0D7416F1082DA28874440BF90AE7F506F3BC676307FCF35DDF90CD74B403D85A28A298828A2AB5FDCFD92CA
5B8DBBB62E71EB40191AEEAAEB32E9B64774F270C47F0E7B55AD27468ACA3DF3289276FBC4F20551F0C69DB95B54B9CB4D392467B73D6BA4A432086D
60824778A30ACE79C54F55EEEEE1B387CC99B03B0EE6B147886667564B61E516E49CE714C2CD9D151581A66BD26A3AC496B15BEDB78D725FA926B7E8
158E6BC6FAAB699A1486319925181CF207152782AF3ED5A0C7190DBA03B093DFBFF5AC1F1079BAF789D34A41FBB8DBE638E807FF005F35D49FB3787F
4D11C2ACC09240F5353D4AE9635A8AE6B4FD5EF13CC9B5053E4904AF0060D68689AA3EA893C86211A23ED5F534EE2B17EE674B5B792794E11064D707
35D3EAF3FDA6E0E50C9848C0EDF4ADBF19DEAA592D927CEF21F9D57A818E299E13D105BC105E4DC91181183D7EA693D58D68AE6B68DA7FD8E269241F
BD90F3EC3B0AD4ACDD63521A7DB8D8374AFC28F4AE6F51D6F574895D1C21DA4E0629DEC166CEDA8AC2D0354967D31A7BF6DBB08F998609CFF3A7AEB8
B25F4712427C96EAEDC1A62B1B545645D6B7124862B6432BFAF6152E957EF76AC9301E6273903191408D2A28A2800A28A2800A6B6369DD8DB8E73D31
4EACCD76ECDAE9EC23602497E5507F5A00C8D02D927D6EEAF231FBA8FE55DDCF26A5D7AEDEEEEA3D2ED0124B62438E056868B64F69A388F204D202C5
BDCF4AA49A1DCC2EF34732F9A790431CD22AFA9B3636A9676B1C1180028E71DCD59AC1F0F5E5C5C4934370ECFB39CB72473EB5A1AADF7D86D77A8CBB
1C28A6265EA2B80D4353D52DEE16E2498E472003D8D5FF00097886E755BD922B9C9E091D87F9EB4AE163B0A4240C64F5ACFD5AFCD9C004454CCE70A0
F6F7C560DD5BDEADB8BEBA9D892F9C0E093EC3D38A2E091B1E22D5534BB02C5B123F0B83C8F7AADE14D31ACECDEE67E67B93BB3DF6F515CF69D6F75A
F6BC5AF4EE8A2C3B0CE460741F5AEDAF6EE1D3ED0CB2630A3E551C67D85037A2B16E8AE0E2D6F57B9BFDE19923671851D07B62BBCA626AC14CDE9BF6
6F5DDE99E6B1F54BE9E5BB1A7590F9D87CEE3B7B552BEB6974B8E39CCA1E42DC72680B1D45158DE1ED46E75182792E57015F6AF18E2B66813560A28A
2800A28A2800A28A2800A28A2800A28A2800A28A2800A28A2800A28A2800A28A2800A28A2800A28A2800A28A2800A28A2800A28A2800A28A2800A28A
2800A28A2800A28A2800A28A2800A28A2800A28A2800A28A2800A28A2800A28A2800A28A2800A28A2800A28A2800A28A2800A28A2800A28A2800A28A
2800A28A2800A28A2800A28A2800A28A2800A28A2800A28A2800A28A2800A28A2800A28A2800A28A2800A28A2800A28A2800A28A2800A28A2800A28A
2800A28A2800A28A2800A28A2800A28A2800A28A2800A28A2800A28A2800A28A2800A28A2800A28A2800A28A2800A28A2800A28A2800A28A2800A28A
2800A28A2800A28A2800A28A2800A28A2800A28A2800A28A2800A28A2800A28A2800A28A2800A28A2800A28A2800A28A2800A28A2800A28A2800AC4F
155F4965A4B0809134C762E3F5ADBAA77BA75BDECB13CE32633903D68031FC25A1A69D64679635F3E63BB27923DFEB5A4935EDB447CF87CEC138656E
D5A200000030076A5A077381D3E6B8BAF88114D2A303E5B12B9E146319AEFAB97B1888F1ADC4847063600FE39AEA29206145145310557BE84CF67244
A012C3A1EFCD58A2803234CB88ED6D7ECD38F2BCAE070692EB5821B659C0D2B640C9071F956B3A2B8C3AAB0F423342A2A0C2A851EC31401CD6AD677B
34AB713A6E51FC23A0A7A59DE5C5B88E28D562FA75FC4D7494503BB39ED1ADE1D2167FB4029231E091C91FFEBA5BED66696368EC226193B7CC6EDEE2
B6E58629B6F9A8AFB4E4645392348D76C68AA3D146280B9CF785B47366D3DF4EC5A69CE324E78AD0D5AD259DE39572E89D631FCEB528A02E73D2E997
97EE7ED056188A7CAABD075E31F8D16BA6EA167035ADAED8D09E5C30E7DFD6BA1A29582E6447A15A1811670CF28392E18E73FE1572DD6786048E4457
DBF28D879DB8F7C55BA2988E6EFB4ED4AF2F7CF281547DD52E0802A0FF00846E7BA532DCBE19B9DA4F3F4AEAE8A561F3339E8344919F6CC0450AE005
5209356758B006D13ECB07EF15B0362F38F7AD8A2985CCBD1F4B5B08B7C986B87E58FA7B0A73C4D697A678A326261F3051D3FCF5AD2A2815CACB7B6C
D1F99E6A81EFD69F6D7315D45E64272B9C5326B1B6990AB46003FDDE2A48218EDE258A25DAAB4012D14514005733AB3B6A3AEDB5844C42447739C7E7
5D35528B4E862BE92EC64BBE7F0CF5A009AEEEA0B2B769EE64091AF524D605EF8816F14DAE9F1BB3C981B8F071EC2AF6B3A4C9AA15FDF0544E8BC8E6
A4D2746874E52F80F31CE5FD3E940D58768BA7FF0067DA95720CCE72E7F90A9352B26BD540ACA3693C1ABD4502393D674768D124793CCE319C7DD3FF
00EACD4DE14D3ECA3125CC077480ED03FBA3FCE6BA1B8823B888C72AE54D32CECE1B287CB817033927B9A077D0C7B9022D61A6BA43228E5401C63B54
86CAE354983DF662850FCB1E319AD92885C39552C3A1C7229F405CCB874E1A7F9AF60819E4C0DAC7007B9F5AC7BCD3752D4655BABB89BE5E5630D803
F0AEB28A02E6346F656E23F32D99665009017383F8D226A77971398ADAD54647058F2BEE6B670339C73EB4B401836D617F6D24B32006593EF1C839FC
E965D166BCDD25E4E0C847CA07201ADDA280B9CE5B5DB69109B77859C2925B8C56FC1279D0249B4AEE19C1A57449176BAAB0F42334E0000001803B50
0C5AC9D57537B7B986D2D86E9E4E480338146B5AB0B00B144A1EE1FEE8F4A8345D365F305FDEE4CCD92A0F519A406DA6E28A5C00D8E40F5A75145310
514514005145140051451400514514005145140051451400514514005145140051451400514514005145140051451400514514005145140051451400
514514005145140051451400514514005145140051451400514514005145140051451400514514005145140051451400514514005145140051451400
514514005145140051451400514514005145140051451400514514005145140051451400514514005145140051451400514514005145140051451400
514514005145140051451400514514005145140051451400514514005145140051451400514514005145140051451400514514005145140051451400
514514005145140051451400514514005145140051451400514514005145140051451400514514005145140051451400514514005145140051451401
08B7845C19C463CD2305AA6A28A0028A28A0028A28A0028A28A0028A28A0028A28A0028A28A0028A8E6952089A59582A20C926B3B46D663D5CCC618D
8246701BB1A00D5A28A2800A28A2800A28A2800A28A2800A28A2800A28A2800A28A2800A28A2800A28A2800A28A2800A28A2800A28A2800A28A2800A
28A2800A28A280293E9D04B7DF6B946F718DA08E0715768A2800A28A2800A28A2800A28A2800A28A2800A28A2800A28A2800A28A2800A28A2800A28A
2800A28A2800A28A2800A28A2800A28A2800A28A2800A28A2800A28A2800A28A2800A28A2800A28A2800A28A2800A28A2800A28A2800A28A2800A28A
2800A28A2800A28A2800A28A2800A28A2800A28A2800A28A2800A28A2800A28A2800A28A2800A28A2800A28A2800A28A2800A28A2800A28A2800A28A
2800A28A2800A28A2800A28A2800A28A2800A28A2800A28A2800A28A2800A28A2800A28A2800A28A2800A28A2800A28A2800A28A2800A28A2800A28A
2800A28A2800A28A2800A28A2800A28A2800A28A2800A28A2800A28A2800A28A2800A28A2800A28A2800A28A2800A28A2800A28A2800A28A2800A28A
2800A28A2800A28A2800A28A2800A28A2800A28A2800A28A2800A28A2800A28A2800A28A2800A28A2800A28A2800A28A2800A28A2800A8E5952189A5
9582A28C927B5495CE6AD74FA85C9B1B594AC318DD33A9EB8E71F4A00C1D7751BCD66ED61B56DB69BF68ED9AEBF44D3134BB110A80198EE6C7AD66E8
1A7AB4ED76CBFBB4F962C8FBDEF49E23F11A5866D2CD83DC9FBC473B3FFAF4BCCADF44743E6C7E66CDEBBFFBB9E7F2A92B8DF0958DCCB71F6EBB7721
73D4F53FE79AECA989AB05145140828A28A0028A2B07C4FAD8D2ACCA41F35D3FDD51DBDE804AE6A1BFB41742D7CF5F38FF0000EB56AB86F08584973A
9C9A8DD339651FC5EA6BB9A486D5828A28A620A28A2800A28A2800A28A2800A28A2800A28A2800A28A2800A28A2800A28A2800A28A2800A28A2800A2
8A2800A28A2800A28A2800A28A2800A28A2800A28A2800A28A2800A28A2800A28A2800A28A2800A28A2800A28A2800A28A2800A28A2800A28A2800A2
8A2800A28A2800A28A2800A28A2800A28A2800A28A2800A28A2800A28A2800A28A2800A28A2800A28A2800A28A2800A28A2800A28A2800A28A2800A2
8A2800A28A2800A28A2800A28A2800A28A2800A28A2800A28A2800A28A2800A28A2800A28A2800A28A2800A28A2800A28A2800A28A2800A28A2800A2
8A2800A28A2800A28A2800A28A2800A28A2800A28A2800A28A2800A28A2800A28A2800A28A2800A28A2800A28A2800A28A2800A28A2800A28A2800A2
8A2800A28A2800A28A2800A28A2800A28A2800A28A2800A28A2800A28A2800A28A2800A28A2800A28A6BB2A296760A07524E05003A8AE7352F15DB5A
965B788CC47F1670055BD075F835947F2D4A491FDE07BFD281D99B14546F34492246EE03BFDD07BD39D9514B3B0551D4938140875148082010720F42
29A5D03842CBB88C85CF34019E2F67FEDD36BB4791B3D3BE339CD69D145001451514F710DB26F9E4545F5268025A2A182EA0B91982557FA1A8F51BC4
B1B379DF1C70A0F734016A8AE7341D724BDBE9A0B861CF29C631ED5B1A85EC7636AF3CA470381EA681D87CB796F0CC90CB2AAC8FF754F7AB15E45049
A8EA5AD4DA989199A27CF0B9001ED5E8DA2EA525EEF8668D95E3507730C1349306AC6BD14514C414554D46FA1D3ECA4B99880A83819C64FA571B79E2
EBB96CE70B0AC7BC6D5C1E568B8D23BA491245DD1BAB0CE32A734FAC7F0A41241A0C0263977CB9FC6B6281193E20BF7B1D3CF90713CA76C7FD4D6258
C524CE2C6D8825C833C8076EF506B2F2EABAF08AC55A46418562701463FF00D75D269F676DA359957946E6F99DD8F2C7DA915B2296BFA9FF0065DBC7
6362816575C020FF00AB1EBF5AE7F49D2DEE6E12375DF331CCB29E70334AD6F71ADEAB2953CC8E3E6278551CFF002AEC615B4D2ED84664440AB9627A
B7BD03D8B914690C6B1C4A15146001DAB2F53BCBD8EFA0B6B38721F059F19AC7BBF189FB5882CADC32E7EFB1EBCD75503F9B04721182EA188FA8A64E
C4B45145020A28A86E966784AC0EA8E7F88F614015F54D461D32CDEE2623E51C2E7935C5E91637BE20BD9350B91B632DB90B71815A3ADE952E21FB44
BE69964EC49AE905B456BA63C10A85458C8E3BF1D696E52D112DB2C0917976C50AAF5DA41E7DEAC571FA15D4B6973300AD2A39E066B3F5FD5B50BABB
10C45A18C1E1437F851715B53D028AE63C1C6EC4534733B344BD0B1E7356F56BE9669C5859E43B1DAEDFD29858D44BDB5791A359D0B2F519A7C73452
961148AE57AE0E715CA6AB61169D658F3F339EBB7A8AB5E0E8CC1A5CD7572C01790F27B014AE3B6973A6A2B234CD49AFAFEE901C451FDD18FD6AB78A
35A1A6D9BC508DD33A91C7F0E698AC6F0656CED20E0E0E0F434EAC7F0CC13C3A4C725D31334DF3B67B56A4D20861791BA2827AF5A04365B982175496
45566E80D4D587A6DAB5EBB5F5DAF2CD9443D3FF00D5506ADACCB25D2E9FA665A53F7DD7A8FA503B1D1D34B2AE37103270327A9AE6E1D5AE74BCC3A8
2C9293CA93D7F3AC98AF6FF57D7EDA425D21593800E00F4A570B1DE5316446665565257EF007A5666AFABA5946D1C1FBCB83C000676FD6A978522911
6EA6B96FDECB201D7EBC530B1D1D645E6A063D62DED14955EAC7B1CFAD6AB30452CC7000C935C4C0EBAADEDD5FDD4BE5420945E79045008ECDA58D58
2B3A824E00279CD4B5C5E97652EA1A9F98378B78D81DE7F8B15DA500C28A28A041451450014514500145145001451450014514500145145001451450
014514500145145001451450014514500145145001451450014514500145145001451450014514500145145001451450014514500145145001451450
014514500145145001451450014514500145145001451450014514500145145001451450014514500145145001451450014514500145145001451450
014514500145145001451450014514500145145001451450014514500145145001451450014514500145145001451450014514500145145001451450
014514500145145001451450014514500145145001451450014514500145145001451450014514500145145001451450014514500145145001451450
01451450014514500145145001451450015CA6A6D71ABEB62CADDCA4110E581E0FAD7575853E8D2ACF2CB6F3222B1CF3D71E940D1CAF8A66B5859349
D3D5649A4203482BA3F0BF878E90A25773B9931B3D3EB58DE1FD1E3BEF104D772FEFA38091BCF73E95DED2436FA1CFC39BFF0011C8EC088ED80D8436
41C1AADE28BC92E116CECDB29B879ACA7A56ABE9651E4FB34BE5ACA7E6F51FE3528D2ADC5A7918F9BA97C739A2C2B8DFB543A7E9B0F98E5888C6D04F
2DC560E852CF77E259AE6E032641D8B9E318ADC8B46B6570F2B3CA4740C78AB896B04731992245908C640A2C04F451453105725E2E76B8956D91B1B1
4139E9EB5D6D676A5A6ADE7CEA42C98C723AD034735E1891ADF51903BED8163C9639E4D5AD465B8D5E6630C44DAC1EBD09F5ABB06832248BE6CE0A77
0BD715B515BC30C6638A35553D463AFD690EE7112DC226AB6F3C09E584F94FBFD2ADDCC77DAD096764648A324283D38ADE1A2596EDC559B9CF27AD68
2C689188D5404031B7B6280E639CF0758AC16573BE15024973C8EB5D1A471C79288AA4F5C0C6688D1634088A155780053E989BB851451408E5FC4893
5EDCADAA2B1450015ECC4F3C545278694696247C99C0DC55B9C5757B5776EDA376319C738A7503B98F06AF6F15A469B5D9D102E00E32054064BED5DC
C4A0DBDB679603AFB56CFD9ADFFE7845FF007C0A9A8118D1E831412192DA668D88E491939F6A9A0D1E1425A666989E99E00AD3A281DD98C34668AE59
EDA7F2A3639E3EF2FB0A8356D1E34D39DD2476651CE4F5AE829A4060430041E083DE80BB38CD1F47B34CDD5E4B105539D83BD74F67A8DADDCCD05B12
4C633D3031D38A8DB46B2624956C139C03C0AB56F676D6A49B7895091824500D9628A28A041451450051D42D9E6F2A58C06685B76D23EF556BC7BEBC
45822B76855DB0EEDE95AF450172B59D9C36710489467BB6393593A9E964EA4B731C2640DCB606715BF4501739D9F55B988B476F6622FAAF39A7C5A7
5FA5D2DD068FCC93963FDDCF5ADFA281DCC8874488C8D35DB99656041E78E6A21A45C11F65F382DAA9E31D4FE15B945017313FB265B505ECE526461B
49FBA71542FF00C3125C2ABB4825949F9B3FE79AEAA8A02E471462285101C8450BF95477B134F692449F798715628A04623C7A898459C286350B82FF
00FD7AA1A6DACFA65D4CEB6CD34EDC138C9C7AD7554503B9876FA4BDCDD35D6A2096CE5133FCE8B9D1E459B36455118E7FDD35B945017332CB4782DC
334DFBE95FEF335676A5A5DEAC80D892431C0E718CF527F4AE92A1B99D2DADE49E43854193405CE7F58796DB4EB4D39667699C0DC73C91556C7C293A
39496509096DC46724D3FC351CFABEA136B176498D4ED807D3BD75B483622B7863B785628976A28E2A5A28A620A28A2800A28A2800A28A2800A28A28
00A28A2800A28A2800A28A2800A28A2800A28A2800A28A2800A28A2800A28A2800A28A2800A28A2800A28A2800A28A2800A28A2800A28A2800A28A28
00A28A2800A28A2800A28A2800A28A2800A28A2800A28A2800A28A2800A28A2800A28A2800A28A2800A28A2800A28A2800A28A2800A28A2800A28A28
00A28A2800A28A2800A28A2800A28A2800A28A2800A28A2800A28A2800A28A2800A28A2800A28A2800A28A2800A28A2800A28A2800A28A2800A28A28
00A28A2800A28A2800A28A2800A28A2800A28A2800A28A2800A28A2800A28A2800A28A2800A28A2800A28A2800A28A2800A28A2800A28A2800A28A28
00A28A2800A28A2800A28A2800A28A2800A28A2800A28A2800A28A2800A28A2800ACCF10DD1B4D1AE24438765D8A7DCD69D73FE2D4924B285541F2C3
E5F1FE7EB40D16FC3B62B61A3C1181F3B8DEE7D49E6B56A1B420DAC247431AE33F4A9A8130A28A2800A28A2800A28A2800A28A2800A28A2800A28A28
00A28A2800A28A2800A28A2800A28A2800A28A2800A28A2800A28A2800A28A2800A28A2800A28A2800A28A2800A28A2800A28A2800A28A2800A28A28
00A28A2800A28A2800ACDD66CE7BE812089F6296CBB671DBFF00D75A5450055D3ACA2D3AC62B483EE463009EF56A8A2800A28A2800A28A2800A28A28
00A28A2800A28A2800A28A2800A28A2800A28A2800A28A2800A28A2800A28A2800A28A2800A28A2800A28A2800A28A2800A28A2800A28A2800A28A28
00A28A2800A28A2800A28A2800A28A2800A28A2800A28A2800A28A2800A28A86EAE22B4B692E276DB1C632C68026A2B157C41084F365825484AEE571
CE6B427BC8E1B13764318C207C0EB8340EC5AA2B94B8F1CE9F6EF0A4904CAD29F9430C56BDBEB76D3DCC7014911A4FB85870695C2C6A51595AC6BF61
A2B20BE675F307CA55739ABB65790DF5B2DC5B9251BA6460D311628AAB6B7F6F74EC90C9965EA3DBD6AA5D78834AB3BB16B7176892938C1EC6803568
AA6754B00501BA886FFBBCF5AB6480092700753400B456747ADE9724C214BE84C8DD14B609ABE08600A9041E411DE801D45473CD1DBC4D2CCE1117AB
1ED4F0410083907A11400B45145001453490A32C401EA68565719560C3D41CD003A8A6E46EDB919C6714EA0028A28A0028A28A0028A28A0028A28A00
28A28A0028A28A0028A28A0028A28A0028A28A0028A28A0028A28A0028A28A0028A28A0028A28A0028A28A0028A28A0028A28A0028A28A0028A28A00
28A28A0028A28A0028A28A0028A28A0028A28A0028A28A0028A28A0028A28A0028A28A0028A28A0028A28A0028A28A0028A28A0028A28A0028A28A00
28A28A0028A28A0028A28A00298E88E30EAAC3D08CD3E8A0028A28A0028A28A0028A28A0028A28A0028A28A0028A28A0028ACABFD72CEC6E52D99B7C
CC40DA0F4AD5A0028A28A0028A28A0028A28A0028A28A0028A28A0028A28A0028A28A0028A28A0028A28A0028A28A0028A28A0028A28A0028A28A002
8A28A0028A28A0028A28A0028A28A0028A28A0028A6BBAA296760AA3A9270053619A39E31242E1D4F714012547E745E7793BD7CCC676F7C541A8DE25
8D9C970FCED1C0F535C7E83AC47FDAD3DFDFCBB565E130B818FF00269363B1DDD154ECF51B4BE245B4A1D80C91E83357298828A6B30452CC70A0649F
4ACAB6D7ADAE35036815978E1C9E0D006BD15936DAF5A5DEAC6C2DB323283B9C740476AD5240EA40A005A28A67991F99B37AEEFEEE79FCA801F51F9D
1F9DE56F1E6633B7BD4959B6DA718B5192F1E404B67007BD006951451400514514005145140051451400514514005145140051451400514514005145
1400514514005145140051451400514514005145140051451400514514005145140051451400514514005676BB62DA8E8D73688DB5A45E3F3CD68D55
D40DCADA3B59FF00AE5E42E01DDED4023CEC788F51D26DE3D3B53B4135AA809B87048F5AF40B416D7FA44210EF82488007D47F9154F514B1BFD2C9BA
8919DD3A01F306C7E756743B3361A4C3032ED603247A7F914914CE3BE26C50DB41A3F951AA0139036AF6C0AEBECA18AF2C2C2E4F2F1C68CAC3AF4E45
737F122D0DFDAD945010D2C336F2A3B022BA4F0EB2B68768A08CA461587A1A3A8BA1CD7C548F76816EE172CB70307D335D1F8618B787ACCB104F97DB
EA6B37E2069F26A1E1A758012F1C8AFC7A7F9C54BE089C3F876DADDC813C00ABA1EA393FE34750E8695900BA9DF2A800650FE24126B86F89F188F55D
26E235559093960393823AD77D6F1B2DFDDCA570AFB003EB81FF00D7AE2BE26452BCFA7C908CB2673F4A4F61ADCD6F1358C173E0995F6AACB15A878D
C0C104007FCFD6A3F873AACFA968263BB7324D6CC177939CA91C7F5AB9A849BFC1443E374B6E23031D4E3155FC09A59D3F4A2EE9B5DF03EB8EBFE7DA
98BA1C9DCDADBC7F155AD91311164C85ED9033FCEBD51542A85500003000E82BCBF578A58FE29C57417E432C673DB0001FCC57A950819C6FC4ABE306
826D61244D2B06E3B283FD6B5FC25A92EA9E1CB4B904160BB1C7A30E31FCAA8DCC33EACDA8B0B4134322182272E3803B807DEB9CF8693CDA76AF7BA3
5CAB26FF009D037F7875FD28EA3E87A6514514C931BC536E27D0AE58332C91A1642A71CD70DF0EF55BA83C452E9B792CA639E3250484F0C318FCC577
FE2204E8577B7AECE3F315C178AAC24B1F12695AAD929DF2C685BD0103149EE52DAC55F156A5A847E382629A58E2F33CB450701B0003FAD7A8D8AC8B
691899F7C98C93F5E6BCC7E215BB59EABA2B9CFCB182CDFED6EC9AF5381B7C11B6739507F4A103D8CBF15DF7F67F86EF6757D927965633FED1E9553C
0BAA36ABE18824958B4D0931499EB91FFD622A5D43CBD47527B3B884CB6D146490ABBB2C47B7D6B8FF0001DDBE91E2DBDD1A72E21B9CB45B863E71CF
EA33F90A2FA8ADA1DDEB31DCCB1430DA4A62679402EBD4753FD2B85D4AF756B3F1DC5A67F69CBB1D179CFA8AF4DAF2AF15BB45F13ADE46518DB1ED3F
850C68EF74D5BEB7D525B7BBB83346620E84FD715C8F8C35DD6347F12416F1DDEDB794EE0AA3B67A57A30AF2FF008971BAF89B4C9460EE5C28FC6860
B73A2F17DF6ADA668116A16B76A1900F306CC6E27A55CF04DFDEEA9A1ADEDF4AB2190E06074C75ACEF1E4376BE0C915E68D9515778D983FCEA5F866D
BBC2317B4AC3F951D43A193E31F10EB5A3788E0B6B7B8410CC37A8DB9C0C91FD2B7EEA3D7D2C1E78AEE376750C502F2BF4F4AE4BE26A63C53A6B16E2
48D473DBE735D869F1DEC7E232F7254A496A301790002280E868406F4E851B2B86BC310219B1C9AE3B4FF1DDD26B51D86A7147B1E5F28C8BC6D3FE4D
7A0A3ABA86460CA7A10722BC93C5F647FE11EB3D5ADD36B25CB2BB01CE72704FE543123B4F1AF89DFC3B0DBADBC6AF34C4E037602B43C31A9DD6ADA6
2DDDCC51A06C6DD9DFD6B8AF174A354F0847AC3AA97611C2BB87209E49FCF35D6780DFCCF08D8B1EBB4FF3A3A8743A2AC4D0B5E4D5EF2FE00154DAC9
B0007923D6AEEAF722DAC9B1FEB25FDDA0F735E6FA3C9FF08BFC40FB24B21F2AE7F77213D327A1FCE86091E9F793FD9AD259C0DDB14903D4D725AA78
E5F4C861927D31C2C871B89201FA715D84B1ACD13C6E32AE306BCFBE2DA22E97A72A8C6D90803DB028608EA2DF5F79AE2CA3365246977C86639AA5AD
78CA0D17555B2BBB4970E7E5603A8F5AD3D1E08AE743D32461CC51A3A9EE0E2B84F8A8DB7C43A5376080FF00E3F402DCF4BB4B85BBB58EE115956419
0187359B0788ACAE6EE682D5659FC9C8768D72011DAA3F11DDC969E10B9B887E593C80171DB381FD6B9BF8464B69BA8BB72CD3824FAF0680B6876961
A8C17FE688438311C30618ABB50436B1433CD320C34C416FA8A8F52BAFB258C92AA96908DB1A8EACC7A0A6212C351B6D40CE2D5F77912189F8E8C2AE
5798780EF67D37C5D7BA45E9C35D65F693D241C91FCFF2AF4FA486D5828A28A620AC9B9F1069D6BE699A56558890CD8C8E2B5AB96F1959C16DE16D46
486350EFF31247A9A068D083C4DA55C41E741333A138C85AD7560CA194E411915C37C3BB482F3C1CD148BF3999BE6C743DB15DCA2044551D14606692
063A8A28A620AC697C4DA3C5742DA4BC51213B7041EB5B35E65F12F4A105CDB6A16719572C0923EEEE0693047A13EA3671DC790F708B27A1FF001AB2
CC154B310001924F6AE5FC2F347E21886B7710A063FBB58FAE081824FEB5D4328652AC01046083DE9819D16BDA4CD742DA3BE88CA780B9EB53C5A9D9
4D37951DC219338DBD39AF35F155BFFC23DE31B7D4E38775BA90E571C11DEBB5F0C5A4735947A9CA8864B9FDEC78FE053D3F1A571D8D4D4354B1D323
DF7D731C23191B8F269F657F6B7F1F996932CAB804E3B66B8FF8ACA83C3B149B14BF9A141239C56AF80228D3C1F60E8A034884B1C724EE228EA16D0D
7BAD5F4EB362B7577146C0E082695B54B0458D9AE50090654F3C8AE53E226891CFA63DEC31A875FBE471CFAD41E0BB95D76DA08258911AC005947F7F
D3145C2C77AAC194329041190477A7535542A85500003000ED4EA620A28AAF7D75159594D7539C4712166FC2802479A28DD11DD55DCE1413C9A7332A
29666000EA49C5729E0E33EAEB2EB77BCF9AE440A7A2807AD69F8AE18A4D02E649037EE9770C1C77140EC6A413C370A5A191645562A4A9C8C8A27B98
2D9035C4AB1A938058F535C57C2790B787EE94FF000DD1C7D081FF00D7A8BC5AB25DF8C74C82590A59C724618678624E7A7E34AE16D4EEA0BA82E33E
44A8F8EB83D2A424282588007249ED5E6EB7AD61F15E3B5B7252097F74E99E395247F4AB9F10F59952F74FD0EDA6318BA71F6865EA14B0007F3FD28B
858EDA3BDB696630C53A3C80670A734EB8BA82D954DC4A9186381B8F5AC2D12E344B7BB5D3AD1365CAA705D705FD71583E2B579FC77A45BDF6EFB0C8
76007856E39FD4D30B1DDC17305C67C89924C75DA7353579F784D2487C71A9C7660FD814B26D072171D2B4FC497B3787AF20BF84FF00A3CAC16443D0
74CD203AEA2A38E459625910E55D4303EC6A4A620A28A2800A28A2800A28A2800A28A2800A28A2800A28A2800A28A2800A28A8E496389434AEA80903
24E39A00928A6B30452CC70A0649F4AE0F58F15DDDD4CF6FA67EEA356DA58FDE3EF8FF00EBD034AE77D4555D34CEDA7406E7FD7141BB356A81051451
4005155FED76E6E3ECFE72F9BFDDCD58A0028A2992388E3676E8A0938F6A007D739E29F107F6642D05AFCD76C3D3EE8ACED67C5D756F6F2B5B5B2A8C
61589C9150F867C3F35F91AB6B0EF234D87018FDEFF014AE3B7717C31E1E926BA4D4EF7748AC3702E7073F4AEE6B03C53AD7F61E9A3ECEA3CF6E235C
70056768DE3113A08EFA222518CB28EB46C369BD4EB5D95119DCE154649F415534ED4EDB52590DA96223201C8C557FED0B0D5A27B386E4092418C639
F5A688ADB41B265B642F24A78527EF1FF0A623628AE32F75CD5B4E9F13805E52008C8E9F4ADED27564BD448E61E5DC15C953DE80B1AB45145020A28A
8AE255820795BA28CFD68033B55D661B12228B12DC138D839C7D6A93789E3B684BDEC5B08E720E07D2A2D2E181DE5D5350705B79280E3F95737ADDC4
BE27F11AD9DA0616F1E14F2067DE91491DDE9DAAC37DA71BE2A628C1C1DDCE3FCE6AEC32A4F0ACB19CAB0C8ACE6B116BA01B383198D064FA91824D4B
A6CA9169513CF2228504139E3A9A622EBBAC685DD82AA8C927B53609E2B9884B0B8743D08AE5B5CBBBAD44982D78B6070587F10F5ABDA75DDA691A1C
4A1B7B2E7E4079CD01636E79E1B74DF3C8A8BEAC69619A39E31242E190F715E5BAC5EDF6AD70D3CA6455DFB625CE140FA57A268366D63A4430B48643
8DC49F7A5706AC8D2A28A29882B3754D66D34C8CB4CFB987F0AD6830254856DA48E0E338AF38F14C6D6F74D1C120925071F32E41A4C695CEC65F10D9
25AC72A92F24A32B18EB4586AF25CDE0B796109B8654D65F85FC3A6DE34BABF01A4FBC8A7B55EBD2916BA93CFB82200463E9FE340686ADF5DC7636AF
7131F9547E66B9AB7F1B5AB23199017CF091E738AABE24BA7D4E65895196D978C93FAD1E16D12DE59E499E30F120C027BFB503495B527D43C4B7C51B
ECB6FE58EC48C939AEA2C25926B18659862464058631595AF5B451C506C8C28079C77C6314B3EAD24D12C1A7C2C6665C671C2FD281323F106B9F6547
B6B239B8E8CDD93FFAF563C2FE79D1D24B87DCD2316073DBA7F4AA975A32C1A69924FDE4E4E5C9ABF657D656FA6C20C8A81571B475C8A00D19658E14
2F2B8451DC9AC8FEDE47D5A2B48955A390EDDC7A935972CF75AEEA8B083E55AC6C339EF525C5A456FE26B648C7CBB94E0D170B1D5D14514C41451450
01451450071FE2E17465D9E7110B81B501EA71FE39AD7F0F42DA768310BB22320176C9E949ACD84B737314F1206318F6CF1FFEBA8574ABDBB7517D31
10A9C8407B7A521F430FC43A8B6AB7D0DBC24A5BAB7CC73CE3D69749D1E3D56E6666DC208FE50E58927D38ADABAF0FA99D7ECA36A91F3331AD7B0B38
AC6D560847CABC93EA7D68B0EFA15F4AD220D35A478C967930093D87A5695145324CCD6DDFECA2041C49F78FA015C46B32244545BA0DDF7010707FFD
75E8B3797E5B34C14A28DC72338C573163689AAEB925D98C0B5858ED18C6E27A7F2A4CA4EC49E0DD19AC2D9AEEE03F9D3818F339602ADEA8CD75AB41
671B91B705B1E9DFF4C56E565DC58CE2EDA6B4654320E589FBA7BD02BDD8ED57544B25F2A21E65C30E14738F7359DE1BB7966B89AFEED9BCECE029ED
9CD6958E930DB1324D89A63D59874A965B1CC86482668588C10A38A6066EB1E238B4E9D608D3CC7DC03139C0ADB89FCC891C0C06506B24F876CD9D5E
46766CE58FAD6C2A8550AA3000C0A007514550D52FC59C20261A673845FEB408BF45456DE6FD9D0CE4798465B152D001451450014514500145145001
451450014514500145145001451450014514500145145001451450014514500145145001451450014514500145145001451450014514500145150DCD
C476B11965CEDCE3819340135159275DB75621E29C73C6133C52FF006FD9FF00767FFBE295C766693A248B8911587A30CD11C71C4BB634545F451815
987C4160A406320CFAA5386BBA798F7891BAE31B79A2E166699008208C83D41AAB1E9D6714FE7476E8AF9CE47AFD2AA378834E4505A56193D369A7B6
BBA62919B9183DF69C0FD28D00D3AAB71616974DBA7815CE319355D75DD29BEEDF4471E9934FFED7D3B19FB5C78FC698822D26CE3545F27784391BD8
9FD3A55FAA2BAC69CDD2EE338F7A7AEA562EBB85E418F7703F9D005497C3FA7CB3098A3AB2F4C374AD19A159ADDA266755618255B9C7D6982F6D48C8
B9808F5F305385CDB9E93C47E8E280196567158C1E4C25B6649F98E6B39BC3B69F6D8AED1E549A362DB811935ACB2C4D9DB221C0C9C30E2859636385
7427D01A00928A28A00A7A959BDF5B9804E624604361739A6CFA6DBDCC36D1DC02FF006620A9E992055EA2803135AF0F45ACCC925C4D811FDD5D99C7
EB5A91C6E96A22F306F55DA1C2FE47153D1401434EB292D0CE659448657DD9C73F8D64EA3E1AFB4EAE9A8C322A4CB22B8392B823F3CFFF005EBA5A29
5877238F7F96BE66DDF81BB6F4CF7C5731A9F8627D47505BEB896269E30BB31C7233ED5D5D14EC2B91C3E68897CF2A5FB95E95CCF88B409F59BF59DE
242B0AE223BC71DF35D551405CC4D634DB8D4FC3E2C19503BA287F9BA63D293C2DA5C9A3E9AD62D16C4572CA77039CF5FE5FAD6E51405CE3FC55E1CB
9D76FE19D6201601B572C3919CFAD6D4E7513A73410D96D90C7B0319578E315AD4503B99FA6C1269FA3C50C80BC91464900E493C9C560A69525E784E
F34E9EDE4594C8CCA1D48CE4E78FC335D751405CE375ED1597C1B6BA4441CB2B29242E791C91FAFE957FC116F2D9F87E3B3991D5E0623E618E09CD74
745017306F82DEEBD6F04A927930FCD9008CBFF86315CDF8E3C3C2EB51B6B8B1597ED18058E770C0E9C9AF42A29582E51D26ECDE5846EE1848140937
0FE2C735CA7C44B3FED34B782105A48492D8ED915DCD14C465786C04D06D221F7A2408C3D0D71BF112C24D4B55B592D86FFB34787DA33CE738AF47A8
9E085C92F0C6C4F52541CD2B0D333A6823D5BC3AD6E8D912C214E3B1C0E2B9BF87D69269126A36772762EF0CBBFE5AEDD2348C62355504E70A3151C9
696F2BEF9218D9B39C95E680B9474ED41EF753BB8D4036F105D8D8EF55AF59EFF5B82DA195516D879849E72DEDEB8AD98E18A2CF9512267AED503348
208449E62C51873CEE0A33F9D303CD7C59A6DC699E2AB7D5ADA6F325DE246083041EF9F622BD262B98A6B61711B6E4DBBB8A492CEDA5984D2428D263
192335224314685238D117BAAA8029581B317C3BE238B5B9AEA21118DE073807BAE719ADEACAD3F4A82CB53B9B982010ABA85007F17726B56840C2B0
7C67189FC3377086C3B27CABDDAB7AA1B8B682E542CF18703A66988E5FE1BC7E5786446C087590860456E6BBA89D334C92E5137B83851EF53C1A7DAD
BB878605461E84D4D3C115C44639903A1EC690FA95B47BFF00ED2D3A3B929B19BEF2FA1ABD51C10C56F108E140883A0152531057257F1C9E23D5EEB4
F59765A5B275C757AEAA58D268CC720CA9EA3355ED34EB5B225A08F0C460B12493401E7BE0AB89B42F11DCE9536E36F3B707B2B57A61200249C01D49
ACF6D134D693CC6B605F39CEE6EBF9D5C9618E684C3202C8460824F3F8D03672BAE599F135ADDB47810DA8611F1FEB5872467D3FC6A87C35D42E2386
6D1AF1587923CCB72C3AAF71F9FF003AECE1B0B782D24B6823F2E2933900E7A8C1EB54E0D02C6DE749E10E8F1B02B82063DBA74A560BE860FC50824B
8F0FC2912E5BCE071F856A78111A3F085846C305438C7A7CED5AD77A7DA5E10D7308908E064914B67650D946D1DB8654273B49C81F4A2DA85F433F5D
94DCA0D2ADF634D7236BE7908BEA6B87B582EFC25E368F7481AD6E70929FEF0F5C7B715E871E956D1DE7DAD4CBE6FA97EB4971A4595D5D35C5CC6656
600618F03E9405CBE0820104107A11597A56B96FAA5CDC410A3AB40792DDEAFDB4096D088A32C547F78E6A2B4D3ED6CE5924B78C219319A622DD735E
3D5B993C31711DB7F1F0F8EA475FE95D2D43756E9756EF04A328E30680461F80D90F842C028C1542ADF5CD5DF14027C377E075311AA3E1C89F49966D
2EE3E552E5E16ECD9F4ADABCB44BB411CAF204EEAA701BEB4BA0FA9C7FC2B85A1D12E92404379F920FD2BA2D5B46B6BF923B966F2A68483E6E33C0FF
003D6A7B2D26DEC1B36CD2A8272577707EB535F59457D108A62E141CFCAD8CD017D4F3ED0B4492FF00C687530CF3436F2963339FBC474AB1E2CD3A56
F1B69B78403189A323238C703F98AEF6DADA1B584436F188D0740292E6D20BB5417118708DB973D8D160B9C25D5ACD2FC52B6B8B74C24670C7A8FBBC
FF005A7EB28FA97C45D3ADEF10FD92DCB044C753B739CFA6715DB45676F15C493A42A2590E59FA93FE149358DACD709712440CC9F75C1208FCA8B05C
E2F4937365F10AE6DED066CE566F310738382739FC0559F8A273E1D8A3500C924E15457516BA65A5ACEF3C31012C8492E4E4F35977D60BADEB10348A
1ECED092493C33FB5006A68F1B43A3D8C4F9DC96E8A73EA1455DA28A620A28A2800A28A2800A28A2800A28A2800A28A2800A28A2800A28A2802B4B79
6F0CCB14D2AA3B0C807FC6B8DD47547D7FC4B169B69201676EF992407A9AD2F1AB4525BC3688A5AEA66C281D76FBFB560DA785EFEDCA242B2C734830
CC060019FF00F55265AB1B7E30D6E15D35ED2D24DF2B9DAC54F0077E6A3F09E8205A2DCDFC6198FDC079247AE6AECBE1887FB3444A434E3E624F4634
F59B537845A41088D9000580DB8A057D3429F89358BB1AA59E9DA7100C9305761EBE95A7E22D50E9BA71F2981B993E58C0EBEE7154E6D0E5B7B882EA
13E6C91B06638E734FB7D3A6D4AFCDE6A2BB5146225F6F6A0342DE83A835DE9A8D74EA265F95B3C1358BAEF8BA112A5869A4C92CAFB0BA8E9EB52EB1
A15E286FECF91C21209DBD6B2B4BF0F6A161A8A5E2DA96913382CBEBDFEB40D5886CED676F12DADB19242FB83CC41E76D7A3D73FA0E9B3437D3DEDD4
7891D76E5873D735D050852614D650E85586558608F6A7514C9393D7F4FB56D42D6223CB84B2B39F4C1ADBB9D52CAC6344DFBF8C2AA1DDD2A7BAB1B7
BB1FBE4C9C60107A54506936703EF11EF61C82E738A06725A94571A8EA2B7B756CCB6CE3E54F5C715734AF0F24F70D74F1F930ED2107727D6BB0A295
87CC71B63A0DD5A789219D5488933B9B1C1FC6BAC7B68A49D2775CBA0C0CD4D4531377396B8B6FB778BE26DAC63B7E5B3D32077AB7ADE95BA237164A
1255E4E2B68471890B8450CDD580E4FE352501733B439EE2E34D47BA044809049FE2F7AD1A28A04151CD124D134520CAB0C1A928A00E72FF004CB7B0
D3EE2E6E6662A832A0700554F05691E4097519A3556989110EFB7D6B43C556F717363147029605F951DCE38FEB5A9A747245A7C11CA143AA00428C01
487D0B559926910CB316676F2C9CF9606315A74531112411243E4A46AB1E30545531A358872C62DC0F627815A345007377FE178EE1C185C281D3713C
7E55BB676E2D6D22801CEC5C66A7A281DC28A28A0415850F87D0DFBDCDDC9E6FCD9500F5E7BD6ED1400556B9B286E590CC09DBE9DFEB5668A00A9269
F6B22AA989405F418CFD6A68208ADA211C08110738152D14011CD1473A6C95432FA1A228638576C48147B549450035943A9561904608AA11E8F668D9
285FD031E0568D1401CBDEE8B7A75069AD30A80E570C07D2A4B2D26FDF518AEEFE5DC6339C920F4EDC5749452B0EE14C92448A33248C15546493DA9F
54756B16BFB4F25582FCD9E7A5311883599B51D7EDE0B2DDF678DFE723B8EF9AEA6A8695A55BE9906C8465CFDE7239357E800A28A2800A28A2800A28
A2800A28A28031BC493B4562A80B059090C57AE315774B8E38F4DB75890A2EC07047278EB52DC5AC172A04F1870BD324D4C0000000003A01400B4514
5001451450014514500158161653DC6AEF7F760855E111BFBD9EDF4ADFA2800A28A2800A28A2800A28A2800A28A2800A28A2800A28A2800A28A2800A
28A2800A28A2800A28A2800A28A2800A28A2800A28A2800A28A2800A28A2800A28A2800A28A2800A28A2800A28A2800A6B2AB7DE507EA29D450042D6
D037DE8233F541510D3AC94605B478E7B55BA280291D2AC09C9B58CD2FF65587FCFAC75728A00CD3A1E9A5CB9B44248C6327148DA0E98CBB45AAAFB8
26B4E8A02E659D074D311416E067F88139A43A069A558791F7801D7A56AD1405CC8FF846F4DDBB42381ECD4D93C37A7B0E048A7B7CD5B345160B98D2
7876CD9308F2AB7625B355A4F0AC4D8DB772211E8BFF00D7AE8A8A561DCE793C3780375FCBD3A28C7F5A917C3EEA7FE4233E00C631FF00D7ADDA28B0
5CC1FF008479B18FED19BF2FFEBD43FF0008FDE8391A8BF279F98D74945160B9CE1D02F71817ED8EF9279A70D1F550AA06AB20C760E71FCABA1A28B0
5CE7FF00B2B581D3566FC49347F666B80E575407EA4FF8574145160B9CD4BA7F88558797A8861EC7FC69E96BE2146045D83C73B829AE8A8A2C17300A
788C37134447A6D5148478901E1A13FF007CD74145160B9CD3CFE248C8531C6C48EAAA0FF2A925B9F11C650FD9A020019006EC9FC0F15D0D14582E73
8FA9788532469B138CF40AC3FAD3BFB535BCF3A50FC03574345160B9CDC9ADEAD13ED934B03233D48A4FF84835003274C27E87A574B45160B9CD2F89
2E8021F4B7CF621B8FE54E5F13CBD1F4B9811FED75FD2BA3A28039D6F14045CB69F3FE7FFD6A72F89E22326CE61FAD6F9008208C83DA98B0C487291A
29F50A053031C7896D3606782E54F71B471FAD2FFC24F658CF93738F5D83FC6B61E28DCE5D158FB8CD33ECD0631E445FF7C0A00CB3E26D3971BBCE19
F54A0789B4D29BC34B8CE3EE569ADADB2FDDB7887D10535EC2D1C61ADA2EB9E140A5A868678F12E98582EF933E9B29C3C47A615CF9AE31FEC1AB674B
B13D6D93F2A3FB2EC324FD963C9EBC51A868551E23D3482448FC7FB0681E22D3081FBE3CFB559FEC9B0CFF00C7AC7F9534E8DA7939FB2A67EA68D434
22FEDED34AE44F9C1E98A135FD2DBADDA29F46A90E8DA71FF9744FCCD1FD8DA7FF00CFAA7E668D43401AD699FF003F910FA9C52FF6C69BDAF613F46C
D35344D350B15B55F9BFDA3FE34CFEC1D3BCBD9E47FC0B71CD1A8684DFDAFA7633F6B8BF3A06AFA79C62EA3E7DEA14D074C4520DB2B739CB52B687A7
3127ECE064638278F7A03427FED4B11FF2F51FE749FDA9A7E71F6B87FEFAAAE7C3FA718C2184F1FC5BB9348FE1FD39954084A85EE1BAFD734C342CFF
006B69FF00F3F717E74BFDAB61B37FDAA3DBEB9AA7FF0008DE9D9CEC6FD3FC29478774D09B7CA6207A9A0342E0D46C49005DC3CF23E714A751B20706
EEDC1CE3FD60AA23C37A681C46DF5CD07C39A7160C51F8F7FF00EB501A1746A36458A8BB8723D5C50BA8D933151770E47AB81544F86B4D2DBBCB607F
DEA0786F4E0C1B63E47FB5FF00D6A0343445EDA9E9730FFDFC147DB6D7FE7E60FF00BF82B34F86B4E2D9DB267D9BFF00AD40F0CE9A0F0927FDF5FF00
D6A0343485EDAB1C2DCC04FA09052FDAEDBFE7E21E3FDB1598FE18D35C92564E4E7A8FF0A3FE117D33046C7E7AF23FC280D0D3FB5DB7FCFC43FF007D
8A4FB65AFF00CFCC3FF7F05667FC233A6FF765FF00BEBFFAD4A3C33A6818C4BFF7DD01A1A7F6AB73D2E22FFBEC51F6AB7C67ED1163FDF1596DE18D34
8C62503D9E93FE116D2F71256520F62FC501A1ABF6AB7CE3ED1167FDF14A6EADC75B88BFEFB15923C2DA686C8128F6DC3FC29FFF0008E69F9C859339
CFDEFF00EB52D434342496CD991A49212C872A4B0E0D3BED76C064DC458F5DE2B2FF00E117D3B24E25E4E7EF0FF0A55F0C69A01CABE49CE723FC28D4
3434FED76C064DC438F5DE28FB5DB633F688B1FEF8ACEFF846F4FDB82B21FF0081527FC237A7ED0A564201CFDEFF00EB530D0D1FB65B7FCFCC3FF7F0
52FDB2DB38FB4439FF007C5677FC237A76FDDB64CFFBD41F0EE9EC7244A7FE07406868FDAEDBFE7E21FF00BEC52FDAADCF49E2FF00BEC566AF8734F1
9F96424F72D4F3A069F8184718EFBA802FFDAADC759E2FFBEC5114B6C311C3245ECA847F2154068161CFCAE73FED559B6D36D2D995E18B0EA300EE26
811728A28A0028A28A0028A28A0028A28A0028A28A0028A28A0028A28A002A1BA9D2D6D64B894E12352C6A6AE6FC68669749FB1DBFDF9C9071E94022
87858BEBDAADC6B972088E36296EB9E07BD767597E1BD38697A25B5A6DDACAB961EF5A940DEE145145020A28A2800A28A2800A28A2800A28A2800A28
A2800A28A2800A28A2800A28A2800A28A2800A28A2800A28A2800A28A2800A28A2800A28A2800A28A2800A28A2800A28A2800A28A2800A28A2800A28
A2800A28A2800A28A2800A28A2800A28A2800A28A2800A28A2800A28A2800A28A2800A28A2800A2B1B55F10DA69CE6200CD30EB1A1E95734BBF5D4AC
85CAC6630491B49E6802ED14514005145140051451400514514005145140051451400514514005145140051451400514514005145140051451400514
514005145140051451400514514005145140051451400514514005145140051451400514514005145140051451400514514005145140051451400514
514005145140051451400514514005145140051451400514514005145140051451400514514005145140051451400514514005145140051451400514
514005145140051451400514514005145140051451400514514005145140051451400514514005145140051451400514514005145140051451400514
51400514514005145140051451400514514005145140051451400514514005145140051451400531911882CAA4AF42474A7D14005145140051451400
514514005145140051451400514514005145140051451400514514005145140051451400514514005145140051451400514514005145140051451400
51451400514514005145140051451400514514005145140051451400514C76091B39E8A09350585D8BC84C813661B18CE6802D514563F893578749D3
99E43F3C99441F51D68044F0EAF6D3EA4D631E4C8A0927B715A35CBF82F4F963B47D42ED489EE7901BAA8AB771AD4926B31E9D628B2107F78DD80EF4
AE36B53768A8E7992085E594E110649AAFA6DF26A16BF688D4AA162A33DF1DE988B945433DCC16C019E4540C70327AD48ACAEA1918303D0839A00751
4514005723E23F11379A74DD29C994F0F2A73B3DA99AFEA5A8CF27D8E0578B7B6DC275FC4D69787BC3B6FA65A869E3592E98E5D9B9E690F6D4C5D0BC
332CCF1DD5E48C571924F563D735DAC3124312C518C2A8C0A92B07C4FACAE9969B5240256E4E0F2A28D82ED9BD50DD4E96D6D24F27DD45C9AF3ED1DF
56BE8E5BD59E45488ED3E593827E95D35BCAFACE93F673281286C16718DC051706AC5BD0F56FED58E47DA06C207CA0E3BF7AD6ACFD1F4C8F4AB4F223
3BB2724D64DDDDCB06BCAD35C6D85189201FE11DA981D3515CC7FC245717D24834C854C2BC798E39CD6CE9D79F6985564204E07CC31FAD0162F51451
40828A28A0028A28A0028A28A0028A28A0028A28A0028A28A0028A28A0028A28A0028A28A0028A28A0028A28A0028A28A0028A28A0028A28A0028A28
A0028A28A0028A28A0028A28A0028A28A0028A28A0028A28A0028A28A0028A28A0028A28A0028A28A0028A28A0028A28A0028A28A0028A28A0028A28
A0028A28A0028A28A0028A28A0028A28A0028A28A0028A28A0028A28A0028A28A0028A28A0028A28A0028A28A0028A28A0028A28A0028A28A0028A28
A0028A28A0028A28A0028A28A0028A28A0028A28A0028A28A0028A28A0028A28A0028A28A0028A28A0028A28A0028A28A0028A28A0028A28A0028A28
A0028A28A0028A28A0028A28A0028A28A0028A2B3E7D4E28355834F2A4C9302411DA8029CDA85E3789E2B08D40B655DCE7B9E335B95188A312B4A114
48C305B1C9152500145507D62C12FF00EC2D3AF9E3AAFA7D6AFD0014515CF5CF89228FC450E93100CCCE15CE0F04D0163A1A28AC6F146AC349D1E595
1B13B82B17D680360107A1A5AE43E1D497B368B34D78CCCAF2E6327D3BD75ACC154B310001924F6A018EAA967A85B5EBCC96D20730B6D6FAD477530B
DD3E51A7CAB23B0DA0A30E2B33C3F603C3DA4CF25F3057790BB9CE73E948674545733A87886758925B48C2A7525B9CD6FDA4E2E6D229C0C798A1A985
89E8AC2D5FC55A5E93722DA79774DDD57F87EB5AB637715F5A477301CA483228158B1451450014514500145155B50BA5B2B09EE9C802242DCF4A006D
9DC4970D36F8F62A36D5CF53F5A7DD5DC3691EF9DC283D0773591E10BE9F52D2E4BB9FEEBCEFE5B7AAF4FE957AF34EFB5DF472C8F8851305475279FF
001A4325B0D421BF5730861B0E0EE1566491235DD23AA8F563815CD696534FBED5242CCE9080B81DF93583A95FDC6AED3CD348B1C11B6DF2F76338E7
A5171A8DCF4249A297FD5CA8FF00EEB03525707F0FD84B7B76C6666655185E8304FEB5DE5313566145145021A4850492001C926B097C59A73DD8B741
2312DB738A9FC4D1C92E8EEB1965F986E23A815C85869F04720B9BA6D8ABD00EA71EB48691DBEADABDAE950EFB86F98FDD407934FD2F508F52B41711
0DA09C633935CE5869F26BDA89BEBF56F2106D41EB53DACB6DA35FDE13910A291827DE80B1D2CD3470466499C220EE6AB5B6A967772F9504C19F9C0C
1E715CD69D1DCF89EE1EF2E1DA3B452551454DA74515A789D614FBA8ACB9C753D280B1D44D2A410BCB21C220C935CB5AEB57FA96B28B6C36DB06E171
D477CD6A788CC8D6D1DB458CCED8627D07FF00AEAC40965A35A2C5BD5401CE7EF35302DDCDC456B6EF3CEC1510649AC87F124011648E091E36E43671
C5626AD7971AFEA30E9F6CED15A9601C63923D4D6DEAF6D6B65A208800BE5A85423A9A43B0D93C5562B1968D249368F9B1818AD1D26FC6A567F6811E
CF988DB9CD6468BA259C3A3192EA12CD22976CF503AD57B5D4E3D2F4B94C085A79643B53B0F4A05637356D561D3200F20DEEC7840D8247AD41A06BB1
6B5148C9198DA3EA09EB5C26BB1DF5F5D99AE15FCC03239E39F4ADBF087D974BBA785E6DD35C10AAABCE07145C76D0EE28A28A6494B54B84B6D3E691
FBA951F522B888350D52DED9A48DD911C9208E95BBAC4D36A7AB45A55B83E521DD31E99E3A5539F4FBEBD97ECEB11548D785C600F634996AC6BDB6B7
145A640D72E24B96405912B9A9ADAEF5FD6A37BFF92D54E4ABF0AA07A5751A4E85059A24932879F6E0E7902AFDDD9437617CDDC36FF74E2815D2D8C9
D4B5958A130E9CBB8A8C6E1FC23DAA3F08E9DE45BCD7938CDC4CE724F615B76F676F6D098A38C6D6E5B233BBEB546F34EB96CC769318E290E5B27EED
02B991E2ABE373247616D28C370C47AD6EC260D234C8E276198D3EEF7635593C3D68A992CCD31393211DFE94F83468CBF9978C6523EE804E00A03439
2B9BD9F55BBB9B894ED5886D41FD05753E13666D1137162439192793D2B3B59F0EBC93036884A311C6EE86B7B49B2FECFD3A2B6CE4A8C9FAD036F42F
514514C930A58D20D6C5CDD608E4AB75C03D3F2AD03AA592900DC2F3EC6AD3A2C8851D4329E08355D74EB4520881723D49340CCEBFD66431B2E9F196
61FF002D1871F8572D736571792C9E6B19A7946DC7D6BD0C2284D81405E98038A863B4B78A532C712AB9EF4813B1434FD296CB433651FDE6049FF7BF
C8AAFA468D35A5F1BA99C0C02020ADFA2985C63B0446739C2824E2B8E3A75CEB5A9485B31C1B8963E99AED29AAAA8A15142A8E800C0A02E57B1B182C
2010DBA800753DCD612CD2C5E29DAA0ED772AC6BA7A80DA426E7ED2501971804F6A009E8A28A04145145001451450014514500145145001451450014
514500145145001451450014514500145145001451450014514500145145001451450014514500145145001451450014514500145145001451450014
514500145145001451450014514500145145001451450014514500145145001451450014514500145145001451450014514500145145001451450014
514500145145001451450014514500145145001451450014514500145145001451450014514500145145001451450014514500145145001451450014
5145001451450014514500145145001451450014514500145145001451450014514500145145001451450014514500145145001451450015C7D893AA
F8EA6BC8D49B7B1431E73C16C62B43C4FAABDAC49636403DE5C9DA141E547AD5DD0B4B4D274F584732B1DD237F798D219A745145311E55AAA5C45A8E
A1772E4CE5FE518C1CF35E89A5DD01A45AC977222C86205B26A1D4F41B6D425F35898DCFDE2067350C7E1B838134F248A06318C521DEE667897C55E4
21B6D34333B0E661C05AC0B08CC3AAE9F30B8F37CC9448E57A927B7B5767A9F872D2EEC05BC2AB1B20C027BFD6B234AF09CF697F0DC1904691481B6E
ECE7E94149AB1D83B2A233B1C281927D057956B73DC788F55BA952476B4B772B10E82BAAF1CEB9F64B41A5DA066BDBB1801467683EB5A1A1787E1D3F
428EC6650CE465D875CFD4D0F512D03C31F66B0F0DD9C4F346A553272C3D6B9BF1A6BF71751B58E9C48808FDE483AB57487C336DE603E749B01FBB8E
DE94FD4FC3D6B7569E55B2240DEBD88A05A1CC7C33D3EED0497B2332DBB2ED45CF079AEB354D3E5BFBD8158FFA2A8CB8CF7CFF00FAAAC68F6034DD32
1B4041F2C6091DEAF509037A9CFF0089A0823D2047142A189C28031C77FE954EF75C169A425B69AACF3AC60123F878ED5D2DD5A4377184997201C83D
C557B6D22CAD9195210DBFA97E4D31A68F234D22EF50D5BCA043393991DBB66BD774AB7B7B3D3E2B4B6915D6150A4820F3EF8ACCB9F0CC524B2CD04D
B1A43C0C7007D7357B46D28699148BE7190C841391D31492B049DCD3A63491AB846750CDD013C9A8EEEE05ADAC93B296083381DEB90F0DB5E6B9E20B
8D4AE58AC1036D8D074CD324EDA8A28A002B91F1BFDB2EA116102308DC02594FDECF15D752100E323A74A0114B46B14D3349B6B18C61614C7E3D4FEA
6AF51450073574B7705C5DC315B16FB4B139F5EE31546CFC2524CE5EF4EC56E4824139F5E3FAD7674503B9CFE8FA00D2F529268F1E5ECC03DD8D7414
5140828A28A006488B221491432B7041ACDFEC1B3EFE6104E482C3FC2B568A006471A451AC71A8555180076AE6B58D2E6D475B08A1D62CA976C718C7
FF0058D7514500436B6D0DA5BAC16E81117A014DFB1C1F6CFB56DFDEE319CD58A2802B5DD9C776815F20AFDD23B5558F44B4570F26F908ECC78AD3A2
8030EFB4C9D2ECDD58108DB7A28C7E1491E953DDCA1F52762AA0606E0735BB4503B8DDABB366D1B718C638C56747A35B25DF9E7900E4263815A74502
33EFF4C8AF4EE6628F8EA07EB55ACFC3F6B6B7A9759DCE8B81C639F5AD9A280B8566EBB7C6C34D92446C48785F5AD2AAB7D6305F44239F760742A79A
00C4F06C3BAD26BE725A599F049EC07FFAEBA5A82D2D62B3B758205C22FEB53D00145145001451450014514500145145001451450014514500145145
001451450014514500145145001451450014514500145145001451450014514500145145001451450014514500145145001451450014514500145145
001451450014514500145145001451450014514500145145001451450014514500145145001451450014514500145145001451450014514500145145
001451450014514500145145001451450014514500145145001451450014514500145145001451450014514500145145001451450014514500145145
001451450014514500145145001451450014514500145145001451450014514500145145001451450014514500145145001451450014514500145145
0014514500145145001451450014514500145145001451450014514500155EFAEE2B1B492E27601106793D7DAAC5739E23D12EB58BA8104C45A8C6E1
9C63FC68020F0C59BDEDD4DADDDA7CD331F2430E83D6BAAA8E1892085228C61114281525001451450014514500155B50BC8EC2CA5BA9BEEC6338CF5A
B3591AEE9B3EAAB15B8755B7CE5CF7A00E7FC2D652EB3AA4DE20D4632183ED807B0E86BB7A86D6DA2B4B748205DA88302A6A06DDC28A28A041451450
0145145001451450041736E9736ED0499D8FD71D7AD32C2C2DB4F83C9B58C22E727D49AB545001451450014514500145145001451450014514500145
14500145145001451450014514500145145001451450014514500145145001451450014514500145145001451450014514C9644863692560A8A3249E
D400FA2B2B4DD662D4AE1E286270AB9C39EF5AB400514514005159FA8EAD6DA7E04A7739E8A28D3757B5D4B70B763B946482280B1A145151CF208617
95B184527938A00928AA7A65FA6A169F6845DABB8AE339E9572800A28A2800A28A2800A28A2800A28A2800A28A2800A28A2800A28A2800A28A2800A2
8A2800A28A2800A28A2800A28A2800A28A2800A28A2800A28A2800A28A2800A28A2800A28A2800A28A2800A28A2800A28A2800A28A2800A28A2800A2
8A2800A28A2800A28A2800A28A2800A28A2800A28A2800A28A2800A28A2800A28A2800A28A2800A28A2800A28A2800A28A2800A28A2800A28A2800A2
8A2800A28A2800A28A2800A28A2800A28A2800A28A2800A28A2800A28A2800A28A2800A28A2800A28A2800A28A2800A28A2800A28A2800A28A2800A2
8A2800A28A2800A28A2800A28A2800A28A2800A28A2800A28A2800A28A2800A28A2800A28A2800A28A8F2FBB1B576FAEEE7F2C500494514500145145
001451450014514500145145001451450014514500145145001451450014514500145145001451450014514500145145001451450014514500145145
001451450014514500145145001451450014514500145145001451450014514500145145001451450062DFEB9F6595E28EDD9D94E324D509F54D6A4B
669A3B611A6DE0AA1FE66BA631C6CE1D914B0E848E4573FE2CD58DADBA58DA9DD777070AA3B0F7A43466F856EAF6FF005D9E5965768E14DB265B826B
B5AC7F0CE93FD93A5A44E14CEFF3CAC3B9AD569111955DD54B70A09C669831F45145021923AC685DD82AA8C926B93BC9AE7C477A2D2DB747671B659B
D7EB5278C35816EA961148A8D211E613D87A56C684B669A720B275718F9987526915B2B96ACACA0B180436E9B57B9EE6ACD14D61B94AE48C8C647514
C900C09201048EA3D2A3BB9D6DAD649DC801149E6B3B4DD34E9F34B713CEB823FBC703EB9ACBF15EB110B58E2B693CCDC72405DC0E281A576675E4E0
DA4D7D7EC4BCB9D8A7AF15A5E06B4916CE6BE990A34ED8553FC2B5CADEFDB352B8863963282460107A0F6F4AF4EB6896DEDA3854001140E2A56E54B4
44D5CDF89AEDDE48F4CB605A597A81EFFE7F5AD9D4AF63B1B379A470BC617EB591E1CB3925964D52E89324B9099F4F5A64A35B4BB31636115BF52A3E
63EA6AE51453105145140051451400514514005145140051451400514514005145140051451400514514005145140051451400514514005145140051
451400514514005145140051451400514514005145140051451400514514005145140051451400514514005145140051451400514514005145140051
451400514514005145140051451400514514005145140051451400514514005145140051451400514514005145140051451400514514005145140051
451400514514005145140051451400514514005145140051451400514514005145140051451400514514005145140051451400514514005145140051
451400514514005145140051451400514514005145140051451400514514005145140051451400514514005145140051451400514514005145140051
451400514514005145140051451400514514005145140051451400514514005145140051451400514514005145140051451400514514005145140051
45140051451400535890A4A8DC40E06719A751401CFDC7885ACE263756CC1C640DBEBF4AE5BC2578758F17C979788C5F0C621D877CD6FF008DEE1C59
A5B471A967C1DC7A8C9C715A9A168D67A5DA45E4463CC2832E7A9A4335AB0EE74DBC9B5F8EED8836F1B290377403DAB728A620A42400493803A934B4
50079CEA7A6CBA86A725DBB8C236E604F18AE87C276F329B89DC058C9DA8A071F8569CFA45BCD36F39543D51475ABF1C69146B1C6A15146001DA8B15
7D2C3E8A28A0928EAB626FAD4C4B2143F5E0FD6B985F0A5D45209108620E402DD2BB5A281A76394B5B5BB935BB779EDCA88CEE1F2E303D6BABA28A01
BB983E28B39AFA2B68625DC1A4C1E2B66DE210DBC71000045038A968A041451450014514500145145001451450014514500145145001451450014514
500145145001451450014514500145145001451450014514500145145001451450014514500145145001451450014514500145145001451450014514
500145145001451450014514500145145001451450014514500145145001451450014514500145145001451450014514500145145001451450014514
500145145001451450014514500145145001451450014514500145145001451450014514500145145001451450014514500145145001451450014514
500145145001451450014514500145145001451450014514500145145001451450014514500145145001451450014514500145145001451450014514
500145145001451450014514500145145001451450014514500145145001451450014514500145145001451450014514500145145001451450014514
50014514500145145001451450014514500145145001451450014514500145145001451450051BED360BE9A0926CFEE8E71FDEF6ABD4514005145140
051451400514514005145140051451400514514005145140051451400514514005145140051451400514514005145140051451400514514005145140
051451400514514005145140051451400514514005145140051451400514514005145140051451400514514005145140051451400514514005145140
051451400514514005145140051451400514514005145140051451400514514005145140051451400514514005145140051451400514514005145140
051451400514514005145140051451400514514005145140051451400514514005145140051451400514514005145140051451400514514005145140
051451400514514005145140051451400514514005145140051451400514514005145140051451400514514005145140051451400514514005145140
051451400514514005145140051451400514514005145140051451400514514005145140051451400514514005145140051451400514514005145140
05145140051451400514514005145140051451400514514005145140051451400514514005145140051451400514553D4B51B7D32D5AE2E9C2A8E83B
93401728ACFD1F56B6D62CC5CDA9F973820F6ABE4850492001C9268021BAB98ED2069A63845FCCD36C6F22BEB65B88376C6E9B860D73BAA4573AAC17
170A716EA76A64F6F5C56CE8102DBE910229278279FAD2B8CD2A28A298828AAF7D77158DA497131015067AF5AC0F0BCB73A95D5C6A9333794C76C59E
8450074F45145001451450014557BCBCB7B284CB7320451FAD73777E2B95DB6D8DB92A4FCAEDFE14AE349B3ACA2A1B6677B589E4C6F6404E3D715353
105145140051451400514514005145140051451400514514005145140051451400514514005145140051451400514514005145140051451400514514
005145140051451400514514005145140051451400514514005145140051451400514514005145140051451400514514005145140051451400514514
005145140051451400514514005145140051451400514514005145140051451400514514005145140051451400514514005145140051451400514514
005145140051451400514514005145140051451400514514005145140051451400514514005145140051451400514514005145140051451400514514
005145140051451400514514005145140051451400514514005145140051451400514514005145140051451400514514005145140051451400514514
005145140051451400514514005145140051451400514514005145140051451400514514005145140051451400514514005145140051451400514514
0051453588505988007249ED400D9A58E089A599C2220CB31ED5E6FAEEAB3EBFA93456C0FD92207693C7D7F1AD4D73529B55BA92DE270B630F2CC3F8
B155ECB4E7B6D0E7BE750A25388F6F042D22E2AC68FC3D564B0B9509B62126149FBC7EB5B1AE4EC562B1849F36E18038ECB51F85A136DA2A99405DEC
CF93E9EF4BA60FB66A93DF942107C91927FCF6FE7409EE5F5B244D37EC699DA1700FBF5FE755F479A668E482584A2C27686C75AD3A299214C9644863
692560A8A3249ED4FAE0FC6BE200E874EB224E09F34AFB7A7A8A0695CCEF126B92EBBAAC3A659026DDE50AC7BE075AF45B2B58ACACE2B68502A46A00
02B89F879A08556D5EE55B7B92210E3B7AD77D490E5D828A28A648550D5353874E8773FCD237DC41D4D4F7B751D9DABCF2FDD51D3D4FA563E9B66750
65D535004124B244DD1476340D1421D3AEF5395AFB5A674B65059632718FC2ACE8F610DDCFF6930EDB784E225619DD4FB8B8975DBB369644AD9C4DFB
D9BFBDEC2B7E1892189628C6150600A436C928A28A6485145355D5890AC091D403D2801D451450014514500145145001451450014514500145145001
451450014514500145145001451450014514500145145001451450014514500145145001451450014514500145145001451450014514500145145001
451450014514500145145001451450014514500145145001451450014514500145145001451450014514500145145001451450014514500145145001
451450014514500145145001451450014514500145145001451450014514500145145001451450014514500145145001451450014514500145145001
451450014514500145145001451450014514500145145001451450014514500145145001451450014514500145145001451450014514500145145001
451450014514500145145001451450014514500145145001451450014514500145145001451450014514500145145001451450014514500145145001
45145001451450014514500145145001451450014514500145145001451450015CBF8975279A65D26C24FDF39FDE30FE11E86B43C41AA7D82DC450E4
DCCC308076F7AE682FD9213146BE76A973C06EEB9A452424566F7D731E936795B68CE6E24EE7D6B5F598BCFBDB2D2AD48091E0B2FB7BFE15A5A55845
A369ACD2316908DF33FA9ACEB39D2DE2BAD6EF70AF28C460FA7B7E940EE5CD66EFC98A3D3ED06669405007F0AD69595B2DA5A4702FF08E4FA9EF591A
05BCD72E756BB63BE607621FE11D2B7E992C28A2A0BBB98ED2D9E799805519FA9F4A0465EBFA8344A2CADC133CFC647F08AE0A4B069F5FB7B58A4323
B3E1829FD49F4AD592FA7B9B9768B7CB757242AE3AA8FE95A167A2CFA3CC97D2A2BB1605806E9ED9FC296E5AD0EBA1892085218861114281ED58BABE
BE2DAE3EC7649E6DC03873D4271FCEA1D6FC442080C764ACD230FBDD31EA3D8FBD605ADFD8E93A3C8AFBEEB539C9076F504F38CFD28B8922E378BEE9
A54B48E35F38C813795E4FE15DAAEED837E37639C74CD703E12B4B73792EA17E151A3398D5DB255BD7DFAFE95D2EA3E21B7821945AEE96551C607028
40D6BA15F52946A1A98B279156DA16CC8738CE3AFE358BAF7899AEB1A569EA6156654326EE719C7159F15ADD5D235ECCCC90312C06719A6D9D893796
F752C6046AD9CB0EA010682D247A1E95A7C5A6D8C76D105C81F31031B8FAD55D6B546B3D96F6DF35CC9F746338ABEB776ED1975990A85DC4EEE82B99
B6D42DA3D4A5BEBD932173B3033C93FE14C843E5B4BCB6B66BFBAB868E5246D04F427D7F1A3C21AE5D6AD3DC893E7823030C7D7D07AD64EA775A9F89
1DA1B789A3B44C93804E71D39A7585ACF6969169965136F672CC57B9EFF854DF52ADA1D26A3ABEE94D9587EF276E0B29FBB5976E6E2CB518813862DF
3E4F02B6744D256C2269655FF4994E5CF5DBEC2B9DF114B2DFEA62D2C5499647D9B97F840EF43B8935B1AD7FE268C5CFD974F512C99C337615A1A3DF
3DEC0E6520BA1C1C0A8B46D12DF4EB601E3579D87CEC46734CB7922B2D4A5461B1189EBDBD29EA2D3A1B5593AA6BF61A6BF973396948385419A8B55D
761B78B65AE6695F81B73F2D70496B23EA32DDDE2EF7E557737CAB9EFF005A1B1C637DCF41D2F5CB7D4E731451BA919FBD5AF5C8786ADA2B5B896F66
C461970B97CEEF702BA28F53B396710C53ABB9E3033D684C525665DA28A2992145145001451450014514500145145001451450014514500145145001
451450014514500145145001451450014514500145145001451450014514500145145001451450014514500145145001451450014514500145145001
451450014514500145145001451450014514500145145001451450014514500145145001451450014514500145145001451450014514500145145001
451450014514500145145001451450014514500145145001451450014514500145145001451450014514500145145001451450014514500145145001
451450014514500145145001451450014514500145145001451450014514500145145001451450014514500145145001451450014514500145145001
451450014514500145145001451450014514500145145001451450014514500145145001451450014514500145145001451450014514500145145001
45145001451450014C91C471B3B745049C53E9AE8AE8C8C32AC3047B50070B73737175AABCC886499F0B180385AE8B44D1CD966E6ED84978FD5BA850
7B0AD282CEDEDF0618514FAE327F3AB140EE60F8A2F7ECF6C90065024E5F3FDD1583A5C171E28BA0D723CBD3A0C055031BB1DAB5F55D0A7D53516338
516E48E41EC3DAB7ED6DE2B4B7482040B1A0C00290EF644888B1A2A20C2A8C01E829D4514C90AC0F11C37572121B704E54E06CCF3DEB7E8A00C5D074
48F4B4691C0699F193D76F1CD6C3AABA95750C0F5046453A8A00C3D5B454B85DF6B0A873904671F8D54D27C2B0C0C27BE02593B27607D735D3D140EE
65BE8564CBB42B2AE73806A55D2ACD207892151B860B1193F5ABF45017673E341926655B8955615E3627719AD8FB241E4AC5E4A144185057A558A280
BB333FB16D79E6419CE791FE148DA1D83B02D19238F973C715A9450176450410DB44238235441CE00A548628C9648D149EE140A928A0435977215C91
918C8EA2A9586990D8C8F229DF23FF00111D055FA2800A82E6DA1B94DB2AE7DC75153D14014EDB4DB6B625922058F196E78A49B4CB499B718429FF00
64015768A02E657F60D963A49FF7D74A9ACB48B3B190C9046779FE226AFD1405C28A28A0028A28A0028A28A0028A28A0028A28A0028A28A0028A28A0
028A28A0028A28A0028A28A0028A28A0028A28A0028A28A0028A28A0028A28A0028A28A0028A28A0028A28A0028A28A0028A28A0028A28A0028A28A0
028A28A0028A28A0028A28A0028A28A0028A28A0028A28A0028A28A0028A28A0028A28A0028A28A0028A28A0028A28A0028A28A0028A28A0028A28A0
028A28A0028A28A0028A28A0028A28A0028A28A0028A28A0028A28A0028A28A0028A28A0028A28A0028A28A0028A28A0028A28A0028A28A0028A28A0
028A28A0028A28A0028A28A0028A28A0028A28A0028A28A0028A28A0028A28A0028A28A0028A28A0028A28A0028A28A0028A28A0028A28A0028A28A0
028A28A0028A28A0028A28A0028A28A0028A28A0028A28A0028A28A0028A28A0028A28A0028A28A0028A28A0028A28A0028A28A0028A28A0028A28A0
028A28A0028A28A0028A28A0028A28A0028A28A0028A28A0028A28A0028A28A0028A28A0028A28A0028A28A0028A28A0028A28A0028A28A0028A28A0
028A28A0028A28A0028A28A0028A28A0028A28A0028A28A0028A28A0028A28A0028A28A0028A28A0028A28A0028A28A0028A28A0028A28A0028A28A0
028A28A0028A28A0028A28A0028A28A0028A28A0028A28A0028A28A0028A28A0028A28A0028A28A0028A28A0028A28A0028A28A0028A28A0028A28A0
028A28A0028A28A0028A28A0028A28A0028A28A0028A28A0028A28A0028A28A0028A28A0028A28A0028A28A0028A28A0028A28A0028A28A0028A28A0
028A28A0028A28A0028A28A0028A28A0028A28A0028A28A0028A28A0028A28A0028A28A0028A28A0028A28A0028A28A0028A28A0028A28A0028A28A0
028A28A0028A28A0028A28A0028A28A0028A28A0028A28A0028A28A0028A28A0028A28A0028A28A0028A28A0028A28A0028A28A0028A28A0028A28A0
028A28A0028A28A0028A28A0028A28A0028A28A0028A28A0028A28A0028A28A0028A28A0028A28A0028A28A0028A28A0028A28A0028A28A0028A28A0
028A28A0028A28A0028A28A0028A28A0028A28A0028A28A0028A28A0028A28A0028A28A0028A28A0028A28A0028A28A0028A28A0028A28A0028A28A0
028A28A0028A28A0028A28A0028A28A0028A28A0028A28A0028A28A0028A28A0028A28A0028A28A0028A28A0028A28A0028A28A0028A28A0028A28A0
028A28A0028A28A0028A28A0028A28A0028A28A0028A28A0028A28A0028A28A0028A28A0028A28A0028A28A0028A28A0028A28A0028A28A0028A28A0
028A28A0028A28A0028A28A0028A28A0028A28A0028A28A0028A28A0028A28A0028A28A0028A28A0028A28A0028A28A0028A28A0028A28A0028A28A0
028A28A0028A28A0028A28A0028A28A0028A28A0028A28A0028A28A0028A28A0028A28A0028A28A0028A28A0028A28A0028A28A0028A28A0028A28A0
028A28A0028A28A0028A28A0028A28A0028A28A0028A28A0028A28A0028A28A0028A28A0028A28A0028A28A0028A28A0028A28A0028A28A0028A28A0
028A28A0028A28A0028A28A0028A28A0028A28A0028A28A0028A28A0028A28A0028A28A0028A28A0028A28A0028A28A0028A28A0028A28A0028A28A0
028A28A0028A28A0028A28A0028A28A0028A28A0028A28A0028A28A0028A28A0028A28A0028A28A0028A28A0028A28A0028A28A0028A28A0028A28A0
028A28A0028A28A0028A28A0028A28A0028A28A0028A28A0028A28A0028A28A0028A28A0028A28A0028A28A0028A28A0028A28A0028A28A0028A28A0
028A28A0028A28A0028A28A0028A28A0028A28A0028A28A0028A28A0028A28A0028A28A0028A28A0028A28A0028A28A0028A28A0028A28A0028A28A0
028A28A0028A28A0028A28A0028A28A0028A28A0028A28A0028A28A0028A28A0028A28A0028A28A0028A28A0028A28A0028A28A0028A28A0028A28A0
028A28A0028A28A0028A28A0028A28A0028A28A0028A28A0028A28A0028A28A0028A28A0028A28A0028A28A0028A28A0028A28A0028A28A0028A28A0
028A28A0028A28A0028A28A0028A28A0028A28A0028A28A0028A28A0028A28A0028A28A0028A28A0028A28A0028A28A0028A28A0028A28A0028A28A0
028A28A0028A28A0028A28A0028A28A0028A28A0028A28A0028A28A0028A28A0028A28A0028A28A0028A28A0028A28A0028A28A0028A28A0028A28A0
028A28A0028A28A0028A28A0028A28A0028A28A0028A28A0028A28A0028A28A0028A28A0028A28A0028A28A0028A28A0028A28A0028A28A0028A28A0
028A28A0028A28A0028A28A0028A28A0028A28A0028A28A0028A28A0028A28A0028A28A0028A28A0028A28A0028A28A0028A28A0028A28A0028A28A0
028A28A0028A28A0028A28A0028A28A0028A28A0028A28A0028A28A0028A28A0028A28A0028A28A0028A28A0028A28A0028A28A0028A28A0028A28A0
028A28A0028A28A0028A28A0028A28A0028A28A0028A28A0028A28A0028A28A0028A28A0028A28A0028A28A0028A28A0028A28A0028A28A0028A28A0
028A28A0028A28A0028A28A0028A28A0028A28A0028A28A0028A28A0028A28A0028A28A0028A28A0028A28A0028A28A0028A28A0028A28A0028A28A0
028A28A0028A28A0028A28A0028A28A0028A28A0028A28A0028A28A0028A28A0028A28A0028A28A0028A28A0028A28A0028A28A0028A28A0028A28A0
028A28A0028A28A0028A28A0028A28A0028A28A0028A28A0028A28A0028A28A0028A28A0028A28A0028A28A0028A28A0028A28A0028A28A0028A28A0
028A28A0028A28A0028A28A0028A28A0028A28A0028A28A0028A28A0028A28A0028A28A0028A28A0028A28A0028A28A0028A28A0028A28A0028A28A0
028A28A0028A28A0028A28A0028A28A0028A28A0028A28A0028A28A0028A28A0028A28A0028A28A0028A28A0028A28A0028A28A0028A28A0028A28A0
028A28A0028A28A0028A28A0028A28A0028A28A0028A28A0028A28A0028A28A0028A28A0028A28A0028A28A0028A28A0028A28A0028A28A0028A28A0
028A28A0028A28A0028A28A0028A28A0028A28A0028A28A0028A28A0028A28A0028A28A0028A28A0028A28A0028A28A0028A28A0028A28A0028A28A0
028A28A0028A28A0028A28A0028A28A0028A28A0028A28A0028A28A0028A28A0028A28A0028A28A0028A28A0028A28A0028A28A0028A28A0028A28A0
028A28A0028A28A0028A28A0028A28A0028A28A0028A28A0028A28A0028A28A0028A28A0028A28A0028A28A0028A28A0028A28A0028A28A0028A28A0
028A28A0028A28A0028A28A0028A28A0028A28A0028A28A0028A28A0028A28A0028A28A0028A28A0028A28A0028A28A0028A28A0028A28A0028A28A0
028A28A0028A28A0028A28A0028A28A0028A28A0028A28A0028A28A0028A28A0028A28A0028A28A0028A28A0028A28A0028A28A0028A28A0028A28A0
028A28A0028A28A0028A28A0028A28A0028A28A0028A28A0028A28A0028A28A0028A28A0028A28A0028A28A0028A28A0028A28A0028A28A0028A28A0
028A28A0028A28A0028A28A0028A28A0028A28A0028A28A0028A28A0028A28A0028A28A0028A28A0028A28A0028A28A0028A28A0028A28A0028A28A0
028A28A0028A28A0028A28A0028A28A0028A28A0028A28A0028A28A0028A28A0028A28A0028A28A0028A28A0028A28A0028A28A0028A28A0028A28A0
028A28A0028A28A0028A28A0028A28A0028A28A0028A28A0028A28A0028A28A0028A28A0028A28A0028A28A0028A28A0028A28A0028A28A0028A28A0
028A28A0028A28A0028A28A0028A28A0028A28A0028A28A0028A28A0028A28A0028A28A0028A28A0028A28A0028A28A0028A28A0028A28A0028A28A0
028A28A0028A28A0028A28A0028A28A0028A28A0028A28A0028A28A0028A28A0028A28A0028A28A0028A28A0028A28A0028A28A0028A28A0028A28A0
028A28A0028A2B91D6BC4A6C3C67A769A250B0B8C4D9E80B74FE9401D7514514005145140051451400515C4EA9E34BCD33C409A54BA6A6E91C047DE7
04138CD76D40051451400514514005145140051451400514514005145140051451400514C9244897748EAA338CB1C53E800A28A2800A28A2800A28A2
800A28A2800A28A2800A28A2800A28A2800A28A2800A28A2800A28A2800A28A2800A28A2800A28A2800A28A2800A28A2800A28A2800A28A2800A28A2
800A28A2800A28A2800A28A2800A28A2800A28A2800A28A2800A28A2800A28A2800A28A2800A28A2800A28A2800A28A2800A28A2800A28A2800A28A2
800A28A2800A28A2800A28A2800A28A2800A28A2800A28A2800A28A2800A28A2800A28A2800A28A2800A28A2800A29AECA8A59D82A8EA49C0AAEDA8D
8A10AD7B6CA4F4065519FD6802D51502DD5B31C2DC44C7D038353D001451450014514500145145001451450014515465D5F4E8665864BB8848E70AB9
CE4D005EA28A2800A28A2800A28A2800A29AEEB1A33B9C2A8C93E82A858EB1A75FDD4B6D6772B2CB128670A0F03A50068D1451400514514005154750
D56CB4D42D79308F033D0F4A934EBFB7D4ACD2EAD18B44FD0918A00B5451450014514500145145001451450014514500145145001451450014514500
145145001451450014514500145145001451450014514500145145001451450014514500145145001451450014514500145145001451450014514500
145145001451450014514500145145001451450014514500472CA90C4F2CAC151016663D80AF20F1DE9B770DEC7AE4AD9378E1A25539DB8191FA735E
99AF2BDD5B2E9D0B859AE4E393D147249FCAB9EF1BE9BA9DD7865A3682D644B5024531160CA1460F5F626931A3A3F0DEA4356D0AD6F47DE9130FCE7E
61C1FE55AB5E6DF097522D0DE697237DCDB2C4B9EC786FE95E934C1851451408CDD4B5BB0D31D23BA9B12BFDD8D4658FE15069DE24D3351BE7B18647
4B95E7CB917693F4AF3CB4BA9EEFE2A96BE61BA3B831AAFA60E0015D8DC784D5FC611EBAB75B30CACC87AE47181EC695C6D1C97C4298C7E3BB29000D
E52A1C7E39AF40BCF12E8F65771DB5DDEA452CA0150C0E39E9CF4AF39F889163C796D9FBB245103F991FD2B6BE2ADAC506816023451E5C9B01C738C5
007A07989E5F99BD76633BB3C63D7359D3788747850BCBA842AA0E0924F5AC5B384EA3F0CE2865B9580BDB60CACDC2ED6EE7F0C572D1FD8EDBE19DDD
BB29BA90C84ACC89F2AB1231F37E0680B1E9F0DF5ADC5A1BB82749200092EA72303AD430EB1A64F2F970EA16CEF8CED128CE3E95C8FC287693C3D731
38CA2CDC03C8C11CD727A4E98973F10AE2C009162F3A456F2DB690BCE79A2E163D86DAEADAED0BDADC4532A9C131B8600FE1435EDA20CB5CC2A3A732
015CBE9BE14B8D1B47D56DACAF886B91FB8258809F5CF7AE5ADADACAD3C15AA5B6AF34535DEF668CA92CCA7B73F51FAD1703D5629639937C5223AFAA
9C8A79200249C01DEBCF3E12DC4A74BD42277668E2756553D890D9FE42B3341D4E4F1578EC9BB27ECBB58AC20E0600E3FA5170B1EAA086195208F515
C94DE358A2F1647A11B3604CBE534A5C607A1C551D1F4BD6B4CF1C5C1862906972B1C64E576E38E7D6B9596C607F8A66CC86589AEB1C37238CF04D01
63D9E982443218C3A97032573C8FC2B8FF001DEAEFE1CF0F436BA7B324B704A239392A072C73EBCD72F3AB4FE0DB3D4EC27BAFED62C0C80127CC0091
9FD3345C2C7ADD2120724815E75E25BDB9D57E1FC3AB3B4D6F796EE16454CAE49201E3F2358D6167E21D7BC2D14B67772C8B6B36DF2F790CC78CE09F
4CF4A2E163B8F1C68779AEE9D047612949629776376D0C31EB5B9A6C325B69D6D04EE5E48E355663DCE2B81F88577ABE9165613C1A95C44F708239E3
47E37051923EB5A1A84BAADDF80B476D3E797ED739883329E58153D7F1C7E5401DCD15E4DE2DD5756D2751B6B0B0D6EEA6B858C79E03701BD3DEBAD1
2EBBA3F81E5B89E4377A928DE377CC54123F3C5170B1D6515E6375A86B0DE035D70EA1770DD09871BF0A46719C7F4AEA3C0FAB5DEB5E195B9BB9019D
6468BCC0A39C63071F8D30B17ED7C45A55DEA874DB7BA57B919F940EB8EB59FE39B5D62EF4B853456904825CB843824579F7842DAE66F884EB1CC126
89A5667DBF81E3F1AEBBC77AD6B5E1CB6B696DAF239167250EE8802081D45203ABD145E0D1ED46A249BB118F349EB9ABF5CA58F894DBF80E0D7751C3
CA50E42F1BDB7151FCAAAE8BAC6B7AEE8379AB41710DB98D9C450F94083B467927EB401DAD15CAF843C507C49A64CA76457F0AE180FBBD386FCEB1BC
3DE31D66FF00C54349B9861645665728B8231DE98AC7A1D15C16A1E2DD5749F1841A5DEC56EF6F33A85F2FAE1B81FAD77B40051451400515CB6B9E29
369AE5AE89608B2DE4EC0312784CD41ACF8BA6F0FEBF0596A90C46D6740C258F395E704E28037F5CD49349D1AEAFDF9F263240F53D00FCEB848C788B
54D0B4ED66C2F2EA59E498992246F954027B568FC5396E478693C968BEC92CA818E4EE279231DB1C0AB1F0C12E17C2B0B3B46607790AF5DC0EEC7D31
C1A433AF8779853CD003ED1B80F5EF52514531051591E26D4A7D23449AFEDE1495A1C1657CE31F87BE2B0FC29E2DBFF10D8DC98ECE2FB54278F98AA6
0F4CD007674570DE1AF1E1D57527D3AF6CFCBB82488C43C86C7507269D73E35BAD37C509A46A7651246E542C91B9270DD0F34AE3B1DBD1583E2BF124
5E1CB049DE26964958AC6A3A647AD61F883C67A8E86B6134B676D2C5771893E52C081E9FAD311DD5158D79AA5CFF00C23C9A9E9B6CB3BB4624F2D89F
BB8C9AE6347F88771A9DD5C40DA74307950B4837CC4648EDD2803D02B8CD3F53F10378EA7B1BA8C7D8467002E0018E0835B5A0EA1A8EA56E2E6F2D23
B68997283712C4FAFD2B9FD17C6977A8F8B868D3D8C5060C88EC092D9504FF004FD690CEE68AE1FC59E34BBF0FEB4B629670C91B2AB2BB139C1EBDFE
B567C51E28D4745B3B5BF8AC226B59480448FF00372323A74A62B1D7D15CA2F8AA7BDF0C36AFA5D8F98621995653B40C7503D6A7F07789D3C4D69348
D12C3342F868C36782383FCE803A4A2B8F83C5571ABF88EEB47D2922516E85BCD933F3E080463F1A83C27E349B58D6A5D2AFE08E399436D743C3153C
8FCB34058CEF889ACEAFA36AF6CB63A83C50DC26760FE1C1C57A059B17B381D9B7168D496F5E3AD7977C5F3FF139D3C0EA2027FF001E35E99A41DDA3
58B1EA6DE33FF8E8A065DA28A2811C878C75CD5F4AD4B4D8B4D83CC8656FDE1DB9DC738DBEDFFD7AEBEB85F1878C6E740BD5B46B089D987991484E70
B9C03F5E2B63C55E245F0EE8F1DDB4425B898848E3CE06EC6493EC3FAD2195BC7F797DA668BF6FB0BD7B774608570086CFF5AB3E06BBB9BEF0BDBDD5
ECED34D23312CC79EBD2B93F19DFEA77FE02B3BABD86155BA955FF00764F00825735BDE02B98ACFC050DD5C36C8A2DEEC7D81A03A1D8515C4D978C6F
B53D3B51D46C6D205B5B33C79A4EE71F81AD6D03C516FACE8736A2B1321B707CD8872460678A623A0A2B8AD3FC55A9EB3A36A9A858416D0AD9B1D8B2
64B10064E7B74AB1E05F154DE248AE92EA148E6B6DB929F7581CFF008500739E20D4B5183E22DB589BD9BECCD711911838001C1C7EB5EA15E45E2573
FF000B6ADD49C81736F8F6C84AF43D475A5B6D5EDF4B8BCA1713A19374AD8551CFEBC521B36A8AC08756D4078922D3AE2180C1246CE2588939C0ADFA
620A28A8AE0CAB6EE6DD55A5032A1BA13401CAF85EFB5FB8F10EA306AA87ECB18250EDC0072303F2CFE55D7D713E16F195DEB9E229F4D9ED2189228D
DB7212492AC07F5AA7AF78EB51D1B5FF00ECD92C6DE408EB9285B2CA4718E7AF34867A1515C3F89FC5BABE82D6734D616EB05C027CB2E4B8C6320F6E
E3F3AEC6CEE16EECE0BA41859A359147B119FEB4C44F45145001456078AFC4B6FE1CD3C4CE0493C8711459FBDEE7DAB2F59F12EA7A069D637B7D1413
FDAB1BA2404797C03F7BBFE5401D9D15CCF88BC5F69A2E8F6F79E59965BB40F0C5EC46727DAA9DFF0088355D3BC3167ADCCB0C8B36C77894636AB0C8
E6803B2A2B1ADB5C8F52F0E36ABA6856210B6C93F848EA0D617843C5F7FE201711359C6B321CABF21147FB5401DB515C168FE35BD3E269745D5E08CC
BE67971984630DE9CD3FC51E2DD4FC3BAD5BC53436F25B4C376D5CE40CE3AD2B8EC7754573BE2ED72EB46D05752B28E37CB28224CE403D2A8787FC43
AB6B7A0497AB0C30B441F3230382473C0FA5311D8D15E7FE13F1B6A3ADDC3E9EF6B13DCF2E24076A220EA5877FC3D69D178BB53B3F1A8D135248648E
49163468C6082D8C1A571D8EFA8AE3FC7DE22D4BC3D05B4B60B094949562EB9391576F757BC5F052EAF6ED12CE2012B647078E7F5A623A3AC8B6F116
9577AA9D32DAE7CCB90092A01C0C75E6B33C03ACDF6B7A34B737F20693CD2A0A80315C2F87A1BCFF00858D347A75C4692069373BA646DEFC521D8F64
A2B96F1A6B37FA1E910B58C625B894EC3215E878E703D6B17C65AC6B1A069DA7490EA0CD2CEB9901403A01FE3FA530B1DFC922C68D248C151416627B
0159DA56BBA76B13CF169F379DE4637B01F2F3D3F9541F6A9752F060BC62164B8B2F31B68F55C902BCFF00E175BDE4F7D7ED697BF678D047E68081B7
F271D7A743F9D2048EC3C63A76B5A84D66BA4C8E91A9CC9B642B9E7BD7510A9482346392AA013F8579F7C42D7757D0E7B686CAFE40B32B336634C819
C6320574D236A373E17B6B8B3BA2939B759189504BFCA0D00CDFA2BCEFC05E21D42FB56BE8758BEDE218B2AAC303AF26B7F419AE45A5D6AF797B24B6
C431891B81B477FD3145C474B4579BE8FE25BAD7AEF51B99F53FB1416C84C112B05DDE99EE6B6BC0DE20BDD5D2E6DB505065B73912631B867D28B8EC
75D4514531114F04571198E640E879C1AF18D22CADEFBE204569751AC90191818C138E01AF6CAF13D205CCFF001110584AB14C667C3B2EE1C039E3F3
A068E8FE25585A585BDA496117937049FF0056C470318AEABC1377757BE17B49AF771948232C392077AE67C756B2DA4B69A8EACC2FE04601605F9003
F4AEA744D5EDF54F0E1BAD26308123655887F0301C0A5D43A1B94579E69FF6E9349D56FB56BD9D2ED188894C98DA3E9527C34D4F50BE7D421BC9DA58
E223696EC4D170B1DFD15E5BE2AD4358B4F165BD9596A1719691422B9E066ADF8D86ABA1DADADF41AACF25C33379849C2F6E83D280B1E8F456359DDD
EDC784E1BB8B0D7925B071C672C47A77AE41E3923F094D7B7B7728D524932479BF30E7FBBDA988F48A4240192401EF5C6FC37D4AEAF3C39732EA170D
29B7B8640CDCE142A9FEA6B22DB519BC5BE38367348F1E9F6C8CE215623781C738F7348763D2ABCB35D444F88F622318CDC005474AE8345B6D5ED3C4
ED1B09C69F83F2B676AE01C73DEB97F12C6D71F10A04B794A4A6E028607EED008F5BA2BCEFC77E76870DB5C5A5F5CB4ACD97DED90DCD775A6CC6E34E
B699BABC4A4E7D71CD3116A8A28A0029A0820104107A11585E2F8AFA6D28AD8C81173FBD3BB6F1C639AE2BC6125BD8D8DA49A4DE486558F7332C84E0
8C52B8CF4EBBFF008F59BFEB9B7F2AF3DF870FBBC4BABE32004191EBF35761A6DC4B3F84E0B8998BCAF681D98F73B6B90F86884EBFADBE4F01063B72
4FF85007A1CB2C70A6F964545F5638151CD796D05BFDA259E3587FBFBB835C45C2C9AB7C415B1D443496B0ABB2C7B885E3A1C0FA5655DD9A6A3E311A
042F2AE9E642CC11CF18CE7F95170B1EA2B22344240C3611B8376C7AD3629E298130CA9201D763038AE07C7BABCD636F169168CF1C28A88E47DE6181
C67E955668E659ED24F0E5B5F42100DF95C6F345C2C751E3AF2FFE11F7DFB776EF97239E87A527822E634F0B5A0B89238D88E0338048EC6ABF8D9B3A
05BCB708F91F3328EB9C7358DA66956EFE04FB7B44D24ECDBB7331CAAE7B7A51D43A1E88D22226F67555F52702A17BFB348D9DAEA1DABD4870715C27
86F507D5FC3FA8D9DD83235A8CC4D9E40FAD61787F423AEDE5FDB2CCE9E5203B77F04E7BD0163D7619A39E31242EAE87A1539A61BCB50DB4DC43BBD3
78CD6341A55C69BE1636115D032A29F9C9C0E7B66B8EBBFEC783C3CD1B4C26D40C9F3BC7927BF19A607A7AB2B28652083D0834C927861204B32213D3
73019AE47E1F4F33F85A66694B1477D84F38C5625DF882D20F323D534F9AE247248979E3E94AE163D2A19A29E3DF0C8AE9EAA7352560F83DED9B4189
6D9CB80496C8E413CE3F2ADEA620A28A2800A28A2800A28A2800A28A2800A28A2800A28A2800A28A2800A28A2800A28A2800A28A2800A28A2800A28A
2800A28A2800A28A2800A28A2800A28A2800A28A2800A28A2800A28A2800A28A2800A28A2800A28A2800A28A2800A28A2800A28A28032BFB233A98BF
3773175040181DFF000FE957AEE13736B2421CA6F1B77019E2A7A280390B0F03C1A5DFC37DA7DE3C73C7C312BC32F7E2BAD50768DC4138E4818A7514
05C28A28A00E3F5FF08BDC6B316B7A4C822BD8DD5D91BA3906B5EEECAF354F262B90B6F046CAEE15B25C8F4F415B3452B05CF3FF0016F85754D6F5E4
D46D955042A154171F36D2483F8D5DF18E8BAAF8874BB4B58ADD124421E46328E0F4C7F5FC6BB3A298EE70F3F87B539BE1F7F63310B730B7CA01E1D4
1CE3F5FD2A968DA2EB83C2177A04D68222E4B47331C770483F957A2D14AC17385F87FA6EAFA2A4F657967B62662E5FBF4E307A1AC8B6D0F59B3F1A36
B66C6630B5CBB954193B0D7A8D14582E61F8AE0D42FBC35711697B96E2451F29E1B6F71F5AE37C2D69A9BF85750D0E4D2D84936EF2E6200193FDE3FC
ABD3A8A623CCBE1C417FA26A9796379692813B2A92070A573CFEB496BA0DD784BC6915F846934C91994BA8CEC0C3807F1C7E55E9D484020820107A83
4AC3B99A2F24BC9E05D3CE610D99A52840C0FE119EF5E71AE5B5DE97F1317537B693ECE6759038524118C1E9F8D7AD535915C61D430F4233408E2FC6
BA5CFE26F0F25CDAC2CB25BC8CF1A13CC89DCFB74E9567C1FA9C16FE1AB6B6BADD1DCDAA14788A9DC79CF02BADA8C471890C811439182D8E48FAD03B
9C6F8CA4BA9BC0B3C7731017570E0A46AB8214480F23D40A8BE15CCB1E85359CBFBB9D2E19B6370482ABCFE86BB9650C30C011E869890C48FBD22456
231B828071405CF3DF8B8CB2DA5844843491C8CCCA3A8040AD0B3D7574CF86F697500DF32C421451D43D76135ADBCFFEBE08A5EDF3A06FE748B696CB
00856DE21103908106DCFAE2981E55E03B3866D42EF5BD6581FB39DC15FEF33FAE3BD6ADEF8A35897C257FA9408559AF8C2994FF00551EDFF3F9D7A0
0B2B51D2DA01FF006CC522D95AADBBDBADBC4217C9640A3073ED4AC173C9F6D9DE7C3A32C6D249A92CC4B29249FBC33C7A62BA4F8577B08D0DEC1C95
B859998295232081FE15D8D9E996162AEB676B1421FEF055EB4FB7B0B3B591E4B6B68A2793EF322804D0173C9FC233A58FC47BA374DE5A9699496FAF
15B7F15668AEB42D3A685B7AB4CC54E08CF15DACBA26973DC1B896C616989DC5F6E093EB525DE976178A8B736B148B1FDD04631405CE10E9926B1F09
6CE0B41E64D062408BD490CD91F5C13563E1DCA20F046A226F95A09652EADC11F20EBFAD76F65656D630F93690AC51E776D5E99AAD2E8BA7C824FDC0
4137FAD084A8907A1F6A02E703F09B4CB949EEB5170CB03279439E19B20FE954FC22E23F8A3741881B9E651EF5EAF6F0436D0AC36F1AC71A0C055180
2B3FFE11ED1C5C8B916318994E438241CFE7405CF3AF1D909F11AC1C9C00D09CFF00C0857AD5635DF86748BDB81713DAEE9B39F30B1273F8D6BA2844
0ABD1460739A621D45145007956AD65358FC56B4B99B2219E75911FB63B8A77C55B596EF5ED363814BBB425703D7757A36A1A65A6A213ED510668DB7
238E194FB1A626936835017CEAD2DC2AE15E439C7D290EE71DE3FB76B6F879656D29669227881DDD7214E6B5FE1B107C15640766933FF7D9AD7D5344
B1D5828BE88B85EC1881FE7934BA5E8F69A4A94B2F3123C7FABDF95FAE3D68034A8A28A6230BC6DFF227EA7C67F73FD45733F0800FEC7BF3819338E7
FE035DD5FD95BEA16AD6D749BE26EA338CD54D2B41D3B4832FF67C262137DE018E2819E6BF0F157FE13FBAF3106552623FD93B87F4CD1F10DD078FAD
49200558437D3766BD06CFC2BA558DE9BCB68E449CF570E7279CD1A8F85348D4EF5EF2EE17799F196DFD31D314AC3BAB9CFF00C4BD6E3B28ED74FF00
B2C534937EF37483210038E2B07E27472476FA32C93F9A4C2718501470BD2BD0B57F0CE97ACC30C77F13398576A386C363EB55AE3C1DA4DCE989633A
3C823E52466CB29C63F2A045BF0CCB1B786AC1D59760817273D3039AF29F17580B7D6A6D5F4E88B58C92FCAE47CACE3EF0FA57ABD96836767A336950
B4BF676EA7773DBA1FC2993F86B4DB8D2A2D325590DBC4DB946FE680B8BE19D66DF5DD1E1BC81550E36C918FE061DBE9E95E77E1F607E30CE4720DCD
CF3F83D7A1691E1BD3F4572D60664041DC85F2ADF514C83C2DA55BEA2BA845132DD07DFE66EE49EFF9D00701F15588F16D8F4F96DD0FFE3ED5D0FC5A
6DBE1CB703A35C018FC0D743AAF8674BD5EEFED57B13BCA136021B181ED4ED47C3B61AA430C37E679A387EEAB49C67D4FBD00739E1F60BF0AA42BC7E
E24CE3F1AC9F846EB0AEAB2BF0AA8A4FD066BB883C3B610695369B1F9C2D66EB1F99F77D853346F0BE9BA2CAF2589986F5DAEACC0861EE31401CBF85
3525D67C59A8CF631C16712C7924282EF93D7F4AC0F01954F888E19B39F3802DDCD7A041E0DD26DB536BE81248D98E4A23607D3E952B784B463A8477
C96ED1CD1B061B1B033405CE07E2F1FF008A86C863FE5D473FF036AF4FD1D4A68D62A7A8B78C7FE3A2B3F57F0AE97ACDDFDAAFD657709B170F80A3DA
B52CAD56CED9604924755E8646C903D2981668A28A047927C5DC1D76C70727ECF823FE046AFF00C56B6925D3F4BBB4C948D595BDB216BB0D67C2DA56
B5722E2FE3769000010D8C5685C69F6D75A79B2B94F360281086E49C77FAD21DCE07C5CC07C2BD157BB2DB81FF007ECD3B4EB59AF7E0F3456EACF2E1
99557A9C49CFE99AEAAE7C29A55D69B1585C2CD24109CC61A43F27B0F4157B47D22DB46B6FB359993C9CE423B642FD280B9C0FC33D2ADB50D1EFE3BA
697FD685788360631DC5746AFA378674ED56E34DB76658F6F98839566C1C01FD6B4C7876C52E679EDCCD6E67C798B0BED56C7B62AD7F665A7F67B589
88181C618776F727D680383F0C4EFA8F85F5BBA6315B40FBBFD1EDD40E71DC9F5E955FE11C8915DEA91C8CA8CC230A091CF2DC5763A6783F48D31E53
024AC92FDE8DDB2A7F0A7D8784F48B0D48DFDBC2E26CE465B81F41401E7DE25E3E2DC0D8E05D5B13F4012BA1F1668F1F897C4BF62B7B8105D5A5B866
2DC6FC92401F9D7477BE17D26FAF1AEE785BCF2C1B706C106AD5CE8D617371F68788ACF803CD5621881D06680B9C1685A4EADE1CF18D85B3CC644BB5
73200DB86D1DFDBFFAD5E9D546DB4BB6B6B97B95DF24CC31E648DB881ED57A9882A399B6C2EC4E30A4E7D38A92A1BAB74BAB678242C11C60ED383401
E4DF0C496F1B5D3F5060939FF812D278A9449F15E043CAB5CDB023FEF815E8BA6F8634AD32F45DD9C2F1CB8C643751EF4B73E19D26E6F16EE4B6FDFA
B06F3371DD90720E4D21DCE43E3095FB3696A7EF6F931F4F96BB6F0EFF00C8B7A67FD7AC5FFA08A8751F0DE9BAA386BF4926DA30A1A424274E9F955E
D3ECA2D3ED56DA0791A35FBBBDB7103D07B5302DD1451408F31F8B763752CD6776885A08E32AC47639AEE9E1B2BBD1216BC8926852257C37B0ABD756
B05E40D05CC4B2C4DD55BA1AA89A3DAA42216F324857EEC6CE76AFD00A0679EFC4EB39AF2D34ED5E085C5B8842B291F733C8FE75A9E2E9037C2EB20A
BCC915B80A3B702BB99EDA0B9B736F3C4AF130C14238C5511A158E205915E58EDCE628DDB2ABF8520B9CB7852CA6D33E1BDF35DC6434D1CB285E8769
4007F2AA7F079716FAA13D77C7FC9ABD0AE6DA2BAB57B699774522ED65E9C554D3F45D3F4C575B283CAF3061886233401E63A437DA3E2BBC8C338B97
3D3BFAD5AF8B526FD66C208F9711723FE05C577D1786F4985C4896B8973B8CBB88627DCD25F786B4AD42E05C5DDB9925030199C93FAD00627C42CFFC
204B9EB98BF952780C84F8759F45989FCCD74B71A469F736B1DB5C5B2C90C5CA231271FAD3A1D2EC60B16B286DD52D9B398D4902803CE3E122AFF6CE
A4C506E10800F7037723F97E5506B2EAFF0014ED484C6CBC8F2E3BE18715E9767A1E99632F9B676690C9FDE527FC6A06F0CE8AF71F686B1432EEDDB8
B3673EBD68B05CE53E2EAC8FA6D8048D8A8958B11DB8156B5ABD863F86088AEA5E4B58D02A9C9ED5D85C59DB5CC1E4DC411C918180A5781F4F4A863D
1F4E8AD5AD52CE210B755C6680391F84E42E85711B32EEF3890B9E718EB583E1693ECBF11EE67BB6F290ACBCB715E9D65A55858396B3B58E162304A8
EA2926D274EB8BA5BA9ACE17997A395E7FFAF405CE27E22EBB7B1DA5ADBD8AC915BDC287794A60E0F6E6B2BC70964745B04D3A57BA911332CA58BF1C
77EC735EA179636B7D188EEA159554E403D8D3069B622D4DAADAC421618281719A02E607876EA3BBF87E915B386962B131B0C746D86B94F8657B1694
7586BB66EB1A845196620B6703F1AF4FB5B582CE110DAC491460E76A8C5578346D36DEE3ED10D9C4B2939DC050173CEFE2D4A5EFECA18803888B3FCB
C8E78AF41D00C6742B158DC304811091EA140352DC695A7DCBB4971670C8CDD4B2E735620821B688456F1AC683A2A8C0A02FA1E65AEF866E97C6D041
A7EE8E0BCE5993A2A7F103FAD7A0DF69E92E87369D6E022987CB41F871573CA8CCA2531A1900C07DA3207A66A5A2C173CD3C01A75BE9FAA6A36FA945
1ABAE0279CA3B75EB5DBD8CEB3DF4C6DAD235B745DBE7000176F4FA55A9AC6D27904935BC6EE3B95E7FF00AF53A22C6A11142A8E800C0140363E8A28
A6220B9B88AD60696660AA3D7B9F4AF21F0A472A78DEDEFAE94C519795C92B80A0E71FA9AF62655752AEA181EA08C8A61B5B761868223F54140CE07E
23DCBEA715B69BA62B5C481CBBEC048FA66B4742D1AF3C37E0BB848B26FA45F30ED1CA938FE5CD75F1C3145931C4884F5DAA0549482E795F87421D07
521736D2DCEAB2920332924D5AF87065D3750BCB4BAB79D659D81036FDDF735E8D1C51459F2A344CF5DAA06694468242E11439E0B6393F8D160B9E5B
AEBCF71F10A0BB11B35BC53200F8E302B47E29092E23B5B6B60F2483259546700D7A0F969FDC5FBDBBA77F5FAD0634670EC8A587424722985CE3AFE7
BEB5F00DB7F6597F352254664192B818C7B7D7FC6B1ACE0137815ACADAD2493517CB48FB49239E5B35E96AAA8B8550A3D00C508888308AAA3D00C520
B9C1FC3A824B6D3EFF0049B8475777793711C0042AFE7D7F2A6E8BA3BF877C5F2DC5D37EE2584A8931903907FA577C00049000CF27DE8655752AEA18
1EA08CD160B998BA84B757D1476281ADD72669581031D80AE066B79A6F8811EA0236FB12DD093794C7E15EA400000030076141008C100FD68B01C0FC
49867BEF22DAD6377298DD85C81CFF0085765A4606956A06788C039F51C1ABB4531051451401C878E61BB985A05595ECC3832247FC479EBFA562F8BF
4E7D5349B7934EB1961B7863C72B8EF8C1AF49A29587739DD3EE1C783B06DE48CC50792AADD5B000C8FCFF004AC2F8736371617FA935DAB06B908C9F
29EDBB3FCC577F45023CCAF7C41A441A83DD5ED8DF36A4E482607DBB46718EBED567C3FE22D2EDEFBF77A45CC524ED86B991F71E4F4EF8EB5B9A668F
A69699F54B2B7378246776913B1231CF4EF566FF0048B3BD8E2B6B6B48E24DE1D9D23D8001E98001CF6A066578C741B9BBBC8AFADD7CC50CA5D40E57
15BF26A0F3C0A9A746649D8609DB858FDCE7F4AD5A2988E63C6363757FA525ADBA492CD839654E3A564A437B078421D25AC6E45C0381C11BABBDA295
80E4740D065D0BC3B74A2232DDDD64BAAF3B72318F5AA3E0BD2EF748D52EE6BAB39C2DD155076FDDE7A9AEF28A63B983E2EB5BCBCD19E3B2723BB85E
A476AE7D2CEFEEBC2B1E9B169CCAC8C048E382DFE735DF51408E2FC2769A8E9FA34FA59B47898239129E3E63D07D7FC4569DAE8B6CF681F54D3F7DC2
A90CC0EEDDF91AE868A00C8F0FE8C9A35BCD1A602CB26F0A3A28C74AD7A28A0028A28A0028A28A0028A28A0028A28A0028A28A0028A28A0028A28A00
28A28A0028A28A0028A28A0028A28A0028A28A0028A28A0028A28A0028A28A0028A28A0028A28A0028A28A0028A28A0028A28A0028A28A0028A28A00
28A28A0028A28A0028A28A0028A28A0028A28A0028A28A0028A28A0028A28A0028A28A0028A28A0028A28A0028A28A0028A28A0028A28A0028A28A00
28A28A0028A28A0028A28A0028A28A0028A28A0028A28A0028A28A0028A28A0028A28A0028A28A0028A28A0028A28A0028A28A0028A28A0028A28A00
28A28A0028A28A0028A28A0028A28A0028A28A0028A28A0028A28A0028A28A0028A28A0028A28A0028A28A0028A28A0028A28A0028A28A0028A28A00
28A28A0028A28A0028A28A0028A28A0028A28A0028A28A0028A28A0028A28A0028A28A0028A28A0028A28A0028A28A0028A28A0028A28A0028A28A00
28A28A0028A28A0028A28A0028A28A0028A28A0028A28A0028A28A0028A28A0028A28A0028A28A0028A28A0028A28A0028A28A0028A28A0028A28A00
28A28A0028A28A0028A28A0028A28A0028A28A0028A28A0028A28A0028A28A0028A28A0028A28A0028A28A0028A28A0028A28A0028A28A0028A28A00
28A28A0028A28A0028A28A0028A28A0028A28A0028A28A0028A28A0028A28A0028A28A0028A28A0028A28A0028A28A0028A28A0028A28A0028A28A00
28A28A0028A28A0028A28A0028A28A0028A28A0028A28A0028A28A0028A28A0028A28A0028A28A0028A28A0028A28A0028A28A0028A28A0028A28A00
28A28A0028A28A0028A28A0028A28A0028A28A0028A28A0028A28A0028A28A0028A28A0028A28A0028A28A0028A28A0028A28A0028A28A0028A28A00
28A28A0028A28A0028A28A0028A28A0028A28A0028A28A0028A28A0028A28A0028A28A0028A28A0028A28A0028A28A0028A28A0028A28A0028A28A00
28A28A0028A28A0028A28A0028A28A0028A28A0028A28A0028A28A0028A28A0028A28A0028A28A0028A28A0028A28A0028A28A0028A28A0028A28A00
28A28A0028A28A0028A28A0028A28A0028A28A0028A28A0028A28A0028A28A0028A28A0028A28A0028A28A0028A28A0028A28A0028A28A0028A28A00
28A28A0028A28A0028A28A0028A28A0028A28A0028A28A0028A28A0028A28A0028A28A0028A28A0028A28A0028A28A0028A28A0028A28A0028A28A00
28A28A0028A28A0028A28A0028A28A0028A28A0028A28A0028A28A0028A28A0028A28A0028A28A0028A28A0028A28A0028A28A0028A28A0028A28A00
28A28A0028A28A0028A28A0028A28A0028A28A0028A28A0028A28A0028A28A0028A28A0028A28A0028A28A0028A28A0028A28A0028A28A0028A28A00
28A28A0028A28A0028A28A0028A28A0028A28A0028A28A0028A28A0028A28A0028A28A0028A28A0028A28A0028A28A0028A28A0028A28A00CCB2D4CD
CEA97963240637B6DA73BB2183569D735631ABF8E3580D920DB42319E0F1517872F26D46E75AD36EB7C96F693F971B9241C73C67DB0295C763AAA2B8
3D3E5B8B7BCD53C3D733DC35D3CAB25AC8D2139889F5ED8FEB5DB5B40B6D6E90A16210632C724FBD098344D451453105145140051451400514514005
145140051451400514514005145140051451400514514005145140051451400514514005145140051451400514514005145140051451400514514005
145140051451400514514005145140051451400514514005145140051451400514514005145140051451400514514005145140051451400514514005
145140051451400514514005145140051451400514514005145140051451400514514005145140051451400514514005145140051451400514514005
145140051451400514514005145140051451400514514005145140051451400514514005145140051451400514514005145140051451400514514005
145140051451400514514005145140051451400514514005145140051451400514514005145140051451400514514005145140051451400514514005
145140051451400514514005145140051451400514514005145140051451400514514005145140051451400514514005145140051451400514514005
1451400514514005145140051451401972E8C9F6D9AF6DA7920B89976BB0E411F4A4B7D3C68BA5DC0D3A1F3EE08321DE7995FDCD6AD1401CFE916B35
CEA4BAB5F5A986E45B081B70C12D9C9207A5741451400514514005145140051451400514514005145140051451400514514005145140051451400514
514005145140051451400514514005145140051451400514514005145140051451400514514005145140051451400514514005145140051451400514
514005145140051451400514514005145140051451400514514005145140051451400514514005145140051451400514514005145140051451400514
51400514514005145140051451400514514005145140051451400514514005145140051451400514514005145140051451400514CC3FF7D7FEF93FE3
45003E8A28A0028A28A0028A28A0028A28A0028A28A0028A28A0028A28A0028A28A0028A28A0028A28A0028A28A0028A28A0028A28A0028A28A0028A
28A0028A28A0028A28A0028A28A0028A28A0028A28A0028A28A0028A28A0028A28A0028A28A0028A28A0028A28A0028A28A0028A28A0028A28A0028A
28A0028A28A0028A28A0028A28A0028A28A0028A28A0028A28A0028A28A0028A28A0028A28A0028A28A0028A28A0028A28A0028A28A0028A28A0028A
28A0028A28A0028A28A0028A28A0028A28A0028A28A0028A28A0028A28A0028A28A0028A28A0028A28A0028A28A0028A28A0028A28A0028A28A0028A
28A0028A28A0028A28A0028A28A0028A28A0028A28A0028A28A0028A28A0028A28A0028A28A0028A28A0028A28A0028A28A0028A28A0028A28A0028A
28A0028A28A0028A28A0028A28A0028A28A0028A28A0028A28A0028A28A0028A28A0028A28A0028A28A0028A28A0028A28A0028A28A0028A28A0028A
28A0028A28A0028A28A0028A28A0028A28A0028A28A002A189651348CEDF213F28A9A8A0028A28A0028A28A0028A28A0028A28A0028A28A0028A28A0
028AAD34CD1DC409B09490905803F29C71FD6ACD00145145001451450014514500145145001451450014514500145145001451450014514500145145
001451450014514500145145001451450014514500145145001451583E2F8A55D0AE6F6DAE2482E2DA32E8CAD81F95006F515CA78474C96482DB5BBB
BE9A79E7872558F009AEAE800A28ACED6ADE69AC647B6BB96DA589199593A1E3B8EF401A345719E0E379AF68115F5F5FCE66599806071C0A9B43D6AF
23F155E787F5093CDF2977C1291F330EB83F81FD295C763ADA28A298828A2B95F1DEABA8E8DA6C37BA7CA8ABE688E456507AF20FE9401D551581AB49
7B6DA20FB3DE31BE970235C0F998F380319A974EB0D523104D7DA93B48066489572A7DA95C66D514561F8BF563A3F87AE2E91CAC9C2C6475C9A62372
8AA5A4DEA6A5A5DB5E4672B3C61BFC6AED00145721E2AD6753D2B5AD32DAD2788437D26CC346094C601E7F1AD5D73FB4ADB41966B4BC1F688119D98C
43E7C73D3B52B8EC6D515CEF82753BCD63408EFAFA45691D8AE15718C1AE8A98828AE73C67AF3F87F4EB7B88B1BA4B855208CE57A9ADF8668E789268
98346EA195877068024A29AD9DA76900F62466B90BAF15DCE8FE271A76B0B00B4940293C60FC99E80D00763456078A2F352D3F4C6BED34C4EA800646
19EA71B81FC69FE1FB5D6ADE59CEAF74B3A3852983D0E39A571D8DCA2B0BC557B7FA56933EA362D115850168DD7DFAE7F1A8740BED6754D12DEFCBDA
66742C14A918E48FE945C2C7474561F86B53BDBF4BD8B518D12E2D6731909D3A66B7298828A2B3AE35064D5ADEC2040EEEA5E53FDC5F5A00D1A2B3B5
AD562D26CC4F20DCEEC238D33F798F4A19353F23CDF3E3F340DDE508FE53FECE7AD0068D1595E1FD6A1D6EC4DC46A51E3731CA87F85856AD00145145
001451450014514500145529EF7C8D420B790288E7042313C971DAAED0014514500145145001451450014561F8AEF2F74CD1AE351B39D54C0A094640
41E71D7F1ABDA35D3DEE8D65772E37CF0248D8F52A0D005EA28A6B6769DA4038E0919A007515C9E81AC6A777E27D434CD4648956D065555002E0F19C
D33C4DADEA1A6F88F4DB286E122B5BC20331404A9CE295C763AFA2B17C45717B67A7A3D84C1676916350E80EE27F0F6AD2B48E78EDD05C4A65971966
200E7D06074A622C51451400552BFD4ED34F7B74BA93635CC9E5C6319C9A8F5BBB92C347BABA81419634CA03EBD05634504F6FA7C5ACF882EC3CB6CA
660AA9C2E470290EC7554572BE08BED5F5682E352D4A41F6699B16F1EC00800F5AE967595A32219046DFDE2B9A6225A2B89F09EA9AAEB77BAA2CF782
3167308D42C60E7EF0E9F85685BDEEAA9E2B86D27B88E7B29A066474503247D295C763A6A29AE094215B6B11C1C6715C7D86A1AD5FF8B750D25AF516
DAD1412EB100C73D2988ECA8A6A8214024920753DE9D4005155750B9167613DC9FF9668580F53D87E758FE0AD68EB5A1AC92BEEB8898A480F5EBC13F
E7B5007454515CC5D7F6C6A3E259ADACEF4DAD85BA2EF655C92C7B0A00E9E8AE2F57D4EFE5F1259683A35D92E8BBEEA53CED19E73EF8AECC0C0C139F
7340585A2A29D5DE1658DCC6E470C06715CA780F50BFBFB1BEBFD4EF0C8893346AA7A285E49FD680B1D8515CDD95B6B179772DC36AB22596EFDD0F2D
72C3FC2BA4A0028AC4D67567B7D42CF4BB36517774DD48CEC41D4E2AFDD9B88ACDA4865DD246A5882A30F8EDED401728ACFD17538B57D323BC886031
2ACBE841C1AD0A0028A2B365BE63AC45616FB4955F326247DD5EC28034A8AE77C66DA8C1A33DDE9972F1490F2C07422ADE9776973E1C8AE8CEE43C05
9A4272C38E4FE14AE3B1AF45711E01D46F350B8D41751BD9E49EDDF6AC2E780A7A1F7AD186C66BAF134ECD7B72D6B66A81537FDE73CF5F6A2E163A6A
2A2B8B886D937DC4A91AFAB3629619E29D37C322BAFAA9CD31125151BCD12384795158F452C0134C8AEEDA699E28AE239244196556048A009E8A88DC
40A4833460AF072C38A7492C7128691D50138CB1C73400FA2A332C62511191048C32149E48FA50F2468555DD54B1C2827193ED4012514C4749065195
874C839A1E48D387755FA9C5003E8A88CD12C823695039E0296193F853F72EEDBB86EC6719E71400EA298F2C71E3CC755CF4DC719A5475719460C3D4
1CD003A8AA135DB5B6A90C52B0F26752138E8C31C55FA0028A28A0028A28A0028A28A0028A28A0028A28A0028A28A0028A28A0028A28A0028A28A002
8A28A0028A28A0028A28A0028A28A0028A28A0028A28A0028A28A0028A28A0028A28A0028A28A0028A28A0028A28A0028A28A0028A28A0028A28A002
8A28A0028A28A0028A28A0028A28A0028A28A0028A28A0028A28A0028A28A0028A28A0028A28A0028A28A0028A28A0028A28A0028A28A0028A28A002
8A28A0028A28A0028A28A0028A28A0028A28A0028A28A0028A28A008E5952289E4760A880B313D80EB581E25BCB6BFF05EA73D9CCB2C5E530DCBEA0F
35B77B6C977653DAC84849A36424760462B060F0E4B0F86868A9242227CF9AFCE5B273C71498D173C1E08F0AE9C1BAF923F99ADAA86DADE3B5B68EDE
15DB1C6A1547B0ACBB6B2D553C473DCCB79BF4F75C2424F43F4A005BCF12699637E6CAEA492394639287073E86B46F413637000C931363F2A74B6D04
C434D0C7211D0BA838AAFA99BC36CD1D944AEEEA41666C6DE2981C77C345BF6F0BBFD9A48553CF6C6F0720F19FE95D2691A04763A85C6A33C9E7DEDC
1F9A423181E82A8782B47D4B41B23617491B42CECFBD5F9048F4FC2BAAA481B0A2B1EC22D5D75EBE92F2406C18016EA08E3F0AD8A620AE43E2764F85
D0019CDD47FD6BAFAE5FC6FA56A5AD584367611AE1645959D9C01919E3140D6E62EAF79A9687AFE95AB6A0564D3DE31090BFF2C7239FC7DEBBE8A58E
6896589C3A38CAB03904564EA3A7BEA7E197B3BAB70667882ECDC0E187420FEB59FE11B7D774AB04D3F52B5578A33B6391645CAAFBF3401D5572DADB
59EA1AC4BA6DFBA2DBC7684FCDC8DEDD0FE02BA2BA964860692281A771D114819AC6D0D6EDA7BB1AAD81596494BAC84061B78C0CFB52608C3F85FA86
6CEF34891F7359CA4C7FEE13FE3FCEBBCAF3C8B4FD56C7C6EFABD86932ADA4A36CB1E40CE7A91FCEBBF8246961491E368D9864A3751420670DF10C48
7C41E171085321B96DB9E99DD1E2B7B5D93564D1B50668ACCC62DDFEEBB6E0369CF6ACCF19699A9EA3A9E91756366645B19FCC3F3AF232A7FA569EB9
3DD5DF872E62874F9CCD7313C423C0F97823268033BE17316F08267F86671FCABB0AE3FC01697FA4697FD9D7D652A13234824C71CE383F95749A85C4
D6D68D25B5BB4F2F4545F5C500F7398D712CB5A9F5AB6B8285AD6D4A41BBA07C1248F70714EF867AA0BFF0D2DAB3E65B2731B03D769E54FF0031F855
BF0CC53C3632FF0068E9F22DD12F2B394077E7B03EB8E2B1B44B0D4F44F156A37906993369F799211700A9CE471F89A0677F5CCEA9A741AD6A3A8E99
70A30F6C8C8F8E50E4F22BA3CB7979DAD9C676F19FA7A57371CBA80F133DF1D2E716CF6E223D0B6436734C48C4D33519D347D47C3BACBEDBCB35FDDB
31FF005880E473DFB577F1B6F8D5BD4035C9F8DBC38DA90B7BFB246FB640EA1827F1A679CFAD759102B12023042818F4A481987E3938F066A9FF005C
7FA8AE6B4B92F97C1FE1E8C451FD99E750EE1CEE0371238C7D7BD753E2DB69EF7C3779676B134B34E9B140FAE727F2AABA043369FE13B4B4BDB29649
ED860C4AB9C9049041E98E94305B1BD1C304534AF122AC9290D263B9E809FCAA7AE6F401A95C6B9A85FEA16B25B46CAB1C28DFDD1FD6BA4A620AE5FC
3331BCF10F886E2403747722DD48FEEA8C7F857515CAF873367E28D7AC25C0334C2EA3FF006830E7FA52634667C43F31B5DF0EA64F9467391DB395C5
77B593E20D21356B244E04F0389616F461FE348353BA31F906C6617BB718C6533D376EF4CD311CB7C3D674F1278820DC4C6250DED9C9A99DF513F11A
5D323BC94D9DC40659549FBA31DBD39205743E1DD11348B694B6D6B9B87324AE0773DAB166B5D422F1E3EB29612BDB080C1F2F56F7FA7148647793EA
7A0F872D7FB5A6925CDEED94C67737918271FA54DA7EA16F7DE2CB26D3751371686090988B67CB238C73CD5FF125AEA371A86952E9F18956095DE546
3852300609FCEA8AE8F6D71E27B4BFD3ACE5B296DD8B5C129B15C1ED8A00B33BDCCFE20BDB7D49DA1B25881B6D8F8DC7B9C773D6A0F87BA8DC5FE997
6B732BCBF67B864467EBB68B3B4D4ADFC47A95DDD5B497124ADB6D24FE08D79FCBB555F0AE89A8585DCF6F7913228BBFB4F9ABF71F83C0A00B3E23B9
BDB1F11E95BEFA68F4DBB90C7204382AD8C0191DB2454D67617727886FA35D4EF3EC7004E3CCE4B9E48CFD31F9D4DE37B35D43C3175E5BE25B722446
53C8653D2AF7872DA5B7D1ADDAE5CBDCCCA2499CF566228EA1D0A7E347FB3E882F940F36D678DD1BD3E600FF003ADF46DE8AD8C6466B9CF1BEEB8D2E
0D322199AFAE11140F407713F862BA35015401D00C51D407565EBAF7696B11B238FDF279AD9C623CF3CD6A56078A2C2F2F8589B5CB450CE1E78C3637
AD0C48CFBBD524B4F881A6E9F04B2986EA06F35598B29E18A919EFC53BE215F5D58F87E2BBB0BA785BCF45DD19EA0F3FD2A2F11E957F278BF4DD5ADE
269A1851A2654EAB9079FF00C7AA2F156897F37852CB49B281AE1E2915DDB3D319C8FD68289BC4CF7BA3E9763A8ADFCEF70B711AC809F95F77518EC2
B4BC65AC49A3787E5B9B73FBE7611C647627BD51F1ADADFEADA3DADB59DA4864F352671FDD033C7D6AF6B1A74BE20F0DF912466DA7243AABFF000B29
E3F3FEB408CBF155BC907C3CB93E748F2491A3CACE49DC4919ADEF0C023C33A603DADA3FFD045676BF677DAC78425B28A064B9DAAA55F8DD8EB8AD6D
12268347B581E2689A28C26D6C67818A03A1A145145311C56B7FF128F1EE9BAA31C417886DE463C007B6693C5907F6A7852EB5145C490CDE6C2DDF62
9C71FA9AD3F1BE8CFAD686208A32F324AAEB8EA0679FD2B63EC512E97F615553188BCB008E3A6290EE61E9B7E9AEBE932C643A45179F2FCD9C3E3033
EF9CFE75D3D731E07D15B45D2658E58D966926624B0EA074FC2B7354378BA6CE74E0A6E827EE8374CD080B74552D27ED9FD976FF00DA38FB56CFDE63
D6AED3110DC5BC57313453A0746C654F7C1CD713E249DFC4BE2287C336A48B58184978EBEDFC35DA5F34E965335A2079C21F2D4F735CF785F48BAD17
4F9259ADFCFD42EE4324E77818F6C9FAD26346CC5716B67756FA5431ECFDCEE40B8C2A8E2B42B9CD1F4ED48F886F354D55635DC8238115B3B57FCFF3
ADD99E448C98E2323F65040FE6684079C784D3517D3BC4D2D81894C933805F3BB2013C76E86BB7D0D238F44D35A708B28811559800724741F5AC7F08
697A9687677D15EDA8904F299804753C9001079F6A75C5B6ABA9788B4C2D6A6D2C2C8F9A46EFBC7B0FD3A50367595C4F877E5F88FE205F58D4F3F515
DA31210955DC40E067AD723A2697ABDA78B352D52EED57CABB5C00922923078EF4091D36A37B169D6135DCF9F2E25C9C75A96095668239933B6450C3
3D70466B98F16DAEB7AB431E9D6369B2DA4606591A419C0EDD6BA7853CB8638F8F9540E3D853031F5CBD8D353D3ECDE40B1973349F2E785E83F139FC
AB97F0EDD47A57C42D434E5205B5FE5E238C027A81FCC5747611EA035EBAB9BCB2FDDCC5523F98108ABDEB3BC59A46A577AFE9DA8E9B6F96B23B89C8
F9F9CE2931A3B3AC5F11EAB0787F4AB8BF28A667C045EF23E302B56DDE4920479A3F2DD864A6738AE535FD06FF0057F13DADCCC825D32D4065883805
9A9890783B4D6D2B4BBAD6F55CFDB2EC19A563D553A815D4D9DC2DD5A457080859543007B66B03C4D6FAD6A9A63D8D8DBC702C980ECD202719ED5D05
AC22DAD21806311A05E3D8520648D9DA71C9C5798F856CF51D43C27AC59DA4DE4C6B74DB0A8CB3B6012BFA0FCEBD22E8CC20616F18776E065B18E3AD
73DE0CD2750D0A0B8B3BA891A39676984CB203C900631F80A01167C21AB8D4B4A58A65F2EEED408A68CF504700D7415C71D1357B4F17CDAB69C224B6
9B89222DFEB3DEBB019C0C800F7C50819C3DE2BA7C55B573928F6D81E80E2BB4B86D96D2B1E8109FD2AA5F582DCDC4176981716E7285BA11DC1A4D42
2BBBC85ADA30218DC6D924241241E085FF00EBE280303E1AC1245A0DC4AFC24F74EF18F45E9FCC1AEC2A0B5B78AD2DA3B78176C71AED5153D30615CF
787879DABEB174DF78CFE5AFA8515D0D73FA4C2F61E20D46DDF88EE584D19FE63FCFA52606C5EDBADDD9CD6EE0112215E7A5713E18949B07F0EC92AF
9F15D10C33CF960E4FE1C577D58565E1E8ADBC4B79AC332B493AED403AAFAD0D023135175D07C6DF6EC08EDEEE06C8ECCCA3F9D74FA242F1E9AAF37F
AD9C9964FAB7FF005B1553C4DA20D6E3B353B4F9138720FF0077BFF2ADA5508A15460018028481B39BF10220D7F4C9AE2746857205AE37348C4F040F
CBF2AA76738B3F883A844488ADC5A798CB9E074E7FCFAD6B5DE8B24BE2387578E552628CA08DF3807D6AA3F8664FEDB9F505BB5737717953EF5C1C64
7DD1F850067FC488639AC6C245DBB9E754DF9EAA7D2B760D12D74B7BCBCB42524922239E8B81D7F4A83C43A1DD6B32DAEC9A18A3B69048808393EC6B
62E2DCDCD8BDBCAD83226D665F5F5A02E79FDE6AD34966969A969B2DABA48A5678C6030CF24D761ACC5606CA0BEBD948B5B41E6819FBDC7153CDA734
9A5A69CCCB2465763BB8E71EC2AA6BBA136AD0DA5A79E23B38594BA7397C741405CCAF0C249757F2788B547F29EF9BCBB38493F2A76FC4FF009EB5A7
786DE2F134324B3B4D24916C5B60BBB6F3D7D073DCD4CFA54F2EB16972F2C62D2D01F2A1507838C0FD295B4A9A3D66E352B5923F3278C467CCC9D98F
4A00C6F0B5D326BDAFC570FE5C30DC2A2216CA8DC4E3F9559F17DBAC02DB540B910CCBE6AE4FCC3A5161E166B3D42E263742586E6559650C08662B93
FCCD6F6A5682FB4F9AD49C798B8CD01D4E2F559187892D7568917FB3E37113E7804E393E98AEA34674BB33EA2A491336109FEE8FF3FA5324D06093C3
EBA4972AA00CBE3393DCD6959DBC7676915B420048942802806CC679E3875E9D54B5DCF24631080331F1EBD2B2FC177CA92EA71DCB98D4DDEC8918E7
1C74CFE35BF67A3AD9DD5EDC24CC5EEC92491CAFE351693E1DB5D377157694997CD05BB350043E2D2D1DAD9CE98DC974833EC6BA0AC9D52DFED9A858
418CA23996418E303A7EB5AD400514514C41451450014514500145145001451450014514500145145001451450014514500145145001451450014514
500145145001451450014514500145145001451450014514500145145001451450014514500145145001451450014514500145145001451450014514
500145145001451450014514500145145001451450014514500145145001451450014514500145145001451450014514500145145001451450014514
500145145001451450014514500145145001451450014514500145145001451450014514500145145001451450014514500145145001451450014514
500145145001451450014514500145145001451450015997BA52DC6AD69A8C5279734070DC7DF43D4569D14005145140051451400514514005145140
1836BA01852E2179CB433CDE6BF2496F6F6ADD0000001803A0A5A280283E9CB26AE9A8C9212628CA468470A4F5357E8A2800A28A2800A28A2800A28A
2800A28A2800A28A2800A28A2800A28A2800A28A2800A28A2800A28A2800A28A2800A28A2800A28A2800A28A2800A28A2800A28A2800A28A2800A28A
2800A28A2800A28A2800A85EDE27B88E765FDE460853F5A9A8A0028A28A0028A28A0028A28A0028A28A0028A28A0028A28A0028A28A0028A28A0028A
28A0028A28A0028A28A0028A28A0028A28A0028A28A0028A28A0028A28A0028A28A0028A28A0028A28A0028A28A0028A28A0028A28A0028A28A0028A
28A0028A28A0028A28A0028A28A0028A28A0028A28A0028A28A0028A28A0028A28A0028A28A0028A28A0028A28A0028A28A0028A28A0028A28A0028A
28A0028A28A0028A28A0028A28A0028A28A0028A28A0028A28A0028A28A0028A28A0028A28A0028A28A0028A28A0028A28A0028A28A0028A28A0028A
28A0028A28A0028A28A0028A28A0028A28A0028A28A0028A28A0028A28A0028A28A0028A28A0028A28A0028A28A0028A28A0028A28A0028A28A0028A
28A0028A28A0028A28A0028A28A0028A28A0028A28A0028A28A0028A28A0028A28A0028A28A0028A28A0028A28A0028A28A0028A28A0028A28A0028A
28A0028A28A0028A28A0028A28A0028A28A0028A28A0028A28A0028A28A0028A28A0028A28A0028A28A0028A28A0028A28A0028A28A0028A28A0028A
28A0028A28A0028A28A0028A28A0028A28A0028A28A0028A28A0028A28A0028A28A0028A28A0028A28A0028A28A0028A28A0028A28A0028A28A0028A
28A0028A28A0028A28A0028A28A0028A28A0028A28A0028A28A0028A28A0028A28A0028A28A0028A28A0028A28A0028A28A0028A28A0028A28A0028A
28A0028A28A0028A28A0028A28A0028A28A0028A28A0028A28A0028A28A0028A28A0028A28A0028A28A0028A28A0028A28A0028A28A0028A28A0028A
28A0028A28A0028A28A0028A28A0028A28A0028A28A0028A28A0028A28A0028A28A0028A28A0028A28A0028A28A0028A28A0028A28A0028A28A0028A
28A0028A28A0028A28A0028A28A0028A28A0028A28A0028A28A0028A28A0028A28A0028A28A0028A28A0028A28A0028A28A0028A28A0028A28A0028A
28A0028A28A0028A28A0028A28A0028A28A0028A28A0028A28A0028A28A0028A28A0028A28A0028A28A0028A28A0028A28A0028A28A0028A28A0028A
28A0028A28A0028A28A0028A28A0028A28A0028A28A0028A28A0028A28A0028A28A0028A28A0028A28A0028A28A0028A28A0028A28A0028A28A0028A
28A0028A28A0028A28A0028A28A0028A28A0028A28A0028A28A0028A28A0028A28A0028A28A0028A28A0028A28A0028A28A0028A28A0028A28A0028A
28A0028A28A0028A28A0028A28A0028A28A0028A28A0028A28A0028A28A0028A28A0028A28A0028A28A0028A28A0028A28A0028A28A0028A28A0028A
28A0028A28A0028A28A0028A28A0028A28A0028A28A0028A28A0028A28A0028A28A00299248912EF91D557D58E05617887C4D0690A638D04F718FB80
F4FAD79C6B5AA5DEAD279B7173315E408D5B007F8D439A4691A6D9EAB71AF6936E7136A16EBFF03CFF002A9AD753B0BDE2D2F2094FA2B827F2AF11F2
61F35B7460A021403E87BD35AD1E197CE88CA98FBA633823D0E697395EC8F7CC515E4DA1F8F352D35D6DF55CDE5B81CBE3E741EE7BD7A6E9BA8DAEA7
6897567289226EE3B7D6AD34CCE516B72E514514C90A28A2800A28A2800A28A2800A28A2800A28A2800A28A2800A28A2800A28A2800A28A2800A28A2
800A28A2800A28A2800A28A2800A28A2800A28A2800A28A2800A28A2800A28A2800A28A2800A28A2800A28A2800A28A2800A28A2800A28A2800A28A2
800A28A2800A28A2800A28A2800A28A2800A28A2800A28A2800A28A2800A28A2800A28A2800A28A2800A28A2800A28A2800A28A2800A28A2800A28A2
800A28A2800A28A2800A28A2800A28A2800A28A2800A28A2800A28A2800A28A2800A28A2800A28A2800A28A2800A28A2800A28A2800A28A2800A28A2
800A28A2800A28A2800A28A2800A28A2800A28A2800A28A2800A28A2800A28A2800A28A2800A28A2800A28A2800A28A2800A28A2800A28A2800A28A2
800A28A2800A28A2800A28A2800A28A2800A28A2800A28A2800A28A2800A28A2800A28A2800A28A2800A28A2800A28A2800A28A2800A28A2800A28A2
800A28A2800A28A2800A28A2800A28A2800A28A2800A28A2800A28A2800A28A2800A28A2800A28A28012B3F5BBF1A66953DD632EAB841EADDAB46B8F
F887201A6DAC479324A405CE335327645415E491C5F9AF7D31B99D8BCADCB127A93FE1CD35E018DCA7201E0761F5A618630A06CE00C1C1EF53C70CA4
A2464927AF702B0B9DA90CFB2BB90A5371E48E2AC8D3988180D9C63039AD8B3B668813202C3A7D3A55A78D083F2E39E091CFE153CC3B1CE5D69F1AA8
F3620393CEDEB56FC0D3CBA4EBE6D59CADBDC641423A1EDF4AD1910C61B6A6E241CB03C9FA8ACAB85104E92C0423A91C74E73D69C6762271E647AB51
50DA4BE7DAC329209740C71EE2A6AEA38828A28A0028A28A0028A28A0028A28A0028A2A3918AA16556623F857193F9D00494573967E2DB7BAD526D34
585E8B987EFA84076FEB5D1D0014514500145150DB5C45730ACB0B6F46CE1BA67071401351451400514554D42FE0D3AD4DC5D36D8C305E3B92702802
DD14514005154351D5F4FD2C037F72B0EE1919079AB514A934292C6DB92450CA7D41E45004B451450014514500145145001451450014514500145145
00145145001451450014514500145145001451450014534B004024024E07BD3A800A28A2800A2A369110AAB3AA9638504E33F4A92800A2A8586AD61A
84F3C36970B2C96E76C8A0118FCFAD5FA0028A28A0028A28A0028A28A0028A28A0028A28A0028A28A0028A28A0028A2AADEDFDA69F0F9B7B7090A7AB
9A00B5454504F15CC09340E1E3719561DC54B4005145140051451400514D765452CEC15477271557FB534EF33CBFB7DAF99FDDF3973F966802E5150C
7756F2B6D8A789DBD15C13535001451450014514500145145001451450014515461D5F4E9EF7EC70DDC725C609D8A72702802F514514005145140051
45672EB5A7BEA7FD9AB72BF6AE7E4C11D3AD0068D14514005145140051515C5C456D1196770883A9350586A765A92C8D6370932C676B15CF06802E51
4514005145140051451400514514005145140051541F54813578F4CC3199D0BE40E055C7748D0BC8CAAA3A963802801F453410CA194820F208EF59B7
BAED8596A30D84D237DA2660AAA173D6803528A28A0028A28A0028A28A0028A28A0028ACEBDD552D352B4B230BBB5C9C6E1D17EB5A04800927007526
80168A8E09A3B8856585B723720FAD49400514514005150DD5C2DADB493BAB32C637108327150693A845AAD8ADDC2ACA8C4801BD8D005DA28A2800A2
8A2800A28A2800A28A2800A28A2800A28A8E2952689648DB723720FAD004945145001451450014514500145145001451450014514500145145001451
450014514500145145001451450014514500145145001451450014514500145145001451450014514500644974DE636D76C64E3E63455F36C849396E
7E9454D99A7322C51451546614514500145145001451450025729F10EDCCBA2C538CE6DE75380339CF1FD4575755354B18B52B09AD25E1655C6E1D41
EC69495D58A8BB34CF264E7000C7A2FAD6969D1969448415DBC609EA3D6A9FD9DECEF25B2BA5613C2D8393C1EC1856D59C247DC55C81D2B91DCEF563
4570C7686FBA3BF1F4348D2C80132856EC368CFE3512C6A406D9BBB104D49B9860041B73DEA52B80D0D9385C03FC5B8565DFC2934A2250482DE9DEB4
A691553E63B7E5EE2A4D1ECE4BEBEDE015857058E38C7D7D6AA316DD852928AD4EA34880DB6956D0B124AA753F9D5EA0515D88F3D851451400514514
00514514005145140051451401C1F85D8FFC2C8F11296EA3A7E22BBCAF3DF0D03FF0B4F5D2071B1BF0E56BB9D42E52CEC27B991C22C485B27B524363
352BE874DB196EAE080918E99C6E3D80AA171AEFD8A08AE6FECDEDED64655F37786D99E8580E83A5713E34BBBABBF0EE82CEE7CAB890B3BB1CE4E78C
FE15D86B5A4EA3AD688DA749716F1248177328249C107F98A02C4FE27D46EB4CD0EE2EECA0333AA1F98300107F7AB33C177F76DA0582C9673B8909FD
F646DC16EBC55CD417C9F025D47E789CC364E8655E8C5548FE9553C377925A7C3AB4BB8516478A2242B743872280E8759454504A26B78E5C637A06C7
A6466A5A620AE03E2735C88B4D95083669703CC00F56CFFF005ABAED5A79E3804166337329C27FB23B9AE3FE23BDE0F0D5B24D6CB1A473202E24DDC8
1F4A4C6B73BE8DB7C6AE380C01A7D56D3DB7E9D6ADFDE850FE82ACD311CF78EA18A5F08EA0D222B3471EE4247DD391D29747D6204F0DD95C309A4458
543B246485C01D68F1D9DBE0DD48838FDD8FFD0854BE108917C29A7A8030D08278EB9A5D47D0D4B4BA82F2DD2E2D6559627E8CB55E5D5AD63B892DD0
BCD344332244BB8A0F7AE27C257B750CDE27B388958EDC3CB0003846F9BA7E952FC33BBB8934ABFB892092E6792E4969432E5BE51807273EBF9D170B
1DBD95E5BDEC026B590491924647623A83EF566B97F0DDA6A567AFEA8D790A436D78DE7448AC0E08383D3EB5D4531329DCEA16F6D28858B3CC79F2E3
1B9B1EB8A8ACB58B3BCBB96D23664B98BEF4520DAD8F5C5655E59C163E2E1AC4B7A774B0F962D54659B8C640AC89E79A5F8B1A73181A256B365C3756
5C39048FAD2B8CEF68A8ADE78AE6059A16DC8FD0E319AA1E24BC3A7787AF6ED73BA388E08EC4F19FD698891F57B15778D66F3191B6B8452DB3EB5721
9A3B8856685D5E36195653C1AE57E1F891FC16924640B998C877B776C9C134EB0B6BDF0CF81EFBED5379B3C292C89B79DB91C0FCF9FC690EC740350B
4DF222CDBDA3FBE114B6DFAE054D6D7305DC226B695648CF46535CD7C3843FF0894372EC5E4B996491D8F527711FD2AA7876E5E0F1FEB9A720C5B9C4
A107456C0E7F1CD0163A25D774D3AA9D34DC8174A09D8CA5781D793528D56C4DE25A19D44D27080820391D403DEB93F1142B2FC4CD0D4F4685B781DF
EF75A93E28EE4D0EC6688EC952F90238EAB956E87F014058EA6EF53B3B2245D4E23C10092A7033D39C62AC4B3245034C72C8ABBBE519C8F6AC7F13DB
AA783B518880425AB11F50339FAE697C1523CBE11D35E5259CC5C93F534C09B4FF0010695A8A4AD6B740F93F7C30208FF1E9DAA6B3D5EC6F1EE238A6
C496DFEB524054A0C67273DBDEB93F08DBC6DE3AF11481576C72FC800E0124F2297598BCDF8A7A642AC5639AD7F7C14E3CC0BBD803EBC85A416474F6
FAF69D3DF2D9C739133825032950FF00427AD59BED42CF4E8C497B3AC2873866CE3FCF35CAFC48FF0046B3D2AEA1C24915EA2A95E08041E3E9C5741E
27546F0DEA1E622B81031C1FA50058FED4B0362B7AB72AD6EC702450483F97D29DA76A565A9C267B0B849E356DA597B1F4FD6B1BC031A2F832C1410C
19093F893C566FC300134BD490745BE7FE43FC280B0ED47C4B66FE34D3F4F0D22A5B97F34B211F363031DCF4AEB2DEF2DEE51DE27DCA870CC41007E7
5C7DF807E2F69A02E716249F63893FFAD5D45E91A8437FA6C12B453088297C74DE0E31F9500C8DB5ED3D559F7B9895B6994212A3F1AD18268AE2149A
07592371956539045733A544B6BE149F4A4B817F2430C88CF10CA8C83819A8BE1DCE22F08D9A5C3E1DA6923407D724E3F4A02C53D63C490BF8D34CB0
F266115B3BB39284166C10303D3835D6C7AADAC96725DFEF1628FEF168C835CA6AA31F17B46C01CD9B67F296BB08EEE292F66B45CF990AAB3647186C
E31F9500CC7F0DCFA05CDD5ECDA2A013337EFC852327E87A558B8F13E8F6D2A453DE2C6EEC5406522B0BC02C5F56F1202071798FFD087F4A7FC49D3D
A5D1E2D52055FB469F28941C755C8FEB8A3A0753A1B8D734CB58925B9BA11249F75995803FA55AB6BDB6BBB6FB44132BC3FDFE83F5AC1BDB9B7F1069
FA5DBA461A3BF2B3B0FEE22E09FF000FCEAAFC47BF7D2FC2862B6509F689041F28C6D520938FCA80B1BA35DB06F2CAC8C6295B62CDB0ECDDE99AD3AE
5658AE6EBC11158D9E9F2169AD1029DC8141201CF5CF5F6AD8D052F22D1AD62D480FB5226D7F9B39C50069555B8D42D2DA5114D708B2119DB9C9FCA9
F7938B5B39EE1864431B391EB819AE2FE19B9BE8354D5AE0EFB89AECAEE3C90300E07FDF54C476F0CD14F1892175753DC1A8AE2F6D6D885B89E3466E
8A4F27F0EB589E1DB3D4ECAF758B9D48A88A69BCC891307E5FFF00562B2BE1F5DB6B177ABEAB75FBC99A709196EA8983803D2803B48678A74DF0C88E
BEAA7354AFB5BD3F4FBA8ADAEE711CB2B055041C64F4E6B989EF9F4EF8A315AC39F26FA10244CF1BB04EEFAF1517C518DAE06896A87066BBDA3EA703
FAD21D8EBAE357B0B761E75C281BB696C12A0FD7A54F71750DBC26595F08A012402D81EBC76E3AD6178EE248BC0FA8246A155117681DBE61573C371E
DF09D8ACECD287B60CE58E490C33FD6988D5B7B886EA049EDE41246E32AC3A1AA435CD2FEDDF62FB743F68C676138FD7A5627C397DDE1C930C4AADDC
AAB9EC32302B3FC476515D7C4CD1124030D1177E3EF6D2C467F20290EC7611EA565248912DC2EF93EE0208DDF4CD6678E1237F086A46450DB62DC33D
8E6B3FE261687C371DDC476CD6F728F1B8EA0F3563C4F782E7E1DDCDE32E3ED1688FB7D376DFF1A00BDE0F3BBC29A69273FB81561B5AB20D26D7678E
23896545CA467DCD73D71752699F0BE39ADF2B20B450A476CF7FD6A4F093489E06B610584933C91312372E24249E4E4D00753F6983ECDF69F317C9DB
BB7E78C555B6D674DBB99A0B5BC8E695412C8872C3F0ACCF065AEA765A37D97568846F1C8DE580C0FCA4E71C561F8B606F0B6AD1789F4D45C4A7C9B9
87B3E79CFE94058EB6DF5CD3AE6EC5ADBDC79B3775553F2FD78E2B4AB2B438ADE4B55D46208D2DF012C8EBEE3A0F61D2B569886BA2BA95750CA7A823
20D79C68F6566DF15B5380C113451C05910AE429F93B7E26BD26BCE34C88BFC5AD482CB221F259B298E7EEF1CD26346F78E2D2DA3F0E4F730A2C5751
63C878C6D6DD91C715AFE1E9EEAE742B49AFD192E5E305C30C1CD67EBD672DB594BA9C52BDD496A0CE21B820A100738000C7156743D660F1268AF756
9BA2DDBA26CF546C7FF5C51D43A1625D62D231232F992A4470EF126E553EE6ADDADD417902CF6D22C91374615CE786C43A5D94BA54731D427477790C
6BF28CF624F154BE1D4C62D2AF7CF936A1BE648D09E0127181F8D0163A19FC41A5DBDFC7633DD84B891B62AB230C9FAE315625D4ECA131092E5144C7
11B67E563E81BA5721E3AB44BDF15F872DE41F23CA4363838C83D7F0AD0F88AAB1F82EE4A285F2DE2D9818DBF381C50075134B1C1134B2B848D06598
9E00AA8FABD84768974F701607E1642ADB4FE38ACED534F9356F0825A49722DDDE08CB48FD0100139ACDD62E233E02BA86CDCDE2456DB0CE176AF1F5
EBF85303B0560CA1948208C823BD50B8D62D21BDFB1A7993DC8E5A38537151EFE9557C27316F08E9F2B92C45B827D4E2B0FE1C5C1BDFED9BE9803712
DE10CDDF18E05023B18A78A607CB704AF0C3A153EE3A8AE3A510D97C5184AAC70C72D99DC78504F3FF00D6AD8D0F47B9D3B57D52F2EEEC4C6F9C3228
C8DAAB9EDF88AC4BFB75BAF8AD64B22EE48AD4BE0FA8071FAD263475516AF632DD2DA8B80266FBAACA54B7D335A15C778B9BCBF13F871D701BED0413
EDC576340318CCA8A59D82AA8C924E00159A75CB310BCE1276810F332C44A7D7356756B4FB7E977369E6F95E7215DFFDDAC2D32582DBC312E9D613BE
A26D607433050A84E0F19271DFDE803A3B6B886EEDD2E2DE4124520CAB0E84572B7857FE168E9E0280DF639371F5E0D3FE1D48B1F832C3CC3B4C9248
AB9CF277371FA556BC6FF8BB3A7AFF00D393FF0026FF000A00EB6F6F20B1B66B8B990471AF527BFB0F7AA8BAED87962499DEDD58654CC85770F6AC2F
88314AF16992070208EE94C8BDCF23F4A678F4FF006869B6369A70F3AE64B8478CC6376D03B9FCC500747A56B165AB4124B67212B136D7DC3054F5A6
43AED94F1C92C7E6982362AD308C94E3AF358DE3265D0FC217AF6512C72DC61247418E5B827FCFAD1A01BB93C0F696B6BA7152F69B0333A80491F7B1
9CF3D7F1A00E9D4C3776E186C96290647704570DE00BA82D6E75B80962C6EC948D149240CF402B7FC1BA7DF697A1AD96A2E1A58A46C61B76D538207E
A7F3AC3F877121D4F5F9CA8F305D6C07DB24E3F4A00EBED353B5BB9DE08DCACF1FDE89C6D61F8553D4FC4961A4CD1C77E278448FB15CC7F293EB9F4A
C1F134AF63E3BD0EE6DF3BAE33148A3AB8CE29BF13625B9B7D26D587CB35D852DE80E07F5A0123A58F5DD3DEF2DED048C24B905A1CAE038F514F9756
B64BC3668259A7500BAC49BB6FD6924D16CA4D4ED750643E75AC7E5C401C2A8E7B7E359761059E97E24BF2B74F3DD5FBEF3022E4A0F7F41401B163AA
5ADFC92C50BB79B11F9E365C32D5EAE234C99E3F88BABBCB9402D03BA039E98FE95D8DBCCB716F1CC808591430C8C1C1A130684BBBA82D21F36E1F62
6719C13FCAB36CBC4BA55FB482D6777F2BEF1113601F4E95B046460D715A64B1681E30D5AD266296D729F698801C03DC7D79A608E874ED7B4FD4E778
6CE5791E3E1C043F29F7F4A96D357B4BDBB9AD6DD9DA484E1F28401F8D725A2CD2681E2CD45353D90C5A8446E93FD923248FCB3F9575BA4C463B3F39
D024B72C66907BB1CE3F01481A292EB36ADE258F4C7B374BA3192B2B28E140CF5ACFF883A9DD58E88D15B5B3B89885797F85013D3EBC54539FF8BA96
A339C59B71E9C1A97E243EDF0DA81F78DCC78A5D07D4DAD1EF1A7B4B746B3B98488972644C0E9EBDEB0FC653DBD9EB3A05DDCB0444B86DCF8ED818FD
6BA14B910DB592302CF3054F4C71C9AE5FC791ADD6B1E1E819415FB496607B805723F4A625B9BF2EBF67024524E93C514A40595D30A7F5AD50430054
820F208EF5CBFC4651FF00088CD8E36C88463B735D158FFC78DBE7AF94BFC85302C56341E22B39B585D2CC5731DCB2960248F6F15B35C7DCFCDF13ED
0600DB68C49EFD293046BDF788AD34EB9860BB8AE63799F62663E09F63DEAFDFDDAD8594B74F1BC8B129621064E2B9AF1C806EF44C60BFDAC6077ED5
D1EACDB349BD61DA073FF8E9A00A1A5F896CF5483CCB38AE2461F7D157253EBCD5FB0BFB7D42376B76398D8ABA30C329F43585F0F6048BC2903840AF
33BBB1C72DF31FE94DF0C6D3E28F117390265C7E3BB3FCA8B858BB73AE35BEB567A7CDA7B86B86C2485D4E07AF149E2FBDBAB3D167FB2444B3C64193
B2FF00F5EA8EBFFF0023B68633C73C7E75A1E33556F0B5EEE38C28239EF9A005F0B5C4CFA1DA24B6B3C6CB18CB3E307DFAD6E565E9322DBE87601B92
CA8839EE6B5298882F2EA1B2B592E6E5C24518CB1359D6FACB4D02DDCB6924368D8226661D0F4247A5617C4F69174284216DA65F9803807A6335A97D
61AA6A5A135830B4B759102E5598E14631DBDA95C7635EEE40DA7CD246720C4CCA7F0AC8F0412DE1E462C4E657EBDB9AD0860FB1E846DDA412795032
96F5C039AC6F09063E0B70BF7889718F5E68035D3564B8BA96DEC223726038918300AA7D33DCD5AB4BC86ED19A16C953B594F553EF5C678013538EC6
E7CB8E031BDC365E427271F4AE8345D1AE34ED46F6F27BAF34DD9DC5003853400C9FC4A906ACDA749673F9BB729820963D863DEAC9D61609A34BFB77
B5591494776041C73838E87DAB2248525F88F1BB853E5DB12A31DF039A97C6A91C91E9892A12A6F501F4C6791F88A00D6B8D512DEDCDCBC137D9C121
9F03A76207A1F5ABD1CD1CB0ACC8E0C6CBB8376C550D74247A05DA818558481ED59B6D238F0124885831B6EDD79A045C835A37924A74FB46B88626DA
650E141C75C0EF57A3BE8EE2D1A7B51E763F801C1CFA5737E0B5BFFF00847A316E6DD62C9C17C96CF73C77AD3F0F69434896F227BB13C9338948C636
8E681900F14C63527D3DEC675B918D91F52C4FF21EF5A29A984BE4B3B985A19241943BB706FC6B27498A36F1A6A92B0CC8A800CF619ED526B6DFF153
E8EA7032C4D1A8177C4971736DA34EF68809D843313F7474A4F0ECD33E9167BA13B0C63E72DDB19E94FF001336DF0FDE93FDCFEA29BA44C2DBC3BA7B
6C66DD1C6800EBCD1D43A1AF45158FE28B992DB4399A10C5DF09F2F519A6223BDD78C7E6B595B1B98A0E259036003E83D6B5ADE513C11CCA080EA180
3DB3585A36892AE9B6D1DEC8A5140628A3EF1F735D080140000007000A431D4514531051451400514514005145140051451400514514005145140051
451400514514005145140051451400514514005145140051451400514514005145140051451400514514005145140051451400514514005145140185
E21F0F43AB849A36115DC5F724C75F635CBFD8EE74F9BECF78B2213C860721BE86BD16992471CABB644575F46191512826690A8E27151B4831C02B8E
A4D35F29C2AB1DDDABAF874FB3854AC76D1004E795CF3F8D4C90C519CC712293DD540A8548D3DBF91CA5B687757A01B8CC709393BBA91EC2BAAB6B78
ED6058615DA8A38A9A8AD23151329CDCB70A28A2A880A28A2800A28A2800A28A2800A28A2800A82E12678C0B79444F9FBC53771F4C8A9E8A00E52D3C
2F79A7EB73EB16DA82CB7371912ACB1E148C83D8FB53B5FD0756D76D7ECD35FC10459C948D4E1BEB5D4D14582E62DCF87EDEF7C3B1691764111A2857
55FBAC3B8A7595A6AF05A2DACB756ECAA3689803BC0FA74CD6C51401465D3A09348974D04AC32C4D1123AE181C9FAF39AC5D2F43BFB3D157449DA16B
5427F7A339652D9C63D79FD2BA8A2801A88B1A2A20C2A8C01E82A2BA6B84849B48D2497B2BB6054F4500739A6596B906AD7179786DE559C05DA2423C
B033D38FA7E55178AB48D4BC41A5369E23B6887981D64F3091C7B62BA8A295877323495D521582DAEE08161862086457C9620606056BD14531189E2D
B29752F0EDCD9408CCF36D181D40DC0E6ABE92BA8E8BA7269D25BB5D9857114D1FDD61D81F4C5747452B0EE61787345FECC8EE669D57ED1792179470
71D78CF7AABA469571E1DBEBD16D119B4EB8712AAA91BE363C118EF5D3D14582E65D9C53CFA8B5EDCC4620ABB21463CE3B93F97EB5A94514C471F3E9
DA95AF8E7FB64C0D7368F118D427253818E2AADDD86B1078C6DBC446D5E684C650C118DCF18C600FD7F9D775452B0EE52D22DE4B5D3208663990025B
EA493FD6A4D42CE2D42C66B39C663990AB559A2988E4FC296B7BE1EB1934BBCB69668E3919A19A11B8329FE55B11D89BCB6BC17F1ED17A0A3479FBA9
8C01F5EA6B528A077399F0E5B5C7876C5B4A9EDE69A18A46682689376F5249F9B1D0D49E1ED2E68752D4757BC4D93DEBFCA9DD5074CFE42BA2A29582
E711A8DBEA977E2EB0D623D325582D10A3A92373649E953F8F6CB50D634E82C6C2D2462B2ACCCE71818078FAF35D8514582E60EB46F2FF00C2F71045
632FDA2EA168B61C7CA48C64F3557C3771A869FA15A58DD69570258536B15C60F27A575145311C5787EDB51D33C43AA5D5C697702DEF9F7A152AC57D
8F34CBD8F5393C736BACC5A65C1B4B788C472BF31C83C81F8D7714500711E38B7D4B5AB7B282C34DB82219C4CE5D40E83181CFB9AE92E11F58D1EF2D
9E292DCCC8F10F3060F2383FAD69D1401CAF8560BFD2FC3E34B9AD645BB88BAA363F767D0EEAA1E028350D20DDD96A56170AF3DC1916554CA74C124D
773452B0EE727AB69B750F8E74DD6E185E5B6585A09446BB99786C1C7A722AEC96B7B7167ABCF0030DC5DA6C801E18055C0CFA1C935BF4530B9C5F87
20D4A0F0AFF647D825B6BBDAEA66917E4C9CFCC4F7F4A6F842C350B3B2B6D36F2CA48CDA5D3CA6661F2B039E87F1AEDA8A560B9C9EB9A65E278C74CD
7A081A782DE2314A91F2E33BB903BFDEAD3B64BA8AE6FB5136CCCD39411C3B8021541E4FA139E95B345160B9C6F8574ED534AD5B537B9B1C45A85C79
A08901F2C649E7D7AD7557D6CB79653DB3805668D90E7DC62AC514C4725E02D0EEB49D3E5FB7EEF33794895BF8501CFEA79AD6F1368D1EBBA24D6121
DACDF346DFDD61D3FC3F1AD7A28030348B8D42D34982CEEAC2537502088145FDDB60601CE6B4B4EB67B6B76F37996576964E72031EC2AED1400C9635
96278DC655D4A91EC6B8EF0B594DE15BABFB0B98A56B59E612412A2E4723073FA57694500674025B8D49A768DD20488C681C60B12412707E82B9DF0D
D89F0B6A7A8DADC2BFD92E641241305C8EFF0029C7D6BB3A280392B4D28EA5E336D7E4474B782211C01D705DB182DF4EB557C6569AA6A3AAE98D65A7
4CE9A7DC094BE461F9078FCABB7A29587739DF15C379AA78566B4B6B394CF74AA36123F77C83CF3ED5368A6F20F0CC16F71632ADC5B5BAC5B32BF390
A064735B945311C97816CEFB48D2FEC17F692C723CCF26F182A0103A9FC2AAF88A1B893E22E86D6814BA40EC437420139FE75DA39654628BB980240C
E327D2B95BAB2D6E7F15D9EAEB668B15BC262311954939EBCD21DCB3E24B097C43043A6AC52C5089D5E791C63E51D87AF5A93C5904F71E1FB8D32C6C
9E533C5B176602A6318ADBB7691E1569A2F29CF55DDBB1F8D4B4C4739A559CD79E184D2354B4920C5B085C9C60F6FE8299A0477FA0E9E34CBBB696E1
2124412C237065EC0FA574D452B0EE67DBA4D04735D5C2334F272D1C783803A01EA715936ED757D777526AFA54CD6E4ED8222A1805C60920F735D351
405CE3BC2F65AC682B716925ABCF652485AD879837463D1BD3FF00AD5D8D14531105CCFF00678C3F952C9938DB1AEE35C469D63A9DBF8EAEF5C9B4E9
96D6746450A32FDB9C7E15DF5140EE73FABDC5EEA3A5CF67A7D94F1C9708D1979D760504609FCB34ED17404D0FC3B269D6AECF232B3339E32E571C7E
95BD4501738CF095A6A9A76893E9AF6AF1DDB4AEC267184C1C0DD9EA4D41E14D2353B056D36EED8958EF04E6E49F95C0F4F5E6BBAA29582E71DE20B4
D46F7C4DA5DEDB584AF0D8484B1240DDF4E69DE368EFF59D08E9F616173BA59159CBA8002839F5EB9C575F4530B9CDEBD6F79ABF8427B6B686586E59
14796DF2938C1207B555B8B3D46FFC1B1E916F6E6DA6FB3AC529946D1C0C607AE715D751405CE7FC290DEDB6916B65776C61482128DBBF89B3FCB155
F49D26E3C39A9DE9B689E7B0BB6F302C78DD1BFD3D2BA8A28119D124D35DFDB2788A88D0AC31923773D49E700F6AE756DB5393C749ABB584C2D843E4
A8380C3DCF6EF5D9D140EE721E2CD3AFEFB5CD267B4B57923B293CC76C8E7A703F2AEAE26668C33C66362395241C7E55251408CAF12437971A15D43A
763ED12465473838EF8AC7D26D2FA5F070D1D6D1ECA710F92F24830BC8E587739FEB5D6D140EE727E12D3EFF004FB0B3D32E6D8C6B632C8CD293C499
DD8DBFF7D9FCA9B75A5EA9278DE0D6E3B353043098B619543375E7F5AEBA8A560B9C978AEFA6D3F56B19A3B74BE2EA516D3396CF3F301FA67DAAB9F1
1EA764A266F0A4B1C207CE6304328FFBE6B66EB499535D3AC5A859A53179462725703D54FAD3EEE4D5AE21F260B1483CC1B5A59260DB01EA703AD002
CC967E27F0E32E3305DC7C67AA9EDF91155B491A8E95A6C3A6C968D70F026C8A54236151C2E7D38AD5D32C63D3AC22B48892B1AE371EAC7B9AB94C0C
ED36D66B582696E1849753B191F1C0CE3851F4AC1F0C6997DA19BCBAB885A417F2995E38D7E68CE4E38FC6BAFA2811CDC5A4CDA8789935ABE431C56E
9B2DA16FBC0F763E954FC51A76AFACDFD8F916412DECAE44A59A55CCA011D0678E9DEBB0A28020B6926962DD716E607CFDC2E1B8F5C8AE4741B1D5B4
9D7758B8B9B56B96BC7531CA0F03AF7F4E47E55DAD1401C4C7A2EAF67ADEA376556E64D4A2D8655E046781F962BAFB484DBD9C1092098E3542477C0C
54F450015CC788BC3EFA9EBDA55E47C2C127EF48F41C8AE9E8A00E73C4DA136B17DA5C88A0A4131339CE0ECF4FE75D09C853B00C81C0E829D4500728
746D5CF8A86B456D3708FCA09BCE31CF3D3DEAFF008A7489759D216088A8951D640A4F048ED5B9452B0EE64DB43777135B4B7508816DD30137062CC4
633C718AC8F10693AB6A7ACD95D410C690D949B94338CBF3D7DABADA28B05CE63C59A76ABAE6982C2DA28E15660CEED203D3B74ADBD34DD0B609790A
44C80280AD9C803AD5CA29882B06F74A9C788E1D62D95662B198DA32DB703D41ADEA280392D7B42D4B51D56C3518FCB2F6CD91116C051F5EE6B6F548
AFAEB48960B658D679A22ADB9B852460E2B4A8A00C5F0CE9D75A4E890D8CE10BC208521F3BB9CFA7155346D2752D3AFF0050BB2B6CCD7D287237B7CB
C13E9EA6BA5A29587739CB8D0EFEE75EB7D5A4B8B70F6EA5523018AE2B5B53B21A8E992DA4857322E338E335768A62326D2CEEBCBB686E4C6B15B1C8
D8725C8E87D856B51450050D634C8757D3DECEE090AC410476229910D463852DC243950079DBB8C63D3D6B4A8A00A17B6724BA649696CCA9BD0A163C
6011C9E3BD56F0F6973691A7AD8CAF1490AE482339393DF35B1450064D9699269B34A6C4A3432B6F68E46C60FB100D5F488890CB260C8C31903EE8F4
FF00EBD4F4500603E8F7A7C4BFDACB2DB81B3CB0A41C85FF001A4F1168D79ABCD6AD14D0C696B28955589F988C75E2BA0A28032B55B1BCD47486B232
448F2A6D95F9FD29BA4E9B7567671D8DCC90CB6B1C2235500E49F535AF450064D969B3E971182C1A230962C04B90549FA75ABB040D1F992B10F349CB
1E838E807B559A28039FB5D2751B6D62EB50592D89B9C064258803F2A2F346BCBBD5EDB5179600D6DFEAD0138FE55D05140EE53BDB25D434F7B5B938
F3170C53B1F6A8ED6C1E2B6B6B77906DB6236ED1F780E99F4AD0A2810541796B15E40619B3B4907838353D14010471BAAAA3B2385031F2E0F1D0F5A8
755B16D46CCDBADC3C1920965EA47A55DA28022B788416F1C418B0450B93D4E2A5A28A0028A28A0028A28A0028A28A0028A28A0028A28A0028A28A00
28A28A0028A28A0028A28A0028A28A0028A28A0028A28A0028A28A0028A28A0028A28A0028A28A0028A28A0028A28A0028A28A0028A28A0028A28A00
28A28A0028A28A0028A28A0028A28A00AB7B7B6D616C6E2EE51144BD588240FCAA3B2D52C75081E7B3B859A3419665071FE78AB173047756D2DBCAA1
92452AC08EC6B8FF0087723DAFF6968531F9EC273B72304A93D7F41F9D219BEDE25D1526F25F508D643D158107F95596D574F4B74B86BB896273B431
6EFE9ED5C6FC4DB2611DA6B16A8864B493121DB92070466BB0B2BDB7BFD1A2BDD9982587CC2A40E98E46280B125A6A9617AFB2D6EE199BD11B3572B2
F42B28AD74E46485237972E7006704E40FCB1534FA9DB4533C2BE64D220CBAC4BB8A8F7A622F5155ACEF2DEF61F36DA40EA0E0FA83E86ACD00145159
D77ADE9F69218E6B81BD4E18282DB7EB8A00D1ACFBFD674FD3AE60B7BCB858A49FEE03DE91B5BD3164810DEC3BAE31E50073BB3589AF9F0B6A7A8C31
6A973B6EEDD8AA00181CE7A74E79A5719D65154AE6FEC6C0C715D5D4506E1F2891B1C0F73FD688F56D324CF97A8DA363AE2653FD6988BB45558EFECE
470915DC0EC7A2AC809FE74B757D6B6607DA6748F3D013C9FC2802CD151413C57110920916443DD4E6A5A0028A6B30552CC40006493DAB923AFDA5DF
8E2DEDA0BB5305BDBB973BF0858D0075F4C49124198DD586719539AAB756F69AAD998D9C4911390D1BF43F51F5AE47E1A9FB3DA6AE2E27C88AEB6177
7E0003DFEB401DDD151433453A6E8254913A6E46047E94E9248E24DF23AA28EEC702801F454714D14C9BE1911D7FBC8C08A92800A28AE3FC6374D26A
FA4E8A6468E1BC76F30A1C13E83F3A011D8515CDF8334CD574AB3B9B7D5640E3CDCC387DDF2D6F99A25916369115DBA296009FC28025A2A1B9812EAD
DE09376C7183B4E08FA1AE23C1C935AF8DB5BB17B99E68EDD004F35CB60120D0077B4579D78B2DE4D37C5BA2CB6D7771B2EAE177C464240F9873FAD7
A2D00145145001451450014515CE1D02E3FE12EFED94D45FCA030D6F927B62803A3A28A2800A28A2800A29B91BB6E467AE29D400514532491235DD23
2A0E99638A007D151B012C44062032E3729C1E7B835C0C02F2C7E28DAE9ED7F733DB340CEA92C99C0D8DC7E9401E85451450014520209201071C1F6A
5A0028A28A0028A62C88CC555D4B0EA01E453E800A28A2800A2B0749D4E6BFF11EA90A48AD676A1235518E1FBD6F500145145001451589E28D2AFB57
D305B69F7DF64937862D92323EA39A00DBA2AB5942F6D630412CA6578A354690F5620609A9D59594329041EE0D003A8A2A3134464318910B8EAA1867
F2A00928A2B989B59F3FC7367A559CE0A430C925CAA9E092381F5E86803A7A29922978D94315241008ED589E16D224D22DAE525BEFB599652DB81CED
F6A00DEA2A39658E14DF2C8A8BD32C7029E082010720F42280168A28A0028A28A0028A2A369110AAB3AA9638504E33F4A00928A2B8FD66E66D4BC50F
A22CAD1AA5A1963556DBBDFDCFB5007614564F872DF51B5D1A18755904974A4E486CE067819FA56B500145359954658800772690C8814317500F439E
0D003E8A6BB2A296760AA3A9270294104020820F42280168A6BB2A296760AA3A927028465750C8C194F420E45003A8A28A002986441208CBA8723217
3C9FC29F5CD4BA24573E2E87581A98CAA8DB6EAC0EEC0C7AF4EF401D2D151BCB1C79F32445C0C9DCC0714AD22232AB3AA96E14138CFD2801F4535995
14B330503A927142B2B28652181E841CD003A8A8E49A288032C8899E9B980CD39595D43230653D083906801D45725F102C94F87AEAFA3B89A29E2504
6D908079E98AD7F0C127C37A7B124930A9249C93401AD45155AFE25B8B29A0694C3E6A15DE0E08C8A009D595D43290C0F420E69D591E1BD223D174C1
6B1DCBDC02E5CBB1EE7D3F2AD2B988CF03C6B2BC45870E87045004B45711236A107C43B0D3E4D4A796DDA16942138E80F5F5E95DBD0014546B2C4EE5
1644675EAA181229E480324E00A005A2A28AE2098910CD1C8475D8C0E28B8823B981E194128E304038A0070910B940EA587519E69F5C3F85EDD2DBC7
3ACC0858A451AEC0CC4E39AECE49E184A896548F77DDDCC06681B25A299BD77ECDC376338CF38F5A7D020A896785E5689258DA44FBCA18123EA29D23
A4685E465551D4B1C0AC5D2B4BD2EDB59BCBAB3B8125C49CBA0707666803768A864B9B789B6CB3C48D9C619C034A6E2012F94668C49FDC2C33F95004
B45145001451450021200249000EA4D3629639937C4EAEBEAA722B99F1CEAD2D8E8B3C1679370E9C9519D8BEF563C0C77784AC5FBB264FD695F51D8E
868A28A620A28A2800A28A2800A28A2800A28A2800A28A2800A28A2800A28A2800A28A2800A28A2800A28A2800A28A2800A28A2800A28A2800A28A28
00A28A2800A28A2800A28A2800A28A2800A28A2800A28A2800A28A2800A28A2800A28A2800A28A2800A28A2800A28A2800A28A2800A28A2800A28A28
00A28A2800A28A2800A28A2800A28A2800A28A2800A28A2800AE135A2BA27C47D3B51E56DF518CC32EDEEC38C9FCD6BBBACCD5F4B8F536B2670B9B5B
959816F41D47E3C526343B56D3E3BED22F6CC83FE908D9E73F3638FE42B90F026A50C9A17F625EB37DA22B8F2B6639DB9CFE5C115E8158DA7E856D63
ADDF6A4889BEE482A7BAF1F35017346F5DE2B29E48F8748D994E33C81C570DF0DEF2F2E6D351BC787CF966B8F98EF0A471EF5E804020820107A835CD
E85A0CFA0EAB79F650B2585DBF998C80D11F4F714021FA0D96A76BAD6A92DD22C76773279912060486E3D2BA2A8D464AB9CAF1F773EB8EBEF525310D
6C95214E1B1C1AE2FC068F1CFAC45A901F6E139126EFE25FF0AEBAF5657B2996DE411CA50EC73D8E2B81B5D7844C60BEF0F5CDE5FC24ACB3C6D91211
DF3DE931A27F09E8F0DD5EEA89710EFB38AEC3DBB74E41E83DAA6F1EC6A7C41E1A628B86BBDA4E393CAD6A7873C4D65A9CED62B69258DC01B840EB8D
C3B91C556F13E99A8EA5A8D85DC167B92C262F8322E5C67A8FCBA51D07D4E8B50B0B6D46D5EDEEA249159481B9738CF715C678096DADD752D12FADA2
3358CACDB9E304B267AFF9F5AEE2D6669EDD25689E1661CA38C15AE4756D1EE97C7F677969B920BD88A5CB2719DA39CFE18A188D0F0FE8F66A67D49E
D510CB319221DD1474AC2F066A69AB6B3ADDFDF2B4815D562CC65F62FCDE80E3A0AEFC2A85DA000A06318E315C6E89612F8575DBF89A191B4DBC6124
7328DDB0FA37EB45807F86E6963F166A71476D343617204909642AA58019C7EB5D9565A3B5FDEC52A2C91C36E490C463CD2463F2AD4A1031ACAACA55
802A460823822BCFF4DD32CD3E27EA969F668CC06D8384C7009099FE75E855C798A5D37E235C6A370856D2F2D04692E32030DBC1FF00BE68608E9ED4
5AC6AF0DA08D44476B2200369C77AE0FC01A7ADCDEEB7F69512C0B74404278DD93C91F4AEB74F296B15F5F5C0D8934E5C37F79380A71EF5CB781B50B
7B2BDD612F0BC1F68B932C7BD0E19727BD005CF0EEDD3BC77AC6956C7166D1ACC9183C2361738F41F31FCAA1D435682FFE2147A45D30FB15AC64B2B1
C2B3EDCF3F4E2A2D1EF6D57E23EA770F21F2E6882C5200769E99C9FC2A7BDB19343F1EAEBA62DD61789E5CCE067C9620727DB81F99A0075E5F43A6F8
C2CE4D306EB59D7CBB95887C8A73C1E38AEDEB0F50BA86F923B6D3E4492795D49755CED507A93F856E50818579978B344B1B4F16E851A2CBE5DD4844
80C873F787439E3AD7A6D79D78C353B39FC53E1D9A094491DB4E4CAE3A2FCCBD7F2A1823A2D6A583C29E1BBBB9B30FBBF8048E5C6E3C0EA6A3F0EE9B
6DA8787ED6EF50569EE6E23F31E57625813E87B62A7F135947E23F0E5D5A594F1C926032EC607E61C807D2A0F0A6A36D6DE1AB7B7BA9045359A7972A
3704107B7AD1D43A15BC09AC497C97FA6DCCA65974E98A2C8C7978C9207F2FE555FC387FE2E5788863F817FF0065A7780B486D3D2FB57BBF90DD3315
0472101CE4D52D0EF6D61F891ACCCF21586E6302391810A4FCB91FCE802D78E973E25F0E608DDF68E9F88ABBE3BD766D2AD2DAD2C5F65DDEC8115FFB
8B9E4FEB591E3ABCB57F10E86E932BA432E6564E420C8EA455CF1F69EDA9E9D65AB69804EF672071B0E772E7FA114017B58D12DF4FF0EDDCF64D2457
56D0B4CB3F987716519C9FAE29FE1DBD1E29F0AC12DCBB2CD92AED1B60ABAF7FC88FCEA6F106A56F2784AF1E3911DAE2D9A3445392599718C7E351F8
17479347F0D410CC3134A7CD71E84F6FCA80323C2F0CF1EB9ABE87A9DD5CCCD1E1E366908DC87A1AADE16437DAF6ADA66A9777339B760D06E948CA64
8EDF51573C636F7565E27D2B56B007CC9F36B2EDEA41E9FA13F9557F16DACFA36BFA56A1A5B057B85FB1B83FC59EE7F4FCA81A35FC2B60244B8BB379
772406E2458627949431F4191DEB2ACAD62B2F8AAF05B33AC46D7798F7123247FF00AABB3B74834CD39232CB1C30200493C0AE12D751B65F8A3757CD
2FFA2CB6E1125C7CB9DAA3F983408B1F13606B7B2B7BDB59EE22B992654212420371E99ABDE23D3E58FC29717EF7533DFC312CAB2EE2369041C01E95
57E25DCC3F65D3E11229952ED19941C951EA456BF8CAEA03E0FBD659908961FDDE1BEFF23A7AD032C697A9CD71E11835478CCB39B5F34A28E5980E9F
98AC1F09DB1D7B41B8D4F529E47B8B8760ADBCFEEC2F4C01E956F4AD40D8FC3A82E6D621712DBDAA96887273DF23F33F8567E9F75A2D8696561F142C
2B7A4CD2A614E19B1BB68032B4010F846DAE35F9E0D5AE6F246B8B1B864773D2541D062BD12B13C2EBA443A52C1A1CE935BA1F9995B24B1EE6B6E842
615C5786EED7C4FAE6A57975F3C167208ADE23CAE39F9BDFA576B5C0783AD9BC31AEEA5A5DF911C770C24B691B857033C67D706811B3E1E96FD7C41A
C59CB6AF158C4E0C05871F87B1EB587AB4525D7C57B586DE578BFD1BF7AE9C32AE09E0F6CF033EF5D87DADAE354861B3903C318669D9705791F28CFA
E79AE45AFAD57E2C2CAD3C7E535A796AFBBE52D8E99FC2819635B1FD81E2CD05B4E69238AFE568678F79656E57079CF3F37E95A9E2E3A995B18B4DB6
9268A49BFD2427F7063AFF009ED58FE3BB98A2F11787242E0FD9AEB7CA073B0129C9F4E9567C7BAC5CDA59E9A963334705EC989668FA84E3A7D73401
0F8C1E1D19F4EBCD2E610DC2DCAC6C88D90C0F504575BA8E9F6FA8DB7917218A67702AC5483EB915C278F0DA4163616DA7C5BE3B6B959269546E0A3F
DA6F5AEB6C2FECDAF6FEE56E62F258C615CB601C2F41401C5782A2D4350B0D7ED0CF24B1C395851DCE1A53BB1CFF00C057F3AD5961F1269FE08B5B75
8E49AF4C989954E5963E78CF6A67C33B8B78EDF5385E68D657BD62A858648C0C63D6B47E21EA37BA6E80AD60CC8D34A23791472AA41E87B501D4CEF1
68B6D23C391EA1A6916D7B1C91B10B212727A82335D6BEA096FA27F68DC11B52DC4AF8EFF2E715C678DE1B3B6F089B5D2A3F3CB32BC92A9DEC17AE58
FBD6A29FEDED16E34DB77DD17F67C6BB9791E66338FD314009E1A7FF0084834E9B5CD45A46323388E28DCA88914F4183C9A97C23793EB9E1DB88EF84
D18499A156DC55CA7040CFAF38A8BE1B3F93E1E3A7CC0C7756B3389236E0F2739FA735BF6770F35F5CEC2A6D23015481F79F92C73DE80389F03E916F
36A7AF4427B94582E7CB50B29048C9E4FAF4AE87C5FAD7FC235A021809699F10C45CE7A0EA4D65783254B3D73C502E5D63DB76580271C65BA7E950F8
CB4BBBD4FC196774AA6492DE469DD0F565627FC4501D4BDAB81A1F86975382F247BD85524F9A4C894923208F4EB5A1AEF93A9F8565D52DA6993168D3
C4C8E578DBBBA7E15358C1A35E69105C98A0788C40127D71C83EF55F51992DFC1B7D15CF91017B4944512363E52A42803F1A00C0D0AD35AD53E1F896
1B8792F1E5CDB97948014360FE3F7AAD7C468648FC2915DB4922DD466347657201F5C81C1E6B4BE1DBA7FC217A747BD778572573CE3CC6AA7F136789
BC34D68AE1A66993F76A72D8EBD28E81D4DBD496E25F09C91D9EE37125A844DBD72401FD6B8CD72CC594361E18D29E59356B9DA659BCD6CC60727BF4
EBF8576126B10D8785D2FD4F9DE5C2A02AF24BE0715CEF864DBE9F6777E27D666592FAEB2FB179645FEE81EB401BDA95BDEE91E1092DB4612CF77146
150F2CCC49E48F7EA6B32F2DADED7C1867BB856CAFD6DCC9966FDE07F5E393FD3348FAC5FDBF806E7595120BCB877745907280BED1F80519150092DD
3C0E65B50F7DA95EDB1DDFC721623E6CFA01CFE54058E87C3570FAB7852CA6BD3E6BDC418949E377506B91D2F42D3E0F89F796714256DE2B4F311448
C0863B7BE73DCD6BF836459749D16246064B6824F3973CA64E066A38E36B5F8A7737330D90CF6402BB70A48DBDFF0003401D1FF66D8BD9496889BA35
6248121255B1EB9CE79AE73E1792FA0DD39766CDD328DCD938007F8D6CD8CF6F696B7FA8CEC5219A766F9BB81C715CFF00C2CB88E2D0E6B398ECB96B
A76119041236AFF81A009B469878A3C47A94D780B5A5837930C393B73CE49F5E9534D77FD89E39B2D3E12DF65D4623BA327215C13823F2A8FC2B6526
83E21D5ECEE3E586E9C4D048DD1C64F19F5E696FAD1B59F881637100CDBE9B1EE924CFCA5B27007BD0055F8990CB143613D8BC8B7535D2C5B439DAFC
7008FAE2B43C49A6B58F86AEAFD6E656D4205128B92C412411C63A018E2AB7C40B88166D155E41BA3BF8DD80E4AAE4726B43C6D7B6DFF0886A012746
6922DA81581272476A00D7D12E9AFB45B2BA7FBD2C2ACDF5C7357EB17C272C3278674F58A456296E81C03CA9C720FA739ADAA6233B5DD40697A35DDF
EDDC608CB01EFD07EA6B13C1D6716A1A0C3A9EA2866BBBC6792477278F98E001D86056D6BDA7FF006A68B77619C19A32A0FBF51FA8ACEF08CDF63F0E
4169A832C571680C6EADC1C0271F5E290FA153C33A939D7756D0277693EC8C5A1624E761C6467F1159175A35B1F895059F99308DED0B1C484377EFF8
56BF85F4D71AE6AFE20B80635BB72B086E3F760F53F90AA1717707FC2CFB5BDDFF00E8C2D0A1971F2EEC1E33F8D033B4B2B48EC6D52DE12E513A6F62
C7F3AB34C8A449A3592270E8DC860720D3E99273DE33B2B6B8F0DDFCF32FEF21B7631B824608E4743EB59DE0AD296F3C25692EA25E66910840588F2D
3240C7F3AD4F1A3A8F0ADFC64FCF244422F763ED4CF07DD40BE11B22F2AAF9108593271B4FBD2EA3E862F83927D4DAFAD2FA5696D2C2E1A28A3248CF
27AFAE2AEE9D33D8F8F6E74A809168F6C2511924853C74AABF0FE685AEB5A1E600D35E33C6A782473C8FCE9F05F5A1F88F752EF4D9F6654F373F2823
B6681B346FECEF6EFC616BE7DB79BA5A407A9F9564EB923BF402B36E6F96C3C7DA759E98CBF66BC8C89E34FBB91BB07D29FADCCD77E36B3D3AFE5316
9B1C5E6105F6ACB276079AAFAB5DB27C40D1A4920F234FB6474495970AC4A91C1F4E571408EB352D2EDF5244139915933B191CA919AE174EDBA26B17
BA46AE66B99A5E6CA46724BE7A0EB5D5E957B0D968EF777D2F96B34D2C8BB8E4B29638FD2B135AB4835FD39B57FB4086F2DB2D6AB9C15039C1EE7340
1D0683A145A65BB7985A59A507CC2EC5860F6E6B9BFB24107C56B58EDD7620B667280F19DA7FC6B73C35E24B7D53498E5BA9121B945C4A8C71C8EFCD
73F26AB643E24A5F190FD985B797E6E0E334022EFC4CB7813C38F7821C5C2CA8A255E08E7BFA8ABB6DA446DE1959EF59DEF5AD7CC6999B956DB9C0F4
C5677C4DBE85FC31F6585C492CD2236D5393B7AE7F956F5DEA369FF08E492A4AACAF032280792DB7A501D0CAF060975BD063B9D57F7DB498D14938C0
EF4BE13BA74F106B7A506260B694188139DA3D29BF0F6EEDA1F09C31CB32C6F1B3175738239A8BC2D3403C57AE4CCDB7ED3301096180E3DA81972DEC
AEA5F146A175AC448D628A16D9A520A8FA0354B49BCF2BE225DE9D65386B236DE618D4E42B0DBD3F3A70B94B9F1A5EA6B6DE5DB59A036B1C9C2367AB
63B9ACF86F9E1F8817FA8DDC12DBDBCB60C906E186940298DB8EE48A00DFF887FF002265F7D17FF4215ABA0284D034F51FF3ED1FFE822B9DF18C861F
02B591DC6EE4863022FBCD9C8CFF00235B7E1ABB86E344B48E2903BC3046B20FEE9DBD3F4A05D0D7ACED6ACADEF34E9FED09B8C713953DD78FFEB568
D65F886F62B3D26E4B9F9DE2708BFDE38C7F5A62333C00E17C1D6D23B1C069092C727018F5A87C337EDE26BABFBF94C8B6904DE4DB22B951C725B83D
7914BE05F2DFC271593B6C9B0E248CFDE5C93DBF1A83E1DD9CBA5DB6A5A75C0DB2C7745F07B8C0191EA3E5A4577229D1E3F8AB68F2BB3EFB67D800C6
060F157FC7BADCDA3E9112DB1DB2DD49E507FEE8EE6B22E2F925F8AD6AC582C10C0C9BCF009DADDFEA6B53C75A3C9AB69704F6C3CC92D6412041C865
E33FCA802AEAEDA7697E1833E9B2A35E20564646DCECD9EF56F5A9755D5FC2968FA6C3FBDB90AD2AEEC1031922B56D2F34E92CE3911610C540F2C28D
DBBD00EF591E37BCBFB3D1208ACD8ACD3C8164302F213BEDFCE8023F187D8AC3C3323DBF970DDDBEC286118656C8F4E95D2E9333DC693653CDFEB258
11DBEA54135C978AFCB1E1196C347B7697013CD75E76804753DCD6EE93748228D8BEDB782D235724E1437D2842E865E8213FE13CD70A8E760CFE74FD
6756B4D2BC461F57858DBC91810CA5772A71CF1F5AABE1ABB0FE30D62E25578E3B92AB01653F3E0D685AB7DABC43A9DA5F5BF9B6CC40469065738E83
3D290C77869209751BFBEB2B933DB49B56324F4EE4574B5CC691A6269BE28B91A6AECB278B32A03F2ABE7B574F4D098D745752AEA194F504641AE574
08E14F19EB622454DAA8303DF9AE9EE675B681E670C55064851935C8786DE75F14EA97771673225E15F29B6F18F7A18219F1034F865FECE3122C73CF
7214B8EB5ABAE585B58F866E9A34CCB0C7B84A7EF96F5CD66F8CEE7CDBED2E1B78DE430DCEF91829C20F7ABDE2DD46297C39750DAEE9669936AAAA9C
F5FA503353C3F3C975A0D8CF29CBBC2A4D69563786268FFB16D2D0BE678615122E0F15B34C961505DDCC5696EF3CEDB51073EFED45D5CC76B1F9926E
C67002AE493582DA81D535C82DBECB38B38B2C5990E1DBB647A0A00A9AC5D41278635195C33DC5C216FB8781D867D055DF011FF8A46C79CFCA6A6F13
DD08F47B9B58219659E48F6A224448FE58AA3E0DBA367A25969F736D70973C82A63C63DFE94BA95D0EA88C8201C7B8AC4D0F48BDD3EF2EA7BBD41EE1
673F2A12485ADCAC7D2355BABFBDBA82E2C1EDD613F2B93D6811B1451453105145140051451400514514005145140051451400514514005145140051
451400514514005145140051451400514514005145140051451400514514005145140051451400514514005145140051451400514514005145140051
451400514514005145140051451400514514005145140051451400514514005145140051451400514514005145140051451400514514005145140051
4514014F55B3FED0D3A7B4F30C7E6AE370EA2A86951DD695670594964D204F97CD8483BBDCE71FAD6DD1401970DAB5CEA31DF5CDA8864814AC7F3062
73C1C915A945140051451400514514005145140051451400514514011B43132046890A8E80A8C539955D4ABA8653D4119069D4500451C10C47314488
718F9540A968A2800A618D0820A2F3D78EB4FA280228A18A1CF931247B8E4ED50334D6B5B6772EF6F1331EA4A024D4F45001510821520889011D30A2
A5A280227B782424C90C6E4F5DCA0E695238D23F2D1155318DA0607E552514015574EB2495A45B4843B724EC156A8A280229218A528D246AC51B7212
3383EB493DB41706333C49218DB726E19DA7D6A6A2801ACA194AB00411820F7A8DADADD942B411151D01418153514015A4B0B39A532CB6B0BC8460B3
2039A6B69D62C141B3B7210E54796300D5BA280307C456B791685723408638AE9F19F2D002477E9DEB374ED17C2515880F6D6FE622832F9C4EFDD8E7
AF7EBD2BB0A8DE18A439922463D32CA0D03B9CBF8774BB583C45797FA54261D3E484478030ACE083900F6C5759480000000003A014B4082A19ADE0B8
50278639429C80EA1B07F1A9A8A008A1821814AC10C7103C908A173F9545FD9D63B81FB1DBE4743E52FF00855AA2802B7D86CF24FD960C9E0FEEC734
E92D6DE584432C31BC4B8C23282A31D38A9E8A0080DB5B985A13045E537DE4D8369FA8A6AD95AAC0B00B687CB53909B0601F5C559A2802A269F651C8
248ED2DD1C1C8658941CFD7156248A3990A4A8AEA7AAB0C8A7D14015E1B2B582030456F1AC47AA05183F5A75B5ADBDA47E5DB429127A22E2A6A28029
5C6976573279B2DBA193FBC3E527EA475AB51469146238915117A2A8C01F853E8A00A52E996335D0B992D6369C0C6FC727EBEBF8D5C2010411907A83
4B4500505D1F4E493CC5B38836777038CFD3A54F73656B778FB55B45363A6F40715628A00AF059DADBB6E82DA188E36E5102F1E9C539ADA0690C8D04
65CFF11419FCEA6A2802258215428B146109C950A304D1F6783FE7847FF7C0A968A00694529B0A82B8C6D238C5436F676D6A0FD9ADE38B3D4AA804D5
8A280238E18A224C71A216EBB540CD24B043363CE8924C74DCA0E2A5A28018D1A347B1914A74DA471F950B14687288AA7A640C53E8A00867B682E001
3C11CA0740E81B1F9D2C30C50A6D86248D739DA8A00FD2A5A2802BC9676B349E64B6D0BBFF0079A304FE7486C6D0800DAC181D079638AB34500436F6
B6F6DBBECF0C71EF396D8B8CD4D4514005466188C9E61890BFF78A8CFE75251400545F67836EDF263C7A6D152D140080050028000E001DA968A2801A
E8AE30EA187A119A62C3122B04891437DE0140CFD6A5A2802348228CE6389109EEAA0535ED6DE4C992089B27277203535140104D6D05C6DF3E08E5DA
72BBD0363E99A74D0C5326D9E249141CE1D4119FC6A5A28022686278C46D12320E8A54103F0A558A3442891A2A9EA00C0352514011A43146DB923452
4632140A76D519C28E793C53A8A008DE28A43978D18F4C9506830C4C815A342A0E429518152514011A410C6ACB1C48A1BA855033408620432C6808E8
428E2A4A2802096D6DE691649A08A475FBACE8091F434E78629191A48D1D90E54B28254FB7A54B45002606738E69A88918211554139C018A7D140051
45140055596C6D26984B2C08641FC58E4FD7D47D6AD5140152E34EB2B903CFB689F68C0254640F4CD588A348A358E350A8A3000ED4FA28022104224F
3445187FEF0519FCEA42010410083D41A5A2802348A38D0A2468AA7A8550050B1468851235553D401806A4A28010000000600ED4800049000CF27DE9
D450014514500145145001451450020001240033D7DE968A2800A28A2800A28A2800A28A2800A28A2800A28A2800A28A2800A28A2800A28A2800A28A
2800A28A2800A28A2800A28A2800A28A2800A28A2800A28A2800A28A2800A28A2800A28A2800A28A2800A28A2800A28A2800A28A2800A28A2800A28A
2800A28A2800A28A2800A28A2800A28A2800A28A2800A28A2800A28A2800A28A2800A28A2800A28A2800A28A2800A28A2800A28A2800A28A2800A28A
2800A28A2800A28A2800A28A2800A28AE6FC6EF77FD8F0C1A7CEF0DC4F7291AB23608073FA5007494570FE31B7B9D1BC1ED345A95EBDCA3A0F31A73D
49E7F0A4BCD1EFFF00E1188353D3355BEFB72DBACA732960F900918A5719DCD1593A46A69268163757B711C6F2C2A59A470B938E4F357E0BBB6B8FF8
F7B88A5E33F2386E3F0A6227A2A97F6AE9A4E06A1699F4F397FC69EDA858AFDEBCB71F5957FC6802D5155CDF5A0009BA8003D3F7829F0CD0CE9BA095
2450719460467F0A00968A6492C7126F95D517D58E05320B882E10B5BCD1CAA0E09460C07E54013515149710C4C1659A3463C80CC01A3CF8766FF363
DB9C6EDC319A00968A648E91A17919554752C7005305D5BB3EC13C45FF00BA1C668026A2A3324624119750E4642E7923E948D3C28E11E54562320160
0E28025A29AAEAE328C187A839A63DC428FB1A68D5FF00BA58034012D14D474910346CACA7A153914EA0028ACFD475096C8FC9A7DCDCAEDCFEE5727E
958FE14F134FE20D47518DAD7ECF0DA155556077E493F7BDF8A00EA28A28A0028A28A0028A2B0FC5BAC0D1740B9BA560B315DB10EE58FF0085006E51
593E19BC9350F0DE9F7733EF92584176C6327A1AD6A0028A28A0028A28A0028A2AA6A11C92D9C8B0DC9B6703225001DB8F5CF6A00B745723E00BED5B
55B0B8D4353B9F323924290AED0300753C7BD75D40051451400C01816CB02BD801C8FF001A7D1450014514500145145001451450014514500145616A
8BAF402E2E6C6E6D9A245DC90345C9C0E99CD66F8335CD53C456B35D4ED6F0AC3218CC6B1E771C67D7DE95C763AFA2B93B5BAF14DE6A5711452588B2
89F6ADC14396FA0CF35D4A065450EDB980009C6327D69887D14514005145140051451400514541797515959CB733B058E252CC49C5004F45733E0AD6
AEF5DB2BBBDB965317DA0A4202E30A3D7F315D350014514500145145001457137977AD6BFF006F9341BA6816CE710C4A081E6E3EF124F6AEC6012082
3131064DA3791DCE39A07625A28A2810514514005145417493C9095B6984527662BBA8027A2B924D5757FF0084DD347134525B2442595847838FE95B
3E217D4A3D1E63A3C7BEECE020E38F7E695C763528ACDD01B506D1E03AB205BBC61C0AD2A620A28A2800A28A2800A28A2800A28A2800A2A0BCF3FEC9
2FD94299F69D9BBA66B33C331EB1169C7FB7640F705B2A320955F7228036A8A28A0028A28A0028A28A0028ACAD4F5A874D57335B5D3ED190638F21B8
EC6A9F84FC467C476F75702DFC98E297CB404E49E33CD0163A1A28A2800A28AE691BC4DFF0951DE919D272476E17B1F5CD0074B45145001451450014
51450014514500145145001451450014514500145145001451589159EB1FF091BDCCB77FF12FC1DB103F971401B74514500145145001451450014572
BE25B8BD8353D3EDAC6F250F732FCE831F2AFF0091FA5752060004E7DCD002D145140051597A5DB6A30DC5CC97F75E6ABB7EED01E00FE95A44641009
19EE3B5003A8AE5F55BABA5F13E9DA6C1A84AA2605E45C2F41938381DF18FC6BA8A0028A28A0028A2A96AD2BC1A6CF3A4C6231A139001A00BB456578
71AE9F478A5BC94CB2CA4BEE6F43D2B56800A28A2800A28A2800A28A2800A28A2800A28A2800A28A2800A28A2800A28A2800A28A2800A28A2800A28A
2800A28A2800A28A2800A28A2800A28A2800A28A2800A28A2800A28A2800A28A2800A28A2800A28A2800A28A2800A28A2800A28A2800A28A2800A28A
2800A28A2800A28A2800A28A2800AAB75630DD4D6F2CBBB75BBEF4C1E33EF56A9ACC154B310001924F6A00E3FE29123C20F8FF009EF1E7F5ADCF0E2A
BF85F4D07E60F691E7F1519AE67E23EA36977E196B6B5944D29954ED50490066B63C29AB69EDE1DD3A06B98D2548150A39DA410307AD21F43427B1B6
B5D026B448C3C3142FB55F07B13583F0DE6821F09A3C92451933B825980EFC0AE9352BAB786D6549E78D19E360AACC016E3B0EF5CCFC36BB8078645B
348A93452B964638603AE7068021D360B197E25EA2B12432C5F64560A006556F979F4EE6AD7C45B3B61E13B899608D5E374219540239C555B1BAB68F
E296A12799188E4B3550F90149F949E7F0AD0F8872C67C237B087532928020393F781E9F4A03A9761934C4B2D36DE486176BA8D5146D19236F5FD2A5
D0F471A2BDE45030FB2CD379B127F7091C8FA702A0D3AD74FBAD2F4BBA97CB592D228D830600A1DBC83F8E69DA3EB5FDA9ADEA304641B7B60A10E3A9
E72681195E17BCFF0084835CD56FAE70F15ACBE4DBA1E428E79FAF1FAD1E2A9C787B56D3754B35318B99FC9BA45E16453DC8F5F7A4F0C5AB787FC43A
AE9F738482EE413DA313C30C9C8CFAF23F2A9FC5B62DAE5FE9BA6C03FD54A6795C72100E307EB9A063FE20C104BE0FBF9CC48D2222149368247CEBD0
D3A0B4D3DBC136FE6C11B092CD3040F999CA8E87AE734BE3BB8B74F09EA16FE62798515420232390471F8555D5EE08F8751DDE9EF896D6289D0A7F0B
29507FAD005FF08595FC3E198ADB5B1BE5C9F92419C2F606B12C2DAD21F89F796AB0A795F640EA8DC80D819C035D3E8BABC5AAE9D05E4663DAD16E98
EFFF0056FC6548FCF9F6AE76DE4B597E28DCB178DB75A2A2F39DC71C8FCA8034351D3EDF5EB7FB5E9C045756738F2264206FD98C8E3B751F80AD6D62
CAD2E6C27375024816263961C8C027AD63DDDDD9785347BD48268F7BCACF044A7254B638C7B1CD68EAB751C1A04E67B85F37ECC4E72016246338F426
8039EF0706D1BC43AB6873B911AE2784B9EABDF1EDFE154EEAC6293E2769AEE59E3BB85A72AC78E8D81FCAAD789ADDE4F14687AAE9CC0A5D836F2B8E
9B7AF3F813F953353BBB61F1474E7966458D2DC80FB863273DE819DC5BDBC36D108A04544049C0F7A9A991BA4881E375653D0A9C834FA6485709F0FD
C37887C4E463E6B953C7D5FF00C6BB679A246DAF2A2B63382C01C5703F0F2E6D935BD794C8AA5A6DC0B30191934868F43A2AADF5F5BD8584B7D70F88
224DECC39E3DAA0D1358B4D6EC05ED91731162BF3AE0823AD311A34514D665452CEC154772702800660AA5988000C927B5721E2E54BFF086A77B201B
767FA386FE1008C91EE79AD0D4358B7B9D62DF42B69D0C9282F3907A20FE1CFA9E9557C65A569E9E13BFD910468A2DE8031EA3A719A432FF00831B77
8434B3FF004EEA2B72B99F02DDDB9F08E991F9F1EFD85369619CE4F18AE9A9882B9AF116B0F69ABE99A779BF6786F1CF993E70401D141ED9F5AE96B9
1987FC251A8EABA36A5671A5ADA10219864481BB1140D16747D42EBFE128BED2C8B896D228C3A4D3724371900F71CFE95475DD5AFA24B2BAB7B83149
737A21B683B4B1E79247F9EA2AC780A6BCFECEBBB3BC7F356CE731452FF7D47FFAAA96AD08BEF88DA6C170563B2D36DFED0A090A376783F985FCA907
52FF008C352BBB7BDD234ED3E668E7BDB8C31538F90633FCEAD78E2FFF00B3FC297D2AB6D774F290FBB71FCB358CF99FE231BDBA702DACF4FF003A15
F50491F99C93F953FC5E2E6FACB43B5BC510A5D5E2FDA0630AA07201F4A00A76F16A3A07C3D86FC5CBC335B462448140D84338E18773F35747AAEB6D
6BE103AD451925ADD24500676EFC73F866B2BC7314DAB0B3F0FD80F9E69049315E91C63D7FC3D855CD4EEA4B6D4B4BF0EC16114F677316C90C99C2C6
A3040F7C0A00A97FAB5DC3FD8F73A54B3DE9BD70651D536F19E3F87AFE95B1ACDD3B69BA80B5B8F28DB40CC6546E4381903F4E7EB5CFE81A75CF87BC
6AFA6594924BA65C40662AED9F279C0FC7231F8D5DF8917CB63E14B8452164BB7118E3AFAFE828020D2750D5357F073EA37570F686281CA489C190A8
3F393E9C7F3AD8F08EAD26AFE1AB5BEB9C09086576EC4A9209FD2B23C453AD87C3FB7B2B007CDBA863B78514F272067F4CFE7538D32EB42F8713D8C2
DBAEA2B490B107B9C96C7D0138FA5003AD758BAD76EAF5F4D9BC9B0B3E1260B93338EA39ED51F82B5BD435EB9BFB8B92A2D622238D40FE2EE6AAE8AC
74AF86D08B38C4B77730B98D2319DCED9E7F018CFD2A4F018369A1E956D0104DC09279F8E40E40FD714058B579ACDC6A9E236D0B4A731ADBFCD7772B
D500FE11EF9A77853569351D43578BCC99EDAD644585A6FBDC839FE5599F0EECC5ACBAD5D5CB1137DA4A39638C28C9C9FCEAF6B37777278735DBBB28
D7C968F6C0506D661C877CF53D491F4A00175B9F593A84F6376B69A6D96545C051995C7A13C63A7E62B5FC377B75A86896F737A9B6671CF18DDEF8AC
2F09E8F6161E11B4BABCCC8244170E19C94CB723E5CE3B81D2BA8B069DAD55AE62489CF445FE11D81F7A00B74514531084060410083C106BCCBC1EF2
5B78AB56F0FC53A8B5791E42CBF79B1D81FA576925E8D5AF27B0B0B8DA9001F689939C13D147E4735C8EB1690681F1234ABBB7096F6F7636B04181E8
73F5C8A4C68F458A3486358E250A8A3000ED4FA28A620AC4D775E5D32686CEDAD9EF2FEE33E5C0871C7A93D856DD71FAF16D1BC530788268249ACDA0
F26428B9688FAE280449FDA7E322991A0DA86F4330FF00E2A96DF50F18493C2971A35AC31990798E25070BF4DD5723F18F87E440DFDA48B91D0AB647
E94ABE2ED09E558E3D411999828C29A9F98FE46F5148082010720F42296A84158DA85BFF0069CEF6EEC86CA253E68C8F9DF1C29F403AFD6975BD6ADF
4C8D231347F699DC4712139E4F73E80525BE81656F6CE25DEF2382647F3580271C9EB480E73E13393A05DA73B16E89538EB9515DE5703F09A441A35E
DB06F9E3BA2C47B15001FD2BBEA10DEE1591E24D6E1D0B4C6BA906F95BE5863EEEFD856BD719E30D325BCF13F87A5DD88567C30F42086FD714C45DBA
B8D534DD09B57BDBACCD1C6B24B6E54041D32BF5AB7AACD35CE82350B3D41ACD05B9989540F905723AD6478FBCFD4ADED7C3F60775C5E4A1A5F448D7
BB7A0CE0FE1563C561746F87D716C84B08EDD6DD4F739C2E690CA7F0DECEF13438EEE5BE2D0CEED2791B06327BE7AF5E6BB6AC8F0A5A9B2F0C69D6E4
E4AC00FE7CFF005AD7A62615CAC1ABCDAFF882E2C6C25786C6C0E27990E1A47FEE83D8707F2AE9A50C62709F78A9C7D6B8FF008710258787EEA4BA3B
26FB4BF9C5CF3C600CD032F7877569A5B7D4AEAF96686D2DE4FDD34E72DB40E79EFDA9DA2DDDE7882092FDA592D6D5988B544E0B01FC4C7BFD293C57
E7DFF83AFF00ECD13E594EC00E0BAE7AFE3E951B5CC7A3781215890893ECDE54483925C8C7F3A4059F0C6BEBAAE8D2DDDC108D6CEEB31C600C739FCA
AAE8DAADFEBA973A9C2E21B089996DD0AFFADC7527DBB562CFA2DD68BF0BEF21E45D4F89AE077192323F203F5AD2949D13E1CC56D6EBFE90D624A803
3F315DCEDFF8F1A0637C0F2B6A7A8EB3AF4AAAAB34A228C83C0551CFF4A9B55F105CDB2D9DCC132ECBABB11450941F3A67939EB4EF04E2C342D22C10
2349344D71311FC2A7241FD547E1546ED12FFC7DE64CC058E8D6E1F66303791D076340753A6D73506B1B68D20C1B9B971142A46793DFF0AC4D0B59BB
BFBBBED2DAE7CC9EDEE5E332ECC6D8D7827EA4F02B5A49C0D246ABA95BC624B6469A351FC3C71F89AE67C128349F096A1E20BA4C4D725E6C9EA547DD
1F89CFE740BA1A9A1EB936A335D69D0CBE6CD6D74F1B4E46408C1EBF53C8156BC2FAACDA8C9A9C32BF9A96972634908C161FFD6AC1F05DABE8DE04BE
D67616BABA8DE7E7D141DBFD4FE35A3E0B85AC74FB2B5465796E15AEEE9BA9F9BEE8CFE54016B54D5E5B8D7A0D074B976CC47997328E7CA41D87B9A3
49D4DD6F7538DE696E2CAD36EC99864B3721941EFC8AC9F0FE931DCF88F5EBBD48BF98936DDBBF6A95C67271DABA7D34DBDC45325BC300B3560A8153
8623A923F2A00C89F5AB8D3F5FB2B6BDB9056E60691E154036738503BE7271F854F3EA97569E28D3F4D9265905D46EEF1AA81B31D39EB8FF000AC8D3
211AEFC45BED4655CDBE94A2187D0BF3CFFE847F115359A799E2DD5F5BB971FE883ECB6C8C7A9FFF005FF3A066BF8B355974AD363962F904B2AC6F2F
FCF307BD546BFB8B6F11E9BA7DA35C5DC13C65A791CEE5031C1CFF009EB497B35DEA5AE26837B6503593DBF993B127AFFB3F8D57F09C771A6EBFA8E8
EB335CD942AAD1B93F70FF0077F0CD0217C4BE2192DB49BDBDB791A15B799618081FEB5FBE7D8735D13DF2DAE922FEFB112A442497FD938E47E75CBF
8A6D86ABE29D2F4B99912C6056BB9FB6704F5FAFF534DF17CF75AAF826E668D36C6D723CBDBC6E8B3804FD4D01634F4DBCBFD674D93551335AC4DB9A
D62001CA8E85F8EF8E952693AF9D43C3EB7A881A7DFE485ECCFF00E79A4D5254D17C2020B6059C5BAC1028EAEC463A7EB4EF09697FD8BE1CB786E462
500CD216FE16239FD3FAD0052B5D66F60D6AF3489655B9B9C2188E300139CFE038AB3A56A5747C517BA45C5C7DA04302CBBF681B58E38E3EB599E098
7EDBA96ADE22997FD74CC90F1D1475C543E1996448751D6822497FA9DD9482327A2838E7DBAFE42819DED14D5DDB46EC6EC738E99A7532486EF8B498
FF00D336FE55C6FC2D5C68F7C7A66EDB8FC2BB1BB655B4999C80A11B24F6E2B90F85EA3FB22F9C3860F74C401DA90FA1DB514514C466EA7772A4B0D9
DB10B34D92CE7FE59A0EAD59BE0ED5A7D56D6EFCE6F316DE731C729EAE3D6AC78AEE8D968377730A0370C9E521C73F31C62A878751F48D3EDB4C8112
478E06B8B96F463D00F7CD2EA31D777DFDA1E25974992EDAD608230D8538695BEBE82A5F086A579A825F2DD2CC63866DB13CA30585605E58CFE25D16
E759BB88585DDB96303A1C6E55F5AEB3C2D717375E1EB39EF00133A64E0633EF40DEC6BD14557BEDFF0061B8F2C65FCA6DBF5C1C5324C18B51B9D775
A9ED74E98C1636476CB32FDE91FD07B0AA3E1CD4F55D4BC537914B3E6CEC54C6540FBEDEBF5A6F85DC68BE1069B1E65E4F23B88C72CCDD31F8556F05
9B9B7B2088DB6F2FAED9E60464A20EF48A3BFA28A299272DAEEA92BCD7715A5D3C11D8C06596441C97CE029F6A8BFB56E2FF00C3726A91CD2DBC76F1
82A41C798C0739F6CD33C7F2A5BE91F61B4454B8D4E608DB472C3B93FE7D6A1F14DA7D8BC2961E1FB4FBF33245F8773F9D4948B336A93DF787A6D57C
F7820822054AF1BDC77FA6715BDA1DD4977A25A5D5C70F24419B35CE789EC0C1A2E97A1DB3048249047231FEE8E4FF005AD5779A6D1AF85A710A45B2
D48EE14609A68454B0BE9FC49AA5C35BCD243A65A36C568CE0CCFEB9F41557C277DA8EA7E22D5259EEA46B2B5610C71FF093EBF5A76821B43F05431C
0A66BB943308D0F3BDB279FA62A978316E6DB4AB4B7B67D973777325C5D640631A038C1F4CE05007592ACB7D725239E486088E19A338677F4CFA0AC4
82E751D5B4EBD6B79E556B5CC70F94D832B0EE4F7E95227899EF34ABD923B52B38768ADE2CE5A4ED9C543A03BE97E14B7853FE3FEE0B623EA43938C9
FA628034ADDF525D36CE1BA2C927945EEA63FC239E01F5A8FC1F7D71A869D34F3CAD2A79CCB1331C9C0A7F8A1A787C312C3192D3CAAB0EECF73D4FF3
A7F87A34B4B58B4FB68B10DB4404921FE290E323DFBFE94074337C45ABCF676C35259DA3822B908888789547DE27D47071525DEA73DAE9D6DABDD4EE
8B2CCBB611F28119E4023B9C73FA551F1746BABF88348D0557302BF9B385ECA3B7E43F5AB5E2589756F10697A3E0985099E7DA3A01D07E940C9758D4
2E20D2E0D47CC92192499447086C02A4F191DCE3AFD6BA19AE12DED5AE2E088D11773E4F4AE77C4500B8D7F4C5B83FE876CAD3B2E38247414FD7D2F2
FF00C2730C61E520E1783B09EF40834CB9B8D52D67D6649248A38F79B6881C0DA39C9F5CD1E0BBBD4350B1BABDD4256712CE7C907A2A8F41516A3249
67E151A7692A1E516E031073E583D73FAD5DF0E0921821B35656860B750CCA3AB939FF001A1019E8CDA97C426C13E4E9D060E38F98FAFEB5D7570FE1
E32DADBDEDCA737D7D77B101EA141EA47A735DA461846A24397006E23B9A10992561F886E27874EBC9924785604CA15382EDF5F6ADCAE73C5B19BB5B
2B1276C52CBBA5278C2AFF00FAE98D1A5A3DCC936896D7374D8768B739358BA4DE5C7892F6EA759E582C607F2E211BE3791DCFAFFF005AADDF79D7DE
1BBB4B6431C4136C231F3328FF003FA557D18AE8DE12B786DD7CCBA68CB08D48C973CE4D2029786DA4D4BC65A9EA129FF8F745857FAD6B6B7AABC5A8
5B6956B2A4534E0B49231FF56BFE26B17C266E2D6D6D208597ED17970D717408C9551C119FC2ADE9FA624DE2FD56E2FD03AA6D6883F403D714017347
9253ACCF1DB492CD62A9F34921CE5FD454B71A93DDEBCBA4DB48C8B1AEE9DD783F407B55FB2B8FB43B886255B68FE4538C163DF03D2B23C35035BCFA
85C5C92D73248783F7881D85022A3DC5F4BE32834EB7BC98C1065E4F9BB7A1F5ABFE339D9347F2226FDEDC38455FEF0EFF00D2B363B6BA8F50BD9E04
717D70E1413C6D1D49F6A9F568A5B8F14698B72479302798C46704E39A0A3A5B38BC8B3862031B100C7E153D55B269258DA6763B64398D4AE0AAF6AB
54C90A28A2800A28A2800A28A2800A28A2800A28A2800A28A2800A28A2800A28A2800A28A2800A28A2800A28A2800A28A2800A28A2800A28A2800A28
A2800A28A2800A28A2800A28A2800A28A2800A28A2800A28A2800A28A2800A28A2800A28A2800A28A2800A28A2800A28A2800A28A2800A28A2800A28
A28021FB2DBFFCFBC5FF007C0A4FB25AEF57FB3C3B94E55B60C8FA54F4500452410CCCAD2451B95FBA5941C7D2912D6DE393CC4B7895FF00BC1003F9
D4D450056FECFB20722CEDF3EBE52FF85492410CBCC9146E7D5941A968A0067971F97E5EC5DBFDDC71F95247147167CB8D573D768C549450057B9B4B
7BA409730A4AA0E46E19C1F6F4A5B7B5B7B552B6F0A460F5DA3AFD6A7A2802B3D95AC92B4B25BC4EEC006664049C5392DADE24644822446186554001
1EF53D14015E1B3B78203043046913672817839EB9F5A68D3AC571B6CAD863A6225FF0AB54500531A5E9E1CB8B1B6DCDD4F94296EB4EB2BC9049756E
92B01B4161DAADD14010ADADBA42B0AC118893EEA6D18151369960E086B1B63918FF0054BFE156E8A008E186382258A1458E35180AA3005494514014
2EB48D3AF6E05C5D5A4724A06031CE71518D034801B6E9F0AEEEB85C66B4E8A00AF2D9DBCD646CE4895ADD93618CF4DBE94CD374EB5D2EC92D2CA311
C28490BF539AB74500150DC5B457303413A078DBEF29EFCE6A6A2803253C39A3C7379C9A7C4B27F78673FCEA7B9D274EBC91DEE6D23959C618B679ED
FD2AFD14019B6FA26976D2C72DBD8C31BC64942AB8DA7E95A54514005539F4FB59E532CB082E460B025491EF8EB5728A0086D6D60B380436D12C71AF
45150DCE9D6577711DC5C408F2C630AC7D3AE3DC55CA2802ACB636B35D4773242AD3C780AFDC739A75CDAC1790986EA25963C83861DC77AB1450055B
5B1B5B5CB410AA33756E4B1FC4F34FB8B682E542CD1AB01D091C8FA1ED53D140156DACADAD599A0882BB0C16C92719CF53F5AE47C4D243ABF8CF43D2
3CC8E48632D712A75048E80FE00D7692C492C4F148A191D4AB03DC1EB545742D255C3AE9D6CAE33871180C33EFD68008746D312F16FA2B541301F23E
E240CFA0CE07E15A2402082320F514C8A28E189628942228C003B5494019F65A559D948CF6B0842738C1242E7A803B73CD4967A7DA589636B02C65FA
9193FCFA0AB9450066DCE8B61732B4B344C19CE5C23B287FA8079ABCB0C4B0790B1A88B6EDD98E31E952514019B6BA258DA8448D5CC71E4C713B9644
CFA035A5451400551BFD2E0D40A99DE6DAA31B5242A0FE557A8A00C4D3BC33A6E97233E9E26B7DE417093361FEB4B79E19D36FA749EF566B89233946
9252767D2B6A8A07733EDB4A8ADAECDCADC5CBB95DBB5E52CB8FA5685145020A4201182320D2D140151F4EB090112595BB67AE6253FD29A349D3BE5F
F41B61B0E5711018ABB450014C74592364719560411ED4FA28031BFE117D1090D2584723AF21DC927F9D5D974CB39ADD2DE4803451FDD5DC78FD6AE5
140198342D2D18B45651464E33E5E533F5C75AD3A28A002A29EDE2B88F64C81D73919EC7D6A5A28028325969B15C5F18F6E1774AFCB3103EBCD72BE2
6BFB7F12D8E8F6360E596FEE43B023958D49049FC6BB5923596368E450C8E0AB03DC1AA963A4D85836FB3B658DB1B7392481E833D290D17238D638D6
3418550001EC29F4514C4159F2E8FA7CD399E5B556909049C9193EB8CE2B428A006ED5DBB7036E318C718AA51E95651DC24E916190614124807D707B
D5FA280239A28E785E2954346E0AB29EE2AA41A4D9C313C7B1A45923319F3189F90FF0FD2AFD14019DA5E8F65A52B2DA464678CB1C903D07B532E343
D3EE2F8DE4B0932903700C42B63A6477AD4A28020BBB582F6D9EDAEA3124320C329EF514BA7D9CBA77F67BC0A6D7684F2BA0C0E83F4AB9450040B6D0
25A0B458D44013CB083A6DC631F954765A7DA58822D6058F2003824F03A75AB74500675EE8D637B2192546566FBC6372BBFEB8EB5762892189628942
A28C003B549450053B0D36CF4D5916CA111095CBBE3B9A6AE956697AD78B08F358EE3CF05BFBD8F5ABD450053BDD3AD6F8A34E877A7DD6562A47B645
3ACAC6D6C6364B68826F3B9CF52C7D49EF56A8A00A179A4D8DF4AB35DC1E63A2ED07730E3D383CD5A92DE196DCDBC91AB4257694C718A968A00CEB6D
1ACEDE65982BC9227DC3236ED9F4ABD222C91B238CAB82A47A834FA2802B59D95BD8DAADB5AC4A9128E1477FAD57B3D1EC2C66696DA0DAE4920924ED
CF5C7A568D1400514514019D2E8F633BCAF346EFE69CB032363F435557C2DA4448C96F0C90231C958A56507F5ADBA28B05C8ADE08EDA0486204220C0
C9CD4B45140105DDAC3796ED05C26E46EDE94CB2B2B7B187CAB68F683C939C927DCD5AA28033E4D1EC6490BBC4704E4A063B49FA55E550AA15400A06
00038029D45001451450054874FB4867F3A281164E791DB3D703B53ADEC6D6DA5925B78123793EF328EB5668A0028A28A00AD7161697534335C40924
909CC6CC39534E96D6DE69E39E5891A58B3B188E56A7A28020BAB4B7BC4097312C814E467B5491A246812355551D028C014FA280208AD6DE191A48A1
4477FBC546334B15ADBC2CCD0C31C6CDF78A2004D4D4500548B4FB486E0CF1C08B29CFCC074CFA7A510E9F6904ED3C5022CAD9CB01EBE9E956E8A006
48892A14911594F50C320D2A22A285450AA3A0030282CA18292013D067AD3A802116F08B86B858904CCA159F1C91E99A510C4263308D048460BED192
3D3352D140114D6F0CFB7CE8D5F69C8C8E94F2030218020F041EF4EA2802082D2DEDD192189555BEF0EB9A95116350A8A15474006053A8A00856DADD
6532AC31890FF10519A9A8A2800A8A68219C013469205391B867152D140080000003007402ABC1656F6EEEF040A8CFD48FF3C7E15668A00AD6D656D6
858DBC2A85CE4914B25A412CA249214671DC8FE7EB5628A006801400A000380076A688E31219022872305B1C9FC6A4A2801B8192703278CD23468ECA
CC8AC5795246714FA2800A28A2800A28A2800A28A2800A28A2800A28A2800A28A2800A28A2800A28A2800A28A2800A28A2800A28A2800A28A2800A28
A2800A28A2800A28A2800A28A2800A28A2800A28A2800A28A2800A28A2800A28A2800A28A2800A28A2800A28A2800A28A2800A28A2800ACB4B4D4175
E92E8DE2B583C6145B9CE5580EA2A7BF925DA90C21D7CCC86957FE5901D4D65F847519B50B3B9F3A4695219D92395872EBEF48674345607887537D2E
E74E74938B8B81098F1C107927F951AE6A53D85C69D240D21FB4CA23FB3EC196CF27DC1A61637E8A28A041451505C5CC706C0E7E691B6A28EAC7DA80
30B5AD7DACFC4DA5E93091BAE5B32E467E52703F91AE92B83F122FFC5CDD00A8E5A3E4FE2D5DE521B0A28A86EAE23B4B67B898911C63270326988C6F
10DF6B7677368349B1173139FDE9DA4EDFF0AD9794456ED34C368442CFDF181935CE59EB5A87FC259169973E5343736E67408398C73804FE1FAD3BC6
F2EA36DA35CDC5A5CC490084A3C4D1E4B13C707B70690EC67DB6AFE24D465B2D434F8A37B09E66431EDFBAA0E324FE15DBD73BE10B0D42C74BB78EE2
E524B530298A309864279E4F7AE8A840C28A6B30552CC405032493C0AE6D75CBCD423BCBBD39235B2B327E7719F3F6F2C07E54C474D4573BE13D7E7D
7ECEE2F65B7582DD64DB11CF240EB9A7CFAECB06A5631B4518B7BD62A99CF9831FC447A52B8EC6FD158D77A9DCD9EB5696D71147F66BC731C6CA4970
C077ED8356759D56DB47D3E4BCBA276AF0AA3ABB760298876A32DEC5133D9C703ED524890904FD3159BE12D7A6F10D94F732DB2C023936280739F5AA
3AE7882FB4EF0E1BFBC8A18C5C4456355CEE491972A0FE19AD1F0658FD83C31651B0C3C8BE6B7B96E7F96290CDDA28AA9A8DFDB6996725DDE4A23863
1C93DCF603DE988B74560B6AF769A3BEACF6F1C76C8A2411B64B98FD73D01C76C55E7D4524B3825B45F31EE941857D723393EC3BD00685158DA56B3F
DA114D1379497914AF114192A597BE7D29FE1FD54EAD6D70D222A4B6D70F6F26C395257B8FCE8032B5DF13DFE89E4FDA34D43F6893CB8B6CB9E73DEB
A8425915986D240247A5717E232353F1EE8BA5FDE4B606E645278F51FCAB775BD56E6CE1B83611C52496B1F9B319890AAB8240E3B9C5219B54573D7D
E22FB3F82C6BAA8AAEF024888C72373638FD6B4F48B8B8BBD22D6E6EE311CF2C61DD476CFF00F5A988BD45154AFEE9E0465816379C2170AEFB4051D4
9F6FF3DA8033343F1049AA6AF7F632594907D91B01DB3F3738AB5ACEADFD9EF6B6D0C6B2DDDDC9B2242703DC9F6A5D23533AB68316A56B1059264665
463C6E191827EA2B8DBCBBD5EF3E20D944FA7DB8B9B280B0884D9539E739C71486757E1BD5AF75337A97F67F667B69BCB180407F7E6B72A384BB448D
220472A0B2839DA7B8CF7A9298882E65FB3DACD3052E63467DA3AB60671597E17D6DF5DD39EE64B47B76490A10DD1B1DC547AE6BD269F1DC359C2938
B45DF72CCD8083B28F56357FED727F674375636A2512AABF97BC26148CE7A521D8996FAD5AF9AC96E23372ABB8C59F980F5C7E22AD579E7815AEAFFC
4BAAEB4F69FB9B890C624320263C761DCF615E874D098514560FF6E9BD9EE23D2CC461B607CFB993EE21C678C75A00DEA2B1344D7EDF51D0DB5494AC
11465964627E5F97B8F6A87FB76EDB4D3ABC1668DA72A17F99C895907F1818C018E714058E868AA4DA8C3F6086EE305D670BE528EAE58640AA167AEB
5D5ACBB6D81BC8EE1ADCC2AD91B877CFA62803728ACAD0F571A98BA8DE3114F69318655072323D2A0BBD6667D61B49D3228A5BB8E312CA657C2A29C7
1C739E45006E5739E28F1543E1B7B6596D649FCF273B0E368F5E9CD59D1B5F83518EF7CDD91C960E5272AD94EFC83F81A86D75B4BED5A1D36F34D319
9A13710B3B07CA83DC63E53480D9B59D2EAD62B88C3049503A861838233C8A9EB2F58D5534C481153CCBABA93CA822CE3737B9F4155ED75B3FDBA347
BE8A38EE9E332A796FB811E87D0F19A60695F6A167A7C6B25EDC242AEDB54B1EA7D2AD579CF8B67BED5FC45A4E9474C21E273398FCD53B87D7A01C57
67AA6A8B6096F184592EAE9C470C5BB019BBF3E8295C7634E8AC2B6D6A45D75347BF5856E6583CD4F29891C75073D0E39ADDA620A28ACA9B598BFB48
E9D671B5CDCA0DD22AB00B18F73EBED401AB599ADEAF0E8F6427951A577711C7129E5D8F614EB0D4E3BD373184649ED5F64B164120F6C1EF9AE2BC45
A8EA977E2FD1AD4E95207B5637020F3D4F99EF91E801A0674BE1CF11BEB3797B67716125A4F6846E0C7239FEB5D0D55B26796DD679AD45BCD28CBA64
123D32475E2AD5020A2B2F55D5858C7288A233CD1446674071B507727FA544DAF466C61B9B5B76B8692013B22B01B131D493401B34555D3AF62D46C2
1BCB7CF9732EE5CD524D685DDE4D6DA643F6936E76CB26FDA8ADFDDCF73401AF4567E99AA41A85A493AFEEFC9768E65623F76CBD4135043ADC4DA97D
92589A20D09991D9BEF2038CE3B501635E8AC6B7D6F7EB11E9D35ABC2D3C6D242C48F9947A8EDDEB66800AC5F12EBF1787ED22B896DA59FCC7D8163E
D5B0CC114B31C00324D61C7AEC571AADAE9F35932ADDA34B03B907217B95EDED402362DE65B8B78E640C16450C030C119F5153553BEBD5B48D888DA4
90297D8BD768EA4FA0AA9FDB91C9A7C37569034C658CCBB030055475268035E8AA5A7EA36FA869B1DFC4DB62752C771FBB8EB9FA62B2C78953FB52C6
DBECEC2DEFF72C1296E58AF538F4FC6803A1A2B9ED13C4A355D56FB4D96CA5B59AD4E7E73F79738CF6C73576F7588EDEF92C2DE17BABB65DC634206C
5F563DA803528ACAB2D6E0BDD38DD471B07593CA3092376FCE3159D278B6DA3D22E6F1E1659ADE5784C20EEF99700F3E9C8A00E9AA082E6DEE431B79
E29769C36C70D83EF8AC4BBD74FF00C220DAB8B59489206621481B3B64E7FF00AF543E1EA3D8F860B5E5BBC19633195F1FBC0DC83FE7D68035B59F11
DA68D77696B7314CEF74C154A282064E39E6A6D77561A2D89BC92DA59A14FBE6323E4E98E0D56B3D66DB50D69F4EBAB168AEA05F362F35431DBEBFEC
9ACBF8997A21D021B151992FEE123DA3A950431C7E2147E3486749A46A0BAA6976F7C91B46B3AEE556EA066A8EB7AF3E8C1E59B4F9E4B540099D1860
67DA963D463D35AC2C64B478209F6C3031619CE3A11DAB3BE235CB45E186B58CFEF2F2548547AE4E4FF2A03A9B9A26A71EB1A5C57F0C6F1A4B9C2BF5
EB5A15474BB44D3749B5B41F2AC112A9C9E981CD644DE2DB6485EEE3B59E5D3E393CB7BA1C2E7206403D464D0074B45351D5D03A9CAB0C83ED4EA620
A2A8EA7A95BE9B02C939CB3B6C8E35FBD237A0AAF6BAB897521A7DC5BB5BDC345E6A2960DB97BFD0D0051D57C5F67A53C697963A8234ADB63CC380C7
D8E7DEBA2460E81802030CE0D711ADAFF6BFC48D2F4E215A1B088DC4808EE7FF00D4B5D76A17F6FA75AB5CDD310808000192C4F400773486CB74562D
96B6D35FC76575652DACD3219220ED9C81EBE94DBEF11DAD95E2DB98E49079AB13C8B8C231E829858DCA2B0E5F10C31DEC31F91235B4B28805C745F3
3D00EE3DE9977E248ADAF5215B495EDCCA217B8E8AADE9CF5A570B1BF452121412C4003924F6AC293C4B0859A782DE596CE06DB25C0E141CF38F5A62
37A8ACCD4B564B28E3F2A16B896552EB1A1C1DA06493E82ABFFC241149636F3DADBC93CB7085C42A46401D7268B858DBAAF35DDB5BBA24D710C4CFF7
43B852DF4CF5A834BD4E1D574C8EFAD7732383F2F19C8E08FCEB8ED4AE5F53F8896102D9DC3C562819D4280413D09ED8E94AE3B1E814565EA9AA8B19
61B686133DD4F9F2E2071F8934CB1D6E39F526D2EE23F26F5137B206DC08F634EE2B1AF4514500145145001451450014514500145145001451450014
51450014514500155EF1A45B66F243173C0DA338AB145006669566F0EF966DC5D8F1BBA8AD3A28A49586DDC28A28A620A28A2800A28A2800A28A2800
A28A2800A28A2800A28A2800A28A2800A28A2800A28A2800A28A2800A28A2800A28A2800A28A2800A28A2800A28A2800A28A2800A28A2800A28A2800
A28A2800A28A2800A28A2800A28A2800A28A2800A28A2800A28A2800A28A2800A28A2800A28A2800A28A2800A28A2800A28A2800A28A2800A28A2800
A28A280327C4F3CB6FE1EBD96DD4B4823C2E3DC81FD6A2F0C5A0D3F46B5B4501CAC7BA4914F1B8F6FAFF009EF5B0CAACA55802A4608238229238E389
76C68A8BD70A30290CE5B57B66D57C6FA6C1B730D8A19A4CF439E9FD29F7F1B5EF8FAC2170562B381A6073F7D89FE9C574C228C4A6511A8908C16C72
47D690C3134EB3B229954150F8E403DA8B05C7A36E40D8232338230453A8A2988826BAB7B7CF9F3C71E064EE6038AC8B46B1F12AB5D87768ADE568E3
0AD8E4756E3D7FA5695E699657CC1AEED9252063E6F4A659E9363624FD8EDC43904108C40FCB3486715E248A0D3BC79E1F951DB04E1B7C84E39C77FA
D7A2563DC786F47B97DF7164B2BF66676247D0E722B46D6DA3B4844508217AF2C58FE66840C9E9922248856455653D430C8A7D14C470DE172F77AC6A
BE2195774B24BF6582327008C8E87F0FE7563E24CCDFD896D668326EEE513DF8E6B60C1A5787E196F243E4C464CE589210B75C0ED599ACC3FDAFE2DD
2204F9EDAD14DD3B2F233D81FAE07E748AEA74D6B17916B0C39CF968A99FA0C54D4514C939DF1E4D3DBF842FA4B6243E14123D0B007F9D64C8D1E85F
0C1110912CD6DC63A96719CFE46BB1BAB68AF2D64B69D774722ED615562D1ED23B316AC1E58D63F297CC6C90B8C60521A662F83E26B0D274FD29EDC9
12DB19E52790371E86A9DE46FAB7C49B3B7898A5BE910891F6F4DC7B7EAA3F035D4E9DA6C1A745B21691CE31BE46DCC00E833E94DB2D26D2CAFAEEF6
056F3AE9B74849CFE54582E6048A750F88524B2FFC7B68F6C0F3D37B0CFF002FE554BC5C926A7E29F0E5B0DCB6D2B79A411DC1CF3F80AE99F458CEAF
26A093BA0982F9D08036C857A66A5D574B8F525858CAF0CD6EFBE3953AA9EF405CE5FE213A6A173A4684012D71768EE47F08E47F226BAAD3EF12E24B
982185962B57F243766207207D3A555D47438EFADEDC19DD2E6DDC48972065F3DEAFD859C763689045921792C7AB13D49A03A16AB8BF88B6F3DE2E8D
68A48B79AF556523A724019FCCD76954754D3A1D52CDAD672C01219594E0A30E845362461F8F676B6F0BBD8DA2666BB2B6F146A3271DF1F80C7E356F
4FB35D0FC336CD7855A6D3EDD8EFC9C02464FF00855AB2D1D609D2E2EAE66BC9E20446F29FBA3E9EBEF56355B08F54D367B195D91265DAC57AE33486
71DE1257D27C1B7FAE4E4F9B701E6407F8476FCCD6878223FECDD16C619416B8D45DEE5FFD8046727F2515BEFA65B3E8E74A2A7ECE62F280EF8C75A8
F4AD2974F1B9A779E6DA103B8C61476028B05CE5B41996E3C57E23D7EE037936ADF668C819CEDE3E5F7E07FDF5563C73231D0A2B7B347497569D11C9
FBD820707F01D2B4ADBC3115BDFDD3A5CBFD8EE67FB435A84006FF005CFA67B52DFCBA5EB3AA368E6674D42C88995947FAB381C83DF83401CE78C631
25C681E13B71889DD0C9B7D071FCB26BD0800A000300702B8CD3AC12F3C72D76AAEF0E9D114F39F9F3253D4E7F3AED6840C2B91F1BC8BA7699733DB0
2D7BA96DB4504E700F0702BAEAC8D7F454D662B71E71865B6944B1C817760FD286088FC3B0C5A659DBE888199EDA00F23EDC2EE27247D7249AC1F0D3
7DBBE237882F31958008013D41185FFD94D74F650DBE9B0CC1EE04937334EE7EF1F723B0AC0F87168E9A65E6A7383E6EA372D26E3DD474FD4B500763
45145311C2F8BE056B9B4F0F58AB799ABDD19EE5B233B0104FE1C13FF01ADAD435382D3C3BA9B5B2622B188C48D9F95885C601F63C54BAB6822FF54B
5D46DEEDED6E6DD1A30EA8189539CE33D0F279F7AB173A3DACFA1CBA511B61910AE4F273D771F539E690EE6378342E95A16976322335C5E069D80FE1
1D727F4AEB2B2749D2DACD564BA984F72B18895C0C0541D001FCEB5A840CC7F15DC5C5B7863509ECF89921241F41DCFE59ACAF0369F6E3C0D025C286
5BA56926DC719C93DFE82BA89A18E785E19903A48A5594F420F5159369A17D96CE3B14BC9458C79DB10500E324E0B0EDCD01D0E63C63147FF08F691A
6E9919B6B1BABA11B0030719C0FCFAFE55BBE36BA4D37C1F730C20069A316B0A019CEEF97007FBB9AD5D474B82FEC56D9BE4119568980C9465E84540
34733DF4375A95C7DA9A1E6141184443EB8C9C9A00AFA258FF0064F86EC1AFDB32D85B33364E02E464FE438AC1F021FB2689A97882F19962B891E554
F450724FE278FC2BADD66C0EA7A45D58ACC62371194DE0671FE7A547168F6F1F8757473FEA441E493F8727F3E68030BC15B6CB485BB9C3F9DAADC3CC
0019217D4FD00FD6B2751692FAEA6F187856490C96CC62B98645C095540C91CF231FCABB0D27496B38A35B99C4EF1C7E526136855F61FD6A945A05E5
959CBA7585F243632B31C98B32206EA01E9F89A006F872D749D4BC362E2D2CFC886FCF9934609F998372393D320FE06B23459D2F3C53ABF8867DCB6D
623ECB0285E4F6E07E1FAD75FA569F0695A745636ABB62886064E4924E49FCCD64E95E193A7CCE8F7024B25B86B88A12BCEE3DD8F7C50052F13C56BA
F6A69A658DDBC3ACE9E3ED10B63E51D3827D7A547E0F92DF57D4EF2FB50B3F2B5CB3C4170DD01E080C07AF047E02B664D1E6835BB8D534E921596E63
092ACAA71C74208AB1A3E949A77DA2566592E6EE4F327902EDDC7D31E839A02E739A748977E3ED6B549083069B08B7539E43639FE4D563C4C2D75DBD
874CB1BC68757B41F6980ED3B47B13DB22ADDA78725B6D52FDBED286C6F27FB43C617E72DFDD27FBB576E34D9D3517D434F9214B8950472095090C07
4E41C8A00C2F084D06A7AADCDCEA56423D72C4086693B11C8C8F7C0E7EB5DA56568FA4AE9CF737124825BBBB7DF349B71F80F615AB420655D4677B6D
3EE278D4B3C713328033920715CB7C358F7E8771A94EC5AE2EEE19A476F41D2BB1601810C0107820F7AE7E1D06EAD6D65D3ECEED61B09189E17F78A1
BA807A7E340191E0879750F12F883561936D2CA228C9EF827F90C7E7496ACD7BF16EE9DB94B1B30AA7D32067F566AEBB4EB1B4D32D63B2B3458E38C7
0A3A9F73EBF5AE63C071FDA6F75CD65811F6BBB2AB9ECAB93FFB30FCA803B3A28A2988E3FC6ED1E93A0EAB76B2B1B8D47642B9C70318C0FC377E7599
A85B49E1BF8722D8EEFED0D43CB818F248661F77F05047D7EB5D3EBFA02EB5756124D36D86D25F31A3DB9DE78FF0A7F88B4897555B27B795239ACEE1
674F3065491EA290EE67CEE744F06DC58DA93F6AB1B205B68240623D7EB9A6F82A18B4CF05437521C9746B8958F527AF5FC2B66DB4D58ED678EE5BCE
92EB267623862463007602AAFF006213A526926E3364A3070B8765CFDDF4C7BD0073DE01B49351D1B549EEC4AB0DFDD175C1C640E723EA78FC28485B
5AF89B73B4E2CF4DB7585C0E8C7AEDFCC9FF00BE6BB8B7822B5B78EDE0509146A15547602B2FC3FA20D20DEC8F28967BC9CCB23818FA0FE74582E73F
6772B79E3ED535494B0B4D2A1FB386EA37771FCFF2AE83547D6E68E09343368AACBB9C5D03DFA74AAF6FE1E686FAE71703EC371319DA155C166F427B
8CD743401936706A377A64906B8205959B04DB31C15EBDEB9BD2A48F52F1C6A5AC48E174FD2A316B037401B1838FCDBF315DD560E97A02E9E668A391
4D9BCE67116DE598FF0078FB1C7E54009E279E3D3B41D435156225683CB4DDDB3C00076CE6B054DC787FE19867CFDB278C22F1C82DC01F97F3AE97C4
3A2FF6E5BDBDBBCC638A399657503EFE3B51E21D23FB5AC238A3658E586559632C3E5C83D0FB50057D16CADE1D02DB43964C4BF650655438383D4E71
DC9AC6B58DF53F89ADE5AEDB2D12DC468074DCCB8C7EA7FEF9AEAEC2C8DB79934CDE65CCC732498FC80F6154468F716DA9DDDD69F3AC62F8A99B70C9
523BAFEBF9D0022CCBFDB9A85E8894C76B6DE5B38182C47CC466B9EF03C3A9DFDA5FEA85E38DEFE63895B25D547A57676F630DBD97D946595948663D
5B3D49ACCD2F4CBDD3AC5B4E8A4892D958F9728C97DA79E9EB40096F65A5E93A719D0B14B03248EC4F2CF8F989CF7F4AC4F877A52CBA2DC5FDEA997E
DD23ED5906404C9C9FC4E7F215D26B3A51D4342974BB794422450BB8AE78CF35734EB48F4FD3EDECE2394823080FAE050072DF111B3A3596930280F7
D729122AF1851E83F215BCA2D9E1FB00958258797E6F1C300B900FB703355FC41A3CDA8CF617769244975632978FCD04A90473D3BF02AEDA69FE4D94
B1CAE1A6B8DCD33A8C6588E71ED401CA7856E8CFAA6ADE23BD6C437130B6B738C9C038E2AD6A10A6AFF112CEDDC6E874C80CCCA471BCE31FCD7F2AD0
D174092C1208AE274920B524C08831C93D5BDEADCDA6491EAB2EA562D1ADC4C823916504A903A1E28039CF124B2EA5E3FD134DB5CB0B2617337A0E41
E7F01FAD3FC447FB4BC79A3699D63B75373203D0FF00902B7F49D162D3EE6E6F647F3AF6EC832CC4638E3E503D2AADFE837327883FB5EC6E638A4783
C9937824E3D57D0D01732FE206AEEBA025B69D265EF27F20BAF4E3A8CFE5CD4973A2C89A25A41ABCC91D85A04CDB5B8C99181E012719AD3D53C3D0DD
E936F671B00F68E2489DFF00BC0E79FAD490DA5FDF5C452EACB1247072B146721DBD4FF85006CA80A0050001C003B52D145311E7B3497DAD7C49923B
73108F494F944B9C0271CE0753935D4D8D95B5BEB73CB2CC67D4268F7331180880E0003B7FF5A98BA2C969AE5C6A9A7C91F9974B89925071C7423157
ECEC85B87925732CF2F3239EFEC3D05219C9F84EE219F58F10F886E19562330811BAE157FC7E5FCA99ACCB73AA7C45B3D3EDF6B45A7A09D91CE17775
CFBF515A1A3786AEAC81B2B89226B05B8338DBF7A427A03EC2ACDFE937D1789175BD2C40F23C0209639491919CE41FC05005882CE18B5D5B9BEB9F36
FE58C8890290B1AF70BFFD7AE5FC7168D15B5A69703992F352BDF30B85C1FCBDB35D869DA73C13CB797922CB752F1B80E117B28AA373A14B79E2D835
5B8954C16A988501390DDE8039DF1A412A7FC23FE1CB094891A4524F7E380DF87CC6A4F17C4B26A7A1F86ACB72FEF04F295E4ED07AFD78635D09D0E4
97C6035B9E6568E287CB823C1CA92304FF003A6EA3A25D4FADBEA5693C51B496BF666DC0E541392463BD017337C69ADB1F0A4ADA779A1EE26FB306C0
18E79FF0A96DF44963D02D6CF549A3B7B28117CD8E3049739FE23F535A3A968314FA026996D85103078B77F781CFF5A16C751BFB888EA8638E08B07C
A88E7CC6F5340117884C3A6E89A96A409DED6DE5463A6C5230001F539AE7AC6DE6F0F7C37B9BF919BEDD7106413D5037000FA039AE9FC53A34BAEE9F
15A473AC4826579030FBC076A76BDA39D53484B28990189D1944832ADB7B1C50172AF84218F4AD034FD364622E5EDCDC38C7DDDC7273F4CE3F0ACDF0
8CF1CF7FAD6B972CAA93DD7910B7FB20F031F95743A6584D019E7BF78E5B99F86283E5541D1477C726A868BA0CD6221826963682DE4691028E5DCFF1
1FA671F850067F887FE26BA8F9FA0DC06D534AFBE98FBCA7B0F5A97C142D3516B8D61E074D459CC73EE270A7DAB5A2D2E4B2B9BA9B4E16EA6E5B7B17
53907DB1DAAC691A5C3A5DBBC71316791CBC8E7F898F7A02FA1A14573B068177178B26D5DF507781D485839E3DBD315D153105145140051451400514
514005145140051451400514514005145140051451400514514005145140051451400514514005145140051451400514514005145140051451400514
514005145140051451400514514005145140051451400514514005145140051451400514514005145140051451400514514005145140051451400514
514005145140051451400514514005145140051451400514514005145140051451400514514005145140051451400514514005145140051451400514
5140051451400514514010DD5B41770186EA149636EAAE322A2B4D3ED2CD9DEDA05467FBCD9249FCEADD140051451400514514005145140051451400
51451400514514005145140051451400563C9A05AB6B6FAB452CD0DCCB1F97298DB1B8718FA7415B1450057B4B482CEDD61B7408839FA9F53EF5628A
2800A28A28039CD47C2905EEA92DFA5EDCDBB4E9B25443C30FE95BB6D6F15A5B476F02048A350AAA3B0A9A8A0028A28A0028A28A0028A28A0028A28A
0028A28A0028A28A0028A28A0028A28A0028A28A0028A28A0028A28A0028A28A00E7AFBC38D71AF3EAB05F3C0F243E4BA85CE063B7BD6B69F636FA75
9C76B6A9B628C71EA7DCD5BA280B85145140051451400514514005145140051451400514514005145140051451400514514005145140051451400514
514005145140051451400514514005145140051451400514514005145140051451400514514005145140051451400514514005145140051451400514
514005145140051451400514514005145140051451400514514005145140051451400514514005145140051451400514514005145140051451400514
514005145140051451400514514005145140051451400514514005145140051451400514514005145140051451400514514005145140051451400514
514005145140051451400514514005145140051451400514514005145140051451400514514005145140051451400514514005145140051451400514
51400535982A96620003249ED4EAE53E236A3269FE1697C972924CE230475C77A00D18756B9D4D5DF468A278958AF9F3121588F4039AAF2F89058CF0
5A6A76AF05D4D20450A728E09E4AB7B7A53FC131AC5E10D376FF00145B89F524D49E26D34EA1A7A18A3DF3DBCC924781CF519FD3F952195BC59E20BD
F0F5BADCC7A7A5C5B9214BF998218FB62B7E17F3608E42305D43633D322B97F89833E127FF00AEF1FF003AD9B8D4ADF4BD056FEE9B11470AB1F52703
00501D0D4A2B9B9751D5A1D2575868E16802095ED80C304C72437AE39ABBA96BF6965A08D5C1DF0BAA98FDF774A2E1635E8AE6F53D5F52D334D4D55E
1B69AD1555A444DC1829C7209E38FEB5B6B706E2C56E6CF6BF991878F71C0391919A622CD54D42E67B5B632C168D74CA798D18038F5E6B9DD3FC4DA8
4FE229F44BAB0892E235C87473B4719C9AD3F0F6B0FABC776248823DB4C62257386F719A571D86785FC429E21B69E78ED9E05864F2FE739C9C64D6ED
713F0D015D2F54DAA377DB9F009C0E82AF4BE2AB9B7D7A0D22E74968A6B8198D8CE0AB7E38A2E0D6A7514562EADAB5EE9F3DB47169DF68FB41DA0A49
D1B1F4E95A364D74D6E0DF471C72927E5439005311668A28A0028AC99757CEB126956688F751C6247F31F68507A74049A8B4DD73ED5ABDC6937907D9
EF604126D0DB96443DD4FE228B858DBA2B1B5FD5E7D16D5AEFEC0F736D1AE6464700AFE1E950D87885B54D3A3BDD3AC24994A16705C2ED23F841EE68
1D8DFA2B2B4BD5E3D5749FB7D9C4CCDF30F24B00438FE1CF4FFF005D65AF8B1C6B31E8F2E91711DF48BB963322E08C673BBF03F950163A9A2B9DD57C
4579A75C5BDB9D1A5924B93B62DB32905BD2B66C25B99AD124BC816099BEF22B6EC517158B3451450014572DADF8B5F43D5ADED6FF004F22DEE1F097
0B2E4019C648C75F6A9F52F11CDA64115EDD69CCB60F2057984996453F758AE2803A2A2B127D6A531493E9B646F2DA21B9A512050DDC85F5ABDA66A3
6DAAE9B15F5B3E619173CF6F506802ED15870EB53DF8966D2AD56E6DA26285BCCDA5D875DB56743D62DB5BB26B9B6575D9218DD1C72AC3A8FD680B1A
745737AE78B23D0A6DB7DA6DE790C76A4E814AB1FCF8A9EF7C426CADED65934DBA7FB510A8232ADC9E9DE80B1B13CCB05BC93BE76C685CE3AE00CD65
F86FC416FE21B596E2DE192258DF610FDEAC5CDE5CC7A6FDA534F7924C12D0170182F3F813D38F7AA7E18D62D759D1DEF2C2CFECE8B232795C0CB003
D38E73401BB4572B0F8C63FEDC9349BDD3A7B69D572A37072E7B0007A8AB6FE248ECF508AD357B66B1F390BC52B48195B1D4123A1F6A2E1637E8AE6E
E3C4EF6B0477D75A5DC47A6C840171B86E00F4253A8078FCEB6EEAFADED2C64BD9A402DD137971CF14058B351C8DE5C6CFB59B68276A8C93EC2B125D
7A78AC46A0DA6C86C8A86DEB202C17FBC57D315B167750DEDA45756CE1E195432B0EE28030ACBC5F6779ACC7A50B5BA8AE1C124489B715D2570D7247
FC2E0B41FF004E247E8D5DCD2436145150DCC92476EEF0C466914655036377E34C4655AF88ED6E7C473E8AB14AB3C2A4EE61F2B63D2B6EB9AD1F5F87
52F10DCD849A4BDADEDBC7BA491F6938E38C8E7B8AD7D4B53B7D3962F383BC933EC8A3419673E8280285EF8B348B0256EE69226DC542B44C0B1071C7
15B48EB222BA1CAB0C83ED5C07C4EBE825F0EC50491490DC8B957449539C60E482323BFAD77DB9523CB10AAA3249E00148649456241E208EEA19AEAC
ED659ACE1621A7040071D4A8EA6B461BEB69AC05F24ABF6729BF7F603BD3116A99236C467C31DA09C28C93F4AC53E2063647504D3E67B00BBFCE0E33
B3FBC17AD6B59DDC17D6915D5AB878655DCAC3B8A00C9B6F15E9573AA2E9A8F2ADD3F4478C8ADDAE22E6285BE2C5A1551BC5A3337B9C1C1AEDE81B0A
28A281051596BAB09279561B59A58A062B24AB8C023D0679A974FD56D6FECA4BB8198431B32B338C7DDEA6802FD15869AFF9B6A2FADEC2E26B2E4F9C
B8CE07521739C569D85EC1A859C5776AE248655DCAC2801D35D5BDBC90C734AA8F336D8C1FE23E82A7AE0354D72E6E3C71A65B1D3EE123B62EFB4A8D
EF918C815D9C776D25B492ADA5C0642408DD42B37D39A571D8A5A5788AC755D4EF2C2D4B996D0E1C95C038383FAD6CD73BE1FD5B4DD4B55D462B5B17
B6BB85B170CC801639C751F4AD09B5402F9ED2DADE49E489434DB081B01E9D7A9F6A00D2ACBD475CB1D366F2AE646F336191951771541FC47D05167A
D59DEC57525B89996D4ED93F76739C670075AE2E1D5AD6F7C69AD41796B74D6D756CB1EC119DE0000741C807AD0D82476577AF58DAE87FDB0CE5ED4A
82A547272718ABD65770DF59C5756CDBA29577292315CF5EDD691A27862286E2C2E64D3D50615E3CE32DC6EF439ADED3DEDFFB36DE4B7458ADCC4AC8
B8C0552334016E8AC44F11DB491C92C36F7535B444AB5C24794C8F4E7247BD6A0BAB76B4FB5ACAA60D9BF78E98F5A6227A2B15BC456B15B35D5DC53D
BDBE32923AF0E3DBDEB56095668239933B6450C33D704668025A28A28021B9B986D2069AE6458E35EACD4F6915633213F2019C8E78AE07E206BD70B1
C5A747A6DC089EE14199D7024DA73841DF3EB5D9595FADD36C5B4BA870B9CCB1155FA034AE3B1521F13E9171A9C7A7C7704DCB9C2A9423DFFA56D570
DACB5ADBFC4BD2E5902463ECEC5DFA64F38FAD74BFDB5683518EC2549E1966CF96648CA87C7A1A6163528A82EEEA0B2B592E6E6458E28C1666359C9A
F5B34115C3C5711DB4C408E574C29CF427B8A046C51451400514551BDD4EDACE68A0937BCD2FDC8E35DCC47AD005EAA97BA95969EBBAF6E6384609F9
8F6A860D62D25BDFB1BF9905C91958E65DA587B52EB56F6F3E9775F68863902C2E46F5071C76CD003B4BD52CB56B6373A7CE268958A16008E47D7EB5
7AB88F877796D6FE168A25579273248EF1C48588F9B193F90AEAB4ED4AD75285A4B67276B6D75618653E8450368BB4567EA3AB5A69ED1473B334B29C
2451AEE66F7C7A52DB6AB6D3DC9B661243381911CCBB4B0F6F5A045FA2B2F54D7AC34A70B7F24910233BFCB62BF98A659788F4BBEB779ED667789064
B089B1FCA8B858D7A2B234EF10E9BA9DD9B6B391DE419CFC84018F7A2EB5C861BD6B2B6865BBB941974886760F73EBED45C2C6BD15534FBF8350B7F3
606E870CA7AA9F4356E800A28A2800A2B2E3D72D679654B58E7B8111C3BC499507EB9ABD6B730DDC0B35BB8746EE3F95004D456436BF6A6F24B5B68A
E2E9E23890C09B827D4E6AE5A5FDBDEC0D2C0E484243AE3E6523B114016E8AC37F1569515FA595CC92DBCCE7004D1951F9D587D72CA39E18A469104E
DB6391A32118FD4D170B1A955EEAF2DACD03DD4C9129E858F5A9C900124E00EA6B286A9A65F334617CE895B6194C798C1F4CD0059B0D52C3512E2CAE
A39CC7F7B69E950EA3AD5869B736F6D7729596E0E100527DAB98F87891A6A5AF884058D6EB6AA8EC326B6EFEFF004B6D5ADA0BCB395AE565D90BB45C
03EA3D452B8EDA9BF51433C53C7E642EAE9923729C8E2B23C53AA4BA5E8F712C16D24CFE59E40F953B649ACFF06EAD1BE85A7C1F64BA1232805BCAF9
4924E5B3E9FD28B858DBD6357B6D1E08E5BADE448E140519357A3903C6AEBF758023F1AC3D7759B1B1789352B195E33280B23460A83EA0D6D991047E
6160131B8B1E001EB40125158F16BF6B3A3CB6D05CCF0A1C79D1C7943EB839ABD05E437369F69B73E6A609C272723B63D6988B54573A3C61A67F68CB
62EB3472C4A490C9C93E8077356E1D7ADE4BF8ECE6867B79261984CAB80E28035E8A28A0028A28A0028A8AE264B785A69490AA39C0ACF935816EBBAE
ED6784160138CEEA00D5A2A0BBBA82CAD9EE2E6411C48396354A2D5D5ED5AE9ED678EDC7DD720723D71DA803528A645224D1AC9138746190C0F069F4
00557BEBB8EC6D24BA98318E3196DA326AC550D56F12CACA49A6B779A255CB60023F1CFF00850049A6EA106A766B756DBFCB6E9BC60D4FE646251117
5DE4642E79C7D2A8E9DA84171A425E5B5BBAC4572B146A09FC00AE5ECB5DB86F175EC975A75E154458E258E32767FBDEFCD2B8EC775454504BE7C2B2
797247BBF8645C30FA8A969882AB5FDEC3A7DAB5C5C1608BC7CA326ACD676AF7F0585A3C9770492C38C9DA8187D0D005D86559A18E54FBB228619F42
3352556B3B88E6B18AE117CB8D90300DC6D15564D66DE18DA59629D625CE2431E15BE9401A745416770B776715CA82AB228600F6A9E800A28A2800A2
8A2800A28A2800A28A2800A28A2800A28A2800A28A2800A28A2800A28A2800A28A2800A28A2800A28A2800A28A2800A28A2800A28A2800A28A2800A2
8A2800A28A2800A28A2800A28A2800A28A2800A28A2800A28A2800A28A2800A28A2800A28A2800A28A2800A28A2800A28A2800AE67C7DA3B6B1E1B95
225DD3407CD41EB8EA3F2AE9A8A00E7FC0D32CDE11D3F69C948F630F4209E2B43589DA0B165864D93CA44717A9627B5556D11EDA7966D22E8D9994EE
78B6068CB7AE3B559B7D3489D2E6F676B9B841F2920055FA0A40617C466F2FC27B5CEE2D3C4A4E3A9CF5FD2A97C4549DBC136ED111E5A18CC83DB1C5
6FEBFE1C835F8D61BCBAB9584306D91900640383D3DCD58B1D256DF4C6D3EE6E24BCB628230B380485C63048EB458772841A4CBA8F8761B57D5AE0DB
4D6CA876220254AE319C55578F43D33C1CF6F7664B9B0B691A3559B82CC0F00631DEB4AC345974D8FECF67A94CB680FCB13206283D031E94FD5FC3F6
7AAE8A74B72D1C60EE565EA1B9E7F53405CC3F131B89FC013CB23476911854C708192578DAA493D6B6BC3B322695A6D98C922CD5F71F418155C786DA
EB461A66AF78D731246123DABB7181804FA9A9F46D125D32DCC6F7CF3BAC7E544E500F2D7B7D7B500636923FE2EA6B60FF00CFAC647E495D7C70C50E
EF2A344DC72DB540C9F53580BE1768B589356B7D4E74BC946D91CA2B065E38C7E02B6AFA0B8B8B731DADD9B590907CC0818E3D306840CE63E1C0DB65
A98038FB73E28F88D68EBA65B6B16A08B9D3A659030EBB49E7F5C569787BC3D2E86F2ECD45E68A572EF1B44A32C7BE4549FDADA4EAF63A844B32C914
2AC93A918C0C11FE34740EA2693A8DBEBD7115FD9CBBEDE18BA7A48DD41FA0FE75B9581E0AD34699E18B58B66D670646E393B8E47E98ADFA620A28A2
803999A4D3A1F1816B48DE7D5668424986F92351DDAB32D9677F8A8C6E244731589E10600071C7BF5AD293C3D79078AA4D6F4EBA897CF40B3452A9C1
FCBE9449E1A9A1D746AF617852E2552972651BB729C72A074E9D29144FE21BC8EE7C0F7B79B08496D0B056EBC8A97C1E889E14D3820007920F0299AE
E8975A868FFD9767751C16E621136F4DCC40C639FC2934BD3F5AD3B474B21736733449B23919194E3DE8174333E1AE7FB1EFB38C1BD930076E83FA53
3E225ABC16F61AFDAAFEFF004D9D4B71D5091D7F1C7E66AFF85B43D4741134124F6F3DBCD2195986438623B0C6319ABD7975A6EAFA66A76CB731CB1C
71BC73EDE76707FC3F4A3A05F521B395359D5E0BF81C3D9C36E1A33D7E77CE7E840C56F573FE0AD3C69BE16B384801DD4CAC7D4B1C8FD31F957414C1
851451408E17E21059B55F0F5BBAE54DDE4823823238FD2B53E21AE7C0FA8803A08CFE522D55D77C35AB6B1AA4176750B7885ACA5EDD4467819C8DDE
F5A1AEE91A86B5A51D3A4BB8A28E44026658CE5981CF1EDC0A4327D0A255F08D8A22850D648703A64A64FEA4D627C310927834A3FDC33481813DB02B
6AC74FD42CFC3D1E9DF6A81A68A3F2525D871B00C0E33D40FE5547C31E1FD4742B436125E433D9B6E270855D491DA8022F0B0B4B659ED3C3D14925B0
B86335C4CDF2EEEE140EBC62B37C09742CB4FD66F2E099776A210941F78B155C81F56ABDA1F87F5CD1A29B4E86F2DCD8CB21612E0F9A80F5C0E99F7E
69746F0ADEE953CD6C9751369AF769720924CA4A90429E318C8140CDED734EB4D4F489ED2FC85819725BFB847423E95C8F81EF665D41B49D60BFDA2D
13FD13CC18DC9EA07735BDAC69BAE6A13C6B15DDA4568920731E1B320072037E9F9553D73C33A96AFAA59EA0B7B05A4D680ED6895893CE79CD023AA9
98242EE7A2A93FA571DF0A811E137CF5FB53FF0025ADD9AD7589ECDA06BCB547652AD20889241F6CE2AA785BC3F77E1DB49AD96F63B985D9A45063DA
439C77C9E38A00C868D5FE2FE5C03B6CC15CFAE29FF102259756F0D2B8CA9BDC30F504AD4E3C3DAE9F130D75AFAC7CED9E598846DB76FA67AD3BC41E
1FD6B58D46C6E3ED967147652092350AD92D91D7F2A00B3F10768F04EA25803854C7FDF6B52FD9EDAF7C0D0C3A84C62B77B28FCD93382A36839A4F12
E91A86B9A08D385C410B4A419DB6920E0E401F8814D9743BEBCF0A4BA2DE5CC25BC858A39A304676FF00787E02802889647F064D0E8F184B08ED9D12
E2E0E4C8B839200FC793567C112C76DE17D1AD4EE67B889DD4F6183939FCC536CF44D4AEFC3ABA36AD2C514288103DB93B9C0E99ED8E9F5A97C35A2E
A3A6430417D3C0F1D9A3C76E63072431072DF950065DF0DBF17F4E391F3591E3F093FC2BB9AE2E5F0CEB9278946BE2FEC85D20D891EC628170463D7B
9AEBE0F37C94FB46CF3768DFB33B77639C67B508192D1454176B70D6B22D9BA24E47C8CE3201F7A62390D1B8F8ABAE0FFA7543FF00A05437C64B9F8B
5650C8E7CBB7B72E8A7A6769CD5EB1F0D6AB69E21B8D6DB53824B9B84D9223407690318E4118E1455ED6B479A4D5ACF59D3973776DF23A16DA2443D4
7D690CC6F8B11A9F0B45211F325D2E0FA64356B78E25953C157F2DBB157F2D391E85973FA567F8C343D5FC4BA6A47104B5F2A4C8859C1DFC7524715B
96515E5E584D6BADD9431A3208F08FB8382307E9DA8030FC21FDA6FE0CB186DA2B7D8F130123B9E32CDC9153E8DA0CB67E19BED2355BD8DA22589788
E0C6A79C9F4F5A9347D2F58D02D4E9F6AD0DEDAAE4C2D29D8C993C823BFAD5EBCD224BBD16F6D4CF8BABC5FDE4BD81F41EC3A5017326C1E187C1F776
BA3A4B35A416F285B9B96C090F39C0F4FCAA5F8793469E0BD3119B0D234AAA09CE48773FC8545A6E91AD3680BA15D886D2DD22685A78DF7BCABCF41D
B39E6A4F0B68FA9E9B05AD8DF242B6FA7BC8D14A8FB8CC58B76FE100350057B8C0F8B569CF26C5BFAD7695C84FA1EB3278B575F57B31E4C662484B31
DCB83DF1EF5D5C45CC4865003951B803C03DE98324A28A2811C36757D2AFAFEFB428E3D434D965669202D86471F7B1FE15A1A1CBA5EB1E17BB4B50F6
D04E6513248C33116EBF853B41D3B53D212F2478D6E3ED772F29815C2F97927904F041A924D119B47D55117C89B502CECAA01DBEDC75348A2B687FBA
F0EB69FA187B98225744BA9FE5424E4F03A900D55F86F72B0784AD84A5DBCEBA68D38CE0919FCB834FF0FE97AD8F0E8D1AF047696FE5B209D0FEF0A9
278C763CD58F0A68FAB68E9169F726D9ACADA491E391492EFBB38C8EDD4D0053BC427E2C58F3C7D918FE40D7550DF2C9A95CD96C20C0A8DBBB36EAC4
D5B46BE1E2BB5D7B4E58E668E231490BBECC8F507F1AD2B3B6BE896EAEA416DF6D9C8C2F3B540E0027A9A0473FE0BE7C5DE2C27AFDA507EB27FF005A
9754FED5B4D7EE752F0DAADE0F962BDB46E3E603823F0F4A9FC3FA0EB3A56B17F7D34B6522EA12879554B7CBC9271C7B9AB76961A8586B9AB6A623F3
21BA74DB02B0C90AA06EFAF5E28189E13D4EDB53FB64EB6E6D6F1A406E206EAA40C66B3EC811F1675120706C5727FEF9AD9B5B2B85D42F3555B748AE
2789512266EB8EEC4563DA691E2183C4F73ADBADA319E311F93BCF0074E7F0A045CF88BFF2246A1CE3FD5FFE86B55B5FBAB8B2F86B1CD6C4ABFD9215
271D01500FF3AB7E28D2F57D73476D3E25B389660A6466918952083C71ED56EDB4FB9B8D0DF4CD5E381A3F2961CC4C4EEC0C6E39031D01A00C6F0C4D
A8C9E0AB482D74E570F032AC8D3800E491923AD49A6E9D77A6781AF2C359915422481190EEC29E9FA9AB5A15A6AFA2587F677D996EE38B3E4CC250BC
7A1079A9750D26E6E7C3D7D6C195AF2F39247001C8E3F0A0672D26A7AAD868C34BF1469A5ED1E3D8B7911CED03A13F4AF41B168DACA030BF991F96BB
5BD463AD624B16A09A2B691736A6EE4783CA130FB8DC63E6F4AD4D12C9B4DD1AD2C9DB734310527DE84265FA28A2988E23E213137DE1C8B1957BF527
D7AA8FEB5D6CB76915F416A558BCCACC08E836E3AFE7589E2FD0EEB56FECFB9B174F3EC27132A3F01F9071FA55FB482EEE3515BEBF896031218E38D5
C3673F7893F8521F439BF10402E7E25688B2282A9116E7BE09FEB563C7C7175A09C907EDCA3F3228BAB0D7EE3C596DABFD860F26D90C6B11986E20F7
CD3FC53A6EB5ABCB62D6D6702ADA4C26F9A719620F4E940CC9F8B6D22DAE98ABBB6348FB8038071B700FF9ED5D06BF6BA96A1A13E9B65A7842E8AA1E
495404C63A60E7B54DAD68E7C49E1F16D7D10B6B9CEF4C36FF002DC671C8EA08FE74FB2B8D616C52DEE2C4FDAD10069B72EC6F7EBFA5022EE8CB3C7A
4DBC3772249711208E42A72370E2AFD57B2B65B4B64886091CB1F563D4D58A620AE6B58B5B3B6F11D9EB335E4892A218D6DD5771933E807B66BA5AE4
356D27518BC696BAEDAC46EA144F2DA10E015C8233CFD6931A2A6BB2DC4DE36F0F3C96E604DE76EE605D87B81D3B575B7F2452697787219562915B1D
B00E4573BAB699ACCFABD86B415585ABF364A41214F539E84D6C24175FD817493A6EB99D24631AF62C0E0500627C31B68A2F0B8B84FBF3CCE58FB024
01FE7D6ABC533587C5992D6DCFEEAFADF7CCBE8C149CFF00E3A3F3AB5E18B2D5BC33A6FD86E2D4DD444EF8CC38CA93D54FF8D5CD1744946BB79AF6A2
A16EAE46C48C107CB418C723BE00A00C6F0F4EF7DF12F587B8DD9B68CA441BF840602ACFC42965B37D22FAD8E268EE768F707B55DB9D2AE34EF141D6
EC2DCDC2CE9E5DC44A4061D3E6193CF4A4BDD3AEB5FD52C27BAB736F6166FE688E5C6F91FB703A5006EDE5AC7A869D25BCF18DB34654A9ED915C7782
7563A7E89A869F7836CBA5CA5157B904F03F3FE75DE571AFE1991BC7FF00DA6148B36412B8CF0CE3A7EBCD006E69B6BFD95A23B381E76D79A5207563
926B91F8737B77345A9DC8B369DE6B8DCD2F9807E1CD7A0BAABA14619561823D4572FA568F75E1CBFBB36307DA6C2E984854300F19F4C1EB40068163
AB5B788EFEEF5110C30DEE4A422404E47D3DABABACCB78279EFF00EDD770F93E5A6C8A3DC0919EA490715A74C1856178CAF64B1F0BDECD167794DA08
38C67BD6ED54D4EC61D4B4F9ACEE07EEE55DA4FA502313C170EFF03DAAC4FB249A263BFB863919A8174F9FC2FE0CD44ADC99EE36BB97EC09E38FA54D
E1FB5D4BC3F65FD9925A9BA891898654700609E873D2B67EC693DB4EB791A9372312A83918C6001486607C3654FF00844609460CB2C92348DD493B8F
53F4C567E993CB0FC55D46DA31FBA9A3DD201D010A083FE7D6B5345B6BAF0E5949A70B296EA1572D0C91639CF63E98F5A97C3DA3CF6FA8DFEB3A8A01
7B7A7EE2907620C617F41F9503323C556B1DF78F3438278C18C2963FED73573E23201E1A8D906D31DC47B48EDD7A554D4535ABAF15D96A8BA4B882CD
4A6D2DCB67BD58F1845ABEB3A6C76767A732FEF56472CC08C0E45202C78D7519B4EF053CB0B6D925548B7771B8735A1E1685478534D4655C3DAA1600
63395E734CBDD3DB5EF0E4B6377035BC8576A799CE180E1B8A874B9750B4D1A1D3BEC32ADE431F961C8FDD71FC5BA988CAF8796A6D2FB5F847DD4BA0
0739FEF77ABDE25CCDE24D12D9A4D91F985FDC9E9806B3FC256DAB681777F15FD84F38B99C30962191DF9FD6AE788B4BD6AF7C47A7DED941135BD91D
C374801634BA0FA9A9E3060BE15D44B74F2BFA8A6787EE16D7C1FA74F2AB10208F21464F381FD6A7D56CA6D5BC3D71672A2C534F1E0AEEC80739EBF8
567E9B6BA8FF0067D86973DAF9315B6D13485F2240BD00C7AE29F51742BFC4594A6896F12AA9125CA0C9ED83527C42B99ECFC1F31B72559CAC6C4765
3D693C5BA3EAFADCD6F1DA1B78EDADE4127CEC72E7E98AD5BED39F57D124B2D4923492418F909600F63DA8024F0F208F40B0518C0857A0C76AE7BC13
34A357D6ED0E3C98A6DCB8E809EC2B56D1753D374A8EC23B4134B1A948E5DE0291EA7D2A7D134B3A4E9F228025B9909924238DCDE99A00C1B78E197E
2A5CB145668ED06DC8FBA78C9FAF26A7F1B02350D00AE39BC00FAF6A659E95AD41E2FBAD6DEDA1293C42311F9A32391FE152F882C35BD52FAC658ADA
18E1B3944BB4CA0990FA5007594565EA905F6A1A1CD0DB916B752A1182D9C7B6456478060D42D74BB9B7D4653298A72A8739E9D79A02C75745145311
9DAE057D2E584CBE534B8446F46CE47F2AE65AF758B645B7F10DA992D5655C5CC439183DEBA1F10D84FA9699F65B67D8ED229DF9C6D03BD45711EA13
D87F67C96E8ECCA15E6661B08CF5C75A43317C7B78EC9A6DBC6BBA09A64666EC4678CD753708A74895368DBE4118EDF76AAEA7A3C5A8697159B361A1
DA637C742BFF00EAA8EEE3D4EE74DFB1470A42D22796F3170428E870280337E1F5C493685246E9B560B874424F51D7FF00AD5D6551D2B4E874BD3E3B
383255072C7AB13D4D5EA10985677883FE407799FF009E66B46B2B5EB4BEBFD3A5B4B330A79CBB59E46208FC8536088FC28D9F0DD903C958F04FEB59
DA19CF8D75AFF757F9D687872CB51D3AC63B2BCFB398E14C2B46C4927F1150269F7BA7EBB7B7D690A4E976A3ABE0AB0F5CF6A43362D6E7ED066F90AF
95298F93D718E7F5AB354748B59ED6C82DD48249DD8BC8474C9ABD4C4158BE2D95A2F0EDD15C65863FCFE55B5595E21D366D574A7B4B79123666072F
9C63F0A0096C9376871248C14341C93D811D6B95BED4355D32D9EC352B2FB4E9E46D1731764EDF5AE896CEF17C3CD66FB0DC797E5AF96DDBEA7F1A8A
D93514D39AC6EEDC4EFB0A89010148C607E2290CD2D3A4825D3EDDED8E61318DA7E9C7E756EA96956A6CB4D82D9B1BA35E71D3AE7FAD5DA620A28A28
00A28A2800A28A2800A28A2800A28A2800A28A2800A28A2800A28A2800A28A2800A28A2800A28A2800A28A2800A28A2800A28A2800A28A2800A28A28
00A28A2800A28A2800A28A2800A28A2800A28A2800A28A2800A28A2800A28A2800A28A2800A28A2800A28A2800A28ACCD4B5BB2D38B89CBB18D3CC90
20CF96BEA7D28034E8AAB0DF5B4DA7ADF24A3ECCC9BC39E98AA363E21B2BDD4458A2CB1CCD1F9881D701D7D45006C5158C3C436C7C4CDA118A41384D
C1F1F29F9777F2AD9A0028A2B1F5FF001059E80B6EF7692B09DF6AF96B9C74C93F98A00D8A2988EB246AE872AC0107D453E800A28A2800A28A8A59A2
876F9B22A6F608BB881B98F403DE8025A28A2800AC183C2D636F2DDF94D224176E1E4887720E7AFA7B56DB32A217760AAA3249380052AB2BA86460CA
C32083904500280140000007000A5A28A0028A28A0028A2AA5EEA1656210DEDD456E1F214C8C141FCE802DD155ED2EEDAF6013DA4F1CD11240746C8A
B1400573967E1486C9EF120BB916DAF25F3258B68C9F6DDE95D1D1400D550AA1540000C003B53A8A2800A28A2800A28A8D9D15D51994337DD04F27E9
40125150DC5CC16B1196E678E18C7F148C147E6699697F677BBBEC97704FB7AF97206C7E5401668A28A0028A28A0028A6AB2BAE51830F50734EA0028
A6B32A296760AAA32493800556B1BFB4D46133D8CE93C618A9643C6476FD6802DD1556F2FECEC541BCB98E10DD3736335345247346B244EAE8C32194
E41FC68024A28A2800A28A2800A28A2800A2A269E15CEE9A35C7072C0629CCCA885DD82AA8C924E001400FA2AAD8DFDA6A30996CA749A30C54953D0D
5AA0028A28A0028A28A0028A63BA4485E475551D4B1C0151C7756D2E7CAB889F1D76B834013D15134D1205DD2A286FBB96033F4A12786472B1CB1BB0
EA1581228025A28A2800A28A2800A28A2800A28A2800A28A2800A28A2800A28A2800A28A2800A28A2800A28A2800A29ACC154B310001924F6A1195D4
3230653D083906801D453249638F1E63AA6E381B8E334FA0028A29A485196200F53400EA28A2800A29AAE8C48565257A807A53A800A28A2800A28A28
00A28A2800A28A2800A28A62491C99D8EAD8EBB4E71400FA28A2800A28A2800A28A2800A28A2800A28A2800A29010C010410790452D0014514500145
14500145145001451450054D422B89AD1A2B4912391F82CC3381DF1EF469D671D859A5B44490BD58F563DCD5BA8A49E1876F9B2A47B8E1773019FA50
04B45145001451450014514500145145001451450014514D7758D4B3B0551D49381400EA28A2800A28A6960B8DC40C9C0C9A00751451400514514005
145140051451400514514005145140051451400514514005145140051451400514514005145140051451400514514005145140051451400514514005
14514005145140051451400514514005145140051451400514514005145140051451400514514015AF6EA3B2B669A4E83000CF524E00A8AFB5182C65
823B859313B8457032B93D8F35078861B7BAD2A4B4B9731ADC32C6AE07DC62783F9D7331DF6B7A24F0DAF89618EF74D322AC77683255B3F29348763A
0F106BA742896796CA59ADCE03488C06D39E8454DA9EACF61A72DFC7672DCC1B0C8FB18028B8CF4EF599F10959BC2B32AF79133CFBD5FBB257C1B312
70469E791D8F97401634DD521D4B494D4ADD2431BA960847CDC76C7AF15C9D9DEA5FF8A35BB7B8D36EE5B696140F1BA7CE31DB19FCBE95B5E00FF912
B4DFF75FB7FB6D54F4838F891AC2719FB2A1E3EA3FC68024F17CE9A6F84EE6D20B29D231000A6351B23E7A139A9FC3BF66BF8EC6FDF4C9A0BB4B5551
3BAE0118C601079FC454BE3938F076A7FF005C87FE842AD7863FE459D33FEBD63FFD0450067586BE973E287D2EE34936D7891963296563B7EA074AD8
D4F53B7D3628DA6CB3CADB228D7ABB7A0AE68647C586F43623FCFE955F5A99DFE28E8F039FDD2464A8ED9C13401D52EA2D1CF0C57B6E6069CED8C870
C09F43E86B33C49ADDB6972449A9696D3C323858E4215973C763D2AE6B3A3C5A8DE58DDCB772402C5CCB853807EA7B74AC8F88A33A5D8FCE01FB647C
7AF34023AA9678A081A799C2468BB999B8C0ACC6D7A05823BA68275B393A5C10368F72339C7BD61FC51BB9EDBC32A906424D3049187A60F1FE7D2B42
E2D752BCF0AFD820B7B44692D844A5A73803180785FA77A00D4D4F538B4EB5173243713444125A08F7ED1EA7DAA9D8789F4FD46C7ED564279F920C51
C7BA41EE40E94EB4B5BAB6F091B4BF9166B98ED9D1D94E41E0E3F4C564FC318228FC2892A280F2CAFBCFAE0E0500743A5EA969AAC0F2DA393E5B9491
186D6461D88ED5C7F88BC4024F18E93626CAE4476F3190E508690E38DA3B8ABFE101B3C51E29453F2FDA95BE84EEA87C4A07FC2C3F0D9F67FF003FAD
01D4EB6CEED6EE36648A68F69C112C654FEB566ABFDA13EDA6D70DE608FCCCF6C671F9D58A6238BF885E20161A6C96096B70D2CC5479BB088C0C83D7
B9E3A5745A76A56B7412085675608081240CBC7D48C573FF00131C8D06DE2C712DDC6335D489961FB34321FDE4A36A8F70B93FCA90C6DDEA16D6B224
52316964FB9120DCCDF8532D753B6B9B86B605E39D464C522ED6C7AD721E19B996FBE24EB924E4910218A356FE10180E3F2A3C7F732586BFA05E4276
C8256527D572383F9D0163BDA6B30552CC40006493DA9D553518A39B4DBA8A694C51BC2EAF203F7410727F0A622A1D7F4F41BE695A1898129348A551
F1E84F53ED49E2186DAF3C397AD2471CABF6591E32CA0E3E424119E95C8A6B0D6DA22E8DE2DB396385C08E2BB44DD1B2E32A73EA38E9F8E2BB0BD446
F0C5C24520910D930571FC4361C1A4339FF016A16767E0AB579261BB7BEE45F9981DC7A81EC335D6D9DE5B5F5BADC5A4C92C4DD194D72DF0C6D120F0
8A4A17E6B895DDC9EF83B47E1814EF016F57D762518B78F52956303A0E7903F4A019D0DE6A96765224534B999FEE44837337E028B1D4ED2FDDE38243
E6C60178DD4AB2E7A641AE4FC1F72FA978CB5EBCB839789BCA8D49CEC5048E3F2A778BE63A7F8CFC3B77012B24F27D9E5F468CB0183FF7D1A00E96EF
5CD32CAE5A0BBBC8A19171912640E467AF4A25F1068F0A8697508114F425B8354BC6BA47F6C7872EA2450678D7CD88E39DCBCE3F1E9581764F89BE1B
43E42C6F751A02467EE94EBFA01F98A02C75726BDA545179B25F44B19E8C7383C66AE59DEDB5FC1E75A4AB2C79C6E1EB5CBE9D796FE21F0CE916F246
A4CCEA255C7DDF2B963EDC81F9D75C88B1A85450AA3A00300500C47758E367738550493E82B80BAF16D95DF8DB4D50D3456B6A92EE778D86F665C0C0
1CE2BD0AB89493FE2EF48ADFF3E185FD0FF8D30475104F63AB5A968C25C41BB690E8719FA11EF5C5F8165B1D36FF00C41E63C7114BBD8A3B9504F005
77515CC735CCF02677C0543E471C8C8C570DE06B085FC5BE23BB9006920BA645CF6CBB73FA5203B6B2D42D2FD1CDA4EB2EC3B5C0C82A7DC1E4536FB5
2B7B128B3176964CEC8A342EEF8F402B9AB599A1F8A7796EB858E6B257603B918E7F9D4D7D6B2D878CE3D72EB528A2B2780C1E4B72EC78E1463D7072
39E3D2803734ED5ECB517923B794F9D17FAC89D4ABA7D41AABE22F1059E8366F2CECC65284C71AA93B8F6E7A019AE72F6E9A6F899A4BDB5BCB0B496E
C199C63CC4E4E71F9F5ADDF18BA4BE0ED4DA36575F24F230470680B199E0AD674C83C3368B35C849E4772EA54E77163EDE9815D9D733E1CB8B6B3F03
E997572A36470273B7241240E2BA6A6233F55BDB1B4B465D4A5F2E1995909DAC78C73D071D6B3FC31168763A3CB26872EEB3DECEEE4B1E4019EBED8A
D8BF0A6C2E03A865F29B208CE78AE5FE19ED93C20BC02AD3CB918E319E9480CBBCD574BD5BC5F7F6D35FC2B6874F30C73EEC047CE4E09EA79ED5D6F8
634F874CD0EDED6DAE8DCC4B9225CF0DCF6AE6618213F172E13CA429F62C91B463385AEEA38D224091A2A28E814600A10D8FAA573A95A5BCBE4CB301
281B8A2A96603D7001AB1733082DA59D8644685C8F5C0CD717F0D6E1F528354D52E5B7DCCD75B4B1ECA14103E9C9A623B1B6BA82EE1F36DA45913246
47623A83E86AB5F6B3A7D84A21BAB954948C88D54B363E801AC9F0FD86AB61AA6B57DA9988457126F89633C1033CFB718AC8F86D70FA95CEB1AA5D28
3732CCA0B7A0E78FA521D8ED6D2F6DEF63325B4AB2283838EA0FB8ED525C5C436D1196E2548D07F131C0AE1E7B8974FF008AF04511221BF8712AF627
69E47E207EB4B3DDBEA9F14A3D36E07FA2D84459623CABB15CEE23FE043F2A00ADE33D2F47FECBB4BDD304609BB44631B93B813DF9AED353BAD3ECAC
4C57F208ADE4431F209C8C608E3DAB90F895651431699790AEC3F6A48DD578561D4647AF15DC5F056B19C3A865F2DB208CF6A00C7F0A41A1DA6972BE
8537996A6562F2313F780E7A81D055EB4D6F4BBC95A1B6BE8249114B3206E540F51DAB03E17C6A9E142CBFC773231FD07F4ACBD4AC61BDF8B5042C9B
6316E24942F1BC804F3EBDA803B8B5D56C2EE7305B5D4724A067603CE298FAD6989722D9AFA01296DBB7777F4AE4FE23A795268B25B318656B9F2B7C
7F29DA40E38ABFE34D36D2DFC1376B040AA6DD55E3603E60C18739F5E4D0075B4566F87A77B9F0EE9B3C872F25AC6CC7D4ED19AD2A62229E086E6168
678D648DBAAB0C83DEBCE3E246956BA64DA66AB69671222CDB26540155CF05723F06AF4CAE37E22AFDB3C35736F6F09B896122672A78840EE7DF19E3
EA6931A2F6A96563AB45FD9567045B5947993A463112FA03EA6B534BD22C7498445636EB1E14296FE26FA9ACCF01CD0CDE0FD3DE18D631B0AB81DD94
9049FAE335D1D001451543546D512343A5476D23E7E713B11C7B62988BF4562A49E24FE3B7D37F091E9525F1173BEDB4F1F491B9A571D8D9A2B32293
5969904D059AC7BBE721989DBEDEF5A74C44370B3340EB6EEA9291F2332E403F4AE0BC6977E22D3A5D3A23A8A08AEA5084C09B1836477FC6BD0EB86F
89470FA11FFA7E1FD2931A3B751B540C93818C9EF4EA28A620A28ACFD762BB9B45BB8B4EFF008FA78CAC7CE393EF400EFED6D3FCE307DB6012019C17
1FCFA55D072322B8F11DADBF835F4DD563B66BE8AD1F7C518123060A7E6F63DEAC7C3DBB9AE7C2164F74E4BEE7442C796009C7D7A1FCA95C0E823BDB
49A411C57303B9E8AB2027F2A72DD5BBCC6159E26947540E0B0FC2B80D7F488A5F88DA65B5AA9B55962323987E5DC0139E9F4ABFE37B68B4E5D26EEC
6358268EED10320C12A4F43EB40EC76F597A858DB6A2E47DA2449D1481E54A411F503DEB3BC75ACCBA1F87249EDF227958451B0FE1241E7F214EB0D1
2DE4F0E5B8542974F6EAC675387DE541CE7EBDA8028FC387B87D2AF05CCEF332DD328676C9E2BB1AE27E16A94D12EFCC24CDF69612739E6BB6A103DC
CCD6E0B5BFD3A6B0B9BA5B7F393AEF00819EB8CF4A341B0834AD1A1B5867F3A38C13E693F7BDEA8F8D6D2DA7F0BEA32CD0A3491DBB32391CA9038E69
3C2304573E09D3A19D03C6F00DCA7BF39A3A874306F6D8789B52D6E117AAB3DA158ED72FF228EA5BEB5D9E916D2DA6956D6D3CDE7491C614C99CEEF7
AE1B44D1AC25F883ADDAC96CA61891591324019C7F8D7A1471AC512C718DA880281E80508192554BEB1B6D422115DC2245072012460FD455BA2988F3
4D7F4E8748F1A6951ACF72B6176DF3C7E6B6D539C71CFD2BA3D6744875658ED2C1DA0549079B3A487E503F8473C9AC2F8A13092DED5ED519E4B29B74
9228E133D89F5AEDB476864D2ADA5B60BE5C9187C8EE48E7F5A431BA4E916BA4C2D1DB073BF05DDDCB1622B468A2988649247126F91D517D58E05245
2C7326F8A4575F553915C759DC8F1078E6EE29D835A69A988E13CAB3670588AD3897518BC5E6382DBCAD2C41F3B018567ED81EB48763A0775452CEC1
547524E0536396394131488F8EBB581AE305DB6B5F10E5D3A566367A7C45BCACF0CF91C9F5EB52F8BEECF87AEB4CD4ACD762493F953C69C075238E3F
0A2E163AD92E2189B6C93468719C3301C50D3449B774A8BBFEEE580DDF4AE5BE2259DB49E19B9BBF2F13A05DB22F07AFAD49A6E816CDE158C6A05EF6
57B7DE6491CE4657202FA62988EAAA2FB443E7791E747E6E33B370DDF9573FE01BA96EFC2F0B4ECCCD1BBC60B75C03C56669BA8E9B0DD5EE95AD2FD9
AF5A6722E24E3CCC9E0AB76382295C763B294C6607DEE15082A5B7631DBAD62F86745B5D28DD3DADF1BB33365BE6042FB704D4FA2E9CA9A0C3677789
D7E63F31CF049C73FE7AD61FC37082DF55D99E2F1873E83A500760D344922C6F2A2BB7DD52C013F414A24421886521490C73D31D735C97F6CE9B69AF
DED86B11B413CB2131DC4832A548C2807F878AD2F0EE95141A55D5BBBFDA21B99E470ECDBB7A1E01FD28B858669BAC2DFF008A750862B806D6D6258F
1B86D2F9E48FE55D06E5DBBB236E339CF18AF3EF0A68BA6CBE27F1044F6AAD1DBCC123058F039FF0AE9F5E8DADBC3B2269E9F2458DC88792A0FCC050
819AE934522EE8E5475CE32AC08CFA54B5C4FDBB48D7ED2D574B9D6DEE52540D6E7E56600E0823BD76D40053249638F1E63AAE7A6E38CD24F2086092
53D114B1FC0571FE0D71E2137BACEA0AB29698C70C6D9658D47A034C47660820104107A11585E2272F3E9D64D2F970DD4DB64C360B00338FA55186FC
E99E351A32716B731798887F85BDAA878F2C613A8E8CCA1835C5D847C31E41233F4EB498D2D4DAF0DE937FA64F7E6F2757866937428A721064F4CF4E
B5D0D54B0B0B6B08D92D50A86396258924FE356E9A131A4818C9033C7341650C149009E833D6B81F883A7A4779A5CF67FB9B89EE76B38E99E0E715A7
E22D322D3BC393DE465CDE5BAABF9E58E59B233DFA7B521D8EB6A192E2084E259A38C9ECCC056647AA3C3E134D56E14348B6C2461EA71593E0AB68B5
5D11756D4105C5C5DBB96320CED01880050075AA4328652083C823BD33CE88C823F353792405DC3271D78AE5746D45AD3C6179A2024DB905E119CEC2
0648ACEF13E9A9FF00098E95158B1B669C1F31A3247193D3D3F0A2E163BB59E1694C4B346641D5430C8FC2A5AE3BC65690697A2477BA7A7917104ABB
5D3A9CFA9EF5D49B88E0B317170E150202CC68B8162B0B58D334ED5AF6D5EE2EC2C903F11AC83E639E98F5AB76DAEE957770B6F6D7B1C92BF4519E6B
9FF16C11FF006FE8C523019E6F98818CD0C11D6CB34508CCB2A463AE5980A68B980C22613C4626E8FBC6D3F8D53F10416D368F72D750095638D980C6
4E7158DE07D3D1B408AE2E984EF283B430C88D7D07A500754AC194329041190477A635C40B9DD34631D72C38AE7FC27232DCEAF665CECB7BA2B1A139
DA2ABF8CF488469FF6FB684092170CEABC065EF45C2C7522689A3322CA8507560C31F9D2C72248BBA375619C654E6B06786C26D0205B7817FD242840
9D727AFF005A4D7B6F87FC317326971889D570081CFD6988DB3756E25F28CF17999C6DDE339FA549248918DD232A02719638AE4347BAD2078683C8DE
6CF2C45E42AACCE5BD8E3AD69785EE24D4F4245D46D8F990B98C89573BB1D0F3EC690EC6C9BAB7085CCF105071BB78C6692E2EA0B6B733CB2011819C
E7AFD3D6B8BF0E582DF5C6B718C796AE5212DD149F6FC2B72D7455D3BC3CD0DCB0B99A1491D5C8E848EDF90A02C4BE1CD562D474F331B805DE560159
86E1CF4C54BAFD95A6A366B6D7779F66C3070C1C29FD6B2FC17696A3428EEE58235983B167618C60F534EF1E411B68826D83CD49502B0EBD68E81D4D
D96CD24B15B5124A8AAA1432390DC0C7515CDF82235B38B52134E4E2E36EE91BAE3EBF5AEA6D73F65873D762E7F2AE37C2DA6C17BA96A53DC132451C
FF00246DD3711C9F7A00ED51D5D43230653D083906B93F1A59297B2B847911E49D55F0C707F0A7DF5C7F60788AD12DFE5B4BCE248FA853EA2A6F19E4
C7A7283D6E97F98A00E8E34F2E244CE7680334D59E1690C6B2C65C7550C323F0AC7F176AAFA4E8D24909DB2B8215BFBA3B9FD6B2EEA7D1A3F0FB3DAD
C24976918606239919FE9F5A623B1A2B2FC3B7D2EA3A2DBDCCE8C9230C303D78AD4A0028A28A0028A28A0028A28A0028A28A0028A28A0028A28A0028
A28A0028A28A0028A28A0028A28A0028A28A0028A28A0028A28A0028A28A0028A28A0028A28A0028A28A0028A28A0028A28A0028A28A0028A28A0028
A28A0028A28A0028A28A0028A28A00A1ABE9C353B316E6668B122BEE5193C556BDD3EEEFB4F7B0B9681E272019493BB6839FBB8C67F1AD8A2803335D
D223D6346974E2E630EA36B75DA474ACE7D2F55BCD0A4D26EE6860511794268B24BE3A7D07AD749450073BE1BD3B56D274B1A6DC3DAC89029F26419F
98939C11C71D6A9DBF87759B5F115CEB31DFDAB4B72A15E368DB6E001FE15D751401CF6B9A66AFABE952D81B8B38526501D95589EB9E29DA1D96B7A7
59C565732D8CB0C316C8DD43EEC81C6474C74ADFA2803918FC3DAD2F884EB8D7D66D394F2FCBF2DB6EDFE79AD1D6F436D49ED2F629560D46D0868E41
9DA7D41F6EB5BB45006544350B95F2AFE0821888C3849092C3D07A0F5F63597E23D0756D72410FDB2DA1B58DC491E1096C8E99AEA68A2C1731EF7497
D6344934FD60C4CCFD1E104608E8DCF7CD4760759B3B54B39ADA19DA2508970B26D523A0257AF1DEB728A00C8BC8AF6DB4DFB3D85BADD4B2061233B8
4E4F53EFD6B33C2163ABE89A22D8DD5946EC8C594A4C3BFAD755452B0EE725A2D8EADA6EBFA85DCB651B41A94CAC76CA331633927D473567C49A4DED
CEABA5EABA7AACB2D8B9DD0960BE629C6704F4E95D2514582E64E9B6D72D7D3EA1790882591151620E1B681EA475AD6A28A6239EF1A68D3EB3A308AD
08F3E09566453D188ED56AD56EEF6E2D6EAEEDC5B88109D85B259C8C1FA015AF45007327469F4CF13CFACE9F089D2ED36DC43B8290783B867AF4E949
79A55C6BBACD85D5F5B1B7B5B12CEB1BB066918FA81D0715D3D14005677881669341BE8ADA1334B242C8A83BE463FAD68D1401CF41048BE1F8B4CD52
C9A7912DD50ED4DE8C40C0E7B11486DAF34CF0AC5A7416CF773790D11DAC06D241C75FAFE95D151401CAF839353D33C3D1D8DEE99379B0676ED74C30
249EED50782ED756D326BFB7D4F4E7892EEE5EE1255911941239070723A57634503B9CB7F67CBA378AEE355B6B6967B6D4102CCB18CB46E0FDEC7706
966D2E4D6FC4F6BA8DCDBBC56960B98965182EFD738ED835D4514082B92F08E873695AAEB2658DD6DDE73F6653F7761249C0FCBF2AEB68A00E37C1BA
04DA56B3AB4D3A388FCD2B6FB871B492495FD057654514030AE42FB4EBCB4F8810EB51DB493DB4B6E617F2864A1F535D7D1401936A66B786E6F27B79
19EE250DE54601655C05191EBC735CDF8623D4B49D675BBBBCD26F7C9D42E3CE88A056206E738233C7DE15DD5140EE70D149A8A78CA7D665D06F4DB4
90AC11B00A645FF68AE69DA9C57D67F1062D527B2B8BAB0106C8CC6BBFCB38EC3B1CFF003AEDE8A02E705AB45AA45E2ED37C48D632496CB1988C318C
C9102081BBEA5F3ED5BB71A6CF79E0FB8B02025C4F13E149C60B12C07EB5D051405CE2749FB6CBA058E82D613C53C2516E1A54C20456CF07BF6AEDA8
A281199ACDD5CC16C52D6C65BB79559708400BC77CFD7F4AE7BC149A9E89A43E9D7FA55C06466915E3C306CF6AED28A00F3D8BFB5C78EDF5C3A0DE0B
7783CA2B91B870067F4E95DCD9DC35CDAA4CD0490336731C830CBCE39FE75668A561B647246B2C6D1C8328E0AB0F506B8CD02CAEFC1D757B6AF65717
361712F990CB6C9BD97B61876ED5DBD14C46458B5C5EEA4F7D3432DBC09118A28E55C336482CC476FBA0562685A7FF00C2257FA8C52452BD8DD38921
9517704EBF2903EBFA57654500727A7E9526A1E317F114F0B450C7108ED92418627182D8EDC1AAFA9E9D2E95E3CB7F1004792D2E10C5705549311DB8
0703B703F5AED28A00F3EF1D4977ADD95A9D26CA7B8B6B5B90CEEB1B658E3F84752073935D15C6B62F34EB85B0B0BD9E62A63F2CC5B0AB11DF76315B
F4500725E06B5D4F48D15AC752D3E488C6EF2075756073CE000735450DE8F888FAB36957A6CCDB792AFE51C8380738FCFF003AEEE8A00E1BC7CD7373
7BA4C56B63713A5BCC2E24288718F4FAD6A78BA49350F0B4F6D6967772CB7710D8AB11F97907E6F4E95D2D140189E1BBC8DAC6DF4F4B3BCB736D6E8A
DE744540C003193D4D6CB36D524E70067819A751401C8EB1E2D3F67F2748B3BB9A773B4C9E490231DCF3D6AF47369B068AF6E04D2ACC8C2402362CC5
8739E2BA0A290CF3CF00DDDD689A65CD96A56378882532C24444E41038FD33F8D77969702EADA39D63923120C859176B0FA8A9E8A626154752B07BE4
8D63BDB9B528D9CC0FB777B1ABD450063C5A2DC47BB3AC5FB64E465C71FA527F62DC73FF00139BEE7FDA15B3452B0EE64268D38642DAB5F3052091BF
1BBEB5AF4514C456BCBB8ACAD24BA9F779718CB6D193D71D2B80F1CEA69AA3698BA7DB5D4DF679C4EE442C00518F5AF48A28033B4ED62DB52256DD66
0C1771DF195C7E75CBEB1E27D66C7C67069E96C7EC0F2C7164C67E72D8E8DF9D7735CFDC2CBAA78922824B765B3D38F9C5987124B8F971F4CE69303A
0ACBF11BDE47A0DE3E9C18DD08CF97B7AE7DBF0CD6A54371047730B4528255BA80C47EA2981C4F87AE62FF008439AD6C6DDE4D5A689D255742199CE7
2598F6A77816571A7699A63C72A4F6524EF3AB290064B01FFA1D77200000030076A5A5619C15FEA70C9F10B4EBE4498DA416EF13CBE5B6D04EEF6F7A
5F1E6A56F756B6096AB2CE63BA495CC71B1DAA3F0AEF28A02E723E33B3FF00849BC2064D34B4AD1B89915472D8C82B8F5C1ABBA6EB569FF08DDB3A4B
BEE16D957C95FF0059BC2F231D7AD7435025A5B24E674B789656EAE10063F8D311C2FC33D420874F9ECAEDBCABB7B92C11D482D9F4FC457A0D43F66B
7DE1FC98B783B836C19CFAD4D401CD78D6F00D0EF34F8A19E6B9B8870891C65B209C76FA1A4F0D6A16D63E15B34BA1340D6D0AA481E16041C7D2BA6A
2803CE344D6ACE0F1B6B1A8DD79D05B5CA848A49226C3118F41ED5E8714A934492C6728EA194FA83525140114F3C56D0B4D33848D7AB1ED5CC6B1E2C
85678ACB4F5998CAC164B9119D9129EA7A75AEAC80C0820107820D342204D81542F4DB8E280395D767D23FE113BCB1B770E1A23B14292CCDEBD3AE6A
97843C476DA7F86EDAD75512DBCD0E502B44791DABB90000000001D00A0807A8A43B8C8DD658D64439570181C7506A4A28A62387B4B37F0EF8EAF2F2
E13FD0B5153B651D10E41C1AEA21D463BBBA48ECF12C6326590745E381EE73579D16452AEA194F504641A48E348D76C68AABE8A302803905D3FF00B1
FC7F2EA2548B4D421219C9E11F209C9EDD3F5A4F15DBFF00C245A8699A5DA36F8239BCEB9953954503819E99EBFA5762EAAEA55D4329EA08C8A44444
1845551D7006290EE727F116FADE3F0E5CD90901B8936ED8D464E33572DB59D397C3119370BF2DA842BDF3B718C56FB4514872D1A31F52A0D34DADB9
183044475FB82985CE4BE1D5E5BC7E1AF224709324923B46DC36DCE738FA549A51B3D6B4316DE20854CA338927E19D7390C18FF9E2BAA58628D8B471
22B1182554038A6C96B6F322A4B046EA9F74328207D2905CE6F42BA8741D0DA2BE9E56812775B576524BA76C63F1AA5F0FE74816F2D668A68669EE1E
48D5E32372F5CE715DB2A8550AA00006001DA9D4582E729A5DD5B6B2BA85A6B16EAD11B87F2CCE3865CF404FA7B54BE13B23A7DC6A56D6D233E9C928
FB3E4E429C7CC07B56FBDB40E851E08D94B6E2A50104FAFD6A48E348902468AAA3A2A8C01405CE33C3B247A5F8A7C431DF3F96D3CE248CB0C6F1C9E3
F3157B5496E74AF0F4535AABB4AD721CA2F258124E3F2ED5D218A3670ED1A975E8C5464539955B1B941C1C8C8E868B05CE37C4DA669779A70D4EC952
1D490ABC4633B58B71C102BAF837FD9E2F37EFEC1BBEB8E693ECB6FE6F9BF678BCCCE776C19CFD6A6A2C1722B8884D6F242DF76442A7F118AE57C136
A742B3BBD36F8889D276752E70194F707F0AEBEA39218A5004B1A3E3A6E50714C47276964FAB78E1F58087EC7691F97139FE36C75154BC67AAD94FAB
68D0C33076B7BB0F230E8A011DFF000AEED40501540007000ED48638CE728BCF5E3AD2B0EE43697B6D7A8CD6B32CA10ED6DBD8D59A64714712958D15
0139C28C53E988E17C777B6926A1A440250CD0DD664DBCEC1C75AD4F195FDB3F84EF5619925695004546049E476ADF6B4B6662CD6F112DD4941CD2FD
9ADC631045C7FB0290EE6368E96FA9F84E1B2F3158FD9D5245CF2A71DEA9783A78B49D24E957E7ECF35ABBFF00ACE032924E41EF5D3C70C511262891
33D76A819A24862971E6C48F8E9B941C5160396D034F6BAF135FF881D196293F776DB860B0EED8F4AABACEA56D2F8D7489222ED1C2089240A768CE6B
B8A8BC8840C0863C7FBA28B05CE47C79A8413E842DED4B4D23CAA4AA29381CF26BA5B3B8B4BDB6D91B24CAAA03A91903EB9FA55AF262C93E5264F53B
4734A88918C22AA8CE70062802BC5A7D9412F9B0DAC4927F7954022B95F115D9B9D774B960B1BC962B49899DD213F28E95DAD14582E62EB7A85BAE85
2BB09099E26F2D3690CC71FA75AA5E13BFB7B6F0EDBC772FE4C8A3E6461C8AE9E8A00E47C2F7B00D53562EB2446EAEF747B90FCFF435D3DD409736B2
4120052452A7233D6A7A284071FE0BB3BD492E4DE13E4DB3B436E0F1DFA8ADCF10E9CDAA68D716898F3197299E9915A94530B9CDF87F5386D3458ADA
F94C1736EBB5E32BC923D2B5A2BBC5A35D4E0AC6C4951B4E42F6CFF9EF570AAB3062A0B2F4247229D408E3FC23791C5757E92452C2D737398F721C1E
B5D36A31493E9F71144C55DA32148EB9AB54500723A05CC69E1F1A770D785993CAEFD7AFF5A3C7B73FF129105BAB493ACCAE405FBA00CFF857581103
970AA18F538E4D3A90EE53D3AF20BEB4496D89283E53918C1C5739E16B88EC6E2EE1B81E52CEFE623B746E715D7D46D144E811A34641D14A8205311C
BEA111D6BC4D686DD19EDACCFEF255E99F635078DB518566B08A360CF0DCAB48472106457628AA8A15142A8E800C0A1955C61D430F4233480E73C496
51F893C3E5AC5CC8F19DC9B78C9EE2AF69FA8D90D32294BC6922A00E9801B77718FAD6B0000C0E0533CA8FCCF33CB4DFFDEC73F9D3022B1F3CDBEEB8
50ACCC485031B47615668A2800A28A2800A28A2800A28A2800A28A2800A28A2800A28A2800A28A2800A28A2800A28A2800A28A2800A28A2800A28A28
00A28A2800A28A2800A28A2800A28A2800A28A2800A28A2800A28A2800A28A2800A28A2800A28A2800A28A2800A28A2800A28A2800A29A48504B1000
E493DAA99D52C43A83700076D8ADB4ED27D37631FAD004B797D6963187BCB88E153D0B9C6696CAF6DAFE0F3ACE649A3CE372F4CD41ABDA4179A6DC24
F0C72E236DBBC03838FD2B07E1A127C21167B4AE2803ADA2B3EEB56B0B494C33DC8120192AAA5B03DF00E3F1AB2B73035A7DA52556876EFDEBC8C7AF
14013D159967AFE937D23C76B7F0BBA0CB0CE303F1A92D357B0BDB86B7B5B9579917714C1071EBC8A00BF4565278834992F8D9ADF45E70EC4E013EC7
A1AD5A0028ACE7D6B4D8E73035E441C100F3C0CFBF4AB93CF0DB4465B8952341D598E05004B45538754B09EDA4B886EE17863FBEE1B85FAD4B6D756F
769BEDA649541C128D9C1F7A009E8AE53C6DE20B5D3ACD6C85D04BA9DD46D53CAAE4649F4AE86D6F6CEE7F776D770CECAB921240C71EA714016A8A84
DC4025F28CD18938F90B8CFE54F5756CED60769C1C1CE0FA5003E8A28A0028A8AE2E21B688CB71224718EACE702B9FF0C6B50EA53EAD37DAFCC812EB
6C5BDF80B8EDED401D2D15CF78CECA39BC3B7B73B9E39EDA0778E4462A460671C7514783AE97FE111D364BA9D43C887976C16F98FAD007434514C791
231991D546719638A007D145303A162A194B0EA01E6801F45315D1890ACA4AF500F4A7D0014514C74591191C65581047B5003E8AE1FC65A53691E1F9
351D1AE6EEDA6B7652DB27621813839CFD6AF412DA4DE16B2B9DD3C93DC440443CF72CD21EB8E7D690CEAA8AC0D03C3EDA72C73DDDF5D5CDC8C93BE5
25467B63DAB7E9886921412C4003924F6A1583286520823208EF59FAD5AC5A8E9B3E9CF742DDAE130181F980CF50334CF0EE98348D1E3B35BA374132
7CC3DFF53401AB453048864281D4B8192B9E47E14FA0028A29ACCAABB988007726801D45355D5C6518301C70734EA0028A2981D19CA0652CBD403C8A
007D1451400514521206327AF4A005A28A2800A28A8E646922744728CCA4061D54FAD00494573BE15D1B51D1A3BC4D4350379E6CBBE3E49DA39CF5EE
6BA0DC376DC8DD8CE3BD003A8A28A0028A28A0028A6B0CA9009048EA3B570DF68D5B4FF88761A449AACD736B710B4BB6400767E3DFEED00777451450
0145610D37551E2B37FF00DA1FF12D31EDFB364F5C63A631D79ADDA0028AE4BC5F79ADE89A5DC6A3697F132070046D6E32809EC73CFE35D1E9F23CDA
75B4B21CBBC48CC718C920134016A8A28A0028A29AC32A4024123A8ED400EA2B83BBBDD6F4BF1B699A6CDAA9BAB6BB3B8A989548192306BBCA0028A2
B9B1A7F883FE12E178D7C9FD96010210C738C1ED8C75C50074945145001451450014514500145717E3DF13C9A4ADB58584DB2F27914B30C1D899FE66
BB4A0028A28A0028A6B862842B6D623838CE0D711E30D47C41E1AD3A2BD8B538AE55A508C1ED957A827B7D2803B9A2B13425D5EE2DADEF6FEFE2759A
30E608E0000C8C8F9B39EF5B740051451400514514005145140051451400515CCF89F56D6F45D32E350861B29608DB8C96DC149C027B7715B3A45CCB
77A459DD4FB7CC9A1576DA3032466802ED145140051451400514551D62FD74CD2AE6F5B6FEE50B00C7009ED4017A8AAD63349716304D32049248D599
41C8048AB34005145140051451400514C9376C6F2F1BF0719E99F7AE4F4CD6359FF84D24D175130B42B11903226370C71401D7D14514005145140051
4514005145140051456643A9F9FAF5D69AAAB8B7895D9B3CE5BA50069D1557503722CE4366E8B328C8DE320D65783F54BCD5B4A92E6FB66F599906D1
8E05006FD145626BC35F32DA8D11A254DFFBF2F8CE3F1FE9401B74514500145145001451450014514500145145001451450014515CC4FAA6AF6FE2AB
2D3E5FB3FD9EE0337CAA79007E7401D3D1451400514535B3B4ED201C7048CD003A8AE72DF53D4D7C59FD9B762030343BC18C1EBDBAFD0D2F8A752D4B
4AB4FB55AF93E5EF0BCA9279F5FD681D8E8A8A86D9CCB6D148C006740C71D32454D40828A28A0028A28A0028A28A0028A28A0028A28A0028A28A0028
A28A0028A28A0028A28A0028A28A0028A28A0028A28A0028A28A0028A28A0028A28A0028A28A0028A28A0028A28A0028A28A0028A28A0028A28A0028
A28A0028A28A0028A28A0028A28A00E4BE235F5D59786C7D9437EFA558DD97B03FE34ED456E2EFC222C2CF4CB913490284C85014F1CE735B1AF6951E
B5A4CD612B141260861FC241041FD28B03776BA75ADBDD44F2DCAA6C2CB8DA481D49ED9A4323B3FB58F0D01A8C612E96DD848A083C80476F6C561F80
66369E023758DC54CB2600EB8FFF00556FEA725CC5A4BA081EE6795193118180483FA7F8562F81A0BBB0D0D749D474E9E260CFF3919420E4F5EDE940
7429F807512FA65C5E4D6D7535C5CCECEF2226ECFA0CD5FF000826A506A1AAC7776325B5A4B319ADC30E064F23F97E551E85697BE18B9BBB06B59EE7
4E924F32DE4853714CF50456DA34E6E24D4248E611A46523800F9986412D8F5E3A50072763656F2FC52BF8962558A3815D900C06385EDF539AB3E368
5DBC41A10B693C99AE2530BBA1C31438FF00EBD54B193508BC7975ACBE8BA82DADCC2B19FDC9CA9014671FF01A9FC4F25EEA1ACE972DA68F7E52CE6D
EF21888C8E3A521B2FF8FEDA187C153986344FB334462C0FB9F3A8E3F3A975ED5E6B2F020BF524CF2DBC60367A3381CFEB4CF183DD6A1E1B92D2CAC2
EA49AE3690A62C6D01B3CFA1E288AD63F11F845B479E19ADA68618E37F3131B5D40C63D4645311474E9AD57C15F626B39E49E6818B01116DF21079CF
D71CD6CF850DEDCF87228B59B5649A3263649467701D0E3F4FC2ABF87EF65D3B4B834CD52DA74BAB75D83646584A01E0A915AA2EDACAC64BBBEF3006
7DDB150B1407000C0A00E47C01A7DBB5F6B6EC84C515DB2471372A393DBB9C62B434F3F65F8917D6B08090CD6AB2145E9BB8E7EBD6ABF826E5ADAE35
73756D750A5CDCB4F133C2DF329CFA0A517691F8F5B5265996CDED8462431B609FA633403DC3E22DBDA14D29E5823323DF2076DBF33260E467BF6AEB
21B5B1B5940860B78A47040D8814B0EFF5AE7FC776B7325BE9D776D09996C6ED269114649008AD182E5753D4EDAE2D958C1046C77B295F99B036F3DF
1407431E4D4A0D1357BF875DB62B6D773EF86EB6EE5C1EC7D315ABE12B25B1D32658AE04F0CB72F24526EDD9438C73F856668D7114F16A165ADC2BF6
696EE4F29E61F2C9CF2326AFF84AC469F6B790C193646E59ADB9CFCBEDF88A1033A1A2B8DF09F8A750D5F5EBDD3EFECC402356923E082A0301839EBD
7AD7654C443756B6F790986EA149633D55C645711E01D334CB8D36FE492D617885E3842EBCAA8E8335DED701E0994E97A46A5A5DDC7225DADC39488A
F2FB86063D7A521A3A6F126D3E13D47C90190DA3EDDBD31B7B7E1587E07D2ADEEBC176A6ED04ED2AB85DFCEC1B9B007A7AE7DEAE6BEF1E97E0A7B29E
5C4EF6C51140C966EE07E74DF08EAFA75B784AC9259D637861C3A370D91D78A07D06FC39BC96E3C2C5AE6567F227910331CE1460FF005AA5E15D5AD7
59D4F55BDD4AE22E25F2ADE29986D58C7A034DF8653C3FD8B75A7CCC5256B976556182CA54723F234EF095B7FC22FA9EA1A5DF8F2E39A4F36DA76FBA
CBD319EC68132DF85F532DE20D4F4A851DAC9099ADA42A718C8DCB93EE78ACD5D3C597C538A1B17F2626B5F35D739CF27207E55D8ADD7DA2E77433AF
D9A253E61C70CC7A60FB004F07B8AE346A56971F15619E29D5A05B531990F0A1B078CD005BF1322E89AFE877BA78311B9B9F267504E1D4E3A8F5E4D7
6F5C27C42BCB64D43438DA650F0DD89641FDC5E393E95DA5BDC41751096DA68E68CF46460C3F4A103D89E8A2B0F56F14E97A5C539927F3268720C283
2C5BD298893C5051FC3F796CCA6492E6178E28D4659D8838C0AE73E15BA4DA14A9292D3DB4C5006E7CB523200F4E735B1A05EDA5EC0355B9B8885CCC
B9DACC0792BD9466B91F09EB567A478C35C82E2658AC6E25778E423E5F95CE39F4C134AE3B1EA34557B3BDB6BE84CD67324D1862BB90E464558A6239
EF1AE9F6979E1BBE9AE2DD1E5860668DC8F994819E0D1E0ADC7C13A779580C603827A6727FAD278C352B7B7D22F2C1848F71736CE238D109CE411D45
51F08DDECF0525A112C3756D6D2165642A57A9041FC452199B617BA35FE9D369F7732D8EB512C80CB2FC8C1F9E41EE3DBD2BAEF0DDBC96BE1FB18655
6590440B0639209E4FF3AE6EDAC747F13787607D62258AF522C3C8DF249903AFB8ADCF068B85F0B592DDBBBC8030DCFD4AEE3B7F4C5086CDCAE56F57
53B9F1B4104966D368EB11DCCC3E40D8CE7DCE463F135D55711ADEA534BE37B7D1AF2630E98C81C85F9448D8CE19BD28624417D729A67C43D320D2A4
0B05E26DB88633F21393CE3A67FC2BBEAF39D5A7B58FC7BA25E08C41A6C08C827DB88D9B0C78FC48AEC7C3EFE769A6E173B679A4914B752A5CE0FE54
03356BCEBE23DB8D3EE74EBED3C347733DC10FB58E1CF046457A2D79E78FF57B3FEDAD0E38E7566B3BD124E07F000CBD7F5A182DCD8F13E9CB65E1CB
8BF86465D42D809FED1939660467F03E94ED2F4F4F116850DFEAA375C5D43952A4A88BD368FD697C6D7F672783AF847731319A20230AC096E474A9FC
3BAAD847E14B399AE63548A001F2DC8207231407421F006A571A9F878B5E486496DE7780C84E4B818209F7E6AB6816D717F15F5CF89EC5A2944C4C6F
2B95010F4039C0C7AFBD45F0DA7B54D0E68C4C8257BB91F616C1C718E3E82ABF84F544D5FF00B4753D6E657961976476E7EEA2E38C2F724E47E14017
FC097F35CB6AB68F3B4F059DCEC8246393B4E78CD2F8BBC3D7179653DD69571711DF03BF025389001CAF5E3DBE958DE1EBE68EEFC430B11697B3DE21
8605E1C8CFF08EFC575BAC6BF65A5C7B0CF13DC3B0448B773B89EF8E8050072FA44D6BE28D0ED2CECAE2E62BE8B69BB937B7C83A31E7825B1D07AFB5
75C344B11A70B10AFE483B87EF1B3BB18CE735C24E6DBC1BACDAEA765A82DD4178C52FE2560793CEE0074E73C7E1DEBBD6D6B4C4B017A6F61FB3B701
F7753E9F5A10339BF86934F2E99A889E6925D978CAA6462C71815977F62D6BF13EC6CB4D99EDE3B8837B80C4E061B7633D090B57BE184AA74CD4D988
03ED6CE49E3008EBFA554BCD46C66F8AFA6DCC575134296EC8F2670A1B6B8C67F1140177C5E3FE11D9B4BBED301595A7F25D598959011FC557BC5969
F63F0D5E6A114920BF8C2C826572086C8E9EDED597F12750B47B5D3638AEA06912F11D94480955C1E4FA0E95B5E34BCB57F06DEB0B984ACF162260E0
87391D3D6802D45A9DD4BE0E8F53862125DBD9ACC100CE5CAE7F9D656996A2EBC29FDA77ED716BA83C4EEF399191948270719E074E29FA6EB10D8780
62B8B69239E6B4B2526356E41C01C8FAD67585D4171E0E9357D4A56D42FE78E4091B9DDB5B9C2AA8181DBB5006F7823559F57F0D5BDD5D1DD3026366
FEF63BD63EAB91F16F442471F64703FEF992AAFC3EB971A1E9B6D0CA158DDCBE620232542F7A7EAD7B683E28E9529BA8B6240CA5B70C29C3704FE340
58EFE9AEBB94A9CF231C1C1FCE8565650CA415232083C114D96448A3324AEA88BC9663803F1A6238CD355ED3E255CD8C77572D6AB67BD62925665427
6F4C9FAD68788B5A9E0D6B4DD0EC2411DC5EB664908CF9683D3DCE0D64596A364FF142EA68EE23313DA2A07CE14B003806975E81ED7E2668DA8DC102
DA506256C701B6B000FE2D48A2C7C41825B6F084EA9349346CCA1FCE6DC472390715BDF6E8F4BF0BC57B2AB34705AA3155EA7E51589F13AF62B7F0BB
40CC3CD9E5408BF43927F4AD2BCB57D5FC0FF66B771BEE2CD361078CED07FA5022B7870DC6BFA39D52F2E268A4B82FE52C4E5442B9C0C7BE052F87BC
46D7DE18BABB60D2DC58EF47E33BCA8C83C7AD53F01D8D85C787628EE6D22FB540CD1CAAC3E6C83DC56AC77165A5D96A97163630436D6A9BB7428144
AC14E7A7A7028033F4759F56F0DB6A973717B05E3ABB6EDDB02119FBA3FBBFFD7ABFE0BD5A7D63C3B15D5D60CAACD1B30FE2DBDEB374AB84D5FC2E75
9D62732F9CAE04209089CE000A3A9E3F5AAFE049CC1E1AB18619151E4BD70C9D4EDE723DBB5030D7CE7E28F87D40C910B1FF00D0FF00C2BBAAE175B9
E25F8A5A1E6451B2DDC31247190FC1AEE4104641C834C4C8E68C4D134659D430C6E462AC3E84571DA5CB756FF122E74C3A85E4F6C967E6049E62C031
DBDBF1AED199514B3B055519249C002B82D3AEEDE6F8B3792C53A3A35A0456520827E5E33F852608B1F104DF69B60354B0D56EA1669913C90FF21C8E
C3B74CD5EBD8B543E177BF96FE68AF92DFCC0B19DA8A719C11DCD52F8A854F86225DC370BB420679E8D5BDAEC822F0ADE3865E2D4E091907E5A004F0
AEA8DAC787AD2F641FBC65DAFEEC3826B66B9AF87CBB7C1D63F329CA93C76E4F5AE969898553D46F63B0B37B8908E3851FDE63D05599258E25CC8EAA
3D58E2B9B596C7C53A94D18904D6760402A0F1239EF91D8521A39AF8816A2D61D1AE64556BB92E0BCD200324F071F41D2BD32BCDBE2669F656361A7B
DB298DC4FD3793C7AF26BD02DAE619E34F2E6476281B0AC09C7AD081EC59A2B2DF5FD3535B4D21AE07DADC642638FA67D6B529882B87F8B07FE29588
7ADD2FF26AEE2B85F8B3222F86A08C91BDAE94819E7015A805B9D568993E1FD3F9209B58F91DBE415C9DCDEEB761E38B5D2E3D45AE219D376245036F
D6BACD04A9D034FD8C197ECD18C8391F74573573C7C58B439EB64DFC8D21A2CEB9797BE1FBDB0BC92FE49EDAE2E3C89A2750000DD0AFA62BACAE2FE2
3906DB488CB6375F21CFF9FAD7694C185731AEEB7347AE58687A7C8B1DC5C9DD24846762E33F99AE9EB83D46DE5B5F8A7637B2A1305CC7E5A3760DB7
14023535EBBB9F0D5BC5A98B89AE6DFCD54B88A56CF073F32FA738E2A7F1199E5D0E4D534ED42683CBB732A04236BF1919AA3F12C97F0C8B58F0669E
74544EEDF4AB9AAC074FF87F3DACCE3741A7F96493D484C7F3A4066473EBF3F846CAFED27927BC665DAB18C82B9E4B7AF4AD5D6B5A9ECE6B0D36DD53
FB46F8E003C88C6396F7EFF953BC0F8FF843F4D0181C45DBEA6B9DD7E3961F8A1A3DCC9C40E36AB13C0E0834016BC76751B4F075E47713ADDA4A554C
823D8CBF303D0718E3AD743A24A90F85EC2690E123B34666F40101AC8F89770B0783AE1091999D107E79FE95705BC979E008EDEDCFEF25D391531DCF
9628020D02FAF3C4705C6A715D35BC1BCC76F128CE36F76F5CD5AD0F546D734C9D566115D41218656419C107AE3DC573FF000E6C629FC3B246F3DD45
2C73324B1A4A540FA8FCEB7B44B0D374BD52EECF4DB6753B03CD2972D9627A7D7BD0064683AB6B3378B2FF0044BBBA4952D9491208C640F5FAF22B4B
C1DA86A57D0DE8D4E2900827290C9226D2EBCF5FA7F5AA1E1931BF8F7C48D19040F2C67DF9CFF2AED2800AE0BE275A5CB68ED75FDA122DB891545B05
0149F5CF7EF5DED71DF1389FF845800091F6842481D07343046E69563770470BCFA94970BE58CA32000F1542F7599AEBC48BA0E9EE23744F32E26C64
A8C741EFC8AD1927641A6451B90D2B0C81FC4A1093FD2B98D2ADE6B2F8A7A91B9CE2EED8BC47B30CAFF200D008D3D4351BAF0EDE59FDAEEBED363732
796CD2001E36F5CF715D30208041C83D08AE23E2640F7B6BA659400B4D35CF0ABD718E4D7676F1F936D1444E4A205CFAE0628064B583F6D92F35FB9D
3E3BC36AF6CAA5536826404649C1EA2B7AB89BCD2AD7C53ABDF833CB697FA7C9E524911FE0C719A6237745B9D427B9BF8B510A0C122AA051818C673F
8F15CC5F3DEAFC520962A85DED42967E88B8193EF5B9E187D4AD6EEEF4BD5665B89210AE9303CBA9E39F7E0567C2A4FC57B87DCB81660019E7A0FF00
0A434695CDFDEE97ADE9D6971702EA0BE263CEC0AC8C0673C76AB3E26D65743D29AEB686919C47129EEC6B2BC4CC9FF097F86D19B1895DB8EBD38FC3
355BE275ABCDA5D94EBBB10CF938E8323A9FCB1F8D03355E1D62DB4B6BD92FDA6B848FCC683CB509D3256AE699A9FF006D68897960EB1C920C7CC321
181E41AB37D322E8F3CDB805F219813D3EED73BF0D6DE787C34649F8F3A667518E828174134BD575BBAF105F6957335BC525BAEE4CC7F7C67A8A874E
D675CBAF145E6892DCDA46D6CBBC49E56778E38EBEFF00A1A4F1748343F11E9DAEAC6EEAD986554EADC703F9557D7ADDB43BFD2BC46532E1D96F0A2F
DFDD9C7F334866D4571E209F5896D126B35B6808DF308CE4E7B633D6BA5ACBD06D9A0D3C49292659D8CAE5BAF3DB35A94D09914E8D244CA8ED1B1180
EB8C8FCEB85F0C69F76BE28D7962D4A6578E440CECA1CC9F7BAE6BBFAE37C24C5BC5BE25DCBB49950E0F5FE2FF00EB5007569148B6A22794C8FB7697
3C127D78AE33C142FAE2D6F2D6298DB4315CB9691402CCC7B0CF618AEA34B9649F4D69A69B3BD9CAB7F757247F4AC5F87C54E9377B4E40BB7192724F
D4D1D43A17742D4EE26D4AFB4BBC7F365B4208976EDDCA78E4557F135F6A1A7EA7A59B7B9C5BDCDD244E9B474C8CF355FC3EC0F8E75EF9B3808393EE
78A3C712C6350D0226701BEDAAF8F60451D03A9BBADDEFD82C3CDDC5033AA17FEE03DEA94E9A92FD8CE99786E2DA590799237CC40F5CFA558F105D5B
DBDAC515E22B5BDCCA229371EC7D3F1C562DEE9B2F87047A86937727D915D7CCB666CAB03E86803B1A2988C1E35719C3004669F4C4723AE6A3ADE93A
8D8C42EADE4B7BB976798D160A74F7C51E23D4B58D1E5B4D9790B45712842ED10057F0FC6B4FC5BA79D47419D130258BF7A87FDDE6B04412F897C317
1752A10CB1A8B7F5F9793498CDDBD3AAC725AC1697D13C9293B8B44380075A6EBFAA4DA5585BC51C9E65E5CC8228DD94004F738E955FC2172FA9DA8B
C991834082DD0B77C753FCAA878F6CE496F749B862FF006649C2C855B6E39E39A02DA9A7A9CB77A1D92DFBDDBCE88CA27593A6D3D76D6F4322CF0A4A
872AE030ACABBD1B4B96C4FDB8CB710200C3CC9D88C0E98E6AFE9A50E9F07970881360DB1839DA3B5022D5635E5C99F551609766D995370200CBB1ED
F4C56CD72BAAD8D9F88B549EC6791E0B8B50A61910E1B1DFF5A606A584B7C753BA82EA55923895761540339EE7DFDAB0BC5135C5BF8AF487B38BCDB8
65291A9E8324824FD055BD020D4749D55F4CBC9CDD42D0878E63F7B8EC7F5A3560A7C75A31CFCCB14991E830690C9350BCD5B4636B3DD5CA5CC12CAB
1CAA230BB33E98E4D696B9A8FF0065E9CD75B4361828CF419EF59BE37745D260566018DD47B46793D6AEF88352B5D3F4BF3EE2159B790238DD4105BB
668115AF65BFB3D37FB516FBCD5545768DD02A953FFEBAD8B0BA4BEB286EA3042CABB803DAB9AD6C3FFC2292CDA9DC8124B10D90AE1547B0C75E2AF7
872E47D974FB5465F2CD98907A93BA802B000FC44F71699FE5FE26A5F1C67FB047240332E7E9CD408CADF11DC6EE52D00C54FE37318D0D448C0033AF
7FAD03EA6E5A716900FF00A66BFCAA7A86D71F65871D362E3F2A9A98828A28A0028A28A0028A28A0028A28A0028A28A0028A28A0028A28A0028A28A0
028A28A0028A28A0028A28A0028A28A0028A28A0028A28A0028A28A0028A28A0028A28A0028A28A0028A28A0028A28A0028A28A0028A28A0028A28A0
028A28A0028A28A0028A28A0028A28A0028A28A0028A28A0028A28A0028A28A0028A28A0028A28A006BA2BA95750CA7A82320D27CA89C0C281D00EDF
4A7D140185A1A497B7F73AC5C4650C83C98108E5630793F89FE55BB4514005376AEFDFB46EC6338E714EA2800A28A2800A6491472AED9515D739C30C
8A7D1400D550AA1540000C003B535628D0E51117E8315251400CF2D0163B172DF78E3AFD69B0C31429B618D235CE76A2803F4A968A002A136B6E58B1
8222C7A9D8326A6A28018B1C6885111557D00C0A3620180AB8C6318A7D1400C8E28E25DB122A2F5C28C0A7D145002606738E7D691943295600A91820
8E0D3A8A00A9369F6571B3CFB5824F2C009BA30768F41ED568000000600E8052D14005413DA5B5C106E2DE194AF42E81B1F9D4F45004325ADBC90885
E089A31D1190151F8548AA1542A80140C00070053A8A002A116D6FBD9FC88B737DE3B064FD6A6A2802B7D82CC800D9DBE07FD331FE148BA7D92A955B
4B70A7A8112E3F955AA280204B3B58DC3C76F0AB8E8CA80114C4D3EC92E0DCADA42B31EAE100356A8A00AAD6366D762ECDAC06E40C098C63781F5EB4
9369D633C865B8B2B6924FEFBC4AC7F322ADD14014FF00B2F4EFFA07DAFF00DF95FF000A93EC76A22F285BC3E5E776CF2C633EB8AB14500451430C20
88A248C1EBB540CD356CED5576ADB4200E7023153D140159EC6D1E4691ED60677186631825BEA69A74DB02814D8DB151D0189703F4AB74500578ED2D
625744B7851641870A80061E87D7AD57B1D174DD3E4F32CACE285F0465474CD68514014AD34AD3ECE779ED6D228A59092CC8B8CE7AD35F46D2A4CEFD
36CC93D4F90B9FE557E8A008E28D228C47122A22F0154600FC28963496331CA8AE8C30558641FC2A4A2802A8D3AC4631676E31D3F74BFE1497F616DA
84021BB8848A0E57B153EA0F6AB74500654BA0E9B3DAB5B5CDB79E8C413E6B16271D39EA29DA5E87A7691BBFB3E0316E1839919B8FA135A74500675D
68BA75D4C679AD10CAC305949527EB82335696D2DD6D7ECA214F236EDF2F1C62A7A28032ACBC3FA669F2192D2D95327705C92AA7D40A96D347B0B4BB
92E60B6449A4258B0F53D703B56851401972F87F4795D9E4D36D59DBAB18C67F03DAAFDBC115AC090408123418551DAA5A280239A28E789A29903A38
C3291C1154D344D2E30A134FB65DA723118CFE75A14500665EE85A5EA12996F6CE39DCE325F27A0C530787348F21616B30C8BD999B9FD6B5A8A00CDB
1D134DD3E533595A242E463E52703E83381F856951450051D434AB1D4B67DB6D966F2F3B7248C7E550DB681A4DA126DAC9222C72769233FAD6A51401
972F87F48B818B8D3E1947FB796FE752D8E8FA769B2349656A90BB8DA4AE791E957E8A00C87F0F696FAD8D65EDF3783186DC7191C038F5AD7A28A006
BA87428C321860D634BE15D1E720DCDA99F1D048E4815B7450065D868561A7481ED1244C2ED0BE61C01EC2AB7FC225A2FDA96E4DBC8665E8E677CFF3
ADDA28B05CC4BAF0B68F7854DCDB3C9B7EE8333FCA7DB9AD1B1B282C61315B86DA58B1DCC5893F5356A8A002AB5D59DBDE22ADCC4240A772E7820FB1
15668A00CF8B48B44B88EE1C493CB17FAB699CB94FA67A52DFE9363A98C5F5B898636E198E319CF406AFD14019F61A3D869AC4D8DB2C248DB852718F
A668D5748B5D56345B9560D19DD1C8870CA6B428A00C4BCF0EDAEA362D6BAACD35D827E566382BF4F7F7A9F47D1D349531C3777124014247148D9540
3D2B528A2C1732A4D0ED8DDC9750493DACB2FF00AC36EFB77FB9E3AD3A3D16CE3B192D144A1246DCEE2421D8FA922B4E8A00C08BC25A540CCF6C2E21
91CE5DD26605FEBEB5BA0600033C7A9CD3A8A002AB5FD8DB6A366F6B771892271C83FCEACD140146CB4CB7B32ACAD2C8C8BB11A56DC547A0A4BFD2ED
AF9E3964DF1CD17DC9A26DAEBED9F4ABF45006741A4DB47742EA5696E275184799B714FA7615A34514005637FC23F6F1DD4D776B34D05CCF21792556
C96CF6C74C56CD140142D34E1024CCD33C97138C493F463C6063D3159B0F8474C82F0DEA3DDFDA89CF9A66F9BF3AE868A2C17306EFC2B63797515D5C
DC5E3CD09CA399B95FD2B57EC50B58FD8E6DD344576B79872587B9AB345160B993FD876CD6C2DDE6B97814F11349F2E3D3A6715A714690C6B1C4A151
4615474029F450054D434EB7D4638E3BA5DC2391645FA8A6EABA6DBEAB682DEEB77961C3FCA7B8ABB45002018181D2968A2800ACC9344B392F9EF409
239641893CB72A241EE2B4E8A00A73E9F0CB62B64A648A2500011B60E076CD54B4F0EE9F651B4767E7408C30C1252335AF4500625BF8634CB6BCFB5C
2265989CB3094FCDDF9AB375A1E9B7938B8BAB512CA0E4333B647D39E2B4A8A00CF9F48B2B88A08668B7450125632720E477A61D16D99638DDE678A3
60C222FF00271D38AD3A2800A28A28018E8B223238CAB0208F515058D8C16168B6B6E0F94B9E1B9CE6AD514015ECED21B2804302E1724F3DCD3AE6DA
1BB84C37312C91B755615351401991E8B6AA115E4B896343958DE4CA8FC2B48000000600E82968A002B365D16D24BB6BB06449D8E77A3608FA569514
014AD2C16DA679E49649A671B4BB9EDE80553B8F0E59DCDD25D4D2DD34F19F964F37047B7D2B668A00C5BFF0DD9EA3B3ED935D4810E5732F43F954B3
E836773A7AD95C19A58D0E5599F2CBF438AD5A2803223D02CFECCD0DCB497395D80CA73B47B7A549A7E8D69A78FDCF98CDB360676C955F418E95A745
00639F0F599BC6BCF32E05C118F304986FCEAC6A7A5DAEA91797761D978C00D8C7BFF9F4AD0A2802B58D9C7636E2089E465073991B26ACD145001451
450014514500145145001451450014514500145145001451450014514500145145001451450014514500145145001451450014514500145145001451
450014514500145145001451450014514500145145001451450014514500145145001451450014D2C01009193D067AD6478A6CE3BAD06F19B21E285D
9181C11C67FA57390E9173AE78334D9ADF68BE554DB234AD9083838EC09C6695C763BCA2B87D735799BC4BA77871E6F2A22A8D732EEC1938E99EDD2A
6D6750B4D0B56D367B29879524863B98D24DC369E871EB45C2C6B78BB509F4ED0A49ADDF63B3AC7BFBA8638C8F7AA1A069BABD86BF3992EE49B4C685
4A191F712DEDFAD53F88FA5C4742B8D404B3093CC8F29BCEC3C85E95D06936F0D869497026B974F2039124A5C28C64E01A3A87435E8AE37C34C3C570
DDEA3A906920F38C7043BC845000E703BD49A65FBE99E2F97C38CC5EDE48FCEB6DCC4B27192B93DB834058EBA8AF3CD7EC67D37C5FA443A55DDC411D
DB10D19998A8E79C0CFA1ADEB0D22FECFC4D77324A7FB32788643484BB49C73EA0F5E680B1D2D15E7BA72DDDA7C48BAD36D2EA536FE4866F31CB1504
02719E339E2AFEBF7B73A1F88AC3EC259A3BFCC5246EC586ECF0DCF7145C2C7674571FE2BB79746D14EAD65753FDB6DA442D23B93E60240208E98E7D
2B71A33ABE916D2FDA67B632C6B2EFB79369195CE3E9CD311A945719E0A96F757D16E92FB50B9768EE4AA48AF87007BD3FC25797F2F88B5CB0BBBC96
E61B47558CC98C8E4FF852B8EC52F105D6B3A478974989754926B7BB9C06428071B86471DABBDAE2BC62A1FC5FE1957FBA673C7E22B4BC59AA5C69E2
C21B74976DD4E239248BEF22F1D3F3FD2803A3A2B8CF13CD2785D2C750B3B9B87125D2C534334A596404125B9EFC75AEA750B537B6725B09E580B8C6
F8C80C2988A3A36AADA9DE6A29B63F2AD27F291D09F98F7CFE95B15E77F0F2C66934ED4A217D73118EF1941461CF03939FA5743E2DD7CE83A743E510
D7570E2288B6383DD88FCBF3A43B1D1D15C8EBD3DE78674A8752FB6DC5CB4722ACE92B6E590375C0ED8AEA6DE78EE6DE39E16DD1C8A194FB1A6225A2
8A2800A2B987D5EE752F14CDA358C82182D23DD71301962C7A28F4E696E7529741D6F4EB1BAB87B8B7D40B223C9F791C63BFA72295C763A6A2B92F1C
DF6AFA3E9C751D3EF95230EA8616855BAF1D4D3FC4B2EAF6DA0BEAB61AA9B7315BABB44D023063C67A8E3AD30B1D551581A1BEA1AB784AD25B9BB686
F274DED3468B9FBC71C74E98AC6F0F6A7ADCDE20D5B46BCBB137D94656508011E98FAE475A047714573DE0FBAD56E34B97FB6A3759E39D91778C12BC
75FCEBA1A0028A2B88F19DFEBDE1ED3D6FADB544951A608636B6518C827AFE1401DBD15CDDDCBAD7F64CDAAC17091158BCF4B46841F942E76B1EBBBA
FE353D95EB7893C3D05DD85D3D9C92724A00C558755340EC6ED15E7BE1CD535FD5749D5C4B70EEF6ACE2368D4077603851F956937892EF4BF03DBEA3
AA42CBA8C80A244E36966DC42923E9834AE163B0AA5A91BF16B9D316DDA7CF49F3B48FC3BF4AC7920D6E3D1FEDEBA9B1BD11095A268D7CB38192B8ED
F5AD2D03558F59D1EDEFE3014C8BF3AFF7587047E74C46278375FD535ABED4E1D4A2B788D9388CA44A7E56C907927D8D75D5E7BE0EFB49F15F8A6184
88D5AED99A5C6E2877BE001EF93F956D45ABDCE99E248B46D525F3A3BB4DD6D3EDC1CF7538A432978C75CD7B439EDDADBECAD6F732F96BF212CBE99F
7AED6B88F890A1DB425248CDF28C7AF22B63C55E2387C3B60B2B2F9971312B0467A3118CE4FA0CD006FD15CD5FDDEABA3E90354BBB85B811056B8844
61005279C1EB919AB1AEC97B3E82D7FA45F7D9CC7099C7EEC379836E71CF4A6163768AE3BC2FA8EB3E22D022B97BBFB2B29646956204CAC0F519E00C
607D734ED07C41A8EB125D69A5122BCB390A4F70172A0038E17D4F3F950163AFA2B98D3B59BD8FC5D3E837EE93621134532AED38F422BA7A04145145
001599ADEB56DA3C08D36649A56D90C2BF7A46F415A75E59A95ECFA8FC5CB5818E63B49D52351D801927F3CD008EED25D716313490DA3038630A960E
07A03D09A8749F100D62EB52B7B381D1ECC2A8F3C6DCB9DDC1C741C0FD6B7AB22CF4E6B5F10DFDEA22AC577147B88C72EB91D3F1A4332B47F116A779
E2CB9D16F2CEDE316C85DDE3624F6C633F515D65711A4F1F16F5B18E0D9A1CFE1156A6A3AF48DE218741D34A7DA1977CF2919F257AF4EE714033A3A2
B02D754B9B3D7D747D4A4595A78CCB6D30014B01D548F5A8E5D6A4BAF13C9A0DA4C2DDE18BCD794A862DD3E500F1DFF4A623A3A2B9AB3D66EEDFC52D
A0EA2C92B3C5E7433A26DC8F423F03D297C597FACE936136A3606CDE0850174951B70E719183EF4058E92B90F10F89757D0AEA212E9D6ED6B3CFE545
279A73ED91F4A82F3C51AA45A069FADADB4612E655436E1492C0E7273DBA71F5A4F89D8163A493DAFD38FC0D21D8ED867033D7BE296B3B579B5182D8
CBA64304CC8096490904F1C631593E18D7751D7EC24BA1159C5B18A18F2D9561EB4C563A7A2B9ED22FF5CBEBA985C41631DB43298D9D0B12F8EBB79A
E86800A28A28033B51B8D4E2751A7D8C770A46599E50B83E98AC3D23C5B77AA6AD26989A498EE2DF3E7EF970179C7A575B5E6FA83DC68BF14049696A
653A84402A16DA198F539F6C521A3A4D4B5FD434EBDB5B37D23ED32DC6EDBE4CC3B639E471D6B6EC9EE64B656BC856198E728ADBB1F8D57D3EC3ECCD
25C5C389AEA6397908E83FBA3D00C9AD1A0028A28A620A28AC1D4F5F4B5D66DB47B6311BCB85DC0CA48451DBA773401BD4573A75DBAB2F105AE93AA4
501FB6026196027191D8835A3AD6A434BD3DAE7CB32396088B9C6589C0C9F4A2E163468AE6F53D6755D262B7375630CCD7132C61E163B533D8E79CFE
95D2500145148480092700753400B4573B61ADDDEB12DD1D252D7C9B690C7BA6639623D874A9F40D7975592EED6685ADEF2CDF64B19E9F51ED4AE3B1
B745145310514D660AA59BA0193580BACEA12DA3EA56F690CD60A18EC572652077F4FC2803A1A2B12D355BDBBF0F5BEA56D60B34F30DDE48902E064F
73EC0557D07C4536B9A64F756F6012686531189A5EA7EB8FA52B8EC7475CBEB3E318B47BA482EB4EBB5F31F6231C056E7A8353683E246D5356BDD367
B26B6B8B3197CBE41E6B1FE26C890A68D24870AB7793C678E334304B53B8072A0E08C8E869D58975AD4F69602FE6D3DFEC9C33B0901744FEF15C7E3D
6AF1D4AD8E907538E4125B795E6AB0FE218A622ED15CCDC7886FAC2D3EDD7DA662D9F051A2932403D377157F51D5A7B2B55BA8F4F96E6131798CD1B0
CAF19E94AE3B1AF45731A4F8ADB5DB45934AB091E40712090E153DB77AD6B695AA45A92CC16378A681F64B13F543D698AC68D158C9AD35DDDDC41A5D
B2DD0B73B657F38200DE83839A934AD660D4E6B8B608F15D5A9C4D0B73B7D39E8680B1AB4566EAFAC5BE9314665CBCB33848A253F33934C7D55AC8C1
FDA91240276D8ACAF901BD0FF8D006AD159DAC6AD0697144645692599F64312F576F4A813569A0B9860D52DC5B99CE23915B7267FBA7D0D170B1B145
14500145159375ADC31EA634CB58DAE2F3019914E022FA93DA8035AAADFDEC16168F7370DB634F4EA4FA0F7AAF06A919BD5B1BA5F22ED94B2465B707
1EA0D739F11AFEE2DF4A445B593CB59E3713020838ED8A571A46CE91E25B3D52FA4B154921BA8C6E31B8E715B95CF6822DAE6F5B505D1E5B4B89A205
A6947DEE9C0AE86840C28A2AB5F5EDBE9F6725D5DC8238A319626988B34561C3ADCED6A6FE6B23159601DCCE0301FDE23D2B524BCB78EC9AF1E5516E
A9BCBE78C50058A2B06DB5EB8BA80DDC7A4DC7D871B9652E03B2FA84EBFAD6A58DEDBEA16A9736B209227E8476F634012CF3470426595B6A2F53827F
95665AF88B4BBCD5069B04CCD7254B6D2847039CF3539D5ED8DDC96D0ACB3C91712089376C3E86B9AB79A09FE28B34041FF433920639E295C763A5D7
3568745D35EF678DDD54850A9D4935634FBC4BFB186EE3565495770561822B33C49ACB68D69E7C9A735D5B606F60E00539E9835A0B70C74D8EE6080B
968D5D62040382071408B945733A678BA1D4A59ADA0B19FED91BEDF232338F527A015763F10DB2CD35BDEA35B5D45FF2C4FCC5C1E9B71D69DC2C6CD1
58AFAE9B7BCB68350B47B55BA3B63919C119F43E86B4350BFB6D36D1EEAF2411C49D4FF4A00B54561AEB5722D96FA7B1D964E01DC1F2CA0FF111E95A
17378C9622EACE1FB5023700AE1723D726802E515CCE8BE2B5D69192D6C25FB4C6E4491EEF9507625F18E7D2B4B4CD5E3BFB8B8B6789A0BAB66C4913
303C7A823A8A0763528AA3717AF1DCADB5B41E74BB7730DE1428A65AEA42E6EA7B710489240A0B038E4F703FC68110DCEBF6F6FAEC1A43452B4B30CE
F03815AF5869AC46FE244D325B0649CA332CCD83F28AD1D42FEDF4F804B72C406608AAA32589EC2802DD1590FAB4D6C8B2EA162F6D031C193CC0FB73
D0903A55EB9B830DB19E289E7180408F0723D6802CD15CEE9DE2AB6D43CD8EDAD2E1AE636C187033F527F3AD0D37578EFE79ADCC5243710FDE89FAE3
8E7F5A02C69514514005472CB1C319799D5107566381525737E2BD1ED27D12EA6C3AC91A1607CC3C93EB9A00DDB6B886EA1135BC8B246D9C329E0E38
A9EB2BC3091A7876C162E13C918F7F7AD5A0028A28A0028A28A0028A28A0028A28A0028A28A0028A28A0028A28A0028A28A0028A28A0028A28A0028A
28A0028A28A0028A28A0028A28A0028A28A0028A28A0028A28A0028A28A0028A28A0028A28A0028A28A0028A28A0028A28A0028A28A0028A28A0028A
28A00C6F14DCC36FA0DDACA72D244C1107563E83F3A83C157113F86ECA056FDEC3161D08C1073E95BCCAAC30CA08F714891A267622AE7AE0629751F4
387F14699F62F17E9FE21922135A12B14E0AE76760D8AE96E24D3BECE0DAC56734AF81128553927A74ED5A8E8B2214750CA7A82320D4505ADB5B9260
B78A227825102E7F2A62393F88F7F6A7C35716AB7319B81226630D96E083D2B7ECA6B6BDD252DA0B889D9AD82955704AE571C8AB5258D9CB379D2DAC
0F274DED182DF9D3A0B3B6B762D6F6D14448C1288173F95219C9F80217D1EDAF749BE2B1CD0CE59431C6E523A8A4B6B16D47E24C9ABC2C3ECB67104D
E0E43B152B81F99FCABAABCD3ACEF8A9BAB7494AF424723F1A9A0822B78845046B1A0E8AA3028B05CE17C47A959CBE36F0FCD15CC6D146583B03C29C
8EB5DD2CB13C5E6AC88D1E33BC30231F5A8C58D98040B58003C91E58E69EB044B1792B1208F046C0A02E0F5E2981C4D84D037C58D424F3A32AF68A10
EEE18E10601F5EB53F8BA5853C57E1D6799479739DE33F77A727F3AEA3FB274C0411A75A641C83E42F5FCA9CFA758C921924B2B6773C96689493F8E2
905CC0F8873C23C21749BD7749B0A00793F30391F956CE92D149A45B42B2C6E52045708C0EDF971F85492695A74854BD95B9D9D3F7638A7DAE9F6766
EEF6B6D1C4CFF78A8C66981CAFC3B096DA5EA2B33796D1DDB07DE71B7D3355FC157515CF8C7C4524470B2BAB283C123279C5754FA1E9B25E1BB7B446
959B731C9C31F523383F9531FC3DA449766EDEC53CF2DB8B8241CFE069582E73BE339235F16785B2C322E5B77238E5319FCEB47C63ACDCE930594766
4235E4E2232919F2C71C806B5E6D234DB894CB3D85BC9231C9668C139FAD2DE697637D6696B736E9242806C07F871D30680B9C3FC47B6B4B6D26CA12
ED2DCBCEA4CD231248E9927A576B6572F36A57B17981E28847B303A657279EFD8FE34A747D3DACCDA496C9242D8DC1FE6271EFD6A6B0B1B7D3EDC5BD
B2954073C9C93405CE5BC09E5DBC1AD06629B2F2466E3A0E79ACDF1D584F3E95A1EA3307B886D0837447076B6DE7F4C66BB59348B092E8DD496CA666
C163920363D46707F1AB92C692C6D1CA8AE8C30558641FC28B05CC8BFD1F45BFB4F3EF21FB4418F301699C8231C11CD5DD1D836956C56DC5BAF96311
0E8A3B5471E8D6311411C45638F948F79D8A7D719AD2A620A28A280385F0BDB4BA778F75B86E14AFDA143C44FF0018CF6FCE93C7968DA9F88FC3B610
963279CCEFB7AAA6572DF929FCABB0BAB0B6BB7492543E647F724462ACBF4229B6FA75BDBCA6701E49CAED32C8C5988F4F6FC290EE73BF1380FF0084
41C9ED3C67F5A9BC51228F87970C180DD6A9B727AF4AE82FACADF50B492D6EE3124320C329AA5FD8362D646CE7F32783CBF2D5657C855F41FA7E5405
C8BC1B93E11D2C9FF9F75AC4F0ABABF8EBC49CA310EA323F2C7E95BD63E1CD3EC2CDAD2113BC0DFC0F33103E9E9490F86747865F3A3B4C4A4E59F7B6
5B9CF3CF3CD006CD14514C415C67C530C7C2A98E40B94C8C75E0D7675C87C46B88E3D0A1466524DE4595CF6CE79A068DDBEB88CF872E6E410D1FD91D
F8E98D84D61FC33B49AD3C25199D4AF9F2B4A80FF748007F2ADA9349B2BB83CADD21B761CC4B210879CF4AB371656F756EB04D1E63520AAAB15C638E
D8A00E57E1C6026B6A08205FB63154BE2E59CB2E9B617883296D2307F6DD8C7FE835D5D9F86F48B090C96969E53960D9591FA8FC6B46E2086EADE4B7
B88D648A452AEAC38229582FA90DFDC469A45C5CEF0B18819C37B6DCD73FF0E2C24B3F0CA4B3392F74E65C11F76B55741B6FB325AF9D72D6A800F21A
4F908CE7078C91F8D6AAAAAA855002818000E00A6238BF0418D3C49E268F703235D06C8EEB96FEA69BE2D85EF3C71E1D82DB719612D2C857F813239F
6E86BA08F40B38A34F28CB14A8CCC2646C39DC49209C7239F4AB365A65BD9CD24EA5E59E51869656DCC4761F4A43399F88BB43E84C4F22FD71FA554F
8B3653CFA6595DC6098EDA461263A8DD8C1FFC76BA7D47C39A7EA9379B7C6E25C1DC8A666DB1B7AA8EC6ADA69D07D81ACEE1A4BA85F2185C36F247A6
68B05CA37BA1DAEAB6262BBBDBB9ED641BF69940523A83903A535D2DADFC173C76A245B74B49113CCFBD8C11CD4F6BA1456B1A4115DDDFD990604264
E31E9F4F6A9B50D22D75040971E688C214D91B9552A7B60500647C3AFF00912EC7803EFF00FE846B3FC0193AEF8A5CE39BEFFD99EBA0B1F0F58E9D6A
D6D6525D4311E8AB70DF2F3DB9A8AD3C2FA7D95D35CDA49770C8E7326C9DB121FF006BD68031DC01F1710EE19361D3F3AEA6FF0053B1D356337F731C
0246DA9BCF5359137842CE5BE37E2F3505BDC604DE7F207A7D29FADF856CB5D8AD16FE4959EDB03783CB8EF9FAD006F8208041041E8452D311163454
418550001E829F4C415E6DAB58BE91F1474ED470C2DEF6503791C0720A919FC6BD26A9EA3A75AEA5686DAF6212464E47182A7B107B11EB40172B3A29
A49B5995639336F043B5C0231E6139FD00FD6ABC7A5DF448214D666F257850D12970BE9BBA9FAD5E82C62B5B46B7B7DF1AB672C1B2D93FC593DE8039
2D1496F8ADAD9FFA7551FF00A056769F6D3C7F162F967B992179559E22A01DEA40207391C01FA574F0F84ED6DF507D421BDBD4BB97FD6CBE60CB8EE3
A55ED5B45B5D51A19642F15C5BB6E8A788ED74F6CFA521DCCDBDD074F8B57B2D56EEE2FA5BA49923889901193D063038F5A86EB518A7F18B699A65B4
51EA0B0667BD6504A271C01DCF22B6A0D2C09D2E2EEE66BB962FB864C055F7000033EF546EBC331C9AE9D62CEF26B4BB64D8FB002AC3E9F850073CAA
B0FC5BB58FED324F20B42256723EF6D63803B0C63815B5E2BB8179E01BF9F6EC0F0E40CE7F8877A75DF83ECA5BF8AFAD669ADAF1010F3A9DCD282307
767B904F356B52F0E5B6A160962D3CD1DAAA04312370D8EFF5A350D05F0BC692F8534B59515F16E870C33CE2B13E260FF43D27FEBFD063F035D1E91A
57F6545E4C7753CB0AA8548E4208403D2A86B5E168759B959AEF50BD0B1B89238D1942A301D4714740EA7435E796F7ABE10F156A967206FB35FAF9F6
A00E0C9D368FC4E3F0AEA352D5AD3C33A74526A5737132BC9B15CA86624F3DB1C556D52C6CF54D7F479CA869200D383B4FDCC7CB91F5A18235B47B4F
B0E976F037DF0B973EAC793FA9ABD4514C414515535082E6E6D8C76777F659720F99B03FE18340116A9A82D8C4A150CD7129DB142BC963FE1EB5C0F8
FF004EBAB11A76BB34B24B7293012E0FC883A803D3D2BA383C31A945AC36A87C4123DCB26CCBDB29017B803381D3B558D63C3773AD59FD92FF005576
8321B6A40ABC8E9DE90F6372D675B8B58A74FBB2A071CE7A8CD25E5EDB585B35CDE4C90C2BD5D8F1595A4E8779A5C36D6E9ABCB25BC1C08DA25E57D3
34EF16E936DAC688F6F793B41144DE7175EDB41FF1A606ADADD417902CF6B2A4B1374653906A6AE6BC0160D61E13B547DDBA52D290DDB278FD00AE96
8130AE4F5FBAD3E3F12D8470E9E979ACB0FDD1270225FEF37D39AEB2B97D5FC3935DF89A0D5AC6E0DBDCC2832EE37237518C7D2931A333C448EBE37F
0E35D5C2B4ACE7F74AB855FA1EBD6B67C453DADDE8F2DADF232433DD0B5DE0FDC39E1FF3151EADE163A85C5BEA0978C9A95BB06498AFCA7DB6F6A9F5
4D01B52D121D3A59F07CD124D2E3963C9247BE4D007337675DF0AA44BA96355D1D5D4976FBD1723193D7AE3DABD061912689258CE51D4329F5079AC7
BCD2EFAF2D1B4F9A781ECDC05772A448541CE3D3F1AD8863586148A31848D42A8F61C508192514535D43A321240604645311CA6853D9BEA17D0F86EC
E28A15936DCDC92482C3FBABDFBF7AC9F0CCCF69ACF8A3506713BC526DDDD031CD6BE81E17BED09EE61B6D454DADC3EF394F9D4FB510F83FECB7D76F
677863B5BCC19A365CB641CF073DE915A1D444FE6448F8C6E5071E99A929AAA1542A8C003029D4C9239995217771B955492319C8AF3CB6B2D7ED2DDB
55F0B5D24D6371B9CD8CBCECE790074FCB15E8370865B69635386742A09EC48AC3D3348D4343D3D6D6C268EE54F27CE272AC7A91EDED498D16BC2D7D
15FE8504B1406DF6E51E2FEEB0EB587F0C72747D458F53A849FF00A0AD74DA5590B0B310960CECC5DC81C163D71597A6E8579A4CD751E9D730ADADC4
C653E6292E84E338EDDA80327C36777C47F10B290576A8C8F5F9693E252ABFF622BE369BD00E7D38AB1A47852FF45F10DC6A1657904B0DD7FAD1329D
DC9C9231EFEF5ADE20D061D74DB25C3158E0666CAB10C188E08A077D4BBAC229D12F9081B4DB4831DB1B4D71BE0631AFC3BB8FED2DDF6532BAFAFC87
68C8FC735D0C963ACDD58BE9F753DAAC2CBE5BCC80EF7523078E80E2AC5EE8D13787A5D22C55218DE3F2D73D00EE7EB408E42EA3F11E8560619C26AD
A31001FEFA271DFAD7633DC4373E189AE2DB985ED19907B6C3C556FB16AB6DA7269B04B14F198FCBF3A4386518C74EF56A7B38EC7C333D9444F9715A
BA027FDD340189F0C618A3F07C4F181BA5964673EA738FE40557D21E55F8A1AC451E4C2615690678070B5278121BB8BC256AD61E4B2C8CEC5662460E
E3C8201ADDD23491612DD5DCEEB35EDDB069A50B8E00C003D80A00C4D2C69D63AFDF5AE8103DC5D4ADBAE6577FDD45CE71C0FD2A1F0F4CD178C3C493
5CB2318514B088606064F7EFC549A2685ADE87A95F7D99AD64B6BD7F30BB31CC679ED8C9EB4EB5F0CEA3A7EB37D359DC2496FA820133CBF781EFC7E2
681953C44C57E26681E6B7EE769DA0F40C777EB9C574DADE8D6DAB3DA3DDC8EAB69279B80700FD7F2A8F5DD062D58DBCEAFE55E5AB6E825C670720E0
FA8E2A416FAA5D288AFA4B78E2CE1FC9C9320FC7A502395F154F72FE3CD0A3B7D800E632E095627BFBD6CF88744D4356860371A9416B1DB3F9C4C711
3C8EFC9ABDAEE869AA2412432082EAD58343201C0F6C7A5235A6AB7900B6D424B54848C4A612C5A41F88E280B9AB0389208DD5F7865043631BB8EB52
D323458D151061540503D00A7D3105703E037925F14F88A49D8B49BD464F5C65BFFAD5DF5627F64B59EB33EA761B375C2E2689B8DD8EE0FAD263460F
8FD9E2D5F4096D891746E76AE3B8C8E2AEFC495CF84DCF7134647E75A834C92F753B7D4B50D8AF6C0886143B8293DC9EE7E98AA5E26D1F56D6ED5ACE
39ECE1B7F3030243162074CF6A02E6E581CE9D6C7FE9927F21566B334A87548008AFA4B5785230A9E5060D91EB9AD3A620AE1BE294EF1E93650827CB
96E06FF4C0C75AEE6B3B5AD26DB58B136B7238CE55B192A681A23BFD321D63435B1791A38A445C14F4C5735F1014E9DE088ECAD37B441D622D9CE00F
5FC6B7A08B5A8208EC82DAB222EC13EE2085C607CBEB566FB47B6BED15B4C9F26365C6EEFBBAEEFCF9A4058D3A248F4DB68902EC1128000E318AE2FC
0B3DC4365AFC9129648E7678549F973F3671FA57470C5AC4162B60B0C0C51046B7064E318C676F5CD59D1747874BD33EC8BF397C991B1F789EB40181
F0D99E5D22EEE2520CB2DCB17F5CD10A4517C52942280EF6659BF4ABDA469775A04B7315A5B8B9B69DFCC5C4814A1F439ACF9F49D7AD7C5BFDB70430
DE6E8BCB68C4823007A0CD032FFC417D9E11B93D8B203F4DC2B6B4A39D2ACC8E8604FF00D04564EBFA56A1AEE9AD66D2436B1BED623059811CE09E87
9F4AB96506A96DA5AC0CF6AF73180A8C0305DA3D7DE8EA2E860781224FED2D767DA37FDA76E40C71D6A4BA8526F89B6864507CBB5675FAD58F0E687A
BE8D7176F2CF69325D49E6381B8107DB8A64FA4EB6DE264D613EC99488C423C9C11EB4BA03DC67C434074DB07CE0ADE21047E355FE27A17D16DC8C80
25E4F61D3AD5EF11E8BABEB496D1896CE1482612F058E71EBC55F7D3EE755D3AE2CF5A8ADC2B1C466124E3FDA39EF4C0AE2CF54BED063B6FB559C6B2
42A37C71B7DDC7B9AB7A558CBA56842D269FCF6891B0E0638A6E9D67AAD95B25A1B9B79628FE549594EFDBEE3A66AE4D0CEB6262B474F348237BE7A9
EA7EB9A00C0F87D02C5A2CEEBFF2D2E5F3F5069DA6AFFC5C2D55B017FD1938F5FBBCD4FE1FD2755D1ED5ADDAE2D664672F921B209A65AE8FABDB6BF7
1AA7DAAD5FED2A15E22180503D0D003B5DB7BDB8BD32E87751C57F6E837A38E1C1E82A3F0F6A1772EAD3C1AC598B7D40C63054FCAEA3D2ADCDA5DDAE
B526A905C0DC54288BA07503A1F7A9E3B6B89F50FB7CF0A432C7198E242DBBAF72476A00C799FF00E2E55B2F3FF1ECDD47B541E2D0F2F8AF42898E20
1206E4E06777F80AB4DA2EB8DE215D63ED565BD633188FE7DB8AD5D534BFED286D9A4655B8B6612230195DD8FE5401278806ED06F41C7FAA3D6A8F82
E59A5F0CDAB4EC59864027AE3B53EEEDB55D52D8DADC086D2273890A36E661E83D2AFB432DB69E20D39230C8BB53CC6200FD3AD0230BC130C6A3509C
20DEF3905B1DA96DCFFC5C4B9E08FF0045FCFEEF34EF0E699ACE911CF14CB6730964326EF3581FA7DD34834ED720F11CDAAEDB59925411AC4243F20E
33C903D2819D4514D5CED1B80071C80734EA620AC8F14B6DF0DDF12BB8797CFE75AF5CFEB569ACEA7A7CD631ADB42B28DA64DC4F1F4A068B7E18057C
39601B19108E95AB58BA1DBEAB636505A5E2DBBAC436974739C7D3149E29B7BFB8D2CFF674EF0BA64B146C122907536E8AA7A48B85D32DC5E3979F60
DEC7B9AB94C4145145001451450014514500145145001451450014514500145145001451450014514500145145001451450014514500145145001451
4500145145001451450014514500145145001451450014514500145145001451450014514500145145001451450014553D57EDA34E9BFB3021BCDBFB
ADF8C67F1A87433A9B698875854177CEED98C11DBA500695145140051451400514514005145140051451400514514005145140051451400514514005
14514005145325916289E47FBA8A58FD05003E8AA9A75F45A8D8C57700611CA095DC30719C55BA0028A2B2F54D5D34F9EDED92233DC5C6E2B1AB0076
A8CB35006A51591E1DD7ADB5FB392E2D9593CB90A32B738AD7A0028A28A0028A28A002A85D68DA65E166BAB18252C72772E726AFD140156CF4FB3B1D
DF64B748B775DA2AD514500145145001451450014514500145145001457371F8BED1B5C83486B2BE8AE26381E6C4100EBCF5E9C55EF116B31E83A61B
E96092750E176C7D79FE9C5006B5155B4FBB5BEB082ED11916640E1587233566800A28A2800A28A2800A28A2800A28A2800A28ACFD5F56B5D1ED3ED1
74C793B51179676F402803428ACBD0F5EB1D76DE496C5D8F94DB5D5860A9AD4A0028A28A0028A28A00A1AB69369ABDA8B7BE8F72AB07520E0AB0EE29
D65A7C368EEE85E495F01A47396C0E83E95768A0028A812EA092E64B64954CD1005D33CA83D2A3D4EFE0D334F9AFAEB77950AEE6DA327F0A00B74552
D2B52B6D5B4E86FECCB186604AEE183C1C1C8FC2AED0014514500154755D39353B7582595D23DE19D57F8C0EC6AF51400D550AA1540000C003B53A8A
AF7973159DABDC5C3158E31962013FCA802C51591A5788F4CD5EEA5B6B09CC9244BB9BE5C00338AD7A0028A2A29E68ADA079E77091C637331E805004
B456268FE28D2B58B96B6B299FCE519D8E84123D6B6E800A28A2800A2B24F8874A170D01BC5DCA705B69DA0FA6EC62AFFDA60DD1012A13367CBC3677
E064E28027A2B3358D6F4ED19623A84FE5F9CDB5001927FF00AD5A08EB246AE872AC0107D45003E8A28A0028A28A0028A28A0028A28A002B3357D3EE
75085A082FDADA2910A3A88836E07AF27A55F9A68E085E699D5234196663C015462D774895B6C7A8DA938CE3CC02802AE83A1CBA1DBC56B0EA0F2DAC
79FDDBC63273EFDB9ADCAABFDA365B03FDAEDF69E84C805471EABA7CD742D61BC8649C8C8446DDEFDA802F5145140051451400514514005145140051
45140051451400514D665452CC42A819249C002ABD8EA167A846CF65711CCA876B1439C1A00B54514500145145001451450014514500145145001451
45001451503DDDBC73A40F346B2BFDD42DC9A009E8A28A0028A28A0028A28A0028A28A0028A28A0028A2A29258A15DD2C891AE7196200CD004B45202
08041041E8453239A2977795223ED386DAC0E0D00494514D6655C6E603270327A9A0075359770C1E9DC7AD3A8A0028A28A0028A28A0028A28A0028A2
8A0028A28A0028A28A0028A28A0028A28A0028A28A0028A28A0028A28A0028A28A0028A28A0028A28A0028A28A0028A28A0028A28A0028A28A0028A2
8A0028A28A0028A28A0028A28A0028A28A0028A28A0028A28A0028A28A00C4F15B6A1068B35D69B79F6692D91A46F9036F0074E41A5F0DEA13DE785A
DB50BB6F326788BB60633827FC29FE2A709E17D4892066D9D467D48C552F07CA63F03594A9866581980F7C938A5D47D08EDEFAFEF34A7D574EBC499C
02CD6A546D18FE1F5CE056BE83717175A25A4F78733C89B9F8C739F4AE421F0ECB77671EB7E19BF96CAEA552EF0939476EE3DABA9F0DEA371AA68F1D
C5E43E54E18A38078247048A10335EB9FF00166BA744B4B7318FDEDCCC2256232147735D05739E30D5EDF48B181A5B24BA9A6936C292282A1BD4D0C1
106BFA86A9E1FB14D49AF52F60122AC91B421386E3208AE960944D047328F96450C3E8466B86F882862F0BC4FA8DC992EA49915150ED4CF53C77E075
FA574DA65CCE6F12CDB6F929671BA903BF4EBFE7A501D0D8AA7AB5CC967A55CDCC5B7CC8E32CBBBA66AE5453BC51C123CE5444AA4B96E98EF9A6238B
D435AD7B4DF0EC37B75776DF6DB920416AB6F92D93EB9F43E95D2585C5EDB68A6E75C7884C886493CB180A00CE3EB5CE688ADE26F10BF882E38B1B36
2968A4F071FC5FD6AD788EF24D67C11ABBD8A1215F626D3CBA2B2963F96EA432C58EA7AAEA5A51D5ADBECF140773470BA925D467AB7A9C55FB4D41F5
AD123BCD265589DFB48BB802382A7FC6B13C25A735E7836C40D4AE16392265654C00064E474FAD6C78734EB5D1ADE7D36CE49A458A4DCC6439E58670
28030F41F11EABAA596AB1CDE4A5CD8336E91532A00EC0773C1AD3B1F1118FC24BACEAD1989829CA918DE7B6076CD667C38C1FEDB7007CD7C4E68F8A
514B2787E068D4B471DC032E3D307AD1D07D6C6CC57BAD3699FDA12456606C3288006DDB719C6ECE338F6ABBA2EA96FACE990DF5AE76483953D54F70
6A48DC1D2124E83ECE1B9FF76B94F85B04D0F87EE1A5CEC7B83E5E7B8C0E7F1CD023B7AABA8DDA5869F3DDC8A596142C40EF56AB375F9ED6DF46BA7B
F52D6A5764A01E76B10A4FE19CD3114EFB55D434DD385D5CDB432798542F944E109FEF7A8F7156F5CD660D16CD66951A4924709144BD6463D8572735
AEBFE148FCC5B81ABE88982F0C8B978D3D47B0FAE3DAA2F1E492DCEABE1CBB8A42968D22BAB91C0624107F2C52B8EC74F36A57FA6A4336B115B081E4
11BBC0CC7CB24F04E474EDF8D55F1DDCEA76BA04F269CB1F97B089DC9F9829E38153EBBA04FACD97D96F35531C00EF6F2E10B9C7A9CD51F144C3FE15
BCCFFBC3BADA300B8F98F2BC9A01177C3125FF00F62E96AB043F6536C9972E770E3D315D0D73763732D8F83F4B788A062B0A7CC32082467F4AE92840
C2BCEF5ABFD4A2F885A3CDFD9A7CE10BA241E70F9C10D939E83FFAD5E895C4788495F897E1D38E0C6E3F46FF001A1823A0B47B7D2B499AE67B08B4D8
933248919523EBC71505AEABA95DDA7DBE0B2845A95DE88D21F3245F51C60565FC52328F0837959DA674F3303F879FEB8AE87C3E43787B4D61D0DA44
7FF1C14011DBEA6FA8E94B79A52C723938F2E438C11D467D6B2348F164DA85FEA9A7C96891DE596EC286CA9DB90493E99C7E7557E1A79DF65D54C99F
2CDE36CE38F7C543E0D8236F1BF8A65232CB36D1F466627F90A00DFF00096BFF00F090E96D72D0793246E51D474CFB56ED416D6D6F6B17976D0C70A6
4B6D8D428C9EA78A9E988CDD63516D3E183CA8D649AE66586256381B8FAFE555A4D5EE2D352B2B0BDB31E65DB10B2C4E7670327A8CFE15078BE0B5BD
B5B5D3A799ADE6B9980B7997F824038FCFA7E358D05DEB3A46A76569E26862BCB633016F7CBD51B040CFFF005E90CE935AD722D2DE0B648CCF7974DB
61814F27DCFA0A825D6AE74E9E15D62D122827904693C4DB9558F40DE9935CDDFB4A7E3169FE6E4C6B0911E7D3CB6FEA6B4FE26B11E11900382D3460
73839CF6A02C6C788751BCD2B4F6BCB5B54B98E205A552FB481ED54749F115CEB5A5A5DE9961BDB67CFE649B103F75048C9FAF4A357F393E1ECDE716
F38588DF9EB9DA334EF014290F8374D0800DF1976F7249A005D37C4F16AD6BFF0012F80BDE2961240CF8F288E32C7D2A7D235C379A95C699796FF67B
DB750CC9BB2194F71585F0EE2437BE21B84006FBE655F61927FAD4B1293F16276EC2C568037754D652C2F2D6CA288CF7976488A2071803AB1F41FE15
1DB6AF3A6A91699A95BAC571323491B46DB9180EA3EB5CB48D792FC5BDB948CC76E444655DC36EDFE1008F527AFAD6A788744B9B9BDB4D5EF7588ADD
34F70CAB1DB1C0CB0E73B89C9E38A00EBE8A28A62389D6807F8A3A2AAF55B6663F9B56B78AF597D174E6B87D3C5D5B701C9718049EEA7A8AC8D5403F
1674703A9B3627F0F33FC2AD7C4C62BE0DB8C779101FCE90CDF6BB917494BAB6B56998C6ACB023004E71C03D38CFE9585A278CA2D65A78ADAC275BA8
9B1E4E739F727181D315D0697FF20AB3FF00AE09FF00A08AE4FE1BC68E35BBD407335FB2F3E839FF00D9A803A1D2F594BEBCB9B19ADE4B5BDB600C90
B907E53D0823A8A6EA7AEDB585EC3611A1B8BE9CE12DD18038F524F005635A8CFC53BE23F86C933FA564A45723E2EDC18638CBF93B97CE2718DA3918
A02C7636BAC17D4974FBDB636B72E85E21BC3AC8075C11E94EBDD53CABD5B0B48C4F7AC9BFCB2E142AFA93FD2B2B57D2750BCD534DD4EE6F2DADE3D3
9CC8422B1DC091919F7029BA98D3AD3C6305DC4B35CEAF343B120460142FF798F6140176CBC41BB59FEC7D46D8DADE326F8F07724A39C907F0A975DD
69F468FED0F613DC5B05CC92C441D9CF71E9EF5CE5EC778FF137447B958F70B676658F384187EE7AF26B775FB98EFF00C1BA9CD0E4235B4B8DC31F77
20FF002A00AD2F8CAD134CB5D492D2E5ECE6203CB8188B271CFA9AE951D648D5D0E558020FA8AE4BC27A4D8EA5E00B0B5BB843C4E0C8C01C7CC18F35
D6A22C71AA20C2A8000F414C18FAF3BF136B84F8834399EC2F6310CEC1A375C07CF1C7A9AF44AE17C7CC06BFE1DDD8C1B83927EA293046BE9D1E9DA2
D8DCDDD8E937D10966CC91888973C6738CFDDA9B4BF1359EADA6CF7D6505D48B6EFB1A311FCE4F1D0679EB5BB5C47C2C18D1B515E38BF7FF00D05680
36F46F1359EB17B71650437705C5B8DD224F18523F5350EA7E2EB1D22F16DF50B6BC855DCA24AD18D8D8EE0E7A564F8787FC5CDF1015208F297247AF
CB51FC4D85269742490655AF3611EA095CD007413F896CA08BED0C939B4DDB0DC04F90138C73E9CF5AD90EA503860548CE73C63D6B9EF1E22A782350
45002AC69803B00EB515F5CCD67F0D45C4676CCBA7C601F4CA81FD6802FC3AF25DB4874FB4B8BB8A172B24A9B42F1FDDC9F9AAD7F6944FA535FDA249
72B8CAA4432CC738C63D73D6B1FC196E24F01DA451486279E17CC8BD431246EFA8FE9573C2BA40D0B4B3A79B91712AC8D23B771BBDBF0A00E63C25E2
19A7D4F58BC9F4ABD769A5507C88B718F031B4F39AEA75DD4E2B2D3BCDBAD36E2E6D9D374A02021071F7813585F0F07FA67889C721AFD803F8B7F8D6
A7886EC5DF81AFEE950A6F81BE52738C1C7F4A06F734B42B8B5BBD1AD6E2C20F22DA44CC71ED0BB467D05412EBB6FF006A92DECA19AF258BFD6792B9
09ED9F5F6ACED22E24B6F86F0DC439F363D3CB263FBC14E0FE7517C3355FF844629B25A59A691E466EA5B763F901408E834CD4ED75381A5B4727636D
74618646F423B564AF8CB4B1AAAE9D731DDDA4ED9C1B88B62F427AE7DAB26C2496CFE29DE5B44C1A1BC87CC700E76903BFBE47EB4CF12DAC575F1274
28A489190C459811D7058F3EBD280B1D09F135826A10DA4C9730FDA18AC32CB1158E423D09AD2BEBD834FB57B9BA902449D49EE7B015CA7C51E343B1
753B5C5FC786F4F95AB63C5BA6DB6A5A23457D7CD691230769876C7B5000DE248601E65FDA5CD9DB310167957E539E99C74ADA1B5D3230CAC3EA08AE
2FC5170F77E02B81670BB5B2448A659FE5660A47207E039ADFF0E4DB744D36DA57267FB1A39073C8C019CD00735E14458FE22F8856350A8074030072
2BBDAE17C301C7C48F11799D768C7D3231FA57754C18579F78EBC4B6773A5AE99672C8649E75493E52BF203CF5EBD057A0D709E29417DF107C3F6215
4840667E3A8CE71F9252608EA6D747B0B5B98EF22B75172B108BCDE412A07715A54514C41597E253703C3D7C6D3779C222576F5F7FD335A959DAEDE9
D3F489EED111BCBC643F4C1201FD0D0064786BC4DA26A56115BC72C304C142C96F2614E71CF5EB4D6D196CFC67617B6CCFF67963947943EE46DB4723
D3350EADE05D17593F6B843DACB280C1E2C6D39E73B6B12DA0D7BC13A95B79F746F74AB895626EA7693D383D0D228E8FC41AA6802EA2B5D6ECE47757
C445E2E09F639ADABFD46CF49B712DE49E5440606109000FA0E2B98F88E46CD15481CDF29FD3FF00AF5D4EAC8AFA4DE2B804181F8233FC26811423F1
4E8D35B24F05E79AB2676AC68CCC71FEC819157F4ED46D753B6FB4594A248F254F041047620F4AE63E18D9C31785E2B8112F9D248FB9F1CE33C537C1
1F27887C4F0A1C469780AA638192D9FE545C2C74DA96AB67A688C5D4B8794E1114659BE828B1D4EDAF599222EB22F54914AB63D79AE4347B93AAFC4D
D45AE1722C62290AB0FBB82013F5E6B766D275093C631EA7F6B02CA38B688C75CF71FD680B1A7A8EAB65A6ED1732E247FB91A8CB37D051A76AB67A96
F16F21DF1FDF8DD4AB2FD41AE274BD464BAF897AA4CF6D2DC2DBC662884401DB82067923DFF3ABF33EADFF0009BDB6A56FA5DC4566E82DA73232AEEE
786233EE31F4A2E163B5A28A298863AABA95750CAC304119045705F12F45B54D0D2FED2DA28A4B7946EF2D02E54F1DABD02B9AF17BB5F69379A559C1
E7DCBC5B981E020F5CFAD0343227B5D4B44B7B6D3E185DE78577B84188863924FAD68E8DE1FD33468916CED904A060CA4659B3D79AC5F8637427F0B2
C271BEDE5646C0FCABB1A4818514566F886F5F4FD06F6F22203C5112A4F63D29884BAD734EB5B836F24FBA553865452C57EB8ABD6F3C575024F03878
DC655877AE4BE1E2349E1492ED406BCB9790B3BF3B9BB67DAA4D0EDF52F0E7853519B502AD346649D1036428C6690EC6F5D6B1A7DACDE4CD74824032
5402C47D71D2AD5B5CC177089ADA559633D194E6B97F8704DCF87DF509FE6B9BA9DDA490F56E6AAE973C9A77C4ABFD32361F67BA8C4DB00C0538CD01
63A88F59D3259D204BD84CB2676A6EC134E8F56B096EBECD1DD46662480B9EA47A7AD71FE2BD3ADAF3C79A2432C6024AAC642BC16C64F26AD7C45822
B5F0F41716E8B14B05CC7E5BA8C151FE45176163A8D4750B4D32DCCF7932C69DB3D4FB0AE5FC21E22D35EC2E659EF3CB796E5C84724E076C57557F1A
4DA7CCB2A2B8F2C9C30CF38AE57E1D410AF84DA59E28CED9A4258A67814750E875B717505B40679E5548B8F989E2ABDB6AD61733F910DD234B8CEC39
527E99EB597AF595EDDDDE9F7F6D7D041676C4492ACD9C303D7B7A7AD6278F2FE39A2D2AE2C165F37ED4A63BA58CA8039E031F7A2E163B2D46E6D2DE
D1DAFDC2C0FF002364120E474E2B3BC2D6DA1DBDACE74171244F26646049F9BF1AD99977C4EB807729183DEB94F8711C91E8D75BB010DDBED03B734C
0E962D42CA69C4115D4324A46762B82714457D673CC6186E627907555604FF00F5EB8BD474A8A7F8936B1440431180C9284F977E3B7F2AB7E3A55D36
0D2EFACA358E786ED554A8C64107834AE163B092448A3692560A8A3249E805528B58D3A684CD15E44D1821721BB9F6EB5078922B8934697EC89E63A1
0E63CE37807245732DABE93E2492CEDC06B1D4629D18C522EDCE0F4CF7F6A1B0B1DED14514C4546D46C925689EEE0574382AD20041A7C9796D14A239
6E61473D15A4009FC2B8FF00895A7DABE9114E88B1DC9B8501D4609CFAD69EA1A45941E10963F294B241BFCC3C9DC075CD2B8EC74808232391552EF5
3B1B270975750C4CDD159B93F8564F852EE56F075BDCCE433C709393DC01597F0F586A96B7BAB5E8135D4972C9BDC64AA803007B7345C2C7650CD15C
46248645910FF129C8A92B8B6BB6D23C7B1595BA916B7CA37C60602BE3A8AED28421ACC154B310001924F6AE76E2D342BBF14DBDD35DAB6A11805630
F9071D2BA4AE2755B78E3F899A4BAC6AA648598903A91BBF5A18D1D5DE6A365618FB65D450EEE46F6C5365D574F8111E5BC8555C65496EA2A878C628
E4F0CDE99141289953E8722A3F0B6956B0786ADA22825F363CBB30C939FAD006DACB1BC4254911A3233BC30231EB9AA475DD2436D3A8401B38C6EAC7
F00CAEDA6DE42DBB64174C881BA818154FC71A44117D8F578608D45ACBBA60171BD723FC28BE970B6B63ACFED0B2DA18DDC00119E64029F6F756F741
8DBCC928538254E706B22E92C752B7B1856DA3297789301402AA064D56F17DCB68FE1990E9EA90066099518C647F9E680B1B8B7F68D70205B988C847
0A187F9CFB54B3DCC16E019E78E207A1770BFCEB964B9D293C34B15BA34B3887E531C0DBF7E3AE48F5EF57F4359359F0F5BB6B7663CEC152245E4F6C
D0163585FD914120BBB7D8780DE62E3F3CD3350D4ADAC2C1EEE5957605257073BBE95C7785B434BAD2F5219C0691A280B0C88FD702B6A5D2574EF0A4
D6F2B09E48A3387DBC8F6145C2C58F0E6A515DE956F24B7B1C934B93B4B8DDC9E07AD3B5F834BBD8638751BB58447206004801CFA62A1F0BE9B66BA2
58CE6D62F382677ED19CE6AAF8FA185B478E668434AB3A2AB6391CE7F2A3A075376EAC60BBB216CFB84600D9B588C6060735CFF832386C6D6F9A6B90
03DC6034920EDE9D2BA98862141FEC8E95C9F84B4EB7B896F2F250260262B1AB8F957D78F5A00EB6391255DD1BABAFAA9C8AE57C5F6D1A5EE9B3AAB6
F96E151B0C79E454934BFD8FE2DB3B6B7045ADF210D1E78561D08A7F8C4959748603205E2E7F31402DCE9AA24B88247D91CD1B38FE15604D6178D354
6D3345CC6DB1E6709BBB01DEA196EEC06881AC833DCC480A148CEFDDF5A2E163A8A2A9693732DE69B0CF7113452B2FCCAC31CD5DA620A28A2800A28A
2800A28A2800A28A2800A28A2800A28A2800A28A2800A28A2800A28A2800A28A2800A28A2800A28A2800A28A2800A28A2800A28A2800A28A2800A28A
2800A28A2800A28A2800A28A2800A28A2800A28A2800A28A2800A28A2800A28A28029DF69D67A82AA5EC0B32AE70AC4E39FF00F5532DF4AB2B482686
CE0581265DAC23E3B11F875ABF450065DA68F0D84021B09668230A176AB641F53CFF0011F5AB9696B0D9C0218176A2E4FD49E4935628A002B2F5AD1E
D75AB78E0BB8C32A3EE072415E3B7BF4AD4A280316E7C35A6DE5A35BDF096E8150AAF2BEE641EC7B55AD334B874C42237965760AA5E56CB600C01F4A
D0A28B05C2A96ABA6C1AAD8BD9DC99044FF7BCB6DA4FE35768A00C48BC31610D8AD9472DDADB2F48C4C40ABFA7E9B69A6D88B3B4882C033F29E739EB
5728A00C3B6D00D8348BA65FCD6B6F23163005565527AEDCF4ABD369E1ED5608EE6E20C36E3244C03B1EF9241CE6AF5140187A57872DF48919AC6EEE
952493CC9119948738FA56ADE5AC37B6B25B5CA078A55DACA7BD4F450064FF006449F62165F6D90DA81B4A951B8AFA6EABF696D0D9DAC56D6EBB6289
76A8F4153D14005676B7A70D5B4D92C5A431A4ACBBC8EA541048AD1A280317FB3AF96C0E9A2E626B531F9424753E604C608F4271DEA7D4345B2D4349
5D3674FDD22858CF74C0C022B4E8A00C8874BBA5B75B59F52925842ED204615987A16AB5A869F6FA8E99369F3A0F2654D8401D3D31F4E2AED140187A
7E852DADBDA5B4DA834F6D68729198C0CE3A64E4E715B9451400573575E179AF35AB5D527D56469ED3FD5010A851F866BA5A280289B26BAD324B3D4D
D2E3CD52AE5576820D50B4D2B51B3D3974D86FA2FB322796923C44C8AB8E07DEC71EB5BB4500665A696BA768EB61A6C9E4141C48503739E491DF3595
A7785AE74ED4AEAFEDF576F3EF1B74F9B71863927A678EA6BA8A280B8514514018BAFE843599EC5DA76856D2432829F7B7718C7E545D585E6A0F1C17
ED6CD691CA24CAEEDEFB79191D07BD6D51405CC9D57458EFEE6DAFA37F2AF6D0930CB8C8E7B11DC555BAD22EF5796D975892D8DADBC82510C20FEF18
0E0927A019ED5D051405CCBD774FB8D4EC1ECE0B9581255292929B89523B7A543A4E9179A568C74F86F91BCB5DB03B45F707BF3CD6D51401CC7867C3
777A03DD01A84734772C6460D0E087C75EBD298340D617589F584D4ED92F268445B05B92800FC735D5514AC3B985A9E83F6DB8B2D4127F2753B45C2C
CAB956E39047A1E7F3A90D85FDF346BAA4B6E608C8631420FEF18720927A73DAB668A62B851451401C949E1AD566F12DBEB72EA76C67B7431A279076
9521B83CFF00B46AE6BFA0DF6BB60F672EA31430BED2CA2DCB1047A1DC2BA1A280326DADB55B4D356DD2E2D279A30A91B3C6C836818E704F3593E18D
0758D023B98BCEB29639E5329FBE086231E9EC2BACA2958773928B45D72D7C4775ADAC9652BDCC6226872C02A8C773F4ABFAD684D7B7F6BAA584AB06
A16BC2B30CABAFF75B15BD4530B9911DB6A776F11D45ADE3891B73450E5BCC23A649ED9E715997FA06A03C630EBFA7CB03FEEFCB922989500631C100
FD6BAAA2811C85FF0086F526D7ADB5DB4BA8DEEE3E1E2932230B8C617DAAF5FE917EFE1F6D26C64B64496131C8F26EC82D9DC463B726BA1A281DCC1F
0CE9FA8E9163069B73F6792DE0421654241249CF4FC6B7A8A28105719ADE81E20D6352B1BB79B4C8BEC32178C2EFCB648EB91ED5D9D140156CBEDDB5
BEDE2DC367E5F24B118F7CD737E1FD1F54F0F4B7D6B0C30CF6D7171E6C729936EDCF5C8EBC01DABAEA280389D1B41D6F47F13EA17E05BDCC178D9663
26D232D927A76E6A6F16E8FAD6B37567F658AD162B39C4C8CF29CBE31C1E38AEC28A561DCE73C456DACEA9A2C9A7C1656C8F70A03C8F3E553041E9B7
9E95369705ECDA5FF666B1A7A471240212EB30712F18E9D456ED14C2E735A0E9FA8E810369C90ADE5A2316825F3023283FC2456BD85AC96FE74B70CA
F3CEFB9CAE700740067B0157A8A04735A36917FA3DE6A496E913DBDE4E668E42F831923B8C73526B5A56A175A4B69360F6D1DBCB0EC9249776E073C9
007AD743452B0EE61F87F4CBAB0D1D34BD41ADE78638FCB564DC0B2F7041FAD54D1AC6F7C3504B610DBCB7967E69781D5D432038F948F4CE79AE9E8A
2C1739ED13447B7D5AF758BDDA6EEECE0228E2341D067D78159F7B63ADDD78B2C7571A7C222B58CA6C3382C41CF23F3AEC68A2C1738DF1AE99ACEBF6
115A59D946822944DBDE71C9008C01F8D4DE26D3753F117864225B1B6BB494482069411201D89E9F9FA575945160B9C76ABA6EB1AFF86CD835B8D3CC
71A8D8CEAFE6B01D38E00F7CD5FD0EDB5186DD65BBB6D93DBDAAC11C46507CCC1C9391D3B0AE8A8A2C17388D1F4DF10D8F89350D526B081D6F48CA89
C65066BB7A292989B16B89B08FFB47E295F5D9F9A3D3ED96243E8C78FEAF5BDAAC9AD457F6D2699047716B82268CB8539EC727B7D293C3FA43E9B1DC
CD732092EEF2532CCC3A03D80FA52193EBA7501A3DC9D280379B3F77D3F4CF7C556F0A1D5CE8A875ECFDB0B1EA003B7B671F8D6DD14C41593E25B19F
53D0AE6CADB1BE6014E4E38C8CD6B5140187A6DBDDE89631DAB07BCB688615D3FD628F423BF3E94B3893589A088DA4B15B432095DA74DA58A9E001FD
6B6E8A00E3FC55A66B3AC4F69F65B18163B1B8F395A49C66503A0C638FC6B5B53BABD3A232FF0066CAF733C6C8628D83042411C9ADAA281DCE3FC182
F748F0F2D8DE69576B2C2CC7E55043E4FD6A0F0A41A959788357B9BCD32E638B509C323100ED196EBCFB8AEDE8A560B9C9DDE9173A678ADB5ED3EDDA
E12E5365C428406CF1C8CF5E95A8B2DD6A57109FB24D6B6F0B8918CE006723A00066B628A2C173919F48B9D1FC56FAE69D035C5BDD214B9823FBC09C
1DC077E83F5AD39D65D624862F22782DA371248D2AEC6723A281FCEB6E8A2C170A28A2988E6F58D6F5648244D2745BB927230AF22E02FBE3BD2E9378
B63A72472697A979CC374A7ECE58B31EA735D1D1401E77E0B5D4344BED456E748BE4B3B973245884B11C9C023E95DEDA5C0BAB649C452C41C6424ABB
587D454F45036EE155353B28F52D3A7B29890932152476AB7450238FF0C79DE1AD364D3350B59CAC52B18A58D0B2C8A4FE95BD14525FC575F6C8DD2D
EE13CB585FA85C10491D89CD69514AC3B9C9F856DDFC396773A65F0758E39DDE097692AE871E9DFF00C6A5D174795FC4B7BAF5DA18CCC3CB8236FBCA
A38C9F4E9D3DEBA7A2985CE1B599A497C67A6EA11D9DDBDAD92B24AE223C139EC7EB478EAE5F57D056D34FB3BB925332B90213C000D773452B05CCCB
4BD4D56CA758A19E1214C67CE8CA72476F5AE6FC30F71A6689268D25ACC2FBCC70A0A1DA41EF9FA576F45160B9C578BA1BB5D43483224D2699015372
A8A4A923D715178C16EF55D3AC67D3ED24FB0DA4E92B7CB86201FE15F402BBAA28B05CC48751BB782E2EA3B2B89239187D9D3033F77BF390322A8783
63D4AC2CE7B6D4B4F921DD2BCBE6021873CE3039AEAA8A623889AEE37F1E5BEA2B15C35AC76ED19944476863478E2F7EDD6F676F610DC4EF0DD24B2E
C85B0A833C9E2BB7A29587731359BB67D3AD6EACA279DBCF464400827AF6AA1E23B4B0D56CF72C4175342BE5154C48AF91C138E40AEAA8A2C1723843
8863590E5C280C7D4E39A928A2988E1BE20DF5BCF69059C1BE6B88EE559D110E540EA6B5B51D734F7D025109695A484A247E5372718C74AE8E8A561D
CE4BC17756F75E1C8B4D90149FCA60F1B2E303A52786618FC34B7BA75E1D88D7066864DA486560076FA575981BB7606718CD2900F519A2C17396874F
9356F172EB051D2D2D536C45C60C8DEA3DABAAA28A621AC700939E3D066B8ABFBD9AE3C5DA7EA31E997C2D2DA365797C93CE41E31F8D76F4940239CF
155E89B43B8B48ADEE1A69E30557CA3C73DFF2A7689AA456FA0DBA4905C2BC318565F2CF5F4AE8A8A5619C6782EEBECA92DB4D6976B2DCDC338CC470
07A935D46A36AB7BA7CF6AE0112C65707F4AB745160B9CC782F4FB9B6B279AF376F27CB8837F0A0EDF9FF2AD0F12E9A754D224863456954EF407B915
AF45160B98BA66A90AD8C50DC0315D46A11A00A4B6471C0F7AD196E961B5F3E647518CED0A58FB0E3BD59A2988E4BC197663827B6B88A78DE4B8778C
346704135BFACC52DC691751400991A32140EB9ABD45219CC787B5341A4595827986F1308EA6327660F24FB6335178DA696E34F4B6B5B69E6649D5A4
2B19C2815D50500920004F538EB4EA2C172859DFA5D58BCF1472208C11895769C819FCAB03C2CEFA57DAA1D47F72933F9B13B7DD606BAEA632ABA957
50CA7A823228039B689B5AF135ADE42AC2D2C4644A471237A0A6F8C27649B4EDB04AEB0CE249595490AB5D3AAAAA8550001D0014EA2C17307C4365FD
B7A1AB5A1DCCA44918C7DEC7B568E997293DAA260A4A8A03C6C3054FD2AED14C41451450014514500145145001451450014514500145145001451450
014514500145145001451450014514500145145001451450014514500145145001451450014514500145145001451450014514500145145001451450
01451450014514500145145001451450014555BFBDB7D3ECE4BABA7D9146324FAFB0F7ACA7D72E6DED5350BBD3CC562D8CB07CC880F42571D3F1EF40
17B57D523D2AD4DC4B6F7332004B1823DDB71EBE94BA2EA916B1A6C77B046F1A392007EBC53F50749748B9753B91ADDC823B82A6B17E1E951E11B6E7
80CFC9FAD2EA3E874F4561C5AE4B7D73711695682E63B6709248650A09EE055DD37538752B433C2183212B24471B9187506988BF4573369E29F3B58B
8D2E7D3E586E21C055F3036F24F1F41EF564F88A3B6D456C754B7364F2465E26670CAF8EA011DE81D8DDA2B9D9BC48D66F0C9A8E9D25B594CDB52E4C
8180CF4C8C6466B7A49A38A169A47558D5771627803D6811251589FDBE9F646BF16939B15CE67046703BEDF4F7ABB757C63B14BAB6824BA47C1022C6
7691D79A00BD45739A678BACB55819ECADAEA49158AF96106723F1E95A5A66AD6FA8B4D12078EE2DC859A190619091FA8F7A02C5ABAB982D22F3AE24
11A02064FA9E054F5C1F8EF57BB8F51D374D1632FD9A4BA42CE304CD823E551F8D75B63A81BC623EC3776E42E733C7B47D3AD005FA2B34EAB16F9B64
52C904390F32005411D47E152E95A95BEAB64B776BBFCA62429618CE0E2802ED1451400C7758E367738550493E82B3343D7EC35D8E57B0773E536D60
EB83572FAE05ADB191ADE69D73829126E38FA7A5657856F749BEB6B89348B316CAB295906C0A49F7C50336649E18A48E392545694E1149E58F5E2A6A
E12F35B7BAF881A7DB4B6170915A248CA1A3CB3165C6E03D3A575F637C2F049B6DE787CB6DA7CE4DB9FA50162E5145472CA9044F2CAC111016663D00
A04494565DA6B96B77199A28AE4418CF9CD110847AE6ADDD5E5BDA5A35CCF2AAC406776739F4C7AD0059A2B25F5CB384A1BB59EDE3908092CC9B5189
E40CFF008D685CCEB6D03CCE1CAA0C908A58FE428026A2B1ACFC4DA35F4B2456D7C8648BEFA302A7F0C8E4FB0AB161AC58DFDC3DBC3232CF18DCD148
851F6FAE0F6A00D1A2B2356F10E9BA3B117D3321181C213D7A7345CF88F4AB655792E731B757552557EA7B5006BD151C52C7344B2C4E1E3719560782
2B3EE75CB4B7B96B65134F2A101D2088BECFAE280352A3791136EF755DC768C9C64FA541637D6D7F1192D64DE14E18104329F420F22B8AF137896C6E
7C41A5E9C12468A2B85966768D8647A018C9A571A47A05154ACB52B6BF0E6D1D9C2100928CBD7EA2AED3115EEAF2D6CD55AEA78E20DD0BB6334DB2D4
2CF50477B2B88E7546DAC50E70699A958DADF5A3C7770A48A14E370FBBC76AE2BE1B5DDA5868775E639DCF72DB5514BB3003D064D21D8F42A2A9E9FA
8DA6A30B49673ACAAADB5B19054FA107914CD4756B2D34C6B7329124A711C6AA599BE8053117E8AA561A9DADFEF581CF991FDF8DC6D65FA8ABB40053
1DD2242F23AAA8EA58E00A7D50D46F34D897ECFA84B0E24FF966FCE7F0A009D6F2D5E5112DCC2D21E8824049FC2AC579E3D9DA5AFC57D2BEC091AC32
5A33E23208276C833FA0AEFE69A38226966758E341966638028024A2B32D75CD3EEAE04104F976194DCA543FD09EB4FBFD674FD3182DF5D2404AEE05
C1C63EB401A14566CFAE6996D6F1CF35E44B14ABBD1B9391EB57A19639E159617578DC65594E41140125159EFAC69F1DEFD91EED04F9C6DE783E84F4
06AD41710DCC65E0912440C5495391907045004D45539753B18A531497512B83820B743EFE956C104020E41E845002D417575059DBB5C5D4AB1449F7
998E00A9EB03C4F71A0CF66FA6EB77A902480372FB48F420D006D41345710A4D048B246E32ACA7208A96B2F4E3A6697A1DB8B79D12C110797233F041
E739FC6A78755D3AE2092682F609238BEFB2B8217EB401768AAB677F677E85ACEE239829C1D8D9C55AA0028A28A002A259A2699A15950CAA32C8186E
03DC5477D74B676AF3B23395FBA8BD5CF602B89F8777936A1ADEBB75791EDB93200C0FF00C9F97F4FD2803BFA28A2800A28A6B10AA5988007249ED40
0EA2AAC7A8D8CAAED1DDC0C10E1889071564900124800752680168AA91EA3632CDE4C77703C9D943824D5BA0028A29A59571B8819F534011C93C513A
24B2C68D21C22B30058FA0F5A9ABCFAD6C60F186B3AC4B7970C16DA658ACD95F063DB924815DF46BB11577336D00658E49FAD003E98AEAD9DAC0ED38
3839C1F4A6CF2C50C7BA6952353C6E66DBFAD60F86346B7D326BE9ADB523782E1C12030213F2279A00E8E8A825BAB782458E6B88A377FBAAEE013F41
53D00145145001454325D5BC4FB249E2473FC2CE01A9A800A2A379628C812488A4F40CC0669E4800924003A93400B4535595BEEB03F4340652C54302
47519E9400EA2AAEA17B6FA7D94B757522A451A924B1FD2B33C1D77717FE1BB6BBBB98CB24A58E4F50371183F95006ED1585E27D0E4D6E1B68E3BF7B
3F2A4DC4A8FBC38E3A8AD9863F2A148F716D8A1727A9C500494514500145145001451450014514500145145001451450014514500145145001451450
01451450014571DE2FD535ED0EC1AFA09AD5A1320503CBE541FAD75568CEF690BC9F7DA352DF5C501627A28A2800A28AC7F125BEAF73628BA25CAC13
87C9278C8A00D8A2A0B45996D214B960F30402460382D8E6A7A0028A28A0028A28A0028A28A0028A28A0028A28A0028A2A18AE219A59638E40CF1101
C0FE134013514C75628C11B6B11C1C6707D6B94F125EF88744B237915D5B5CC6180C343B48FD6803AEA2B26CAE2EDF461793DCC45A48448B88F01491
9F5E6AAE923C4371224D7B710C306ECF96B18DCCBFD295C763A0A28A298828A28A0028A28A0028A28A0028A28A0028A28A0028A28A0028A28A0028A2
8A0028A28A0028A28A0028A28A0028A28A0028A2A18651286214850C541FEF7BD004D451450014514500145145001451450014514500145145001451
4500145145001451450014514500145145001451450014514500145145001451450014514500717F1405C9F0EC22DF3E5FDA17CDF4C76CFE38ABF79A
66AFAA685FD9ED7367145342AACE88C5B181EF8ADDBEB38350B392D6E937C520C30AA9616B7D69125B34F1490A7CA1CA9DFB7B7B52B0EE56B780DA78
424B66B9170D05B488651DF00D63783437FC2B89369258C5363F235D16A3A7CD3697258D8CB1C0922B231752DF2B6738E7DEA9786B45BCD0EC56C24B
B867B54071FBA218927EB4018DF0ED6FBFE1170D6AF6A15A6727CC56273C75C1F4C56B787344B8D1F50BD96EB51FB4497CC6631AC7B4039E4FFE3C29
F6BA2DD693713368D25BADB4EDBDADE60D856F5523FC2B42DED674796E66647BC74DABD7620EC077C67AD00D9CB6930997E296A9363221840CE7A671
4EF1AC69378AFC311BAEF1E7B657D895FF000ABDA5E85AAD9F886EB5696EAD18DDA8591155B0318E9F9549A9E87A95FEB767A89B9B55FB192635D8DC
E7D681F521F89033E119BDA58CFF00E3D55BC7734F1F8010DB120482257C7F708FFF005569788749D535AD35EC7CEB38A372096DAC4E41CD4F1E9975
79A34DA66B7F66789904686DF703803A9CF7E9D281152DED7559BC2F6F656A2CA357B458B73163C15C67007A55BF0A59DC69BA32585DDDC57135BB15
26339DA3B035069F6DADE936B1E9F1456D770C4BB219DA42855474DE3BFE1565E0D434ED34FF0067C715D5E4B2EF97CC7DAA491C91F90E280313E1AA
E34BD45B600CD7CFCFA8E3FF00AF4FD2588F89BADC79E05B4671F821FEB56FC2BA4EA7A2D8CF15C2DB486595A5DAAE4104F6E98AA961A56BB6BE28BD
D6A4B6B57FB5204F2C4D8200C77C7FB22818CF1D61B59F0D7A8BDCFE1F29FE95D43DDB2EA915988F21E2690B67A6081FD6B33C4BA34BA91B3BAB3641
77652F9918738561DC55CB286EA4BA6BCBE8E249BCBF2912372C146727271D738FCA811CB0B9D6B48FB45EE9B6ABA9E8D7123C862CE248CE4EEFC383
5D1F862EF4FBCD1A2974B5290166CC67AC6D9C95355343B3D4F47B3314B10B85965790A2B00D112DC0C93823BD5DD2B4E6B386F248E35826BA95A4D8
30429EDD38A101AF4571DE0AD3F55D3352D52DB529C4C84AC8A436704E7F2E2BB1A62615C6FC3A5261D6646EADA848A47A63FF00D75D5DC1B8109168
2232F6F3490BFA0AC2F0C68F7FA21B8824682486698CC58392D9239EDEB8A0667CED9F8B16CB8E962DF8F5AEA60BAF3AF6EADFCB2BF67DBF313F7B70
CD62EA7A4DE278AAD75DB1459D92130C90B305E3D413F5AD6B0826479AE2E555259CAE514E428030067BD202FD719F142F9ED3C2FE4C4C54DCCA11B1
DD7A91FCABB3ACDD7B48B6D734C92CAE97218651BBA37634C456D4F4D9355F0CFD82D2616AD2C4803A8E00E091C551D57438E6F0843A7CBA908D2D76
97B838C1DA79CFF2ABDA61D56DB4E86D2E6D03CF10118983A842070188CE7F0C56778BB43BCD43C302D2C46FB95984CEA085129392DD71DCD21943C6
7A84979E09B8FB25ACB25B011837328D990187CCA0F3D40E7DEBA7D16743616568CECD32D945212475046339F5C8AC0D734FD5F5FF000DFD86183EC6
2345DC8F8CC8C0741E82B534DFB6C503DF5C5832CCB04502C0AC092149C9FD7F4A03A187E1DB481BE246BD36D04C406CF452719A9B55668FE2B68C55
880F68EAC3FBC30FFD71F951A1DB6A76DE23D4F569B4A9156F982A26F00A81DCFE54CD4EDF58B8F19586B11693234169198CA9750C739E7AFBD032CF
C4E00F83A6CE33E6C78FFBEAAF2D9C10F819EDE28D550D812463A9299CFE759FE37B4D575AD216C2CB4A7662EB2176990018078EBD7B55EBAB8D4468
06CE2D1AE4CCF6FE480648C8076E39218D02297842E9ECBE1C4576C4CAD6F04D2803AF05881CFE550FC2E266D06EEF2525EE27BB6F31C9E5B007F89A
BDE0FB7B883415D1752D3E588C28C8ECC014901272011EC69BE1ED366F0B25DD90B796E2CE5B832C0F08DCC320021867B607340142D677B3F8AF756B
16443796C1E451D37019DDFA1FCE9DAEEC93E27E811B64910BB63B74723F956B693A4CC75CBBD72FE3549E6511C51E726341EA7D6AB6BDA55E378AB4
BD76D61F3D2D54C72C4A406C1C8C8CF5FBDFA5033A15BB89B506B21BBCD48C4A78E30491FD2ACD66595BCEFA8DCDFDC45E5F988B1C6858121464E4E3
8CE4FAD69D3246B00CA54F43C570BF0B2CE38F4ABDBB03F7925CB47D3A05C1FEBFA577123148D9C2962A09DA3A9F6AE4FC29677BE1DD3520B8B79A68
EE7F7EC634CB44E40CA95FCA90D6C50FB449A6FC5A3042A043A8423CC51DCEDCEEFAE57F5A8348BE92F7E286AB2CF0C927D8E1912145C12BB5957232
7B824FE35BDA7E9135D78AA6D7EF61F2709E5411B7DEC74DC7D2A3BCD1AE74DF15AEBFA5DB79E93C662BB84300C7241DC33EE178F6A00A8D1EB5278E
ADB548B4E7B6B364104ECEC32C39E48CFAE315DBD61CA973AC4D6EAF6D35AD9C32095FCDC06908E8B807815B94C4324711C6CED9C2824E2B88F87976
759BAD6752BA1E64D24C154B73B5083851ED8C5772402082320F6AE2FC2F60FE17D5B53B39A293ECB752892DA4542C08E7E5E3BF22802BAE9CB65F15
EDDA250B14B6AD22AF65E0838FC7F9D3BE225F3ADEE8FA6FCE60B99B74A88325C02303F5A5BE92FAD3C7306AF79653B588B7648CC6858A0F71D8FB56
878934897C4165A7EA5A70F2AF2D241342265DA48C8383E9D01A451178B5F50BDD3618B47D32EBED504EB246E50204DBE9CD4BE3756B8F01DDBDC461
65F2A36208E51B72E455F5D52EEEE110C3A7DD4374EB82658CAC719F5DDDFF000ACFF174730F0BCBA3DAD9DE5DCD242A8AF1C5B97823A9EDD28122C7
8574AB587C2764822573716A8F2171924B2838CFA0CD54F86AE5BC30C858B08AE6445C9E8063FC6AFE837125B7866CA2B8B5BB4B886DD6368FC96272
A31E98ACEF87D0DD69BA39B0D42D27827795E51B90E08E3BF6A00CC5D6E3F0EDE5E697E20B369ACE799DC5DAAEE0F924FCC3B9CF1ED572FA73E1FF00
86C65D1AE964553FBA9D79F95A4E3F100E3F0AB1A3335BE91359F886D716C6573119632C0A9278231C1EB8A4D33C3A93F85B54D259248EDAE6E246B6
126728B85287E9919A00A5E1F9EC5FC0E167B69AE2EAE227326206767639C1CE3E95ABE02B8BF93C3A90EA76D3C32DBB18D4CCA54B2F51D7EB8AAFE1
8BE6D0F468F4AD5EDE786E2D772A95899D6519246D201F5ADED2A19F7DCDDDC828D72C19633D5140C007DE8034AB9AF8811452783AFCCAAA59101424
720EE1D2BA5AE63C6A6FAEF47B9D32C34E9AE1EE100F3171B473FF00D6A608BDE1EB757F0AE9D05CC68EA6D903230C83C0AE43C09A4DB3F897C41B97
30DADC18E38BF83966C71EC1715D66893CD6DE1D87ED965710BDB46A8D1EDDC5B1DC01D6B03C18F25A6B5AE49736579125EDC79B0B340DC8CB1EC3DE
900B022D97C57305A811C53D9EE78D38527D71F857755C1DC5C797F1117547B3BDFB2476DE4F9DF676C6EE7DB35DCA3AC91ABAF2AC011918E2840C7D
145727ADF8A25489ADB4AD3EFA5B92C519BC92047CF26988B736A56D26B3BA7F37C9B504464465959FB9E076AC1F8792ACDE23F134899DAF70186460
F2CFDABACB5BFB78B4B8E5115C2A28DBB7C86DD91D78C5723E1F92E747F136B57575A55FADBEA0FE644CB0138C12707D3EF5219E834567EA97373068
B7175650B3DC2C45E38C8E738E98AC8F036B97BAD699336A716CB8824D8C76EDDD9F6A62B1D3D65F88A29E6D0AF23B55DD3797955F5C1CE3F4AD4AC1
F193DC2787E4FB196F3CC918455382C770E2804734753D27C4963158C69F61D522641E538D9BB046573DC77E6AFF00C45D4AE60B5B0D32CD8A4BA84D
B0B0383B41031F8961478862D2F57D19DDA248B54DAACA146D9164F4F7A3C5FA35EDEE91A5DF47BE5BED399247551CB8F97771EA0807F3A4517FC43A
3D8DB7862E3ECD0471496B0978A451865239CE6A6F056AB36B1E1BB6BBB8C998651D8FF111DEA2D5F57B4D43C3D3C764C679EEA231A40A3E70CC3BAF
518ABBE17D27FB1742B7B266DCEA3739FF0068F5A3A8BA1AF593E22D3EC6F74D964BE8838B78DDD0962369C75E3E95AD5CF78EEE45A7843507CFCCF1
F96BF5620531185F0DB42B26D121D527B60D78D2B32CA49C800F15DF564F85AD3EC3E1AD3EDCF558413F8F3FD6B5A920643736D05DC262B989258CF3
B58645719F0EA0482FB5F8E32424777B1573DB9AEC2F2EE2B2B66B89CB045EA554B1FC85719E059CDBEABAB2CF6F7118BDB932C0CD11C32F3DF1EF47
51AD87DFEA3A758F8A6EEDFC490011DC05FB35C30CAA28EDC7DD3EF5AFE12B316D1DF491DE1BB827B8DD0CC5F7975DA31CFB74FC2AADAEDD4356D574
FD66D04D6A67FF00479265C83C00541EDD2A4F0A68E746D4B5682DE466B0DE9E4A939DAD825867DB205086751591E28D53FB1F40BBBD5FBE8B84FF00
78F02B5EB27C4FA5FF006CE817562A70F22E539FE21C8A6498BA0E9F13F8346A13DBFDBAF2E2069FE6F999988C851E9D00E2974AD4F51D33C0525F6A
B13ADC5BA37968E30C57F8739FAFE42A6F08EA305AF87E0B1D41D6DEEEC90C72C6FC1F949C11EBC7A55BBFB59F5BD0B5184E425D2158118018007073
EE79E690CCDF0A6976DAAF87A3BFD5635BAB8BDCBB3C9C90338C0F4A8FC2B7C6E351D63C377AC678ED1CAC5E61C931E718FC38A9FC197D0D9F8762B0
D4245B7B9B32C8F1C8769E0E7F1EB4CF0B694D06AFAB7886E9822DDB318D4F554CE493F90A0664C9A67D93E24DBD8E953496B6EF6E249155CF039CE3
3D09C55FF152AE87A8E9177A66E85E7BA114A01243A9F5CD575D46CA7F8A315CC5728D0FD8B67999C286E78CD58F8813C3F69D1A2320DF1DE2C8E00C
955F53401ADE35B1B6BCF0CDEBDCC619A0859E3393F29C75FD2A9782347B06F0CE9578F6CBE784DE1F2783B89CD6D6B917F687872F63B5657335BB04
20F078AC2D0A44BAF03D9E9E8E05CEC588C59F98156E78F4C0CE68174297C45D3ED626D36F638825C4B7888EEA48DC3FC8AE9BC4BA7C17DA45C194BA
BC313BC6C8E5482067B7D05737F11AEE39A3B0B683334B05DAC92AA0276003BFE75D36A5A8DA3E877330B88D52485C29638C9DBD31F8D0072FE00D3E
5D4BC3625D4AE2596032B88A2DE71E9B89EB9CE6AEF82AE668F54D6F4991DDE2B3B8FDC963921493C67F2A6FC39BEB48FC270452DCC48EB23E55DC02
32C6A2F075C43278B3C48CB2AB7993831723E61F3671FA521BEA43AD599D2FC63A2AD8DDDD46979291247E6923F5A5F8936D258E98BA9D95D5D45334
E11C2CC769041EDF855AF182987C4BE1EBE75C411CE55DF3D33507C4FBD8BFB023B745F318CCAE4AF21319EBF9D305D0E8F59D4C68BE1C96FE4CBB43
10C027EF31C019FC4D65F8774FFED5D162D4F5092492EEE54B87DC4797E800AB7E26B1FF008487C2B3DB58CA8C640AC841C82410714784EF62FF0084
5ED04F2A2496F1F972A9382BB78E47D3140BA1078675893531A86937AEC6EACD8C6F20E0B29E01FAD67E96D776DF11AE74C4BD9E4B44B6F34A48D9E7
8FFE2AA6F0458DC7F686ADACDC02AB7D31F2BFDA504F35142A61F8B570EE30B3D8E109EE46DE9F950323D7BED3A5F8BF418EDEFEEBC8BB9C89226909
07919FE75B1E30BFFB069F03CBE60B579D56E1E2FBC13DAB03C5F7F049E36F0EC7BC6DB79B73BFF0F257FC2B7FC4F7F15A369E2E5165B3B894A4CA57
70208E0FE14019A6D8CBA868F2E877EF369667DD2C6929700804E493CE3D8D7695C1DF78793C3FA9DA6ABE1F95D1259D565B60D95753E9ED5B7E2BF1
441E1A8ED9A685A669D880AA71803193FAD084F53A1A2A0B5B84BAB486E63FB9322B8CFA119A9E988E53C437F7DFF092699A54104ED67700B4ED09C1
2391D47403A9F5AADA9B3F873C41A3C361712082F6531CB6EEE5C1E9F30CF4EBFA559F126AB7B1F8874CD1EDE436B05D9CC9700727AFCA0F63C75F7A
CCF10C5A4D8EBFA246B08DE9700CB363249C8FBCDF99A96523B1D4AD4DE593C4B7335B36322489B691C7F2AE43C12FAAEAD6572B77A9DD18619CA870
DF3BFD1BB0AEA6DAF126B7BC91E553189991589E08DA3A5607C350BFD8D765581CDD3F00F419A7D45D0B5A1DFDD45E29D4F43B99E49E38116685E4E5
8038E09EFD4549AE6AD336B767A0D849E5CF700BCB28EB1A0F4F7AABA6CA8FF1335555032B66809F5FBB54E7B39AD7E2ADB5E4BFEA6E61608DEFB718
A008BE21DBDFDAF867CB3706E6D4C8BB9A4C6F53CF7EF5DC59FF00C7941FF5CD7F95723F14AE922F0DA404FEF2599481EC3AFF003AEB6C8E6C6DC839
FDDAFF002A3A83D8B1595E23BB96C74979E12570CA1DC0C94527935AB591E25BD8AC34BF36E22496079163951BBA93CD3624674CF3182CEE746D4A4B
985AE104ABBF79209E7E9F4A8BC7AF7961A50BEB0D427B7612AAB22B70C0FA7A551BCF0C47A395D6FC393C8850AB983394914F5AB3F12E551E154DFF
002B34E876F7EF9A92BA9D2DCDB497562B125CCB049807CC8CF39C572DE109F59D56DEF62BBBE944314C5166C0F31B1E9E82BB14911A1591594A6DCE
E078C7D6B9BF0110DA5DDB646E377212053EA227D22F2EADF5DBAD1EF276B808A248A461F360F626A3D7F559CEB561A1D94FE4C974499A45E59140CF
1E9F5A8AD0EFF88D7FB48212D533CF7C552D52DE4B2F897A76A53716D3A1883760DB48FE64503356FAF8687A869D671CB71706F25DACB21DE40F5CFE
755B5BF104EBAFDB68362A639A761BE66E30B8ED9AE86FA4860459E68D5D94E23CAE4EE3E87B5713E27D357FE13AD3EF6F91CD94B84DDB880A71EA3A
50C48DDD4EF24F0F4D672497324D693C9E5CBE7B6E653D8835D18208041C83DEB0B50D0746117DA2E6C8CCD1FDD0F2BB64F6EA6B7140550AA0003800
76A100EAC9D6F561A74714310DF7770764283D7D4FB56B57117BA749AB7C421E7CD2470D9C41A3087072450C11A367777F6BAFC16173762E8CF19791
428FDD1FC3A574D542CED2CECE774B68B12B0CC8FD49FA9F5357E8426636AFE24D2F47252F27FDEEDDDE5A02CC3EA074AAFE0F996E74E96E7ECD3C0F
34A5DBCD046ECF422B8D58B5A3F6C115EDA19AF2EBC8759533201D8EEF4C57A4E9F14D0D84115CBABCC88159947048A10FA16AB37C436F1DCE857914
8514189BE67E80D5E9E68ADE179A77548D065998F005508D57545F3AE14FD94FFAB8C9E187F79BFA5311CF7800C97BA446D3DD7991D9B98E38978031
D09F5F6AED2B85F05AAE99E23D5B48DC08DDBD403C753DBF1AEEA921BDC28ACA975FD3A1D5D34B9262B72FD38E33E99AD5A620AE53C4B7FAD6973C52
C37510B3964DA498816407FC2BABACFD6ECD6FF4AB881B1CAE413D88A4C68C3F106A1AB693696D751DF452C3232ABB1840233DC7B55DBD7D5963B516
77C923DC11F7A35C018CE47AD6769F6FFF00090787EE2D2E148F254C519F7F5FD0559F0B89EE5217B8271608605057196E84FE540CB97F7B3E89A3F9
97338B9B82D852CBB739F61D714972751B3B1FB7ADD99595433C2CA3691ED8031557C71632DDE949247232085B2C17A9048AB76DA358BD8A99249E68
993244B316038A046ADACE9756B1DC47F76450C2A6AA9A7A451D8C29044628946114F5033C7E7D6ADD31051451400514514005145140051451400514
5140051451400514514005145140114FBFC86110CB9181CE28B7884302463F8473F5EF52D14005145140051451400514514005145140051451400514
514005145140051451400514514005145140051451400514514005145140051451400514514005145140051451400514514005145140051451400514
5140051451400514514005145567BDB446DAF730AB038C190039F4A00B351C9BFCB6F2F6EFC1DBBBA67B66A4A28029E9F6AD6D1379AC1E691B7CAE06
32DEDED5728A2800A28A2800A2992491C485E475451DD8E0511C892A0789D5D0F4653906801F45151BC91C6543BAA9638504E327DA8024A28A2800A2
8A2800A28A2800A28A2800A28A2800A28A2800A28A2800A28A2800A28A2800A28A2800A28A2800A28A2800A28A2800A28A2800A28A2800A28A2800A2
8A2800A28A2800A28A2800A28A2800A28A28020BB99ADED64992179990644683258FA0AADA4DBBC5685E750B3CEC6497EA7B7E55A1450014D2A18618
02339E69D450033CB4DFBF62EFFEF639A7D145001451450045333AC2ED128670A4AA938C9EC2B92B817BE2786C2C6FAC5EDD527F36F0152170B9C283
EF9AECA8A431000A000000380052D14531051451400535542A85500003000ED4EA2800A28A28022782191C3C9123B0E0165048A968A2802192D6DE47
DF25BC4EFF00DE6404D4A4060410083C106968A00896DE05C6D86318391851C50D0C4EFBDE2466C6371504E3D2A5A2801A0050028000E001DA991DBC
113978A18D18F52AA01352D1400C0881CB8550C7A9C7269B2410CB1F972451BA8E8ACA081F854B450056FB0DA73FE8B073D7F76284B1B38DC3C76B02
32F20AC60115668A0086E6DA0BA88C57312C919E76B0CD4274BB16B436CD691184F542B9CD5CA28021B6B682D61115B44B1460E76A8C544DA7593CC6
66B588C84E4B15EA6ADD1400D00280140007000ED55EEF4FB3BC656BAB78E564FBACC391F8D5AA28028DCE91A75D2C6B71670B88BEE657EED4C2CAD4
45147E446521FB80AE76D58A2802A41A7DA40C1E281430E9924E3E99E9595E2B8ECAF2D22B0B886296E2E5F65B875C942782DED8AE8298514B872016
5CE0E3919A006C112C36F1C29F723508BF403152D1450055BBB0B5BE54175024A50EE42DD54FB1EA2993E9763716C2DA6B589E15390A477F5ABB4500
546D36C9AD16D4DB47E429C84C719A4B5D32C6CCB7D96DA28B72ED60A3823E95728A00CA6F0F690F22C86C220EA7218641CFD41AB979676F7A8AB3A6
4A9CA30E0A9F506ACD14019979A169F7F686DEF22338273B9DB2D9FAD4BA76996FA6A3ADB19087C677B96E9D2AF51400543716F0DCA04B88D6450770
0C38CD4D45005682CEDEDD0A471FCBBB70527214FB67A536F34EB2BFC7DB2D629B6F4DEB9C55BA280292E9766B6A2D5200B003BB6292067F3A659E8D
A7D8861670793BBAEC7619FD6B428A00CBFEC1D385D2DD240527539F315CEE3F539ABD716D0DD45E55C44B22673861DFD6A6A280295B697676D26F8E
01BBB1625B1F4CF4A9EEAD60BC80C373189233D8D4D450053B4D36DAD4AB468CCEA301DD8B11F9D5CA28A002AA5CD841712898EF8E61C7991B6D6C7A
55BA280218218EDE311C4B81DFD49F5353514500733AA689A7EABAFA064647893CC95E33B493DBF1ADFB3B58ECAD63B78776C8C606E3935360649C0C
9E334B401977DA0D8DFEFF00B579D22B9C9532102A4B4D22D6CA3D9019C2E00E677E00FC6B428A00C56F0C6926F8DEAC322DC139DEB3303FCEB4AD6D
92D6111C6CE541272EE58FEB5628A00C5BCD06CEE75CB5D44C604B0659883D7D38FF003D2B6A9A1543160A016EA40EB4EA00290804104020F506968A
0086DADA1B48BCAB7411A649C0F53496D6D15AC7E5C09B1492D8C93C9A9E8A0064889221491432B0C107BD476F6D0C1198E24DA87AA9248FD6A7A280
0A28A2800A28A2800A28A2800A28A2800A28A2800A28A2800A28A2800A28A2800A28A2800A28A2800A28A2800A28A2800A28A2800A28A2800A28A280
0A28A2800A28A2800A28A2800A28A2800A28A2800A28A2800A28A2800A28A2800A28A2800A28A2800A28A2802B5F5EC1616C67B962B183824293FCAA
9C5E20D326B15BC8A72F1392176A924E3DAAD6AA71A4DE1C671039C7FC04D739F0D6D561F09C321F99A577E4FA6718FD2901BDFDAF6074FF00B7ADCA
35BF4DC3B9F4C75CFB52D86AF63A8BBC76D366541978D810CA3DC5725E11B2DBE2DD6FCC60D1C13131A0270AC4F5C7AD58BB668FE2B5884E04B66C1F
1DF018F3F90A2E3B1D4DF5FDB5846AF75284DC70A304963EC054569ABD85E4FE44370A66C6EF2D8156C7D0F5ACBD66C593C4361AC497F1416F6CA51A
393AB64FF0D657882F3CEF16787E48ADE68D0CD8F38A63CC1C70075C5170B1D8DDDEDAD92AB5D4E910638058E3355135FD2242447A840C475C367156
BFD1AFA392365495118A38619C30EBFCEB8BF0615D23C49AB787A78D172E6580E39653DBE98A623A66F12E8A93080EA111918E00505B27F015AF581A
0D9DBCD2DCEACD04625B994F96760F9517818FC89CD6FD0015C578ABC4168FAA58E8F15D147374BF686C602007A64F5FC2BB5AE0FC6022FF0084DBC3
51AA26E33167F97EF02C3AFE46931A3B3B5BDB5BCDFF00659D25D98DDB4E71FE71505C6AFA7DB48E935C8564FBF85276FD481C54C1E08E56B484AC72
94F336AAF404E33E9D6B07414FEC4B6B8B0D42EA2BBB9925690476E85D8A9F51EB4C474704F15CC2B2C122491B0C8653906B81F1E693671EADA135BC
2227B9BC09214E376597FC6AE7806EC8975F8E4260B686F498D1F0A23DCCDC7B7402A7F1C2A9D4BC3A59B18BE0724E0763FD293D86B73A86BCB652C1
AE21050E181907CBF5A963912540F13ABA1E8CA720D72EEBA4F89B5E9610913AE9CC1A52061A562081C8E7031567C5D75368FE1699F4C0B13AED442B
C6C04F5A2E163625BFB48582C9731A9DDB7EF743EFE9F8D590410083907BD72BA035B4BE0B8E231CD399606326D899BCC63D4E718CD5BF05BDFB786E
04D4E395278999079A0862A0F19CFB71F85170B1D0535DD514B3B0551D4938029D50DD5B43776EF6F7318922718653DE988E4B58B8B4D5BC636DA3DC
B29B3B781A693E701598F406B7B40D160D074EFB15B16652ED23337727FF00D42B94F0BE9F6379E34D7E5FB2C661B47582342B90B8E0FF00E826BD02
921B0AE77C4163A56B52DB4573A8A412DACDB955655049E38C57455C17C43B7B78AF745B94891267BD50CE060B0E0F268608EE249A28B02595133D37
3019A8E4BDB48E159A4B985636FBAED2000FD0D50F135AC371A05F34D123BC56D2B46CC3953B4F23D3B573DF0FB48B7B9F0B4171A820B81233F94927
2224DC4600FA8273EF4C0ED91D6450E8C194F420E41A7D71DE0091D4EB3645D8C36978C912939DABCF1FA5763409854125CDBC52AC524F1248DF7519
C027E82A7AE22E356B0D3BC437F6BE2484849E40D6D70E9B976607CB91C8A4C68ED37A1DD861F2F079E9F5A50430C8208F6AE7BC37A6431E9B7D133C
7736B7572EF1B06DC1A3206075EDCD73DE0CB2B78BC5BE22D3D109B44F90231C8C6EA02C7A1646EDB919C6714A480324E00AF3EB2B48ACBE2C7D9ED9
4C70ADB170818E395F7AB9E31D5646F11E91E1F8DD922BB756B82A70594B602FE8680B1D7C3756F3B32C371148CBD4238247E55C77C4AB43068ADAA5
B5C5C453AC88A424842B03C74AB3E2582F6CF50D27FB0F4EC133012CD12F451D41FC33C9A67C527D9E0F718CEF9D17E9D4FF004A011D2690E068F61B
DBE66B78FEF1E49DA2AF570B71A3C2BF0F9AF6E773DF8B11389892191828200F61D2BA0F08DDCB7DE16D3EE6E1F7CAF17CCD8EB8247F4A606D514547
3F99E449E4E3CCDA7667D71C5022392F2D6297CA96E6157EBB59C03F9558AE3F4AD2E18340993C510C0B752BBEF76219D8678208E7F2A77C36BC96EF
C34FE74CF2986E1E352C724280081FAD2B8EC75D45145310515C178F6C65D3608F58B0B8BA8D45C0375124CC1581EFEDCFF3A83C77E5DB5B69BAA585
C5D2453B8694248C43C670771E7DFF005A43B1E8945727A8E936D7B71A5C16973731093F7B21599B2D1819E79E3248A7F8962D4621A55969B6D35C5A
34E16E82B9CEC18C02DD8727F2A6163A7465750C8C181EE0E69D5C478ADE2D02FB489B4A65B6924B9092C28C42C887FBCA3AFD6BB0BBB64BB81A190B
853DD1CA9FCC50227A424282588007249ED5E7FE1B8CC5E31D5746D56E6EA665F9ED7CC9DB94EBD8F5C11F91AD083C3905DEAF7ACF7B7AD651304111
98E0BF56E4F6148763B1A29A88B1A2A20C2A8C01ED4EA6230FC4FAABE99611ADB3AADE5CCA91401B07927D3D2B69812A4024123A8ED5E7BF10B4BB58
F55D22F97CDF3EE2F151CF98718C8E9E9F8576565A6C36376EEB777523CAB811CD3971C63900D0072464D5EC7E20D9E94756BAB8B5917CED8E79C739
071DB8AEFEBCEFC44B76FF001434E8ECE411C925B005C7555CB64FE55A1AFCD3787359D267B5B9B8786F27F2678A47DC1BDFEBCD219DA515CBF8E351
D46C34D80E990CCED2CE1246881DCABD703D33EB547C5F15BF8774E8B54D3E6782EE29142A34C4F9A33C8209E78A6163B6AC0D3B549EFBC53A9D9A48
AD69651C63000FBEC33D7AF635A4506A5A5A87796217112926262ACB900F07B5705E07D1609358F10426E6E94413AC6AC929566196E491D7A52047A5
515059DBADA5AC70249248B18C0691B731FA9A9E98828A28A0028A28A0028A28A0028A28A0028A28A0028A28A0028A28A0028A28A0028A28A0028A28
A0028A28A0028A28A0028A28A0028A28A0028A28A0028A28A0028A28A0028A28A0028A28A0028A28A0028A28A0028A28A0028A28A0028A28A0028A28
A0028A28A0028A28A0028A28A0028A28A0028A28A0028A28A0028A28A0028A28A0028A28A0028A28A0028A28A0028A28A0028A28A0028A28A0028A28
A0028A28A0028A28A0028A28A0028A28A0028A28A0028A28A0028A28A0028A28A0028A28A0028A28A0028A28A0028A28A0028A28A0028A28A0028A28
A0028A28A0028A28A0028A28A0028A28A0028A28A0028A28A0028A28A0028A28A0028A28A0028A28A0028A28A0028A28A00A1AB35D7D8258ED2D8CD2
C88CA3E6550B91D4E4D63F85EDF50D0FC3E967796323BC3B9B313AB6E25B38033EF5D3D1401C47852D751B1D6F53BABBD3E644BE937A003EE6493839
FAD1787546F1A5BEA9169172D6F0C4D19C81939C8CF5AEDE8A561DCE335EB6BF93C51A46A925A4B2584518F322037796C739247B647E551F886DF58B
AD6B4BD5EDEC9DEDAD25FF005007EF3AF2C4576F45160B99FA444E904B2C91B46D71334BB1BAAE7B1FCAB98F1AE917336B9A45FE9F27953C927D9A46
538254F3CFA8C6E1F95745E22D326D5B4B6B5B6B936F2170C1C7B76A5B7B091A6B396E483F64882A0CEE2588C1627F01408BF6F0476D6F1C112858E3
50AA2A5A28A6015C8F8AB4CBD7F10E8FAC5B40D3C366F89523197C67A81DEBAEA28032E08E69EFA5BE58DE2CDB88A3120C12724E48EDCD739E108754
D32DB5017DA74CF7D34E4AC8C38933EADE95DBD140EE707E1AD2354B53AB5A6A30481EEAE525138195386C939FA558F1C413DF5E69891D85CCD05A5C
0966744E36F1D3D6BB4A29582E713AC5A5C69BE23B3D6344B09E4DE9B2EA18E120153DFEBEDED5BF7F6F0788347B8B39E19E159171FBD42841EC7F3A
D7A28B05CC4F0DA5CD969D06997D0BACB02ED122AE51C0E8723807DAB5204923F30CB26E2F2165F455EC3F2153D14C414C91C471BC8DD14127F0A7D6
25DB6AD72BAA5A0B7458D82A5B4B9EA187CC4FD280287C3FB478B479EFA75C4DA85CBCEC4F5209E3FA9FC6BAAA82CEDA3B3B386DA218489028FC2A7A
002B8BF1B5B5E6A1A96931DAD85C4F1DACFE6CAEA9F2E38E07AF7AED28A01195AECE8340BA2564FDFC0C8AA1096CB29C0C76ACAF004861F07DBC3708
F1496DBF7AB210402CCDD3E86BAAA280388F053C96BA9EB5F6CB79E01797864837C4C03024E39C7B8AEDE8A280615C95B7D92E6EB56B4D6ADE3368D7
8C20698641247CD83DABADA63C6922ED911587A30CD00739E12D3174A9F53B5B5777B0F355ADF2720657E600FE554B4545D3BC71AE1BB0B12DCED78A
466C023B8FAD7620050028000E001DA9D48773CFA7BA5B4F8AE6E2ED4C30FD9762B919078EBC55AF1869329D734BF11DA21992D5D3CF542325036411
F99AED76A960C546E1DF1CD3A80B99377A945756CD6FA65CA3DDCA87CBD87EEFB9F4AE73E25F3E154B30CD2DC2BC6E70A7240C8DD5DB2A22962A8016
EA40EB4E2010411906988E4754D634FB8F02DCC704DB9DECCC4B1853BB76DC631563C017313F85ECED41659A0421D1D707EF1E7E95D1082100011460
0390028E0D3D5554614003D00A063AA0BA6912DA5681774AA8C517D5B1C0A9E8A0470BE0DBAF3AC6EEE6F3F79AD492C8047264B01D8007A2FF004ACF
F0943716F11D32E232B769AB899E351F7542F278EDC57A46D5DDBB03774CE39A40881CB8550C782D8E4D2B0EE3E8A28A62296AB6316A7A6DCD8CFF00
EAE74287DBD0FE75C9F85B4FBAD4FC0D3D8EA19CED92D61CF50149E79E9CF1FF000115DCD4688A8308A00C9381EA4E49FCE80394F00437C74BF3F512
BBE21F66894754553CE4FD6A3F1DEAD796375A5D9C52496D6975301717287042E40201EDC126BB1550A30A0019278A64D0433A859E249141C80EA08C
FE340CE0FC75040B61A67F67402586DEE95A59A3C305FAB7AD763A6DD3DD4D76C306012811306CE46D19FD7F9D5A5B6B7588C2B044233C9408307F0A
A7AC1BBB6D16E4E9110FB4A266145504673E9484729E3D8EE34CD6F48F1158C65A48DFC8900FE207903F105C7E55D969B6CD6B669148C1A524BC8DEA
C4E4FF003AC6B196E356B7D3E1BE4579E1659EE98260238E540FE55D2508028A28A60715F1251FECDA4DC007CB86F55A4603EE8AD46686F3C5D6735B
4824582D9CBB21CA8DDD327F3ADD9A18AE2168A7459237186561906996B696D671F976B047127A22E280384BED46CDBE2CE9D225C44634B531BBEE18
0C43F19FC47E7527C42BCB63A8680A278CECBAF31B07381C735D79D1B4A672EDA659963D49814FF4A1B46D29860E9967FF007E17FC290CC5F1BEBF73
A3E990358ED06E9F67DA3EF08978F980EFD6B9CF1F5A69763E1F8E3817ED377248A5EE198BBEDEA493DB35E833E9D657165F629ADA36B6C63CBDBC0A
8E3D274F8AD64B64B38BC99461D48CEE1EF9A00CED06E7ED172AB0CE65B78ACE25C2B654363F9D60781E68AD35CF12FDB1BCA904C1D83F1F2E5ABB4B
1D3ED34E84C5636E90464E4841D4D473691A74F75F699AC6079BBB94193F5F5A00CDD2EEE3D2FC1E2F2E4485215790E7EF11BC91FCC54DE19F115AF8
8EC9EE6D91E331B6D746E48F4AD0BF7B5834F99EF157ECC8877A9191B7D31597E12B086CB4D924B7B7104775334C898E421FBA3FCFAD006FD14514C4
145145001451450014514500145145001451450014514500145145001451450014514500145145001451450014514500145145001451450014514500
145145001451450014514500145145001451450014514500145145001451450014514500145145001451450014514500145145001451450014514500
145145001451450014514500145145001451450014514500145145001451450014514500145145001451450014514500145145001451450014514500
145145001451450014514500145145001451450014514500145145001451450014514500145145001451450014514500145145001451450014514500
145145001451450014514500145145001451450014514500145145001451450014514500145145001451450014514500145145001451450014514500
145145001451450014514500145145001451450014514500145145001451450014514500145145001451450014514500145145001451450014514500
14514500145145001451450014514500145145001451450014514500145145003248D24429222B29EAAC320D38000000600E82968A0028A28A0028A2
8A0028A28A0028A28A0028A28A0028A28A0028A28A0028A28A0028A28A0028A28A0028A28A0028A28A0028A28A0028A28A0028A28A0028A28A0028A2
8A0028A28A0028A28A0028A28A0028A28A0028A28A0028A28A0028A28A0028A28A0028A28A0028A28A0028A28A0028A28A0028A28A0028A28A0028A2
8A0028A28A0028A28A0028A28A0028A28A0028A28A0028A28A0028A28A0028A28A0028A28A0028A28A0028A28A0028A28A0028A28A0028A28A0028A2
8A0028A28A0028A28A0028A28A0028A28A0028A28A0028A28A0028A28A0028A28A0028A28A0028A28A0028A28A0028A28A0028A28A0028A28A0028A2
8A0028A28A0028A28A0028A28A0028A28A0028A28A0029BB86E2B9190324524922C71B48E70A80B13EC2B1BC34D25E4136AB2939BC7263523EEA0240
A00DCA2B3358B6BCB8B29D6D2FDAD58C4402AA0E0FAE7AD735F0C1E79B4CBE9EE6779A46B920B31CE481C9FD690EC76CCCA8BB998281DC9C53AAA6A5
651EA1632DAC9D1C707D0F6355B40BD6BDD3879873342C6293EA3BD311A9451450014514500145145001451450014514500145145001451450014514
500145145001451450014514500145145001451450014514500145145001451450014514500145145001451450014514500145145001451450014514
50014514500145145001451450014514500145145001599A5CF7867B9B6BE525E27CC72EDC09233D0F1C67B7E15A745001451581E28D566B14B4B2B3
38BBBE94448C064A0EED8A00DFA2B9DF0CE9FADE9F737C356BEFB542EC3C824E4F7C93E9DA9754D42EE5F1159E95A748119479D72D8CE13D295C763A
1A28A298828AC4D7755B8D26E6CA6280D8C9279731032549E87E95B4082010720D002D1451400514514005145140051514D20862690ABB0519C22962
7E80572FA978DADACEF2DACE3B1BA334F22AFEF93CB0013D73DE80B1D6D1451400514514005145140051451400514514005145140051451400514514
005145140051451400514514005145140051451400514514005145140051451400514514005145140051451400514514005145140051451400514514
005145140051451400514514005145140051451400514514005145140051451400514514005145140051451400514514005145140051451400514514
005145140051451400514514005145140051451400514514005145140051451400514514005145140051451400514514005145140051451400514514
005145140051451400514514005145140051451400514514005145140193E2695E0F0EDFC919C3088FEBC1AB3A3A18F47B342A14AC08081FEE8A5D56
D7EDBA65C5B77923207D6A2D06569B46B567C6F081587A11C7F4A5D47D0B579CD9CFFF005CDBF95721F0B1C368576076BA63F9815D1DEDEDC23CD6F0
E9F34DF2F0E3014E45739E0BB3D4F40B29ED6F34D9087977068C83C1A3A8743B6AE7340902F88B5DB65538132C84FB915D1D6078694CB71AA6A1FC17
572447EEA9F2823DA811BF45145300A28A2800A28A2800A28A2800A28A2800A28A2800A28A2800A28A2800A28A2800A28A2800A28A2800A28A2800A2
8A2800A28A2800A28A2800A28A2800A28A2800A28A2800A28A2800A28A2800A28A2800A28A2800A28A2800A28A2800A28A2800A28A2800A28A2800A2
8A2800A28A2800AF38D534CB9BCF8896762FA9DC811C0668E4046F4EB5E8F5C7E96AD79F12356BADBFBBB3812053EE403FE3498D1D4DAC2D05B242F3
493320C192423737D715CB7821FF00B4751D73577C3196ECC319F4451D3F515D81AE37E1B2F9163AA593FF00AD82FDF70F4040C7F23407436FC45ADD
BE83A5BDE5C727EEC69DD9AB06FAFF0051B68F4BBB4BA924D42EA542F648772796473F28E98F5A83E21D8BDFEA9A1C1212B6924E524627E50491D6B6
6E62B7D3A682C7468225BE9081BCAEF31277249F6ED4319A1E21B35D4341BDB66C02D1315246704720FE62AAF836F8EA1E17B299DF7C813639F71C56
9EA122C5A7DCC8E4054898927B0C1AC2F877098BC216CC47FAD667FD71FD28EA2E874F55EEEEA3B383CD955CA0383B1738AB1556FDEF12006C218E59
770CABB6D18A62290F10D81247EFF8FF00A62D51FF00C24DA79206CBAE7BF90D4B25C788BFE59D859FE331A8CDDF89463FE25568DEB8B8C5218F87C4
DA6CD7296E9E7F98EC157311EA6B6EB045CF895A44074DB3442C371336E2077ADEA1030AE1BE23605E6807F8BED7D71F4AEE6B83F89F288934760C37
ADDEE03DB14305B9DE5145417B77058DA4B7575208E1886E763D853113D155ECAF2DEFED12E6CE55961719575E86AC50014514500145145001451450
014514500145145001451450014514500145145001451450014514500145145001451450014514500145145001451450014514500145145001451450
014514500145145001451450014514500145145001451450014514500145145001451450014514500145145001451450014514500145145001451450
014514500145145001451450014514500145145001451450014514500145145001451450014514500145145001451450014514500145145001451450
014514500145145001451450014514500145145001451450014514500145145001451450014514500154ED6D5EDA79B6B8F2646DEA9B79563D79F4AB
945001451450055BD4B896031DABA23B1C166CF0BDF1EF525B5BC76B6D1DBC23091AED15351400514514005145140051451400514514005145140051
451400514514005145140051451400514514005145140051451400514514005145140051451400514514005145140051451400514514005145140051
45140051451400514514005145140051451400514514005145140051451400514514019936B76306AA74D999D2611F9B92BF291F5F5ACEF06D95C416
B777F7A08B8D46E1A620F65FE1FD2B765B6B799C3CD047232F42C809153D218572A623A1F8C5EEB6E2CB5601588E8B28E84FD727F3AEAA98C8AF8DCA
A7072323383EB4C457BFB0B6D46D8DB5E4424889CE0F1823A114DB0D32D34E5616B0842DF79892CCDF526AED140183E28792E6D5347B57DB717CC109
C7DD8FAB1FC8115AD636B1D8D8C16908C470A045FC0548638CC82428A5C0C06C723F1A92800A28A2800A28A2800A28A28021B9B78EEADDE098128E30
707159573E14D1AED225B9B6694C5F759A46247EB5B7450055B5B382CD0AC0AC01C677396CE3EA6A4BAB686EEDA4B7B98C4914830C87A1153514015E
CACEDEC2D52DACE258A14FBA8BD0558A28A0028A28A0028A28A0028A28A0028A28A0028A28A0028A28A0028A28A0028A28A0028A28A0028A28A0028A
28A0028A28A0028A28A0028A28A0028A28A0028A28A0028A28A0028A28A0028A28A0028A28A0028A28A0028A28A0028A28A0028A28A0028A28A0028A
28A0028A28A0028A28A0028A28A0028A28A0028A28A0028A28A0028A28A0028A28A0028A28A0028A28A0028A28A0028A28A0028A28A0028A28A0028A
28A0028A28A0028A28A0028A28A0028A28A0028A28A0028A28A0028A28A0028A28A0028A28A0028A28A0028A28A0028A28A0028A28A0028A28A0028A
28A0028A28A0028A28A0028A28A0028A28A0028A28A0028A28A0028A28A0028A8D6589DCA2488CCBD54302453B70DC572320648A0075145140051451
400514D775452CEC1547524E00A6C5347326F864591738CAB0228024A28A2800A28A2800A28A2800A28A2800A2985D03842CA18F419E4D3E800A28A2
800A28A2800A28A2800A28A2800A2932338CF34B40051451400514514005145140051451400514514005145140051451400514514005145140051451
400514514005145140051451400514514005145140051451400514514005145140051451400514514005145140051451400514514005145140051451
400514514005145140051451400514514005145140051451400514514005145140051451400514514005145140051451400514514005145140051451
400514514005145140051451400514514005145140051451400514514005145140051451400514514005145140051451400514514005145140051451
400514514005145140051451400514514005145140051451400514514005145140051451400514514005145140051451400514514005145140051451
400514514005145140051451400514514005145140051451400514514005145140051451400514514005145140051451400514514005145140051451
400561F8A2FC69B696D3CC5C5A99D56E1973C21EB9C76ADCACAD724B668E0B2BE844905EC9E41FF6491C7EB402112D34ABFB6592C5A2E00D935B91B9
7D0E477FAD66F8523D42DF57D66DF5195E765910C7337F1260EDFD2B02EFE1F5FD85C1B8F0EEAB24441CAA331523D811C1FC6B47C27AE6ACBAAB685E
22422E8A6F81C8E5C0CE73EBC0A43B1D7ADD5B3CBE525C44CE7A28704FE550C9AAE9D15DFD964BD8167C6761700D71FAB5AC36FF0013346FB246B1EF
8D9E40BC03D6A3F88F656C2FB449522449A5BADACCAA016195EBEB45C2C7A01200C9E2A826B1A6BCDE525EC25F3B71BB8CFD7A560FC47D524D37C36D
1C190F74FE56EFEE8EA69F7DA7DD2F81A1B3D0EDE233491270C3A6E1F330F7E7AD303675EB386F746BB8675DCBE5330E718201C561FC3DB88A1F055A
BDC4A91A89245DCE428FBC7B9AD8459EDBC2BE5DDE0CF15911260E7E609CD733F0CAC619BC34F3DC20977CAC815C642A8F4A41D0EE229639A31242EA
E8DD194E41FC6A4AE37E1DCECF6BAA5B64F936F78CB129FE107B57654C4151CB2C70A6E95D517A658E05495CFDC59DD8F148BEB99E0FECDF23CA11C9
FDE3F5EF9FD28037239A29B3E54A9263AED6071525723E7C507C4986D2D97CB592C98CAAA3018F5071F4AEBA800A86E255862676206012327D067FA5
4D54F50B3B6BBB675B985240AAC5770E9C638A00E2B43B27D7E1D2B568AEC248B72EF73F37CCD83900577DE6464E03A939C633DEB82F85BA75ACBA09
BC920569D2E98239EA000B4CF1369B05B78EF403683C933CDB9F693824367A5229EE7A092063240CF1CD29200C93802B84F89F6F0C1A6DA5FC4BB2E4
DE229707048DADFF00C48AEB352B486F34A921B9467511938C9072053117448853786529FDE078A11D2540F1B2B29E854E41AE2FE1ADB33F87676B96
324524EEAB1B124000F3C7D697C1D70D078A75DD262622D60937C51939D9CE0E29058ECA496385374AEA8B9C658E052C6E922078D9594F42A720D717
35C2EB3F118E993E24B4B280B8407037F19CFAF5AB8B2EA70F8F23B3B7B278B4816E77322E2366C6777A673814C2C5FD77568ED2FB4CB349F64F7174
80A83D539CE7DBA56D06049008247519E95E7DE30D32D25F1D6848F0A85B963E7638DD835D8DB69BA65A4D2C76F122C9347F3A6F24B274E84F4A4052
B2D1043E25BAD623D41A4132ED300FBA3F5F6ADD765452CEC140EE4E2B84F0D411DA7C4AD6ED2DC32C0900654DC7009D84FF0033536A17DFDB5E3EB7
D0FCC22CED14CD2A8C8F31C0E87E9914033B4474906637561EA0E69C480092700570FE3794F864D8EB1A6010969FCB9D178590104F23F03577C751DB
DF783E5BF491C18E312C0EAE5721B1FD2988EB28AE4AC3455D53C276F26A7248F74D6A36C8AC418F8E3152780755B8D4FC2492CD992E20678B731FBF
8E47E840FC2803A9A2B91F0DC3757963772F896DDE299663B6495CAFCBDB1CF18A6FC3DD4EE6FECEFE0B895A68ED2E3CB865639257B0FC31FAD0055F
1225DC3E35D1521D467115D4C0B43B8E0007271ED5DD5713E3071178C7C33231C2F9CC0F1F4FF1AEDA90D8514555BEB71730ED37134014EEDF13ED3C
7AFB5311C9EA125D789F59D4F4AB2B8F26DEC220A1C315FDF927938E703047E15D56956D359E996F6D71399E58902BC87F88D70DF0E74F179F6ED525
B9B8137DB483B64C07C007E6F5EA6BD16921B0A28AE17E28EFB6D2ED6FADE59639C5C2C7F2B9008C31C114C47754572DAD5BDDA785A5BF96EE74BF8A
1128649085561838DBD29FA13CFE25F0DDB5DEA2D2C2F229C792DB01EDBBEBD6803A6A2B98F03EA773A869F7905E4AD2CD6574F06F6FBCCA3A13EFD6
9DA7DCDCEA9ADEAB6B7D6F7096D6EE12121991187E1D4F7A02C74B4571BE15D56E1FC4FAC68F24EF716F6C7740ECDB8A0071B73DFA8FCAA7F19C3ABD
B69D3EA5A46A57113C4A19A10032903A9191C52B8EC7574572FA65CCB7FE0B8AF9352996768CCA6625490E38DA46318CF6AD1D06C6F6D6D43EA57D35
D5CCAA0B873F2A1EE00E94C46BD1451401C745AAEB7078E6DB46BB9EDE4B79626972B1E09001FD78AEC6BCEBC442F0FC52D352C9C2C8F6E06E3FC2B9
6DC7F215A7E289EFBC34B6FAB5BDF5C4F6E2511DC5BCEDB8303DC7A521D8ECA8AE5BC59E274D2BC391DF59BA99AEF6FD9F233C1C1271EC3F5350EBAF
7BE1ED263D546A37371242EBE7472B029203D401DA9858EBE8AAF637515ED9437701CC732075FC6AC5020A2A86B3A8C5A5695717D311B614240F53D8
7E75CFE9F35D5CF86DF5ED4B549E012C6654584E1225EC318E6803A3D423BC96D88B0B858260721990303ED5CFF81759D47588B503A9C88EF6F3F94A
1502E319CD697853577D6F4186FA44D8EC591BDC838CD733E00F35EEBC436B197894DD3112A8E54E48EF4867A0515C9786F57BDFF849352D0B5098CE
D6E37C3232E0B2E7BFE629FA8EA92DDF8B97C3C934B6F1983CC32463966C138CF618A6163AAA2B9386FEFB4AF16DBE8D7574D776D7711789DC0DE8C3
D4F7AB3E3117D6FA44FA8D86A12DB496D11628002AFF0081EF4058E8E8AF35FED6F1349E12B1D6A1BB67FDE856892304C80311B89C67B74ADED5F5D9
86ADA6E8EAD259CB7B186790264A93D00FC734AE163ACA2B8F9B53D43C3FE23D374EBBBA1796BA83145675C488D90073DC722B775FD5E1D0F4996FE7
058270AA3F898F414C2C69D15C8493EA90680DE20BAD4191D63F385BECFDDED3D148FC7AD68DCCD71ACF8556EACAE5ACDE6844BB82E4818C95FF00EB
D0162CE87A9B6AB6F71318D5163B8789707390BEB5AB5E79F0EAD3539B4032DAEADE4402E1C7946DD5F278C9DC4E6BA2F11F8923D265B7B2B7512DFD
D3848D0F41938C9A4163A1A2B95D5354D47C3B2DB5CEA3751DD58CF208A50230AD131C9C8C751C77A6F8D350D5F48D39F53D36F6210A951E53440F5E
E0F7A6163ACAC8D6A7D4E196C974B844824971312B90ABC7E55CCDCF8975A8B48D23548E2328BC7C340230723D73DB35378A757D7B46D56C043750FD
96F660810C40B274C8CF7EBD69058EDE8AC6F146A373A668EF7164AAD705D52356190493581AEF883C41A3DB5821FB1CF7F76D8FB3AC472BF8EEFE94
C2C77145635DEAEBA2684979AD488270BF32A7F13F5DA2A8C9A9EB10E943596481EDB6095AD80C3043DF3EB8E680B1D3D56BCB9FB25BB4C62965DBD5
625C9E95126A968DA42EA86402D5A3126EF6ACAB1D6B51BBB64D492C1069CE0B0F9FF7A57FBD8E98EFEB8A045AF0E6BF06BF6D34D043245E4C8632AF
D6B66B88F86B22CD6FAB4A8728F78CCBC63839ED5DBD0361581E24F11FF6045E6CB613CF171FBC520283F5ADFAE57E23A96F075D60818653CFD68047
4B049E7411CB8C6F50D8CE719152D72FA5EABA95D6896F75A7E9E8D02C400595C879001FC35A3A46B50EB5A49BDB08F738CA9899B69561FC24D2B858
D7A2B905F16DF7FC240346934331DC11B86EBA0011EB9C735A579AB6A306A10D9C7A4B3990F12897298EFDBF9D170B1BB453573B46E001C7201CD3A9
8828A28A0028A28A0028A28A0028A28A0028A28A0028A28A0028A28A0028A28A0028A28A0028A28A0028A28A0028A28A0028A28A0028A28A0028A28A
0028A28A0028A28A0028A28A0028A28A0028A28A0028A28A0028A28A0028A28A0028A28A0028A28A0028A28A0028A28A0028A28A0028A28A0028A28A
0028A28A0028A28A0028A28A0028A28A0028A28A0028A28A0028A28A0028A28A0028A28A0028A28A0028A28A0028A28A0028A28A0028A28A0028A28A
0028A28A0028A28A0028A28A0028A28A0028A28A0028A28A0028A28A002B03C53A5DEEAD0D9C3632885A2B8131949FBB8071C7E35BF4500528AE6549
120B8864DE703CC8D4B21E3AE7B73EB55A5D356E3C470EA0E9836B194463D58B67A7B004FE35AD4500725AEDA4D178CB49D58C45ADD11A295979D879
C1E9EF543C7CB7B35EE8D71059CB25B5ADC091D821DD9C8E31F415DE51487739CF146909E29F0E08E062AE0F9B0EE18F9802307F3356346BF68F4AB7
B7BAB5B85BB823547896224F000C83D31F8D6DD14C463EA01ED343B9431DC5CCB2AC80AC2858EE604F03B0E71593E015B9D3FC31F66BEB1BA82681D8
95688E5F278C7AD75D4500719E038AE6CE5D560BBB4B985A5BA32A17888041F7C7B57674514030AE52F45CFF00C27109D423924D37C9FF0045C21655
978EB8EFD7AFB5757450071B7F0DEC7F102C3536B395ACBC8685591724120F51DABA9B08E486C618E6C7988801C559A2800AA5A95D2DB59C8C55D98A
30555527271D2AED1401C6FC39DD69A00B4BA89E199EE1D951908C8C039FD293C656D247E20D0B5528CD6D6D36252A33B791835D9D359430C30041EC
452B0EFA9E7BF132F0DF68B6EB6713CB6F1DC077942FCB9C10067F1AEC5F53B39B466BC49D7CA78F00E79048E9F5AD1DABB76ED18F4C53A80B9C4FC3
9BD8A0F0D982EE4114B1CAEEC1F8E09CE6A0F094F1BF8DF5FB92C0453BAAC2C4603F27A7AF4AEF0804104020F5068000180300501738A3A7B68DF101
F55704DADFC643381911B71D7F2AE9269DAEA68A2B09810922BCCEA72028E76E7DFF00C6B4A988AA8A15142A8E800C0A60723E2EB598F897C3DA8A46
C60B79CACAE067602475FC8D6EDB85B8D6A6BC8C82890AC21BD4EE2C71FA56A51408E17429E07F88DAD5C2F11CD1224721E0330DA0807DF1535CE96D
A6FC438359DB8B4BB8CC4EE3F864C77F40702BB4A29587738BF1D5B7FC241FD9DA3593AC8CF71E6CACA72235008C9C7FBC6A7F1E1B6B6F05DCD9ABAA
B08D1624CF2402074FA0AEA2382188968E24427A95503344B045363CD8924C74DCA0E29818FA45DDB8F08DB4AD322A8B60A496030DB7A7D6B1BE1A5C
42BE19782360678E691DA2FE2C76E3F2AEBFECB6FE5795E445B339DBB0633EB8A12DE08E432470468E460B2A0048F4CD2038CF096A306A62FAFF0057
903DF473305B790FFAB4C70154F5CD56F075D3226A505D46F05D4FAA46E21285480594FE5806BBC1696CB706716F1098F59020DDF9F5A779311944A6
24F3071BF68CFE7405CE2BC6D2C6BE29F0DB1F99619D9A5C73B4657AFEB5DB433473C4B2C2E1D1BA30E869AD6F0B3167863663D495049A7A2AA28545
0AA3A00302980FACCF10DC8B3F0FDFDC119D903719C76C568B1DAA5B04E06703A9AE2350D657C55A2C9A65B432417535D085A1639211482CC7D06293
0468FC3CB092C3C296E274D924CCD2907AF3D3F4AEA2A2B78560B78E14FBB1A041F403152D3105707F16278BFE11E8210EA65174ADB01E47CADCE2BB
CAAF2D9DB4EDBE6B685DBD5E304D008C9F13CF0BF846F9D2446592DC84208F98E3A0A8BC117108F06E9F99517CA8B6BE481B4E4F5F4ADB3676C62F28
DBC3E583BB6796319F5C7AD22D95AA46F1A5B40A927DF511801BEA3BD0071BF0F2EADCDC6B999A30F2DFB328DC06E0738C55AB0BE9757F13EAB6D7D7
6D15B58308D2D836D0E33F798F7AE9469960B324EB6702CA872ACB18041A6CDA569F35D0B996D226987F1E393F5F5FC690EE705A25D4116B9E2B96DC
2DBE613F67DBF213F7B0477C938AF41863F3B4E8E29C97F3210AF9EA7239A8E6D274E9EED2EA6B489A78C615C8E40ACED675BB8D3B5CD36C9608FECD
745BCC9DDB0131DA8D83739BF0AC33DAEABA8F8564466B5867F3B7E38D9C1DBF8F15E87591A4F97777977A924785958451B9041755EFF42735AF4C18
51451408E17579123F8B3A3BB1C2FD94A138E371F3001F991573E24032F86C5A460BCF713A2C683A9F5AE866D36CE7B869E6815E560A0B1CE70A7231
E9CD034DB3174B75E56E9D7EEBBB3311F4C9A43B9C3F8DB42913C31A2B847922D376ACE8064ED21727F0DB8FC6BAA934AD0F50B01712DBA4F6D22F9B
96766046339C66B61955D4ABA865618208C822A945A3E9F0AA88ED82AA1DCAA18E01FA671405C76911AC3A5DBC7142218C2FC880E70B9E3F4ABD4514
C473BE3AB09350F09DE4302EE91409001DF07349E06921B9F06D8C7957D9118E4523A1C9C822BA3ACD5D16C639DE5815E1F30E64589CAABFD40A004D
2A542D7304104715A4126C87CB180DC65BDBAE6B99F87247DB7C443393F6D39FCDABAF9ACA096C9AD08290918C21C11CE7AD65DA784B48B2791EDA29
A3790E599676049F5EB48664690E9FF0B3B570D82C6D9769C7A119156350D56EAF7C64341B265B511C3E64F71805D9700ED5FCC7EB57CF84F49FB57D
AD5275B9EF289DB71FA9CF352EA3E1CD3B50D4A3D46412C57483689627DA48E983406872B2456F6FF14ACA182591DD6026569642E49C1EE7A56B6AF7
335DFC3FD465BA61E618DD738C679E2B52FBC37A65F490CB2C4C9343F7658DB6B9FA9EF4FBCF0FE9B7B6E904F0B18D176ED0C4647BFAD0172A780C2A
F8334D031FEACE71EBB8D54F106AADFF00091D868D651C497537CC6E5D0131AFFB3EF5AFA568561A448CF628F1865DBB3792A0673C0A8B56F0DD96A9
7F05FBBCF05D43F76585F6B6280EA725E22B68A1F1EF87A23732CF3F98AF2191F3FC591C741D2AC78F85E6A1E13B874F985ADF307083F8149009FA57
4979E19D3AF6DE34984BE744DBD6E43FEF437AEEABF6DA7DBDB587D8D54BC7821B7F25F3D49F53405CADA6C76DA8F872CE395126824B740CAC320E00
FEA2963789F46B85B5804514692244A060100119151C3A0C36E0C56B75730DB1C93023FCA33D87A0F6AD48E28E284451A058D4602F6C50070DE0094C
3F0FEE6559363AC929DDE8D81543C656D29F10F872F267748D92347957A86CE49CFAF35D85AF86B4FB39241079AB6F23976B7DFF00BB2DEB8FF3D2AF
EA7A6DAEAB68D6B7B1EF8DB907BA9F507B1A02FA98FAF786EC6FF4E6FED1BCBD7487E705A5CE08F6E3354FC66236F874DFBB6551145B55BA8E98AD88
7446408936A5793C09F761919707D031039A7EABA15A6AC9E5DDBCFE490018924C21C7B50171FE1C454F0F69EA830A204C63E95CC7C45C0D4BC364FF
00CFE7E1D56BA9D2F4B8B4B8BCAB79AE1E3C0555964DC100F4155351F0C69FAA5C2CF7AF732B231640653843EC3B501D4BFA9DCDAD9D8C9757DB3C98
46F3BC03C8E98F7AE5BC316CDA85EDCF8B756DAAAE0FD951F9F2A31FC55D0EA1A159EA5A62E9F7AD3CB0A9CE4C8431FA9EF551BC2B6AF69F626BDBE3
6806D1009BE503D3A50073DE3EF3F57F04DB6A31458549BCD20738424807F97E75BB1E8EFAAF866CA07D4EEA38A4B74DDE5ED1B8151C74AD91A7DB0D
2FFB3B666DBCAF2B69FEEE315434ED124D314C16DA8CFF0063C10B0B004A7D18F4A02E727E25B7165F0E5ECAC24B8961B7B9F2D9A5EA54127DB8CE2B
B1F0D146F0BE978036FD922E3DB68AB771616D7160F632440DBBAED283D2B36DF49BCB7D3A3D322BC516B18081C291204C741DBDBE940189F0D84223
D63C9E9F6D6E7B11CE315DBD729A1785EEB41BCB93617D1FD927943989E3C903D335D5D08185733F10A233783AF0019DB86FC8D74D54358D322D5AC0
D9CEC5626605B1DC0ED4D8909A1AC4BA1D888142C7E42600EDC0AE37E18ABA6A3AF040DF67F3D4293D3396FE98AEAADF4DBAB2B05B0B4BA558146D47
64F9D07A0EC6ACE95A5DB69364B6B68985CE589EAC4F52690CE5BE215BBD94FA77886DB892CA50AF81D549FF00F5D749A4DE45AAE750B690B5B9511C
63B7624E3EBC7E1505CDDE9DAEE8FA85BC3289234568E4E3A376EBEF57344D3E3D2B49B7B2897023419FAF7A03A17E8A28A620A28A2800A28A2800A2
8A2800A28A2800A28A2800A28A2800A28A2800A28A2800A28A2800A28A2800A28A2800A28A2800A28A2800A28A2800A28A2800A28A2800A28A2800A2
8A2800A28A2800A28A2800A28A2800A28A2800A28A2800A28A2800A28A2800A28A2800A28A2800A28A2800A28A2800A28A2800A28A2800A28A2800A2
8A2800A28A2800A28A2800A28A2800A28A2800A28A2800A28A2800A28A2800A28A2800A28A2800A28A2800A28A2800A28A2800A28A2800A28A2800A2
8A2800A28A2800A28A2800A28A2800A28A2800A28A2800A28A2800A28A2800A28A2800A28A2800A28A2800A28A2800A28A2800A28A2800A28A2800A2
8A2800A28A2800A28A2800A28A2800A28A2800A28A2800A28A2800A28A2800A28A2800A896086395A548A3591BEF305009FA9A968A0028A28A0028A2
8A0028A28A0028A28A002A39628A64D9346922F5C3A822A4A2801AAA1542A8000180076A751450014514500145145001451450014514500145145001
451450014514500145145001451450014514500145145001451450014514500145145001451450014514500145145001451450014514500145145006
5D9E8969633C925B34A8B2486568B77C85BD718AD4A28A0028A28A0028A28A0028A28A0028A28A0028A28A0028A28A0028A28A0028A28A0028A28A00
28A28A0028A28A0028A28A0028A28A0028A28A0028A28A0028A28A0028A28A0028A28A0028A28A0028A28A0028A28A0028A28A0028A28A0028A28A00
28A28A0028A28A0028A28A0028A28A0028A28A0028A28A0028A28A0028A28A0028A28A0028A28A0028A28A0028A28A0028A28A0028A28A0028A28A00
28A28A0028A28A0028A28A0028A28A0028A28A0028A28A0028A28A0028A28A0028A28A0028A28A0028A86EAE12D6D64B894E12352C6B3BC3A1E5D3FE
DF39DD35E1F3189ECBFC23F01401AF4572D69ABDA6AFE24BDD32596656B7251210C555F6E77138F7AA77DA8CBE1DF17585847349259DF81B9246DDE5
9CE323F4A571D8ED68ACDD721926D31DA0729343FBD8C83DD7FF00AD9A9347BF5D4F4AB7BC423F7880B0F43DC53117A8A28A0028A28A0028A28A0028
A28A0028A28A0028A28A0028A28A0028A28A0028A28A0028A28A0028A28A0028A28A0028A28A0028A28A0028A28A0028A28A0028A28A0028A28A0028
A28A0028A28A0028A28A0028A28A0028A28A0028A28A0028A28A0028A28A0028A28A0028A28A0028A28A0028A2B1FC4FAAB68FA3BDC440199D962881
FEF3702803628AE6F4783C4306B932DFCCB35818832B36321CF61F4E6AC788F53B9B48EDECF4CD8DA8DE49B610C32140FBCC7D852B81B94557B15BA5
B3896FA48E4B803F78D18C293ED562980514561683AB5D5D5EDE69DA94423BBB56C8C74743D0D006ED1451400514514005145140051589A9EB0F6DAC
D8E996912CD3DC1DD2027EE2773ED5B7400514531D9823145DCC070338C9F4A007D159BA26A8BAAD979DB3CB9518A4B193CAB0AD2A0028A28A0028A2
8A0028A28A0028A28A0028A28A0028A28A0028A28A0028A28A0028A28A0028A28A0028A28A0028A28A0028A28A0028A28A0028A28A0028A28A0028A2
8A0028A28A0028A28A0028A28A0028A28A0028A28A0028A28A0028A28A0028A28A0028A28A0028A28A0028A28A0028A28A0028A28A0028A28A0028A2
8A0028A28A0028A28A0028A28A0028A28A0028A28A0028A28A0028A28A0028A28A0028A28A0028A28A0028A28A0028A28A0028A28A0028A28A0028A2
8A0028A28A0028A28A0028A28A0028A28A0028A28A0028A28A0028A28A0028A28A0028A28A0028A28A0028A28A00C4F18BB27862F4A8CE42A9FA1602
B52D13CBB4853A6D8D47E42B2FC5F6D2DDF86AF2285983850F8519CE0E71FA5686993FDA74DB69DB83244AC47A1C734BA8FA1564B5D334C9AE354785
2291F9925E493F41EA7DAB223D21B51D6A2D7F57611A43816D075DA33C163F5355F5BD535D1AD116DA0CB75696E4ECC8203374DDEFED53E9B7FE22D5
B51852F34B3A7DA44E1DCB672F8E839A067584023079AE6FC160C706A76E3EE437F2AAFE75D0CB22C513C8E7088A598FA01589E0E8593477B87E1AF2
E24B823D89C0FD00A3A8BA1BF45145310514514005145140051451400514514005145140051451400514514005145140051451400514514005145140
051451400514514005145140051451400514514005145140051451400514514005145140051451400514514005145140051451400514514005145140
05145140051451400579AF8E6D3539B5CD1F4D6D48C8973396406351E59DC307DF19AF4AAE12F47F687C5AB28B194B0B52EC3D0E0907FF001E5A4C68
EAF49B7BEB6B531EA37A2F25DC48904610E3D0E2B9FD249D4BE216AB75210574E896DE15F427393FA1FCEBB0AE3FC267CBF16F896065F98CEB2027AE
0E6803A4D4AF7EC368D288CC92121238C75763D0572F7BAC6AFA36AF696F72E9772EA2BB521C6044F90074EDEB5B3AD6A76365AAE956D756E659AE65
2B091D233C0C9FCFF4AC4691752F8A31A91BE2B0B66DA71C2BF7FE7430449E1CD4B573E2FD474BBF91E7862883EF2B80ADC74F6E4FE552F888FF0067
78BB45D4970AB70CD6B31F507A56CE9F3CF73A8DE49E603688C2388051F780F98E7AF5FE7585E3F469A4D0E08BFD6BDFA9503AE00E68E81D4EC6B275
6F10699A448B15E4E7CE7FBB146A5988FA0AD6AF2BD2EDB50BDF1B6A971FDA105A6A10CA42A4F186CA678C0FCA86091E87A56AD06AB1C925BA4CAA87
07CD4DA4FD29757D5AD748B5135D31CB1DB1C6A32D237A0153D889D6D505D4C934C33B9D0601E7D2B92B8B77BCF8A700B86630DA5B79B12F6DD8C7F5
27F0A0075EF8E66B0BA4B59FC3F7AB34A9BE34DC0961F402BA0D02FEEB51D2D2EAF6D4DAC8EC7119CE40EDD6B9DF17DDBAF8B7C3B6F689BEE5663236
3AEC240C13E9D6BB53920E0E0FAD00CE4BC1E5350D675BD60FCE5AE3ECF0B1EA1140E07B74FCAB7F53BE36502F969E64F2B6C8A3CFDE63FD2B9BF867
95D1AFA26EA97D203F90AD796FB4A7F154164E59F5148895FEEA8EBCFBE3F9D1D01EE505D7753B7D4FFB1AEA2866BF940689E2076AA9EA587B54BE12
D7EEB569F51B3BE8D04F632EC3220C2BF247F4AA3E1D6FED4F1CEB1A96E0F15B016F1301FE7D2BA1D26792E1EEE530C51C26622274183201D58FAFFF
00AE840CC984BE9DE3D92D933E46A5019B19E8EBD4D7535CAEAEBBFC7DA22AB6196095987B7F9CD75542061451453105145140051451400514514005
145140051451400514514005145140051451400514514005145140051451400514514005145140051451400514514005145140051451400514514005
145140051451400514514005145140051451400514514005145140051451400514514005145140051451400514514005145140051451400514514005
145140051451400514514005145140051451400514514005145140051451400514514005145140051451400514514005145140051451400514514005
1451400514514005145140051451400514514005145140051451400514514005145140084060410083C106A969967FD9F0BC0B2178FCC2D1827EEA9C
71F9E6AF51400514514014F52B692F2C9EDA29047E6FCAEC73C2F7C7BD4F04490411C318C246A1547B018152D1400514514005145140051451400514
514005145140051451400514514005145140051451400514514005145140051451400514514005145140051451400514514005145140051451400514
514005145140051451400514514005145140051451400514514005145140051451400514514005145140142F356B4B2BEB5B3B890ACB74488F8E38F5
35CEF83E06BBD7F5CD7240764F2F93031FEE2F5FE4BF9574B7FA6DA6A0A8B7902CA236DC992460FE153C10456D108A08D6341D1546052025AE535455
D1FC616BAA8056DEF93ECF704740DFC24D75750DCDB4177098AE6259233D558530454D4F48B4D50426E03892062D1491B6D6427AE0FBD453683632A4
0AA2489A0042BC6F8620F504F7AD400280140007000ED4EA0086DADE1B584436F1848C7402B9A0A759F1BA5C261ECF4B8CA86EA1A53E9F4CFE95D5D5
7B4B4B7B287CAB689624C962147527A9A432C565EA7E1FD2B547F32F6CE379718120F95BF315A94531152C2C6DB4EB45B6B58F644BD016279FA9A835
0D261BE9E2B912490DD420849A338201EDEE2B4A8A00C8B0D0AD6D2FDB5095A4B9BE71B4CF29C903A600EC2B5E8A280394D1A3FEC5F15EA3A7BE161D
41BED56E71819E8CBFCAB4EEFC3D6777AB2EA65A58AE827965E36C6E5F43C5685C595B5CCD0CD3C2AF240DBA363D54D58A561DCCAFEC4B5494BDBB49
6C1805748480AFF5E2B462892189628942A28C003B54951C88B2C6D1B8CAB82A4671C1A6239BD2A1FED3F15DDEB18FDC5BA7D9A03FDEFEF1FE75D454
36D6F0DAC0B0DBC6238D7A28A9A806145145001451450014514500145145001451450014514500145145001451450014514500145145001451450014
514500145145001451450014514500145145001451450014514500145145001451450014514500145145001451450014514500145145001451450014
514500145145001451450014514500145145001451450014514500145145001451450014514500145145001451450014514500145145001451450014
514500145145001451450014514500145145001451450014514500145145001451450014514500145145001451450014514500145145001451450014
514500145145001451512CF0C92B4492C6D22FDE40C091F514012D145140051451400514514005145140051451400514514005145140051451400514
514005145140051451400514514005145140051451400514514005145140051451400514514005145140051451400514514005145140051451400514
514005145140051451400514514005145140051451400514514005145140051451400514514005145140051451400514514005145140051451400514
514005145140051451400514514005145140051451400514514005145140051451400514514005145140051451400514514005145140051451400514
514005145140051451400514514005145140051451400514514005145140051451400514514005145140051451400514514005145140051451400514
514005145140051451400514514005145140051451400514514005145140051451400514514005145140051451400514514005145140051451400514
514005145140051451400514514005145140051451400514514005145140051451400514514005145140051451400514514005145140051451400514
5140191E26BDFB0683793248A92889BCB04E093EDF9D62E83A44536A3A7EAF6B7AA116D47990AB7CCCC7A93ED5A1E36B4827F0DDECD244AD3450931B
91CAF22A5F09DAC10787AC248E048E47B752C42F27233D69751F435A7B9B7B750D713C5103D0C8E173F9D2C33C33A6F8254953A6E460C3F4AE4BC32E
35BF106AF7F76AB20B7710408DC845E7B74A76B178347F1A6962222382F50C72A01C31CF071F8D170B1D23EA7A7C527952DF5B249BB6EC69541CFA63
34F17B6A6E7ECA2E2333E09F2C30CF1D78AE37E25D98921D3668424731B90A5F0327381F8D6ECFE1E5935AD3AFE39F62D926D29B797F7268037A8A2B
03C6DA9DCE93E19B9BBB375498155563DB271C7BD311A8752B113BC06EE11227DE52E062AC3BA4485E475551D4B1C015C4F872F74D5F08234F04D3DC
4E8CD2E2DD98BB9CF46C63F5AD4F061BEB8F0E7D9F58B770D1BB4404C3E664EA323F1C7E14AE3B1BCB736EC14ACD11DDC2E1C73F4A719E112796D2A0
7FEE9619FCAB86B4B1B5D13E22984C4A2DEF63DD06E1908FE8BF88ADDFECDB5D47C4F35EBC2AC2D50459238793AF3EB818145C1A37DDD6352CEC1547
524E00A6433C330261952403A94607F9571BAAEB0B73E3DB6D225320B5B44F36455190EFB72320750011FAD4DAAEA170DE26D2AE349B2BC92324C777
28B760A632475C81D39345C2C75535D5B40C126B88A366E81DC027F3A779F109842D2A0948C842C3711EB8AE47E25AC4346B593016E3ED71AC4F8E47
5AD5B8D1279F58D2EFFED4116CA3DAC9B725CE39E6981B4D2C4B9DD220DBD72C38A5791131BDD573D32719AE23C7BA6C56B2D8EBD0411E6DE702E46D
E1D49EA7D79FE75B7AB59D8EB30E9B032AB4723F99195E30806481FA0A41636E3963954989D5C038254E79A92A1B6B682D61115B44B1C63A2A8C54D4
C414D242A96620003249ED556E217B99D23900FB328DCC33F7DBB0FA0EBF9572DA579FABEB1AEDB2CDBAC6D5DA282224ED0E41073ED9CD2B8EC7668E
B2286460CA7A10720D3AB17C2BA5CDA3E891D9DCB8794333120E7AD43ADF888586A36FA65A46935F4E3761DF6AC6BEAC680B1D051597A41D4584AF7F
796970A4E13ECEBC2FD4D6A5310D660AA59880A06493D05719A7EB126ABF112581370B7B5B7213B6ECFF0017EB5B5ACEA7611CC961777291ABFCF303
9E107638E99AE4AC2F639BE2FCBF66915E0787CB5DB80302207F1E41A4C68F48A28A298828A2ABDDDBA5D5B3C0E5C2BE3250E0F5CF5A00B145703E15
5913C51ADE8FA85C5CC8D183E4EE908FDD93F5F42B4BE076D9AD6AFA7DD5CDC4B227FA9324879889EDEFD3914AE3B1DED71DA7CD776FF112E34D7BEB
99AD7EC9E6A47239214E47F89A9BC29A480F2EA925E5D4C1A6905B2C8E7688F24038EF9ACABB8A597E2D2C714AD106B40642A704AE3919FC28047A05
15C36AAC740F19E8CB60CEB05F9F26688BB303CE33C93EBFA559F1EEB92E9CB63616F2796F7B285771D4264038FCE8B858EC28AE17C5B7767A158DAD
EE913C62EA3B8505639371913077023357FC616F6FA8784E7D4D43A4A96C2589C3104023383F9D170B1D5D15C059695A96A5E06D1E4B1B9FF4E4904A
B23B11B1493919FC055FD6B567B9F18E9FE1D8E6F2E22BE6DCB29C16E0909FA7EA28B858EC28AE36FF00508344F136991584A9F66BD2639E246DC17D
1BDABB2A620A28AC1D575E961D563D274CB71737ACBE63866C2C69EA68037A8AE26E358F177F6B4DA45BD9E9E6758FCD4977100A138C804F63C57497
975358787E6BAB86533C36E59881C170BFE34AE3B11595F4B7FAB5CAC0EBF64B53E5B60677BF7E7DAB5EB9FF00045BFD9FC2964CDCC93A99DCFA9639
FE58A9EFAE2E2EF521A6D9C8D108D77DC4CA3951D947B9E280366AAEA2D72B6133598CDC2AE5063A9F4AE36DAF759D50EB1A6D95D48174D2E23B8539
9256C9DAB9FC0FE95D0783DF5493C3F03EB21C5D127EFF00DEDBDB345C2C5DD23518B54D3A2BB88FDEE197FBAC3823F3ABF5CAF864B59F88F5BD2C2E
D85245B88D7D03FF00FAABAAA60C2ABDEADCB5B30B39238E6FE1691722A3D52FA1D334E9EF6E0FEEE142C467AFB564692FA8EADA3DBEA8D78609E606
48A3500C614F40C3BFD7B66810CF07EB37FAA49A9C3A9795E659CE221E5AE077FF000AE9AB84F86D34D35DEBED729B263740BAFA37CD9FE55A371ADC
9A878B9741B398C51451B3DCC8A30C48C7CA0F6EB4AE3B6A7554572C750FEC5F14DA692D70F2DADF4676091B7347203D3E873557C7B7DA9E969653E9
B7A50CD3AC7E49030C7EBE87BD170B1D9D15C97880EB1A7E8926AC7522B716CAAED02A8119E4647BF5EB5BB6FA945268B1EA721D91340266F6E32453
0B1A14572BA1DC5E789AC1F5237B35AC2EECB0470F1800F53EB5378775B96EB50D4347BEDA6F2C1F1BC71E62766C7AF4FCC52B88D0BFD53ECBAB5858
2A876BA66DC3BAA81D7F3AD4AF3CBCD3EECFC4CB541A9CE59A07915881941CFCBE95D95BA3E99673CB7D7CF70880C85DD402A00E681B3468AE4B4A9A
FF00C4F6526A71DF4F651B3B2DB4719E303BB7AD5CF0EEBEDA9E873DCCA81EEAD0BA4C91FF00132FA7D68B858E868AE334BBFD4350D025D5B4ED4BED
171B18B5B3A8C2119E3DBA715D0E81733DE689677176C1A6923DCE40C73F4A61634AB335FD4BFB23499AF02AB32602AB74249AD3AF3CF8976F7A5AC5
DEFBFD125B8541084C6D3EB9CF3498257677D03996DE391861990311E848AE36E75CD7F4EF15E9FA55E7D99A0BC93E5745E4AE79FC6BA1B0B3D52DAE
57ED3AA2DC5B2A63CBFB38539EDC8AE57C6F70F6FE33F0FBC3179F32EED91671925863FCFB5008F40A2B91D5756D6742BDB3B8BE96DE7B3BA9844D12
26D3113E87BD6DEBDAB43A2E9335F4E46106154FF131E829858D3A2B92B8D435B83408F5E696320462692D3600021C1C67AE715D16997D16A5A7C17B
067CB9D032E695C562DD14564DFDEC9FDAB06976F32C12CB1993CC650DC0E3001A606B5158565A9DFC9E21934BBAB78C2456E25332E7E72481C7B75A
CED7BC4DA9787EEAC96F6D6D6582E5CA97899B2318E307BF34AE3B1D7515CDEA5AA6ABA5DBB6A5710C0D66982F08277C6A4E339EE6B46FB5BB1B2D18
6A93CB8B6640C84756C8C803DE8B858D3A2B9CFED3D5D74E1AB986236ACA24FB3807CC58FAE73EB8AD5FED4B3FEC81AA19716BE5F99BCF1C53B8AC5E
A2B9E8F55D4AE74E1AAC305B436663F31526277B2E339C8E066AE0D45B50D0FEDDA36C94C884C7BCE3FC9A2E162EDB5D4572D308F3986431B67D47FF
00AEAC5707E06BBD725D1E7996DA09BCDBA762EF260EEE01E3F0ADFF0012EA5A969162F7B656B05C4112EE903B90C39EDEB4AE3B1BB4572BFF00094C
EFA4D96A71596609E458E41B8EEC9FEE8FCEA6F12F8964F0F6C96E34D696DA49046B32CC072467A633DA8B858E928AE7F5FD7AF346B7376BA535CDA2
286795660A467DB157EC3517D4B438B50B483E79A2DE913B639F42698AC68D15C8693E3096FB52BBD366D3992F606C2C4AF9DDCF393DB15A30F88638
B516D3B54885A5C040E9F3EE5907B63BFB52B8EC6F5158577E206B0BEB68EFEC9E0B6B97091CFBC1C31E81876CD6ED310514551D5354B5D26DD66BA2
D877088A832CCC7B014017A8AC87D68DAAC6FA959C969148C15642C18027A671D2AD6A37DF61B56944134DF2923CA4DD8C0EF4019F2F8A74E8B518EC
9BCE2D24BE4ACA13E42FE99CD6ED79B784EE74DD5B4CB1B7BFB1BC966B6BA768E558CEC5666CF2C3EA3F2AF49A486C28A2B3EEF57B4B49764A642010
1DD50954FF0078F6A623428AC6BDF1269D677705B4AF2319D82ABA26E404F626B66800A28A280285DEAD6364CCB77729115193B81C0FC6A183C43A4D
CCC90DBDF4724921DAA132727F2A678B115FC2FA90740D8818807D7158BF0CEDA11E1382630AF98F2392C4024E0E290FA1D9515897DE28D274EB9FB3
DDCEF1BEEDB931B633F5A927F11E996EC9E64CDE5B903CD08760FA9A62B1AF4535583286520823208EF4EA002AA43A958CF74D6B15CC6F3AF54079AC
EF136AF6FA769574B2CC239DE16F2C60F2707BE2B94F0D8D1B559B47BC37D243A940A2392044387209EA71FD695C763D1E8AA773A85ADBB6D9A5F9BB
AAA9723EA0038A5B2D42DB5088BDACA1C29C30C60A9F420D3116E8A64B22431B492B05451924F6AA316B1632DCADBA4A773FDC2CA42B1F404F7A00D1
A6A3ABA86460CA7A10722B23C4DACDAE8FA54F25C3FCED190883A93D3FAD67F84F59D25748B2B58EF91AE250328492779EA3F3A0763A9A28A2810541
3DD5BDB141713246646DA9B8E326A7AC1D60E8B7B716F05FDD0596297F76BBB6FCD91ED401BD51F9D1799E57989BFF00B9B867F2A8EEAE63B681E492
48D30A48DEC00245717A7DBA6B11E9BA87F69416F791484CE378DCE73C7D78A571D8EF28A28A620A2A333C2090654047505850B344E7092A31F40C0D
004945355D5C651830F50734330552CC40503249E82801D454714D14C331489201DD581A92800A28AE4F5581E1F19E992A4F2797331DC9B8919C5007
59452646719E7D296800A28A2800A28A2800A28A2800A28A2800A28A2800A28A2800A28A2800A28A2800A28A2800A28A2800A28A2800A28A2800A28A
2800A28A2800A28A2800A28A2800A28A2800A28A2800A28A2800A28A2800A28A2800A28A2800A28A2800A28A2800A28A2800A28A2800A28A2800A28A
2800A28A28039CF1A4F28D0E7B382D6E2696E142AB469900E7BFE5563C3333AE9369653C13473DBC0AB2074C00401C7F9F4ADBA290EE71DA4DA7FC23
3ADEA06E98FD96F0F991CAA848073D0E29EF62BE21F15DAEA8BBBEC3609846231E6499278CF6E95D7514582E713E36796FE5B182DEC2F5DADEE43BB0
818AED07A8C75E95D5D8DEC77B1B3C6932056DB89632A7F5AB7451600ACBF10E969AD68775A796DA655F95B1D181C8FD456A514C472BE19BBFEC8D1A
2D2F538E486E2D32B8085838C9208238ADEB213B0926B9CA991B2B193F71474FC7D6AE51486721F10B4E9A7D3ADB50B2DC2EECA505368E70481FCF15
D0E956EF6FA744933179986F9588C1663C9CFF002FC2AF514C2E71BAD69735878C2DBC4B6F19961DBE5DCAA8CB28C63701F4C56CDC5F8D42DC43A548
ECCEC034A159422F739E306B668A02E713E3F326A56D6D656505C493C176923E2238DA01CF3F88AEBACEEA3BC804B12C8AA78C3A153FAD58A28114B5
5B18B52D36E2CA6194990A9F63D8FE75CCFC3CB5BEFECF69B5266DF066DA243FC2A0F3F8E6ADC0FAD59EA1A949722E2E04A716B1A2931A7BE7B0ADCD
32D8DA58C50B9CC9CB39F562727F534865CA28A2988E6D755D4D757D4A19610C8985B48557E663FDE27D3BD66F86629F43D22EEC678DDF5496773840
497C8E1B3E9DEBB6A29587729E9D6DF62D3E1B62DB9A35C16F53D4D727AFE8D247E2E5D6A7D38DFD8BC3E5CA89CB2103AE3B8AEE28A2C1732F445B31
667EC3632DA44CD9DB226C2C7D719CD43ACF886CF4B8E656F326B98C710C6849248E39C62B6A8A6239EF0C5D5BCFA735ECB201772832DCEE0415EBC6
3D0572525D449F16175201FEC9B36B4BB0E01F28AFF3C57A7521008208C83D41A43B956D6FA1BCB43736FBDE3E71F2105B1E80F5AE6BC27E2CB9D735
7BBB2BBB45B7F2D4BA601048C81CE6BAE00280140007000ED58BA2A35E5FDDEAF24663136228411826353D4FD4FF002A00DCA28A2988E0FC656D7969
E2DD1F51D2FE59AE736CF81D791C9FC18FE54BE2BB19EC7C47A4DEE991806E14D948064E54FF00F5BF9576EF1A3952EAADB4EE5C8CE0FA8A478A3919
5A48D58A1CA9600ED3EA3D2958771B0C71DA5AA449848A140064F40077AE17FB46C4FC55138BA84C5F64F2FCCDE36EEE78CF4AEFD943295600823041
EF55DB4FB26186B4B723D0C4BFE14C11C5F8BEFAC4F8BFC3D27DA6222194B48CAC084191D6A4F883A73DFDBE9DAC58E2E16CA5CB88F0D9424648C75C
15FD6BAE3A5E9C4B13616A4B7DE2615E7F4AB11471C5188E245445E8AA3007E14AC1732DB50D324B51756E6DE59644C44A1416638E17039F6ACEF13D
C259782EE6D2F2757BA92D4A05EEC4F1C01E9FD2BA08ECED2197CD8ADA1490FF001AC601FCEA49218A5C79B1A3E3A6E5071401CFF81AE616F0C585B0
954CF1C3978F3CA8DC7B563EB56274AF88569E20B842D6330F2E493191136CDA09F6E9FAD7729145164C71A267AED50334B22248852455653D430C83
4582E64DFB5AC7120D3E0B692EE46021081720FF007BE8315B355EDECAD6D8B35B5B43096EA5102E7F2AB14C415C6786217B4F16F8824D42602E2475
2B9E018F9208FCC576754EEF4DB3BDDA6EADD242A3009E0E3D38A4072964D3EADF12A6BEB7E6C6C20FB3990747383C67BF2C7F2AE8FC496DF6BF0EEA
16E3ABC0D8FC067FA55EB5B582CE010DAC491463A2A8C0A94804104020F5068B0EE61F832559BC25A6329CED8021F6238C7E954EDBC4AAFA96B686CD
638EC4EC12670D3480602E3A9CE38A97C3AAFA55FDDE8B32111EF69ED5C0E1918E48FC0D6A9D274F3A88D44DA47F6B03FD66DE7FC33EF401C7784AE8
697E13BFBD9CECBE9AE5C989BEF073C28C75F7AECB4A86E21D32DD2F24696E0266463D4B1E4FF8523E8DA6C97BF6D7B389AE339DE4739A76A97F0E9B
612DD4E70B18E07763D80A01EA60687BE7F1E6BD70BFEAD123849FF680AEB2B9EF0869B2D9698F717431757B21B8947A67902BA1A10339BF1E594D7F
E13BC8EDF25D009303A903922A5F054C8FE0FD3983821210ADCF423AE6B78804104641ED5976FA159DADD34D6CD3C4AC773429291131F52B401CBFC3
99D27D63C4CE87E57BB0EB9F42CF50369E9A77C506B8D442FD9AFE3631337DD2DC7CA7F2EF5D33784F49174F730452DBC921CB98642A1BEA2B56F2CA
DAF62F2EEA049541C80C3A1F63DA8B05CC4D574AD3ED7C99ACECA01A8B4A82DD98648208C9E4F4033FA566FC48C6CD14E4717CB5D55B69D696AFBE18
42B631B892C47D09E959F75E16D22F2E9AE6F2DDE77662D8791B009F41DA8B05C8BC76C078335339EB1003FEFA150DB5A4BA8FC37B7B484812CDA7A2
AFD768ABF71E1AD22EA1F267B5668FBA995F07DBAD4FA668D63A56F1611B448E00D9E631518F404F14018DF0E46CF0941132ED78E474753D41CD6778
76D649BE246B9A844736E8BE5161D198EDE3FF001D35D3CBA25A493C93C4D35BC937FAD303ECDFF5156ECACADEC60105B4611073C7527D49EE6981CA
DD3A8F8AB640F1FE84E327B9E6AEEAB0DD6A56DE21B38277E62448829CE1B6E4A8FAF4FC6B4F54D0ECB539E19EE15D6783FD5CB1B6D65FC6AD5959C3
65114841F98E5998E598FA93482E727E05B1B5BCF0A5BC26E2E43C4592789652A15B713823F2ABF610D96883575D2EDC2A41179AE4BE4348149C1AD1
9341B36B892788CD6EF28FDE0824281FEB532E976A9A649608ACB0CAA558EE258E463393DE8B05CE2C78620D534D8F5CF0BDDBD8DC4E9BDA246F919B
BAFB73C5751E13D4A5D5740B7B99E3F2E65CC720C606E5E2A34F0D456CBB34DBBB9B289942C91C672AFF00ED63B31EE6B52C6CA1D3ED52DAD936C699
C77249E493F8D082E5BAE17E28BE2CF495EE6F01FC8577559DACE916BAD589B4BC5257219594E0A30E8453042DD4D22EAD63046D847F31A45F50178F
D4D72DE254DDF123C399C636B63EA326BA5D2F464D3DCCB25D4F7736DD8B24CD92ABE82AA6ADE14B3D57548F519AE6F22B888611A1902EDFA641A406
47C4CFF8F0D2C7ADEA0A5F8A76B3CFE1657854B2433ABC807618233F99FD6B4B52F095B6A71A25EEA5A94A2360E9BA55F948EFF76B52DB4E48F4F6B2
B9965BB8D810C6720920F6A02E65AE896FABE8B1236A77EF693C29F2AC8A032E0718DB56FC3565069DA59B3B5F30C3148EA8D236770CF51ED505B787
459A182D351BA8ACC927C8041001EA01EC2B6ADE14B7812188611060508096B98F10D8DA6B5AAC1A7492C96B7B0C5E7C171191BBA90463F0AE9EB1AF
74186EB581A9F9F2473A43E5C657F80F3CFBF5E9430463E8D71ABE95E218F4ED75D2E52E232B6F7607CCC473B4FEB50F8FCC5FDADE1B4B923C96BC25
C1E98CAFF8D745069931BE8AEEFAED6E1A15222511040A4F53D4F35CEF8FE18AE359F0D413AEE492F0AB0CE382533CD03EA6D78D403E0FD4F774F27F
A8AE23598EF0FC23D30C8A084752E08E4265B69FE43F1AEF6F34637D6FF64B8BB77B338DD16D00B01D016F4FC2AF4F656D3D8358CB1836ED1F9653B6
DC62815CC78F4D9751F0FDBC11EAB30B696D906E58D41652BFA718AE7BC53A7B597C3B96C2CDE6992D2E363B3F5DA09FD3915D358E8971A7DB7D8ECF
53912D14111A34619901EC1AB40D85B1D3CD8347BE0652ACADCE73D7F1EF4582E50D22D21BFF00075859DC8DF149671A301C7451FE152E9F1D958E95
3D969EACB159864E477C64F3DFAD4163A35EE9D17D96CB532B6AB9F2D24843B47EC093C8AD2B3B18AD2D3ECE32E1B25D9BAB93D49A00E5BE1FCCD178
3AEAE1704ACF33A83D38E6AFF896692E7C03773B285925B35765EC0900914DB1F0DDDE9B0CF61697CA34E99CB6C64F9D01EA01AD0D5B499750B26B28
AF5ADADDA3F2D91630DC7D4F3E940752A7846DA197C1FA6473C31C8AB18601D4100E4E0F3DEB1FE2BA83A159B16C15B918FC8D74BA26953E956B1DA9
BE6B8B78936A2B46148F7C8A87C51A0C7E21D345AB4C6178DC491BE33823D47A51D03A91F8D70BE0DD481C0C418FD453FC19FF00229699918FDC0AA9
AB68FAD6B1A24B61797B6913B2E330A31121FF006B3D3F0F5AB9A169DA869BE1D8B4F92780DCC4A56390296503B641C13DE8EA1D0E7FC210A378F3C4
F3B7FAC49022FD0B127FF411537882059BE23685BC65562623EA09AB7A3F87752D2F5BBCD405EDB4BF6E7DD3A9888C727A73EF4CD4341D76EF5FB6D5
45F5921B504449B1B183D73407519F13067C3B010395BC8CFE86BAE43B9031EE335CDF89344D5B5CB48ED85DDA40892090E1189247E35AFA5AEA6A25
FED392173BBF77E50C71EF4750E85FAE0BC56CEFE3EF0F432022059432FBB673FCF15DED61F887425D4E4B5BC80AA5ED948B242CDF74E08383F95304
41E3E553E0EBF2DFC2AAC3EB9153F86AE1EE7C21673DC36E76B6F989F6C8A8B50B3D4B5BB6365776F1DA5B315F3584BBD9C03D063A55EBD86EE0B05B
6D2E18594218F6C8C460630290181F0C5B3E1EB800702EDF1F90AECAB8FF000CE9BAEE816C6D3ECB69344F299199652A46719E0D76142D81EE4729C4
4E4304E0FCC7A0F7AE6BC51B74DF055E4624F35E71B7793CBB39EBF97F2ADAD6ECE5D4348BAB382411C9326D563DAB9DBAF0D5FEA5A25BD95CBAC3F6
30A6141216F318752C7F950C1165744D3D345D16D6FC48248593CBD8D83E611939F519AEA2B26DA0BCBBBE8EEB5089215801F2A257DD9623963FD2B3
B5C8BC4527882C1B4B94258823CDC1E3AF3B850074F451453118FE2C19F0B6A4338CC0DCD67FC39C7FC2156214E7064CFF00DF6D5378864D4EEAC2EA
C2D749697CD52AB279CA0633D7FF00AD59FE0E8F5AD1748B7D32FB4872AAE712C72A1C0273C8CFBD21F41BF11E05B9B4D36171F2C974AA4F7C1AD4F1
65AC69E0CBFB78C0548EDFE4CF6C631FCAB3BC5F06A7A9C96D0596953C8B6D70B2990BA057031D39ABFADDC5CDF6893DA7F64DE07B88CAE3E5217EA4
1A065AF09B33F85F4E6724B184724F5AD8AC2F0A1B8834A82C2E6CA681EDA30BB9C7CADF4ADDA10999DE2048E4D06FD645047D9DC8CFAED38AC9F87D
1463C1FA7C9B54B957F9B68CFDE3DEAFF8826BA1A7CF6F6BA7CB746685941461804823BD65782DF51D3F49B3D2AF349B888C7B81972368C927FAD1D4
3A16EDED4691AEDF5E4D786617B831DBA8CB839EC3D2B334899DBE25EA03C978164B30C51B8C9CAF3F5A75AAEA3A778BB52B9B8D3EEAE21B81FB8789
7705E7D7B55592DB58B2F16DE6A373612DD47796CD1466D86EF2F38C027B74EB40CBBE32BF93ED9A4D8447F757132BBB03C30CF0335A7E28D36FB51D
3A2B7D34C48CB2063BB8C01E954BC47A04FA969761240717962158293C3600C8FAF157CEAD3CF0F950E9D76B74CB82ACBB4464F7DD408678BD3FE28F
BE5930C4420138CF391CD3F4496DACBC31A6CB200ABE4C6A085E7240A5D5F4C9A6F0A4DA742C649BC808093D48C7F856569F3CF77A3D96906CAE63B8
882ACA658CAAA85EF9A3A8743AFA2B95F1B6A9AC6976D0CDA55BB3C60E65755DDB707BFB5741A7DC3DDE9F6F71247E5BCB1862BE9914C562D5723E37
863371A44A5407FB5A8C8EFC8AEBAB93F157DAEF25B25B5B1BA916DAE43B911F500F6F5A4C68E8351B3B7BCB4923B882394053B77A86C1C7BD721E01
D16C66D1DE6B9B48E4996E1B1211F371EFFD2BABB8BF44D384F2433A991481198CEECFA1C74AE7BC17782CF4EFB15CDB5CC72BDC12331103E6C77A3A
8743B1A6B28652AC01046083DE9D45311C26A1A65B69DE34B49258C3DA5E0D9B198E11A9352D3E2D33C636723C2DF60BAFDD90492A1AB67C6B60F79A
2192004DC5BB078CE40C1F5A8351D3E77F06C5E7BEEBA831386C7F1673DBEB48A2E7F6069F2EACD7022F92250BE5862177F5CFE5B7F3ACDD5F5685FC
5D6FA6DDDC7976B12EF65ECED8E86BA0D144C74B824B93BA79143C8DEA7FFD58AE7B5BB3934DF165A6B8B197B63F24C557EE678CFF005A1890FD4B50
B3B0D56CEEB4E04EE3B2648D4852A7BF4AEBAB2A5BD8352B6782C66124920E0A838519EA6B500C000927DCF7A10316B91F154724DE22D1E189F633B1
04FB575D5C7788AF2D97C55A4C8678C25BB91292DF77FCE686087F8AB4F8F4FD3A2BFD399A0B88645F9F71F9C7BD6E6AD2DD8D25A4B08CBCEC06157A
E0F5C565F8B6EADAE341FDD5C4727EF11B6AB0248CFA55AD435645F0D4B7B612F9AC881432F386E3AD0054D4ACED534069E48D2D2F521DE0EEF9837A
F1C9AD3F0E5D4B7BA0DA5C4E774AE9F337F788E33FA563BDE5ACBE19985B48B737B241F3EDF99C93D726A0F0ADDC32DAE856F11CC9147289477538EF
45C0AFE2EB29F4CB8B6BDB5B8B85B37902DC2090F033D8F6ADBBED2A1961B71672CAAF21001F34FCCB8C93F95686B964B7FA45CDBB75284A9F423915
8BE0F5BAB9896E2F508FB3298236DDF7B9EB8A3A8742DEAF31D03C3D2FD9E592493A2191F2413FD3FC6AB27D8A5D0E3BA178AB7CB0890C8B365B7E33
83CFAD5AF17D8CB7DA23ADBA96950EE00771FE7151F87CE94DA442CF0DA47304559832AAFCC060F5A03A17BC3DA836A5A3C37121CC9CAB9C7F10AD4A
CED26311A5C18D42C2D2931A85C0C7AD68D3105145140051451400514514005145140051451400514514005145140051451400514514005145140051
451400514514005145140051451400514514005145140051451400514514005145140051451400514514005145140051451400514514005145140051
4514005151AC819DD76B02A40C918078CF07BD49400515C36ADA86B77F7FAAA68D23AB69CC91AC0879933D4D75DA6BDCBE9D6EF7C812E4C60CAA3B36
39A0762DD145140828A28A0028A28A0028AA5AB5EA69DA5DCDE48C144519209F5EDFAD3745B89AEB47B5B8BA004D2C4ACD818EB4017E8AE7B5EB4D7E
E351B19348BC8E1B68CE664638DDCFD3918AE86800A28A2800A28A2800A28A2800A2B98F1B5C6AB61A3CFA8E9B7C2148117744620C5B2D8C83D4751F
955DF095DDCDF7866C2EAF18B4F2A16663D4F271FA501636A8A2B9FF00126B12D9DD69FA659102EEFE5D81B19F2D7BB50074145607872DF5BB79AF86
B338963DE05B907395E79F6ED5BF4005145140051454534A904324D2B0548D4B313D801934012D1599A05FB6ABA2DBDF48AAAD30270BD386207F2ACF
F157FC243BECBFE11F0BB7CCFDF1247E19CF6A00E8E8A0671CF5A2800A2B1FC4BABFF6369C93AA2BCB2CCB0C61BA6E6F5F6C03543C3FABEAD3EBB7DA
56AB6EBFE8E372CE8B8520E303F5A00E9E8A28A0028AE7FC5579AE59DB5BB685662E5DA4C483192076FC3DEB7909280B2ED2472339C5003A8A28A002
8A28A0028A29ACC154B310140C927A0A00751589A3EAD77AA482E23B351A7B97F2E7DD82C01C03B7AF38A96E7582B24915959CD76D11DB214E003E99
F5A57035AA29618A601658D2400E40650707D6A8695AD5AEA724D0C5BA3B8B7389A09061D0FD3D3DE9750D5534FBFB482688F9574DB04BBB856EC08A
6069D1451400514551D47518EC24B44910B1BA9D615C1E84F7A00BD4547334890BB449BDC292A99C6E3D867B573BA378A8EA9AFCFA44BA74B6B34119
77DEE0E3047F88A00E9A8A28A0028A28A0028A28A0028A2B017C42CDE2E6D0BEC528511EEF3C8E338CFE5DA8037E8AAF7D790D859CB7772C56289773
10334EB7992E6DA29E30424A81D7230704668026A28A2800A28AC2F10788E1D070D756772F1119F3635CAFD09ED401BB4552D23504D574BB7BF89191
275DCAADD40CD5DA0028A28A0028A28A0028A28A002B02F7C2961A84F0CD79717D2BC0DBA32D39F90FB7E55BF450056B4B6FB34453CE9A5C9CEE99F7
1FA67D2ACD4125D4114D14324AA924C488D49FBD8EB8A9890A0924003924D002D1595A3EBDA76B325C2584A59ADDB6B8231F97B56AD00145159B7DAD
E97617090DDDEC7148DD14E49FC71D3F1A00D2A2A9DF6A367A7D9FDAEF2E122B7E3F78791CF4E953DBCF15CC093C0E1E291432B0E845004B45145001
4514500145151096332988488645192B9E40F5C5004B451450014531DD571B980C9C0C9C64FA53E800A28A2800A28A2800A28A2800A28A2800A28A28
00A28A68652C5430257A807A5003A8A28A0028A28A0028A28A0028A28A0028A28A0028A28A00CBD62392EE34B048DFCBB8389651D15475FC6B4638D6
289235FBA8A147D053E8A0028A28A0028A28A0028A28A00290804104020F506968A00290804104641A5A2801AAAA8308A147A018A7514500155DECED
6424C96D0B92724B460E6AC5140108B5B700016F10C74F90539618910A2468A87AA85001A928A0082DECEDAD4B1B7B78E22E72C5140CD10DB5B40ECF
05BC51B37DE288013F5C54F450014C8E38E25D91A2A2FA28C0A7D14005566B1B4694CAD6D1173C93B473F5AB34500145145001451450014514500145
145001451450014514500145145001451450014514500145145001451450014514500145145001451450014514500145145001451450014514500145
145001451450014514500145145001451450014514500145145001451450060789F50BCB15B38ED124617136C91A35CB2AFB7A75EBED58FE2CBB9FC3
4F6779677B34AF2CA124867977875FA568F8C755BBB0161696442497B36C3291FEAC0239FAF3FA5739F10A1B183FB3ADB0F35DBCAA5E6762CE147A9E
DD6A5948F4746DE8AD8C640355EFADDEEED5E18EE25B676C7EF23C6E15534D98CDAB6A212632409E52A2E7843B4E71FA569BBAA296760AA3A92702A8
93CDFC396975278E75BB58F529D020F9E5503739040E7F3ADAF1E3DFE9BA49D4AC7529E2742A853236B7BFB1AA9E139221E36F11B99230A5F20E40CE
49EF57FE23A997C1F318FE61BD1B239E334BA15D496FE2BC93C2897F16A5750DCC562AE4A30DAEC1724918EA6A7F084D757FE16825BF99DE6941CC80
E1B07DEA0D62F22B4F01B34840325884403B9298FD293E1EC81BC1766CCC38DFB8E7A7CC68174316C66D6BFE13FBCD1E3D4EE1ED628C39695B715421
7A76CE5B19AD6B9BBB8D23C59A658A5DDC4F6F7AAC19667DD83EA0D53D0248E5F895AD4892061E4050720E795E9F952F89E5897C75E1CCB0243303CF
4E40E681B34BC5FADCDA6C7696565FF1F7A84BE546DFDC19009FD6AAF89263E1BB1826B47BB9AE669C027796DFEB91FE7AD52F8836332EABA46B0A19
A0B59544A00FBA036735D85C5DDB2590BB66478C0DF19FEF1C718F7A0472BF126C56E7C2F26A0F24F1C9088F6C41F0877300723D79FD056C786F4D10
E95A74E2EEEDBFD1D1B63CB94395F4FC6B3BC7AF34BE0590CB1ED924F2CC8A33F2F39FE95720907F606876B14ADE648B6E3F7679DA14127E9C501D0C
BD7E4BED37C53A2C56FA8DCFD9AF2E00785DC91C30E33D7183DEB5FC5DAEB689A7466000DD5CC822841E80F73585E2E9E13E39F0F032A1114997F987
CBF377A9BE24D84973A758EA102990594E19F1CE10E39FD07E740CD3BDB1BBD3B479AF63BFB96BB82333485DF2AFB4648C7E78A9B4ABE1E27F0EC175
04F25AC927DF30B60AB8E08FA7FF005AACEB57301F0DDE4E5C08A5B56DADEBB94E3F9D67F80B4D934BF0A5BC738DB24A4CCC0F6CF4FD00A0463782AE
357D5DF52B7BDBF99ADA098C7BF3FBC63E99EC38ED5A5A16A3730F8BF53D02591E782DE25961790E59410BC67D3E6FD2AB7C3960F1EB0C0AE0DEB636
FA73CD374C913FE16CEAA415F9ACD578EE7F77FE06802FDEDE5CDE78C17466F3E2B3583CC2F171BDBAE09A8A2BFB9D3FC750E8BF6A92E2DAEAD8CBFB
E6DCD1B0DDC03F414EBAD627BDF188F0FC32B410C517992C8870EE719C03DAB219EC7FE169D80B52A8B0DBBABB16C876DADD09F4CD033A2F1D82DE0D
D4C0FF009E43FF0042152F83542F8474B03FE7DD4FE7CD66EBB76F71F0F2F67B97F9E489F048033F39C71F4C56A78478F0A69A320E2051C502E86ACA
ACD13AC6FB1C8215B19C1F5C579E4B69A8DD7C4D8E03A9FEFED2D0B2CC615F9411D97A7F157A3D70BE131F6BF1E788AFC65D508855FD39E9FF008ED0
08ECAD239A2B548EE67FB44AA3E697604DC7E83A564F8A350934FB5B60B71F668EE271149718CF96A41E7EBEF5BD5CE78AEEAD91AC2C7508D1ACEFA5
30CA587238F9483DB9A6223B96D4ACAFF484B6BD6B8B1966D9233619C9209FBDDC75FCBDA8BCD52E2E3C5F16848D2C10887CE69631CB9EB8CF615932
68D77E11BFB5BBD36EE59F4F92758DED64F98AEE38F94D696A1AB5C4FE328741B364B73E47992DC6D0CE475DA3D291423EA975A4F8BED3489AE5AE6D
AF632C85F05E36C9C73DC7152F8FE19E5F0BDDBC376F02C684BAA8FF0058090304FA562EA1F6593E2BE8F1A4865686DCA4993921955C824F73C8AD9F
104D25DF81F559198386128420632A1C81FA0A008FC1B69A8FFC237A549FDA2562F2C3187C8520AE7819EA38AADE2EBED5B47D4F4E9ADF516FB2DDDD
2C6D098D7E5E9DFD2A7B29E587E1D589B694C7318A24420E0E770C81F866A9FC4B60A744258645E83B7B9F7A05D4E8FC43A9FF006558C5316F2D649D
2269719110279623BF4FD6A95D4FA9595E69A20B817565733849242A09507DFD0D4BE2BB8B54B6B5B4BF8964B4BD9D609491F773D083DB9C573B7BA4
6A1E10B98B50D2EF5E6D304CA24B594925431C71EBD680B0CF8A51EA2BA6452BDC46D682E94C71AA61836D38C9FCEBAFB5875936F28B9BC80C8C8BE5
B2458DADDF20D73FF14C03E188972326ED303D786FF1AECA239894FAA8A03A1C8691ABEB53F8A6F347D46E2089A051247B61FF005AB9FAD5C79B5F7F
121B186F2D3ECAB1099DFC9F99416202E33D78ACBF1CABE8FAC699E25B70408A4115CED1F790FAFE1915D1E80C2EADA5D50823EDCFE62E7A88C0C283
F80CFE34018FE33D4F57D1E4B39AC6EA310DC4EB118DA2048CE3A1AD6F123EA5069B2DE699751C4D6D1BC8E8F186120033D7B6306B0FE25102DF4704
E01BF4CD747E21FF00916F54FF00AF497FF403401CFF0085F59D73C47A2998182D9958AF9E533BCFB2F4AD2F09EB536AD05DC57A145DD94C6294A8F9
5BD08FD6AAFC36503C156647566909FF00BEC8FE9553E1F8C6A1E221E979FE34017AEB58BDBEF11CBA268F2450B5B46249E775DFB7A7CA07E22A85CF
8875C86E6E7468A08A6D56D904C1D57E59E3F40BD8F3EBDA8D3ADE66F19DE6B9A45B89ECE78CC32E5C29320232573FEE8FD6B474DD16E5FC5371AFDF
858A468C430C48DBB0B8EA4D006E594934D636F2DD47E5CEF12B4883F8588191F81AC4F1C5DCB0787DADEDDC24D7B2A5B29CE080C707F4FE75D25719
F10FE53A1C9D96FD734305B9D5585AC763616F691001218D5001EC2B9BF06D96B76179A9C5A983F643333C2491C927A8F6C575B5997974D3DCFF0067
D94C8B3E332BE798D7DBDCD0071F6C935E7C5D9AE6CB982DE3DB70E3EE9F936E0FBE7F95759E28B2FB7F87EF2218F31633246C7F8597906A2DD069B7
765A5E9C11249A42F30182C5002496F73C735AD758FB2CD9E9B1B3F950066F85F52FED6F0ED9DE336E91902C9FEF0E0D6C572BF0E2331F84A1CFF148
EC3E99ADDB0D56C3526956C6EA39CC2712053F74D081976BCE3C4526B6DE3CD1E095ED8AEFF32DE304EDE3A93DF3C57A3D711AF64FC4FF000F827E5F
25FBF7C3D0C11D558B6A07CDFED04B75C1023F249391DC9CD7162E9ED3E2B6A221B769A79AD552355E99DAA72C7B0E2BB04BA94EBD2D9FCA625B6593
8EA18B11FCB1F95735A7C41FE2B6A7328004768A0FD4ED1401A96BAEDCC7E204D1F57B68A29A742F04B0B128E0751CFD2AF6B9AC43A3DB46F2299259
A411C310382EC7B573BAF65FE26E8080F485CFFE847FA555F1B8B8FF0084D7C3C6391634DDF2BB8CAABEEFFF00550163A29B58BAD366B71AC430C70D
CC8234921727631CF0D9FA7514BE27D62F344B16BE86D12E2145F9F2F82A73FCAAB78AB41BAD6AC024FA92431C1209C6D870015079CE49EE6A3F1BC8
0F80EE9924F3164893E720E58641CD004F1EB9A8DF69106A3A669EB246D1F9922C8FB49C750BFF00D7AD1D0F56B7D6F4B8EFED77047E0AB75561D453
74A4587C3368B120005AA90A38E76E7F9D607C2C663E12DA7F82E1D47D383FD6981D7CBE6794FE485326D3B7774CF6CFB5731A7F892F26F16FF625E6
9B14132C459E557DD918C8C71D2BABAE10B9FF0085C614F00D9E07BFCB9A0427C4DBED4E0D21ADE1B755B299955E7DFC93D76E3B7FF5ABA0D3EF3543
73696F2E9422B568FE6984A1B185E381D3B0AC7F8A049F0FDB20380F7495B379A8CD63368F6B1468E2EDC46E5B3950141C8A43E86DD62DE6B63FB50E
95A7469717CA9BDD59F6AA0F7F5ADAAE3EEE1D2EC3C66B75671CF71ACCF19CDBABE10023EF124714C469E99AF3CDAB49A46A56EB6D7CA8244557DCB2
AF3CA9FC2A2F1E807C17A9E7FE79A9FF00C785614E6E7FE16C6926EC4424366DF2C44903E593A935ABE27BC5D47E1F5FDD221559212402738C3E3FA5
21973C292C70782F4E9A421238ED1598FA00324D4137899EDECDB52974F97FB39977472AB02C47A95EC0D1E1F7817E1EDAB5D92201627CCC75DBB4E7
1F8572AF2F88341D1CC4601A9E832C5947C65923603D3918CF7E28047A141A842FA526A13E6085A3123799C6D18ACAFF0084A336DF6D8F4CBC7B0EBF
680074F5DBD715CDF8F750DFE04D3BEC259ADAE3602FDF017807F1FE55DCDBC31A68D1C31A8F2C5B85000E31B68025B3BA86F6D63B9B770F14ABB94D
58AE03E13CF2BE9DA840EC4C50CE3CB07A0C839C7E55DFD30661B788A35F3A51693BD942486B94E46475E3AD5CB1D62CAF74AFED38A5DB6B824BB8C6
00AE26DAEF5AD1E292EF4D83FB4F419DDE4543CBC6B9391FCEB6F4FB6D275AF04C96B612C96F67306DE58FCD11DDB883DB8FE54930B16DFC4F125BFD
B1AC6EFEC18CFDA36718F5C7A56DDB5C43756F1CF6F22C9148372B2F422B92B699078367B1D26392F2DE081E33713E111C7392BEBDFF004AB1E04B85
8FC2DA5C12365E6DE23C7230324D00CC5D43C4265F883669358DDAC36513E23111677247DEC7A74FCABB1FED78FF00B3CDEFD8EFB606DBE5F9077FD7
6FA7BD73401FF85BC48FF9F1E7F215D6C17B14F7B7368A1BCCB6D9BC91C7CC32314018FE14B9D0AF4DE4DA25AF92C1C2CE4C7B493CD68CFABDBC378D
6A239E56400C8D1C7B9631FED1AE67E1B000EB98FF009FE3FD696E6F353D0B59BBBEB0B66D4B49BA90B4A211B9E19070DFE7A500F73A98355B1B8B59
2EA19C3C11314770A78231C74F7AE06DE6D3757D63C4769A95CB2433488639954E4004E0671D3A715D7784E7D36E74E9A7D25D8C32CED23230C346C7
1907F2AC6F09A85F1D78980181BD3F99A011A1ABDD68561A35A58EAE924D67E526C76462AD8181CFAFB7BD6FD9A409690ADA2AAC01079614606DC715
CD7C4B38F07CFEF2C63FF1EA5F12EA6FA5F8084D01224782389587F0EE51CFE59A00D89B5DD3E16F9E473103869963631A9F42C062B404B198BCD0EB
B36EEDD9E31EB9F4AE4F4399BFE10BB6B48F4CB999A4B5E781B1CB0C939CFBD2E8B65A8A78166D3F5626DE740F1A316E8BD579FAF145C2C6CAF8874C
C485EE0C6B1F3B9D480C3D41EF5A51BAC91AC8872AC010718E0D79D5EEB8B2692DE1FF00145949673ED0B0CE172848E87DBFFAF5E89063C88F0FBC6D
187CE7771D680682E268EDADE49E53848D4B31F61589E14F32F6DE6D66E0624BD62631FDD8C7007E9478E6530F84AFC8CFCC81723B722B5349852DF4
AB48620022428063E94750E853D6A0D4FEC973369FA9790CB1964430AB0181EA6B1FE1B4F7575A04F75792BCB24D72C4B31CE4F00D7517DFF1E371FF
005C9BF91AE63E18FF00C8A4BFF5F127F3A03A1BFADE9FFDA5A64B6E0ED971BA26CFDD71D0D56F0C6A8754D252494017109314C3D1978CFE35B35C8F
865C41E2EF115928010CAB30F62473FCE80E875D4514531114F3456F1192791510752C71515ADF5ADE2E6DAE22978CE15B903DC751591AD59DFCDAEE
9D751C918D3E0FF5EAEC00E4F5C1FC2B1BC597B08D6F40B8D39C79AD78119D4603AB150413DE95C763B59668A14DF34891AE71976005352E6DE488CA
93C4C8382EAE081F8D1756B05E40D0DCC6B246DD430AE1FE1EE9E92A5FC9312F14570522889F9547AE3BD0077714B1CC9BE191645E99560453C90012
4E00EA6B8FD10B59FC40D5F4E8895B56856658F3C06F97247A7DE3526B17CD7DE31B3D0439584219A503F8B0338A2E162878EADE4B0D04EA1A6EA375
CC815FF7BB83039CD76B66BB2CE04C93B635193DF8AE1FE25580B5F0D6EB363142665F32107E5627BFD78AEE6D0E6D213EB1AFF2A103D89490A09620
01C927B573FA46866C75ABDD4175492E44C0FEE48184C9CE49CF278F6ADD9A28E789A29903A30C153D0D71DE0E845BF8AFC430212228DD02293F7465
A8608D6F146AC2C2DE08219C25CDC4C88A030040CF35BC391C5707E38D2ACFFB634798C2337173B65C13F3720D75F69A7D8D94B9B58846ECB8C6E272
38EC4FD28EA1D0B4658C48232EA1D864293C91F4A92B8C5D4B4D6D62FB4DD70FD9EF8CA7C9998950CA7EEED3DB15AFE16B59ECF4E9A2B894CADF6872
B2139DEBC6083E94058DCA28AE2AEAD5346B2BFD57567924B8B872B0C5E61C0CFDDE87AD311DAD15CCF84B4BB8D3EC9EFB5095DEE275DC5324845EB8
E7BD56D1B52B7F11DE5DCD78EA9040DE5C703498FA9EB48763AFA2B98D0B5545D7EF344FB489D225125BBEFDC769EAB9EF8FE95937D0DCDA78EED2CE
C6EEE2382E23CBAEF2700F5C67B74A2E163BDA2B9CD2ACEFF4ED5B5092E3234DD9987326E3C7E39E99A5D0AFBFE1224B8BC2F2C70C5318E1547C6303
A9C77E68B858E8A8AE6AD7569AC7C4EDA35E4A6549977DBC8DD47B1AE9698828AC2D6B52B5B4BEB6B7D458C76B3A9C3F214BFA1228B286F60D69145C
4B71A7B444C6DBB214FA13DE95C763768AE5ED35B9358F10CBA746B25BDBDBA96660DB5DCE703F0AB31DECBA76BD1E9B732BCB05D26E81DCE4AB0EAA
4F7A2E162E6B57F258456FE4ED324B32A00DDC77AD3AE2FC55A6CB71AFE99E55F4F1F9928C8DDC2F27A5749636135ADC492C9A85C5C2BAE0248785F7
1401A145158975A9B4FAE2693692846085E571C903D05311B74573ED78FA76BD6F626E659E2B90410E77346DDB9AB3A86ABB353874BB538B897966C6
7CB5F5C77A07635E8AE7F55BDBED1163BB9E7FB55A9609202814AE7B8C55FD4E5B8FECD7BBB19C0D911907C808718C8EB45C563468AE7342BCD5757D
292E1DD2D9B270FB0317E7A63B0AB3A06AB35F1BBB6BB50B716B26C6206370F5C52B8EC6D5159FF6B5BABBB8B2866F22780AF3804B0233D0F6A8B4FB
ABDB97BDB7B944865818223A0241C8E1B93CD3111C571AC7FC240D0CB6CBFD9E54ED9011C7A56CD733A66A5AA1F145C69B7B2C72431C7B8109B4FD6B
52F2E0B5F456315D0B79990C99281B70CE3033F8FE9401A555EE2EA2B79214933999F62003BD67E9F75A836B12D9DEAA058E00CACA38739C66B13C44
7531E28D3A3B7BC8D779261531F0BC77A00ECE8AA3A6C7A8C6AFFDA33C32938D9E526DC7AE6A8FF685F4F6ED7965F67951376E80E43601F5F5FD2803
728ACFD12FA5D434E5B99E2113B311B076AD0A0028A28A0028A28A0028A28A0028A28A0028A28A0028A28A0028A28A0028A28A0028A28A0028A28A00
28A28A0028A28A0028A28A0028A28A0028A28A0028A28A0028A28A0028A28A0028A28A0028A28A0028A2A2B89A3B681E79982471A96627B0A008750D
42D74DB47BABC95628907527AFB0AE0B52F1FDFCE5D748B658611D259002E3FE03D3F9D63EB7AD4BAE6A26E6542B6A9FEA63CE463DFEB54C2A3AE188
EBC6DFE1AC5D43A2149752DCDE22D7E66DD2EA77098E7720D83F2000AB967E2BD76DC122533AF1912AE73F8D52B4B3170CCA49DAA7B8EBED57D6CD90
848D480DFC2DD3F3A8E736F6313ACD17C5B6F7AEB05DC7F679DBA127E463FD2BA7AF2568E3951036E8580E377B77AED3C25AB9BA83EC572CDF688865
4B7F12FB56B09DF430AB4B975474D45145686057BCB3B6BE84C3770A4A87B30AAE746D34DA35A359C6D0B7553C93F8F5AD0A2802AD9D95B58C3E55A4
2234CE7039CFE269F756B05E40D05CC4B246DD54D4F4500642F86F4542C534E854B7048C8CD68BDB4125A9B57895A12BB7611C63D2A6A28032A0D034
CB749123B5056452843B161B4F503278A6DBF87348B48E48EDEC9634917632866C11F9F1F5AD7A280B98F17863448583C5A7448C3F88139FCF341F0C
E88CDB9B4F8CB673924E7F9D6C514015E2B4B78AD05B242BE4818087E618FC6ABDBE8FA7DB4892436CAAC9F779242FD013C568514011CD1473C4D14C
8AF1B8C32B0C822AAD8E9561A767EC76A90E460EDF4F4FA55EA280335B41D21A532B69F6C6427258A0CE739ABCB144B1794B1A08F18D814631F4A928
A00A0349B1D91A791B9233944676655FA027153DDDA4177108EE23DE81830192391D3A558A28032A0F0EE8F6DBBC8B08E3DDD769233FAD44DE15D0DA
432FF67C6B21392EA4863F539ADAA28032EF742D36FEED6F2E2DFF00D21703CC4628C47A120F2292EB40D2AE961592CD0080E63D9942BEBC8AD5A280
33EF346D3AF6248AEAD239511762839185F4E2A4B0D36D34E8CC76500891B190093D3EB5728A00A571A95A402E0493207B68FCC9149C10319FE9581F
0F2C25B6D124BCB9044FA84CD3B7D3B7F5AD7BCD074BBFBB5BBBBB45927181B8B11903A6403835A4AAA8A15400A060003000A063EB3EFF0048B3D42E
609EEE3321B7DDB149F979F515A145023263D1E159E379A7B89D2021A28E593215BD7DFB63349A8787B4FD42F63BD991D2EA31812C4DB588F7AD7A28
0316E7C33A55CC30234051A06DC92A1C3E7BE5BBE6B4A3B5B78AD3EC8910106D2BB3A8C1EA3F5AB1450063D8787AC6C1A33079CC9112D147249B9232
7B814DBAF0BE8D793BCF7564B2BBB16259DBA9EA473C56D51405CC897C3FA6C904101858C30486458CB9605BDF39351FF61B3DC45F68BE9A6B589C3A
DBB818C8E993DF1C56DD1405CC4D57C31A76AF7026BF37126D20AA79A42AE3D054D2ADBE83A7DDDE192E2489177B2B396C63D3D2B5698EAB2214750C
AC3041190450072777A9DAF89BC1CAC6D245FB748228A27C64B06EA3D40C75AEAADE1582DE2853EEC68107D00C552B0D0F4DD34A3595AAC45010BF33
1DA0F5C64F15A5401857DE15D2F5199A5BE49E562E5D41998043EAA01E29D79E1E4B9B15B45D42FA38F6946FDF16DEA7B1CD6DD1405CE7F49F0C45A5
5A1B5B5D42F920CE5555C0DA4F5ED4B65E19B5D32F5AF2D2EAF15A46DF32EFCF9A7DF02B7E8A00F34D3743B8F115DDF6A5A66A336956B24E47909C9D
C3A923231CF6AB53E9BAD785EFAC24B5D666BF8EE27113DBCCBD41F4E4D7512F87AC5AE8DC5B3CF692B72E6DA4D9B8FA91535A68D6D6D73F697927B8
9C0C2C97126F2BF4A561DCD3AC2F1969CFA97872E22846678889A2FF00794E7F966B768A6233746BE4D5345B7BA8E4E648C063C655B183F8E6B9AB8F
8776D717AF7526AB7BE6C8773907058FAE6BACB3B1B7B1F385B47E5ACCE6465078C9EB81DAADD2B0EE731A07836D744D49AFD6F2E2E2664299971C66
AC78CAF8D97872E162CFDA2EBFD1E103A967E3F964D6FD51B8D36DEEB50B6BD98BB3DB67CB5DDF2827BE3D680B91695A5C763A041A5F21521F2D883C
E48F988FC49AE67C1BE1F4D1FC51AA7D96EA49ADA38D6362D8E5C9C907E98FD6BB5954BC4E8AE5199480C3A83EB5574AD361D2ECFC8889762C59E46E
AE4F7345845EAC3D73C3E9A9DE59DFC131B7BDB3398E403391E87DAB728A6066D9E9EF6E2E666B8DF7973F7A52BC2E060003D0562DA785AFEDBC412E
AFFDAFFBE9D76CA04230C3D31F8575945017395BFF000BDF5E6B76FAB7F6C94B8B604423ECC30011D0F3EE6B4B5AD0A2D6B4E8ADEEE5613C243C73A2
8055C7F101FD2B628A02E644169A9C91A41A85DC124214ABF968434A3DCE78F7C543E24D16E75AD3CD8417496F6CE30EBE5E4F078C1F4ADDA28039E8
B4BD721D2D6CA3D4ADCED8FCB0ED11CE318A83C2BE1FD4BC3D18B537B04D665D9D944643648EC7F015D451405C2B989BC2B2BF88E4D723D5658EE48D
A83C952117A63DF8AE9E8A00C4F11685FDB9A3ADA49395B88CABC72E30378EE47A1F4A65AE99A94D7D6775ABCF6D2359AB796210C373118DC735BD45
0015C8CFE1ED46D7C5EFAEE9CD0CA258CA3472B11B73DEBAEA28038BBDF0BEA635FB7D7ED2F229AFD321E3981118054AFCB8E78CD5DD4B42D4A5F0DA
68B633DA0478CA4EF286CF273F2E3DF3D6BA7A295877395B1D1B59B6F0B5CE8F3CB672816CD15BB4790791C672314BA7DAEBBA4E949A6ADBC77B188F
6A4A650AC99ECC0F5C7B575345160B9CE49E1782E3C26BA1CEF9DA32241FC2F9CE47B7245353FE1228B4DFECF16504922A794B75E700A5718DDB7AE6
BA5A28B05CC8F0D6891683A5259C6DBDC9DF2C98FBCC7AD6AB6EDA76E37638CF4CD3A8A62391D02D758F0F69696F2DA7DAA37767DB0B8DF112738F42
3E9EB4F9342B8B9F0BEA96D1C315B5D5FBB4A101E1790403E878FD6BABA295877390B5D375AB8F0C268535B4766560F25EE7CD0F9038E147723B9F7A
4F09E99AAD8A5A58EA16C91C5A7799B2657C897767B7E35D8514C2E721A9697A8DAF8DA1D76CAD7ED513C5E4C88AC15978EBCFD2B46DE0D52DD2EEF2
3B688DDDDCC18C2F270880607CC3BD6F5140AE71FE12D1F5AD12EAEC5D4568F0DE4C65768E5394EBD0639AB1A426A5A5DC6A335CDBC9259DC5D3C914
718DD247CFA7706BA8A29587739CD22D65B5B9D5750B7B075FB54A863818F964E0609C1E07AD666956BADE9BE21D5B539F4B2F15E91B523954B0C67F
C6BB6A28B05CE3FC510EAFAFE95269D6DA53444BA3F992CA00201CE3EB57BEC136B7E1B7D3753B436AC2354525837CC07DE15D1514582E72FA04DA8E
91A647A5DFE997123DB2EC8A58006491474C9CF06A5D6A1B96F0CEA4EC1CDC5CA8658D7F809C055FA8E326BA3A28B05CE46E274BFF000F9D2F5AD3E4
5BE306C5431960CF8C02ADF5C56F6876D2DA68B656F367CC8A0446C9E981D2B428A02E65789AD1AFFC3B7F6A840692138CFB73FD29BE16BC17DE1BB1
9B2322208DEC578FE95AF58BA469B36977B7714417EC13379B18CF31B1EAB8F4EF400BA96A7E5C37300B1BE77DAC8A560255C9181835CEF822EA6D13
44365A9699A8C32F9ACE0FD998A9CFB8AEEE8A02E3148650C338233C8C1AE67C2C827D6F5ED442E6392E4451BFF7828C1FE95AFAF49789A5CA34E84C
B70FF2281FC39EF4ED1AC174CD2A0B45C0655CB9F563D4FE7401CCCFE24D620F1DC3A63DB6DB0964F2D4ECFBDFED66BB6AC229737FE235696D9A3B5B
1C98DDD71E63918C8F515BB420671FE2EBB9A0D77498E6120D30BE6721495279EBFA7E7557C6172667D1AF6DA077B1B3BC8DDE55438C64741E8315DC
32AB2ED60083D88A36AEDDB81B718C638C5160B98D67AADB0B4B8D46690A5B4F71FBA628DF30DAA0718CF3B4D73DE03D4ACED20D4219A531B35D3328
653CAFA8AEF28A62383B2D56CEDFE206A1752CB8B69ADC049883B73F2F1FA1FCA9DACDA1B0F1FD86B4FCDACE3CB693A84380066BB5682175DAD146CA
3B150452CB14732149A35910F556008A561DCE33E24DF44FE19315B959F7CABB990EE0A067B8EF5D3691A8D8DEDAC71DA5CC72B471AEE553C8E2AD25
A5B4700812DE25887F00418FCA9B0D8DA5BCAD2DBDB431C8E30CC88012298892E278ADA169A77091AF563DAB87F0EEA7670F8BF5A9259B64774EBE4B
30203727BFE35DED3594363201C1CF34AC3391F1D662BCD0EE9D4FD9E1BBCCAF8E17A633F91ADC86686EB58F320712224182E872B92DD322B4258A39
A331CC8AE8DD558641FC2921822B78FCB82348D3D14605160B9C9D941A7F892D27B0D6A2F3278A565495B8665CF055AADF83ADA7B1FED1D3CCC66B3B
69C0B6727276919C67DB8ADB974FB29902C96D19519C00B8C67AF4A96DADE1B584456F1AC718E8050171EEEB1C6CEE70AA0927D00AE212FEC7C43E26
32DDDE451D9581FDCC323ED12B91D79C5775509B5B73D608BFEF814C119563AB45AA6A9756B6EE5ADE28B0580E1989C120D60785B4EB3D3B59D4F4ED
4ADE16B8797CC85E54043A9F4CD76E91A267622AE7AE062A1B9B2B5BBDA6E60490AF424723F1A4172086DECE0BC096D6B6F19084B32460153918191E
B935CFEA7756ABE3FD3C3C91A797036F72C00079E0D7590C31C1188E145451D9462ABCDA6584E312D95BB7D631405C7CD8B9B19042CAE248D8290720
E47AD72FF0EE392D2CF50D3EE415B886E8B32918C02063EBD2BAD8A2486258A245445185551802A9DD6916379309AE60065181E62B15623D0918CD00
60DC589D47C7D0DC4673159202EC3B37A575D505ADA5BD9C663B68963527240EF53D311CF6B7359DC6AB1E91AA468F6B730E5770C61F27BF6E954F4F
D326F0DEB7041673C92E9D784AF9521C98987A1AE86EF4FB4BC2CD71087629B3764E40EBC536D34BB6B460F1991D94614C8E5B14AC3B9CAC16D0278F
AE22D4906F9A2CC3216DA1F9CE38FE55D2CDA7590B98A636A249C481958B1CAFAB727B7F854D7FA75A6A2AA2EA10E50E55B24153EA08A75969F6F621
BC956CB756762C7F334582E60789A511F88346DF909E6E491C56D37CFAE46048DFBB81894078E48A96FEC2D751844377109101DC392083EA08A5B3B3
86CA329086E4E4B331663F526802D57171446C7E23B3CED85BA8D8C671D7DABB4AA77FA6DAEA08A2E23CB2FDD753865FA1A1821F22C225499A356972
111B68DC33E87D3A9FA66B97BC8A483E235ACED809347853EBF2E31FA1AE96D6C22B62A5649A42A30A6590B6D1ED9A75ED8C17A8AB329DC8728EBC32
1F50698191E365693C3D2448373C922AAAFA9AB90C0D6DE1816EE3E68ED0A91FF01E952A697189124967B89DA3E53CD7C853EBC639A96F2CD6ED423C
D322E082236C06CFAD20297855A36F0F5A98BEEFCDFF00A11AAFA1A37F6D6B1212BB7CEDB81D6ACD96830D85B882DAF2F5221D17CC1C7E95141E1D8E
D75317B6B7D75196606542D91263B1A00A9AA69567AFDEDC086692D6FED08559A323278CE71F89152F87AE35786FAE34AD65D6778904915C28C6F527
1CFBD5D934788DCBDD4134B15CB36EF3383F863B8A9D34FDDE735C4F24B24ABB378F90AAFA0C501731ADF27C7D744723C85FE553EBD6169ACDCAE9F2
48F0DCA45E6C32A360F5208F5ED4E1E15D3D6E85CA4D78938EB22CE413F5AB975A525C5E0BBF35D25450232A07C8467F31CF4A00E7746935ED135A83
4DD5E417969700A41703965C763FA55DD71BFE2ACD15467A9278AD78EC19AE239EF26172D17FAACC61421EE7EB4CD53494BF9229D66782E213F2489C
E3F0A0099AE645D4DE1EB12C1E674E73B88FF3F4AE4E4D32EEFB76B3E1ABE36F3B13E65BB1F949CF4AEB2C6C8DAF9864B879E4931B99FDBDBF1ACFB7
D064D3F0DA6DD98D8FDF0EB956FC3B50037C27A8DC5F69F245796FE45C5B3947006013C9C8ADFAAD656A2D63605F7C8EDB9DC8C16356698828A28A00
28A28A0028A28A0028A28A0028A28A0028A28A0028A28A0028A28A0028A28A0028A28A0028A28A0028A28A0028A28A0028A28A0028A28A0028A28A00
28A28A0028A28A0028A28A0028A28A002B93F88776F1686B6F131537326D247703A8AEB2B8AF894A7FB3EC640BB8ACE71ED91532D8B87C48E1C02A08
F9B6A11B31C0356AC628649115FCC58C72770CE4FD6A38D779118273D4E7B1AD7B381046AA573900631D7D6B9A4CED8AB96522011444FC600519EF4E
4DD226EDCA76F38F7A6795B50AC442F5191CF34D7902B7967A161C8EE6B277341B342A41E0F4E01FE9496171F60D66D6E1490AA70C0B7507B0A6CD39
CB09548C700F707DAABC882E750804409008C64735A536EE44E3EEEA7ABD14515D879C636BDAC4FA3DBB5C2E9D2DD42885DDD1C0DB8F5AB1A1EA4356
D22DB50588C4275DC10B6EC738EB557C62DB7C27A9F5E60238F7E2B17C2177A95D786ACC6996D1430409B4B4F92656079DB8E8339A5D476D0ED68AC4
D0F5A3ACE9F33448915EC0CD1C913E70AE3A67BE2B12EBC51ACD96BF6DA34FA6DB09A7C6C944876367B8F4E878A2E163B6A2B95D6356F106993DB469
6B65726E5F64688581CE3D49ADBD217501679D55E36B86627118E147A5170B17EAAEA17B0E9F652DD5C30548D7273DFD055AAE53C5AC6F757D1B4452
76DC4A66971FDD5E7FC69B1235F46BFBBBEB33717B67F64070532D9DCB8CE7DAA31AD8943BD9D94F710A67F7AB80A48EB8ABBA842F2E977305BFCB23
C2E91E38C12A40AC7F05C7A9DA787922D695626849540700841EB4866A697AADAEAB6DE75ABE4A9DAE87EF21F422924D4922D5E3B092265F363DC921
230C7FBB5C77817ED175E2BD6F51452B6523B05201DAC77718FCAB77C73148BA1FDBEDF227B19566561D40CE1BF439FC280B6A74945416738B9B282E
0631346AFC7B8CD4F4C4145155751BE834DB196F2E9B6C312EE62066802D515CF596A7ABEA1A50D4ADADAD84720DF142E1B7B2FD738CD5AD3B585D63
497BBD3154CCA7698A538DAE3A838A02C67EBBE287D3EE6EA1B6B7594594425B87738033D00F7AD9D23518B55D2EDEFE146559D77056EA39C1FD4570
3A749A93F8EB5985B4DB69E79225F3A077CC60614E79EBD457A1D8A347650A3C11C2C106638FEEA1F41486CB34515C96BBE27BFD0752B482EEC525B7
B9936ACB1124919E817D7914C475B45735ABEBDA86936B1EA1716310B32EA1D7CC3E64609C027B7E15BAD796E965F6D92554B7D9E6176E005C673401
628AC0D3757BDD56D5EFACAD6316796F27CC621E603B8C70338356B4CD623D5B4E96E2C93F7D1EE4685CE0AC83F84FB67BD006AD15C7DB78C67FEDD9
F46BCD2E45BB45F91626DDB8E33F80C11CD5E8B5F9EDF5B834BD5EDA3824BA5DD6EF1B6E0C7B8340EC745451514E65581CDBA2BCA07CAAC7009FAD02
30B4EF1235EF89EF7473612A2DBE7137507071CFD7B57455CCE87E219EFBC477FA4DCD825BCB6C81D9D5F76EE40E78F7A4D6BC58745D421B6BCD3A41
1CEFB526128DA47AF4FD290EC74F45739A8F89A4D3AD96FAE74D9974F2E10CA586E00F46DBE95AD7178EB64975636ED761F6B2AA3052548EBCFE14C4
5DA2B97D1FC5A757D426B3B7D32757B738977381B7F3FA54AFE269D355FECD6D1AECDCEDDE02B29057A6ECFA6681D8E8E8AC8D435A5B4BCB5B186069
EF2E79584301B540E4B1E714967AD097556D2EEEDCDB5D84F3117787122FA83408D8A2B92D4BC67FD937F1DB6ABA5CD6EB2B10B289030201EA38ABB6
7E268AE3C4034796DA486578FCD899981DE319E47638E68B858D77BCB78EF22B479409E552C898392075A926956182495812B1A962075C019AE01B55
D55FE23289B4D626DEDD8470A38CED3FC59AEBA4D46F8697F6B4D2A43382775BB4AAA4019E73DFA74A0762AF86FC4D6FE208AE6486DE58BECEDB5837
CD9FA607B74A2D7C5DA4DDEAD1699134E2EA4C80AD1118C0279FCAA3F076B767ADDA4F35969EB642370AEA00E4FE00561EB13C16DF15F4E9A50462DC
8E17EF3156031EA79C520B1E8158FADF882CF4350F7D1DC796573BD23DCBF4CFAD22EBF147771C17D6D3D9099B6C2F30015CFA7B1F6A83C76FB3C19A
99033FBA03F3602988D3D2B518355D3A1BFB5DFE4CC32BBC60F5C74FC2AED617844470F83F4B390A82D5598938032324FF003A17C4314F0CD75636B3
5CDA5B92249D08EDD768FE2C0A00DDA2ABD95DC17D691DD5AC82486519561DEAC50014C7758E367738550493E805675F6B56F69762CE38E5B9BA2BBC
C310C955F53E94FD2B57B7D504A220E934242CB14830C86802B41E28D1AE6FE2B182F03DCCA70A9B181E99EE2B6AB82D4D00F8B7A5155033031271D7
E56AECAFF50B6D3A112DD49B4336D50392C7D00A06CB74565C3AC42D7CB657114B6D3480B4425007983DBDFDAA0D67C4B67A24805FC376A8C42ACAB1
65189EC0E7AD0236E8AC79FC4363046653E7344807992226563CFA9AB171AADAC30C52279970251B91605DE4AFAE3D2803428AA3A5EA769AADAFDA2C
E4DCB9DAC0F0CA7D08ED45E6A96D6B32DB9DF2DC30DC2189773E3D714017A8AA5A6EA76BA9C7235AC993136C91186191BD08AA8BE23D20EA034F377B
2EC92044F1B29E3EA302803628AC98BC41A53DEAD92DD8F3DCE1159586E3E8091CD4BACEAB068DA6CB7B71B8AA0C285524B36381401792447CEC753B
4ED3839C1F4A7D71BE09D66CC787237B89C2CD34D23B0DA4F258E3B7D2BA6BDBFB7B15433B366438455524B1F6A571D8B94567DB6AF617570B6F0DC0
33B296F2C821801D73E95A14C41451505DDC47696935CCCC163850BB13E8066801D34D1429BE6912341FC4EC00FD6A1B7D42C6EA431DB5E5BCCE3A88
E40C7F4AE43C233CBE2AB9BBD5B5505E08A4096D6E798D3DF1DCF4AE835AD06D350B195628521BA087C99A31B591B1C72290EC6D515C7F81BC52757D
2A717F228B9B31991FFBE9FDEFD2BA08358D32E6531417D03B85DE543F41EB45C2C685154AD355D3EF2568AD6F219645EAAADCFE557698828A8A79E2
B685A6B89163890659D8E00AE43C37E22B1BCD675ABA92F424024458B7BE032E0F2050163B2665519660A3D49C522B2BA86460C0F420E4562F8A6C2D
353F0F5D4922890C56EF2C2EA7A1DB90463AF4159DF0EE48E0F04DA49348B1877909676C64EF3EBEC2819D7514C8DD2440F1B2B29E854E41A65CDCC1
69099AE6558E35EACC681135151C13457110960916443D194E41A92800A28A8D6589DCA2C88CEBD54302450049454734D14119926916341FC4C702B0
3C2DA9BEAF3EAB7426325A0B9F2EDF27A0039C7B1C8A00E8E8AE63C736CE343BAD46DAEAE20B8B78BE5F2DF008CF71F8D5AF06BCD278534F92E64692
478CB166EA72C48FD3140EC6ED145140828A28A0028A28A0028AE175C9F56D2FC5FA35AC3A95C496B7930DF1B609C6EE47D315DD5001451450014514
50014514500145646B571AAD9DACD75609692470C65CA4BB831C727907155BC17AC5D6B9A17DBAF1516432B280A303031401D05145140051557517B9
8F4F9DEC63125C8426353D09AA3E1CB8D56E74CDFACDB082E439007723B122803628AE675DF115C58CD78B6514520B2895E6F333925B90063DB1F9D6
8F87755FEDAD1A0BE3179464CEE4F420D03B1AB451591AEEA92587D9ADED555EEAEA4D91A9EDEA6811AF456568371AA5C4128D5ED56196390AAB29E2
41EB8AD5A0028A28A0028A28A0028AAB0DDC735EDC5B203BA00BB8FF00BD9FF0A3509A5B6B09A6822334A8A4AA0CFCC7F0A00B545666857F71A96989
73756CD6D292414208E9DC66B4E800A28AE6351F16C760F1C8F6ACD6B24A635901E4E3A9C7A50074F45616BBAF49A5496822B392E1273CB283C0FF00
1ADB5395070464743400EA28A2800A28A2800A28A8C4D1998C2181703711E828024A28AAF0DDDBCF2CB1432ABBC58DE179DB9E9401628A8E5916289A
4738440598FA014B1C892C6B246C19186411DE801F451450014514500145145001451450014515465D574F8236792F6DC2AF5C4809FC85005EA2988E
B222BA1CAB0041F514FA0028A28A0028A28A0028A28A0028A28A0028A28A0028A28A0028A28A0028A28A0028A28A0028A28A0028A28A0028A28A0028
A28A0028A28A0028A28A0028A28A0028A28A0028A28A0028A28A0028A28A004ED587E2FD39F53F0F5C43182644C4880772BFE4D6ED149AB8D3B3B9E3
1680B00C842B1C020FAD6B44848077807DCFF2ADBF12786A481A4D47488F73B3EE96103B7723FC2B0AC2E619D82BED470790C306B96A45A3BA9C93D8
B21C60194027F871DEA558F217620041CF342AAE300AEDC9CAFA5421E38332175CFBB633591A8928448CC92B700E46FEA6A5F0B58BDFEB715CA64450
1DECC3B9078155ED20B9D72E16DEC91844A7966FBA83DCE39FC2BD0349D32DF4AB35B6B65E3AB31EAC7D6BA29C7A9CD5AA595917E8A28AE8390C3F1A
7FC8A3A975FF0053DBEA29BE0A554F08E9A0003F75D877C9AD0D5AC5753D327B27728B3000B019E320D672E93A95940969A4DF4715B8CE04A9B9907B
52EA3E873BE0FF003A7F1F788278D765BAC8C8C074243607E3C66B43E2358B49A541AA43C4FA7CA1F70EBB0900FEB8ADFD1F49B4D1AD0C36E39625E5
95BEF48DDC9355175AD375AD37538E162F0C08D1CA58603641E94740BEB72BE93790F886FADB52879B782DC1033C095BA8C7A803F9574B585E11D2D3
49F0FDBDBAA90CD995B3D727A7E981F856ED08185721AA3797F12F4766FBAD6AE83EBF37FF005ABAFAE6BC5B6AE8F61ABC03E7B09833E3BA1EBFE7DE
8608E96B1AE71AD33DAC32E2D236DB3BA9E5CFF747B7AD5EBA8FEDFA6BA4170D109D3E5953A807B8FC2B931F0E6CD570BAB6A2A49C92AE064F73D280
46FDBDC4116AB1695A788847042CD32201F27202FE3D699E329162F096A6CC320C057F3E3FAD37C37E1AB4F0EC738B796599E7605E4948278E82AAF8
C90EA31DA6870E7CCBC9433E3F8635EA680EA5DF0E38B7F0969D25C38555B5476663D06335A1637F69A8C265B2B849A3070590F7A8752D320BFD1A6D
31CB470C91797B97AA81D3F9561F80B4D8B4CB4BE4B698CD6CD3E23908C16C0E4D007595C7FC50131F08BF939C0990C98FEEF3FD715D8541776B0DED
AC96D7281E29576B29EE2989147C36CBFF0008C6967200167103ED8419AE5FE19C771E66AF313FE8D25C9D9EED939FD2B64681A85BDA7D834FD57CAB
2C150AF1067453D406AD2B2D261D3B48FECFB17680052048002771FE2E7BD219CE68FB8FC54D74E3E5FB347CFF00C063FF00035DAD73167E17B9B2D4
2E350875A9CDD5C802477855B2074E2BA3855D624591F7B85019F18DC7B9C76A60C92B87F1880FE32F0CA37CCA67CED3F51CD7715C9EA3E15BED4755
87509F593E6DA9DD6FB600029F7E7DA8041F12DB1E0ABAF7923C7FDF62A0F16C73FF00C2B50B0824ADBC3E660FF080335735EF0D5FEBB62B6575AC01
0EE0CDB6D86491FF0002AD2D3F4EB8874E6B1D4AE12F21D82251E56CF900C73C9C9348087C1C41F08E97B718FB3A8E2B9EF0024C75EF114A8DFE886E
9D54762DB8F3F956D5A691A96996AD63A75DC3F63C9F2BCD525E207B02383D7BD5CD27484D1B497B4B39332B6E732B8FBCE7B9FC68039DD1D037C53D
66471B996D97693DBEE8FE54DF1926EF1B786392374AC323AF55FF001AB969E1DD66DB5D9B571A95B35C4E81245F28ED207FFAA8D4BC35AA6A3ACD9E
A536A30092CDB30A88CE3A83CFE540CDBD65B565810E8C903C993BC4A71C638C51A2BEAEF0BFF6CC36F1C831B7C96CE7D73572D05C0B68C5E346D381
F398C1DA4FB669B7B1DC4B6924769702DE761F24A503ED3F43D6992721E1A25BE23788891F75547EA2A3F8980C93E851052DBAEC9FFD07FC6B474CF0
A5E69FADCFAAFF006C9926B9C79CA6D800FEDF7AA0F10F87359D6F56B5BAFB65A471594A5EDD30D9EA0FCDC75E0521F52C7C48207836EB20F2C838FA
D6BF86F9F0DE99C118B48860FB28AC8F11E8FAE6BBA60D3DE6B1891F064750C7915A7E1FB3D46C6C85AEA13432AC2891C2630470077CFE1F9501D0E5
2F1BFE11BF89715D3395B4D55089093C6EFF00EB1C7E75D4E8ABF699EEF55208FB53058B3DA34E07D32726A9F8C7C3ADE22B4B58E274478270E4B775
3C30CFF9E95D05BC31DB5BC70443091A8551EC280679CDA4F792FC5BBD0A63F3234644121380BB47A7B56EDF68BA94BE27D3F5ABBD42D6016CCB12C7
1AB7CE18E0AE4F739356F58F0F4936AD06B5A5C890EA108DADBC7CB2AE3A1F7F7AB76F6FA85E4E926AAB6E91C2C1D228493971D189F6F4A00E5FE21A
ABF88BC3B1C8328D29C8C67382B5D89D2AC0EAC35336CBF6D09B04B939C74FA5737ACF86F5DD5F56B4BD6BEB288D93978362B67920F3F95751622F45
B81A8181A607AC20ED23F1EF401C9DA1F33E2DDF6481E559AE3DF85FF1AEA2CAF85E5C5E5B988AFD9A4F2C9273BB8CE6B12F344D4AD7C5B26BDA5186
53710F9734329DBD00039FC055F8AC352834B905ACD026A171219257704A827B0C7A71FAD0073FF099366837DEA2F197F255FF001A2FA249FE2DD82C
83223B62EBEC4038ABFE15F0E6ABE1E2F09BEB79AD6594CB22F9677648C71F90A6EA3A0DFDC78B5B5AB62A8D6D0AAC01986253CE411D71CE3B5001F1
30EDF084B20E19268D948EC7775A4D7269EEBE18CB3DC7FAD92CD1DBDFA1CD4DAA69BA8F8905BDA6A16CB65631C824987981DA523A28C741567C4B61
A85F694FA669D15B88268BCB667620A7D07D28030753697FE14EC2F1B32B2DAC20ED38F977283FA54DE0E6D513C136D159D8DBB8911F64866C672C46
4AE3FAD6AE81A5DEC3A1FF00636B10C125BA45E5AB236778F71FD6A2D1F4FD5FC3B6A74FB6822BDB4566681FCDD8C809CE181EBF85001E06B0BAD234
EB8D36FA58DE68A6DFB10E760700FF003CD74CC095215B6923838CE2B3F4AB2920F3EE2ED835CDCB07900E8B8E8A3E95A54C4CF30F8772DDFF006DEB
32F926EA6C852CF26D230C7D7F0AE8ECB4DD5E2F1A4BAB4A96F059DD47E53C424CB120707A7278A48F42BAD0FC4171A9E9112CD6D76BFBFB72DB594F
5CA9ADAB55B8BCB88AEAEEDFC810E7CA8C9CB648C12690CE6754567F8ABA56CFE1B7627E983506B777249F14B4BB56591E2863DC918230CD8639E7FC
F153DFE99E21B8F16DB6B515840A96E9B046D38258735A7AE6873DFDDE9FACD92AC5A959952519B875EE84FE279A0657F145AEB3A84B6535A5B0816C
A6F38BBCCA0903AF43E9517C532A3C2277A824DC2053E879E7F2CD6ECC6FEF6116E6D4DAAC9C4AECEAD81DC0C1E73FD6B1FC75A6EABAE587F66D858A
B207593CF69D546403C63AF7A046A1B54B7F07B5B10084B220F1D4ECE4FE7CD667C3507FE1118246E599D8127AE01C0157AEA4D5E5D00C234A53752C
66268C5C280B95C673F9F1553C216BAAE8FA0A58DEE9F978B71531CA84364F00F3FAD0053F030D9E22F14C4B808B78303D397FF0AA1E11D4AE2EBC5D
AEDD4D6F2DC480F96047B7E450C463E6238E2ADF876DF57D275FD52EAE74995A3D52E15C14707CA1963CFF00DF5FA55F8B44BBD23C5171A9E9B0ACF6
D7ABFBF8770428D9CE4678340C874EB3D5FF00E13EBAD4BEC6D6FA7DCC411F7B0C921460E01EB9154BC5B6905D7C41D0629A30C8EA4B8C7DEC13D7DA
BAB84CF7B7914D2412DBC56ECC555CF3212319C0E98C9FCEB97D663D62E3C65A7EA91E8D726D2C8146C32166C93C819FA5021DF145161F0F594D1054
682ED021518DA36B703D3A0AEAF5AC7F625FE4023ECF2707FDD35C8F8EFF00B4359D2E1B2B3D1AFCED9D6566645C60023030DEF5D5095B56D22E1041
3DB34B1B47B674D8412BD7E9CD0062F80E58AD3C0F6D34CDB2350CCCD8CE39F6A7F8A64BEB7BFB0D4F4B8D2EDAD55CCB6C1BE768DB03728FEB543455
D461F0EAF876E6C674B9463109367EEC47BBEF6EFCEAEEB2DA847E31B0B8B3B669A282D5CCC17ABA96C607B8E0E281F52B6817FA4EB9E265D4ED37C1
7896EC935BC8BB49391C8F5EF5DA57233436FAAF8AB4CBFD3D2449200C6790C65414C700FA9C9A4D6F5DD6ACFC5D63A7DA5934965295DEC133B813C9
CF6C502DCEBEB03C6F1B4BE0ED4C21C110E78F404123F4ADFA8A7863B8B7920946E4950A30F50460D311C6FC266DDE179C63A5DBFF00E82B5DBD719A
0D94DE10BABBB4923966D32671243322963193C156FC856A6A9AD06B4960D2E396E6EA442A8110E0123A93ED4AE37B9C2FC2C42FE24D472B98FC860D
E9CB8E0FEB5A6FA35A5CFC54922F282C31DBACAF1A9C066C7A7A5747E13D013C39A4BEF1E6DD4B9927651924FF0074560DBEA0ADF119F5016B762D64
B711798606C06028193F8BD56C3C55E199ECE3485E6B9F29CA2E3729641838EBD4D7755C078CEE85D6BDA0C904174D1D8DD892793ECEF841B90FA7B1
AED6CAF60BF83CEB472F1E7192A579FA114213D89A58E39A331CA8AE8DD558641FC2B82F87B696B35EF8815A089D12EF0819010065BB57A0D79F7852
75F0EDF6B76FA9C724734970248F084F980E7183F8D0C11D65ECD6F73E1CBD6B56568BECD228DBC01852318ED5CE7C38D3619BC271CB7404E2469151
5C6422E48207D4E4FE35A0F22E97E0EFB35E865B99EDE50115092CCC0F1C77E4551F016A56DA7F8520B6BE77826899CB23C6C0805891DBDE80E83FE1
E4F232EAF665C986D6F196207F841278FD2A9787B5B4D5FC45AA4DAB483CAB7C4504441283248271EBC0FCE97E1F5E429AAEB36EFBD1EE6E8CB0ABA1
0597E639FCAAC6956A3C2BE26D44DDFC9617E77C3363214E49218F6EB40C4D1AF7ECFE3BBAD3AC629069D3C7B97E52115C0C92BF8E6BB8AC13A87F6A
6AB6D169AE5E0B772F71328F97A70A0F7CD6F5084CE47C6DABCF6D3E9DA3D9CBE54DA8CA11DC7554240E3F13567C49A6DBD978767B9B14105C59A799
14ABC3647A9EF59FE3DD3267BBD2B59B78DE53613A9955464EC0C1B3FA55FF0013EA76977E17BA8ED2749A4B98FCB8D1082C49F6A009F4B9AD3C53E1
9B6B8BC896559172E9E8E320FF005AE5BE1D787F4BBED0A6B8BDB6F3655BA740C5D970001D81F735D6785B4B7D27C316D644112842CC0F663CD735E0
8BDB6B0F0EEA365773247751DC499858E18E400303BE48A06743E2D5897C177EB1806216E0260F18E31CD66F85B4A37DE0DB46D4269647783F758720
44BC85C01DF18A97C4724161E0396CEE258E398DAA811BB80C5B8CE07D7356BC297D689E0FD3CB5C423CBB650C0C8382077A05D08BC07A9CF7FE1AF3
AF642ED0C8F1EF3D76A81D7F3AA1A0DFC1E24D4750BABF9D7ECF6D2F97040D260018FBC79A8FE1ACB6D3786AE2C0DC279AD3C836646EDA40E40A83C1
7A7DA6977DA9693ABDB43F6859B7C2D32021D31D89FF003CD0334BC37A9C536B1AAE8067FB4410730B97C9287EF2E7BE3359B3FDA2C7E24DA69D673C
EB6B34797432B100724E39F6AEB2D52C17542B616B6EA523265962880EB8C2E477EA7F0AE627BEB21F14ADE76BB87CB168543EF1B771078CF4A00D4F
0E596AD63E20D5239958696C77DBEF7DDC93DBF0CFE95D4D31195D4323065619041C8229F4C9384F1EBDC45E21F0EBD9902733944246402C40FEB577
C5CF73A35A43AB5A5DC8254991250E72AEA7AF1DB9E71EF54FC653C23C59E1C2F2C7B23B9CBFCC3E5391D6ADFC4C78D7C2A4391933C7B73EC6915D8D
4F11EAD3E9FE1E9350B185A57DAA546DCED07B9159DA82CB17860EB167777705CA422E3F7D2121F81F2B29E3DAAE6B1ACFF647849351813CD222409E
9923A9ACBF110857C0F3DCEA131BAB99A02C8DB8E016C7DD038006450246FE8F7DFDB3A0DBDE61A1371173B4E0A9E871FD2B988AEB5BB7F1D3E8B0DF
C925B1877EE980728BD73EE6B43C217463D2342B4465292DB3B3A819C10720E7F3A8623FF175670481FE8431EFD28193EA77977A2788B49856E649ED
75090C2C929C956C8E41FF00810FCAB57C43AB8D1B4EFB4041248EEB1C684F526B17C68C175DF0C96200FB7753F55AAFF136297EC3A7DCA6ED90DC02
F8E8BD39FD295C2DB1AFAAAEA76BA35D5D093ED6E603BEDD86D00639DA474354FE196EFF00843E12C31995F1F9D6EEA53C6341B9B893FD59B66638F7
5AC2F868ECDE10855971B259147BF3D7F32698BA1D6D15952F8834D87598F4979FFD2A4E800E33E99F5AD5A623335BFB70B0926D3EED6DE48519CEE8
C307C0CE39E9543C29A95FEB1E191773C917DA9DDD436CF9460E0702B67501BB4EBA03A989C7E86B9EF87231E148C1EA2693FF0042A5D47D0C5B76D6
57E226A3122DACB2BDAAF981F211970315DE58A4896712CF1C31C807CCB08C203ED5CAC585F8AD3600F9EC793F9576542064726FF2DBCBDBBF076E7A
67B66BCF752935FBBF1E595A83602E6DA32E840631E0F3CF7CE2BD1AB89D10FDB7E23EAF739245B47E50E3A1040C7E868608ECA3DFE5AF99B77E06ED
BD33DF1516A1769616335D48095897240AB359DAEDCC169A35CCD751096054C4884E32A4807F9D3114AFF51D42CF4D5BD8041751B80C4818D99FC791
FAD4DACEB70E996F6F249856B960A85B951F5FCEB9A9B4BD67C3D6ED75A1DEFDB74CDBE635A4C33F21EB8F6C56DEA9A8E9B75A3D8DD5DDB24CD72CA6
DE294670E78FD33487624D6354BCD184577208EE2CD98249B17694CF43D4F15AD39964B56366F18919728CE095FAD72BE393749E18517170A249275D
D1A0C023AE3D78C57536AC0224207DC890E7EB91FD28038EF094BAD4D7FACB2BDA330B8DAED26E3920718C574FAB5C5F5A696D716C9049245197915B
38381CE3F5EB58DE0A0A2EF5BC0209BBC9CFD0D6CDD4CF71E1FBA7751B9A194617DB2285B03DC4F0EEA526AFA3417D2C6B1BC80E557A531750BA7B46
BF8A3865B500B044CF9847F2CD54F05385F0ADBB30DAABBB27D81EB59274CD5EC47F69784EF567B598F98D653F419E4E3DFDB8FC6803A28B51B9B9D1
E1BFB5B4591E44DCC85F6EDFF1E95CB784E6BCD4206867D2A1BAB68AE988792400C47BF18E6BA9F0E6A516A9A4473C507D9CA9292438C6C61D456578
13021D4C28C2FDADB028027D635DBED1F52B4824B58A4B7BA97623A920A8E383EFCD6D6A17B1585AB5C4D9DA38007527D2B9FF0019EC379A203F7BED
836FE62A1F88571716D676524076813727B678C67DA8B858D7B8D52E2D2D12F6EED945B3637046CB46090013D8F5AD0B8BB8E0B659FEF07C0403F889
E959777A7DF6A9A3B5BCB7E83CF8C6EC4231CF3C1CD426DE38FC336DA7DDDDB190BF931CE38F9C31C1FA7140172F3577B07B78EEED8933B840D1B71C
D6BD70335F6B5A2CD05B788ADD2FF4FF00306CBA4E0A9CF04D7788CAE81D482AC3208EE2840C49241144F237DD452C71E82B2BC3C649EDE7BF989DF7
5296033C051C003F5AB9AB161A5DD95EA226FCB1CD56F0CFFC8BB679EBB0FF0033408B37BA7DA5E8C5CC5B881804120E3F0AE5FC056E6D6FB5787797
0B228049CF6AECEB94F08A6DD5F5A1B8362600E0D1D46B63AB2030208041E08358FA4C86DB51BCD299B72C5896224F215BB7E07F9D6CD73F32ECF1BD
B327064B46DFEE0138FE9408E828A28A60155A4B98E39960C33CCC37045C671EBCE076AB35857D6D6D6DAEC7A919E5372536A5BA73BFB5005CB6D5ED
E7BD6B274921B903708E500161EA304D59BCBA8ECEDDA7943945EBB5726B9AD4D9A5F15E8B24B0F95273C16C9C735D04D34573A6DC3C6DBD363A9EDC
804114864169AED85EDB99ACD9E723AA4684B0E71D2ACD85FDBDFC6CF6EC4EC6DAEAC30CA7DC561F812D6DE0D15A58A35579266DC71CD2E89188FC5B
ACE09CB6D6C76A00DF9EE61B64DF3C81476CF53F41D4D71FE348B4CBBF0FB5C41E4C52891413B7691EC455C92FCCDE3E82D1E4063863628BFED1519F
D68F1F58C73E8DE770195C03FED0A1823A2B018D3AD80ED127F21566AB580C585B0FFA64BFC8559A620A28A2800A28A2800A28A2800A28A2800A28A2
800A28A2800A28A2800A28A2800A28A2800A28A2800A28A2800A28A2800A28A2800A28A2800A28A2800A28A2800A28A2800A28A2800A28A2800A28A2
800A28A2800ACBD4341D37505C4F6CA1BB3A1DAC0FE15A94949AB8D36B639593C15003FB8D42E51474070C7F3E2A587C17A6028D72D35C32F50ED853
F80FF1AE9A8A5C91EC57B49F7228208ADE158A08D638D4602A8C0152D1455101451450014551D5354B3D22D7ED37D288E3CE07724D58B5B88AEEDA3B
881B74722EE53EA28025201182320D6047E15B38667F2A5952D647F31ED978563EE7AE2B663BA825B896DE3915A58B1BD41E466A7A0028A28A002A39
23596368E41B91C1561EA0D494500416B6F0DA5BA416E9B224E15724E3F3A9E8A2800AA3169B6D1EA736A1B59AE2550996390A00E8076ABD450055D4
6D0DED8CD6A26787CD5DBBD3A8A7DA5B45676B15B40BB638976A8A9E8A0028A28A0028A28A0028A28A0028A28A0028A865B88619238E59511E538452
7058FB54D4005145140051451400514514005145140051451400514552D3F55B1D4FCEFB0DD24FE4B6C9369FBA6802ED145140051451400514514005
14514005145140051514D345026E9E548D738DCEC00FD69EAC1943290548C820F06801D45145001451450014514500145145001451450014514D5656
5DCA4107B83400EA28A2800A28A2800A28A2800A28A067BD00145145001451450014514500145145001484038C8E9D2968A0028A28A006ED5DDBB037
6319C738A52030208041E0834B4500323448D76A22AAFA28C0A7D14500155D2CED5241225B42AE3A308C03F9D58A2800AAED676AD3099ADA1328390E
631B81FAD58A2802B5C595A5D106E6D60988E019230D8FCE98BA569AAA5574FB40A7A810AE3F955CA28029C5A6D8412A4D0D9C114899DAD1C6148CF1
DA9D77616978145D5BC72EDE858723F1AB54500436F6F0DAC422B78D6341CE16AB49A36952FF00ACD36D1B1FF4C57FC2AFD1400C4458D151142AA8C2
803000F4A7D1450063CDE18D1279DA7974E85E46392C73D7F3A925F0F69338025B189C01801B247F3AD4A280299D2EC0D89B236B11B63D63C7150596
89A75940D0C506E465DA4484B617D067A0AD3A28032EC343D3B4FB7920B3B5F2E39410C3793807A8073903E9500F0A688B742E8599138FF969E7499F
FD0AB6E8A02E63DE786F4ABD60F7303BBA9CAB199C9539CE473D6AE369D6F269A6C270F34057691236E247D6AE51401971E8F025B9B769AE6584A796
237932A07E03F9D43A3F86AC746903584B74A8093E59972873ED5B54500725AE689A75E789F4F9A288FDB4482494A363083B915D6D5486C6086F26BB
40DE74C0062589E07A55BA00AB7B6697B1795249322F7F29CAE7D8D50B1F0E59E9F6D2DBD94D751249D76CC7E539CE47A1AD9A28B05CE71FC23686FF
00FB423BEBF5BDD9B04C66C9C7E55BB6F1B456E91BCAF2B2AE0BBE32DEE71535140151B50B459E680CEBE6C11F99227755F5AE7BC0568EB6377A94C3
F797D70CE09EEA0F1FD6B46FBC376579A84976649A39265092846C075F4FD2B5A18A3821486150B1A00AAA3B0A4325AA7A9D92EA36135A3B6D594004
E33C6735728A62311349BBB780DA5BDF66D08DBB655DCEABD0807BF14DD5FC3906A1A7DB5B453496C6D0EE85D79C1F715BB4503B9897FA02EABA50B5
D4AE5E59872B32AED287D855CD32C1ECA26F3AE5EE25700348C319C0C0E2AFD1408C23A0343A8DCDDD8DF496C2E71E646A808CE7A83DAADDD69AD3DA
2DAC37735BC41769080127F122B4A8A560B989A4E86FA65BCB6A97F34D6EE855639141D84F7CFF004A34ED26EF4AB7586CEE22910F2CB2A9C03EA315
B74531DCCB4D324B780C76776616766795FCB0DBD8F53ED5474DF0FDE697237D8F54263964DF2AC90839FA57454501739ED63C3B36AB242F2EA92A34
0E1E3DB10F94E6AECDA5B5EE92F63A95CFDA59BFE5A88C211E8702B528A560B9896763AC5B42B6A6F6DE4813E512146F336FE78CD4BAAE91F6DD322B
385C46237560C7A8C77FAD6B5145857312E6CF50BDB43A7DE242F0380AF36EC923D71EB5B1146B14491A7DD450A3E829F4530239A31341244DD1D4A9
FC45627850B416971A74A479B6B2B29E7B1E95BF540E9C83545BE89DA362A5644503127A6690125C1BFF0033FD192D8A63FE5A3B039FC0573963A5EB
FA5EA1753C02CA64B823E52C4631F97AD75D45161DC8E0329850CEAAB211F3053900D64588179E21BBBC00B470208636CF19EF8FF3DEB5AE11E481D2
293CB76180E06715158D9C5636AB0459201C927B9EE6988E77C6D0EA88905EE9F3ED4848054360E73D6BA6B7676B689A4FBE5016E3BE39A82EED5EEE
68B7BAFD9D3E66403976EDCFA55DA430AE6B51B4D4ADFC44BA9D9DB2DCC413694DE011C63BD74B45311CA6A5A76AF26AB63AAA431CAD093BA00C06C1
8F5279AD28A2BDB7D24C46DCCD34DE63480381B0B1CFE3D6B668A561DCE63C38BA9E9D65F63974D63F3960FE60039A900D52D7589EFDB4E691270B1E
C49149503BD7474530B9CB6B1A5CD6FE23B5D6ED61695101132275FAE2A2F11CB7BAC69863D3ECA708AD93E646416FA0AEBA8A560B995A26A2D79008
9ED6685A2400975C035AB4514C4145145001451450014514500145145001451450014514500145145001451450014514500145145001451450014514
500145145001451450014514500145145001451450014514500145145001451450014514500145145001451450014514500145145006078A750D22D2
C923D720796091B200424647BF1EB56EE351B6D3F468AEE2825780A2F971C4993823238AC7F89633E109B8CFEF63FE75B9664A6810120A95B55E0F6F
929751F438EF056BB1BA6A5A85D4170D3DCDD1C88E32C1075033F8FF002AED6EAFA1B4B54B89C380E42AA6DF98B1E831EB5CA7C33709A05F48412BF6
A77CE393C0FF000AB7AF4F26A9A5E9975A733437CF2892DA37FE2207209E945F406B5344788ED16F20B2B986782EA7608B13A8CF3DF39E95B75C0DAE
BD1EA5AF585B6BF612585F5B484C6E7EEB9C74CFD715DF5081A0A86EAE23B5B692E2624471A966C0CF153556D42D16FEC66B591995655DA4AF514C46
40F16E9C74D3A8986F16D07FCB568B00FEB5A5A66A506A960B7B02C8B0BE7699176E47AFD2B8DBC44F136BB068365B868DA6106E5D4FCB291FC39FF3
DCD69F8CEFC5A783AF7FB2A4F2FECEEB6EDE58C797C8047E44521D8D68F5C82613C904334B6D086DF3AA8DB91D40F5AB33EA3147649790C72DCC4F82
3C85DC71EBF4AE7BC3126A09E0EB28ECB4E47DD01C3B4C002C73C9156BC13657BA6E93258EA2F119A2973B51F710080467D28027D3BC53A76A5693DC
DAF9E52DDB6B868F047BF5C01F5ABD65ABD95EE963518A50B6D824BBF1B71EB5C77C3CB28EE2DFC4114DCC734ED0B28E3E52181FE7FA5278EEDE3D0F
C2967A5E9CA63B79EE02392793DF04FF009E945C2DAD8EA53C4768F6CD7420BBFB328CF9FE57C98F5CE6B5619639E1496170F1B8CAB0E8455230C51E
8061650B10B5DA476036F35CD7C2DBB9EE3C38F14C4B2412958C9EC3AE2803B6A28A2988A375A9436D7496BB24926742E1235C90A3BD36CB56B2BEBB
96D209099E100C88548287D0FB8AE7F5C6D42C7C46DA968F10BA74B60B736E588F973C11F952F85F54D3B54D52FF0050B78E4B6BA68879F6CE3E6E3F
8BF90A571D8DB9F59B586F8D9224D3DC2A17648577151EFCF5A9EC352B6D420796DD89F2D8AC8847CC847623D6B88F87F7B777336B37E968D732CD71
CB3481081CE0735ABA25B6AF69E2CBBBABD862B7B6D47EEC2250C4328EBC7B03F9D170B1917BE2582FFC7BA7235BDCA4166AE76B4477B391D40F4E07
E55DCE9BA943A9C0D2C09322AB6D2258CA9AE4D001F181C91C9B1C8FC8575D6B7D1DD5CDCC11A386B660AE5860124678A0196E8A28A6229EA1A8DA69
D0896EE50818E14752C7D00EF4DB4D422B92AA525864604AC7326D620771EB5C92DE0BFF008B1F67705A3B3B72155C701B00E47E75ADA9691A95D78C
74ED462B909656D190C99C1C9CE47BE46290EC6D5FEA16BA744B25D4A1031DAA3A963E805436BAC595C5D7D983B4731E55245DA5BE9EB5C94D72751F
8B30DAC84F95610928A7FBFB739FD6A6F89F33D9E9BA75FDBE52E20BB1B1C7040DA491FA0A02C74DA8EB9A6692EA9A8DDA40CC370DC0F23F014C97C4
1A4C31472C97D18495048AC0123691904E071F8D64F8C196FBE1F5C5DCB12867B64940FEE1254F153784AC2187C1D6A9B4399E0DD21700EE24743EC2
811D0C1345710ACD048B246E32ACA720D547D5AC52492379C031B0473B4ED563D8B6315CDFC2C95A4F096189212E1D57D8601FEB59765AFA68097161
ACDA4971A6CD23986E917709013C83EBFD3145C763B6B8D4F4EFECE6B992E15AD1C98CBAE48EE08E3A564F83AC740B48AEDB4099A652E1647639E719
001C0F5ABDE1BB6B24F0EDBC5672ADCDAB0255F6E370249E47A8CE3F0AC0F85C00D3F5400617EDAD8FC85007552EA9631DCFD9DEE54480E08C1217EA
7A0FC695B51B35B49AECCEBE4C24ABB8C9008FE7D6B906D723F0EEA9A9D9EB168EF657770F22DCC637021BAAB7D338ADBF06DBDB43E1F096932CF6AF
2C8D19FF00649E011EB4018FE1AF11D8EA1AE6B1A84B74B142BB5220E71B9077C7E15D91BA805A8B9322884A870E7A60F4AE13E1E436F71AA7890C90
C6EBF6B1B43202002CF5B7E2A95AE3C3B0DC698A2711DC23AC6063CC0A4E540FC3F4A00D88755B09550C77484BB8450720963DB079ABF5C041AC693E
29D674B7B70D6BA85B5C091E3946372807233DCF4AEFE840D0514563F89F516D27C3B7B7A99F3123C263D4F19FEB4C45A9F55D3ADE7104D790A4A4E3
696E87DFD3F1AB4248CC5E6075298CEE078C7D6B9CF045944BE18B7B89944B3DEA996777192E49E87DBA566F82EF5CEB7AE685362482DE5731A9E405
DC415FA521D88352B8875AF1B5EE9B2DCC3E4C76045BB170516438E7EBC9FCABA9F0C584DA66856F693DDA5D3A67F789F7719E82B8CFEC7B03F154DA
7D923FB31B6F30C78F973B7A8AF44B5B686CEDD2DEDA311C4830AA3B50819353599514B310AA064927000A75473451CD13452A0746182A7BD310D4B9
B7933B2789B68DC70E0E07AD3A39A294131488E075DAC0D79DE8FA4C4DF11B54B68898ACE28F2F0AB70E091C1F6CF6AD2D44C7A3F8FB468EC2310C77
713452A270AC074247A8A43B1DB5422E6DDA6F256788C83F8038DDF95727E3ED666B5FB0E936CE637D42409238EA10900E3F3AB5E25D26D6D3C31753
59C621B8B388CB14ABC382BCF27BD311D3332A8DCC401EA4D2175085CB00A06739E31EB5CFE89711789FC276F3DFA02CEA439E986191B87F3AC3F873
672DD6957AD7D3493DA19DA38E1762578C64E281D8EF410CA194820F208EF4EAE33C2F33DAF8C75CD22263F64876CB1A124EC240C819FAD74D73AAE9
D693793737B0452633B59C038A02C3B50896E2C66B7698C2665281C1C1048ED59BE13D164D0F4C7B696F0DD33CA5F7F61C018FD2A5D66CB4FD634792
49563B88C44ED13839C71D411F4AC5F868A27F0604972CAF2C8AC09EDD290743B2A2BCD22B5BB87E25CBA5D85F5C4368220EC3CC2484DA09033EF5A9
AAA8D0BC5BA3B594B22457AFE4CE8CE4AB0F5E680B1DBD158FAEE8A9AADBB159658AE51088A4472307B64743581E04986ABA1DCC3A8DC5C49716D2B2
4A5A77181D8E338EC6988EDE8AE5FC33A638924D4E7BEBC963691BECF1CB2E5563E40247AD65E8B7F6BE2AD5AFEEB50980B3B6711DBC0EFB54E41CB7
5EBC52B8EC77945717E18D4A39F58D5FC3B34CD3C11126025F398CF05723EB5936F6D7107C489B47B5BCBA1646105D4C84E1768381F8D170B1E95581
A2EA926A5AF6B089206B3B4748A318FE2C7CDF5E7354F44D1B5049356B4D5001613C998152539C7A8C1C8E82B9FF0000E8F6D791EACB2B4EAB1DD6C4
29291C63F53D280B1DDEB56D7379A3DD5BD94DE4DC489847CE307FCF154B4AD375483C331D85E6A27EDA063ED0A3715E72073D78E2AB78C2C635F0AC
CF1CB346F650E6174908618C0E7D7A558F04CB24DE11D3E495DDDDA324B39C93F31EF407431FC2B7DAB49E2ED534ED4AF8DC2D9C602E0601248E71F4
AEDABCF34849A7F89BAE451C8D1C4C9FBC64386C617183F5AD0D5B519BC29ABE9ABF68925D36F1CC6EB336E319C8E41F4E7340342F8FEEB57D2B4EFE
D2D3F5131441D50C3E58EFDF35D4D8BB4B636F2C872CF12B31F5240AE67E2873E0F71FF4DE3FE757B51D623D0FC1F15F3E0BADBA2C6A7F89CA8C0FEB
F8501D0E868AE59C5EA7853FB59AFE7179E48B9FBC366719D9B7A63B56AF87B554D6745B7BF5C6645F9C0ECC3AD170B1A945145310515C7D96BD2F88
7C493E9F632C96F67689BA575E1DDB38C7B55CB3D51ED3C5726837530943C226B676397C7753F913F414AE3B19F77AD6B963E2DD3349B87B7922BB39
251392BCE7E9D2BB4AE23C4C3FE2E2786CE7B30FE75AFE26D7BFB2CDAD95AA86BFBD9047082385C9C6E3401D0515CAEA7797BE1D92D2EAEEFA4BBB69
E5114E1D146C27A15C6303AD4FE2FBBD4B4FD1E6D4F4DBC48D60404C6D1860F938C83F88A2E163A3A2B9AB1BDD5758D061D46D2716CCD08211A20DE6
30EA7D81ED49E1FD6AE7C45A607B722D6789B6CCC57760F3C0068B858E9A8AE6F40D66F27D6750D2351DAD3DA1DC922AE03A7D3D795FCE9FA8EB7236
B50E89A618CDD382D2C8DC8880FE668B858E868AC18B5592CB5C8F48D466591AE23325BCDB426E23AA91D33547C5FE20D4BC3CD6D3C515B4B6B2CA10
821B78E3EBF5A6163ACA2B0B59D52FAC34C6D4E18EDDAD9116468D836FC1C679E9DEAFC77325F6951DD5814569E3578CC80E0679E6810906A51CDAC5
CE9C118496E8AE5B3C1DD57EBCF7C1F26B53EABAE4F9B47984FB2432E40C8F4C76E2BB0BEBF6B0B289EE4C4B34AEB10C13B039EF9EB8A4368D2A2B9F
D4B54BFD1E4B649E15BB4B898462641B36E7FBC39FC3E95D0531051451400554BBBE82D25B78E62775C49E5C600EF56EBCF7C417DAD1F1B6976C60B7
0158BC29BCED7E48C938EB498D2B9DECCE6389DD50BB2A92147563E9589A3789E1D53589F4B36935BDC40859C49EC40FEB5A7A7CB7D28905F5B47015
236147DC18571D6F70B65F14F521E5BC92DC5B2AC48B8C93853C9EC38EB401DF5158B6DADEED5BFB3350B636972C9BE205C3AC83D8E3AD5CD5754B6D
2EDC4B72DF339DB1C63EF3B7A0A622F5158EFACBDB4F6F1EA168D6C97076A49BC300DE8DE94FD77591A35A1BA92D65961519775230B4058D5A2B0EDB
C431DF69AB7B61693DC0D9B9957036FA8CF73ED5774AD52D757B317366FB973B581EAA7B8345C2C5FA28A2802392458A3692460A8A0B313D00159367
E26D2AF7533A743330B8C12A1D0A87FA67AD57F195FA5AE87771324DB9E2CABAA12B9CF427FCF5ACBD061D2F599F4BD425B5BA4BE82052AFB488CE3B
FF009F5A571D8EDA8AA371A8C30DCADB0592598AEEF2E200903DF9A6DA6A905D4ED6C43C372A32D0C830D8F51D88A623428AC8D47C4165A65C2C37E2
784330024319287F11534BAC59C3790DB4CCC867C794E47C8F919E0D0163468A28A0028A8E5952189A49582A28C926B3DF5CB389A3172268164FBAF2
C7853F8D006A514C6915633213F2819C8E78AA36DADE997370D6F15E466541965395C7E7401A345518354B2B8B8FB3C570A65C642105491ED9EB525D
DFDAD960DD4C2207A16071F9D005AA2A947AAD84B19922BB89D477073496FABE9D73288A0BB89DC9C000D005EA28A2800A28A892E219257892456913
EF283C8A00968A822BBB79A578A29A3774196556C9153D00145145001451597E22BB92CF47B89214767D87EE8CE0773401A28E922078D9594F42A720
D3EB07C17219BC336D237562DD7EB8ADEA0028A28A0028AC49B5385F5EB7B68EE6311C68CD27CF80491800F63DAB4E4582F6D9E3DE248DB8251BF1EA
2802C515C8F86A2920F13EA5034CEF1C49B577367F8BF4AEB1595D772B0607B839A018EA28A2800A28A2800A28A2800A28A2800A28A2800A28A2800A
28A2800A28A2800A28A2800A28A2800A28A2800A28A2800A28A2800A28A2800A28A2800A28A2800A28A2800A28A2800A28A2800A28A2800A28A2800A
28A2800A28A2800A28A2800A28A2800A28A2800A28A28039DF15E997FADD93E9D008E380B2B190B72D8EA315A9651DC1B0105D46B19541182ADBB3C6
3357A8A2C17389F0FE9DABE896375A5259193CD7252E048022823AFAD686A9A4DE33E870D9702CE4DEEFDB851D7EA735D3514AC3B9CAEA70CFAF4B6D
6EFA54D6E60995DE79718500F3B4E79AEAA8A29882B335F8F51974B922D20A2DCC9F2EE738DAA7A91EF5A74500727A5697ABE87A08B3D3EDAD9EEDC9
6926793F88F723B9AB3A5F87D97C2F3699A9B8966BB2EF3B8E72CC7AFE181F957474503B9CD6930EADA26989A77D8D6ED22CA412248071938DC0D5A5
8351B3B59AEA0B786E350B99034A864DAA00180013E9C56DD1408E37C1FA66B7A2B5C477963132DD4FE6BC8B38F933D78EF5B5E25D1935CD21ED0955
9410F13B74571D0FF3AD8A281DCE724935793493A60B02B72D179267DC3CBC631BBF2FD6B4741D221D174986C20E4264B37F798F535A545020A28A28
039816DAAC7E2FBBD41223F64F2950A83CC807A7BD4B6BA71B8F141D6840F6EA20F2B6BAED69093D48F6AE8A8A561DCE56CB46BCD0759BDBCD3D05C5
9DE3177B70DB591BAE467AF535AB6B6F7371A88BEBE8562F294A4116EDC573F79891C64F4AD5A28B05CE4F55D36FAD7C6106BD656C6EA3F27C992352
011D79ADAD22DEE2249E6BB5549EE65323229C8518000CFE15A545160B8514514C4731AA6873C3E25B7F10E9A8249954A5C439037AE31907D6B6C17B
C8D77432C015D5B0E172D839C704E3902AE51401CCEA7A24B1F8A2DFC41608249150C771167059718C8F7C7F4A4D734D93C4B2595AC904B15941379B
29906D2C4020003F135D3D140EE733E318EE6E3419F49D3EC2599A68828651F2A8046067F0A7E9B717D63E18B7B7934CB8375143E5EC40A46E038EFD
3A57474520388F87B0DF691A51D3AFF4EBB8A4799A4572A368181C139E3EEFEB56BC390DDDB690D67AE5817B667778CBA6FDAA493865E715D6D14582
E727E1EB59F41D12E3ECB652CCB3DDBC90400F29191F2E73DB8A87C0561A9E9515D5B6A164F1F9D299449B81033DABB2A28B05CE4344496D74FBBB7D
7E02D6971752187CC889C027F8876C9E95A1E16D2CE9515DA223C56D2CE5E089FAA2E3D3B56FD14C2E713E16B0BED0B54D6A292CA561732892091465
1FAE327B75AB1A959EA36BA668B05A46CF3C5721DD5791DC9C9FC4D75D452B05CE2BC4BA7DA6A5796571A642C9AAADCA12E8A5182E7E6DE7DABB5A28
A620AA1AD69C9AB691736321C0993683E87B1ABF4500739E1A69F4DD163D3AF6DA559AD14AAEC8CB09141EA08FF3D2ABF85B43974FBDD4B58BD4DB71
7F23482319251092D83EFD38AEAE8A43B9C0A5DCADF109B5636379F61FB2F96B2184F5E39FA57730C89344B2A64AB0C8C8238FA1A92969898547248B
146D239C2A02C7033C0A928A00E0FC3F7903F8F753BA58E7482EA24489E48C80ED919ED4DF135CC0BE3FD1672CC61B5056690292A84E7009F5AEFA8A
43B9C6F8E3489B528B4FD5B4F5F364B27126D1D5D320F1573C51A8A4FE18B986C7F7F73791F931C4BCB65B83903A60135D3514C2E73DA2590D0BC236
F677322A48B11DE7FDB6C9C7EBFA567FC379557C39346E7698AE242C5B8E09EBFA1AEC69891A46088D154139C018E680B9C27876FEDA5F88DACC892A
98A68D5637E81C8DA303F5AEB2FF0042D2B519FCEBDB1866971B7730E71579228E3FF571A2FF00BAA05494920B9957AB67A4E872428A218163648D14
13C90781587F0C55A3F0B79122B24893392AC30403D3FAD76345311E756D7B683E2E5C4DF684F2E4B708AE1B2A5B6AF19FC0D59F1D5E5AFF006FE811
F9C9BE2B92D21047CA38EBFAD765FD9D639CFD8ADB3EBE52FF00850DA75839CBD95B31F53129FE94877278E449503C6EAE87A329C835E7F7D6777A5F
8EE6B5D390793AD45F3127FD5FF788FD7F3AF408E348902448A883A2A8C014C7B785E6499E246923C8472A0B2E7AE0F6A6208A08E2B75811711AAED0
3DAB81F04D841A46B1AA68DAAC51F98D22BDB79A8089579E549FC2BD12A09AD6DEE3699E08A52BF74BA06C7D334AC3B9969F67875B51650411C5144F
F699510281D3009FC2B99B6BDB38FE2BDE5C3DCC5E54968155F77CBBB0A319FC0D77515B41145E545122C67AA81C1FAFAD30E9D624E4D95B67FEB92F
F8501727560CA19482A464107822B83F01DCC5A6B6BF05E48B13C3725CA39C1C60FF0085778AAAAA15400A060003802A8CDA2E993DD1BA9ACA079C90
4B94E4E3D68039BBA97ECDF0D5E3BC7092CB6EFB55CE19B2C4F1EBC115A1E03BAB79BC2BA7C314C8D245161D037CCBC9EA2B7A5B6B7982896189C28C
2EE4071F4A486D2D6190C905B431BB0C16440091F51405CE27449638BE29EB2AE76F9E9B1091F79804247E59A77C43B67D62FF0049D26D54C92B4A5A
40BCEC438E4FA7435D5C9A358CB1CB1BC21BCD90CACD93BB71EE0F6A7D969765632192DE10242305D892C47D4D1A85CE77E26103C2854C814F9F1FDE
EFCD47E35B19B53F02442CFF007A6111CA42739014838FCFF4AEAAF34FB2BE5DB796D14C3FDB5CD2D958DAD84663B38122563921680B983A1E9BA35F
F86AD7117990185438791B8200C83CF1D2AEF8662862B0956DAD8416DF687F254775E003FCEA57F0FE95248643680162490AECAA4FD01C569471AC68
B1A00AAA00503B014582E3E8A28A623CEFC29606C7C73AC59DEB3A34CBBE2DAEC9BC6ECE46315D0CBA4E9165AEDA5E476CD26A134840769D9982ED39
6E49E00FE75ABA869567A8323DC45FBD8FEE4AA76BAFD0D3ADB4FB7B798CE37C93118F32572CC07A64D2B0EE72FE25DBFF000B03C3993F37CF543C7D
6CCBE2BD06F67252D3CC58DA4071B0EFCE7DBA8FCABA9BCF0BE937D7C2F6E2191A75395612B0DA7DB9AB971A65A5D69A74FB98CCB015DB8762C7EB93
DE8B05CCBD73C3FA75DD93BEA77177343082E15E7E322AAF8B769F8757388D917ECF1E158E48E571CF7AD58343862D81AF6FA78E320AC734DB978E83
18E451AC68169ACAB25E4B746360018D25C2FE54009E1BC2F84B4C29C62CE3231EBB0561FC3260DA2DD90727ED4F915BD63A2416160D656F7576226C
01997254018C038E053743F0F58E8465FB03CFB65E595E4DC33EBF5A00C2D241FF0085A9AC73C7D91703B7F07F85504B3913E2A5DAC934B6C6E21DF0
3A63E6E067AFE35D30F0B5926A52EA115CDE477730DAF2ACBC91E9D3DAAFEA5A4DAEA5E5BCC192684EE8A68CED743EC68B05CC7D4FC3F049A9596A77
B7B3CD716F3C622070072E3B0FC6A97C4CC7F676999FF9FE4FE46BA4B5B192370F77772DD329CA6F002AFA1C0EFEF59DE20F0D7F6F320B8D46E228A3
70E91A2AE0103AD01727F178C7853520300792697C2049F09E9649CFFA327F2A6DF689737DA63D84FAACC6374D8E7CA5CB0A9F43D2DB48B21682EE49
E1400461D40283D334750E860F8131FDA9E24C023FD3BBFE357B527B5D57C37145ABFEED6F24F295E3E046F93B5B93C74AB36FA0FD9355BCBDB5BC96
35BD6DD2C414633EA0F6352EABA35BEA7611593FEEE18E40FB40EB8078EBEFD68039699B5FF0E3C235578F53D244AA048C32F1F3C139E6BB89EEADED
EDFCF9E548E2E3E66381CD664DA5DD5CDB9B1B9B98A5B2CAFF0007EF0A8C719E9D475AAFE30D0E3D6345F25AE5AD841F382064100743401BEACAEA19
1832B0C820E4114FACCF0E5B3D9F87EC6DE524BA440127A9AD3A620AE2BC42BFF1717403CFDC6FE66BB5AC2D73416D4AFEC6FEDAE3C8BAB3625588C8
61E86931A3456F376A925A607C90AC99EE092463F95725A5A093E2C6A92380596D015CF55FB8BFCBF9D7516963344669E6995EEA55DBE604C00074E2
B160F0DEA9078925D6D75280CF2A796E8613B4AF1EFF00EC8A00ABE24FF9281E1FC75C1FCB26B3BC7BE71F19E80225524BA6DDF9DA5849D0D6D6ABE1
DD56FF005BB2D516EED125B3FB8851B69AD4D5F448B59B2856ED847770E1D2788731B8EE33DB3405CA3E25D2357D6F4F5B649ED2029209032EE27233
8ED52F88833F81AE7CC3E6335AAEE63FC478E6B5614BA788C57C216DC082D1B1191F4C7F5ACFF11D85E5E68ADA7E9AB085917CB6DEC46D5E318FCA80
17C1F0476BE14D352340BBA0576C7724649AC8F87BC26ACA0FC82ECED1E9D6B6744B4D4AC3468ECE616E64810246D938207AD51F0DE8BA9E8B3DCAC9
25BCB0DC4BE61209047B0E297603A7A28A2A846478A806F0D6A01BA188D47E0EC7FC227A60072040066A7F105A5D5FE93736968230D2A6DDCEC4541E
1AB7D46C34DB6B0BDB7895608C209164CE71ED4BA8FA18FE11B892FBC4DAE5CCEC19D196318E817FC8A678BE692C7C55A25DDBB1576711BE0E03296E
41AD25D225D2B5FB8D4F4F87CD8EED3134218290C39C8A925D264D535AB6D42FA211C76A331444862C7D4D0333FE2644B3786E2CE722E531F8835A73
787AD351874D7B967DF6689B30783802A1F18699A86AFA58B6B258B72C8B27CCF82D8078E981D7D6B5F4E7B86B44FB5DB7D9E45F9766F0FC01D72281
742E514514C464F8934D3AA692D6E2E45B6183F987A0C7AD657886F21B8F07DDA42924F1C7184F3B01464630C33EF8AD1F15D85CEA3A14B059E7CE04
305071BB1DAA85FDBEA1AA786A5D3A0B436ADE4856F307DE231C0FCA931A347C3D72AFA2E9D1C9203335B2B007A90062B08D8DBCFF00131D9C731C02
4C6782703AD68F86EDAEC25A9BCB46B636B6C21193F78FAD572B7B0F8C9B50FECE9DA078BCA2CA0123DFF4A061E304D9AA68B72876C82E36647D4574
B776B0DE5AC96F3A068E404115CCF89A5B8BA9EC1EDB4FBD7FB2CE2471E49E471D0D74F6D37DA2DD25314916E19D922E187D450B713D8E47C1520D36
EF52D0AE786B7732A337431FFF005B8ADFD274EB6843DD88504970E641C0F941E807E1581E2CD1EE67F10E9F7563298DAE3F73363BA719FE95D8C682
38D517EEA8007D2840C7D14514C4472AB3290A704F19F41DF1EF5CBF8850A6ADA759D9011BCFF2BEDE094EF5B1ACADF96B5360ECAAB2665DA3248FF3
9ACB9C5C2EBF6FAADDC121842B471A2A9253DC8A4C68B30689E478A06A115C208841B3C8FE2F4CFD2B62EEEA1B3B67B8B870B1A0C93505A4664B896F
248DD19C04457EAAA3DBB64D335DB06D474A9ADA3204846533D334019A2FB53BC28F6B75630C529C46BBC33119F7EA6BA15CE06E209EE40C573BA2C8
1160B29B45786480E37840429FEF035D1310AA49CE00CF033420624922451B492305551924D60EADA9DA4FA25DA34CE1DE3242942A40F4FF0026AC5C
EAF1C57C9049148B08F99E4643D7B718FA1A9F587B793489D256F9268885C0E4F1C714C0CDF01107C2F085FBA1DC0E73DEBA4AE3FC037220D39B4B99
244962918AEE5201079AEC29207B8514514C47193693653F8E96292DE3787C82E508E8DC574F616B69648F059A24601CB2A9C9071DEB22ED7EC9E2D8
EF670560684A87C679AD4B10867BBBB42047330209E3EE8C13490D981696DF6BF176A1139FF470B974CFDE39A7184E95E31B486D1996DEEA33BE3C92
01E7FC29BA5DF5A41E27D4A496644493EEB104679A7EA97D6ABE29D3E733A7948843B86C81C1FF001A065EF11EA6D646DE08E4F2DA671B9BA6D5CFAF
E755BC417B059C4BA869F7319B80E37AC7203E62F7C8A678A6CCCEF63AA42A258A2652E31905339CFEB5AA21D36E6067B582D1A5742546C5CE71C502
2F5ACEB736B14EA302450D8F4A9AABD944D059C51391B95403818AB14C41451450014514500145145001451450014514500145145001451450014514
500145145001451450014514500145145001451450014514500145145001451450014514500145145001451450014514500145145001451450014514
5001451450014514500145150C53C331610CA9215386DAC0E3EB4013515C47C47D3E04D166D4D249A2B952AB9472037D456FF8557678674E1B8B6610
724E739E681D8D8A2A2371089444D346243FC05867F2A748E91A1791D554752C702810FA2A04BBB6915D92785953EF10E085FAFA53E29639937C322B
AF4CA9C8A00928A862BAB79D99219E29197EF047048A9A800A2A11710B49E509633272368619E3DA9D24B1C2A0CB22A02700B1039A00928A85EE6DD2
212BCF12A1E0397001FC6A4565650CA415232083C11400EA2B93D675C85FC4DA4E99697796131798A364000743CFF9CD74EB2C5233AC7223B21C3056
04A9F7F4A00968A28A0028A28A0028A28A0028A28A0028A28A0028A28A0028A28A0028A28A0028A28A0028A28A0028A28A0028A28A0028A28A0028A2
8A0028A28A0028A28A0028A28A0028A28A0028A28A0028A28A0028A28A0028A28A0028A28A0028A28A0028A28A0028A28A0028A28A0028A28A0028A2
8A0028A28A0028A28A0028A28A0028A28A0028A28A0028A28A0028A28A0028A28A0028A28A0028A28A0028A28A002AB5D5BADD2AC7231081C3328FE2
03B1F6CE3F2AB34500145145001451450014514500145145001451450014514500145145001451450014514500145145001451450014514500145145
001599AF5A5D5EE9A61B290A49B81386DB91E99AD3A28033AC16EE4F9EFE05478D42C7860C7DCFD7A568D14500145145001451450014514500145145
0014514500145145001451450014514500466189BEF44879CF2A2A36B2B56FBD6B09FAC62AC51400C08823F2C2AECC636E38C7A629B1DBC113978A14
462304AA8152D14005145140051451400514514005145140051451400514514005145140051451400514514005145140051451400514514005145140
051451400514514005145140051451400514514005145140051451400514514005145140051451400514514005145140051451401C5F8D3589A2D5B4
FD12291624BC21A67CE085DD8FCB8A6F89B52B0D0E2B19B4B78CCF1CCA19226C974EE1B15278CF47925D574BD6A084CBF649009D40C9281811C7E7F9
D740D7369E5A1B648649587EED0019E7D7D291462FC446593C137120CE18C6C3F12287D5C68BF0FA0BE451E60B7548875CB1E07F8D45F11AFA28FC37
3583067B99C2950887B30E7F4A7BE9F16BFE044B085F13244B81D0AB8EC4501D0A9A2CBA6C7E13FB4EA1790CD7D731BCCF2160CFB8F618E7D2B43C25
771F88BC2E91EA30F9A636F2DC480FCD8E86A1F0BDDD859E830DB6A51A5B5CD9298E412C783D7B7AD6A5B6A31DBDB5CDFDD9582C8BE215098381C670
067248FD29206723E03D3A3BBBDD62194EFB3B6BA2162CF0E724024F7000A9E48E78BE20CBA258CAD6D617312C93221C70064E3D33D3F1A3E1A4D1A5
C6AB1B2BA1B89F7441948DCA3752A6A1049F134DF26F36CB6DE519361C6EC51A0752C78A2DE1D0350D22EF4B4FB3C92DC79526D3C32E3BD4FF0011B5
A9F4BD261B6B46D935E3952D90308073F9E7F9D52F8817514FA968F6B16F90DBDCF9936C52760E3A9FCEAE78DB4B5F1068F0EA1A6482596CD8BA8519
DE3B8FD2987A8BE26B8D3B49F0EA0B09233716ECBE494219F20F73D6AEEAAD0EAFE079AF2E2DC193EC6F2A061828FB4F3EDCD5A8B59B09B4D8EED941
95900F24AFCF9C7DDC7F915435ABE167E149A1D52706F2EA06091A273C8E14003B668119FE02D1E2BCF0D5BDD6A5FE95E616F295F908B923F3E0D58F
005C4CC757B37FF536B75B621FDD07391F4E053FC11A8DA5AF842DE19E611C96AADE6A3021972C48E3BF5ED54FC097714377ABACEC50DCDEE610CA7E
60738EDF4A00ABA869B6D6FF001434C860B78D23687CC2B8FBC72D927D7A57716B058DADC4D15A2471CAD879117DF3838FCEB9AD7D0D9F8F348D5A64
616AB03C2F2052429F9B19C7D6B774C02E351BCBF507CA902244DFDE00727F3A019AB451453105145140051451400514514005145140051451400514
514005145140051451400514514005145140051451400514514005145140051451400514514005145140051451400514514005145140051451400514
514005145140051451400514514005145140051451400514514005145140051451400514514005145140051451400514514005145140051451400514
514005145140051451400514514005145140051451400514514005145140051451400514514005145140051451400514514005145140051451400514
514005145140051451400514514005145140051451400514514005145140051451400514514005145140051451400514514005145140051451400514
514005145140051451400514514005145140051451400514514005145140051451400514514005145140051451400514514005145140051451400514
51400514514005145140051451400514514005145140051451400514514005145140051451400514514005342A862C14066EA71C9A751400C78A3931
E622B63FBC334E0000001803A0A5A2802296DE099834B0C7211D0B2838A7B2ABAED650C3D08CD3A8A004000000180290285CED0064E4E075A751400C
F2E3DCCDB17730C31C72452A22C6A15142A8E800C0A75140118862127982240FFDEDA33F9D29452C18A82CB9C1C7229F45002600CF1D7AD2D1450021
008208C83D41A5A28A0028A28A0028A28A0028A28A0028A28A0028A28A0028A28A0028A28A0028A28A0028A28A0028A28A0028A28A0028A28A0028A2
8A0028A28A0028A28A0028A28A0028A28A0028A28A0028A28A0028A28A0028A28A0028A28A0028A8E47F2E267DACDB413B54649F603D6B89D7249975
FF000D90CE2F6E66DD24672024790718ED8CE3F03401DD515C66B57725F7C40D274A809D9680DC4D838EDC7F4AECE800A28AC2BCD7AE6CE3BA9A5D26
E161B60CCD233800A8EE3D6803768AE5ACFC61FDA4AF2697A45E5D469C175C019F4AB526BBA8A11FF14F5E30C672AC295C7637E8ACED1F509B5081E4
9F4FB8B2656C049BA91EB5A34C4145606A4F2EA7AA8D2619248A08D77DD4887079FBAA0D72D1DC6A5A9C3AC68B6B7131B7D383B4732B7EF2420E557D
E95C763D228AE63C01FDA0DE1A46D50CDE7191B6F9B9C85C003AFE35D3D31051451400514D560C32A411EA0D50D5E299EC9A482EE4B63086909400EE
C03C1A00E66F2EF56D7B50D620D16F1E15D3C08E308DB44921CE467F0FD6BA7D345C5B68D01D4E40D711C40CCFD7040E7EB8AE2BE1C585DDD5B4FAA1
BF9D37DD92EA3EECB8C6735D2F8DAF8D8785EEE44243C8044A7DD8E3F9669798FC89346BEBDD55C5F23C4B625D9563C658803839FAD6E567E896EB69
A359C08B8090AF1EE464FEB54AE27B9D4EE2786C2E8DB5B5B8224B8500EE7EB807DBBD006ED666BD7F3699A735EC312C8B13032839E13B918EFD2B92
3E24D4F50F0CDDDF5B178FEC4047BD17E695FA6EFA63F535D5E8D2C9AAF86AD64D414192E6DC79A36E3248E78A2F70B58BF6F3C7736F1CF0B068E450
CA7D8D4D5C8FC3CB990E9573A6CCC4C9A7DC3439273C678FE46B6F5DD620D134C92F2E39DBC220EAEDE828B85B534E8AE524D5B51B6B9D29EE48F32F
DC2C968146625E4E7D735D5D310514514005145140051451400514514005145140051451400514514005145140051451400514514005145140051451
400514514005145140051451400514514005145140051451400514514005145140051451400514514005145140051451400514514005145140051451
400514514005145140051451400514514005145140051451400514514005145140051451400514514005145140051451400514514005145140051451
400514514005145140051451400514514005145140051451400514514005145140051451400514514005145140051451400514514005145140051451
400514514005145140051451400514514005145140051451400514514005145140051451400514514005145140051451400514514005145140051451
40051451400514514005145140051451400514514005145140051451400514514005145140051451400514514005145140051451400514514005723F
65DDF121AF2F5982C768AB69BBEE963C360FE27F3AEAA62E2173128670A4A83DCF6AE66C75EBE9B4273736E63D5B73471C0508CB67009E31EFF85034
66781E36D43C4FAE6AD703F7AB298541FE1E4F4FCABBDAC7F0CE90345D223B666DF339324CE7AB39EB5B14030AE6757D27C45A809EDE3D5AD63B4981
56536F9214F6AE9A8A0471DA2F8735ED12CA4B6B5D5AD8ABC9E665E227B73FCAB49ADFC4E08D9A85811FED447FA56FD14AC3B9434A8F538E171AACD6
F2C85B2A6152001E9CD5FA28A62393B6F11BB6A5AEC66D6145B12111870F349C8507D738AC4F095C9D27C33AADDCAF9D467B865F249F9C3F403F526B
BB3A6D91BD3786D63370719908E78E9509D134D3A97F681B54FB564379993D477C74A5663BA1743B79ED747B586EE4679C265CB75C939C7E19C566DE
783B4CBBBC9AEE49AF165998B36C970327D062BA3A29D8573998BC13A642BB63B9BF519DDC4D8E7F2A57F05E9CFC3DD6A0C3D0DC1C7F2AE968A07766
5E8DA2DA68CB2A5A34CC2520B79AFBBA7A7E75178BAECD9F85B529D182B0819549F53C7F5AD59655861791F3B514B1C7A0AE27C55AADBF887C2D0DBE
944C925F5CAC411861C60924E3F0FD6906E6A7C3DB616DE0EB2E49328694E47424D41F12E367F09BB2824473C6EDF4CE3FAD747A7DAAD958416A98C4
31AA6477C0EB49A95945A8E9F3D9CE331CC850FB7BD1D03A89633236976F3E7086057CFA0DA0D73A7C436ADE0FBED42C6D8C2BBE48A2403FD639EE3F
5FCAAE7842E65FECD6D36EC8177A7B985C67A81D1BE9576CF41D3AD277961839662C158E5549CE703A77A00E7F43DDA47C3EB5FB3616EEE46E418077
3B1F43ED5D841E67911F9D8F3368DF8F5C73542C742B1B1984902C876E7623B9658FFDD07A52EBD7A9A7E917131FBE54A46BDD988C00280DCC0F0182
DA87886E14E6296F89461D0F2DFE22A978EACE6D5BC53A1E96CEC96B364B1F71927F1C0AE8FC25A5B693A05BC1303F6871E6CD9EBB9B923F0E957354
D2ADF538E2131749216DF148870C87DA8B6817D4AC6DADEC6F6DD2DA112DD4CFF3492B16654FE2393D3DAB66A9D958C76BB9F7BC933801E590E59B1F
D2AE53105145140051451400514514005145140051451400514514005145140051451400514514005145140051451400514514005145140051451400
514514005145140051451400514514005145140051451400514514005145140051451400514514005145140051451400514514005145140051451400
514514005145140051451400514514005145140051451400514514005145140051451400514514005145140051451400514514005145140051451400
514514005145140051451400514514005145140051451400514514005145140051451400514514005145140051451400514514005145140051451400
514514005145140051451400514514005145140051451400514514005145140051451400514514005145140051451400514514005145140051451400
514514005145140051451400514514005145140051451400514514005145140051451400514514005145140051451400514514005145140051451400
5145140051451400514514005145140051451400D650CA5580208C107BD50B4D0F4BB29C4F6B630C528CE19474CF5AD1A2800A28A280204B7852E1EE
1624134802B381CB01D326A7A28A002A19ED60B93199E2593CB6DC9B86706A6A2800A28A2800A28A2800A28A2800A28A2800A28A2800A28A2800A28A
2800A28A2800A28A2800A28A2800A28A2800A28A2800A28A2800A28A2800A28A2800A28A2800A28A2800A28A2800A28A2800A28A2800A28A2800A28A
2800A28A2800A28A2800A28A2800A28A2800A28A2800A28A2800A28A2800A28A2800A28A2800A28A2800A28A2800A28A2800A28A2800A28A2800A28A
2800A28A2800A28A2800A28A2800A28A2800A28A2800A28A2800A28A2800A28A2800A28A2800A28A2800A28A2800A28A2800A28A2800A28A2800A28A
2800A28A2800A28A2800A28A2800A28A2800A28A2800A28A280296A1A959E9891BDECE22591822E7B9AB80820104107A11581E25974564860D6CB05D
E369D8719FAE39FC2B63CE820B4131754811010C780171C5219628ACEFED8B210472CB2B42B29C2091082DF4F5FC2B4698828A2AB5EDDC767019E557
280E0955CE3DCFB50059A2B26DFC45A6DD588BAB698CAA5B6ED4197CFD2AD586A16DA844CF6CE4EC6DAEA46194FB8A00B94551D4754B4D3553ED327C
F21C471A8CBB9F6151C7ABDBBCEB6F2A4D6B349F71275DBBFE8791401A545676ABACD9E931892F4C8A98CEE54240FC688F56B49AD23B9B72F3A3A86F
DD2EE2A319F9BD3F1A00D1A2ABD95E417D6A9716D2078DFA11FCAAC50056BFBC874FB196F2E091144BB9B0327149A7DEC3A8D945776C498A51952460
D52D7F54B4D32C1E4D42DA596DD970C55032FD0F3C54D637B68DA2457B690BADAF97BD1153E6DBEC050068D158765E2AD2AF566F2A4915E1FBD1BA61
B3E98ABB67AA5B5ECF25BC65D278802F148BB5941E86802FD159D7FAC5869CF8BE9FC918CEE643B7F3C7F9C539754B57B16BC43235BA80438420303D
C67A8A00BF456769BACD96A8CE2C99E444EB26C217F3A886BD66F35C47024D3FD9CED95E24DCAA7D33401AD45531A840DA6FDBA20F2C2577011AE588
FA555D33C43A76AB2DC45692B192D8664564208A00D6A623A3E76B2B6D383839C1F4ACA9BC43689A349AA40935C4319219625CB023AE7D2B33C2FE28
B7D4B4CFB4490C8924970C36C7196EA78271F5C52B8EC741A96A16DA659BDDDE3EC890804FD69C6F6D45A2DD35C46B030DC1CB6011595E28BFD26D6C
8C1ADC523DACC06488F2B9FAFAF15A369159DC6930C7044BF649215288471B48C8E2801F657D697F1B3D9DC473A2B6D6287201AB55C0FC36BA82D74E
BCB576C3B5E398D154B12001D80F635D8E9FA9DA6A224FB2BE5A36DB2230C329F714C18CB9D5F4EB5BA4B5B8BA48E7760AA8D9E49FFF005D685709F1
03E4F10F85A5001FF4A2A78F568EBB4BABAB7B380CF7532451AF5663814089E8ACD8358B29E548D5D93CD38899D4A893FDD3F856950014D665452CC4
2A819249C0029D54F52BAB4B5B2924BF9025B905589048C11ED400EB1BFB5D42DFCFB29D668B715DCBD323AD5AAC3F0D3691068664D1E426C559D8B3
67823AD5FB6D574FBB9FC8B6BC865971BB62B82C07D2802ED154E0D46CEE676822B98DE55EA80F3F87A8AB94005145140053199510BBB055519249C0
029F5C97C49BFBAB2F0C39B51813388E47F453DBF1A00EAD595D4323065619041C8229D5574C0069768074F2531FF7C8AB5400514514005148480092
4003A935045796D323B457313AA1C310E3E5FAD00492C91C31992575445EACC7007E34E560CA19482A46411D0D64F88ADF4FD434F6D3EFEFD2D7CCC1
07CD556FD6AD58476FA76996F08B8530C68156476186F7CD005EA2A149E19233224B1BA8EACAC081F8D3E3923957744EAEBD32A722801F4535885059
88007249ED4C8E7866CF952A498EBB581C5004B4514D4759177232B0F50722801D4514500148481C938A5AF3B7B7BAF19CDAB5C5B5D341F649D62B4E
4ED38EB9FAF1F9D0348F44A2A9893FB3F4B592F26DDE4420CB267EF10393CFA9ACEF0FFDB6EED9753BAB8900B82CE90672AAA71B47E9FAD023768A28
A0028AC7F105EDD69B6F0DF42736F13FFA426DC9287BE7B62B4E19A3B885268583C6EA195877068025A28A2800A28A2800A28A2800A28AC69758927B
F9EC34B8A39A7B75CC8CED8553E9401B3456568BAB26A4B346C047776CE639E2CFDD61DC7B1AD5A0028A28A0028A2B9E4D66EF56BCB9B6D152211DB9
D8F73272377A01401D0D15CCEA7AEDF787F4F797568566FE18E6847CA58F40C3B5745096686367003328240E99C50049451450014514500145145001
45145001451450014514500145145001451450014514500145429711492C9146C19E220381DB355756D5ED3488566BD6758CE7E60B9031EB401A1455
5D3EF61D42C62BCB624C528CAE460E338AB5400514514005145140051451400514514005145140051451400514514005145140051451400514514005
145140051451400514514005145140051451400514514005145140051451400514514005145140051451400514514005145140051451400514514005
14514005145140051451400514514005145140051451400514514005145140051451400514514005145140051451400514514005145140054734D1C1
134B2B6D45EA6A4AC7D6CB3CB676DE6158E7930E07F10C8FF1A4DD90D2B9B145311163454418550001ED4FA620A28A2800A28A2800A28A2800A28A28
00A28A2800A28A28038DF89048D3F4DDA9B98DEA6DFAE0D6C78B21FB47876EA10E11DC284C9E0B6E18159BE2FD3B56D5AE2CA0B1B6022B7984C66675
C1C0E98CE6AFF896D6EAFF004658ADA222E0CC8429C10A41EF8ED487D8E78EB71EAB241A3EBB66FA7DDAB2B47263E5C8F4FC8F7AEFAB9BD6614D5AC7
EC37763235E6015655C8539FBC1BD2BA18D0471AA0CE14003342063EA2B800DB4A1BA1439FCAA5ACBD7E6BE8F4F64D3ACFED324A0A11BB1B723AD311
CC7C2AB458F44B9BB641E6C970C9BBBE001FD6A4F0D163F10BC42338500719E3391CD4FE0BB7D4F46D1A5B4BDD3E62EACD2A61979CE38CE7DAA0D160
D5ED7C4DA96A52695288AF71B50BAE57A753525772BC33BDEFC5B9A39C868ECE2DB12B745F901C8F7C9A9FE29BB41A25A5C46C5258EE46D61D4706B4
751D12487C510788AC23F3245431CF082017182323D4FF0080A8755D32E7C51796B15E5A496B616ADE6BF9982D2B7F740F4A604DE2C667F87F76F390
6436CA58F5F9B8AB1E0BB78E0F0869EAAA06F8433E7B93D4D47E2FB4D42FF459B4BD36CD5C5C205F30C814260E718FC2A7D322D434DF0DC16E2D0497
5026CD824183C75CFF004A04657C317DFE1C9989E7ED726467A74ED5A575AFDD41772DBA6897929562A8CA3871EA38E9599E07D3B54D0ED9EC6F2CB8
9A53319448308081C7B9E3A57654207B9CEF8CD8CBE08D41DD0A334018A9EAA7238AB3E1039F0A697ED6C83F4AABE2CB4D6753B19B4ED3E1B536F711
ED7924721873E94FD161D5349F0F2DBCD6D0CD25AC61628E29797C75C9238A03A181E0A4597C69E22B8755DCB290BC7DDF98E71566EA46FF0085B564
8A485364770F5E1AA4F0A695AAE97ADEA53DED9288EFE532064943797CE71FAFE94DB9D3F591E334D79AC15A08613188D2405981C8FEB401D0788F4D
5D5742BBB22BB99E3253FDE1C8FD6B96F0C6B523F821ACDB1F6DB77FB0A29E7713C2FE40FE95DCC659A35674D8C4025739C1F4AE5F4DF0C2D978C6F3
5058B16AC81E219E379EBC7E74022EDF67C3BE109CD92297B5838C776EE7F5CD6078027BA83C221EDB4F6B979649199F7A80C73800E7935DBDDDB457
9692DB4EBBA2954AB0F635CD68B67A878621934F8AC9EF6CBCC2F0491B80CB9EC47F5A180BE07B4D474F86F2DB540237797CE8E30F9DAA739C7B66B2
7C25096F1DF8A1189D8CCC0FFC09CD7530C7790C773A83DA2C97B2E02C0B280028E8371EFDCD60E83A6EBDA6EB7AA6A371A744FF006F70DB56E00D98
24FF005A00D78B47B6D07C35A84369B9B314B23331E58ED3595F0DA4FB3F81FCF941D91C92BF1D768FF26BADBC83ED3673C19C79B1B267EA315C668B
A7EB163A03F878D81CB3BAFDA4B7C8118F27EB4013FC469D67F0319E30764CD1B0C8E403CD74BA2A797A1D827F76DA31FF008E8AC6F1968979AAF87A
2D334C11E15973E6363E551C56CE922E574E863BD8522991421556DC381D680E8721F0CE08776AD7400F3CDC9427D1724E3FCFA53358BC974CF8A760
B6E36A5EC489281FC5962327F21537846C752F0EC523DDDA4D325EB1924112E5E2604E38F422AE5868D71A9F8B5BC43A8C0D0242823B68587CDDFE63
F99FCE81B297C4027FB7BC35E8B760E3FE04BFE1543E284F32EB3A442E58DA37CC514FDE60C33FA62B4BC4B65ADEAFAEE9F2DB69E56DB4FB82D96914
799C8F9867B60568F8ABC3EDE20B1B69923115EDAB6F8D59863B65491F4140883C5715FEA5A6DBDB693A65CF9B14E8EAEFB50201F539AEA61691E08D
A640921505D41CED38E45670BBB9BA85206B09E27994AC8E46163C8EBEF5A3044B0C11C299DB1A8519EB80314012D54D4D049A5DDA3004185C608F63
56EB3B599EE61B265B4B492E649414014801720F26988C2F8659FF0084361CF5F364FE75996B671CBF16EF954048E2B759191780E76AF51F535B7E05
B0BDD27C36B6BA84062911D9B19DC4827D0566E9D15D9F88D77AA3E9D769693DBAC4B23C4461BE519F61C1A432C7893F73E34F0F4B09DB2C9232391C
6E5C7435D85723E2386F26F14E8D736F65732C36521695D138C11DBD6BAC560CA18670467918A1031D45155EF277B6B732C76F24EC08F923EB4C445A
A5EAD8584B72C32557E4503259BB0C5709E33D5567F03082E0CC6E9A44662F19037672467DAB761BDBBBBF14C26FB49BD86DA34C5B314C8DE4F2CD8E
071507C4299EEB429F4AB6B3B89AE6478CA6C88918041C83F81148674DA436FD1EC9BFBD6F19FF00C7455CAE63C2BABC66C6C34A92DEED2E921C3F99
1101703B9AA9E33F13DFE817F671DA5AAC9048374AECA4F7C63345C2C7654503A514C463F8ABED3FF08D5FFD8F266F2F0A0753C8C8FCB35C8BDFE8BE
21D32DEDF4E90586A51145580FC85F070549EE315D3F8CA4963D0246B7244DE626CC776CE403F9561789B4FD1AFB4692FC225A6AA8A19367CAFE60C1
C6DEF9A4C68B9F12EDA17F0BCD7262433C4CBB24DBF328CF635ABA15B417DE11D321BA89258DAD22CA9FF74560F8E6EE44F08C7A64AB2CDA84D1465C
2213C8EA49FA835B1A06A9656FE1CB512B3C0B6D14713F9A857076FD3A70680E872DF0FB4A5B9BDD5C4EECD6905C14580B1C1393C9F5E2B4F4295AC7
E22EABA4C19168D0ACC1324853853C7FDF554FE1D5FDB4773ACA4F3A234F74648F71C6E1F372296CF53B04F8A57D766E63104966235933C16F938FD0
D0162CEB3A99D4FC7363A147362D62CBCE14FDE6C1E0D5AD624D46D3C5FA3C3A658EDB597FD7C88985209F981C71C019ACDD7EC5B46F1DDA78816126
CA621676419D8DD093FA57517FA947716DF67D2EE164B99C6232873B47724F6E33EF401CCEAFABCB7BE3E8743BB9A3874E870EEBBBFD69DB9C31FC7A
559D4352B3D2BC63A4C7A61431DD8305C470F2A391B4F1DF24D56D6B4F3A578FADB5D9E1DF61703CB95802C236DB8C91E9DEB7753BCB691218B4636D
2DE48EBB1A350DB173CB6474A0674545145324C8F1158ADD69D2CBF689E092DE2774689F6F38EFEBD2B98F865A3C4BA445AB34B3F9B248F84DF84FEE
E48EE6B7BC73742D3C21A83838668F62FD4902A6F085A9B3F0AE9D0B2ED610863C7AF3FD690FA14BC7F398FC366057DA6EA68E1F7C13FF00D6AE8225
4B4B3456202431804FA003FF00AD5CC7C478FF00E24304FCE6DEEE27E3A6338E6BAA8A459624910E55D4303EC69874382F116AB75FF127BE86E248E7
BABADB0C28C4011670723B93EB573C6B7B72FAC68DA369D3BC535C4E1E431B60AAE7AFF33F8526A96C2E3E24E9AB704476B676BE6C79C28DD92303F1
C7E555B4261AD7C48BED45D1912C62D91AB0EE78CFB719A433B5D46D96F34EB9B561912C6C98FA8AC3F005C4B3F85A1498E5EDDDE1E9D94F15D14D20
8A1791BA22963F85733F0EF737871E6208135D4B22FD09A3A8BA1D5556BEBB5B1B56B878E591571F2C6BB98D59A8A7F3BCA6FB3797E6FF000F999DBF
A531190BE26B3240FB35F03EF01E29C3C4967DE0BC1F580D299BC4073B6D2C97078CC84E691E6F11E3E4B5D3F3EF23521914FE2CD3E0902490DE0258
2FFA838E6B7EB9E92EFC50AB84D2ECD9B3D44DC63E848AE86840C86E65105ACB39E91A173F80CD70DF0B26FB545ABDCBE4BC9701893D7041AEEA7856
6B79216FBB22143F42315C37C33B56D32EF5AD2E7204D0CC9C67A8C1E7F9501D08219A4B2F8C72C113623BB4C483D7F77BBF9815E8B5C3DA581BEF8A
375A8E088ACA20A091D5B6ED3FCCD75573AB69F69791DA5CDDC51CF2FDC463C9A1032F514514C46278BEF1ECBC2F7F3C4583F94554AF6CF19AC5F856
A47859DDB9325D39CE724F0BFF00D7ADEF13E9C755F0FDE59A67CC78C94C7F78723FC2B07E17109E1B9AD59764F05D38950F04120633FE7B52EA3E87
47AFE9CBAA689776454334919D99FEF0E47EB55B50D5C68BE1E82E6F17370634411F769081C7E75A57F72B6B6724CCE14853B73DDBB0AE0BE267DA7F
B3F44965076AC999B8C61B0BD7F5A1823A2BA9F5CD3F4E7D4AE6E2094469BE4B610E303D036739FF000ABF15FCBAA68D0DEE8F2424CA0302F9207A8F
AE78A5F11FFC8B97FCE3F70DCFE158BF0D22961F07C46604092577407FBB9EDF91A03A0DF0E789AFF5E7BCB65B48A19E062A6419D89DBF1357F47D66
E66D6EEF45D4962FB4DB20903C678914FB7E22B23E1D28175AFF0018FF004D231F8B53EDC11F16AE8E3AD88FFD96819D1CD7CFFDAF1E9F0BC48ED099
89901248CE3819151D8EA925C6B575A6CD6E237B7891CB86C86CFA7B562788EC9B58D74416174F65AA5944248A40DC4887B63EB4FF000FEA3AAA6B6D
A67882DA2178D0EF86E2303E7407919FC68158D5875B8AF7569B4DD3CA4924281E4909CAAE7B6075EBEB52586AA9717F71A75C011DDDBF257B3A9E8C
B5C5780EDEE97C41AFC4B7260951D410C9BCE32DEBE98AE922D10DAF8B62D565BE79669D1A228530301734058D2BFD592DAFEDAC218FCEBAB82488C3
EDDAA064926A2BDD761D3AE3C9D4623116C79451B70932718EDCD731A7BC927C5DBBFB57052D9844318E38C7E99ABBF11233245A4ECC645E2F53DB8A
02C7557F3496F652CF143E7346BBB66EDB91DF9FA565E93AF36B1A28D46C6CD9CEF28626900231D79C56C5C0CDB4A3FD83FCAB91F85B9FF8462507B5
DC9FC9680E86BF87BC409AE35DA2DB35BBDAB847566C9C9CFB7B558BDD5161D421D36DD04B7732960A5B0AAA3B9FF0AE77C20A13C69E270380655C0F
4E5AA9D8C97B37C4DD43CB112CD14640F3B2414E00C63BE314058EAECB58F3754934CBC83C8BB54DEA036E5917D41FE955B59F145BE8BA8456D7D6B3
2C731C24CA32A7A67F9D53D4F4BBF7F11699AB5C5D428229561F2E353F3024F7354FE238432E87E60057ED78E9F4A02C6F5F6BF1D8C10DD4F6D32DA4
AC079C7036E7B91E95A37570D0D9BCF042D70C172A88797FA564F8DA257F086A0A147CB182BEDC8AB7E197693C37A73BB6E636C993F850072FE14F10
DD5DDC6A9336993C9BA7C91160943D36904D6FF8A713F843502C85775B16DAE395E3BFBD61FC3BCC7FDBB2B03B7ED64F1DF19ADDD5EE05E7836FAE02
ECF32D6438CE71C1A4B61BDCA5E07B8F2FC0B63348AEE1158611771C6F3DAB46D7C49A5DE4F2DBDACD24B3C432D1AC4DBBF2C553F87B8FF842B4FC1C
F0FF00FA1B5636B100F0FF008F6C7558C98EDB503E4CFC7CA18F19FCF0698BA9D445AE5B4BA82D879572972DCEC78F181EA7DAB56A9C516FD467B878
97E555489F1C918C9FD4FE957298828A28A0028A28A0028A28A0028A28A0028A28A0028A28A0028A28A0028A28A0028A28A0028A28A0028A28A0028A
28A0028A28A0028A28A0028A28A0028A28A0028A28A0028A28A0028A28A0028A28A0028A28A0028A28A0028A28A0028A28A0028A28A0028A28A0028A
28A0028A28A0028A28A0028A28A0028A28A0028A28A0028A28A0028A28A0028A28A0028A28A0028A28A0028A28A002B9EF10B86BBB78C85F906F04F5
E4FF00F5ABA1AC1BA40DE23883CE07DD651D718EDED93FCEA27B151DCDEA28A2AC90A28A2800A28A2800A28A2800A28A2800A28A2800A28A2800A28A
2800A28A2800A28A2800A28A2800A28A2800A28A2800A28A2800A28A2800A28A2800A28A2800A28A2800A28A2800A28A2800A28A2800A28A2800A28A
2800A28A2800A28A2800A28A2800A28A2800A28A2800A28A2800AC6D62DA5D4AE6DF4F3130B5C89A6931C1DA785FAE4035B345001451450034A86186
008CE79A63410348246863320E431519FCEA5A2800A6491A4A852445653D558641A7D140102D9DB29CADB420FB4629BF61B3FF009F583D7FD58AB345
003595594AB0041EA08A643043082208638C1E4EC5033F954B450035D15D4ABA8653D411906A282D6DA024C16F1444F04A205CFE553D14005437132D
BDB4B3B8256242EC1464E00CF14E9A58E089A59582A20CB13D85391D64457439561907D4500709AD6A51F8B74CD2EC2C832C97573BE68F3FEAD1339D
C7F235DCA22C71AA20C2A8000F402A382CEDAD99DEDEDE2899CE58A2019AB1490CA3ABE9F1EABA55CD8CB8DB346541FEE9EC7F03547C2B7534DA42DB
5DA95BAB33E44A08C74E01FCAB728A62393D17C476FAA5B5D5DEA56D6F1CF673148931973C03C03CE73E956BC21A4CDA759DC5CDE2EDBBBE98CD22E7
EE83D17F0CD6DFD92DBCEF3BECF0F9A4E77EC1BB3F5A9E9587739EF175F3C5A61D3ED32D797DFBA8917A853F79BF2CD69E8F611E97A55B58467E5810
2E7D4F53FAE6AC98626984C6343228C2B951903D8D4B4C4145145001451450014514500159979A2DA5D5E2DE832C174A36F9D0B6D661E87D6B4E8A00
A961636FA7C063B75203316666392C7D49ACED4FC31A76A7AC5BEA774AE6680600070AD8E99ADCA2800A28A2800AC9B8D1A337CD7F6533DA5CBF1234
632B27FBCBDCD6B51401989A5192759AFAEE4BA2872AA542A0FC054FA9E9F6DAAD84B657B1EF8651823A11E841F5AB9450073E740BA96C96C6EB579A
5B418053CB01997D0B75AD392CDD2C63B6B09FEC823C056081B81DB06AED1405CE6B48F0C4FA4DC5D4B6DABC9BAEA4324A0C2A416FF269ADE17B91AE
36B31EB128BB65087F74BB76F1C62BA7A29587739FBBD0259F5E7D5E2BBF2E758C24000C81C721BD41AB1069B7326AD1EA57F243E6C3198E34841C0C
F5249E6B628A2C17302F7C3C4EAFFDADA55D1B2BC65D937C819261FED0F5F7ABF6769702413DF4EB34C176A854DAA9EB8F7F7AD0A298AE61EA7A17DA
753B7D5ACA5582FE01B4315CAC8BE8DDFA5739E3E82EFECDA5DC5D4AA5C5E28DB182117F3E4E6BBFAE7B5DF0D3EB676CDA9DC470870E91A28F948F7A
4C699D01008C11906B9ED1F44BDD1259EDAC6584D84D29947980F991938CE3B1FF00EB56DDAC6F0DB4714B3199D170646182DEF53D311C7691E1CD5F
45D7AEEEEDAE2D26B6BB705C4A58385049E303AF35A1AAE80F2EB10EB5A64A90DF4636BEFF00BB2AE31835D0D14AC3B9950417F7370926A29024711D
C91C6C4E5BD4D62F8A346D635A96D4C096C8B6739950B31CB8CF1FA015D7D14582E616B169A9EABE1E9ACCC76F14F32E1BE7240E7E9FD6A5F0DDBDF5
9E91059DF431A35BA2C6AC8FB83003F4AD8A28B0AE725A568BA8E8D79A8C36D14735ADE3EF595A4DA63CE7391DFAD6A6A96174FE1F934EB0F28B4911
8899091D46091FA9AD9A28B0EE737E0FB1D4F48D2974DD4A284456FBBCB9524CEE04E791F89FCA9FE30D3E3D63C3732C4C8CEA43C4E39C3038E0D6E4
D1ACD0BC6F9DAEA54E3D0D65E83A54FA45B9B533ACB0EF660483BB9FE5405FA97B4CB67B4D36DEDE46DCF1C6031F53DEADD145310514514005145140
051451400514514005145140051451400514514005145140051451400514514005145140051451400514514005145140051451400514514005145140
051451400514514005145140051451400514514005145140051451400514514005145140051451400514514005145140051451400514514005145140
0514514005145140051451400514514005145140051451400958171E537896206224AE3A7738E0FE1C7E55A3ABDD7D8EC1E4DAC4B7C8BB4E0827BE6B
03488C4BA8C0598FDE240FA0CF5ACE72D523482D1B3AEA28A2B4330A28A2800A28A2800A28A2800A28A2800A28A2800A28A2800A28A2800A28A2800A
28A2800A28A2800A28A2800A28A2800A28A2800A28A2800A28A2800A28A2800A28A2800A28A2800A28A2800A28A2800A28A2800A28A2800A28A2800A
28A2800A28A2800A28A2800A28A2800A28A2800A28A2800A28A2800A28A2800A28A2800A28A2800A28A2800A28A28022B8852E207864194914A9A58A
31144B1AE76A28519F6A928A0028A28A0028A28A0028A28A0028A28A0028A28A0028A28A0028A28A0028A28A0028A28A0028A28A0028A28A0028A28A
0028A28A0028A28A0028A28A0028A28A0028A28A0028A28A0028A28A0028A28A0028A28A0028A28A0028A28A0028A28A0028A28A0028A28A0028A28A
0028A28A0028A28A0028A28A0028A28A0028A28A0028A28A0028A28A0028A28A0028A28A0028A28A0028A28A0028A28A0028A28A0028A28A0028A28A
0028A28A0028A28A0028A28A0028A28A0028A28A0028A28A0028A28A0028A28A0028A28A0028A28A0028A28A0028A28A0028A28A0028A28A0028A28A
0028A28A0028A28A0028A28A0028A28A0028A28A00C4F1317FECE0A881816E79C60E381FE7D2B3BC3EA45F5B9DC70370C7E06B57C499FECD1819FDE0
FE46A8F8763837AB48E7CEC131A9E063B91EA7FCFD307FC4364FF7674B451456E6214514500145145001451450014514500145145001451450014514
500145145001451450014514500145145001451450014514500145145001451450014514500145145001451450014514500145145001451450014514
50014514500145145001451450014514500145145001451450014514500637893585D1EC11C00D3CF208A1527AB1EFF85666A7A86A3E1CB48750D42E
7ED30C92049E3DA008F39E548AC8F8852193C55E1FB572445BF71E78C9603A7E15AFF12937F83AE4EDCED746FA734868E9E291268525898347228656
1DC1E41A96B07C117125D783F4C965186F2767E0A4A8FD00A93C41E21B6D1162478DA7B99CE22810E0B5311B54561E977BAEDD5D235EE9D05BDA3A93
90F971E9919AB1AFEB36FA169725EDC9CE38441D5DBB01401A94570DA96B1E2AD3F4C4D52E52D1219186E8B196841E99F5EB5B5E1B4D71966B9D6EE2
1612E0C51C5CA818EB9A571D87F88F5296D56D6CACDB1777B2F96847551DDAB4AC2DE4B5B45866B992E1C67323F535CCDB1FB7FC4DBCF30164D3AD55
63CF40CD83FC89AE8758D4574CD3A4B92373F0B1A7F798F414017E8AE0B58B9D47C3925A6A0D70F717FA8B794F037DC07B607B647E55774EB6D72DFC
648B334F2589B7CCD2939477C74FC0D1711D8514514C028A2A359A27601654627A00C0D00495CE78B3579B4F8ECAD2D1956EAFEE16056EA5149C1603
F1AE81C3346C11B6B104038CE0FAD7996AF61A85C78FF4FD38EAB2C92471991266404C7C13D3FC6931A3AFD02CB59B4D5351FED1B9F3ACDDF36E59B2
DFFD6E2A5D5752B81AB59699A7BA89A46F3272403B231D7AD68E9F04D6B6890DC5D3DD48B9CCAEA013CFA0AE77C2927DBBC45AFEA058B62716E99ECA
B401D442AE9105964F31C672DB40CF3E952D67EA9A80B0B74210C934CE238A31FC4C7FA5738757D4F4DF100D1AE25FB54B78A1A1936802127AFE007F
2A6163B17DDB1B6637638CF4CD66E8BA9B6A30CCB3A08EE6DE431CA801C03DBAD64F86F57D46EBC43A9E9B740C905A31D93326D24EEC63F9D477B3B6
99F10ACC213E56A70949131C6E5E8695C2C75D453490A0962001C927B572CDAF5C6A66FA6D2E64B7B2B20C1AE1C03E6B0EC33DA98AC6FEA9712DA69D
3DCC1109648977043DF1D7F4A34CBE8B52B086F20FB92AE719E87B8AAFA06A0DAB6896F79345B1A553B948E0F38C8F63589E109BEC9ADEB5A216CADB
CDE6C43D15BFFD6290CEBE8A2A396448A369246DA880B31F402988928ACB5F1168CD8C6A3073D3E6A736BDA4A46D235FC1B54E09DD401A54566AEBFA
4B7DDBF84FE357A19A2B88566824592371956539068024A28AAB7F7F6BA75B99EF2558D074CF73E82802D515CDB78D74642048D3A123203444123D7E
95ABA2EA71EAFA6C57B0A322499C03EC7140588B56D6134F9A0B7485EE2E273F2C487071EB57ADE4696DE391E2689994128DD57DAB9CD0D86A9E29D5
35024325AB0B6887A11D48ADEBFBD8AC2D1EE25E40E157BB13D00A432DD31DB6217DA5B68270A324FD2B9D87C50639A4B5D52CDADEECE0C310393203
D3E86A6D0BC469AA6A579A6CB108EEAD09DC14E5480719FCE8B8AC69699A8C1A9DAF9D0123076BA1EA87D0D5DAE4E063A4F8EDED0102DB5388CA013D
1D7D2BACA636145145020A28A2800AC0BFF1569D65709136F75693CA322FDD0DE99A9F5ED622D2ADA532473FFAA2C2444CA83D3935CCF8165B5D4F49
1657D632C8C6569FCC74F909CF5CD26C691DE8208041C83D08A8ADEE61B957681C38462848F5159BE25BFF00ECED1A43160492E22880E304F156F4AB
5FB1E9D0C27EF0505CFAB1EB4C45DA2B3AFF0057B4B0915276624901B62E7603DDBD0532F35ED3ACDD1669F8660BBD465549E993401A25D43AA12033
67009E4E29F5CEF8B5251A42EA764E44F64C2542A786538DC3E98FE55ADA55EA6A3A65B5E478DB3461B8F5EFFAD005CA28A2800A28A2800A28A2800A
28A2800A28A2800A28A2800A28A2800A28A2800A28A2800A28A2800A28A2800A28A2800A28A2800A28A2800A28A2800A28A2800A28A2800A28A2800A
28A2800A28A2800A28A2800A28A2800A28A2800A28A2800A28A2800A28A2800A28A2800A28A2800A28A2800A28A2800A28A2800A28A2800A28A2800A
28A2800A28A2800A28A2800A28A2800A28A2800A28A2800A28A2800A28A2803375A0CDA7305552370C923A0F51FE7BD64684A9F6E8771C1556DA08CE
7DBFAD6DEB1C6973F5E83A7D456168CE89A8C05B0072A3EA456335EFA66B1F819D5D14515B1905145140051451400514514005145140051451400514
514005145140051451400514514005145140051451400514514005145140051451400514514005145140051451400514514005145140051451400514
514005145140051451400514514005145140051451400514514005145140051451401C878F7477BCB7B4D4EDE3DF71612872146498F3CFE5C1FCEA6F
1513AEF8792CF4C6124978E831DD57A927D315D4D471C31459F2A344CF5DAA06690EE47656B158D9C56B02858E250AA00AE17C6BA604F1359EAD7D04
D71A76D1149E513BA23FDEE2BD0A90804104641EC6982673DE11874B8A0B93A4DC5DCCACE0B9B9CE475C632071D6A878FED24B93A3BEC2D6B15D833F
19006475FD6BAF55555DAA0003B018A4655752AC032918208C8229582FA9C77C44BE5B8D161D2AC584F757D322AA467276839CF1EE05755A7C2F6DA7
5B5BC84178A244623B900034D8B4DB18AE0DC456B12CBFDEDBC8FA7A7E15728038ED309B7F89FABC2C30B736892A93DF6ED07F99AD1F15EB074786C9
96DA39DA6B9541E60C84F7F63EF4DF1059BC5A8D8EBB6C8CF259656645EAD111CFE59CD6C4F05AEA36CAB3C493C2F8701D723D8D00715A9DE26A5F11
B4788E4D8DB6F292E7E4690296E0FB10A3F0AE834F924D47C4377791CCC6D2D9440814FCAEDD58FBE0FF004AD29B4BB09AD12DA4B488C11FDC403017
E98A9EDADA0B484436D12C718E8AB405CA5ABE8F6DABAC4B74F32F94495313EDEBEBF95678F06698176F9D7B8FFAEE6BA3A28B0AE73A7C1FA690A0CD
7B85E83CF3C536CFC1DA7596A115E4135D6E898305793209AE928A2C3B85709A491A8FC53D4EE9384B287C9E7B9E01FD735D8DCDF5ADA1C5C4F1C676
97C3360902B9BF0159C82D2F757B8044BA94ED2A83D4264E3F99FCA80475D5C7780418E7D760704491DFB16CFBE715D8D72E13FB23C6AD2642DB6AA8
011DBCD5E3F33C7E74C11A175AADB43E26B1D2DE02F71344F224BC7EEC60E7F3DB5CE68B22EA9F12F51BE0731DB5B8485BB30E013F9E6BA9BED1EC6F
6F61BC9E36F3E1528ACAC572A7A838ED51CFE1FD3E796293CB68CC69E59113950E9E8D8EB4B500D1EEA7BD9AFA7931F66F3B65BE063200C13EBCD61F
8991A4F1CF86C46B96567638ECA319AEBA08A382148A14091A0C2A8EC2B98D3639356F195C6B047FA25A47F67B66FEF13F78FEA68045FF0019CD3C1E
13D464B5DDE6F9581B464E0900FE84D62784742B4B6F08C536A0F24B1CE9E7C9119309CF23818C9C63AFAD7672C492C4F148BB91D4AB0F507AD64597
876DED3CB8DAE279A085B30C2EDF2A0F4F7A02E6869859B4F819E158495C88D7A28EDFA62B9AD2029F899AD3479DA2DA30C7FDAF96BAAB8992DADE49
E53848D4B13585E12D3E7885F6A979C5C6A5379BB48FBA9FC23F2A00E8EA394031B029BC1072BC7CDEDCD4945311CBEE8219231FF08BCA4F24158C36
33D73504B7DA6F9859BC3775BBEE9C5BE39FA0AEBE8A561DCE45EE34511ACAFE1EBB249C736BC67F135D1696221A7C3E45BB5B458F9626182BCD5CA2
8B005703ACB4BAC7C49B3D2DC32DBDA279AC3B3719CFF4AEFAB0B5AD09AF750B6D4F4FB85B4BFB7CA894A6E0EA47423BD0C118BE3396183C45A02450
C6F334A55D0A8E63381F975AECE28A386311C48A88BD154600FC2B0ED3C3CEDACAEADAADC2DCDD4636C61176AA7BFD6BA0A103392F87BB4E9D7EC7FD
69BD937FE7C7F5AD6D462B49F5AB08AE2EB12A8774B72B91274E7F0ACBD0E21A478B754B03C477B8BA873EBCE40FF3DAACEB9A15D5EEB167AA585DA4
33C0AC844809183C647A1E68E81D4CAC0D47E2970A1E3D3EDB93FDD63FFD7AD7D122B07D7B57B9B5B3F2661288E49B77FAC3DF8FA8A8A2F0D4D61A93
5FE9974AB2CF184B832AE779CE770F435B5A758C7A7DA082325B92CCC7AB31EA6840CE67C54B8F187869D47CDE730FC2BAD937F96DE5EDDF83B7774C
FBD72F3A7F6B78F2DF6730695192E41EB2376FE55D650819CE7F69F88958AB68B19C7195938A41A9F88F70CE8B1EDEFF00BCE6BA4A280B981FDA1E21
326D1A247C752675FD39AD5D3E5BA9AD55EF6DC4131EA81B356A8A62399F8853187C1F7B838326D4FAE4D5EF0A5B1B4F0CE9F091822104FE3CD647C4
0496F23D3749857735DDC648F65C7F8D757146B144B1A0C2A2851F414BA8FA1CBF8D181BFD0207388E4BD1BBF0C7F89AEA5F76C6D98DD8E33D335CDF
8EA16FECBB7BF8D773D85CA4E47FB39C1AD8BF492F748996CA40AF3C5FBB63D391FE140193AD46BA7F847529A7915A6B88CEF7539059B818FCEB1174
2B69FC09651DF5CCD13C8CAF95E4B163C023E9564E89AAEA5E144D2E6416CD6F9C6E3932B0248FC3FC6B6AD23B9D485B7DBAD0DAC36DB4F94C73BDC0
FE42819A4F6A89A535A758D60F2B9EE36E2B0BE1E4ACFE158A36E90CAF1AFD01CFF5AD7D7AE05AE8B77293CF94C171D492303155FC2960DA7787AD20
906242BBDFEA79A3A8BA1B34514531051451400514514005145140051451400514514005145140051451400514514005145140051451400514514005
145140051451400514514005145140051451400514514005145140051451400514514005145140051451400514514005145140051451400514514005
1451400514514005145140051451400514514005145140051451400514514005145140051451400514514005145140051451400514514014F548CCDA
74E8B8CEDCF3EC73FD2B03488C1D460590038CB03EE066BA874591191C6558608AE6B4B318D5D0A12230582EE3C8E0E056535EF23483D19D4514515A
998514514005145140051451400514514005145140051451400514514005145140051451400514514005145140051451400514514005145140051451
400514514005145140051451400514514005145140051451400514514005145140051451400514514005145140051451400514514005145140051451
4005145140051451400514514005145140051451400514514005145140152FB4EB3D41516F6DD260872A1BB1AB2AAA8A15142AA8C00060014EA2800A
8DE24902F9881B690CBB86704743F5A928A0028A28A006B28652AC01046083DE9904115BC42282358D17A2A8C54B45001451450043716F15CC2D0CE8
1E36C654F7C1CD4A0050000001C002968A0028A28A0028A28A0028A28A0028A28A0028A28A00A577A6DBDE5D5ADCC9BD65B67DC8C8719F63EA2AED14
500151B8668D8236D62080D8CE0FAD49450065E85A58D2AD644693CD9A590BC927F7BD3F4AD4A28A0028A28A002AA6A96D2DE69D3DB41398249170B2
0ED56E8A00C0D374ABC7BDB7BFD6240F7169198A20A739C8C163EE6B7E8A280229E18EE20782650F1C8A5594F706AA68B612E9B602D649FCE5463E5B
1ECBD856851400514514019BA969EDA84D6CB23AFD9A36DEF1E0E588E9F8569514500145145001451450014514500145145001451450014514500145
145001451450014514500145145001451450014514500145145001451450014514500145145001451450014514500145145001451450014514500145
145001451450014514500145145001451450014514500145145001451450014514500145145001451450014514500145145001451450014514500145
145001451450014514500145145001451450015CDE99B25D76E77C4F9DC5C678DA73919AE92B96F0DCAD26B377BE52E4A679F4C8A896E8B8ECCEA68A
28AB2028A28A0028A28A0028A28A0028A28A0028A28A0028A28A0028AE66F2EAEAEBC5834991254B1306EDF19C6E6EBC9FD2AAA6A33697E358B47496
496D678C38566DC509CFE43BFE34AE3B1D851585E2E8FF00E2417574B34D1BDB44D22EC6C0240EE3BD6378492F7C43E1F8A7D52EA410AE52254620B1
1FC4C7BFB0A6163B6A2B91F006A535EDA5F5BCE598DADC10ACC49C824FF87EB51FF68C9E21F15DD68F14F24763649FBE31B60C8D91DFD334AE161FAB
DCDFE97AEE95136A4CF15D5C00EAC368C67A575F5E6FE2B8A5B2F15E8704D3C9359993746246C90720119FCABD22840C28A28A62182443218C3A9703
2573C8FC29F5E77A6DE69B7515F417F76F69AC2CAD89246C61BB6D3DEB6E3BE7F0EF80C5DDC0633C51938639DCECC707F12734AE3B1D2B4B1A384791
159BA02704D495C4783F464D4F4A1AAEAC649AE6F199C65CFC8338E2A7F046B53DEDC5F69972FE69B3722394F52A0E307F3145C2C7614514531185E2
CD5DB49D2730737570E21814752C7FFAD567CFFEC6D07CFBF99E536F16E91D8E4B1FAFD4E2B9AD5C9D43E26E9566E498AD23336DCF1B80241FE55D56
B360BAAE937362CE504E9B770EDDE90CCBD224D4B55B44D446A31A24AA5A28634042E4701A97C2FE21FED97BAB69D552EAD1B6B85E8C3A67F315158D
BDBF827C32D1CD70D3ED2597E5E5DC800281F8551F01E8D75633DDDF6A7B22BABDF9D61CFCE1339C9FC4D006CDFEA3269BAFDAA4F231B4BDFDD8CF44
71D31F5ADBAE5BE232E3C272DC2F125BCB1BA37707701FD6B7749B9379A4DA5CB72D2C28C4FB91CD01D0BB55EF6D85E5A496E6478C4831BD0E08FA55
8A2988E0FC177778DE2CD66C2E6F26B84B51E5AF98739DAD8CD7795E5DE1D867BDF885AEDBAC8F1446691A628704A87C6DFC735B96171369FF0010DF
47824736725BF9855DB710719E29219DAD471CD14B9F2A447C75DAC0E2B8FF00196AD336B1A7787ACE5F29EF594CD20EAAA4E00FD0D2F8A2DE1F0DE8
8353D283433C322027713E60CF20834C2C69F8C2D5DF479AF61B99609AD636752A786F6228F035CBDD7846C26958B4855B731EE771A8B52D463D57C0
171A827CAB35A9623D08E08FCC1150780264B7F02DBCD2B623895DD89EC0124D20E874B7513CF6ED1C72B44CD8C3AF51CD727A46A37FFF0009F5FE91
757B24B0C10068D580F9890A7D3DE93C3334FE29377AADCDC4F1C0B318EDA38DB6ECC639FAF4AC7F0FC7703E2B5EA5CCAD249142C37B705D405033F8
6281D8F4CA2B98F136A3159EA96106A2CD169B306F31C12017C8C02476A9AD1669BC4693DB5CF99A62DA91188DC9556F9460FE1CFE14C93A1A28A280
0A28A2800A28A2800A28A2800A28A2800A28A2800A28A2800A28A2800A28A2800A28A2800A28A2800A28A2800A28A2800A28A2800A28A2800A28A280
0A28A2800A28A2800A28A2800A28A2800A28A2800A28A2800A28A2800A28A2800A28A2800A28A2800A28A2800A28A2800A28A2800A28A2800A28A280
0A28A2800A28A2800A28A2800A28A2800A28A2800A28A2800A28A2800A28A2800A28A2800A28A2800A28A2800A28A2800A28A2800A28A2800A28A280
0A28A2800A28A2800A28A2800A28A2800A28A2800A28A2800A28A2800A28A2800A28A2800A28A2800A28A2800A28A2800A28A2800A28A2800A28A280
0A28A2800A28A2800A28A2800A28A2800A28A2800A28A2800A28A2800A28A2800A28A2800A28A2800A28A2800A28A2800A28A2800A28A28033B59902
59A8640E8F2286527191D7FA573FE199656D6A7C6420257017F871FE38AD7F12C64DA45206C157DA07D7FF00D55434149535105118A952B23E38E991
CFD4564DFBE6B1F80EA68A28AD4C828A28A0028A28A0028A28A0028A28A0028A28A0028A28A00E5AE753BABCF187F6245235BDBC31892474E1A4E338
07B0AC1492DEDFE2A5B436A8238A284C2CEE4FCCDB493C9EFCE2BBE3656A6F85EF92BF680BB7CCEF8F4A8EE34CB1B99E39A7B489E48DB72B15E41F5F
7A561DCC3D6EE377812FA5B89066449065CE3397381F954BE0831DBF8334F32B2C7B50EF2C40C12C7AFE62B5EF34CB1BE58D2EED92558B3B148E07E1
4C3A369C6C9ED0DA47E43905939E48E9CF5A02E731F0E0AE7595CF3F6B38F71CF355FC27653697E3CD662BA017ED21A484FF007C17CF1F8575F6BA2E
9B66EAF6D671C6CA7208CF5F5A96EAC2DAEC869E305D410AE0E1941F7A2C17384F1C5D413F8C344B50DB8C4FFBC23A292C38FD2BD14104641C83596F
A0693204F32CA36D8495249C83F5CD5DB4B482CE1F2ADA211A6738049E7F1A01962A0BC629673B062A5636391DB8A9E98E8B246C8E32AC0823D41A62
382D2B44D27C4BA042FA89DB7D0862F283B588C9C13EA314E8F4FBDD67C0379A68713C90C9B6DE4FF9E814835D74DA469D3C290C96919445DA00C8E3
D322AE410C56F0AC3046B1C6830AAA3000A561DCE6BC197D6F6FE14B78AE6E63492D43ACC19B0570C4FF00222B23C01A74CFACEA5AD3214826665841
CFCD939CFE9FAD76175A369D78FE65CD9C6EE7A9C633F5C75ABB1A2C6812355555180AA3005160B8FA28A2988E22EDBECBF152D6490616E2DB629C75
3823F9D75D7D796F616925D5D48B1C518C924FE83DEB07C69A74D35BDB6A76485AE74F904985EAC9DC0AD554B2D774A85EEAD84904A03F9720E87DE9
0CE42CFC5FA0DE5F9BFD48C9E6E088622BB9635F5FA9AD3F0D6AABE23D7EF35181596D2D631045BB19624925B8AD85F0D684A085D26CC03E910AB967
636961118ACADA28109C958D428268B05D1CEFC4B9563F06DCAB67F78E8ABF5CE7FA56D787EDA4B4D02C2DE5E5E3B740DF5C722B9ED723FF008493C4
B6BA52296B3B16F36E9BB16FEEFE9FAD764000300600A03A0B51CD3456F134B3BAC71A8CB331C0152541756D05E5BBDBDCC62489F8643D0F7A62383F
045E5AB78B7C47399A30934BB918B0C11B989C1FC4538DDDB3FC5B4992E2268C5A6D2E1C601C1E335D7A683A446BB534FB7518C7094C1E1DD1412469
76B93D4F9633F9D2D477472DE28D3E587C73A3EB982F685923771D23209393ED839FC2AFFC469E37F0D359C6C1E7B9741146BC96E7AFD2BAA6862780
C0D1A988AEDD98E31E9552DF47D3A0B817115A209474639247D33D28B05CE7B52B78345F870FA6DC4A91CBF63650ACD82CF8C9C7E26A0F08A0BEF877
2584332F9F2412A6DCF23208E9E9CD7597BA6D95F94FB6DB4737979DBBC6719EB4B6BA7D9D97FC7ADB451718CAAE0E3D33405CE4FE15831787EEADA4
1B648AEDB20F5E5571FCAB3F43BD8EE3E2BEA0E5970626443D338C63F4AED6E345D36EA57966B552F27DF6562BBBEB8350CBE19D1A5BBFB53D8A194A
04C824703A74A02E64F88A482E7C4F65A46A0824B2B9B762EA7A2B64E1BDBA543A668CDE19F135ADB585C4AF617CAFBE173B8A301907F9735D436996
2EC19EDE366587C8048CE13D2920D36D2DA613430E240BB433316207A0C9E280B956F75EB4B2D6AD34A90399EEBEE951C2F619AD7AAB258DA4B78977
2408F3C636A48C325475E2AD531051451400514514005145140051451400514514005145140051451400514514005145140051451400514514005145
140051451400514514005145140051451400514514005145140051451400514514005145140051451400514514005145140051451400514514005145
140051451400514514005145140051451400514514005145140051451400514514005145140051451400514514005145140051451400514514005145
140051451400514514005145140051451400514514005145140051451400514514005145140051451400514514005145140051451400514514005145
140051451400514514005145140051451400514514005145140051451400514514005145140051451400514514005145140051451400514514005145
140051451400514514015AF6D92EE0F29FA6E0411DBFC8CD3ADEDE2B642B126D52724649E6A7A29595EE3BF40A28A298828A28A0028A28A0028A28A0
028A28A0028A28A0028A28A0028A28A0028A28A0028A28A0028A28A0028A2B9CF117895B46D42D2CA2B46B896E412A17AF5C5007474572BAB78B1F48
BAB2B7BAD3595EEFEE813038E9C74E3A8A7EAFE27B8D1EF6186FB4C6115C3058A4590119EE0D2B8EC74F4573FE2BD6AEB41B21791DBC52DBE42B1627
2A49F4EE2B534ABB6BFD2EDAEDD0234D1872A3B66988B9451450014D00280140007000ED4EA2800A28A2802248638D9DA38D119CE58AA8058FA9F5A9
68A2800A2B3F49D49353825963428892B44A49CEEC77ACEF156B777A141F6A8ECD26B6006E62D821B3FF00EAA2E163A1A2A9E9378751D2AD6F4A0433
C4B26D1DB22AE500145145001451450014514500145145001451450014514500145145001454325D5BC4489278908EA19C0C557D2F55B3D5A2965B29
3CC48E43196C7048F4A00BD45145001451585A1F88A0D62FAF2D62B796336AD8DCDD1B9C7E1F4A00DDA28A2800A28A2800A28A2800A28A2800A28A28
00A28AA3AA6A76DA4D935DDE3308C1030A32493D85005EA2A0B3BA8AF6D22BAB76DD14C81D4E3B1A9E800A28A2800A28AA5A86A763A6846BEB85843E
7696070718CFF31401768A86D6E61BBB74B8B6904914832AC3A115350014514500145145001451450014514500145145001451512CF0BC8D124A8D22
FDE50C091F514012D1451400514514005145140051451400514514005145140051451400514514005145140051451400514514005145140051451400
514514005145140051451400514514005145140051451400514514005145140051451400514514005145140051451400514514005145140051451400
514514005145140051451400514514005145140051451400514514005145140051451400514514005145140051451400514514005145140051451400
514514005145140051451400514514005145140051451400514514005145140051451400514514005145140051451400514514005145140051451400
514514005145140051451400514514005145140051451400514514005604BA74B75E328AFE588791676FB63661D59BD2B7E8A00F35F88E54F8AF4656
C9000F97B72E2BAFF1569FFDA1A415588C9342EB2461464E73DAB92F884B9F18685D39C75E9F7EBD1E914FA1CA7C4B19F07CF9381E6C79FA6EA75A6A
3A90F0C417D630C090416C0849892D2055C678E83FC2A1F89E09F08B804F33C638FA9ABD78A62F00B2C6761160A3E5E3F846681743434ED44EABA2C5
A85922869D3722C87807A72476C8AC8D1FC457B7DADDFE933C1045716B1920862417E303E9CD49F0F9C37832C0018D8187FE3C4FF5ACBF0B81FF000B
13C4673CF03F0C8A00B67C47AA5BF8AAD744BCB6B5DD38DC5E32DF77279193ED5A5AF6A535A5D5858D9944B8BD90A091F90807278F7AE7F53292FC5D
D3109E52D482075C80EC2ADF8C205D4F57D374B69FECCCE1A48A603243FA7D38A065E4D5354B7D7ADB48B982293CD05BED0380547278F5AE8EBCFEC2
E75BD175FB5B5F10E2EA0918A5BDC673B58F048EFD0D7A05084C2B1BC4C351FEC999B4C9A288A23348CE093B40EDFAD6CD67788377FC23DA96DEBF65
97FF004134C48E5FE1F49AA49E1A8DAD8DA18BCE7C994B163CF3D2AF7C493B7C1D727FDB4FC39AA3E079A5B6F878D343812A190AE477A9FE2033C9F0
FA5924FBEC212DF52CB4BA15D4769F7F7F078420BCD3A18C5B5A5B06FDEE4B4A14738F41C56F687AA45ABE8F0EA318D8B203B81FE120E0FF002ACE10
C717C3C3120DAA74C38FA98F3FCCD65F8384ADF0C8AC00194C5385CF7396A0468E9FAC5E6BB7370FA4BC315ADB49E5932024CA7FA0FA568E85AAFF00
69C53AC8823B9B694C5320EC45725F0CE3BA9BC3922C1A80876CCC0C622562A78E4E6B574DD25347BAD5628EF67B8BBBBB769D8B00141E4741DF27F9
D0058835CBAD6B52B9B5D10C6B05A1DB25C48A4866F41FE7B558D2B5991B539747D4B62DF4681C15E8EA7FAD73DF0A3234CD455CE64173F39F7C5335
1595FE2DDA7D9C93B6053281D8739CFE18A02C5CD6BC51AA685ACD9DB5FC16A6DEE0921A2DC4B0F4C1E872455FB0D72F5FC52FA3DEC4996844A36023
CBE3A67BD63F8D915FC6DE1A471F29940391C1CB8E2BAFD46255B2BC9E0897ED260701957E63F2F033F950066DBEB373ABDEDDC1A3B4022B5608F2C8
09CB7B63B55AD27567D445DDBB2471DF5A394913395CF63EB8AE3FE18ADDCBA1EA31DBDC2C32F9F9CB47BB071FFD6ADFD1B484D1B5C9659AF65BABCD
45599FE401463073D7F014032B2F8A6FE0F149D12FAD60690A82861623793D393D2AE5C6BD3E93A9DBC3AD2C11C372A4AC9193F2103A1F5FF1358C91
A49F18A53C131DA86C7BEC1FE352FC464592F3C3D1BFDD6BDC1FCD6819B3ACEB379A65AA5FB5A466CB23CC05FF0078013D7D3F0E6B66D2E62BCB48AE
6024C72A8652463835CFFC44FF009136F73FECFF00315A3E17429E1AD3D4F5102D1D45D0D6A28A2988C9F10D85B5CE8F7ED25BC6D27D9DC86DA37642
9C73581F0A430F091DDD0DCBEDFA617FAE6BAAD573FD93798383E43E0FFC04D72FF0B6477F0961D70AB70E14FAF43FCCD2EA3E87674514531105DB4E
96D235AC4B2CC07C88CDB413F5AC0F0B6BA357BBD4616B24B49AD640B26DC659B241CFE22BA6AE17C09B7FE129F1494FBBF6AE3FEFA7A43E875536A0
C2F5AD2D6D8DC4B18065F9C2EC07A75EB49A66AB06A525CADBAB8FB3C9E5BEE18E7B8FC0E6B99D45754FF848EFB53D01F7B40163B9B671FEB3038C0F
A568783B55B6D562BB97EC82D2FF00CCFF004A8B9193D9B068B858E9E8A28A62296A5A8DB6996DE7DD3ED5C80AA3EF31F402A9C1AEA9D423B3BCB39E
D1E6FF005464C61CFA7D6B97F19DC5CCBE36D06C51731ABACAAB9C6E6CFBFD2AFF008B6CF5BD520B56B7B68ED85ACBE7339983118A571D8EB2E6E22B
4B692E27609144A59893D00ACB8FC4566C207786E61867FB934B1E10F1C739AE57E26EA332E91A6DA472E56E98348EBD1C003F4E6BA9F144312F84EF
A30836476E768F4C0E280B1720D4ADEE2D67B98BCC68E1764242E4B15EB8C751556C3C49A6EA1732DADABC8D71129668CA61B8EDF5ACEF86FBBFE10D
B52C724BB9C9FF0078D66686557E2A6B8A47262057D8E1680B1D0C3E28D3A6D462D3C79EB7331F951D31F8FD2A0F14EA3A6DB5B341ADD94F259B007C
C54DCBBBF3EB58FAD2A9F8B1A1E7BDBB679F41211FD2AFFC4A00F83AE0938C3A11F9D006DC779696BA2DBDC428EB6A224F2D510921481818FA554B7F
14E97753886D9A695B7053B232769271CD5BD023F2BC3FA72039C5B47C9FF74571D6A3FE11EF89F2DB8056D3544DCA00E371FF00EB83401D30F13E92
2E9AD2696586756DBB2489813DB8C0AD0B9D42DEDE61012D24E46E114637311EB8ACAD2A04BED76FB5696352636FB35BB0391B57A9FC4F7ACFD4208B
49F19A6AF73A8315B88CA25AAA967638C703D280B1D0586B1677B7325AA33C7731F2D0CABB5C0F5C7A56778F554F83351DCB9C2291EC770AC0B99A69
3E2BE98CD1490892D8FCA4F2576B9E71F4E95B9E30B98AE7C0FA8CD6EFBD1A2E0E319C3007AD01D4B1E0A42BE0FD2F775F201FCEAE36B1641A401D9D
22E249110B221F42457371DF3E97F0A61BA42778B35553E858E3FAD3BC1772F1F83E126CEE266943B16400873F9E680B1D635D40B6A6E7CC0D081BB7
AFCC31EBC556835BD32E51DE1BC85963FBDF3723F0AC7F03477F6FA54D63A9C451A19328AC7276B7383F8E6B1BC116A8DE2CD7A62A3114BB517B039E
B45C2C76965A9D95FBCA96770B2B4470E06415FCEB0BC73E20B6D2B4A96D44A45D4E0280BFC209E49FC33556599ADBE2A5BC517CB1DD5A3798071920
1393F90A8BE2C408DE1B8A7C0DE970AA0E39C10680B6A75563A9D8DD6C86DAED257080E33C9ABF59C26B3B3366AEA91C93811C6427278E99FCAB4698
82A39258A2C79B2226EE9B980CD49552FECA1BFB59219A347254852C33B491D47A5004A6EADC1C19E207FDF14E4963700A3AB02700839E6BCF3C1B0C
3FD97ADE99A8C0AD7968CC58B0CB85DBC60F5E08FD69FE07B8B6B8D1751B5D4E1479EC5893B87CDB7A8E7B72314AE3B1DF34F02A33B4D1845E598B0C
0FAD73DA2E97630F892FAFACF5213B48833029C8404FAE791C54DE17D0E1B1D2E369D0493CCBB9CB72307A0FCAB17C1C8B0F8DFC47046A1234202A81
C01B8D006978AF5D4B3B8B1B2B7B954966B85594ABE0A2FB9ED5D0096D6ED1E149A2943290CAAE0F07E95C178E6C6D878C3422B6A84DCCB89703EF9D
C3A8EFD4D76F6763A7D95C3ADA430C533AEE60BF788CF5FA66807B1C8E98874FF89B2D8ACAED035BEF8D58E76E467FC6BBA5656CED607070707A1AF3
BD42D3EDBF16560323A218559F61C120274FA1AB7AF31D0BC61A3C9A7AB450DD379532063B5F271D3FAD2D86D5CEF28A28AA242B0BFB577F8B974B8E
5F952D8BBAF182D9FF003FAD6D3A2C91B238CAB0208F6AF3AD3F4AB63F12EEAD5C4862484B28DE7DBFC6931A3D02EE16B8B6789669212C389233865A
E67C077579727554BD9DE6682E44619CE4F03FFD55D35A5B456902C308608BD37313FCEB83F07412DFEB7ADC533B25A09F7346A71BDBA73EDC500B63
D128AE4B419E5B4F17EA7A287736D1C4278C3B6E2B9DBD3DB9A4D4F5CF3FC550E86B27930805A67DFB77F1D3F5A2E163AEA2B92BED420D0FC4D630C3
266DEF7F752C7BF22361D0FF00F5ABADA620A28AAF796E6E6031ACF2C27390F1360D0072571AB6AD7D6B737BA507692DAEBC93029E154752477AEC21
66782369176BB282C3D0E2BCF3C176535E5F6B312EA3750982E73BA36C6ECE7AFF00DF35D36BFAB3C7A9D9691692EC9EEDBE6607951FD3A1A48A6745
58BE2992F6DF4796EAC6E3C968017618FBC2AA6B667F0F697F6EB39E792385C79B1CADBC15279FA7FF005EACEAB7516A1E11B8BA84831CB6FB87F850
223F08EA171AB7856DEEE697FD224DE0BE33821881F5ED59FE1ED62FE5F15DFE937B309842A4EE03001CF6AB1F0E71FF00086DA63FBD27FE846A9689
12FF00C2C8D6188C111061F9E3FAD03EE76B45739AA6AF24FAF5BE87632796F265A6954F2AA0678A8F57BEB9F0E4B0DD3CCF73632B88E4129CB467D4
1FCE98AC74F45626BB73A8DB597F68696F1CB1C69BDA164CEF1D720F5E9FCAA8699ADEA3A87871F56B792DA42BB8989D08DA17AF20F5A570B1D5515C
C697AEEA1A9E84DAA40B67B103165F9B3F2F5EF5734EBAD5B51D396E4ADB5B1917283058F5A61627D226D4A49EF17518822AB8F248180579AD4AE73C
2FAC5E6A775A8437613FD19C2A95E3D7B552D4BC43A8E93AEDAD95E24524771C8F294F2338C7E7482C7614573DAEEB179A298AEE58627B47708C39DC
9F8F4C9FE95B573750DADA3DD4AD889177134C44F4573BFDB3A81B0FED34B48DACC7CC57277EDCF51D8FE54699E209753BEB7105B91693AB1DE41CA9
00F19FA8A571D8E8A8A28A620A28A2800A28A2800A28A2800A28A2800A28A2800A28A2800A28A2800A28A2800A28A2800A28A2800A28A2800A28A280
0A28A2800A28A2800A28A2800A28A2800A28A2800A28A2800A28A2800A28A2800A28A2800A28A2800A28A2800A28A2800A28A2800A28A2800A28A280
0A28A2800A28A2800A28A2800A28A2800A28A2800A28A2800A28A2800A28A2800A28A2800A28A2800A28A2800A28A2800A28A2800A28A2800A28A280
0A28A2800A28A2800A28A2800A28A2800A28A2800A82E639A58C2C171E43E73BB606C8F4C1A9E8A00E5753F089D57538EFAFB5395A4870220B180147
D3EB5B96F6F791CCAD2DF9963031E59880CFBE7AD5EA29587739CF11F86A6D7DD525D49E2B5041F24460E08F7A2EBC3D7B73A48D30EB322DBF97E59C
40B92B8C75CD74745160B98BE1FD124D0ED45A25F3CD6EA0ED578C0209EF9FE955E7F0FCF0EB73EABA4DD2C135D2859D645DC0E31C8F4E95D1514C2E
71CBE0C9A0D73FB5EDF5267B9C7DE9973F36304FD3DAB564D0DEE6FEC6F6EEE7CC9AD6361B95002589C823B56E51405CC39B49BABDD4AD67BFBA4782
CDC491246B82CE3B9FF0ADCA28A041514F124D0490C832922946FA118A968A00E6B4FF000B0D3F4FFECF8AFE436864DED1ECE4F7C673527893C3D26B
D125A9BF6B7B3006F89101C90720FF002AE868A561DCE7AEB41BA93418F4987522A817CB662832530063FCFAD4DE1ED124D0ADBEC8B77E6DAA03B10A
00412724E6B6E8A2C1739F87C3CDA6DCDC4FA2DC0B6170C0BC2CB94FA8F7ABFA769C6D1A59E798CF73311BE4231C7A0F415A3453B0AE60BE87259DDC
F75A2DC25B3DC106589E3DD1B1F5F6A7E8FA10B2BDB8D46F26FB4EA1738124A460281D1547A74FCAB6E8A561DCE5F52F09CDA86A505FCBAC4E24B672
F06225F939CFE35BF670CB040127B86B89324972A17AFB0AB345315CE761F0EC9A76A17175A2DC2DBADC906585D72A48EE2AFDA69F3C724F737173E6
5D4ABB15C2FCB1AF6007D79AD3A28B05CE4E3F0ADDAF88FF00B6DB5522E586D7C46082B8C631F953BC41E18BFD66F609CEAAB1476D2799127939DA7F
3AEAA8A561DCE63C41A16ABAD6922C5F538555C0F37F718DC41C8EFF00E15A1A1596A5636FF67BFBB8678A34548B6260800639AD7A28B05C28A28A62
30F50B2D66EED6EEDA3BBB58D26DCA8C5092A87B7E5599E1CF0EEB7A058359C17F6AF199378CA1E33D6BAFA2958770A28A2988AF76B70D6CEB68E893
1FBACE320735CDE83E19D4347D4AEAF05F42FF006C903CEBB0F201278FCEBACA28B05CE66D34CD4ECF5FD4B528D51A3B8700425C7CEA3B8F43F5AB7A
669D2FF6CDD6AF3C42069D046B18E4E38E4FBF15B7452B0EE145145311CDF89FC38757B8B2BEB690477964E1909E8C339C7E7566ED753D42D1AD3ECA
96E24003CC650C31DF007AD6DD140EE735E29F0C47AC6851595B111CB6B86809E8703A1FAD55F10DCEA973E13BD8CD8790CB6E7CE791B8C01CED1DFA
575F589E23B7D56F2CE4B2D3E2B6786E226491A572ACB9F4C5204CA3F0E063C1B69C6325BF9D523A6EA1A578EEEF588ACE4BBB6BC8B1FBBC650E07F8
55EF08E9BACE8F670E9F78B6AD6D16E3E62312C73D05751401E7DAC596B7FF0009769BAEC3A5C930897E68D1C7CA0E4639EF835A5E2EB6D6B5AD0CD9
DA698AAB2E19BCC99432E0E7A74AEBE8A2C17323C3BFDA29A64506A76896EF0A2A2EC9036E0063F0E9585F10B44BAD4EDECAE34D85DAEE0971B909CA
A9EBFAD76945160B956C2DC5AD8C16E001E5A0071EBDFF005AE3AEECF53B3F88A7566B29AE6D1E2D89E5738F971D3F0AEEE8A2C17380D42D756B7F1B
5978865B2796DF618D628B974CA1001F7E6AEEB16DA91F068D2E2D3A59E7BA46DE51862325B7739AECA8A2C17392D16C6EEF7C2A740D5F4F92D912D8
44652E0EE3D881EDD7F0A3C36B7FE1EB1FECABEB19A75858F9335B26E5753EBD306BADA28B05CC6B5696D92F753BAB7943CEC9FB98C6E608381C7AF2
73F4AE77C2125DDAEBBACBDC69B771C17B3F988FE5121793D7F035DDD14582E79FEA2F709F126CF53FB05D1B386031B3F9679CAB0C81F88AD4F1DDA4
FAC784C1B289A460E9384C7CDB71E9F8D75945160B9CB095B58D43486B682558ED0F992BBA150A76FDDAEA68A298828A28A00E12F74D96DFE255BBDA
262DEFE1DD720AE55B69C9CFBF02A2D534A961F8810C768A52D752406E5472182F5FA577F8190703238CD30A233AB9552CB9DA48E467D2958771B34B
1DBC26473B51073804E3F015C0787F504B4F1CEB1713C37020BD602293C938EBDEBD168A00E13C740C1E23F0F6A122916D0CFF00BC7ECB86079F4ADC
B69E2BDF123DDC2C1ADEDEDB67999F94B139EBF4CFE55B534314C9B668D245EBB5D411FAD3A38E38976C68A8BD70A3028B05CF3CFED4B11F154DCB5C
2883ECFB3CC3C2EEDA3BFE953F8CEF6C7FE121D0AE05C24AB0CD99154EE50BEA48AEDCD95A93936D0127BF96291AC6D1BEF5AC07EB18A2C3B92C3347
3C2B2C2EB246E32ACA720D494C8D1234091AAAA8E814600A7D3242B8789BECDF146EA5BC951164B502324E0118E9F5AEE2A9DE69B657ACAD756D1CAC
9F74B0E47E340D14348B88E35D42EA69D12092E98C6CCD804607AD737E01BDB34D475DDF730AB497594DCE0165E7A67F0AEDE5B3B59ADBECD2DB44F0
FF00CF3640547E15561D0747818B43A65A231E09588026905CE66D6EEDD7E2A5DB79C9B64B10A1B70C120AF1FA1A8355B1B7B1F8876D7BA800F657B1
B282EB90AFD319AEA64F0EE8D2B067D3E12C0E776083F9D5DB9B2B6BAB6FB3DC429243C7C847031D31E94582E64F882C6D5EC14416F00B97917C8210
7DECF5E3DAB6E20E22412905C28DC47427BD57B5D32CED5FCC8205460300E49C7E757298829924891AEE9195074CB1C53EABDDD9DBDEC223BA884880
EEDA7D68038DF00BC2BAB6BD9740EF7236FCDD7AF4FCE9BE26B76B5F1FE8BA948716EEC232C470A79EFF00957483C33A2094CABA744AE7A95C8CFE46
AFDDD9DBDF4060BB85658C9076B7AFAFB52B0EFA989E3B9963F0B5D4479926DA88A3924EE07FA536C6C1ECFC00B67295DEB6A49C9E3272D8FD6B561D
1ACE264665798C7F73CD6DC13E82A7BCB1B6BC50B731F9817381923AFD280B981F0E88FF00844E050465647E076F98D55D1B0BF127580581DD6E31FF
007D0E2B774FF0F69BA679BF6389916652AEBBCE307AD553E11D1CDF7DB3CA9BCECE722661F875CE280B986D0BD87C53864989F2AEE2711B11C64827
1F9F157BE2490DE19F2002D2CB328451D73CD741A969969A9C4B1DD47B8A36E471C321F506A34D1AD7CF49AE1E6B974FB8677DDB7E9405C974D85A3D
2ADA0B8F9D9615570C3DBA115C46831BE9FE22D47C2F2A116D732996327B28E4E3EBC57A2565368B09F10A6B1B8F9AB118F6E3F5CD0D026717A3B1D2
75CD47C2AD14AD15DDC2B44C38DA84E58FE55E8B1C6B146B1C6308802A8F402B39B45B76F1047AC038992231918FBD9EF9AD0B8885C40D117740DC16
43822840D9C8F825186B3ADB13D260B49E288C3F8D7416C636C8BCFAF2DC7E95B165E1AB5B0BE6BBB4B8B98D9F9906FCEF3EF55756F0C5C6A1711DC7
F6ACA2485F7C5941F29CFAD2D6C3D2E33E21853E1693776954FF003AD0D496293C3262998A89215518EA4E0600A8B53D06E354D3DACAF35491A23839
58541C8F7ACBD62DE4D33C31F63D42FA5BAB86947D95946D6520003F01FD69887E96FA9A7873FB2EE7499C4A636891F82A41E84F3C75AD6F0BE92FA3
68B1D9CAE1E4DCCED8E80939C573D6DA2F8C4C79935A55F9461598E7F4AD3F0A5F6A8D7377A66B05649ED7044AA7A83D8D086CE9E8A28A6485145140
051451400514514005145140051451400514514005145140051451400514514005145140051451400514514005145140051451400514514005145140
051451400514514005145140051451400514514005145140051451400514514005145140051451400514514005145140051451400514514005145140
051451400514514005145140051451400514514005145140051451400514514005145140051451400514514005145140051451400514514005145140
051451400514555B9D46CAD1D52EAEA185986407700914016A8ACEFEDBD337229BD8817385C9C6691B5BD3167F21AF6212EEDBB73DFF00C3DFA50163
4A8A28A0028A28A0028A28A0028A28A0028A28A0028A28A0028A28A0028A28A0028A28A0028A2A9CFAA58DBC9E5CD751ABFA6738A00B945355832865
20A919041E08A7500145452CD142A1A69123527196600669EEEB1A96760AA3A92702801D45416F756F75BFECF2A49E5B6D6DA738354ED759B6B9D62E
74C0AC93C033F37F10F6FCE8034E8A28A0028A28A0028A6B30552CC40503249E82B3E2D734A7202DF45C9C0C923F9D006951551352B1791228EF2077
90E142B839FCAADD00145145001451450014514500145155ED2F2DEF6232DA4AB2A062A597D476A00B1451450014542F3C31B849268D5C8C85660091
528218020820F208A005A28A2800A28A2800A28A2800A28A2800A28A2800A28A2800A28A2800A28A2800A28A2800A28A2800A28A2800A28A2800A28A
2800A28A2800A28A2800A28A2800A28A2800A28A2800A28A2800A28A2800A28A2800A28A2800A28A2800A28A2800A28A2800AA7A969D6DA95B182E90
95EA181C329F507B55CA2803163D1AF23528BAD5D95EC08071F8D5ED3F4E834F8DD61DC5A43B9DD8E4B1AB9450014514500145145001451450014514
500145145001451450014514500145145001451450014514500145145001451450014514500145145001451450014514500145145001451450014514
500145145001451450014514500145145001451450014514500145145001451450014514500145145001451450014514500145145001451450014514
500145145001451450014514500145145001451450014514500145145001451450014514500145145001451450014514500145145001451450015C8D
E69A357F1DA9B9C35B585BAB04F56273CD75D50A411473CB3220124B8DEDEB8E05034EC79EF8D59DBC77A259C230046BB540C0E5C8FE42B63C6BA3C1
24B67AB0F95E19904A41FBCB918FCB158DE2E3FF00172B453C8FF54BB876F9CF23F3AF42B9B682EE130DC46B24648255BDAA6C36F63275DD79B4289E
E2EACD9ADC305491641F313DB1DBFF00AD5752F2E25D3E3BB8ACF7338DC23F33042E3AE71FA573BF141D57C2CAA4F2D731803D7A9FE95B71DC5EC3A5
DB98AC84A7C804FEF40C1C7A639EDF9D3114B46F13CDAD2CAF61A6978A27D8CED381CFE549A978A25D32FAD6D2EB4C226BB60B10130209CE3AE3DEB3
BE15A04F0DDC70437DADC11F455AAFE3B253C57E1D942190A4A0AA03D48607FA51D076D4E96FF52D4AD2C9EE7FB354F9632C04DB8E3D8015734DBE5B
ED320BDDBE52CA9BCA96CEDFA9ACAB5B8BAB8F17CF15DDAF931259FEEC121B782C327FA545E239E1D3FC15A97F658DAB02B45B573F292C0375FF0078
D022F5B6AB36A5231D29227B7462AD3484E18FB0153D96AB1DC5E4D652A793750F2632D9DCBEA0D66F801547836C0AF56562DF5DC47F4AC9D42E241F
15ECA38067FD182B81D87CC4FE84501635F55F11CFA3DDDB5BDDD8AB35E4A6385925E3A8033C7B8A76AFAD5E68B08BABDB589AD0B05631392CB9E9D6
B1FE2265750F0F480062B79C29E84E571FCAAD78A91EFF0050D3F46BE222B2BB7DDB9392CCB8F9493F5A4337E4BC9E482DE7D3ADD6E1251B8E640840
238FF3ED59FA17889F5C7BA5B6B1F2FECAFE5B9925C7CDE9C03E95B500863416F0ED0B0A85D83F846381F9579F783AF6EAD63F11BDADA34F225D3B1C
3818396EDDE988EB06B1710EB7069B7968231386314AAF9538CF1FCAA3D6F58BED22379E5B48A4B7DE151958E4E7D452E986D65D1B4AB9BF1FBD5DA2
2662490E78038FCBD2A7F10429710D9C520055EED01C8CF63400ED0F5AB7D67491A8423CB8F2C1831FBB8EBCD645978C8EA3A8CB61A7E9AF2CC99E4C
802E3D7A562F86C3DA2EAFE18DECB299C2C201C1DA7866E3B6D19A9FC1D6C2DBC77AF44189589422E4F6C8A571D8EEE22CD12191363900B2E7383E99
EF525145512676BD73259E877973082648E262B8EB9AC6F065C5A6B5E1B4925B45DE84C5297504B90012738F7AE9D955D4AB00CA460823208AC6BDCD
A5BFF65E870225C38246DE16207AB13EBE94866078135277D7356D2D19DAD6DD8B45BBF879C62BBBAE7FC39A459E82E6C622D25DCC8679A66FE3E71F
D6BA0A10338ED40B6BFE328F4D56C59E9CBE6CE47F139E83FCFA574FA85AC17968F1DCC62441F363DC572DE0076B89F5CBA939924BD2327AE0574DA9
35F2D9B1D35227B8CF0B21C0C7F8D036735F0C948D02E189C8376F8F5E30293C66C349D4F4BD6A1521C4DE5CA47753EBFAD5CF06E95AA6936D3417FE
408DDDA4010E48626ABFC4E1FF00148BB7759E323F3A5D03A9D6A32BA07539561907DA9D5534A3BB49B3620826043CF5FBA2ADD5121451450015E67F
13423DF697A658C2893CCDFC000CE4800715E995E66B9D6BE2E9DC4F97627A751941FE26931A3BCB1D2ACACA28520B788794A02B6DC9EDCE7F0AD0A8
6E678ED6DE49E660B1C6A5989EC0571DE35F11CD6F61616964CD0CDA8852C5465910FA7BD30DCEAA1D4ECA7B936F15C46F2F3F283D71D71EB585E37D
7FFB33495163748B76F32A00AC091CF39F4A643A74D36A3A5A59DABC16765F3C93C8A15A5381F8F3597E22B28355F889A4E9E91285854CF71851CF3B
B9FAE07E74811DCD949E6D942FE62C859065D7A31C727F3AB151C30C704290C28A91A0C2AA8E00A9298828A2A39A68E085E6998246837331EC280303
C69AD8D274799621BAE668D8228EAA31CB567FC2D72FE1424B6E3F697C9F5E0552BDD4ED2EAC358BCBAB7B9696785A384F90C55171C738C75C73563E
131CF8565F6BB71FF8EA52EA5743B8A28A299271DF1205BC1A03DD98FF00D2CB2C51480E0AE493FE35BBA0A436BA55B592CE924B0C404803824377CF
E39AE53C7EA759D6F49F0FC048691CC92B67855FFF005026AD6AED168FE29D0A1B140AD3379720079753C73EBF5A5D4AE876B5147345296114A8E578
6DAC0E3EB5CA78D75E5B49ED34886E56192E9809DFFB91FD7B6699A45B87F15C32E976EF1E9F05BB2BC84E04CC7F8B1DFB517158ECE8A2A39A458627
964384452CC7D8734C4494571BA35DC9E2737BA9DC6E4B38494B54CE002072D9EF553E1BDCDDEA73EA3A85D4CF20DC2352C4F3DE95C763BDA28A2988
29A480402402781EF55354BF8B4CD3A6BC9FEE44B9C0EE7D2B91866BBD57C337DAECD6ED713B83F658173F22838E31F9FE740EC7774571ADAF5C685E
09B19AF94FDBE65D91A3E73924E33F418A35B8E6D2741FEDB9AEDBFB4632AECD9F95893F740F4C7F2A570B1D9515CEEA5ACACDA469EF69234726A453
CBC1C30538248FA74FC6B1BC33A8DC6B2977A3ADC4BE5C12B869B765C479C05068B88EEE8AE2BC21A8CDAAD95E58DBDDB620B99144CDCBA45FC38F73
CF3EC6B43C15A95C6A16776973234BF66B8689247E5987B9EF4C0E968AE0F58F13CB0DF9BB82E8A5BDBDDADB7947EEC83077135775AD6A4D3DF4DD4D
AE1E382EA508616FBBE5FAFD4F5FC695C2C75F45723AE6AB75A7C9A7EA06E0A4771388D2DC7468CF427DF1CD6E6B77735A58836D1969E691618F1D8B
77FC39A2E3B1A545727E17D4263AD6A3A41B8FB443660624F46E847E79E3DAA7D0B57B8D53C41AA22BEEB2B62238F0BC6EFAFF009FD29858D3D73588
346B11733A97DCE11501C124D5E89C4B12480101D43007A8CD79F78B22BFBFF10699A1BEA08EEE7CEF31610A63E78C8CF3D2BAFD43501A1E8C27BC91
EE1A30A99C006473F4E9482C6B515C5EBBAFEA9A2E8D6FA94B22335D8C087601E51232307BFE35D1433DD9D122BA7921137922472CA769F973EA314C
2C169AC5B5DEA973A74224F3AD865C95F97AE3839AAFE27BBD46C749371A5C7BE657191B41F979CFF4AE6FC051DEDCDDDFEBB3CF1ADBCF2BABA1183C
77CFA568CFE269D358D3515556DEFA6F292361F3919C6FA416D4E87499E7B9D2EDE7BB8CC733A02EA46306AE565EAD7D25B88EDED76FDA67C842DF75
07F78D60C3E21D4EF7C397971A7A2DCDDDA5C18D98281B900CEE028B858ECA8AC7D33588EE3C370EAD77F2298F749C7420E0F1F515977FAF6A369650
EAAC90ADA48E1560232EC0F7CFAD170B1D6514D53B943608C8CE0F514EA620A2B2354D63ECD73158D9A09EF653C27641EADE955F48D74DE6AD75A65C
983ED102861E492411DFA9EA38A571D8DFA2B1A4D567B9BA96D7498A399A2E2495DBE453E9C75A9B45D506A76D233208E785CC52A039C30FE94EE2B1
A7456347AEC4D79E5B205B76DC2394367257AF1527F6BA2EA71DACA8112588C91C8587207527D05170B1AB4563FF006C04D5A0B49A108972A4C126FC
EFFC3DF8AB97D766D63558E332CF2644718EE40FE5401729AC70A4804903A0EF589A5EBC753B6B95821437F6EC55A0DF8190719CD5DD27548754B679
62528F1B98E48C9E55876A2E1633E6F14DB5A388EFAD6E6090FDD5DB927F956F83900F3CFA8C571DAFC0754F1AE976609F2EDD7CD94038F7AEAAEAE6
1B480CB3B8551C0F527B01EF490D9628ACC8354125F8B39EDE481DD4B45BFF008C0EBF4352EA1A825881989E4C0DEFB47DD5EE6988BD4563DE7886CA
D51641BE58C91B9D070A0D6AA3AC91ABA1CAB0041F506801F556F6F63B28D5E549183363E45DD8AB5484800924003A93401989AF696EE105D7EF3FB9
B1B23F0C54F16A56934A228E71E61E8AC0A93F4C8AE3BC210B5D78C355D4198E2262AA3B7CC6AFF8EC4F20B1874F5637A64CA141C814AFA0EDA9D7D1
5146CCB6EAD39018282E7B038E6A8B6B9609E6169880833B8A9C37D2988D3A2A24B885EDC5C238688AEE0C3A63D6A2B6BFB6BA491E193708890FC118
A00B54555B3BEB7BD573031250E1948C15A9E4748D0BC8CAAA3A963802801F45513AA58884C86E54286DA720E41FA75ABD40054534F0C2019A58E307
A176033F9D44D7D6CB2F94D3286CE39E99F4CF4A875D609A3DCB94572A9950DEB4016E39629977452248B9C654823352D73FE1BB9B5B6D36D6D249C0
B8906FD8C30727B66BA0A0028A28A0028A88DC400B0334794FBDF30F97EB4E2A92C7860AE8C3A1E41A0055756242B0257A807A53AB9DF0FA06D63539
623FB957D883B7E1F95740EEB1A96760AA3A92702801D45337AEFD9B86EC676E79C7AD2BBAA296760AA3A92702801D4535583286520823208EF58365
247FF0915F4AEC228D06C1B9B00B77EF401D0514D0430054820F208EF4EA0028A28A0028A28A0028A28A0028A28A0028AA3AA4379342BF619BCB756C
91D370FAD5B6611C65E4230A32C40A007D154F4DBE17F6C6758CA2EE2A32739C77ACC6B9BD5F1243642E8BC401660540E319C1C0A00DFA28A2800A28
A2800A28A2800A28A2800A28A2800A28A2800A28A2800A28A2800A28A2800A28A2800A28A2800A28A2800A28A2800A28A2800A28A2800A28A2800A28
A2800A28A2800A28A2800A28A2800A28A2800A28A2800A28A2800A28A2800A28A2800A28A2800A28A2800A28A2800A28A2800A28A2800AA9A85DC969
12BC5692DC9271B63EA2ADD1401E69A8C3ACEABE2DB3D61747B85B5B62BB558024E0E7F2AEE6C351B8BAB8F2A5D3AE2DD42E4BBF4CFA569514AC3B9C
4F8DE2D4F5BB3FECCB3D2A6212E03F9CD8DA40C8E3F3ADF7B8BB8746411584AF36CF2F6647040C67E95AF45160B9C87836DF51D1F4E4B1B8D35C3493
B48F26F00283DFF4AABE21B1D5B57F10E9F3DBE9AE91D83962CEE06FE7B76AEE68A2C1730AF2E7541149716DA32FDA950A464CEA78247A76C8A4D0B4
9B84F0F4963AC04792E0B99429CF0DEF5BD45160B9CDE876579E1EB296C45BBDC5AC44B4322302C73D8AFF00851A368F39D76F35DD4A358EE2701218
81CF96806393EB5D2514582E719E22B2D575BD4F4FF2B4FD90D94C65DD24806EC1FF00EB0A93C57A4EABACC367716D0A433594C6458CC8096E9CE7F0
AEBE8A2C1730629351B3B56BC7D3CCF7974419638DC0098185ACAF0869BABE8C97F2DC58076BC9FCD004CA0AF5E0E7EB5D9D14582E7337967AB5EEA7
A6A476B15AD85ACA25900901C91C8181FE79A9B549F539645FB3E985D2DAE15F2641FBC03A11F9D7414530B9CB45A64EBE20B8F113589129B7F2C41B
86E27B907E83DAB37478358B3F146ABAAC9A4CA62BBE15378C8E95DDD14AC17238CB346ACE9B18804AE7383E95251453114F529AE2DF4F9E5B480CF3
AAFEEE31D49FF3CD70FA6C9E36B4799CE9AB234CFB98B119AF44A293434CE67C310EB335F5DDFEBD10864651145182080BD4F4FC2BA6A28A6238CF02
29B5D475DB0906D78EECB85F5524F3FCABB3AE72F2C2E2CBC4F0EAD630F98972A21BA45EBECD5D1D243615C778FCB5E2E99A343F34D7772AC57FD95E
E7DABB1AC1B7F0FB2F89E5D66EAE4CC429581318F2C1EB43044BE22B9BCD3742966D2A20F2C41405C670BEC3BD49E1CBBBDBED160B9D46311CEF9246
31C7D2B568A620A28A28031753D7E1D3E49A136B7524B1A6E1B23CAB1C640CE6B92F01A5CD95D6A7AA6AB6374B25C370FE593DF247F2AF47A2958773
85F10EA771AEDDD8E8F636972B6D3CCAD3CCC846141E47F5A97C45A69B7F17E8FAC796E6CE01E5BED19F2C8CE381DABB5A28B05CC8B795B54BD59D01
1696CC769618F35F1D707B0CF06B97265B2F19EAF726176BB9E358AC9403CE7BE7A62BBFA618D0B872AA5C0C06C72280B9069D68B63630DAA12446B8
C939E7BD5AA28A6221B99C5B40F3323B841921064FE55C57883C466F67B5B2B7B3BB168D329B993CB39283B62BBBA2931A39DD4F57B31A3C900B6BB5
8E685A38F101C74C0AE63E1F6A3FD8BA4DCDB5F59DE0669CBAED849CF007F4AF49A280B91C6DBE357DACBB803861823EB525145311C481158FC47BBB
ED524112BDA816CEDC29EC467D78A8B4EB2B8D7BC70DAF3AB2E9D6A365BEEE37900F207D79FC6BB79228E5004B1A381D372834E550AA1540000C003B
52B0EE7037905BDB7C4F7BCD5768B67B706279170BB80000F7E86BA0B39A4D575C4BBB7256C6D15910E0812311CF1EDC56D4D6F04EA04F0C7281C80E
A1B1F9D3E3448D0246AAAA3A0518028B05C7D62F8B209EE7C337F15AB3094C591B47271C915B545311C0E9527D97E1E45A7E98C24BE9E17FDDE79527
25B3E98150780251FD9DA758D93307F31E7BCDA48DB838507EA3B577D0DA5B5BBBBC16F146CFF78A2004D36DACAD6D1A46B682389A46DCE557058D21
DCB345145311CCFC41B59AEFC2376B6E0B3261C81D48079A8F4CD42DECBC1764B68EAD29B6558E346F9B791CF1DB0735D4101810C0107820F7AA9069
9636F319A1B5892427390BDFDBD290EE71BE2AD3E48EF3C2FF006B95DAD6399527691B215B83CFD791F80ABBE3979356B68342D34AC93DCB8672BC84
51EBE95D6CF04373118AE2249233D558645456961696608B5B78E2CF04A8E7F3A2C17332FE0B6D1FC305DA38C9B0B622266E769DB8E09F7AE63C2A9F
D8DE01BED5E4E26B95775EC7D063F1AEEAFECADF51B392D2EE312432F0CA7BE0E47EA2A29B48B19F4B4D365B756B5450AA87B63A500719E1D864D07E
1D5EEA4EA7CFBB532F039C37CABFCF3ED9AD1F06A496767A7D944CAFBE26B8BB6EBCB7DD1FA8AEA4D9DB1B2FB11853ECFB367978E36FA545A669969A
55B7916716C5272493927EA680B9C2F8B6C2DAE3C4BA4683651AC71C9319E7553DCF527F0156FC5D12EA9E29D134588068E13E74898E001D8FE02BAB
3A3D836B23563003781360727B7D3D69C34AB21AA9D4C45FE965366FC9E9D28B05CE4FC4089A8F8DEC2CA41E5DA69B01B99318008ED81F90AE92DAEA
E65D0DEF2E235598ABCB100BCA8C1DBF8E3F4A5BFD034ED42FA3BDBA80B4C8BB7218AE4673838EB5A7B1766CDA36E318C718F4A2C173CF7C2C64D3BC
2D79A993E65F6A2CEEBB7AAA8C8DDEC320FE95ABE0155B5D16DD586EB8BD779E420F2ABD891F90FC6B6B4FF0FE9DA69B8FB2C2409C10CA589014F551
E82A5D2F48B3D29196D11BE63D58E481E83DA8481B399B7991FC77AC6A53460269F0089493FC47000FC79FCEA2F14F9DA8F8ABC3F6171118EDDD84A4
67B81923A7B62BA593C3D6526AEDA8B07323E0BC79F918818048A9F54D22D353F28DCAB89216DD1C91B6D65FC68B05CE4BC7F225F6A9A4E97CED3720
CAD8C81C74FCB9FC6ADDE6AD2CBE04D467930AFF00BCB78F6AE323381FA66B66FF00C37A75FD8C7693A481227DEAE8F87DDEA4D4D71A2DA4FA21D2B6
95B7D9B47A8F7FAD1A85D19BE1B8A2B2D1F4FD1CC3B8BDAF9B382DF7437AFD4922B0767DA3C7D79A8DC845B4D1E3C46BD06EC6140AEBF46D222D26DF
6095E799B87964EA4761EC2AACDE19B49B556BC334C2391C492DB83F23BFA9FF000A2C172BF88650BE1D9B5636A52EBECFB570D9D809F5E9DFD3BD41
E1481745F064524C0799365C8C7259BA0FF3DABA4BEB486FECA5B4B85262994AB00706B334CD04D97922E2FA5BA4B7FF0052AE31B7FC68EA17D0E6FC
576F7169A0E8FA2C6CC8D338591C1C007A9FAF24D6BDD699069B6103DC799A8DD2958E08E56C8DD91D00F4E6B5B5BD221D62DA38E4768A489C3C522F
546A5B4D2FCA9927BBB892EA74185670005FA0A2C1734467032003DF14B4514C4707E0A0F7DE26D6351987CD1908149E8CDC9FE4455DBEB6B4D2E3D7
2F34CB7649D6150EEA720927E603E800ABBFF08ECF6DAC5C5F6977DF664BB399E2F2C119F51EF5AB169D0258C96AF9956607CD67E4B923049A561DCC
4F00C263F0C477529FDE5CB34AEC7EB8FE86A8F836DEE2E65D7AE0B34715DCA551B078383C8FA66B4ED7C3F79043F603A893A7072551461C2E73B73E
95BF0C51C312C5128445180A3B5017380B9B261E2FD2B4582590C56716F91C0E5B3CF356F5749353F8836769149FBA8213E76C3FC3D483F9815BBA7E
866DB5DBDD526915DA73FBB033F28A7699A2B59EB37FA8BCAAE6E8FCA00E54673D695877312FA45B8F1A24870969A445972A3B91C0FE95BF37948ADA
CBB9E2D818D64C009C13F99CE2A8CDE1E9A4D5AEE65B955B3BC6469931F39DA3A03E8697C6B73F65F0CCD1C6AB99B112AF4FCBF2A62333E1FC223D2A
F3569C9792E65762C7A955F4FF003DA97C03E64975AC5C6710C930DAB9CF3CE4FEB56B4DD0B52B2D34E971DDC7F64624872089101E4AFBF5ADFD3EC2
DF4DB55B6B54DB1A9271EA4D1606CE6F439A2BBF12EB1AB48EBE4C38891FB6077FC8543AECD797DE31B6D3ED423FD997CDD8FF0074F19C9FD2ACE99E
1EBDB49A5B392646D3CDCFDA377F1B9E30A7D3A558D474BD420F107F6CE94B0CCEF188DE191B6E7DF3F850327B6B1906B50DFEA5701AE9E329145183
B1063279FC69DA8482CFC3B792B4E276915BF79C6199BE518FA7F4AB5F63964B59CBC805D4F1952FD9323851E807EBD6B0EDB42D465D226D36F3CB8A
31964656DC5DCF7F6A0453B7B4B6B4F055AADD9919EE25126E8BEF127FFAD5D95AC6B15A451A6FDAA800DE7E6E9DFDEB074FD3B519E5B15D4ADE1861
B24C05460DBDB0067F4AE96840C2B375EBB165A2DD4E71C21032719CF15A55CC789A1D5F52D25ED60B00B92093E7024D0C1107843478A4D10DDDC238
96F09901121E14F4E94C8AFDF43F1545A5366782EC811BB637A93EFDC56869377AAC50476B268A63589151089400703BF1445A0C973AD47ABEA2EAB3
44731C71F41F534011F8AEFA6592D34BB53892EDC0638CE17FC8359DE3C84AE99A7E9D6D8DF23F96B8E3D066B5B5ED3EE1EFECB53B28BCE96D5B0D16
EC6E5F6FA64D52D4B4ED46EA7B3D50DB799342C49B60FF00746303A9EBDE863469DA98AFEC6086C6446B6B77F2A507209D98181FE7D2B1F439A7BCF1
3DFF0090CCB6A92672063A7F0FF2AD7D1EC25D2B44942AEEB9937CC57FDB2338A8F40D367B5D0E75954A5CDC6F62BDC139C50055F0FCAD73E22D52F1
72B6C76C6A4F42D9FE7CE2A7D5A56BCD76D74B8D8841FBC9403D475AADA0D85CA41059CD6AD045139965671FEB5B3C0FD6A7B8B59B4FF10CBA943672
5CACF1ED013AAB77A00CEF13812F8934EB188AC71B01249838C632791F415A9E20D5952C608ECA447378DB15D5BA0E2A9DDD95FAEAEBA89B633BCB6C
50041FEAD9BD7E829DAAE8924563A6C9650799259B0678D7AB7AE3F1A00B7AC69976FA3C16B6934482320BB487009CF5CFE7F9D33C5533C5A4456AA4
3CF2151EA4E3BD5D8667D5E35125B4B0DB83B9B7F05C83C0FA7AD54BEB4B9BBF14DAB18985B40BBB7EDF949EB8CD0C466EA919B51A3D8425BED6D2AB
3B280587D6BB0491245CA3AB80719539E6B9ED4ED274F10A6A11DBC936D8B11141908C46326B674CB5367651C4C72FF79CFAB1EB42065BAC6D7B5092
D562B7B6204F31001FEE8CD6CD606B30491EAD69A8085A58621B5827241CE41C7A734C1153C556E63B5B3B6B651E6484A640C6EE075ADA5996DB46DC
2542608429653B8060318FCEB12FC5DCDA858EA92C521B68DF708C2F2A3B67EA71525DD95CC3E1894AA3892593CF963039507923F0A40687872231E9
6AEC1774CC5F23BFD69BA6DC26AB753CCE374709DB1A91C0CF5C8F5A9AC49BBD3822EE484C2225CAE0938C13F4E98FC6A8699249A4E8D7027B795658
9D9CEE53B5B3DF3D2801DA3465B5FD4A60C76447CA098E0739E3F2A8E3D556F3559C18E49A287291468A58337F78FEBF4A9BC35089741DFD0DC963B8
8C920F1553409A0D2629ED2E564178D31CA0524C9E98ED401A9A72DC5A59DD4F760282CD2AC63A2AE338AA7A259A4B6F71A85EE1C5C92F8619017D6A
E6A76B3CDA3DD2460B4D28DDB78E3A703F018AA46E84BA345A7D807FB4BC62320023671C93FAD003FC337724FA732AA99365C14E5B1B53D7FF00AD51
594C67D46FB4E81952239395C9DA071FAF15A56568349D20C698678D19D8FF0079B1553C2F6E56CE4BA917125C393C8E83FCFF002A009ED8368BA3B9
BD984BE592C08F43D073559A49A6D1DF52959A39972F1E1B851D3A54FE27B792E34691220C483B9B6F5C60D50D42EC49E1F36F623788E21E71C1C28C
74FAE6802D785269AE74D7B99DD98C92123273C7AD6ED735A02EE9E016CD9B782008C4640663FE7F4AE96840C639207CA327238CE38CD73BA44B38F1
2DE5B2C8D24289F3963FC59EB5D2D72FA3B471C52C70BE6FEE5986EC7DC03D7D3BFE54022FDADD4979AEDD441DD63814A800F19E99A8ECE79F4C4BE4
BC91AE044EA54A8249DDEBE9DAB3FC333476097CB7D37FA419B0548F99B1D38AE92DA1D86499862598867F6E3007E140190B7F28F13DB5AAC8C56580
99933F2AB019E076ED4EF1519A1D3DAE21B992338F2F603F2B67D6B234C9DA1BCB9D43027BD9A4F2E38B1CFB9E3FCF35A9E21612DF5959C9F2C4CDBD
989E3AFA52E805BD02C5ECB4D45927924690072188C264741C562812EABAA5F4B6F7461318C861D401D2B557517FECBBDB966511AB32C1FC2718E2B3
6CD1AC3C2F737E233F68B80496001E0D302F785350B9BEB39D2ED8BBC126C1211F7C56FD729E1C91EDA5B5B38E4DEB22B492F7EC48E6BABA1030A28A
298828A28A0028A28A0028A28A0028A28A0028A28A0028A28A0028A28A0028A28A0028A28A0028A28A0028A28A0028A28A0028A28A0028A28A0028A2
8A0028A28A0028A28A0028A28A0028A28A0028A28A0028A28A0028A28A0028A28A0028A28A0028A28A0028A28A0028A2A1B8591EDE4585B6C8C84231
EC71C1A009A8AC6F0E5BEAF6F6922EB53ACB296CAE0E703EB51DBEB8F73E279B4A862568E01991F9CAF1F975A571D8DDA28A298828A28A0028AF3FF1
449A8D86B9A6699A6EAB79BEF5B90F2E7032075FCEBBD45DAA1724E06324E49A0761F45145020A28A2800A28A2800A28A2800A28AC0BBD7258FC5B6D
A2DBC28C248BCD91DBAA8E7A7E028037E8A28A0028A2A0BCB84B4B39AE652024485CE4E3A0A009E8AE33431A8CBA45CEB3AB5FDC244C8D2C510948D8
A3279A8BC032DFEB363757B7F7D76CA642910121000C7F314AE3B1DC515CAF834DDDDD86A02F2EE69545CBC51C9E61CE14F507B75AA7A119AE3C75A8
C62E266B6B118C3B13B98F1CFEA68B858EDA8A28A620A28A2800A28A2800A28A2800A28A2800A28A2800A2B99F106A5336B763A1D94ED0CD70A64774
3CA81D3F91AB7E1A8F5586CA55D624DF27984264E4ED14AE3B1B745145310514514005145140051599AE6B56DA1D98B9BA0ECACDB405157A1904D0A4
80101D43007A8C8EF4012D15CBE8DA96BD71E22BAB5BFB331DA479F98A8007A60F7CD751400514514005145140051454374671039B5119940F94499D
A7F2A009A8AE4EC3C4DA84DE28FEC3B9B3804A8374AD139214633DFEA2BACA0028A2AB6A379169F6135E4E711C2858D0059A2B90D37C4DAA5CEB1676
D71A6A2DB5E2092391037DD232093D2BAFA0028A28A0028A2B0FC41AFA6913D940A8249AEA508149C706803728A2A392458A3692460A880B313D80A0
0928A8E09A2B8852681D5E3719561D08A92800A28A86EA65B6B5967604AC685881EC28026A2B2740D6E1D76D1EE208DA308FB486E735AD4005145140
051451400514514005145140051451400514514005145579EF2DADE458E79E38D9812033638A00B145416F7505D0736F2AC811B6B6D3D0D4F4005151
4F3C36D11967916341D598E29D14892C62489D5D1BA329C83F8D003EB9BF10685A96B134452FA18A28241244BE59EA3D7D6BA4A2802B58A5D47688B7
D2A4B38CEE645DA0FE1566A06BBB652C1AE21057EF02E38FAD4F40051514534532968A449141C12AC08CD4B400514514005145140051451400514514
00514514005145140051451400514514005145140051451400514514005145140115C799F6797C8FF59B0ECFAE38AC2DF7D7FA3AD8496D309E45D92C
B2A6D5033D7DEBA2A2802B585A47636515AC44948976826A72AA58315048E871D29D4500145145001451450014C4448C62345504E70A314FA2801888
91AED8D5547A018A7D14500151086247691234576EAC14027F1A968A0085ADE06944AD0C6641C872A33F9D4D451401592C6D63B96B94811666EAE053
6FB4EB6D41152E50B04395218822ADD14015A5B1B69AD05A4B12B4000010FB53CDAC06D7ECA625306DDBB31C62A6A2802ADA585A592116B02C60F5C7
3FCEAD514500145145001451450014514500145145001451450014514500145145001451450014514500145145001451450014514500145145001451
45001451450014514500145145001451450014514500145145001451450014514500145145001451450014514500735E23D66E74EB792EA1658E28DB
CB5C804CAF9C1033D8533C5BAF3E9FE1FB792DC91757851620BD79C138FE5F8D67F8C209B55F15E89A4E08B607CF94E71919391F90FD6A96A2EFAC7C
40B12536E9B63F3A39E14EDE49FA678A92AC76773713D9E9F0C79F36EE40B129E397C72C7DBBD63695A9DDB5FEA3A4ACFF0069BA80AFCEC3E5427EF1
3EDED5A76F7E5F4D9F53B948D523DEF0E472100F5F5241AE6BC07198348D4FC4573B9A4B82EF8FF61413D2988DBD0F519E6D7753D3E49DA74B62A77B
6010C7A80076AC9F136B32BDADE5E453BC16F6338862DBD6597B9FA0A83C206EAD74459D30DAA6AF725C6EE76A03CB11D8633F98A93C5B07F6978BF4
8D2CED8ED220D773B64283CF39FC0019FF006A9741F53A6B9135E6831CA2EDED6430895A54E08F9726B8DF86767757525D6B53DCCA0C926D3FF4D7EB
FA56C7883509A4F086AB7718DB048A23B65C60EDE84FE3CD37C29E659D8E91A55AC801488CD79F2838DDC81EDC9C501D0EC68A28AA24E57C43AF37F6
B5BE83A7CE23B99B99A5DD8F297EBEB4FD2270FE209A0D3D9E4D3A1870EE4E55A4C8E41EE783FD2B0346D2D2E7E246B33EA280884978D5FA104F07E9
8AEB74DBA1777134B62ABF61890A20518123E7923DB8C5228E7AD11756F89F75704662D362083273F3F4E3F3AD2F146BAF6B7969A458CA897978E14C
848FDCA773585E0DB96B7D2AFEFB6E751D46F592341827D7207A0C9FCAA48EC04DF13AE24D502491C76C1A32F80AC4000123A7AD0069E90CA9AFEDD2
A5964B08A122EA576F94BF38C67BD65EABE209837F6B89658EDD6EC436A1780EA3EF1C1EB9AEA2D6E16EF5664B17096B6A0F9BB140591CF6FC3D6B98
F11C69AC78EF48D1628C791640CD32AF0B8386E9F403F3A00BBAE6A6FA7FF65EA17324914D7772BFBB07858C0FBB8F7C8CFD699AE5EDDDAC9617AD2B
C77B737188A147071176047F3A8B5A8C6B5F11F4FB203743611F9B2F3919EB8F63D28D62E239BE21235DB85B6D26D44E14F1B9CF231F98FCA901D56B
BAB41A2E9735F5CB61506147F798F415C8FDA64BBD2A1BBF31E5D66FD834491BFF00A919E3E980299E32B7D4750F0E697F6E6E24BA0D2FCA00453C2E
71EC4D7513CD6D62F0D8697147F6B9542A6C51FBB4FEF3103A7D7AD3021D6EE7CC86F526664B7B3837CDB0E373919033E83AD50B4D4C6A1E1D37B792
C82DAD2DFE770C479CFDB3EBD07E269DF10EE9E1F0CAD8C659E7BE9520418E5B9C93C7D074F5AA9E2CD3DACBC1B61A45B7EED649A28A5751C0F527F1
E7F0A011D07872E25D47C336B3DCBBF9934672E0E1BA900E7D718E6B90F0B58B6B1E2ED67527BA9CC30C8228E457C3373C0CFB015D4C3394B4BC4B2E
2CAC6D8C51903EF3AAF507DBA562781C9B4D12CE08195EEEEE7334E71F7533CE7DFB7E7401DD5158975E23B3B6F1043A33A399A5C7CC3A027A56DD32
42A86AFA70D52C5AD1A57891C8DC547503B55FA8A59E1876F9B2A47BB85DCC067E9401CA7C47BE5D37C20F6D0B6C69CAC28A3FBBDFF415A3E09B0FEC
EF0A58C4CB8778FCC6FAB73FCB15C8FC42B9B7D5BC49A4692B2A1895F32306E06E233FA0AF468A5B7CF910C91928BF71581207D290FA146C2CA2D034
69638DCBAC41E5677E327926B87F093CB2E8BAAEAD713341033C8F23A13976C7033E838AE8FE226A02C3C2B711AFFADBB221403DCF3FA67F3AC4F105
A9D27E1FE9BA4C798CDCBC6B3B83D33CB67F1FE540D1D47846F27BCF0C5A5DDEB12EEA4EF63C9504E09FC2B20EBB6BAA6A37A67B9F2F4BB0181B090D
2BFF00335A3BF7E897B0698AA60B5B468612A30247D9D41F6FEA6B1BC130E95A778416FEE96179A42CD26FE58B0380B83DFF00C6811A9A66A7763408
B21BED171334768241CECCF04FD07F4AAFA7EAF72753D4B4CB59DEEE44991227939F2F23E763EC3B56D59CD25BE8ED7B7E530A1A655C00225C70A0FD
3F9E2B95F0046D1E9DAAF88E64264B82ECA3FD95C9C7E740CDBF0F6A52CFAE6A7A799CCF0DA9003B9CB6EE87F0AE96B81F0596B3D262B95092EA5ABC
C5FE6E7080F2C71D873FA577D42130AC5F13EB91E89A779BD6E263B204FEF356D579F789EDA7BEF88FA25ACE035AED122023838C9607DFE5FD698227
F13EABA9F87B4282E9AE59AEEECED915B90876E7E51DB15D5E8E974BA4DA8BE666B9318326EEA09E707E9D2B89F194A9AAF8BF47D3E4C8B38A6FDE39
1F2161CB0FC862BAAF0EEA173A9A5E5CCA7FD1CCE56DC600C28E0FD6903D8DBA82EA1927B69228A7681DC6048A012BF4A9E9AECA88CEC70AA324FA0A
623CE34CB5BAD4FE22EA28D7F2B358C1E5ADC0500839031FAB7E55D66AD7AD1E6C61B878FC980CB71393F32281C73EA4D607C35517175AEEA5963E7D
C85049EA064FF5AD0F887235BF862716C804B772242C55796CF4E7F0A43EA3F43D76597C12DABDEB65A347C3E396DA4804FE34EF0CEA772DE14FED5D
4A42E58B48B93D57A28FCFF9D731E2A2F69E11B0D02C9831468A3B8207F130242F5EB9C93F8575DA32C6F1FF0064A4314967651246EC070D20C13F86
450329FF006EDD59EB10477AE185CDA99BC8551F211D003EA6A4B7D57508FC5169A7DD3A3A5DC0D2988260C0719C13F8639AC6D3A21AE7C4EBEBBCE6
DF4C01570720BF4FE609FC2A4D3AE049E2FD735CB9398AC57ECD08E9B8FA7E38A00E9BC43AD45A2D9095977CD2B6C863FEF37F9358575ACDF695AE69
96D713F9D36A0C9E6DBECE23527036FD39A835A8A4BDF88DA2C37710F2960F342672370C93F5C1AE875268CEBDA7471DBC52DCE5999D9798E3F5CFE7
8F7A04723E39FED2BEF1069BA24571138965F3906CC18F9E3773CE0576D77F6D8349056E2213C69996665E3007240F5AE460944BF11357D4EE57F75A
640234C742C7000FAF26AFF8D351BB1E0D89563DB75A832425179C6E1CE3F4FCE819ABE11D5A5D67424BC9C82FBD90B05C6EC1EB8A8BFB5AE355D4AE
AC74791112D7896E19770DDCFCA3B557D347F66685368F6A8247B3B26696553C09483F2F1DFF00C2ABFC35856DFC266EE4243CF3492CACDEDC7F4CFE
3408D4D17586B8B8BEB0BD950DCD8B7EF1C7CA197D7DAABBEBED1EA963E7E23B5BD56D8857E6007463F5FE46B1FC0F6D36A5A96B7ABDDC58B7BE631A
A9046E5CF38F6ED55F55B71AC7C4AB1D36118B6D3A005C03C01D48FD40A00DEB9D7AEAD35DD3A09C43F66D41B6C4A3EF807A13FA56FDEDDC1616925C
DD3848906493FCAB89B822F7E24CD793B9FB1E8F6FBDB3D9B1D07E2734BE3A92EB53B6D0AC910C42FEE3E64DD9F4C73F439A02C5FBCF125F58E969AD
4D1C46CA66D91C3D1803F7589FC0E6B7B44BB96FF47B5BCB88D639274DE557A007A7E98AE47E25C88BA5E9BA141C3CF32051E8A3E51FCFF4AD7BCD48
D9E83AB471A22C1A7C4B6F1B0CF2DB0023F0C8A00C9F03AB6A5E2CD7F5B6DDB0C9E4479FAFF40ABF9D77B5C6782E43A5E81A5DB3C464B8D424795B69
FBAA79DC7F002BB3A1030AE0BE255E6A10E8FF0064290F937530446427763D0FE95DED705E39DF77E2BF0FE9EADC34A2420F4E1BFF00AD4304747E1C
5D561B3FB36A96F042B022470989B3B801839F4ED5AEEE91A17919554752C70053EB87F1EDDDE1D434CD32D36E27258866DAAE47627FCF5A61B9D447
AD69B2DE0B48AF6279CF4507AFE3D2B42B074082F9493A8E9563666350B1B4243337E35BD4088AE678EDADDE794908832715E65E2F86E17C55A0DCDC
643DC4C8CCBCE17E7181F502BB1D445C6B374D6B6170214B3706490AEE064FEEE3BE2B98F1B25EC5ACF8756F2E239FFD2948658F69CEE5CF7A4CA47A
45705E37D55AC2FCC69A84C4491E24B45524608C66BBDAA1A86996B7D04F1BC48B24D198CCCAA378047AF5A624719E02D4DAEAEE2B692FE44091931D
A91C30E7907F5AF42ACAD2742B2D32DADA28A3577B652A92B0F9B06B569206158BAE6BB6BA5C2E258DA6217748A3A2A9F5FF000ADAAF3EF1F44D6560
6D2DA5692EF55BAE9DC27F747B6481430476BA725A2D8C4D63124704837A8518EBFD6A8C9ADB7D9E5BAB7B4335AC44832F981738EBC547697761A7E8
8F697339096312C5331E3923A0F7AE345CCDA4AB6997C5E6D07547CC172BC3264E7AFAF038A02C7A4595DC37D6715D5B36E8A55DCA6AC553D32CA2D3
B4E82CE17678E25C0663C9EF9AB94C415977DACDBD9DC1B7D924D328DCE912E4A2FA9FC2ACEA5762C74E9EE88CF948580F53DAB92F07B39D0352D76E
98B3DD972A5B39DAB9007E7498D23A2D335EB1D4EF24B6B467678D37B12B818CE3F3AD6AF38F8732C16B1CB7371BCCBA85C98A1C0CF03AE7D2BD1E84
0D0553D4350B6D3A0335D3903B281927E82AE571BE2D8AEA29AEF529D6336905B148013D5C8EE3EB9A1823A39757B28B4F4BE6909864195DA324FE14
EBAD4AD2D2D96E2697113E0A9009C83DF1E95C66932CD69F0F27D4A73891E0D9103C855E9C7F3FC2A6D3EF27B7F876F7D7A7E7684C7083CFCBD07E27
AFE028B858EDA09A2B8852681D648DC65594E411555B53B25BC4B469D44CE4803B6476CF4CD67785C3A785ECED8488B76B6F9DA482549CE091E9CD65
5F9F37C61A56996E4B35A627B873DCF5C9FAF5A2E163AD9EEA0B668967902199C2203D598F6AA9AE185348BB9A68D58242D8DCB9C123FF00D558B772
FF006978E6CED95FF75A7A9918762D8FF3F954FE3799BFB105A46C55EEA458F81DBAFF008501620F042C7A7786219EEDC46D77296CB1E39E83DB8152
687AB9BEF13EA109BBDF0A0C451EEE0F3D40FA56A5CF95A578764DE5563B7B720F1919C7F8D637826D268BC2CB708882EEE4BBAB3F18C9C0E40E9C67
1401AFABDBE9DAA6DD2EF2522463BD511B0DC03DFE99AB9696D6FA6D8476F1B6D8215C02EDD07B9AE66FA57B9F186936518DD2403CCB99140C1217D7
F015D15D44F25C2F9C545A46BBDC1E77B7BFB0C6680265BAB7781E78E68DE24CEE646DC063AF4A7C13C5711092075743DC1AE4BC0F1BDCDCEA3A892C
B119DA389070B8E093EF4BE11D49AE75CD62DA35616E2532A67B64F4C741C5170B15BC6D0AC9ABE93656A8229AEA6FDE3A0C165CF7FE75DA492C36F1
832C89120E0162147D2B96842EA7F10A4933BA3D3A22006ECC78E3F1CD1E2DB99AD6DAFE6BB50233188ECFE6041638CF1D73401ADA2E8D169935CCF0
DD49325C90541E8A3F0EBD7AD69DC5CC36C9BE7915076CF7ACED1257B4F0BDA4B7A583456E0BEEEBC0ACFD0659F55B7BCD558069A4DD1DB06E420038
C0E9C93FCE803A48E449503C4EAE87A329C8348CE8AEA8CCA19B3B413C9FA5735A3C93681E1AB8B9D51807F35D828E80FA7E60D452E5FC2773ABDDE5
EE255F354E79419C003D28B858EBAB0F51D4E787C4163A7C0C8124F9A5C8E719E950F82125FF00847639E73F35C3990739E3A0FE554ACECA1D6BC577
F75728AF15A95488AB771FFEAA00EAAE10496F2465B6EF5299F4CF159DA069736976B2453DC99CBBEE1E8BF4ACDF1ACDF63D1E7B979C96628B0479C6
D39C138EE79AD6F3AE9347B7644692E5E34078CE1881927F5A03A1A01D59994302CBD403C8A010464107B715CAD94BBBC6975676AD20896DCACCF9C9
DDC739F5A9740126953EA96F733B496D0157576E31C73F8FF85170B1D3D15CA1D4AD7528AEA69DCB30052DE0527AFAF1DF38AD9D08DC9D220FB68C4B
8EFD719E3345C2C6954722EF8D93732EE04654E08FA5494D62141662001C927B5311CBD8DC5EB78B5AC96E6692D600CCE19B3DB001FC4FE95BFA8473
CB632A5A3EC988F91B38E73EB58FE1C4125EEA7A932B0134985EFC0E78A7EB7A94B6874F922F311E5976F9248FBA7FBC3E9F952E853DCD0D1A2BD874
F44D4A5125C64E4839C0F4CD68546C37298F2CA594F2BDBF1F5AC6D12F4CF737F64D33BB432B6C727242F007269926ED158BA15E3DC4B7D6D24AD28B
59B62C8DD581CF7A9AF6FD52FA2B1491519C6E91CB01B17DBDCD006A5158FA64F236A9770A4AF35AA2A95776DDF37700D4E6F45D5FCB65012046A7CD
914F20FA0F7A00D1A2B9BD3EF6FA5F11BD9B4CCF0DBAB2B657AFA13EF5A1ACDEB588B690310A64C30FEF0A0763528AC566D53FB761641BAC587504E0
0C771EB5B54082A39E54820799CE1514B1A92B03C572CB0696EE93108E4465028393D739A0058CEAB733DA5DC457C890FEF230701573CF07AF7ADEAC
8D2E1B8B1B3696EAE9A58BCB0CA8460AFB54175A95D41AB58447212E4E1A1201283A673D7DE90CDEA2B2EF6FA45D56D2CEDDBE663BA50003F2FF009C
D5F98910B90C14853862338E3AD3112D159D617FF6ED38CA92209B69C900ED53CE3AFD29FA4DE9BFB113300183156007191E94017A8AA66F164BB6B4
8082E8B976EA17DB1DCD56D2753FB735D421A391ED885DCB901BAFF85005BFB741F6FF00B16E3E76DDD8C71525D492456D2490C7E63A8C85F5AE7B4D
7BABFD7E6BC448E348FE470DCFB71EF8ABDA8EB2B691ACF16C92013794E73C93C74FF1FF00255C762DE957925FD824F2C2D1392415231D0F5157A98E
EB1C6CEE70AA0927D0566C97D7B259B5D5A430F94A09024277301DF8E94C46AD15534FBC8EFED23B887EEBE723D083822A0835306FA4B4B98C432295
DA036ECE7A0A00D2A2B3E1D433AA4B612A05751BD0839DCBEFEF57269A382269656DA8BD4D0049514F29863DE22924C7F0A004D65DFEB6D61179B359
BEC619421BA8F7E38AD3B49C5D5AC73AA95122E707B50052B4D6AD6EEECDAC4B28946721940C7EB5A75CD7876269B55D42FDC0C33ED527AF5AE96921
B0A28A298828A28A0028A28A0028A28A0028A28A0028A28A0028A28A0028A28A0028A28A0028A28A0028A28A0028A28A0028A28A0028A28A0028A28A
0028A28A0028A28A0028A28A0028A28A0028A28A0028A28A0028A28A0028A28A0028A28A0028A28A00A37BA5DADECF1CD32B6F8C60156C647A1F6A65
D68B61762DC4B0FCB6E308AA7008F43EA2B468A00A77D636F7DA7C9613A91048A14AAF1C7A0FCAA4B5B382D2CA3B385008234D814F3C558A2803374E
D12C34D9A49ADA221DC9C166CED07B0F414BA8E8F63A94D14D7711678F80436323AE0FB568D14015EE6D20BAB37B4990185D76951C7151E9BA6DA697
6FE459C7B5739249C93F535728A0028A28A00A577A6595DC9E65C5BABBE319C9071F855A8A348A358E250A8A3000ED4FA28028DB695636B74D7105B4
71CCD9CB01EBD71E94B7BA6D9DF953756E9295E0139047E22AED14010DB5BC36B0AC36F188E35E805471D85A457B25EC76E8B7320C3C8072455AA280
2BC5676D0DCCB731428934D8F31C0E5B150DCE97A7DD5D25C5C5A4524C9D1D979FFEBF4ABD4500473451CF1345322BC6C30CAC3835059E9D676258DA
DBA46580048EF56E8A00AF3DA5B5CCB0C93C2B23C0DBA3247DD3EA3F2A2EED2DEF60305D44B2C648255BD6AC514011450450C0B0451AAC4A30100E31
50D9E9B6760D235A5BA44D21CB11DEADD140145F4AB093524D45ED91AED1768908E7157A8A2800AAB77616B7C105DC2B284395CE78AB54500650F0EE
8C24327F6743BCF56209356ED74FB3B4777B681237718661D4D5AA2802B5DD95ADE84177047288D8328719C1A75DDA5BDEC061BA85268C9CED71919A
9E8A008ADE08ADA0482DE358E24185551802A92E87A525C0996C621206DE08CE01F5C74AD2A28020BBB586F6D25B5B85DD14AA5587A834CB5B1B6B4B
24B2B7882408A5420E983D7F9D5AA2803334ED0B4ED36779ACE0D8EC300939DA3D07A569D14500154EF74DB3BE286EA0590A7DD39208FC455CA28033
2EF40D32F2CD2D66B5531C67298E08FC6AFC30C76F0A430A048D061547402A4A2800AE775BD7ED134BD652372B359C7B18371966040C7AD74554EE34
CB1B95904F690389082F9419623A64D0062F8034C6D33C2D6CB2AED967CCCDEBCF4FD2B6753D36DB55B26B4BB42D13107838208E8455B550AA154000
0C003B53A803165F0C6992E9434E68DFCA0E24DC1BE72C3A1CD685858DBE9D6A96D6B1848D7F327D49EE6AD5140193A1E876BA24770B6CCCED712992
476C649F4FE7F9D451F872C93569750DD231964F34C44FC9BFFBD5B74501733B53D26DB51314926E8E780EE8A68CE190FD7D29D61A6C764F24A64927
9E5FBF2C872C40E83E957E8A00C297C316726B12EA1E64AA272AD3420FCAECBD09AB9AD6930EB1642DA7764D922C88EBD5587435A3450172969FA6DB
E9F6CD0C40B1725A477E4B93D49AA47C3965E5490A49711DBCA77342926109ADAA280228218ADA148608D638D061554600ACED3B43B5D3F56BFD4A32
CD35EB02C5BF87D40F6CD6B51401812F862D5F599750F3A45139065857EEC8474C9EB5A3A96976DA9087ED01C3C1209227438646F506AF514018BA97
876D351B310C8F20956412ACF9CB861D0E6A45D0ED46872692DBDA2954891C9F9998F25BF3AD6A280B98BA2787E2D24AB1B992E5913CB8DA403E45F4
15B545140113CD146C8B24A88D21C206600B1F41EB5C7A42752F89CF3ED0D0E9D005273C0623FC49FCABA3D5F488F5336EED2BC52DB3EF8DD3A834BA
569506982668C99269DCBCB2B7DE7348668D66EB5A2D9EB56E915E2B6633BA37538643EA2B4A8A6232748D17FB325964FB75CDC79800DB23640C559D
42D2E2E9156DEF5ED719DC51412D5768A00E7F41D0EEF479641F6FF36DE462ED194C92C7DEA9F893C2973AEDE4139D4FC9FB3B168F10F29CE477E6BA
CA281DCCDD32D2FEDCB3DFDF0B966500AAC7B4023B8E6B4A8A2810514514005739E26D12EB50BDB1D434E6845D5A13812F420FA5747450073F71E1C5
B8D026B0966DD3CEDE649291F7A4FF000AAD3E8779ABDA5AD86AB15B436D6D22B9F218E5F00818F4F7AEA68A0770A28A28115EF6D92EECE6B67242CA
85491DB3DEB9ED3743BD5D2D346BC31A5944AC3745D66E78CFA7AFBD75345160B9C77867C3179617100BE31793625BC809FC64FF0011AEC68A28481B
B85737E36B1BED4B498ECEC622E24941948EA1457494500729E26D22EA4D02CF4AD3E0124219525C71803BE3F3A9FC4DA54F3F87E0B4B2877AC0E85A
04E3728EC2BA4A29587730FC3F6532B4FA85E42619E7C2AC67FE59A0E82A93E9B796DE26BDBCB5837FDB5511662788863E63F5E98AEA68A2C1738DB2
D1EEECB5BBF8961668AF36A9B91D931F30E7BF5AB3E28B49D67D2EEA085E6B6B2932D120C9E831D7A8E2BA9A28B05CE6F54B3BDBBF0A5FA346ED7573
FBCF2BBAF230BF90A9F458A599E39A4B496D62B7884304520C1F738CFB0FF22B768A61738B363A841AE6A5F67824F3EFDCAA4C4663488F53F5E95BDA
A95D3BC357013188A1D99C7AF19FD735AD5CA78B5B51BCB33676DA64F22F9A0B32904328F4FAFA520DCA9E1E9EF34CF0F9D3D6D64378FB9E02AB9560
7BE6B67C2FA08D16D65DEE5EE2E183C87D31D055FD26669AC222F6B2DB15013CB90608C0ABD403670D616F76979AAE9DE4C82EEF27CB4E54ED11773B
BDFF00AD6878C6D525934C92E23964B4825669420DC4F0319FC8D75345160B9CF436335FE917CF3298A5BC1F22765551F2002AAF862FA2D3F48874B6
B79FED70EE0D16CE4F279AEAE9BB577EFDA376319C738A2C17399F14E9334BE1C8EDADB74A2170F22F7619C938FCF8F7A8359964D57C35241A6444C5
0C4BB89183C0E807B575F45160B9CB786849333DD5BC5B2186D560894820330E49FCEA8F86A593F7567136DB992E4CF72CB90368E8B9F7C1AED55551
42A28503A00314D5862476748D15DBAB05009A2C17393F1618975CD31EFF0077D8A3064C1076961DBEBD2B734549FEC8D7172CC5EE1CC8149E1476C7
A71FD2AFC914732ED963475CE70C0119A968B05F438AF0F5CDBE9736A72DE16FB749360447EF3019C01F9D685FE9F3BF86AF1C83F699FF007CE31FA7
E55BED6D034C266822328E8E506E1F8D4D4582E72BA45D58D9E896C90C71CF78140D817E62D5B9A55B4B69A7C514EECF20C9624E793DAA74B5B78E53
2C7046B21182C140353D09036154759791348BA6886584678F6EFF00A55EA6901810C0107820F7A6239CD026D90D859D9C81800D25CFF101907033F5
C543A9B47FF09319EE48D96B10F290FF00131E471F535D2C16D05BEEF22248F71C9DAB8CD249696F2DC473C90A34B1FDC723914AC3B950DCCD6BA1B5
CDEBED99232C495E87B703F0AC7D0ADCD9E8175A96C63713AB36475C76FF001AE9E48D25429222B29EAAC320D22C51AC5E52A208F18D80718FA5160B
9CC7872396DE0B786DDB74B3BF9D72C4745F4FAFFF005A8B38A0B9F11EA7FDA4A3744C3CBF31B1F29E98FC315D25ADA5BDA2B2DB44B1863938EF55AF
349B2BE9966B88419171F3038247A1A2C17225799ADAE9AC62458B6E210ABB771EEDFE7D2B1BC23710C369793CC7649E6ECD84E5CE3FC735D622AA22
A28C2A8C01E82AA9D32C4DDFDA8DB47E7E776F03073EB4582E735A4CB37DA2FE68C32DECF3840BE83A938357F5E89B51D5AD34F460A1079AE7DB3FFD
6ADA4B2B58EE5AE5208D667FBCE07269DF6683ED3F69F297CEDBB77E39C5160B922A8550ABC003029D4514C415CEF8B03CA2C6D909CCB3741EDFFEBA
DBB9B986D23124EFB109C6704F359655750D7A39A365686D57EF03904FF9FE5498D1B2AA15028E80605738A525F145C5DDC3652D57644A79E4F1C7E2
0D74B555EC2D9EE7ED2D1033641DD93DBA50239AD32E5DF5DBDBBB856DEA024683F898F61DF15BFAC5C7D9F4899DB01D9760C73C9E2A41A6DA2DEB5E
0817CF6EADEFEB52DD5AC377108E75DCA1836338E450339D58DF46F081C12F2CA07CD9C1F9BFCFEB5A5A529B5B04B68A3DD2087CD723A6E3D07D7FC2
AFDC5AC3736E6095018C8E838C7D292D2CE1B342B086E7A966249A00C6F09179E1BBBB9FFD6C931520F25715A5726DB4EB2BAB9B78910E09628A3961
C0CFE3FD698748845CB4B0CB2C4921CC91A3615EAF4B0472DBB40EA3CB65DA40F4A00C1F0AB886C234600C9752C9264761EF506AD14336B16B62A891
C16E7CF978C633C922B674BD263D354859A494E36A96FE15F4145E69305DDC2CEF24A8C06D608400C3DF8A02E67EBD7AD3E82B2C2ACAB2CC139EEB93
CFE95A376CB67A249E6363CB876647AE31FCEA7B9B586E6D5ADA541E591800718F4C557934C595A3134F34B146B811311824742718CD3029F852CA6B
6D1916E472D2B4A833D8F4FF001AADA7C3F6BF135CDE03BE38F807048F6FC6BA50000000001D00AA9A7584761132A317676CB31EF4AC1739FD36749B
57BED5E5CEC8BF74AABD4F6039EBC54DAB5D1BDD62D74EC3A28657607BE4771ED5A367A347693656666883EF58F18C37624F7A96F34E1717515DC329
86E23C8DE17391E845160B985E3374B87B4B28CFEF8BF66E80F1D3F5AD4BBBE822F0FC8F1306544F27BF5C608FCA96E3438EE3ECF23CC7ED30B64CC5
725FD88FF3D2A5BDD1E0B8D28D8464C6A3953D79F53406847E1C410E9912B6034B9900F6E07F856BD54D3609ADAC921B8916474180546001D855BA62
0A28A2800A28A2800A28A2800A28A2800A28A2800A28A2800A28A2800A28A2800A28A2800A28A2800A28A2800A28A2800A28A2800A28A2800A28A280
0A28A2800A28A2800A28A2800A28A2800A28A2800A28A2800A28A2800A28A2800A28A2800A28A280296ABA8C3A558497771B8A270028C924F414FD3E
F22BFB18AEE10424A3201EBE9593A978874E8EE0D9DC4465859846EFD573E9EF5AF3136B66E6D6D8398D729126173EC290EC67F88F5B1A35BC3B2133
5C5C3F970A7A9A8FC37AE49AC2CE92DA3C324070C718527D2B86BBD4350D5BC7F6CE9A7B34F658DB6864E148CF24F4182457A6D979A6D51A78160998
66445208DDDF91403562CD14521200249000EA4D310B45737278A519EE65B4B633595A67CE9F7607E156F43D7175AB69EE6381A28226DAB23370FEA7
A718A571D8D9A2B061F12C325DC31BC263827DE62999BEF05EF8C74A647E240DA8D9C0F6AC90DE92B0313F3360F523B0A6163A1ACBD4AE7528250D67
671CD084DCC59C020E4E7BFA559D46FED74CB292EEF2411C318E49FE42B96D6BC4D736BA04B757767E4C179195B670DF30DC0E370F5C5008D4F0C788
0F8816E24107971C4C0290739EBC1AE82B97F87961F60F09DB961FBCB92677FC781FA015D4500C28A2B062D727D46FA5B7D1E28E48E038926909DB9F
418A046F5159BA7EAB15DDB4F24988E4B66649D7AEC23AFE154E1F10A1D40417317911CB099E3663CEC1DC8F7C1FCBBD01637AB2B55D46EB4E592E16
CC4F6B1C5BDD8481581CFA1F6A82CF5C336A915A4D01896E6332C19FBC547A8F7EB557E205E0B3F085E75DD3E215C7A93FE00D032EF86B5A6D76C1EF
3C8F2A2F30AC7CE77015B35CF6941340F09D846633249E5A2A46839791B9C7F9F4A9746D665BCD4AFAC2E92313DA11931676F3DB93D7A5006E51581A
0EBAFACEA1A8469120B7B49362B8CE58E7FF00ADFCAB7E810564F88B556D1F4B37690999B7840A067AF7AD4660AA5988000C927B572D2F89D9F51B38
05AC6D67793F928EC72C70719C7D71498D1BFA5DE35F69B05D344626917250F6E71581AB78A2EADEEAF134DB25B98AC595276E73B9BB0C56C6B57975
61A749736B0C52F96A59F7C9B7000EA38E6B91F86B26A525ACD73E446F6D7772CF24CCF86538E98EFCD008EEEDE4F3ADE294A942E81B69ED91D2A6A2
8A6232B5DD5E3D22D1642BBE69582451E7EF13FE157A59D2DED8CF707CB5500B77C7B71D6B9504EB1F10DD5C86B6D26205571FF2D0FF009FD2AE78E9
6F1FC3E4E9EAED324C8D841938E7FF00AD48762CBEBE21559EE6CAE20B566C79EE3A0EC48F4AD78278AE6049E0916489C65594E411587AA5D3C1E099
AE3535549BEC9F3ABFF7F1C0FAE7F5AA7F0E05E7FC232AF77908EE4C408C1DBFE7F9501637B4FD462BD9EF21446592D2531B83DFD0FE38ABF5C9EB73
B693E30D32F14ED82F54DBCDCF048C609FC0FE95D65310552D46EAE6D96336D66D73B8E182BE36D5DAA7AADD0B1D2AEEE8F1E4C2CFF881C50066E81E
221ADDCDCC2967243F66E1D99B23767A74ADEAF38F04EA474EF0D5CDFBC0F3992532DC387C6D19C0FC6BD0629925B749D4E23740E09E38233490DA26
AAD7D70F6B6AD3476F25C3023F771F5359ABAF8B969CE9F6AD731419F324DE1471D71EB5A1A75F41A8D9C7756C498DF38CF51838C1A62322CFC5315D
6B4BA50B39927E770246531EA2B7A6631C2EEA85D95490A3A9F6AE1BC128351F166BBACB28C093C98FDB9E7F90AEA353D663B3B98ACE189AE6F6604A
42871C0EE4F6A436676A1E2D5D36289AF74CBA85A5242AB639E9DFF1AE8A2904B12480101943007A8CD79E6BB38F1278CF46B088B8106E7B889BFE59
B039391F8015DAEA9AA45A7F9510432DCCEDB62841C163F5EC28034A8ACBB3D6126BF3A7DCC46DEF026FF2F70618F622B5298828A28A0028A28A0028
A28A0028A28A0028A28A0028A28A0028A28A0028A28A0028A28A0028A28A0028A28A0028A28A0028A28A0028A28A0028A28A0028A28A0028A28A0028
A28A0028A28A0028A28A0028A28A0028A28A0028A28A0028A28A0028A28A0028A28A0028A28A0028A28A0028A28A0028A28A0028A28A0028A28A0028
A28A0028A28A0028A28A0028A28A0028A28A0028A28A0028A28A0028A28A008E58A39A3292A2BA9ECC334DB7B682D90A5BC6B1A9392054D450014514
500145145001451450014514500145145001451450014514500145145001451450014514500145145001451450014514500145145001451450014514
500145145001451450014514500145145001451450014514500145145001451450014514500145145001451450014514500145145001451450014514
500145145001451450014514500151C8A5E265562A58101876F7A928A00F33D4AC966F12E99E19B491A44B67F3EE25EE58F249F7C7F3AF46B7B886E5
19E1903AAB152474C8EB5CC37872F62F135E6A16CF16CBD5DAD2B13BA3078200FC056BDE4F6DA3E8775E532E2D2127696C9048E33F534914F5391F04
037FE35D6F54DB98C332AB8E992D81FA0AF44AE27E1769EF69E1D7BA90306BB937007FBA3807F535DB5084F70AE3BE23EAD3E9BA02C36CA77DE3F97B
C7F08EBF99AEC6B2B5FD1A1D6EC05BCA76BA3078DFFBAC2982396D5923F0F7C2E16AC57CE9A20A4673B9DCE58FE44D5ED02CE35F0DDAE841E4496E6D
4CB2CB19CEC0C47E591C558BDF0CCDAA5BCABA9DD23BF93E5401170A87A83F5AD2D17497D3C34B7532CD72EA1199576A851D0014877391BDB74D4BC7
F63A4DAFCB6BA6403CCDBC7B91FCBF3AB4646D5BE2391F72CF468C9627006EC7F8FE82B67C3BE1D6D2B50D4350B9B9171737AF92C0630BD7155DBC2F
72756BD945EAAD95ECA24990290EDFECE7D280B9CEF8B2F8F88B5DD1F4C84BA594ED9258637F383F87153FC4B962B9934AD0A160249250CC07F0AF41
FD7F2AE9B55D0FCFBBB3BFB0F2E3B9B3CF96AC3E4208FE7D2B2B50F07DC5CA4778B7A8755171E6B4E578C70303D8014B50BA36BC3BA9457D0CD6D0C6
512C585BE49FBD818FE95B359DA26951691A78B68D8BB925E490F5763D4D68D52259CFF8DEF66D3FC277D3DBEE12ED0808ED92013F9553F04C69A5F8
263BA9C81B95AE1D8F5FC4FE15D25E5AC37B672DADCA6F86552AEBEA2B024F0E5DCDA7C7A5CB7FFF0012F4E0AAAE199474078A435B193E05B39B53D2
F52BBBE6221D467DE02F04E18927E9DAA2488EBBF122608C3EC5A74611C0E8D8FE1FA66BBA82D63B5B45B6B5511A22ED40074AC9F0C787FF00B0E3BA
7966F3EE6EA53249263B7614587730AC2E16FF00C79A9EA8E716DA4C4601839CB720FE279A4F1ABB6A9AFF0087F470AC239A413CAA46081EFF00401A
B5AC3C332596A37322DD66D27B8370D101F33B1E70DD88CD4BAD6877177ADD9EAD61711C3736F1B47FBC5C8C1CF4F7E4D0059BA9ACE6325D0CB9D277
9DA07CBBB67F4AE3FC3B7074CF046A5E21B820DC5F3B94CF5EA401F9E7F215D65DBE9DA26902D6F256D97398D9B8DCECC304D735E26D28DA68367E1F
82E9AE1EE2E15608F6005173F31E3D39FCCD008D6F87564F6BE178E595712DD48D331CF5CF4FE55D5D57B2B65B4B282D90E5618D501C633818CD58A6
4B21BAB75BAB592076655917692A70715C0DBC11DE7C405821678EC3428B27272370FF003FA57A257343C305355BBB88EECADB5EC8259E3C61891DB3
E9498D078B2FE36F045F5C82D18961C20718249381F9D4DE08B616BE12D3D026C668F7B0F726B33E22869B4AB0D26DC0DF7B7691851C7CA3FF00AFB6
BACB68560B68A04FBB12041F40314750E84D45145311C6784379F15F88DA51B64338F97DBD6BB1660AA59880A06493D0572D1A1D33C7B2BB0220D4A1
00311C6F1DAB53C43A65DEAB63F65B4BDFB286389085CEE1F5A486CAF2597F6FCA92DF2674F8D83C5093FEB08FE26FF0AD4B09AD66848B32A62898C6
36F4C8EB8AE462F056A912954D75C7181F2FB574FA06949A3693158A4864D992CC7B93D6840CE6FE243158F47F2C66537A368F5E3FFD55DB571DABA1
D63C75A6D9C7CC3A6FFA44C7D1BB0FE55D8D00C2B90F8977A6DBC2B2409FEB2EE458940EA46727F90FCEBAFAE77C45E1F935BD4F4D95E5516B6B26F7
8CF527FCE2982397D66D5B49F04691A145C5E5F48BBF07A93CB7F303F0ADCF1A5ECDA37831E3B66224012DF780460606483F4FEB5A1E21D127D4EF74
EBCB49A3866B276605C673903FC2AD6B1A526B1A2C9A7DD3E4BA8FDE018F98739C76A41732B4945D0BC00B2CFF00785B199F3FDE619FEA0564F85679
F4BF87579A85DB101CC92C608C601C01F99ADBB8D1752D4ACEDEC352B8B6169191E688436E940E839E957F5AD1E2D4BC3D3E91091024918542A385C1
047F2A02E61F81228F4AF07C77370CBE6DD179CE4FCCDEDEFD3F5AA3F0E77DFCFAAEBD7646E925F2C16FE1C004FD3A8AE8745D2EEA182DD75131E6DE
DC5BC71C44EDC63058FB9C551D2F41BFD32CA6D222311B19A76769B3F36C6EAB8F5E31F8D0330FC192C13F89B5FF0010DC3AA40ADE5AC9D86E6FFEB0
AB7E2B6FED4D6603A15CEED574CCB98B1C3AE7A03DC8C1E2AD68FE16BBD3EE65B4FDD8D34DCFDA3767E7931D14FB55B4D16EB4CD7AF351D321B7985E
E09121C188F7C1F439A00A9E0EFB2EB5793EBD2ACB1EA287C99632DC2F1D71DB3E9ED5D9D63787746FEC886E0BBABCD75299642830A0FA0AAB6FA76A
CBE2F9AF669D7EC0508540DD78C018FF003D28133A3A28A298828A28A0028A28A0028A28A0028A28A0028A28A0028A28A0028A28A0028A28A0028A28
A0028A28A0028A28A0028A28A0028A28A0028A28A0028A28A0028A28A0028A28A0028A28A0028A28A0028A28A0028A28A0028A28A0028A28A0028A28
A0028A28A0028A28A0028A28A0028A28A0028A28A0028A28A0028A28A0028A28A0028A28A0028A28A0028A28A0028A28A0028A28A0028A28A0028A28
A0028A28A0028A28A0028A28A0028A28A0028A28A0028A28A0028A28A0028A28A0028A28A0028A28A0028A28A0028A28A0028A28A0028A28A0028A28
A0028A28A0028A28A0028A28A0028A28A0028A28A0028A28A0028A28A0028A28A0028A28A0028A28A0028A28A0028A28A0028A28A0028A28A0028A28
A0028A28A0028A28A0028A28A0028A28A0028A28A0028A28A002B8DBFF0006DCDEEA9752B6AAE96574DBA4842F24673B49EE335D95140D3B1059DAC3
656915ADBAED8A25DAA3DAA7A28A04145145001451450014514500145145001451450014514500145145001451450064EBBA15AEB70C4972595A07DF
1B2F63EE3BD169A32C5A8FF685D5C497372A9B119C6020EF815AD450170A28A2800A82F6DC5DD9CB6E5D904AA5772F5153D1401CEE9BE1B682FA1BAD
42E9AE9AD176DB8EC83FC6BA2A28A0028A28A00A7A869F06A11A24FB818DC3A321C3291E86AE514500145145006568DA2C5A5B5C4BE6B4F7170FBA49
5C727DAB568A2800A28A2800A28A2800A28A2800A28A2800A28A2800A28A2800A28A2800A28A2800A28A2800A28A2800A28A2800A28A2800A28A2800
A28A2800A28A2800A28A2800A28A2800A28A2800A28A2800A28A2800A28A2800A28A2800A28A2800A28A2800A28A2800A28A2800A28A2800A28A2800
A28A2800A28A2800A28A2800A28A2800A28A2800A28A2800A28A2800A28A2800A28A2800A28A2800A28A2800A28A2800A28A2800A28A2800A28A2800
A28A2800A28A2800A28A2800A28A2800A28A2800A28A2800A28A2800A28A2800A28A2800A28A2800A28A2800A28A2800A28A2800A28A2800A28A2800
A28A2800A28A2800A28A2800A28A2800A28A2800A28A2800A28A2800A28A2800A28A2800A28A2800A28A2800A28A2800A28A2800A28A2800A28A2800
A28A2800A28A2800A28A2800A28A2800A28A2800A28A2800A28A2800A28A2800A28A2800A2B33FE279CFFC8387A7DFA6EED7329F2D8F23E6E1BE5FD6
803568ACC55D6C7DE7B03C63EEBF5F5A6C31EB6A4F993DA38EC0A9FE98A00D5A2B2E58B5B76CC7716918F40A4FF3A67D9F5DC7FC7EDB67FDCFFEB500
6BD158E2DF5DC737B6D9FF0073FF00AD4DFB3F883CCCFDBED767A797CFF2A43B1B54562ADB6BE33BAFED8FA7EEFF00FAD47D9B5F07FE3FED88F78FFF
00AD4058DAA2B196DFC419F9AFAD7F08FF00FAD4E16FAE83FF001FD6C47FD73FFEB50235E8AC416BE20DA01D46D89EE7CAFF00EB53BECDAF1007F685
B0C753E4E7345C66CD158BF66F1067FE4216D8CFFCF2FF00EB529B6D7B68C6A36C0E79FDCF6A00D9A2B196DB5F1D751B73FF006CBFFAD41B6F106EC8
D42DB1E9E5FF00F5A811B345639B6D778C5FDBFD3CBFFEB53A2B7D686EF36F603E8447FF00D6A606B5159715A6AA6547B9D4176A386291C63E61E84F
15A940051451400514514005145140051451400514514005145140051451400514514005145140051451400514514005145140051451400514514005
145140051451400514514005145140051451400514514005145140051451400514514005145140051451400514514005145140051451400514514005
145140051451400514514005145140051451400514514005145140051451400514514005145140051451400514514005145140051451400514514005
145140051451400514514005145140051451400514514005145140051451400514514005145140051451400514514005145140051451400514514005
145140051451400514514005145140051451400514514005145140051451400514514005145140051451400514514005145140051451400514514005
145140051451400514514005145140051451400514514005145140051451400514514005145140051451400514514005145140051451400514514005
145140051451400514514005145140051451400514514005145140051451400514514005145140051451400514514005145140051451400514514005
145140051451400514514005145140051451400514514005145140051451400514514005145140051451400514514005145140051451400514514005
145140051451400514514005145140051451400514514005145140051451400514514005145140051451400514514005145140051451400514514005
145140051451400514514005145140051451400514514005145140051451400514514005145140051451400514514005145140051451400514514005
145140051451400514514005145140051451400514514005145140051451400514514005145140051451400514514005145140051451400514514005
145140051451400514514005145140051451400514514005145140051451400514514005145140051451400514514005145140051451400514514005
14514005145140051451400514514005145140051451400514514005145140051451400514514005145140051451400514514005145140051454733F
950492633B14B63D70280303FB5A7BDF167F665A16482D06F9DC2E431C7DD3E955FC57A96A3A75ED9A432AC704F2AAA951CE73CE7DB9AADF0D97CED3
EFF5073996EAE493939200E83F5A4F18D96B3A86A968B6DA779D6F6CE24470E064F19CFE549EC575367C512DF5AE9AB7763215FB3B879147F12FF9ED
5A763751DED94573111B2550C3073F853260F2E97209D02C8F01DEA0E4025791585F0F5CBF86F6B387F2E774C8FC28EA2E8753591E2AB992D3C377D7
10C8D1C91C795653820E456BD41776B05E5BBDBDCC62485F1B90F43CE6988E2DFF00B56CBC0EFAC5CEA5726E96212AA9738C13C03EBD6A745D42E7C1
50EAB05FDDFDB045E690253B5B9E463E956FE20848FC19751A0DABF2A80A38033FCB8AB9E11543E13B05C028D1742382093C522BA191AF6A17507836
CF58B7B99D6511461C0623796C024FBE73F9D753A633B69768D2B1691A142C4F52768CD731F12716BE0E09080AA9346AA3FBA0671FCAA86AD6A9A3F8
1EC2FAD9DC5DA985FCEDC724B7247EB406E75FAD412DCE9734105CADBCB20DAAECD807D8FE19E94ED1ADA4B4D22DADE69565911797539079CF06B27C
530ADDF841EEE4244D0C0265707183819ACF7BC9E2F86B6D711C844AE8A31DE4CB72A3EA2811D8473C3280D14A8E09C02AC0F3D7153579C5A6B3A6DE
6B1A6C76914D633472AF9D0CBC2F5FE75E8F4206AC158BE22D59B4D16690B2F9D3DC2A6D23395EF5B55E71E32B0B78FC69A32ABC8A2F65CCA37138F9
8723D28608EFEFA196E2D5E382768243F75D6B9DF036A57BA9417EF7B2B4BE5CDE5A93D38CE715BD696D6D6EF28858B39DBBF7396230303AF4AE2FE1
F5A7DB2C358F364630BDD322A863807AE7DFA8A03A1E814C7755C6E603270327193E95CD78275396EB43B86BB76736B33A02473B074FEB583A56AB6B
AEAEA177A94B31943EDB7588390A3B638A2E163D1A9BB977ECDC376338CF38AE3B4FF105CDAF822E6F2EC32DCDA31886F1F31E9B491F43FA525BE9FE
77851F5769E537ED09B813072082A0903F4A2E163B4A2B94B5BD7D77C16F74D2BC53451B1668DBEF32AE79F6AC6F0FDBEB1AF6856D7C97ECAF1DCE10
331E10601C7F9ED45C2C7A2515C36AFAD4573E2D6D2AF1E55B38137148872ED8EA7DAAC786B58C6BB71A5A79CD68F97B5691483C0C91F4A2E163B1A2
B8483ED9FF0009FCFA4477F3ADBAC22773BBE66E9C67D7E6EB5344D769E2F9B404BC9459BC425EB9651D700F6A2E163B5A2B93D3AEE5B2F1A4FA22CA
CF6BE489115C92578CF07EB9FCEBACA620AC59E0D5CEBF0C90CC069F8F9D7238F6C56D570E4CB6DF13A1B48A67F2258DE764662464A9CD0347714570
BE22D46D62F13FD93559E78AD9A306278E4C046C7391DAB63C1ED2CB69733BDD8BA8E49CF972862770F5C76A4163A2A28A82F0CCB6731B550D3043B0
1F5C714C44F595E223769A25CCB6339866890C99033900648F6AE3AF6792D3C35FDA52DD347AA4726597CCDDDFBFA1E6BB2B894DC786259A419692C8
BB7D4A6690CA1E06D527D5B40FB45CC864916668CB1F6C7F8D7475C3781A46B3F878D71083E62995874FBD9C679FA554B4D52C6E61B61A5EA722DD49
322CD048D80E739C8F4C5007A2515CF78835678353D3B4BB790C725E483738EA17DAAAEAF7B3F87B55D398DD4D3DB5DC9E53C7236E23DC7E74C2C757
45145023335DD44E996026455691E458D0374249AD3AE1BE224727DAB466499C2BDD2A6C07001CF5FAF35D6ADAC89672402EA72EE0FEF5882CA71DA9
0CB945707A56A1ADEA965A9C76D39924B2982A36EC173DC7D315A3E27BDD4B4DB0B0BD5B931C85E38E68C636927A9A2E163ABA2B94F19DDDF69BF64B
AB5BA758E69D2168C718CE79157BC53AD368D631345C4B3C823562385F5269858DDA2B8F9B573A5DCD8491EA46F22B860B32336EDB9EE3D3BFE55D85
2B8AC158FE20D59F4C8A05823596E2E1F6468DD0FAFF00315B15C078F52F3FB6B438D6E7E49EE488CE00319CAFE7D698D1D6E85A8C9A9E9FE7CF6E60
90394643EA2B4AB1355D4FFE11FD1165B994DCCC3E4566014BB7A9154353BBD4B45D393547BBF3E32CA6589D460027B1EDD69058EAA8AE6FC45AB5C5
B787D757D32E1563501995901DE0FF00234BA55F6B176FA6DDB46AD67750EE90003E4E3839F7EB4C2C74745145020A2B97B2D76FE6F154FA3C91DBB6
C8F7EF5CE13FC7B7FF005AADE8DAB4F3EAD79A5DD8469AD86EF3106030CFA7E34058DDA2B0EFB5699B5B5D1EC762DC18FCC777190A38E945BEB1245A
D0D23501179CE9BA3910E03F193906803728AE6AFB5FBBB1F10C1A6BDB4522DCF316D24123A0C9E9561359B88BC4716917D1C20CE85E278F3CE067BF
D0D006ED14553D4EFE1D36C24BBB83F220E9EA7B0A00B94572B7FE21BBD392C6EEE561F22F080230082071CE7D79AD7D6B554D2A18CED0F34EFB2242
7193FE71401A74573EFABDE69977026B11C4B0DCBEC49233F718F407DAA6F136AD3E8B64B791451C9106C49BB391E98C7E3401B54573B37882E2DF4B
83539AD116D24C6EF9CEE19EFD38A2E7C4335BDA47A935983A73E312799F3609E0E2803A2A2AA5D5CBAE9AF736AAB2909BD01E030EBFCAB2B45F10FF
006B6993DD436A5A58490C81B68C8ED93D0E3F95007414561D8F8806A509934FB396511E44A0B05D84761EA6ADE99AA2EA9A69BAB68DB702CBE5B1C7
CC3B6680B1A34562695AF8D46F2EED3EC924535A81B949CE4F3C738F4A9748D6E3D4EE6EADBC89209AD480E8FD680B1A4D222C8A85C067CED527938E
B8A92B91D5B58B887C5D616A6C24631A48F185604C9904574525F25BE9C6F2ED1E1554DCCAD8247B71401728AC46D792316B2DCC063B7BA388E50FB8
7B6462AEEB17C74ED367BB113C9E5A938419C71D6802F554FED0B61A88D3F7FF00A414DE171DAB17C27AACF77A75AAB5ACA636DC7CFCFCBDCD58FB7D
8B78A96CA4B52B7A232525F55C7F8503B1BB45666A1ACDB69D2A25D47300FF00718282AC7D01CF5A4BCD6EDAC5A117115C013005084C83EDD7AD0235
28AAB717D05B5AADC4858A3E36E07273D38A8EDF52826BB36AC258A6C6424AB82C3D45005EA28A2800A2AB5C5DC503AC6DB9A46E5635196349657D05
E23B40DCA36D753D54FBD005AA2AAC97D6F1C8D19932C83730552DB47BE3A5496D730DDC226B7904919E845004D4567C7ACE9D2CDE4C7701A4CE366D
6CFE58A9E0BCB7B83B62995981C15E87F23CD0059A2A0B7BA82E43791207D870C3A11F8546750B449FC879D565FEEB82BFCE802DD44D344922C6F2A2
BB7452C013F8546B7709BAFB2EE225DBB80652370F6CF5AC2F17048A4D3EE390DE7AA923B8CFFF005E901D3514C92448A3692460AAA3249AAF6F7F69
71218E19D59C7F09C83FAF5A605BA2A296E208789A68D09FEFB014867855914CB1867C1505865BE9EB40135145140051557508527B2915F236A96047
5040AC8F08BB5D69724B24AECE64642379E07AE3B1A00E868AC1D22F271AEDEE9D24AD2C510DE8CE72C390319ADB9244890BC8EAAA3AB31C01400FA2
A38A68A65DD148922E719560453B23246464738A00751451400514DC8C919191CE29D400514562EB97135ADD58BC32B2892508C99E18645006D51451
40051451400514566EBB3CF69A649736AC03C6413919C8E98FD68034A8A86DA5F3AD6294E01910371EE33535001451450014567EB3733D9D91B8B7DA
7CB605C30EA2AE41279D0472E31BD4363D3228024A28A2800A28A2800A28A2800A28A2800A28A2800A28A2800A28A2800A28A2800A28A2800A28A280
0A28A2800A28A2800A28A2800A28A2800A28A2800A28A2800A28A2800A28A2800A28A2800A28A2800A28A2800A28A2800A28A2800A28A2800A28A280
0A28A2800A28A2800A28A2800A28A2800A28A2800A28A2800A28A2800A28A2800A8E68C490BC64E37A95CFD6A4A280391F035BB69ADA8697330F3219
77000755F5FE55D7550934D89F5586FC31596352AC00FBE318E7F3ABF40D99FAE5EA5868D7774EC17644DB73DDB1C0FCEB2BC0166F67E14B6F3410F3
B34A41F73C7E80569EB3A5AEAD0C76F34ACB08903BA019DF8ED9ED57E3458D151061540007A0A5D43A0FA28A2988E4BC5D6FAC6AB05C6936B628D049
B4A4E5F1F5CFEB5A3E1F8F50B2B0B5B1BAB2545893699166047B715B94503B9C978E6C753D5ACE1D3AC6C96589DC3BC85F1823B7FF005EABF882CB59
D47C3F6DA52E9CAADF26E7120206DC63E95DAD14AC1739BD66D354BDF0A8D3ADED144F343E549BA5188F0319F7CD569B44D4A6F06C1A6058E2BA8194
03BF20807A8FCFA7B575B4530B9C7DDE9D71AA4B15ADCE90B1DC4657CCBC046DDA08E47A9FF0AEC28A28105727E2CD1EF2EB58D2756B38FCEFB0C9BA
48C752320F1EB5D6514018892DDC705DEA42C2532CDB425BFF001600C64FE678AC0F0AC7AAE85A55F44DA44ECF2CEF2C6AA7A6400073F4AEEA8A0773
8FF045A6A367697367A8D8CB12CCED26F24719EDFAD56F0EDADD7847ED965269F3DCC2CFBE3B8897208F43FA5773452B05CE563D0EE353D3F5737D9B
77D4C00B19FE0C0E091EB9A8D06AD6DE1CFEC45D399EE4C4605947FAA0A46324FE35D7514582E7302CA4F0F78316C60B692EE4F28A4822EB96072DF4
155BE1F477765A38D3EEACA68B616904AC30A724703DFAD76145160B9C55E69B77A478CA4D7ADED5EEADA68F6C810E594F73F418AD1B1375ABF88135
07B67B7B3B446584B8C348C78271E9D6BA4A2985CE1F4E4BE6F1E5CEAD2E95771C12C3E4A12983C6393CE074A75BADE9F1DCDAAB6977AB6C61110253
9240EBD6BB6A29582E7116A978FF00109F55974CBC4B6784408CC9D0F0327DABB7A28A62615C35C9B8FF00858B15F2E9D72D1C51344485EBDB70FC2B
B9A28038AD42DDD3C41A8CB7FA5BDDE953AA6F2C9921800323DAA4F0469E2D6EF519ED639E1D3E471E44721EBEF8FF003D6BB1A2958770ACFD77ED27
44BCFB0FFC7C794DB3EB5A145311E677B657979F0FE3B2B1D3E66B92C1EE895C1620E7AF7EB5BB05F5DC5E1B9A692C2EDD3C84B78620BF370A41623B
73DEBAFA295877383F0F5BDEFF00C2BEB9B03657315C056640CB82E5989F97E954EFB4EB5B9B482DAEB499E0D646D0B34280027D49FE95E9145160B9
C9F88346BC92FF0046D56DC79F3E9E409901399071C8F5E73F9D26B3672788F53D282C3347656F219647652A49F4FD2BADA28B05C28A28A62393F1DD
85D5DC7A6DC5AC6D27D92E448E17AE38AD5B6BA24DEEA2527309DAA89B0EEC28E481E9927F235AF4503B9C27C3F9648AF7548A4B4B98C5CCE6646742
0018EF5A1F102D2E6E74589ED63321827576500938CD757452B6817D6E701E36BC9353D274DB8B182594C570B2380A40C8EA3DEB43C69653EB3E1FB5
BBB0899A582459C46CA41C77E2BAE0000001803B52D160B9C76A13DBEA5A4C16DA75AFFA64E554E22DA63C7DEC9FCEBAD894A448858B9500163D4FBD
49450015C178EA756D77452914CCD6572B248446482B953DBE9FAD77B4C2A1B1B803839191D0D304731E2EB197C41E1A0FA7073323F991A30C1700E0
8C553D76ECEB3E144D3ED2094DD4FB11E32847978C6739FA76AEDA8A560B9C6789A24D37E1FF00F66CB9694C0A8A1549E4104D6C7842713F862C0608
68A148D811DC28AD9755752AEA194F504641A1115142A28551D00180280B8EA28A2988E034ABDB44F8917F3B4B88E784246C4100B647F854BA7EA56E
9F12AF903A849E258D4E7219B8E9F88AED3ECD6FE6F9BE445E6673BB60CE7EB4D3676AD3ACC6DE3F354E436D19CFAD2B0EE727344DA4FC4317F7036D
95EC251653F751F8E09ED9344F09D4FE20DADC5AB6F86D63CBB85F973E99AEBE7821B988C571124887AAB0C8A58A18A11886248C1ECAA05160B9C4F8
86F6DA3F885A2B3CD1ED8D191F9FBA4938CD26BB7B6FFF000B134491678D9110A390C30A4EEC73F88AEB67D2F4FB89D679AD217954E4315E73FD6897
4AD3E5911E4B280B21DC0EC039FEB4582E88359B0BDBC484D85E9B578D8938E8D9FF0027F3AC6F15E997D378224B6691AEAEA23E63329E58024FF235
D7514C2E7270EA9A65E7856D59CC72CC912A22372CB20007FF005EB3FC611DE44743D4EEBE68EDE48CDC10B8DA7B922BB38EC6CE26DF1DAC28D9DD90
8073534B147344D1CA8AE8C30558641A560B9C6F8E5BFB56D74EB6D39D6679270D94392A3D6A7F888DE57842488CA0C995032705F1D4D74B6D636B6B
FF001EF6F1C7EE0734DBFD3ACF518D63BDB749950E5430E868B01C9F88EEE39BE1CA35BC88EDE44590A4123819A7EBEF1DCFC3788C18C79316154F70
0715BE3C3DA42A141631043C15E71FCE9DFD83A57D8DAD459C6B1372401FD680D0AFA65E5BDCF852299645118B6D8C4B7DD2171835CEFC3E9E11E1CD
42369515CCF2705864F1E95D6D9E916365A77F67C3029B6E728DCE6A28BC3BA34258C5A7C2A58E4900F340181F0EA68869D7ABBD15BED2C48C81C53F
E1D4DBF4EBD8D9C6E5BA63B33D2B653C33A34723C8963182FD7191DB1DAA5D3F42D334CB99AE6C6D961926187209A02E73DE1F263F883AEC6C483220
6C1EFCF1FA53FC35711BF8C75D456032E300F7F71ED5B77BA069D7D7E97B3C27ED098F99588DC07407F2A6AF8774D5D5CEA6B1309F20E01C283EB8A0
2E646AEA3FE16468ACDC28B6939F53CE07E78ABDAC6A50A685772DFC2B2C2D3186351C0719E0E7F0FD2B4352D1ECB53689EEA33E642729229C32FB67
D29F75A5DA5D69AD612429E415DA171D3DFEB401C4F8B259E0D1F480F768DE64CACB0A0002AF1803E9C0AEB27BC6BED175031228645923C1390D85E7
F9D35BC33A5C9A62D8CF09915540DE4FCDC771E957AC74FB6B0B316B6E9FBAE7218E739EB9A2C1739AF085DFD8FC23A4A26374D3F9455BA8CB1CD3AE
6544F8936A1C819B7651F522B5B4EF0ED8E9F72668439009644639084F53F5A65E786AC2EAFD6F58CAB286C9C39C73D7E9405D12789EC0EA1A2CE918
CCD18F32320E0E4567F87AF5F5D36D72F1945B34D8C186732743FCABA8AA7A6E9D6DA640D0DA2155772ED939249A02E739E329A58F52D2915F623C98
0C46406FF38A9355B5D421D434DBF9A78E668A711EC45C17073FE15BDAA69B6FAA5A1B7B95F7561D50FA8AAD65A39B6785EE2F25B9F207EEC49D14F4
CD160B9AD45145311CCDA4C478E6EE09D816F24345F4E3FA55CB0D32DF4CD4A736CD23C93C6642AC781823BFBFF4A9B52D16DEFA78EE43BC1751F093
47C11F5F5AB163666D51B7CAF3C8DF7A473C9C74A4330FC1970F73FDA524A4993CFC367AF4A8BC31338F12EB76C9836C24DC8474CE7A0FCEB4A5D000
D45EEECAF25B4F347EF1231C31F5ABF61A7DBE9F118ED931B8E598F2CC7DCD0073D0C51FFC2C69C607CB6FBC0F72073524615BC7F2C7B570B087C7BE
3AD587F0DCADAA36A2BAA4E970CBB321074A69F0E5C7F6A26A2BAA482E400ACFE50F994763EB400C82729E3FB8B75E15ED7730F53F2D33C50AA75AD2
3CC2555A4DBBB1571340917594D51B5091A61C37EEC0DCBE9EC2A9F8C638E4934C5924D83CEC120E08CE391E9D293D83A8CF1CCED6B1D8DC45C48B21
0A475CF18E6A5F18366DB4F2C76E675247E55A434A13342D79726EA384EE8C328EBEA4F7A835AD0E5D5DC6FBF78A2520AA08C1DA7D739A6055F185CB
40B629B884925C119C0272319FCEA1D6DF5081ED2FE6B58634B5940CA37241238FA569DD6871DF68ABA7DF4EF3BA0F967230C0E720FF002FCA996FA5
5FBC7141A9DEACF0C44300A30CC474C9A00B7ACD9457FA64C921DA421657C72BC5627862E52F6E9E0BE60D7764A123465C1DA3F8BDEB6AF6C2E6E832
7DB99226607688C647B67D2B22E7C337736AFF00DA29A8849800372A6D27EB401D4D151A0658D43B6E60002D8C64FAD494C4472A9785D0632CA40CF4
AE47C2360648AEA55BDB841E71DC11B01BDEBAE955DE2658DF631180D8CE2B9FD2B47D5B4C12243756AD1BC9BF054E4D263447A75DBE9BE259B4BB82
25FB400EB363E627D0D6BDDA81A8DA3991CE588110190783F37B6335143A3AAEB2DAA4D2F993918501768518C7A9A351B1BC975086F6C6689248D0A6
D941239EFC7D680336D2FA3B7F176A31499862585491CE09E0038FCFF3A93C4E192F34E92172923CBB0B038E3AFF005A8AEBC3F7926ACD76B70B27DA
21F2E52FC6DE41E3F2E2AE6BBA7DF5F4F0340B09581F7A64E0E78EB4006B4C2DB55D32553B0C936C7607191C75ABDAB826C8E2E3C9C3024FF787F77F
1ACFD6EC351D420B668920596160F8DE4F3DFB549A85AEA375616AEA2317313877889F95FF00A7FF00AE8029DFDFECF10E9B2449244B302AFB971BD7
FAD686B77CD6F15B242F83712280EA7B6474FCEB3357D2F55BD9ED2FB6C46785B1E403C2839E73572EB4496EF46B6B669F6DD5B90EAFDB39CE0D0059
BE4BB4D46DAE209152DD7E59F73638CFBFD6A978AC6E6D3D7382D3E3F5156E24D4AE9561BE8A18D14832321CF998EC3D39AADE20B2D42F26B736B045
22C2E1D58B6083DF39FA50C08F5441A5EA3A74969B879AFB2450DF7C71D6BA4AE6759B2D57538EDF1691C7240DBC1F30633F9D74511731219000E546
E03A03DE840C92AADF9996D246B6744917E6CB0C8C0EBFA55AACFD66DA7B9B30B6C033A387DA4E0363B53114350D4CA6B1A7C7099556570A49CED707
DBA7E3567C5048F0EDE11D768FFD08566EB163A95E8B1D463B702780E5EDF7751EBFFD6A7EA30EA727875EC4DA349232805964527AE7007B70290CBE
6E1ACFC331CF18CBA4084679EA07F8D55D5AF268749B5BF8247F324D8AC15BE521873C74CFBD35F509E0F0F442EECB6CB9109864EB20C75158D16AF2
59C2BBB47965853F82504841DF6E46280375A69F43D11E59DC4A776630C0FCB9E707F5AA973AAA416715EDBDEC92BE41922639520F5E2AEDF247E23F
0EEEB4720C837A0C81861D8D416D7720D2469D3D8DC0B831940BE5FCA7DFFCFA5005DD6A55B9F0E4F34472AF1AB29FC41ACF93ED767A1C5AA7DA2432
44819A36E54A9C718AB93D84B1785FEC28A64916355C0EE7233552E6EE47D15B4F7B69D6E4A089404E0F4C1A00D58AF7EDBA30BDB56D9BA3DEB919C6
3A8FD08AC98F50D62EF4B5BAB50ACC9294655505986473FF00EAABF6164DA7681F667DCEE236C851B8E4F600552F08DD29B592C9D2459A37673B9703
191401B12DDB47A9456A623E5C884F987A67D3FCFAD56D3AF259351BAB29E4495A1C32B2AE38E98355B53BA68B57852F23636606E4DA3219BDEB22D3
50365E28BB95AD8859A30100CF3F37D28B858EAE5BD822BA104AE11CAEE05B81DFBFE14CD3EE9EE4CCB2AA0789F6E50E411D8FF3ACBD57ECF3EA8F6F
A8C6CB06C5F2E55F5E7AFE6699E19B696CEF6FA04B9FB4DA020A393920FA7E9401D25154ED35186EAE25811645922FBC1D7157298828A28A0028A28A
0028A28A0028A28A0028A28A0028A28A0028A28A0028A28A0028A28A0028A28A0028A28A0028A28A0028A28A0028A28A0028A28A0028A28A0028A28A
0028A28A0028A28A0028A28A0028A28A0028A28A0028A28A0028A28A0028A29ACC154B310001924F6A007535982A96620003249ED5CE2F899E78AE2F
ADA043A7DB36C7959B973FECFE628D675A927F0BCB7DA643E6A491B872CDB4C7C52B8ECCD3D335AB3D4CB0B666183C6F18DDF4AD2AE2FE1D19A4D0ED
CB5A2794A1B170D8DCC738E3DBA8AED2840C28A28A62336F75BD3EC67305C5C0130C7EEC0249CD59B2BB8EF6DD67843846240DCB83C5723730A4BF14
E1126D75FB2EEDBD70403D6BB1BAB88AD2DA4B89DB6C71A9663EC2921B23B9BEB6B59628A6902C9330544EEC49C549713C56D034D336D8D7A9C138ED
DABCEAEAF6EEEFC79A34D74A638E6DAF0C5CE55771009CFAE335E90E892215915594F50C32284EE0D58CFB2D6F4FBD9921B79F32B82421520E07E95A
75C732AC5F12EDE38D022FD949C28C0E86BB1A1030AA7A86A16FA74025B8620330550A32589F4AB95C8C929D5FC7A2D09DD6DA6C5BD97B173FE45304
75D4566EB3AB41A45A89A605DDDB6C7183F339F4158CFE29B98750B5D3EE6C563BBB89142A16FE13FC8D0163ABAC4BFF00125BE9B34AB796D731A21C
2C9B461FE9922ACEB3AC5BE9102493067795B6471AF5635C4F8AAF6EB57D4B47D12EAD0C17124C259141C8287A7D3A1CD2608F4382559A08E64CED91
430CF5C119ACDD4B5EB5D2DF6DE4770AA4E15F60DAC719C039A9353D4EDB47B6895D4967FDDC11275620703F95727E28BD7D575ED23449ADDA09BED0
24906E0E3663AF1F8D0091DE46EB246AE872AC0107DA9F59DAAEAB0E9A91865696695C2450A7DE624E3F2A6E95ABA6A135CDBB44629ED982BA16CF5A
622D5DDD7D920F35A19A419C1112EE23DFE95574BD72CB559658AD4C81E2FBCAEB8C5537F124463BAB8B6B769ACED01F367DE14023D0771EF543E1DA
B4BA45D6A2F9DD7B72EFC8EC3A521DB43AEA28AA7A95F45A6D8C975372A83851D589E805310CFED5B3FED51A679BFE945776DFC33FCB9ABF5CAF872F
ECB59D5EEAE26B148751B6C02C09395231F4AD11E21B417BE4382B1966549B39562BD681D8D9A2B22C35E82F7507B3F2A48A40BB937FF18FE95AF408
28A2B16FFC436F6578D0344EE8842C9283C213DA8036A8AC493C431473C05A071693BEC4B827009F5C7A54BAA6BB6BA6CD14528662EC0311D107B9F5
F6A00D6A2B2E0D6ECEE3563A7445DA4D9BC301F2D457BE20B6B5BB68238E4B8317FAE317223FAFBD006CD150417505C5AADCC322B42C3707E831EB59
C9E23B1372909F395646DAB33262327EB401B158D1F88B4E93526B10E77A92BB8FDDC8CF7FC3AD33C57AAAE93A24F3146264528AC3A29231C9ACBF08
E936D3E9F6B7979A784B88D485959B3E706EA48FC690EC75B1C892A078DD594F46539069F5CCF872592D757D4B4797EEC4FE745FEE13FF00EAAE9A98
828A28A0028A28A0028A28A0028A28A0028A28A0028A28A0028A28A0028A28A0028A28A0028A28A0028A28A0028A28A0028A28A0028A28A0028A28A0
028A28A0028A28A0028A28A0028A28A0028A28A0028A28A0028A28A0028A28A0028A28A0028A28A0028A28A0028A28A0028A28A0028A28A0028A28A0
028A28A0028A28A0028A28A002B3EFB48B3BE95659E33BD7B8EFF5AD0A280228208ADE158614088BD1476A968A2800A28A2800A28A2800A28A2800A2
8A2800A28A2800A28A2800A28A2800A28A2800A28A2800A28A2800A28A2800A28A2803335AD35AFE08CC2FB2685B7A7A1F6AAF3AEB92C2B06DB61E62
E1E453D3B1C83FD2B6E8A00A5A558A69D611DB260E39623B935768A2800A28A2800A28A2800A69504824024743E94EA2801A4060430041E083DE9B1C
31459F2A344CF5DAA066A4A280218E154B8966EAD26074E800A9A8A2800A28A2800A28A2800A28A2800A28A2800A28A2800A28A2800A28A2800A28A2
800A28A2800A28A2800A28A2800A28A2800A28A2800A28A2800A28A2800A28A2800A28A2800A28A2800A28A2800A28A2800A28A2800A28A2800A28A2
800A28A2800AABA94724DA6DCC50B1591E260A40CF38AB5450070DE0C86CC78425B1D55E2C2CCEF22BB60819E0FE62AC35E423C1FA8875486D033C56
EBEA319FD79AE81F46D3A46766B4425FEF7500FE14EBBD22CAF2DE28268479717DD51D052B1574627C3A1B3C276F1B6038662573C8C9ADDB6D4ED6EA
F6E2CE172D35BFDF18A75869F6BA741E4DA42234EF8EA69D058DADB5C4D3C10AA4B39CC8C3AB1A04CB345145311C79551F13C360E5AD0FD0F15ADE28
D36F354D2C4161308E5590373D0E2A94B6537FC2C182EC21F2BECA41603EA2BA281E49220D3446273D50B06C7E2290D9E55AB4978BE2BD2E37BF8E69
E354559D7EEC44B10735E9FA789C5A837372970E493BD00031E9C5667FC223A279A65366A5CBEFC93C8F6CFA569CCD1E9F60CD1424A42BC4683FCFE7
425606EE7333301F14605EFF006427E9C37F857635C4E9D23EA5E3C5BF48648E34B521830E87A0CFB1EB5DB50818571BE1C430F8F3C408E4967DAEA7
DB3D3F515D95735AB5AC9A7EBD6DAC5A44D2997F73346A3AE7A11430461EB111BDF89F670CEC05BDB46B2FCE7E5E39FC3935467D5E2BCF1DDCEA0B8F
2EC61DB6FBFA3B74CFE39AEF350D0F4FD4A5596EA006418F9C70481D8D44DE1BD25B5086F7ECA04B0AED500FCBF97E345877397BC2F7DF11F4DFB6E1
2086CD66504E173824FEBC5525D62D6F3C7771ABE77476B07976B19E4C9211B4607E27F3AEEF51D134FD5278A6BC877BC43039EA3D0D43FF0008DE90
3508AF45AA89621F2FA7D6815CE6D649EE7E24C2351548C5A5A07C6FF9724751F8E6A9D96A36F7BE39D4F5A7C343630F956E07FCB46CE01FE7F9D767
A9681A76A9731DC5DC24CA831B94E323D0FAD25BE81616FA91BE8A2C3EC5454FE15C77C7AD01738ECEA1AAF8FD92EA5FB0B5AC21860F452BCE33D2B4
2E02D87877C42DA62B607DDB83CB4A7A373DF19FD6BA6BBD1ECAF2769E5888998056914E091E9F4AB09676C967F64589441B4AECEA3068B05CE0B53B
D8B48F86705A4254DCDD441594727E6FBCC7D2BA8F0AB2C1A5DBE9F12964B68177499E0B1E700528F0AE92209A1F21889576125B240CE78AD6B4B68A
CED62B684111C4A154139A009EB2BC410C0DA7B5D5C79ACB680CC110E03301C67D6B56A39A28EE21786640D1B82ACA7B8A623CF3C07298AD10801AE3
549599CAE7746832339ED526A96A97BE37D2F418F77D92CA3F32419E4FF1649F7AEA743F0ED8E8924CF6C09691B82DFC0BD9455B8B4BB58B569F5354
CDD4EAA8CC7B0031C7E54AC55CE174DD40B789756D5B6063037D9ADE1DD8F98F1FA0FEB5E8EB9DA37000E3900E6B2ADFC3DA6DB6AB26A31C47CE91B7
104FCA1BD40F5AD7A1098579DF8DACD74ED36D6C16E19E4D42F774B29386FF003CD7A1332A296621540C924E0015E7BAE4F69AC7C40B0B6B9646B282
2386DDF2B31F7A1822CF8DCFEFF46D174F45337981B601D147FF00A8D27C4251158C5676E375E6A93468ED9C13B70071D8574F65A158DA5F9BE4F364
B8DBB43C8FBB68F4145FE89697DAA5A6A13E4C96BF717B13DA8B014EE6C61D134396F6DE056BD82D7679B8F98E056168D3A695F0DEE751948F36E95D
B39C9627200FAF5AEEA48D648DA3719560548F506B022F07E988163632C96C8FE6240CDF2A9A00E62E8EA3A67C39B3B70A44D7B3018EEAADC81FCAB4
3C6482CF4CD2347B2C34C645545F503BE7EB5D5EABA5DB6AB60D6570088CE0A95E0A11D08AA2DA6596965B56BE796E66B75E247E4A8E070280B987F1
21A66D074ED3F0AF3DCDC22100FDE2076FC4D7696F1086DA2857A4681463D862B8FBF921F11F8AF4716CEB25ADA0FB43FA86EA3F90AEC9DD63467738
550493E82840CE523264F8992AA0388ED14B91D3A77FCC5754258CCA621221900C94CF207D2B9CF0AC0D3DFEA5ACCA39BA93647FEE0FF3FA56E7D86D
FF00B47EDDB4F9FB76E73C7F9C508196E8A28A620A28A2800A28A2800A28A2800A28A2800A28A2800A28A2800A28A2800A28A2800A28A2800A28A280
0A28A2800A28A2800A28A2800A28A2800A28A2800A28A2800A28A2800A28A2800A28A2800A28A2800A28A2800A28A2800A28A2800A28A2800A28A280
0A28A2800A28A2800A28A2800A28A2800A28A2800A28A2800A28A2800A28A2800A28A2800A28A2800A28A2800A28A2800A28A2800A28A2800A28A280
0A28A2800A28A2800A28A2800A28A2800A28A2800A28A2800A28A2800A28A2800A28A2800A28A2800A28A2800A28A2800A28A2800A28A2800A28A280
0A28A2800A28A2800A28A2800A28A2800A28A2800A28A2800A28A2800A28A2800A28A2800A28A2800A28A2800A28A2800A28A2800A28A2800A28A280
0A28A2800A28A2800A28A2800A28A2800A28A2800A28A2800A28A2800A28A2800A28A2800A28A2800A28A2800A28A2800A28A280130339C0C9E334B4
51400514514005145140051451400514514005145140051451400514514005145140051451400C74592364719560411EA0D6547E1CD2A2D40DEADA83
2150A013951CE7A56C51400514514005145140053258D268DA39543230C107BD3E8A00A561A658E9FBFEC56C90EFFBD8CF3F9D5DA28A006222448123
45551D028C014FA28A0028A28A0028A28A0028A28A0028A28A0028A28A0028A28A0028A28A0028A28A0028A28A0028A28A0028A28A0028A28A0028A2
8A0028A28A0028A28A0028A28A0028A28A0028A28A0028A28A0028A28A0028A28A0028A28A0028A28A0028A28A0028A28A0028A28A0028A28A0028A2
8A0028A28A0028A28A0028A28A0028A28A0028A28A0028A28A0028A28A0028A28A0028A28A0028A28A0028A28A0028A28A0028A28A0028A28A0028A2
8A0028A28A0028A28A0028A28A0028A28A0028A28A0028A28A0028A28A0028A28A0028A28A0028A28A0028A28A0028A28A0028A28A0028A28A0028A2
8A0028A28A0028A28A0028A28A0028A28A0028A28A0028A28A0028A28A0028A28A0028A28A0028A28A0028A28A0028A28A0028A28A0028A28A0028A2
8A0028A28A0028A28A0028A28A0028A28A0028A28A0028A28A0028A28A0028A28A0028A28A0028A28A0028A28A0028A28A0028A28A0028A28A0028A2
8A0028A28A0028A28A0028A28A0028A28A0028A28A0028A28A0028A28A0028A28A0028A28A0028A28A0028A28A0028A28A0028A28A0028A28A0028A2
8A0028A28A0028A28A0028A28A0028A28A0028A28A0028A28A0028A28A0028A28A0028A28A0028A28A0028A28A0028A28A0028A28A0028A28A0028A2
8A0028A28A0028A28A0028A28A0028A28A0028A28A0028A28A0028A28A0028A28A0028A28A0028A28A0028A28A0028A28A0028A28A0028A28A0028A2
8A0028A28A0028A28A0028A28A0028A28A0028A28A0028A28A0028A28A0028A28A0028A28A0028A28A0028A28A0028A28A0028A28A0028A28A0028A2
8A0028A28A0028A28A0028A28A0028A28A0028A28A0028A28A0028A28A0028A28A0028A28A0028A28A0028A28A0028A28A0028A28A0028A28A0028A2
8A0028A28A0028A28A0028A28A0028A28A0028A28A0028A28A0028A28A0028A28A0028A28A0028A28A0028A28A0028A28A0028A28A0028A28A0028A2
8A0028A28A0028A28A0028A28A0028A28A0028A28A0028A28A0028A28A0028A28A0028A28A0028A28A0028A28A0028A28A0028A28A0028A28A0028A2
8A0028A28A0028A28A0028A28A0028A28A0028A28A0028A28A0028A28A0028A28A0028A28A0028A28A0028A28A0028A28A0028A28A0028A28A0028A2
8A0028A28A0028A28A0028A28A0028A28A0028A28A0028A28A0028A28A0028A28A0028A28A0028A28A0028A28A0028A28A0028A28A0028A28A0028A2
8A0028A28A0028A28A0028A28A0028A28A0028A28A0028A28A0028A28A0028A28A0028A28A0028A28A0028A28A0028A28A0028A28A0028A28A0028A2
8A0028A28A0028A28A0028A28A0028A28A0028A28A0028A28A0028A28A0028A28A0028A28A0028A28A0028A28A0028A28A0028A28A0028A28A0028A2
8A0028A28A0028A28A0028A28A0028A28A0028A28A0028A28A0028A28A0028A28A0028A28A0028A28A0028A28A0028A28A0028A28A0028A28A0028A2
8A0028A28A0028A28A0028A28A0028A28A0063C889B77B2AEE3B464E327D29F451400514514005145140051451400514514005145140051451400514
514005145140051451400514514005145140051451400514514005145140051451400514514005145140051451400514514005145140051451400514
514005145140051451400514514005145140051451400514514005145140051451400514514005145140051451400514514005145140051451400514
514005145140051451400514514005145140051451400514514005145140051451400514514005145140051451400514514005145140051451400514
514005145140051451400514514005145140051451400514514005145140051451400514514005145140051451400514514005145140051451400514
514005145140051451400514514005145140051451400514514005145140051451400514514005145140051451400514514005145140051451400514
514005145140051451400514514005145140051451400514514005145140051451400514514005145140051451400514514005145140051451400514
514005145140051451400514514005145140051451400514514005145140051451400514514005145140051451400514514005145140051451400514
514005145140051451400514514005145140051451400514514005145140051451400514514005145140051451400514514005145140051451400514
514005145140051451400514514005145140051451400514514005145140051451400514514005145140051451400514514005145140051451400514
514005145140051451400514514005145140051451400514514005145140051451400514514005145140051451400514514005145140051451400514
514005145140051451400514514005145140051451400514514005145140051451400514514005145140051451400514514005145140051451400514
514005145140051451400514514005145140051451400514514005145140051451400514514005145140051451400514514005145140051451400514
514005145140051451400514514005145140051451400514514005145140051451400514514005145140051451400514514005145140051451400514
514005145140051451400514514005145140051451400514514005145140051451400514514005145140051451400514514005145140051451400514
514005145140051451400514514005145140051451400514514005145140051451400514514005145140051451400514514005145140051451400514
514005145140051451400514514005145140051451400514514005145140051451400514514005145140051451400514514005145140051451400514
51400514514005145140051451400514531DB646CD8CED04E2800791131BDD573D32714DF3E1FF009EA9FF007D0AC179DE690991896E7BF4148AC47C
A79F6AC1D746AA91D18208C8390696B22C2E5D655889C83C609E95AD5AC64A488946CEC2D14515448514514005145140051451400514514005145140
051451400514514005145140051451400514514005145140051451400514514005145140051451400514514005145140051451400514514005145140
051451400514514005145140051451400514514005145140051451400514514005145140051451400514514005145140051451400514514005145140
051451400514514005145140051451400514514005145140051451400514514005145140051451400514514005145140051451400514514005145140
051451400514514005145140051451400514514005145140051451400514514005145140051451400514514005145140051451400514514005145140
051451400514514005145140051451400514514005145140051451400514514005145140051451400514514005145140051451400514514005145140
051451400514514005145140051451400514514005145140051451400514514005145140051451400514514005145140051451400514514005145140
051451400547329785D1480594819A928A00E602ED936B8C3838618E69E586D04E7DEB4AEF4EDF219A22431392BEB5047A6CC5B2D85C7192DFE15C8E
94AE7529C6C374E8CCB721D54EC5EA6B6EA0B5B75B68B60393DCE3AD4F5BD3872A39E72E661451456848514514005145140051451400514514005145
140051451400514514005145140051451400514514005145140051451400514514005145140051451400514514005145140051451400514514005145
140051451400514514005145140051451400514514005145140051451400514514005145140051451400514514005145140051451400514514005145
140051451400514514005145140051451400514514005145140051451400514514005145140051451400514514005145140051451400514514005145
140051451400514514005145140051451400514514005145140051451400514514005145140051451400514514005145140051451400514514005145
140051451400514514005145140051451400514514005145140051451400514514005145140051451400514514005145140051451400514514005145
140051451400514514005145140051451400514514005145140051451400514514005145140051451400514514005145140051451400514514005145
140051451400514514005145140051451400514514005145140051451400514514005145140051451400514514005145140051451400514514005145
140051451400514514005145140051451400514514005145140051451400514514005145140051451400514514005145140051451400514514005145
140051451400514514005145140051451400514514005145140051451400514514005145140051451400514514005145140051451400514514005145
140051451400514514005145140051451400514514005145140051451400514514005145140051451400514514005145140051451400514514005145
140051451400514514005145140051451400514514005145140051451400514514005145140051451400514514005145140051451400514514005145
140051451400514514005145140051451400514514005145140051451400514514005145140051451400514514005145140051451400514514005145
140051451400514514005145140051451400514514005145140051451400514514005145140051451400514514005145140051451400514514005145
140051451400514514005145140051451400514514005145140051451400514514005145140051451400514514005145140051451400514514005145
140051451400514514005145140051451400514514005145140051451400514514005145140051451400514514005145140051451400514514005145
140051451400514514005145140051451400514514005145140051451400514514005145140051451400514514005145140051451400514514005145
140051451400514514005145140051451400514514005145140051451400514514005145140051451400514514005145140051451400514514005145
140051451400514514005145140051451400514514005145140051451400514514005145140051451400514514005145140051451400514514005145
140051451400514514005145140051451400514514005145140051451400514514005145140051451400514514005145140051451400514514005145
140051451400514514005145140051451400514514005145140051451400514514005145140051451400514514005145140051451400514514005145
140051451400514514005145140051451400514514005145140051451400514514005145140051451400514514005145140051451400514514005145
140051451400514514005145140051451400514514005145140051451400514514005145140051451400514514005145140051451400514514005145
140051451400514514005145140051451400514514005145140051451400514514005145140051451400514514005145140051451400514514005145
140051451400514514005145140051451400514514005145140051451400514514005145140051451400514514005145140051451400514514005145
140051451400514514005145140051451400514514005145140051451400514514005145140051451400514514005145140051451400514514005145
140051451400514514005145140051451400514514005145140051451400514514005145140051451400514514005145140051451400514514005145
140051451400514514005145140051451400514514005145140051451400514514005145140051451400514514005145140051451400514514005145
140051451400514514005145140051451400514514005145140051451400514514005145140051451400514514005145140051451400514514005145
140051451400514514005145140051451400514514005145140051451400514514005145140051451400514514005145140051451400514514005145
140051451400514514005145140051451400514514005145140051451400514514005145140051451400514514005145140051451400514514005145
140051451400514514005145140051451400514514005145140051451400514514005145140051451400514514005145140051451400514514005145
14005145140051451400514514005145140051451400514514005145140051451400514514005145140051451400DDC3715C8C819229D45140051451
400514514005145140051451400514514005145140051451400514514005145140051451400514514005145140051451400514514005145140051451
400514514005145140051451400514514005145140051451400514514005145140051451400514514005145140051451400514514005145140051451
40051451400514514005145140051451400521200249C01D49A28A008C5C404E04D19EFF007853D59594329041EE0D145240C752120724814514C05A
28A2800A28A2800A28A2800A28A2800A28A2800A28A2800A28A2800A28A2800A28A2800A28A2800A28A2800A28A2800A28A2800A28A2800A28A2800A
28A2800A28A2800A28A2800A28A2800A28A2800A28A2800A28A2800A28A2803FFFD9>
}
\xdef\GKZDedicationImageObject{\the\pdflastobj}
\endgroup
\newsavebox{\GKZDedicationBox}
\setbox\GKZDedicationBox=\hbox to 158.6bp{%
  \pdfliteral{q 158.6 0 0 226.4 0 0 cm /GKZDedication Do Q}\hss}
\ht\GKZDedicationBox=226.4bp
\dp\GKZDedicationBox=0pt
\immediate\pdfxform resources{
  /XObject << /GKZDedication \GKZDedicationImageObject\space 0 R >>
}\GKZDedicationBox
\edef\GKZDedicationFormObject{\the\pdflastxform}
\newcommand{\GKZDedicationScan}{\pdfrefxform\GKZDedicationFormObject}

\begin{document}
\maketitle

\begin{abstract}
After some personal recollections of Jean-Pierre Kahane, this paper offers a gentle introduction to the \GKZ\ theorem. 
The theorem asserts that a unital linear functional on a complex Banach algebra is multiplicative whenever it does not vanish on invertible elements. 
Neither continuity of the functional
nor commutativity of the algebra is assumed. We give a detailed account
of three proofs: the original argument of Gleason and, independently,
Kahane and \.{Z}elazko; the polynomial proof of Roitman and Sternfeld;
and a proof by the author. We explain the complex analysis behind each:
a special case of Hadamard's factorisation theorem, a polynomial argument
using the fundamental theorem of algebra, and the P\'olya--Saxer theorem on
zero-free entire functions and their derivatives. The elementary
entire-function lemma needed for the first proof is proved in full;
the P\'olya--Saxer theorem is used without proof.
We give full details of the algebraic argument common to the three proofs,
which deduces multiplicativity from the preservation of squares, including
the non-commutative case. We then discuss applications and counterexamples,
and present an extension to Ransford's generalised spectra using the
Roitman--Sternfeld argument. We conclude with a brief guide to the more
recent literature.
\end{abstract}

\setcounter{tocdepth}{2}

\section{Introduction}
\label{sec:personal-introduction}

\subsection{Conversations with Kahane}
\label{sec:personal-preamble}

I am very proud to have been Jean-Pierre Kahane's PhD student---his
last \emph{th\'esard}, as we say in colloquial French.

I first met him at one of the Monday-afternoon harmonic analysis
seminars in Orsay. I was then studying for the \emph{Dipl\^ome
d'\'etudes approfondies de Math\'ematiques Pures}, a postgraduate
research qualification broadly equivalent to today's French
Master~2. The programme was jointly run by Universit\'e Paris-XI
(now Paris-Saclay) based in Orsay, the \'Ecole normale sup\'erieure de Cachan and the \'Ecole
polytechnique. When I told Kahane that I hoped to go on to a PhD in
analysis, he made it clear that he would not be taking on any more
doctoral students. He nevertheless found time for occasional
conversations with me, despite his busy schedule.

We spoke first about Paul Erd\H{o}s.
Although Kahane had never written a joint paper with
Erd\H{o}s,\footnote{Kahane's
Erd\H{o}s number is two, through either Jean-Louis Nicolas or Jacques
Dixmier. Both paths are documented by the joint papers
\cite{Dixmier-Erdos-Nicolas,Kahane-Dixmier-Nicolas}, published in 1987
and 2007 respectively. For further information, see the
Erd\H{o}s Number Project~\cite{ErdosNumberProject}; MathSciNet also
provides a collaboration-distance tool~\cite{MathSciNetDistance}.} he
regarded him as a friend and shared with me some stories about their
encounters. Readers unfamiliar with Erd\H{o}s will enjoy the delightful book \cite{hoffman}. Kahane also wrote about Erd\H{o}s~\cite{Kahane-Erdos}.
He also explained to me several results that Erd\H{o}s had proved using the probabilistic method, and he later introduced me to some of his own applications of probability to analysis.

My own connection with Erd\H{o}s went back to one of his open
problems, which I had solved when I was a teenager. He had asked about the
arithmetic nature of the sum
\[
    \theta=\sum_{n=1}^{\infty}\frac{1}{F_{2^n+1}},
\]
where $(F_n)$ is the Fibonacci sequence; the question appears, with
related problems, in his book with Ronald Graham~\cite{erdos-graham}.
I found an elementary proof that $\theta$ is irrational. The result was
subsequently published in~\cite{Glasgow}.

Erd\H{o}s himself sent me a very kind letter, even suggesting where
I might submit the paper. His advice probably helped it find its
way into print: my exposition was very poor. He also suggested
further problems, which I could not solve at the time, and invited me to
visit him in Budapest,\footnote{It was difficult for a Romanian citizen
to travel to Western countries before 1989, but travel to Hungary was
permitted. Paul Erd\H{o}s died in September 1996.} an invitation that,
unfortunately, I never took up. He encouraged
me to pursue an academic career, which, with some good fortune, I
eventually did.

Elementary though my proof was, I think it helped Kahane to see
me as an \emph{interlocuteur math\'ematique}---someone with whom he
could discuss mathematics. To my surprise, he wanted to hear the
details of my irrationality argument. By then, however, I had
forgotten my original proof. Trying to reconstruct it at the blackboard
in his office gave me the idea that led to my later paper~\cite{ActaArith}.
I return to these irrationality questions in Section~\ref{sec:fibonacci-memories}.

It was during our early mathematical conversations that Kahane
introduced me to the Gleason--Kahane--\.{Z}elazko theorem. I managed
to prove a few modest results, some of which were announced,
without proofs, in a short note in the \emph{Comptes rendus de
l'Acad\'emie des sciences}~\cite{Badea1991}, presented to the Academy
by Jean-Pierre Kahane. A later account of the
third proof presented here appeared in~\cite{Badea}; further papers
related to GKZ include~\cite{BadeaPalermo1993,BadeaRasa1993}. I believe
that this early work on GKZ also helped persuade him to reconsider
his decision about taking on another doctoral student.

We also spoke about the history of mathematics. I had read several
books and articles about the Lw\'ow school and the Scottish Caf\'e \cite{Mauldin},
and knew a few stories about Stefan Banach that were new to him. Kahane was preparing
a talk on the Polish mathematical school before the Second World
War, and some of the material we discussed found its way into
his lecture at the Polish Academy of Sciences' Scientific Centre
in Paris on 9 April 1992~\cite{Kahane-Polish}.

Kahane kindly invited me to attend his talk. My recollection of the
lecture itself is rather hazy. What I do remember is sitting between two
relatively young Polish guests whose interest in the history of
mathematics seemed rather limited. They appeared to be waiting for the lecture
to end and the reception to begin. There was a generous supply
of Polish vodka and French wine.

\subsection{A souvenir from our first meetings}

Figure~\ref{fig:kahane-dedication} reproduces the first page of an
offprint of Kahane and \.{Z}elazko's paper in \emph{Studia Mathematica},
with a handwritten dedication from Kahane to me.

\begin{figure}[t]
\centering
\resizebox{0.78\linewidth}{!}{\GKZDedicationScan}
\caption{The first page of Kahane and \.{Z}elazko's paper
\cite{KZ}, with Jean-Pierre Kahane's handwritten dedication
to the author.}
\label{fig:kahane-dedication}
\end{figure}

The dedication takes me back to the period before I began my PhD.
Mme Josette Dumas, the secretary of the harmonic analysis group in Orsay,
had given me the offprint. Kahane added the dedication in his
office, writing very quickly, just before leaving for a meeting. I
remember the scene as vividly as if it were yesterday. Soon afterwards,
he agreed to supervise my PhD, just in time for me to renew my
scholarship and begin doctoral work as a \emph{Boursier du Gouvernement
Fran\c{c}ais}, a French government scholarship holder. I remain deeply
grateful to Jean-Pierre Kahane for that opportunity and for all that he
taught me, both during my PhD and in the years that followed.

My thoughts also turn to Wies\l{}aw \.{Z}elazko, who died in
November 2025. I fondly remember our conversations many years ago in
Orsay, Warsaw and B\k{e}dlewo---about Banach algebras, Solidarno\'s\'c
and much else---and the kindness he always showed me. One evening,
he invited me to dinner at his home in Warsaw and, at the end of
the meal, asked each of his guests to write a few words in his
visitors' book. As I leafed through its pages, I recognised several
famous names. Having only recently completed my PhD, I was surprised
to find myself adding my own alongside them.

\subsection{A digression on the irrationality of some Fibonacci series}
\label{sec:fibonacci-memories}

Erd\H{o}s's problem mentioned above belongs to a family of questions
about the irrationality or transcendence of the sums of certain infinite
series of rational numbers. To recall the background, write
\[
    F_0=0,\qquad F_1=1,\qquad
    F_{n+2}=F_{n+1}+F_n\quad(n\geq0).
\]
Good~\cite{Good}, Hoggatt and Bicknell~\cite{HoggattBicknell}, and
probably many others, rediscovered Lucas's identity~\cite{Lucas}
\[
    \sum_{n=0}^{\infty}\frac{1}{F_{2^n}}
    =\frac{7-\sqrt{5}}{2}.
\]
This sum is therefore algebraic and irrational. By contrast,
the sum
\[
    \sum_{n=0}^{\infty}\frac{1}{n!\,F_{2^n}}
\]
is transcendental. This was proved by Mignotte~\cite{Mignotte}
and subsequently obtained by Mahler using a different
method~\cite{Mahler}.

The series $\theta=\sum_{n\ge1}1/F_{2^n+1}$ of
Section~\ref{sec:personal-preamble} has a general term of the same order
of growth as $1/F_{2^n}$, but no closed form comparable to Lucas's
identity is available for its sum. Erd\H{o}s asked
about the arithmetic nature of $\theta$, together with related problems,
in~\cite{erdos-graham}. I proved in 1987 that $\theta$ is
irrational~\cite{Glasgow}. In 1994, Becker and T\"opfer~\cite{BT} proved
that $\theta$ is transcendental, which gives a complete answer to
Erd\H{o}s's question.

One of the further questions in Erd\H{o}s's letter, also recorded
in~\cite{erdos-graham}, concerned the more general series
\begin{equation}\label{eq:fibonacci-lacunary}
    \sum_{k=1}^{\infty}\frac{1}{F_{n_k}},
\end{equation}
where $(n_k)$ is a strictly increasing sequence of positive
integers. Must its sum be irrational whenever
$n_{k+1}/n_k\geq c$ for some fixed $c>1$? In my 1993 paper in
\emph{Acta Arithmetica}~\cite{ActaArith}, I answered this
affirmatively for $c\geq2$.
In 2022, Nguyen~\cite{nguyen} proved that the sum
in~\eqref{eq:fibonacci-lacunary} is transcendental under the
assumption $n_{k+1}/n_k\geq c>2$. The strict inequality is necessary:
for $n_k=2^{k-1}$, the successive ratios are exactly $2$, yet the sum is
algebraic by Lucas's identity displayed above. Thus
$c>2$ cannot be replaced by $c\geq2$ in a general transcendence
statement.

The general irrationality question is listed as Problem~267 on
the Erd\H{o}s Problems website~\cite{Erdos267}. More recently,
the Star Fleet project has made available an AI-assisted proof
covering the full range $c>1$, together with a Lean
formalisation~\cite{StarFleet2026}. As of 25 October 2026, the
Formal Conjectures repository records this case as solved and
links to the formal proof~\cite{FormalConjectures267}.

I also discussed irrationality questions with Richard
Andr\'e-Jeannin when he visited me in Orsay. His 1989 proof of the
irrationality of
\[
    \sum_{k=1}^{\infty}\frac{1}{F_k},
\]
using Ap\'ery's method, is a particularly striking
result~\cite{AndreJeannin}. Whether this sum is transcendental appears
to be an open problem.

\subsection{What this paper is about}

The aim of these notes is to give a pedagogical account of the
Gleason--Kahane--\.{Z}elazko (GKZ) theorem. We present three proofs
and explain precisely which results from complex analysis each
requires. We also include comments and observations, together
with material from the work underlying the 1991 \emph{Comptes
rendus} note~\cite{Badea1991} that did not appear in the published
version.

The exposition is intended to be largely self-contained. After the
proofs, we give applications and examples showing the role of the
hypotheses. Section~\ref{sec:generalised} returns to one of the results
announced in the 1991 note: the extension of GKZ to Ransford's generalised
spectra. We present an elementary proof based on the Roitman--Sternfeld
argument. A short guide to further
developments, including recent results for reproducing kernel Hilbert
spaces, concludes the paper (Section~\ref{sec:further}).

\subsection{Acknowledgements}
In the current period of transition brought about by artificial intelligence, the need for surveys and gentle introductions may well increase. 
I would like to thank Herv\'e Queff\'elec for the opportunity to
write these notes, and for his help, comments and patience. 
I am also grateful to Persi Diaconis for a conversation about Gleason,
Kahane and the \GKZ\ theorem, which prompted me to look at possible applications
of GKZ to quantum mechanics. Although I usually write mathematics in English, I
first wrote this article in French as a tribute to Jean-Pierre Kahane,
before translating it into English. The French version is available from
the author on request. The author acknowledges the support of the CDP
C2EMPI and its institutional partners for the project R-CDP-24-004-C2EMPI,
and of the EU COST Actions.

\section{The theorem and the three routes to it}

\textbf{Conventions.} A \emph{unital Banach algebra} is a complex
algebra with a non-zero identity $\one$, complete for a
submultiplicative norm. We assume that $\norm{\one}=1$.
The main characters in these notes are, literally,
\emph{characters}: a character on a unital complex algebra $A$
is a non-zero complex-linear map $\varphi:A\to\C$ satisfying
\[
\varphi(ab)=\varphi(a)\varphi(b)\qquad(a,b\in A).
\]
A character necessarily satisfies $\varphi(\one)=1$, since
$\varphi(a)=\varphi(a)\varphi(\one)$ and $\varphi(a)\ne0$ for
some $a\in A$. Moreover, if $u$ is invertible, then
\[
\varphi(u)\varphi(u^{-1})=\varphi(\one)=1.
\]
Thus a character cannot vanish on an invertible element. The
remarkable fact is that, for a normalised (unital) linear functional on
a complex Banach algebra, this necessary condition is also
sufficient for multiplicativity.

\begin{theorem}[\GKZ]\label{thm:gkz}
Let $A$ be a unital complex Banach algebra and let $\varphi:A\to\C$
be complex-linear. Suppose that
\begin{equation}\label{eq:hyp}
\varphi(\one)=1,
\qquad
\varphi(u)\ne0\quad(u\in\inv).
\end{equation}
Then
\begin{equation}\label{eq:mult}
\varphi(ab)=\varphi(a)\varphi(b)\qquad(a,b\in A).
\end{equation}
The functional is automatically continuous. With the normalisation
$\norm{\one}=1$, it satisfies $\norm{\varphi}=1$.
\end{theorem}

The normalisation $\varphi(\one)=1$ of the functional is needed for this precise
conclusion. If a linear functional $\Lambda$ is merely assumed
to be non-zero on every invertible element, then
$\Lambda(\one)\ne0$, and the theorem applies to
\[
\varphi=\frac{\Lambda}{\Lambda(\one)}.
\]
The conclusion in this case is
\begin{equation}\label{eq:unnormalized}
\Lambda(ab)\Lambda(\one)=\Lambda(a)\Lambda(b).
\end{equation}

The theorem was proved independently by Gleason~\cite{Gleason}
and by Kahane and \.{Z}elazko~\cite{KZ}. Gleason's paper appeared
in 1967, and that of Kahane and \.{Z}elazko a year later.
Wies\l{}aw \.{Z}elazko once told me how their proof came about,
and Kahane later confirmed the story. During a seminar talk in Orsay,
\.{Z}elazko mentioned the result as a conjecture. Kahane immediately saw
how to apply Hadamard's factorisation theorem. Gleason independently
found the same approach. Sometimes an idea is in the air, and several
mathematicians arrive at it at much the same time.

In a recent conversation, Persi Diaconis remarked that the \GKZ\ theorem is not a result one would immediately associate with the names of Kahane or Gleason. I have explained above how Kahane became interested in the problem. In Gleason's case, I do not know for certain what drew him to maximal ideals (or characters). One conjecture, possibly unfounded, is that his interest came from quantum mechanics. The \GKZ\ theorem does at least have a connection with quantum mechanics: Bach~\cite{quantum} used it to examine the theoretical assumptions underlying Bell's theorem on local hidden variables.

Andrew M.~Gleason is also famous for a 1957 result, answering a question posed by George W.~Mackey, which played a historically significant role in showing that wide classes of hidden-variable theories are incompatible with quantum physics. Gleason's theorem~\cite{Gleason1957} shows that the Born rule, which is used to calculate probabilities in quantum physics, can be derived from the usual mathematical representation of measurements together with the assumption of non-contextuality.

\subsection{How the proofs fit together}

Linearity relates sums and scalar multiples; multiplicativity
relates products. The hypothesis in~\eqref{eq:hyp} says nothing
directly about a pair of elements $a,b$: it tests only whether
the functional vanishes on the group of invertible elements.
Nevertheless, this condition determines its behaviour on every
product.

Each of the three proofs first establishes the apparently weaker
identity
\begin{equation}\label{eq:squares-intro}
\varphi(a^2)=\varphi(a)^2\qquad(a\in A).
\end{equation}
For a commutative algebra, expanding $(a+b)^2$ then gives
multiplicativity. In a non-commutative algebra, this expansion
initially provides information only about the Jordan product
$\tfrac12(ab+ba)$, so a further argument is needed. We give that argument once, in full, in
Section~\ref{sec:algebra}.

The three routes to~\eqref{eq:squares-intro} use different
properties of complex functions.
\begin{enumerate}
\item \textbf{Gleason--Kahane--\.{Z}elazko: growth of a zero-free
entire function.} For fixed $a$, the function
$z\mapsto\varphi(e^{za})$ is entire, zero-free and of at most
exponential growth. Such a function is the exponential of an
affine function. Comparing its first two derivatives at the
origin gives~\eqref{eq:squares-intro}.
\item \textbf{Roitman--Sternfeld: roots of polynomials.}
For $\varphi(x)=0$, the roots of the monic polynomial
$\varphi((\lambda\one-x)^n)$ lie in $\spA(x)$, the spectrum of $x$.
Its first two coefficients show that $\varphi(x^2)\ne0$ would force some of
these roots to escape to infinity as $n\to\infty$. Boundedness
of $\spA(x)$ rules this out.
\item \textbf{Non-vanishing of derivatives~\cite{Badea}.}
If $x$ is invertible, every derivative of
$z\mapsto\varphi(e^{zx})$ is zero-free. The P\'olya--Saxer
theorem identifies this function as an exponential. A scalar
shift then transfers the square identity from invertible
elements to arbitrary elements.
\end{enumerate}
The remark following
Lemma~\ref{lem:defect} explains why the three different starting points
all lead to the square identity for every element of $A$.

The original proofs of Gleason and of Kahane and \.{Z}elazko
use Hadamard's factorisation theorem. Only a special case is
needed, and it can be proved by elementary complex analysis.
Siddiqi~\cite{Siddiqi} gave an early proof along these lines;
the first proof below follows the same principle. We establish
the required entire-function lemma directly, using a holomorphic
logarithm and a coefficient estimate for its real part.

The polarisation argument in Kahane and \.{Z}elazko's paper
assumes commutativity. A footnote added in proof records the
extension to non-commutative algebras, referring to
\.{Z}elazko~\cite{Zelazko}. Combining their analytic argument
with the algebraic lemma in Section~\ref{sec:algebra} gives
the theorem in its full non-commutative form. Gleason's
note~\cite{Gleason} likewise treats only the commutative case.

\section{Preliminaries: invertibility, spectra, and continuity}
\label{sec:prelim}

Unless an assertion is explicitly stated for more general
algebras, in what follows $A$ is a unital complex Banach algebra with a submultiplicative
norm satisfying $\norm{\one}=1$.

\subsection{The spectrum needs no topology in its definition}

For a unital complex algebra $A$ and $a\in A$, define
\[
\spA(a)=\{\lambda\in\C:a-\lambda\one\notin\inv\}.
\]
This definition is purely algebraic. A norm and completeness enter the
standard Banach-algebra proofs that spectra are bounded, closed and
non-empty. The usual proof of non-emptiness of the spectrum uses Liouville's theorem.
The GKZ proofs below do not require this general non-emptiness result.

\begin{lemma}[The spectral reformulation]\label{lem:selection}
For a complex-linear map $\varphi:A\to\C$ on a unital complex algebra,
the following conditions are equivalent:
\begin{enumerate}[label=\textup{(\roman*)}]
\item $\varphi(\one)=1$ and $\varphi(u)\ne0$ for every $u\in\inv$;
\item $\varphi(a)\in\spA(a)$ for every $a\in A$.
\end{enumerate}
\end{lemma}

\begin{proof}
Assume (i). For any $a\in A$,
\[
\varphi\bigl(a-\varphi(a)\one\bigr)=0.
\]
The element $a-\varphi(a)\one$ therefore cannot be invertible. By
definition this is exactly $\varphi(a)\in\spA(a)$.

Conversely, $\spA(\one)=\{1\}$, so (ii) gives
$\varphi(\one)=1$. If $u$ is invertible, then $0\notin\spA(u)$;
hence $\varphi(u)\ne0$.
\end{proof}

A functional satisfying these equivalent conditions may be called a
\emph{linear spectral selection}. The theorem says that, on a complex
Banach algebra, linear spectral selection, that is, a linear choice of a spectral value must be
multiplicative (a character). The spectral-selection point of
view leads to a different interpretation of the \GKZ\ theorem, as a kind
of local-to-global principle. When the Banach algebra $A$ is commutative,
classical Gelfand theory shows that the spectrum $\spA(a)$ of $a$ is
exactly the set of values $\gamma(a)$, where $\gamma$ ranges over the
characters of $A$. Consequently, if a normalised linear functional
$\varphi$ on a commutative Banach algebra $A$ has the local property that, for
every $a\in A$, there is a character $\gamma_a$, depending on $a$, with
$\varphi(a)=\gamma_a(a)$, then $\varphi$ is itself a (global) character.

A remarkable non-linear extension of the \GKZ\ theorem is the
Kowalski--S\l{}odkowski theorem~\cite{KS}. If $\varphi$ is a
complex-valued function on a unital Banach algebra $A$, not necessarily
linear, satisfying $\varphi(0)=0$ and the spectral condition
\[
\varphi(x)-\varphi(y)\in\spA(x-y)\qquad(x,y\in A),
\]
then $\varphi$ is automatically linear and multiplicative. Thus
linearity is now part of the conclusion. For a linear $\varphi$, the
spectral condition reduces to $\varphi(a)\in\spA(a)$ for every $a\in A$,
so the statement contains the classical \GKZ\ theorem.

The Kowalski--S\l{}odkowski theorem implies the following
`$2$-local-to-global' principle. If a complex-valued function $\varphi$
on a commutative Banach algebra $A$ has the property that, for every pair $a,b\in A$, there is a
character $\gamma_{a,b}$, depending on $a$ and $b$, with
$\varphi(a)=\gamma_{a,b}(a)$ and $\varphi(b)=\gamma_{a,b}(b)$, then
$\varphi$ is itself a character, that is, linear and multiplicative.
Indeed, under this hypothesis,
\[
\varphi(a)-\varphi(b)=\gamma_{a,b}(a-b)\in\spA(a-b)\qquad(a,b\in A),
\]
and $\varphi(0)=\gamma_{0,0}(0)=0$, so the Kowalski--S\l{}odkowski
theorem applies.

\subsection{The Neumann series and the spectral bound}

We first recall exactly where completeness is used to bound the spectrum.
The series in the following lemma is named after Carl Neumann.

\begin{lemma}[Neumann series]\label{lem:neumann}
If $b\in A$ and $\norm{b}<1$, then $\one-b$ is invertible, with
\[
(\one-b)^{-1}=\sum_{n=0}^{\infty}b^n,
\qquad
\norm{(\one-b)^{-1}}\le\frac{1}{1-\norm{b}}.
\]
Consequently,
\begin{equation}\label{eq:specbound}
\spA(a)\subset\{\lambda\in\C:|\lambda|\le\norm{a}\}.
\end{equation}
\end{lemma}

\begin{proof}
The series is absolutely convergent because
$\norm{b^n}\le\norm{b}^n$ for $n\ge0$, using $\norm{\one}=1$ when $n=0$.
Completeness ensures that its sum belongs
to $A$. For its partial sums,
\[
(\one-b)\sum_{n=0}^{N}b^n
=\sum_{n=0}^{N}b^n(\one-b)=\one-b^{N+1}.
\]
Since $\norm{b^{N+1}}\to0$, continuity of multiplication gives a
two-sided inverse after passage to the limit. The norm estimate follows
from the geometric series.

If $|\lambda|>\norm{a}$, apply this to $b=a/\lambda$ and use
$\lambda\one-a=\lambda(\one-a/\lambda)$. Thus $\lambda$ lies outside
the spectrum, which proves \eqref{eq:specbound}.
\end{proof}

\begin{proposition}[Automatic continuity before multiplicativity]
\label{prop:continuous}
Let $\varphi:A\to\C$ be complex-linear and satisfy \eqref{eq:hyp}. Then
\[
|\varphi(a)|\le\norm{a}\qquad(a\in A).
\]
In particular, $\varphi$ is continuous and, with our normalisation
$\norm{\one}=1$, satisfies $\norm{\varphi}=1$.
\end{proposition}

\begin{proof}
Lemma~\ref{lem:selection} places $\varphi(a)$ in the spectrum, and
Lemma~\ref{lem:neumann} bounds its modulus. This proves
$\norm{\varphi}\le1$. The equality $\varphi(\one)=1$ gives the reverse
inequality when $\norm{\one}=1$.
\end{proof}

This order of reasoning matters. We shall apply $\varphi$ to an infinite
power series. An arbitrary algebraic linear functional need not commute
with a norm limit. Proposition~\ref{prop:continuous} supplies precisely
the missing justification, without using the conclusion of GKZ.

Notice also that we have not used the general theorem asserting that
every Banach-algebra spectrum is non-empty. Under the GKZ hypothesis,
the explicit scalar $\varphi(a)$ already belongs to $\spA(a)$. For the
arguments below, the upper bound provided by the Neumann series is
sufficient.

\subsection{Exponentials in a Banach algebra}

For $a\in A$ and $z\in\C$, define
\[
e^{za}=\sum_{n=0}^{\infty}\frac{z^na^n}{n!}.
\]
The series converges absolutely and uniformly on bounded subsets of the
$z$-plane, since $\norm{a^n}\le\norm{a}^n$ for all $n\ge0$ and
\[
\sum_{n=0}^{\infty}\frac{|z|^n\norm{a}^n}{n!}
=e^{|z|\norm{a}}.
\]
It follows in particular that
\begin{equation}\label{eq:expnorm}
\norm{e^{za}}\le e^{|z|\norm{a}}.
\end{equation}
We may multiply the two absolutely convergent series for $e^{za}$ and
$e^{-za}$ and collect terms of equal degree. Their constant term is
$\one$, while for $n\ge1$ the coefficient multiplying $z^na^n$ is
\[
\sum_{k=0}^{n}\frac{(-1)^{n-k}}{k!(n-k)!}
=\frac{(1-1)^n}{n!}=0.
\]
Thus
\begin{equation}\label{eq:expinverse}
e^{za}e^{-za}=e^{-za}e^{za}=\one.
\end{equation}
This calculation is valid in a non-commutative algebra: all terms here
are powers of the \emph{same} element, so they commute with one another.

By continuity of $\varphi$,
\begin{equation}\label{eq:Fa-series}
F_a(z):=\varphi(e^{za})
=\sum_{n=0}^{\infty}\frac{\varphi(a^n)}{n!}z^n.
\end{equation}
Since $|\varphi(a^n)|\le\norm{a}^n$, the series has infinite radius of
convergence, and it can be differentiated term by term. In particular,
\begin{equation}\label{eq:derivatives}
F_a^{(k)}(z)=\varphi(a^ke^{za}),\qquad k\ge0,
\end{equation}
and
\begin{equation}\label{eq:firstderivatives}
F_a(0)=1,\qquad F_a'(0)=\varphi(a),\qquad
F_a''(0)=\varphi(a^2).
\end{equation}
These elementary observations support both analytic proofs.

\section{The common algebraic conclusion}
\label{sec:algebra}

\subsection{The square defect}

For a complex-linear map $\varphi:A\to\C$ with $\varphi(\one)=1$, the \emph{square
defect} of $\varphi$ is
\[
d_\varphi(a)=\varphi(a^2)-\varphi(a)^2\qquad(a\in A).
\]
Each of the three proofs is a way of showing that $d_\varphi$ vanishes. The
following elementary lemma explains how the three arguments fit together.

\begin{lemma}[Properties of the square defect]\label{lem:defect}
Let $A$ be a unital complex algebra and let $\varphi:A\to\C$ be
complex-linear with $\varphi(\one)=1$. For $a,b\in A$ and
$\lambda,\mu\in\C$,
\begin{enumerate}[label=\textup{(\roman*)}]
\item $d_\varphi(a+\lambda\one)=d_\varphi(a)$;
\item $d_\varphi(\mu a)=\mu^2d_\varphi(a)$;
\item $d_\varphi(a+b)=d_\varphi(a)+d_\varphi(b)+\varphi(ab+ba)-2\varphi(a)\varphi(b)$.
\end{enumerate}
\end{lemma}

\begin{proof}
(i) Expanding both terms,
\begin{align*}
d_\varphi(a+\lambda\one)
&=\varphi(a^2)+2\lambda\varphi(a)+\lambda^2-(\varphi(a)+\lambda)^2\\
&=\varphi(a^2)-\varphi(a)^2.
\end{align*}
(ii) is immediate from the definition. For (iii), expand
$(a+b)^2=a^2+ab+ba+b^2$ and $(\varphi(a)+\varphi(b))^2$.
\end{proof}

\begin{remark}
Part~(iii) shows that, when $d_\varphi\equiv0$, one obtains the
polarisation identity $\varphi(ab+ba)=2\varphi(a)\varphi(b)$, stated
as~\eqref{eq:jordan} below. Part~(i)
shows that $d_\varphi\equiv0$ as soon as $d_\varphi$ vanishes on a set that
meets every coset $a+\C\one$. This explains the different starting points
of the three proofs. Proof~I shows that
$d_\varphi(a)=0$ for every $a\in A$. Proof~II works on
$\ker\varphi$, which meets every coset because
$a-\varphi(a)\one\in\ker\varphi$; there $d_\varphi(x)=\varphi(x^2)$.
Proof~III works on the invertible elements, which meet every coset in a
Banach algebra because $a-\lambda\one\in\inv$ whenever $|\lambda|>\norm{a}$.
\end{remark}

\subsection{Why polarisation alone is insufficient in the non-commutative case}

Suppose that a complex-linear map $\varphi:A\to\C$ preserves squares.
Expansion of
$(a+b)^2$ gives the polarisation identity
\begin{equation}\label{eq:jordan}
\varphi(ab+ba)=2\varphi(a)\varphi(b).
\end{equation}
If $ab=ba$, this already gives multiplicativity for the pair $a,b$.
Otherwise, \eqref{eq:jordan} determines only
$\varphi(ab)+\varphi(ba)$; one cannot replace $ab+ba$ by $2ab$.

For scalar-valued maps, however, square preservation still implies
multiplicativity. The following purely algebraic lemma appears in
Roitman--Sternfeld~\cite[Lemma~1]{RS}; we spell out its
non-commutative step.

\begin{lemma}[Scalar-valued Jordan maps are multiplicative]
\label{lem:jordan}
Let $A$ be a unital complex algebra and let $\varphi:A\to\C$ be
complex-linear with $\varphi(\one)=1$. The following are equivalent:
\begin{enumerate}[label=\textup{(\roman*)}]
\item If $\varphi(x)=0$, then $\varphi(x^2)=0$.
\item $\varphi(x^2)=\varphi(x)^2$ for every $x\in A$.
\item If $\varphi(x)=0$, then $\varphi(xy)=0$ for every $y\in A$.
\item $\varphi(xy)=\varphi(x)\varphi(y)$ for all $x,y\in A$.
\end{enumerate}
\end{lemma}

\begin{proof}
Assume (i). For arbitrary $a\in A$, put
$x=a-\varphi(a)\one$. Then $\varphi(x)=0$, and hence
\begin{align*}
0=\varphi(x^2)
&=\varphi(a^2)-2\varphi(a)\varphi(a)
  +\varphi(a)^2\varphi(\one)\\
&=\varphi(a^2)-\varphi(a)^2.
\end{align*}
This proves (ii).

Assume (ii). Applying (ii) to $x=a+b$ gives \eqref{eq:jordan} for all $a,b$. The following elegant identity may be viewed as a non-commutative analogue of the parallelogram identity:
\begin{equation}\label{eq:commutator-identity}
(xy+yx)^2 + (xy-yx)^2 =2\bigl(x(yxy)+(yxy)x\bigr).
\end{equation}
The identity is verified by expanding 
$$ (xy+yx)^2 + (xy-yx)^2 = (xy+yx)(xy+yx) + (xy-yx)(xy-yx).$$

Fix $x, y \in A$ with $\varphi(x) = 0$, and set
\[
s = xy + yx, \qquad c = xy - yx,
\]
so that $s$ is the sum of $xy$ and $yx$, and $c$ is the commutator of $x$ and $y$. 
 Equation~\eqref{eq:jordan} gives $\varphi(s)=\varphi(xy+yx)= 2\varphi(x)\varphi(y)=0$, and (ii) then gives
$\varphi(s^2)=0$. Then, the non-commutative parallelogram identity \eqref{eq:commutator-identity} and the polarisation identity \eqref{eq:jordan} imply 
$$\varphi(c^2)= \varphi(c^2+s^2) = 2\varphi\bigl(\bigl(x(yxy)+(yxy)x\bigr)\bigr)= 4\varphi(x)\varphi(yxy) =0.$$
By (ii), $\varphi(c)^2=0$. The range is the scalar field, so
$\varphi(c)=0$. Therefore
\[
2\varphi(xy)=\varphi(s)+\varphi(c)=0,
\]
which proves (iii).

Assume (iii), and take arbitrary $a,b\in A$. Since
$a-\varphi(a)\one$ belongs to the kernel,
\[
0=\varphi\bigl((a-\varphi(a)\one)b\bigr)
=\varphi(ab)-\varphi(a)\varphi(b).
\]
Thus (iv) holds. Finally, (iv) immediately implies (i).
\end{proof}

\begin{remark}[Where the scalar target is used]
The implication $\varphi(c)^2=0\Rightarrow\varphi(c)=0$ uses the fact that
a complex number whose square is zero must itself be zero. Polarisation
also uses the commutativity of scalar multiplication. The analogous
assertion that square-preserving maps into an arbitrary algebra are
multiplicative is false. For instance,
matrix transposition is linear and preserves squares, but
$(ab)^{\mathsf T}=b^{\mathsf T}a^{\mathsf T}$.
\end{remark}

\section{First proof: zero-free entire functions of exponential type}
\label{sec:first}

The basic observation is that $F_a(z)=\varphi(e^{za})$ is zero-free:
$e^{za}$ is invertible, so the hypothesis \eqref{eq:hyp} excludes
$F_a(z)=0$. On the other hand,
\begin{equation}\label{eq:Fgrowth}
|F_a(z)|\le\norm{\varphi}\norm{e^{za}}
\le e^{|z|\norm{a}}.
\end{equation}
The analytic content of the proof is the assertion that these two facts
force $F_a$ to be the exponential of a linear function.

\begin{lemma}[Zero-free entire functions of exponential type]\label{lem:entire}
Let $F$ be an entire function with $F(0)=1$ and no zeros. Suppose that
there exist constants $C\ge1$ and $T\ge0$ such that
\begin{equation}\label{eq:etype}
|F(z)|\le C e^{T|z|}\qquad(z\in\C).
\end{equation}
Then there is $\alpha\in\C$ such that
\[
F(z)=e^{\alpha z}\qquad(z\in\C).
\]
Moreover, $\alpha=F'(0)$ and $|\alpha|\le T$.
\end{lemma}

This is the special zero-free case of the factorisation theory of entire
functions that is needed in the Kahane--\.{Z}elazko argument. Their paper
invokes a standard result on entire functions, taken from
Titchmarsh~\cite{Titchmarsh}. The full Hadamard factorisation theorem is
considerably more general: it describes an entire function in terms of its
zeros and its growth. Here there are no zeros and the growth is at most
exponential, so we can prove the required special case directly.

\subsection{Step 1: constructing a global logarithm}

Since $F$ has no zeros, $F'/F$ is entire, so it is the sum of a power
series with infinite radius of convergence,
$F'(z)/F(z)=\sum_{n\ge0}c_nz^n$. Integrating term by term gives an entire
function
\begin{equation}\label{eq:logprimitive}
h(z)=\sum_{n\ge0}\frac{c_n}{n+1}\,z^{n+1},\qquad h(0)=0,\qquad
h'=\frac{F'}{F}.
\end{equation}
Differentiating $Fe^{-h}$ gives
\[
\frac{d}{dz}\bigl(F(z)e^{-h(z)}\bigr)
=e^{-h(z)}\bigl(F'(z)-F(z)h'(z)\bigr)=0.
\]
Hence $Fe^{-h}$ is constant; its value at $0$ is $1$, and consequently
\begin{equation}\label{eq:global-log}
F(z)=e^{h(z)},\qquad h(0)=0.
\end{equation}
Thus $h$ is an entire logarithm of $F$, normalised by $h(0)=0$.
From \eqref{eq:etype},
\begin{equation}\label{eq:real-growth}
\Rea h(z)=\log|F(z)|\le\log C+T|z|.
\end{equation}

\subsection{Step 2: a coefficient estimate from the real part}

We prove the elementary form of Carath\'eodory's coefficient estimate
needed here.

\begin{lemma}[A real-part coefficient estimate]\label{lem:carath}
Let $r>0$, and let $h$ be holomorphic on a neighbourhood of the closed
disc $\{|z|\le r\}$, with $h(0)=0$. Write
$h(z)=\sum_{n\ge1}a_nz^n$. If $M\in\R$ and $\Rea h(z)\le M$ on
$|z|=r$, then
\begin{equation}\label{eq:carath}
|a_n|\le\frac{2M}{r^n}\qquad(n\ge1).
\end{equation}
\end{lemma}

\begin{proof}
Put $u(t)=\Rea h(re^{it})$. The power series converges uniformly on the
circle, so the Fourier coefficients of $u$ can be computed by termwise
integration.
For $n\ge1$,
\begin{equation}\label{eq:fourier}
\frac1{2\pi}\int_0^{2\pi}u(t)e^{-int}\,dt
=\frac{a_nr^n}{2}.
\end{equation}
Indeed, in
\[
u(t)=\frac12\left(
\sum_{k\ge1}a_kr^ke^{ikt}
+\sum_{k\ge1}\overline{a_k}r^ke^{-ikt}\right),
\]
only the term with $k=n$ in the first sum has a non-zero integral after
multiplication by $e^{-int}$. Also,
\[
\frac1{2\pi}\int_0^{2\pi}u(t)\,dt=\Rea h(0)=0.
\]
The function $M-u(t)$ is non-negative and has mean $M$; in particular,
$M\ge0$. The constant $M$
has zero $n$th Fourier coefficient for $n\ge1$, hence \eqref{eq:fourier}
implies
\begin{align*}
\frac{|a_n|r^n}{2}
&=\left|\frac1{2\pi}\int_0^{2\pi}(M-u(t))e^{-int}\,dt\right|\\
&\le\frac1{2\pi}\int_0^{2\pi}(M-u(t))\,dt=M.
\end{align*}
This proves \eqref{eq:carath}.
\end{proof}

The point of this estimate is that it controls \emph{all} coefficients
from an upper bound on the real part. A direct Cauchy estimate on $h$
would require a bound on $|h|$, which \eqref{eq:real-growth} does not
immediately provide.

\subsection{Step 3: the logarithm is linear}

\begin{proof}[Proof of Lemma~\ref{lem:entire}]
Construct $h$ as in \eqref{eq:global-log}. Apply
Lemma~\ref{lem:carath} on every circle of radius $r>0$, with
$M=\log C+Tr$. For the Taylor coefficients of $h$ we obtain
\[
|a_n|\le\frac{2(\log C+Tr)}{r^n}.
\]
For each $n\ge2$, the right-hand side tends to zero as $r\to\infty$. Thus
$a_n=0$ for all $n\ge2$, and $h(z)=a_1z$. Taking
$\alpha=a_1$ gives $F(z)=e^{\alpha z}$. Differentiation at $0$ gives
$\alpha=F'(0)$.

If $\alpha\ne0$, choose $z=r\overline\alpha/|\alpha|$. Then
$\Rea(\alpha z)=r|\alpha|$, and \eqref{eq:etype} gives
$r|\alpha|\le\log C+Tr$. Divide by $r$ and let $r\to\infty$ to obtain
$|\alpha|\le T$. The case $\alpha=0$ is immediate.
\end{proof}

\subsection{Completion of the first proof of GKZ}

\begin{proof}[Proof I of Theorem~\ref{thm:gkz}]
Proposition~\ref{prop:continuous} gives automatic continuity. Fix any
$a\in A$ and define $F_a$ by \eqref{eq:Fa-series}. It is entire,
$F_a(0)=1$, it has no zeros, and it satisfies \eqref{eq:Fgrowth}.
Lemma~\ref{lem:entire} gives
\[
F_a(z)=e^{\alpha z}.
\]
By \eqref{eq:firstderivatives}, $\alpha=F_a'(0)=\varphi(a)$, and
therefore
\begin{equation}\label{eq:expidentity}
\varphi(e^{za})=e^{z\varphi(a)}\qquad(z\in\C).
\end{equation}
Comparing the second derivatives at zero yields
\[
\varphi(a^2)=F_a''(0)=\varphi(a)^2.
\]
Since $a$ was arbitrary, Lemma~\ref{lem:jordan} implies
multiplicativity.
\end{proof}

In fact, comparison of all coefficients gives
\begin{equation}\label{eq:powers}
\varphi(a^n)=\varphi(a)^n\qquad(n\ge0).
\end{equation}
Only the case $n=2$ was needed to finish the proof.

\begin{remark}[Normalising the logarithm]
If $F(z)=e^{\alpha z+\beta}$ and $F(0)=1$, the exact
conclusion is $e^\beta=1$, or $\beta\in2\pi i\mathbb Z$. It follows that
$F(z)=e^{\alpha z}$, and one may \emph{choose} the logarithm with constant
term zero. It does not follow that an arbitrarily chosen $\beta$ equals
zero. Our construction $h(0)=0$ removes this ambiguity from the start.
\end{remark}

\subsection{A finite-radius version of the coefficient estimate}

Lemma~\ref{lem:carath} does not require an entire function. Let $R>0$, and
suppose that $F$ is holomorphic and zero-free on the disc $|z|<R$, with
$F(0)=1$. Assume that $C\ge1$ and $T\ge0$ satisfy
\[
|F(z)|\le C e^{T|z|}\qquad(|z|<R).
\]
Since $F'/F$ is holomorphic on the disc, integrating its Taylor series
term by term gives a holomorphic function $h$ on the disc with $F=e^h$
and $h(0)=0$, exactly as in \eqref{eq:logprimitive}--\eqref{eq:global-log}.
Write $h(z)=c_1z+c_2z^2+\cdots$. For every $0<r<R$,
Lemma~\ref{lem:carath} gives
\[
|c_2|\le\frac{2(\log C+Tr)}{r^2}.
\]
Since
\[
h''(0)=F''(0)-F'(0)^2=2c_2,
\]
we obtain
\begin{equation}\label{eq:finite-radius}
|F''(0)-F'(0)^2|
\le\frac{4(\log C+Tr)}{r^2}\qquad(0<r<R).
\end{equation}
For $F(z)=\varphi(e^{za})$, the left-hand side is
$|d_\varphi(a)|$, the modulus of the square defect.
Under the GKZ hypotheses, $F_a$ is entire and zero-free
on the whole plane, and the same growth bound holds on every disc. We may
therefore let $r\to\infty$ and obtain $d_\varphi(a)=0$; this recovers
Proof~I. The finite-radius inequality is relevant when only local
zero-freeness or local growth information is available.

One might try to weaken the GKZ hypothesis by assuming only local
non-vanishing near the identity. More precisely, suppose that $\varphi$
is complex-linear, $\varphi(\one)=\norm{\varphi}=1$, and
\[
\varphi(u)\ne0\quad\text{whenever}\quad\norm{u-\one}<\delta,
\qquad \delta>0.
\]
We retain $\norm{\one}=1$. If $a=0$, then $d_\varphi(a)=0$, so assume
$a\ne0$. The exponential estimate gives
\[
\norm{e^{za}-\one}\le e^{|z|\norm{a}}-1<\delta
\quad\text{for}\quad
|z|<R:=\frac{\log(1+\delta)}{\norm{a}}.
\]
Thus $F_a$ is zero-free on this disc. Applying
\eqref{eq:finite-radius} with $C=1$ and $T=\norm{a}$, then letting
$r\uparrow R$, yields
\[
|d_\varphi(a)|\le\frac{4\norm{a}^2}{\log(1+\delta)}.
\]
However, since $\varphi(0)=0$, the open ball cannot contain $0$, so
$\delta\le1$. The right-hand side is therefore at least
$4\norm{a}^2/\log2>5\norm{a}^2$, whereas the elementary estimate
\[
|d_\varphi(a)|
\le |\varphi(a^2)|+|\varphi(a)|^2
\le2\norm{a}^2
\]
is already better. Indeed, the local non-vanishing assumption adds no
information when $\delta\le1$: every normalised norm-one functional
satisfies
\[
|\varphi(u)-1|=|\varphi(u-\one)|
\le\norm{u-\one}<1
\]
on the open unit ball centred at $\one$.

This calculation therefore does not yield a quantitative version of GKZ
under a local non-vanishing hypothesis. What it does show is exactly what
the real-part estimate gives on a fixed disc, and why letting $r\to\infty$
is decisive in Proof~I.

\subsection{Why the exponential function?}

The preceding proof yields a stronger formulation: it suffices to assume
that $\varphi$ does not vanish on
\[
    \exp(A)=\{e^a:a\in A\},
\]
the set of exponentials in $A$.
Indeed, once continuity is known, Proof~I uses the non-vanishing hypothesis only
at the exponentials $e^{za}$; continuity itself, through
Proposition~\ref{prop:continuous}, uses it only at the elements
$\lambda\one-a$ with $|\lambda|>\norm{a}$. These elements are exponentials
too, since $\one-b$ has the logarithm $-\sum_{n\ge1}b^n/n$ whenever
$\norm{b}<1$.

In 1987, Richard Arens~\cite{Arens} asked whether the exponential function
could be replaced by other non-constant entire functions.
This problem has been studied, and solved, by
Krzysztof Jarosz~\cite{Jarosz}. 

\section{Second proof: polynomial roots and bounded spectra}
\label{sec:second}

This proof does not construct an exponential and does not need the
functional to be continuous. Its key device is a polynomial of increasing
degree whose roots are forced into one fixed spectrum.

\subsection{The theorem in its natural algebraic setting}

\begin{theorem}[Roitman--Sternfeld, bounded-spectrum form]
\label{thm:RS}
Let $A$ be a unital complex algebra. Suppose that $\spA(x)$ is bounded
for every $x\in A$. If a complex-linear functional $\varphi:A\to\C$
satisfies \eqref{eq:hyp}, then $\varphi$ is multiplicative.
\end{theorem}

This is the corollary to Lemmas~1 and~2 of
Roitman--Sternfeld~\cite[p.~113]{RS}. There is no topology in the
statement: the algebra need not be normed or complete. Boundedness is a
property of the subsets $\spA(x)\subset\C$,
and the bound may depend on $x$. In a Banach algebra it follows from the
Neumann series, so Theorem~\ref{thm:RS} immediately implies GKZ.
Although boundedness alone allows an empty spectrum in an arbitrary
algebra, \eqref{eq:hyp} rules this out: Lemma~\ref{lem:selection} gives
$\varphi(x)\in\spA(x)$ for every $x$.

\subsection{The polynomials and the location of their roots}

Fix $x\in\ker\varphi$, and for each integer $n\ge2$ consider
\begin{equation}\label{eq:pn}
P_n(\lambda)=\varphi\bigl((\lambda\one-x)^n\bigr).
\end{equation}
Since $\lambda\one$ commutes with $x$, the binomial formula is valid even
when $A$ is non-commutative. Therefore
\begin{equation}\label{eq:pncoef}
P_n(\lambda)
=\lambda^n-n\varphi(x)\lambda^{n-1}
+\binom n2\varphi(x^2)\lambda^{n-2}+\cdots.
\end{equation}
The leading coefficient is $\varphi(\one)=1$, so $P_n$ is monic of
degree exactly $n$.

By the fundamental theorem of algebra, write
\begin{equation}\label{eq:rootfactor}
P_n(\lambda)=\prod_{j=1}^{n}(\lambda-\lambda_{n,j}),
\end{equation}
where the roots are counted with algebraic multiplicity. If $P_n(\lambda)=0$
and $\lambda\one-x$ were invertible, then its power
$(\lambda\one-x)^n$ would be invertible, contradicting
\eqref{eq:hyp}. Hence
\begin{equation}\label{eq:rootlocation}
\lambda_{n,j}\in\spA(x)\qquad(1\le j\le n).
\end{equation}
This step uses only the fact that powers of invertible elements are
invertible; it does not assume multiplicativity of $\varphi$.

\subsection{The first two Newton identities}

Comparing coefficients in \eqref{eq:pncoef} and \eqref{eq:rootfactor},
and using $\varphi(x)=0$, gives
\begin{align}
\sum_{j=1}^{n}\lambda_{n,j}&=0,\label{eq:sumroots}\\
\sum_{1\le i<j\le n}\lambda_{n,i}\lambda_{n,j}
&=\binom n2\varphi(x^2).\label{eq:pairroots}
\end{align}
Squaring \eqref{eq:sumroots} yields
\[
0=\sum_{j=1}^{n}\lambda_{n,j}^2
  +2\sum_{1\le i<j\le n}\lambda_{n,i}\lambda_{n,j}.
\]
Thus
\begin{equation}\label{eq:secondmoment}
\sum_{j=1}^{n}\lambda_{n,j}^2
=-n(n-1)\varphi(x^2).
\end{equation}
These are the first two Newton relations in elementary form. We do not
need the general Newton identities.

Choose a finite number $R_x\ge0$ such that $|\lambda|\le R_x$ for every
$\lambda\in\spA(x)$. By \eqref{eq:rootlocation} and
\eqref{eq:secondmoment},
\begin{align}
n(n-1)|\varphi(x^2)|
&=\left|\sum_{j=1}^{n}\lambda_{n,j}^2\right|\notag\\
&\le\sum_{j=1}^{n}|\lambda_{n,j}|^2
\le nR_x^2.\label{eq:RSestimate}
\end{align}
Consequently,
\begin{equation}\label{eq:RSlimit}
|\varphi(x^2)|\le\frac{R_x^2}{n-1}
\qquad(n\ge2).
\end{equation}
The element $x$ and the bound $R_x$ are fixed, while $n$ is arbitrary.
Letting $n\to\infty$ gives $\varphi(x^2)=0$.

\begin{proof}[Proof of Theorem~\ref{thm:RS}, and Proof II of GKZ]
The preceding argument proves condition (i) of
Lemma~\ref{lem:jordan}: $\varphi(x)=0$ implies $\varphi(x^2)=0$.
The lemma therefore gives multiplicativity. In a Banach algebra every
spectrum is bounded by \eqref{eq:specbound}; this proves
Theorem~\ref{thm:gkz}. Automatic continuity is supplied separately by
Proposition~\ref{prop:continuous}.
\end{proof}

Within this proof, the only non-trivial fact which is used is the
fundamental theorem of algebra: a monic polynomial of degree $n$ has
$n$ complex roots counted with multiplicity. The coefficient comparisons
and estimates are elementary algebra.

The fundamental theorem of algebra is often proved using complex
analysis, for example by Liouville's theorem. Other proofs exist. It is
therefore accurate to say that the Roitman--Sternfeld argument needs no
complex-analytic apparatus beyond the fundamental theorem of algebra.
Completeness enters only if one first derives the bounded-spectrum
hypothesis from the Banach-algebra axioms.

\section{Third proof: P\'olya--Saxer and zero-free derivatives}
\label{sec:third}

The first proof controls an entire function by a bound on its growth.
The proof given in~\cite{Badea} instead controls it through the absence of
zeros of its derivatives. The starting exponential is the same, but the
rigidity theorem is different.

\subsection{The exponential function and absence of zeros}

\begin{theorem}[P\'olya--Saxer]
\label{thm:polya}
Let $F$ be entire. If $F$, $F'$, and $F''$ have no zeros in $\C$, then
\[
F(z)=e^{\alpha z+\beta}
\]
for constants $\alpha\in\C\setminus\{0\}$ and $\beta\in\C$.
\end{theorem}

The statement and the historical account below are recalled
in~\cite[Section~2]{Badea}, which refers to
Burckel~\cite[Chapter~12]{Burckel} for the history and complex-analytic
background. The theorem was stated by P\'olya in the 1920s and proved by
Saxer, Varopoulos and Csillag. It is considerably deeper than the
coefficient estimate of Section~\ref{sec:first}. We use it without proof; no growth hypothesis appears in its statement.

In the application below, \emph{every} derivative of $F$ is zero-free.
It would therefore suffice to quote the weaker version of the theorem, in
which the derivatives of all orders, including order zero, are assumed to
be zero-free; this is the version used explicitly in the proof
in~\cite{Badea}. The stronger statement, which involves only $F$ and its
first two derivatives, makes clear that the remaining derivatives are not
needed for this route to GKZ.

In fact, fewer non-vanishing conditions suffice. The following theorem of
Tumura and Clunie is also stated in~\cite[Section~2]{Badea}; for $n=2$,
a short proof was given by Yang~\cite{Yang}.

\begin{theorem}[Tumura--Clunie]
\label{thm:tumura}
Let $F$ be entire and let $n\ge2$. If $F$ and $F^{(n)}$ have no zeros in
$\C$, then
\[
F(z)=e^{\alpha z+\beta}
\]
for constants $\alpha\in\C\setminus\{0\}$ and $\beta\in\C$.
\end{theorem}

The case $n=2$ strengthens Theorem~\ref{thm:polya}, since the hypothesis
on $F'$ is unnecessary. Proof~III therefore only requires that $F_x$ and
$F_x''$ be zero-free. The example below shows that a condition on a higher
derivative cannot simply be omitted: it is not enough that $F$ and $F'$ be zero-free.

\subsection{Why non-vanishing of \texorpdfstring{$F$ and $F'$}{F and F'} alone is insufficient}

Consider
\begin{equation}\label{eq:doubleexp}
F(z)=\exp(e^z-1).
\end{equation}
Then $F(0)=1$, $F$ has no zeros, and
\[
F'(z)=e^z\exp(e^z-1)
\]
also has no zeros. Nevertheless, $F$ is not $e^{\alpha z}$ for a
constant $\alpha$, since $F'/F=e^z$ is not constant. Its second
derivative is
\[
F''(z)=e^z(1+e^z)\exp(e^z-1),
\]
which vanishes at $z=(2k+1)\pi i$, $k\in\mathbb Z$.
Thus the second-derivative condition imposes a genuine restriction.

The example also has much faster growth than allowed by
Lemma~\ref{lem:entire}: along the positive real axis,
$\log|F(r)|=e^r-1$, which cannot be bounded by $\log C+Tr$.
The two analytic lemmas therefore exclude it for different reasons.

\subsection{Step 1: start with an invertible algebra element}

Assume \eqref{eq:hyp}, and let $x\in\inv$. Set
\[
F_x(z)=\varphi(e^{zx}).
\]
Automatic continuity and the identities \eqref{eq:Fa-series} and
\eqref{eq:derivatives} show that this function is entire and that
\begin{equation}\label{eq:allnonzero}
F_x^{(k)}(z)=\varphi(x^ke^{zx})\ne0
\qquad(k\ge0,\ z\in\C).
\end{equation}
Indeed, both $x^k$ and $e^{zx}$ are invertible. Their product is invertible
as well, with inverse $e^{-zx}x^{-k}$. The hypothesis on $\varphi$ now
gives \eqref{eq:allnonzero}.

The restriction $x\in\inv$ is essential at this step. For arbitrary
$a\in A$, $e^{za}$ is still invertible, but $ae^{za}$ need not be.
Thus the hypothesis guarantees $F_a\ne0$, but it does not automatically
guarantee $F_a'\ne0$.

\subsection{Step 2: apply P\'olya--Saxer}

By \eqref{eq:allnonzero}, Theorem~\ref{thm:polya} gives
$F_x(z)=e^{\alpha z+\beta}$. Since $F_x(0)=1$, we have $e^\beta=1$,
and hence $F_x(z)=e^{\alpha z}$. Differentiation at zero yields
\[
\alpha=F_x'(0)=\varphi(x),\qquad
\varphi(x^2)=F_x''(0)=\varphi(x)^2.
\]
We have therefore proved
\begin{equation}\label{eq:invsquares}
\varphi(x^2)=\varphi(x)^2\qquad(x\in\inv).
\end{equation}
The non-vanishing hypothesis also gives $\varphi(x)\ne0$, consistently
with $\alpha\ne0$ in Theorem~\ref{thm:polya}. If one prefers to quote
Theorem~\ref{thm:tumura} with $n=2$, only the cases $k=0$ and $k=2$ of
\eqref{eq:allnonzero} are used here.

\subsection{Step 3: a scalar shift removes the invertibility restriction}

Take any $a\in A$. Choose $\lambda\in\C$ with
$|\lambda|>\norm{a}$. The Neumann series shows that
$x=a-\lambda\one$ is invertible. Apply \eqref{eq:invsquares} to $x$:
\[
\varphi\bigl((a-\lambda\one)^2\bigr)
=\varphi(a-\lambda\one)^2.
\]
The left side is
\[
\varphi(a^2)-2\lambda\varphi(a)+\lambda^2,
\]
and the right side is
\[
(\varphi(a)-\lambda)^2
=\varphi(a)^2-2\lambda\varphi(a)+\lambda^2.
\]
Cancelling the last two terms on each side gives
\begin{equation}\label{eq:allsquaresagain}
\varphi(a^2)=\varphi(a)^2.
\end{equation}
In the language of Lemma~\ref{lem:defect}(i), this says that
$d_\varphi(a)=d_\varphi(a-\lambda\one)=0$.

\begin{proof}[Proof III of Theorem~\ref{thm:gkz}]
Proposition~\ref{prop:continuous} gives continuity. For invertible $x$,
the derivative argument and P\'olya--Saxer give \eqref{eq:invsquares}.
For arbitrary $a$, the scalar-shift argument gives
\eqref{eq:allsquaresagain}. Lemma~\ref{lem:jordan} now proves
multiplicativity.
\end{proof}

\section{What the theorem is useful for}
\label{sec:applications}

\subsection{Hyperplanes consisting of non-invertible elements}

\begin{corollary}\label{cor:hyperplane}
Let $A$ be a unital complex Banach algebra, and let $H\subset A$ be a
linear subspace of algebraic codimension one. The following are
equivalent:
\begin{enumerate}[label=\textup{(\roman*)}]
\item $H\cap\inv=\varnothing$;
\item $H=\ker\varphi$ for a character $\varphi$;
\item $H$ is a proper two-sided ideal.
\end{enumerate}
In this case $H$ is closed, $A/H$ is isomorphic to $\C$, and $H$ is a
maximal proper two-sided ideal. It is also a maximal proper left ideal
and a maximal proper right ideal.
\end{corollary}

\begin{proof}
Assume (i). Since $\one$ is invertible, $\one\notin H$. Codimension one
gives the algebraic direct sum $A=H\oplus\C\one$. Define
$\varphi(h+\lambda\one)=\lambda$. This is linear, normalised, and
has kernel $H$. If $u$ is invertible, then $u\notin H$, so
$\varphi(u)\ne0$. GKZ gives (ii), including continuity of $\varphi$;
its kernel is therefore closed.

A character's kernel is a proper two-sided ideal, so (ii) implies (iii).
A proper two-sided ideal cannot contain an invertible element, since
$u\in H\cap\inv$ would imply $\one=u^{-1}u\in H$. Thus (iii) implies
(i). The character induces an algebra isomorphism $A/H\to\C$.
Finally, any linear subspace properly containing a codimension-one
subspace equals $A$; this proves all the maximality assertions.
\end{proof}

For commutative Banach algebras, this is the hyperplane characterisation
of maximal ideals referred to in the title of Kahane and
\.{Z}elazko's paper~\cite{KZ}. The codimension-one hypothesis is essential. For example, in
$A=\C^2$ with coordinatewise operations, the zero subspace $H=\{0\}$
satisfies (i) and (iii) but not (ii): it is not the kernel of a character,
and it is not a maximal ideal, since it is properly contained in
$\C\times\{0\}$.

\subsection{A non-vanishing test for point evaluation on \texorpdfstring{$C(X)$}{C(X)}}

Let $X$ be a non-empty compact Hausdorff space. In $C(X)=C_{\C}(X)$, with the
supremum norm and pointwise multiplication, the invertible elements are
exactly the nowhere-zero functions. This leads to the following characterisation of Dirac measures.

\begin{corollary}\label{cor:CX}
Let $\Lambda:C(X)\to\C$ be complex-linear, with
$\Lambda(\one)=1$. If $\Lambda(f)\ne0$ whenever $f$ is nowhere zero,
then there is a unique $x_0\in X$ such that
\[
\Lambda(f)=f(x_0)\qquad(f\in C(X)).
\]
\end{corollary}

\begin{proof}
By GKZ, $\Lambda$ is a character. To show that it is a point evaluation,
put $H=\ker\Lambda$ and suppose, for a contradiction, that the functions
in $H$ have no common zero. For each $x\in X$, choose
$f_x\in H$ with $f_x(x)\ne0$. The open sets where these functions are
non-zero cover $X$. Compactness supplies finitely many
$f_1,\ldots,f_m\in H$ with no common zero. Then
\[
g=\sum_{j=1}^{m}\overline{f_j}f_j
\]
is strictly positive everywhere, hence invertible in $C(X)$. Since
$H$ is an ideal, $g\in H$, a contradiction. There is therefore a point
$x_0$ at which every element of $H$ vanishes.

For any $f\in C(X)$, the function $f-\Lambda(f)\one$ lies in $H$.
Evaluation at $x_0$ gives $f(x_0)=\Lambda(f)$. Uniqueness follows
because continuous functions separate distinct points of a compact
Hausdorff space.
\end{proof}

Thus the non-vanishing condition singles out point evaluations among
normalised (unital) linear functionals. In particular, a non-trivial convex
combination of evaluations at distinct points must vanish on some
nowhere-zero complex-valued function.

The following reformulation of Corollary \ref{cor:CX} uses the 'selection of the spectrum' point of view.

\begin{corollary}\label{cor:measure}
Let $X$ be a compact Hausdorff space and let $\mu$
be a complex regular Borel measure on $X$, such that the mean value property
\[
  \int_X f\,d\mu \in \image(f)
\]
holds for all $f$ in $C_{\mathbb{C}}(X)$. Then $\mu$ is a Dirac measure;
i.e.\ there is a point $x_0$ in $X$, such that
\[
  \int_X f\,d\mu = f(x_0) \qquad \forall f \in C_{\mathbb{C}}(X).
\]
\end{corollary}

Following \cite{Wille}, we give an independent proof of this corollary that does not rely on the Gleason--Kahane--\.Zelazko theorem.

\begin{proof}[Proof of Corollary \ref{cor:measure}]
The mean value property implies at once that $\mu$ is positive and that
$\mu(X)=1$. Therefore it suffices to show that the support of $\mu$ contains
only one point.

Suppose that there are two distinct points $x_1, x_2$ in the support of $\mu$.
Then there are also two open neighbourhoods $U_1$ of $x_1$ and $U_2$ of $x_2$,
which have no common points.

Furthermore, since $X$ is a compact Hausdorff space, Urysohn's lemma guarantees
the existence of two continuous functions $g_1\colon X\to[0,1]$,
$g_2\colon X\to[0,1]$ with the following properties ($i=1,2$):
\begin{align*}
  g_i(x) &= 0 \qquad \forall x\in X\setminus U_i,\\
  g_i(x_i) &= 1.
\end{align*}

Now regard $g_1$ and $g_2$ as the real resp.\ imaginary part of one
complex-valued continuous function on $X$. The mean value property implies that
there exists a point $y_0\in X$ such that
\[
  \int_X (g_1 + i g_2)\,d\mu = g_1(y_0) + i g_2(y_0).
\]
Therefore $\int_X g_1\,d\mu = g_1(y_0)$ and $\int_X g_2\,d\mu = g_2(y_0)$.

Since $\mu$ is positive and $x_1, x_2$ are in the support of $\mu$, the
integrals of $g_1$ and $g_2$ are strictly positive. Consequently $g_1(y_0)>0$
and $g_2(y_0)>0$ hold.

Now $g_1$ vanishes in $X\setminus U_1$, hence $y_0$ must be an element of
$U_1$. But $g_2$ vanishes in $X\setminus U_2$, so $y_0$ is also in $U_2$. This
contradicts the assumption that $U_1$ and $U_2$ have no common points. Thus the
support of $\mu$ contains only one point, which means that $\mu$ is the Dirac measure at that point.
\end{proof}

\subsection{Weighted composition on spaces of continuous functions}

Applying Corollary~\ref{cor:CX} at each point of the target space gives
a representation for linear maps that preserve nowhere-zero functions.

\begin{corollary}\label{cor:weighted}
Let $X,Y$ be non-empty compact Hausdorff spaces. Suppose a complex-linear
map $T:C(X)\to C(Y)$ takes every nowhere-zero function to a
nowhere-zero function. Then there are a nowhere-zero function
$\psi\in C(Y)$ and a continuous map $\theta:Y\to X$ such that
\[
(Tf)(y)=\psi(y)f(\theta(y))\qquad(f\in C(X),\ y\in Y).
\]
\end{corollary}

\begin{proof}
Set $\psi=T\one$, which is nowhere zero. For each $y\in Y$, the map
\[
\Lambda_y(f)=\frac{(Tf)(y)}{\psi(y)}
\]
is linear, normalised, and non-zero on every nowhere-zero function.
Corollary~\ref{cor:CX} gives a unique point $\theta(y)\in X$ with
$\Lambda_y(f)=f(\theta(y))$. This proves the displayed formula.
For every $f\in C(X)$, the function $f\circ\theta=Tf/\psi$ is
continuous. Since the topology of the compact Hausdorff space $X$ is the
weakest topology making every $f\in C(X)$ continuous, $\theta$ is
continuous.
\end{proof}

Continuity of $T$ follows from the representation:
$\norm{Tf}_\infty\le\norm{\psi}_\infty\norm{f}_\infty$. In fact,
$\norm{T}=\norm{\psi}_\infty$, since equality is attained at $f=\one$.

The study of composition operators and their weighted counterparts,
particularly on spaces of holomorphic functions, is an active area
of research in operator theory. For a recent authoritative monograph,
see~\cite{LQ}.

\subsection{Linear maps preserving invertibility}

\begin{proposition}\label{prop:preserver}
Let $A,B$ be unital complex Banach algebras, with $B$ commutative and
semisimple. Suppose $T:A\to B$ is complex-linear, $T(\one_A)=\one_B$,
and $T(u)$ is invertible in $B$ whenever $u$ is invertible in $A$.
Then $T$ is multiplicative.
\end{proposition}

\begin{proof}
For each character $\chi$ of $B$, the functional $\chi\circ T$ is
linear, normalised, and non-zero on every invertible element of $A$.
By GKZ,
\[
\chi(T(ab))=\chi(T(a))\chi(T(b))
=\chi(T(a)T(b)).
\]
Thus every character of $B$ vanishes on $T(ab)-T(a)T(b)$. In a
commutative semisimple Banach algebra the characters separate points,
so this difference is zero.
\end{proof}

The proof illustrates a useful method: compose an algebra-valued map
with enough characters to separate points in its range. If $B$ is
commutative but not semisimple, the same argument shows only that
$T(ab)-T(a)T(b)$ belongs to its Jacobson radical, the intersection of
the kernels of its characters. For a non-commutative target, characters
need not separate points and may not exist at all. The argument therefore
does not settle the corresponding questions for non-commutative target
algebras.

The \emph{Kaplansky--Aupetit problem} asks for the following
extension~\cite{BresarSemrl2008}. Let $A,B$ be unital complex Banach
algebras, with $B$ semisimple, and let $T:A\to B$ be a surjective
complex-linear map satisfying $T(\one_A)=\one_B$. If $T(u)$ is
invertible whenever $u$ is invertible, must $T$ be a \emph{Jordan
homomorphism}, that is, must
\[
T(a^2)=T(a)^2\qquad(a\in A)?
\]
When $B$ is commutative, Proposition~\ref{prop:preserver} gives the
stronger conclusion of multiplicativity without requiring surjectivity.
For a non-commutative target, matrix transposition shows why
multiplicativity itself cannot be expected. The Kaplansky--Aupetit problem is 
still open.

\section{Why the hypotheses matter}
\label{sec:limits}

\subsection{The real-field version fails}

\begin{example}\label{ex:real}
Take the real Banach algebra $A=C([0,1],\R)$ and define
\[
\varphi(f)=\int_0^1 f(t)\,dt.
\]
Then $\varphi(\one)=1$. Every invertible $f\in A$ is continuous and
nowhere zero; connectedness of $[0,1]$ forces it to have a fixed sign.
Consequently $\varphi(f)\ne0$. Nevertheless, for $a(t)=t$,
\[
\varphi(a^2)=\frac13\ne\frac14=\varphi(a)^2.
\]
Thus the real analogue of GKZ is false.
\end{example}

This example appears in the Kahane--\.{Z}elazko paper~\cite{KZ}.
On the complex algebra $C([0,1],\C)$ the same integral does \emph{not}
satisfy the hypothesis: $f(t)=e^{2\pi it}$ is nowhere zero but has
integral zero. The values of a complex-valued nowhere-zero function can therefore
cancel under integration, whereas those of a continuous real-valued
function of fixed sign cannot.

\subsection{An incomplete normed algebra where the conclusion fails}

\begin{example}\label{ex:polynomials}
Let $A=\C[t]$, the polynomial algebra, equipped with
\[
\norm{p}=\max_{|z|\le1}|p(z)|.
\]
This is a unital normed algebra, but it is not complete. Its only
invertible elements are the non-zero constant polynomials. Define
\[
\varphi(p)=\frac{p(1)+p(-1)}2.
\]
Then $\varphi(\one)=1$, $\norm{\varphi}=1$, and $\varphi$ is non-zero on
all invertible elements of $A$. Yet
\[
\varphi(t)=0,\qquad\varphi(t^2)=1,
\]
so $\varphi$ is not multiplicative.
\end{example}

For a non-constant polynomial $p$, every $p-\lambda$ remains non-constant
and is non-invertible in $\C[t]$. Hence
\[
\sigma_{\C[t]}(p)=\C.
\]
Thus the bounded-spectrum hypothesis of Theorem~\ref{thm:RS} fails.
Completeness can be replaced by a suitable spectral hypothesis, as
Roitman and Sternfeld showed; see also Corollary~\ref{cor:Q} below.
The example shows why some such replacement is needed.

\subsection{Why passage to the completion needs care}

The completion in Example~\ref{ex:polynomials} is the disc algebra,
consisting of the functions that are continuous on the closed unit disc
and holomorphic in its interior. The bounded functional $\varphi$
extends to it by the same averaging formula. However, the completion
contains many new invertible elements. For example,
\[
f(z)=e^{i\pi z/2}
\]
belongs to the disc algebra and has an inverse there, whereas
\[
\frac{f(1)+f(-1)}2=\frac{i-i}{2}=0.
\]
Thus the extended functional does not satisfy the GKZ non-vanishing
hypothesis. Continuity of a functional is not enough to transfer that
hypothesis from an incomplete algebra to its completion. An argument using the completion must therefore account for the
relation between the two sets of invertible elements, or between their
spectra.

\subsection{A positive counterpart: normed Q-algebras and spectral norms}

Theorem~\ref{thm:RS} applies whenever the spectra are bounded, even if
completeness is unavailable. The following two conditions ensure this
boundedness.

\begin{corollary}\label{cor:Q}
Let $A$ be a unital complex algebra and let $\varphi:A\to\C$ be a linear
functional satisfying \eqref{eq:hyp}. Suppose that one of the following
holds.
\begin{enumerate}[label=\textup{(\alph*)}]
\item $A$ carries a norm, not necessarily complete or submultiplicative,
such that $\inv$ contains a ball $\{u:\norm{u-\one}<\varepsilon\}$ for some
$\varepsilon>0$. This is the case for every normed Q-algebra, that is,
whenever $\inv$ is open, and for every Banach algebra
(with $\varepsilon=1$, by Lemma~\ref{lem:neumann}).
\item $A$ carries a function $p:A\to[0,\infty)$ with $|\lambda|\le p(x)$
for all $x\in A$ and $\lambda\in\spA(x)$; for example, a spectral seminorm
in the sense of~\cite{Palmer}.
\end{enumerate}
Then $\varphi$ is multiplicative.
\end{corollary}

\begin{proof}
In case (b), every spectrum is bounded by hypothesis. In case (a),
let $x\in A$ and
$|\lambda|>\norm{x}/\varepsilon$. Then
$\norm{(\one-x/\lambda)-\one}=\norm{x}/|\lambda|<\varepsilon$, so
$\one-x/\lambda\in\inv$, and therefore
$\lambda\one-x=\lambda(\one-x/\lambda)\in\inv$. Hence
$\spA(x)\subset\{\lambda:|\lambda|\le\norm{x}/\varepsilon\}$.
Theorem~\ref{thm:RS} therefore applies in either case.
\end{proof}

The GKZ theorem for normed Q-algebras is discussed by Roitman and
Sternfeld~\cite{RS}, and its extension to spectral algebras by
Palmer~\cite{Palmer}. A version for \emph{bounded} functionals on algebras
with a spectral norm is proved in~\cite[Section~3]{Badea} by passing to
the completion: the algebra is
a full subalgebra of its completion, and the upper semicontinuity of the
spectrum~\cite{AupetitLNM} transfers the relation $\varphi(x)\in\spA(x)$ to
all elements of the completion. In Corollary~\ref{cor:Q} boundedness of
$\varphi$ is not assumed, because Theorem~\ref{thm:RS} does not need it.

\subsection{Zero-free does not by itself mean exponential-linear}

The entire function $F(z)=e^{z^2}$ is zero-free but cannot be written
as $e^{\alpha z+\beta}$. It violates the
exponential-type bound, because $\log|F(r)|=r^2$ on the positive real
axis. Since $F'(0)=0$, it also fails the hypotheses of the P\'olya--Saxer theorem.
Likewise, the function in \eqref{eq:doubleexp} is zero-free with
zero-free first derivative, but its second derivative has zeros.
These examples explain the roles of the growth and derivative hypotheses
in the two analytic proofs.

\subsection{A scalar selection need not exist}

For $n\ge2$, the matrix algebra $M_n(\C)$ has no characters. Here is a
direct proof. Let $E_{ij}$ be the usual matrix units, and suppose that
$\chi$ were a character. Each $\chi(E_{ii})$ is idempotent in $\C$,
hence is $0$ or $1$, and their sum is $\chi(I)=1$. There is an index
$i$ with $\chi(E_{ii})=1$, and an index $j\ne i$ with
$\chi(E_{jj})=0$. But
\[
E_{ij}E_{ji}=E_{ii},\qquad E_{ji}E_{ij}=E_{jj},
\]
and multiplicativity would give both
\[
\chi(E_{ij})\chi(E_{ji})=1
\quad\text{and}\quad
\chi(E_{ji})\chi(E_{ij})=0,
\]
a contradiction. By GKZ, no normalised linear functional on $M_n(\C)$
can avoid zero on all invertible matrices. Thus, although each matrix has an eigenvalue, there is no linear
functional that selects an eigenvalue from the spectrum of every matrix.

\section{GKZ for Ransford's generalised spectra}
\label{sec:generalised}

Of the three proofs, the polynomial argument of Roitman and Sternfeld
offers the most elementary route to extensions involving Ransford's
generalised spectra. The key is to retain the two ingredients used in
Proof~II: suitable powers stay in the set on which the functional does
not vanish, and the polynomial roots remain bounded. This leads to the
extension announced in my 1991 note~\cite{Badea1991}. A proof is presented
here in the algebraic form given by Kulkarni and Sukumar~\cite{KulkarniSukumar2005}.
Throughout this section, $A$ is an arbitrary unital complex algebra with
$\one\ne0$; a norm will be introduced only where explicitly stated.

\subsection{Replacing the invertible elements by a distinguished set}

Let $\Omega\subset A$ satisfy the following three conditions:
\begin{enumerate}[label=\textup{(\roman*)}]
\item $0\notin\Omega$;
\item $\one\in\Omega$;
\item $\lambda\Omega\subset\Omega$ for every non-zero $\lambda\in\C$.
\end{enumerate}
The third condition actually gives $\lambda\Omega=\Omega$: apply it
also to $\lambda^{-1}$. In particular, $-\Omega=\Omega$. Define the
\emph{$\Omega$-spectrum} of $a\in A$ by
\begin{equation}\label{eq:omega-spectrum}
\sigma_\Omega(a)
:=\{\lambda\in\C:a-\lambda\one\notin\Omega\}.
\end{equation}
This abstract construction is due to Tom Ransford~\cite{Ransford1984}. The
definition and the first three conditions make sense even in a complex
vector space with a distinguished non-zero vector $\one$; multiplication
is not (yet) needed. They imply
\begin{equation}\label{eq:omega-scalars}
\sigma_\Omega(c\one)=\{c\}\qquad(c\in\C),
\end{equation}
because $(c-\lambda)\one$ belongs to $\Omega$ exactly when
$\lambda\ne c$.

Choosing $\Omega=\inv$ recovers the ordinary spectrum $\spA(a)$.
For a general $\Omega$, however, membership in $\Omega$ is a separate
condition from invertibility. We assume no inclusion between $\Omega$
and $\inv$, and shall not use the ordinary spectrum to locate the roots
in the argument below.

By imposing geometric conditions on $\Omega$, notably pseudoconvexity,
Ransford recovers a surprisingly large part of classical spectral theory
in the setting of a complex normed space. Much of the usual theory thus
survives without multiplication or completeness. As Ransford puts
it~\cite[p.~306]{Ransford1984} about 'Banach algebra theory without multiplication':
\begin{quote} 
`Still included in the ``grin after the Cheshire cat has disappeared''
are generalisations of many of the classical results of spectral theory,
including non-emptiness of the spectrum, the Beurling-Gelfand spectral
radius formula, Newburgh's continuity theorem, and even an analogue of
B.\ E.\ Johnson's uniqueness-of-norm theorem in which all mention of
representation theory has been eliminated.'
\end{quote}
In a related direction, Zalduendo~\cite{Zalduendo} found a geometric
condition equivalent to commutativity in Banach algebras.

When multiplication is available, we impose one further condition:
\begin{enumerate}[label=\textup{(\roman*)},start=4]
\item there is a strictly increasing sequence of positive integers
$n_1<n_2<\cdots$ such that
\begin{equation}\label{eq:omega-powers}
u\in\Omega\quad\Longrightarrow\quad u^{n_j}\in\Omega
\qquad(j\ge1).
\end{equation}
The same sequence is required to work for every $u\in\Omega$.
\end{enumerate}
Only powers along this unbounded sequence are required; closure of
$\Omega$ under arbitrary products is a stronger sufficient condition.

\begin{lemma}[A sufficient condition for the power property]
\label{lem:omega-powers}
If $\Omega$ is closed under multiplication, then $u^n\in\Omega$ for
every $u\in\Omega$ and every positive integer $n$. In particular,
\eqref{eq:omega-powers} holds with $n_j=j$.
\end{lemma}

\begin{proof}
Starting with $u\in\Omega$, induction gives
$u^{n+1}=u^nu\in\Omega$ whenever $u^n\in\Omega$.
\end{proof}

Thus the ordinary group of invertible elements satisfies the required
power property. The proof below only uses the implication in
\eqref{eq:omega-powers}; it does not require its converse or closure
under products of distinct elements.

In fact, the power condition is equivalent to the existence of a single
integer $m\ge2$ such that $u\in\Omega$ implies $u^m\in\Omega$ for every
$u$. One direction follows by choosing any $n_j\ge2$; for the other,
iteration gives $u^{m^j}\in\Omega$ for all $j\ge1$. The exponent $m$,
like the sequence in \eqref{eq:omega-powers}, must work for every $u$.

\subsection{Boundedness and non-emptiness}

Boundedness of the $\Omega$-spectra will be a separate hypothesis in the
algebraic theorem. A simple sufficient condition is available when $A$
has a norm, even without completeness or a submultiplicative norm.

\begin{proposition}[An elementary spectral bound]
\label{prop:omega-bound}
Assume the first three conditions on $\Omega$. If a norm on $A$ and a
number $\varepsilon>0$ satisfy
\[
\{u\in A:\norm{u-\one}<\varepsilon\}\subset\Omega,
\]
then
\begin{equation}\label{eq:omega-bound}
\sigma_\Omega(a)
\subset\{\lambda\in\C:|\lambda|\le\norm{a}/\varepsilon\}
\qquad(a\in A).
\end{equation}
If $\Omega$ is open, every $\Omega$-spectrum is also closed, and hence
compact, possibly empty.
\end{proposition}

\begin{proof}
For $|\lambda|>\norm{a}/\varepsilon$,
\[
\norm{\bigl(\one-a/\lambda\bigr)-\one}
=\frac{\norm{a}}{|\lambda|}<\varepsilon.
\]
Thus $\one-a/\lambda\in\Omega$. Scalar invariance gives
$a-\lambda\one=-\lambda(\one-a/\lambda)\in\Omega$, so
$\lambda\notin\sigma_\Omega(a)$. This proves the bound. If $\Omega$
is open, then $A\setminus\Omega$ is closed, and $\sigma_\Omega(a)$ is
its inverse image under the continuous map $\lambda\mapsto a-\lambda\one$;
hence $\sigma_\Omega(a)$ is closed. Moreover, openness and
$\one\in\Omega$ supply a suitable $\varepsilon>0$, so this closed set is
also bounded.
\end{proof}

Openness does \emph{not} ensure non-emptiness. For example, give
$A=\C^2$ coordinatewise multiplication, with $\one=(1,1)$, and put
$\Omega=A\setminus\{0\}$. This set is open and satisfies all four
conditions, since a non-zero coordinate remains non-zero under positive
integer powers. Nevertheless,
\[
(1,0)-\lambda(1,1)=(1-\lambda,-\lambda)\ne(0,0)
\qquad(\lambda\in\C),
\]
so $\sigma_\Omega((1,0))=\varnothing$. More generally, the
$\Omega$-spectrum of $(z_1,z_2)$ is empty if $z_1\ne z_2$ and is
$\{z_1\}$ otherwise. Thus all these spectra are bounded. This example
also shows that the power property need not imply multiplicative
closure: $(1,0),(0,1)\in\Omega$, but their product is zero.

\subsection{Non-vanishing and spectral selection}

\begin{lemma}[The $\Omega$-spectral reformulation]
\label{lem:omega-selection}
Assume the first three conditions on $\Omega$. For a complex-linear
functional $\varphi:A\to\C$ with $\varphi(\one)=1$, the following
conditions are equivalent:
\begin{enumerate}[label=\textup{(\roman*)}]
\item $\varphi(u)\ne0$ for every $u\in\Omega$;
\item $\varphi(a)\in\sigma_\Omega(a)$ for every $a\in A$.
\end{enumerate}
\end{lemma}

\begin{proof}
If (i) holds, then
$\varphi(a-\varphi(a)\one)=0$ forces
$a-\varphi(a)\one\notin\Omega$. By definition,
$\varphi(a)\in\sigma_\Omega(a)$.

Conversely, if $u\in\Omega$, then $0\notin\sigma_\Omega(u)$.
Condition (ii) places $\varphi(u)$ in this spectrum, so
$\varphi(u)\ne0$.
\end{proof}

In particular, the existence of such a functional forces every
$\Omega$-spectrum to be non-empty. It also rules out the preceding
example. Notice that spectral selection itself forces the normalisation:
applying it to $\one$ and using \eqref{eq:omega-scalars} gives
$\varphi(\one)=1$.

\subsection{The theorem and its polynomial proof}

\begin{theorem}[{\cite{Badea1991,KulkarniSukumar2005}}]
\label{thm:omega-GKZ}
Let $A$ be a unital complex algebra, and let $\Omega\subset A$ satisfy
the four conditions above. Suppose that $\sigma_\Omega(a)$ is bounded
for every $a\in A$. If a complex-linear functional $\varphi:A\to\C$
satisfies
\[
\varphi(\one)=1,
\qquad \varphi(u)\ne0\quad(u\in\Omega),
\]
then $\varphi$ is multiplicative. Equivalently, every complex-linear
selection of the $\Omega$-spectra is a character.
\end{theorem}

\begin{proof}
By Lemma~\ref{lem:jordan}, it suffices to prove $\varphi(x^2)=0$
for $x\in\ker\varphi$. Fix such an $x$ and use the polynomials $P_n$
from~\eqref{eq:pn}, with $n=n_j\ge2$. If
$\lambda\notin\sigma_\Omega(x)$, then
$\lambda\one-x\in\Omega$, since $-\Omega=\Omega$. The power
condition~\eqref{eq:omega-powers} gives
$(\lambda\one-x)^{n_j}\in\Omega$, and non-vanishing of $\varphi$ on
$\Omega$ yields $P_{n_j}(\lambda)\ne0$. Thus every root belongs to
$\sigma_\Omega(x)$, a bounded set.

Choose a finite bound $M_x$ for the moduli of these spectral values.
The coefficient identity~\eqref{eq:secondmoment} and the
estimate~\eqref{eq:RSestimate} from Proof~II now apply with $M_x$ in
place of $R_x$:
\[
|\varphi(x^2)|\le\frac{M_x^2}{n_j-1}
\longrightarrow0\qquad(j\to\infty).\qedhere
\]
\end{proof}

The Roitman--Sternfeld argument therefore extends by replacing ordinary
spectra with $\Omega$-spectra and using only the exponents $n_j$.
It requires no continuity, completeness or analytic growth theorem;
the role of the fundamental theorem of algebra is exactly as in
Section~\ref{sec:second}. For a fixed $\varphi$, boundedness of
$\sigma_\Omega(x)$ is needed only for $x\in\ker\varphi$.

\begin{corollary}[The normed case and automatic continuity]
\label{cor:omega-open}
Suppose $A$ has a norm and $\Omega$ satisfies the four conditions above.
If
\[
\{u:\norm{u-\one}<\varepsilon\}\subset\Omega
\qquad\text{for some }\varepsilon>0,
\]
then every normalised complex-linear functional that is non-zero on
$\Omega$ is multiplicative and continuous, with
\[
|\varphi(a)|\le\frac{\norm{a}}{\varepsilon},
\qquad \norm{\varphi}\le\frac1\varepsilon.
\]
\end{corollary}

\begin{proof}
Proposition~\ref{prop:omega-bound} supplies bounded spectra, so
Theorem~\ref{thm:omega-GKZ} gives multiplicativity. The same bound,
together with $\varphi(a)\in\sigma_\Omega(a)$ from
Lemma~\ref{lem:omega-selection}, gives the stated norm estimate. Thus
continuity could also be established before multiplicativity.
If $\Omega$ is open, a suitable $\varepsilon$ exists because
$\one\in\Omega$.
\end{proof}

\subsection{Why the converse need not hold}

For the ordinary spectrum, every character is a spectral selection.
For an arbitrary $\Omega$, this converse can fail. On $A=\C^2$ with
coordinatewise multiplication and $\one=(1,1)$, put
\[
\Omega=\{(z_1,z_2):z_1\ne0\}.
\]
This set is open, satisfies the first three conditions and is closed
under multiplication, so it satisfies the power condition as well.
Directly from the definition,
\[
\sigma_\Omega((z_1,z_2))=\{z_1\}.
\]
The second coordinate functional $\psi(z_1,z_2)=z_2$ is a character,
but $(1,0)\in\Omega$ and $\psi(1,0)=0$. Thus $\psi$ neither avoids
zero on $\Omega$ nor selects its spectra. The first coordinate
functional does both. In this example,
\[
\inv=\{(z_1,z_2):z_1z_2\ne0\}\subsetneq\Omega,
\qquad
\spA((z_1,z_2))=\{z_1,z_2\},
\]
which makes the distinction between the two spectra explicit.

\section{Extensions and further reading}
\label{sec:further}

The following references provide starting points for further reading.
\begin{itemize}
\item \emph{Other proofs.} In addition to Siddiqi's elementary analytic
proof~\cite{Siddiqi}, discussed earlier, Wille~\cite{Wille} gives two
proofs for commutative Banach algebras with a symmetric involution.
Palmer's book~\cite{Palmer} provides a broader account of the theorem.
\item \emph{Beyond Banach algebras.} For the extensions to normed
Q-algebras and spectral algebras illustrated by Corollary~\ref{cor:Q},
see~\cite{RS,Palmer} and~\cite[Section~3]{Badea}.
\item \emph{Maps between algebras.} Proposition~\ref{prop:preserver}
introduces the study of linear maps that preserve invertibility. Related
generalisations of GKZ were proved by Aupetit~\cite{AupetitPJM} and
Christensen~\cite{Christensen}. See also \cite{BresarSemrl2008}.
\item \emph{Weaker hypotheses.} The article~\cite{EoM} in the
\emph{Encyclopaedia of Mathematics} lists further directions, including
versions in which linearity of the functional is weakened and extensions to
maps between topological algebras.
\item  \emph{Recent work on reproducing kernel Hilbert spaces.} Chu, Hartz, Mashreghi and Ransford~\cite{CHMR} extend GKZ to a
reproducing kernel Hilbert space with a normalised complete
Pick kernel. More recently, Bhattacharyya, Jaikishan, Kumar and Sahu~\cite{BJKS}
give tests involving only reproducing kernels.
\end{itemize}
For further references on the Gleason--Kahane--\.{Z}elazko theorem,
we encourage the reader to consult MathSciNet~\cite{MathSciNetDistance}
or its freely accessible European counterpart, zbMATH Open.

\end{document}